\documentclass[leqno,11pt]{amsart}
\usepackage{amssymb,mathrsfs,color}
\usepackage{enumerate}
\usepackage{graphicx}
\usepackage{xcolor} 
\usepackage{geometry}
\usepackage{amsfonts,amssymb,amsmath}
\usepackage{slashed}
\usepackage{mathrsfs }
\usepackage{txfonts}
\usepackage[titletoc,title]{appendix}
\usepackage{wrapfig}

\usepackage{amsmath}

\usepackage{cancel}
\usepackage{cite}
\usepackage{bm}
\usepackage{nth}
\usepackage[T1]{fontenc}
\usepackage{hyperref}

\usepackage{marginnote}
\usepackage{cancel}
\usepackage{pgfplots}
\usepackage{upgreek}
\usepackage{slashed}

\allowdisplaybreaks[1]
\usepackage{enumerate}

\definecolor{green}{rgb}{0,0.8,0} 

\definecolor{deepgreen}{cmyk}{1,0,1,0.5}
\newcommand{\Blue}[1]{{\color{blue} #1}}

\newcommand{\Del}[1]{}

\numberwithin{equation}{section}

\newtheorem{theorem}{Theorem}[section]
\newtheorem{lemma}[theorem]{Lemma}
\newtheorem{proposition}[theorem]{Proposition}
\newtheorem{remark}[theorem]{Remark}
\newtheorem{definition}[theorem]{Definition}

\newcommand{\pa}{\partial}

\newcommand{\ep}{\varepsilon}

\renewcommand{\Im}{\mathrm{Im}}

\newcommand{\sgn}{{\mathrm{sgn}}}

\newcommand{\bbR}{\mathbb R}

\newcommand{\calA}{\mathcal A}

\newcommand{\calL}{\mathcal L}

\let\al=\alpha

\let\f=\frac

\let\om=\omega

\let\pa=\partial

\allowdisplaybreaks

\newcommand{\beq}{\begin{equation}}
\newcommand{\eeq}{\end{equation}}
\newcommand{\ben}{\begin{eqnarray}}
\newcommand{\een}{\end{eqnarray}}
\newcommand{\beno}{\begin{eqnarray*}}
\newcommand{\eeno}{\end{eqnarray*}}

\begin{document}
	\title[Asymptotic linear stability of Columnar vortices] {Asymptotic linear stability of Columnar vortices driven by Coriolis force}
	
	\author{Shuang Miao}
	
	\address{School of Mathematics and Statistics, Wuhan University, 430072, Wuhan, P. R. China}
	\email{shuang.m@whu.edu.cn}
	
	\author{Siqi Ren}
	\address{School of Mathematical Sciences, Zhejiang University of Technology,  310032, Hangzhou, P. R. China}
	\email{sirrenmath@zjut.edu.cn}

	\author{Zhifei Zhang}
	\address{School of Mathematical Sciences, Peking University, 100871, Beijing, P. R. China}
	\email{zfzhang@math.pku.edu.cn}
	
\begin{abstract}
In this paper, we establish the asymptotic linear stability of a class of Coriolis-driven columnar vortices for the 3-D axisymmetric Euler equations. This result represents a critical step toward proving the nonlinear asymptotic stability of such vortices. The key and widely applicable strategy is to construct a distorted Fourier basis, which is achieved by solving a two-parameter $(c, \xi)$-dependent Schr\"odinger equation associated with the linearized operator of the system.  To capture the precise asymptotic behavior of the solution, we decompose the $c-\xi$ plane into distinct regions, with the partitioning guided by the leading-order profiles of the Schr\"odinger equation across different parameter regimes.
\end{abstract}	

 \date{\today}

	\maketitle
	
	\tableofcontents

	\section{Introduction}
In this paper, we consider the 3D incompressible Euler equations in $\bbR^+\times \bbR^3$:
\begin{align}\label{eq:Euler}
	\partial_t{\textbf{v}}+{\textbf{v}}\cdot
	\nabla {\textbf{v}}+\nabla {p}=0,\quad  \nabla
	\cdot {\textbf{v}}=0,
\end{align}
where $\textbf{v}(t,x,y,z)$ is the fluid velocity and ${p}(t,x,y,z)$ is the pressure.	There exists a wealth of research findings on the local well-posedness and ill-posedness of the 3D Euler equations in various function spaces (see \cite{KP, BL1,BL2, BKM, CMZ, Chae, Che, Maj, MP, Vis}). For the 2D Euler equations, they are globally well-posed for smooth initial data (see \cite{Yud} for the vortex patch), and their long-time dynamical behavior constitutes a crucial current research field, with several significant advances achieved in recent years \cite{BM, Bou, WZZ1,WZZ2, IJ, MZ, LZ, EMS, KS, Den}.  Compared with the abundant research findings on the well-posedness of the 3D Euler equations, the mathematical theoretical research on their stability is relatively scarce (in contrast, a wealth of mature results have been accumulated for the 2D case \cite{MP}). Nevertheless, the stability problem of the 3D Euler equations is not only a core topic in the field of hydrodynamic stability but also of great significance for revealing the intrinsic mechanism of turbulence formation.

We are concerned with the stability of  the columnar vortices, which are steady solutions of \eqref{eq:Euler} and take the form of
\begin{align}
	\textbf{v}=U(r) \textbf{e}_{\theta},\quad \omega=\nabla\times \textbf{v}=W(r)\textbf{e}_z=\frac 1 r\frac {\partial{(rU)}} {\partial r}\textbf{e}_z,
\end{align}
where $\textbf{e}_{r}=\begin{pmatrix}\cos\theta\\
	\sin\theta\\0\end{pmatrix}$, $\textbf{e}_{\theta}=
\begin{pmatrix}-\sin\theta\\ \cos\theta\\0\end{pmatrix}$, $\textbf{e}_{z}=
\begin{pmatrix}0\\ 0\\1\end{pmatrix}$, $r=\sqrt{x^2+y^2}$ and $\theta=\arctan \frac{y}{x}$.
Columnar vortices represent a fundamental and ubiquitous form of organized rotational motion in fluids. In both natural environments and engineering applications ranging from aircraft wingtip vortices and atmospheric tornadoes to oceanic mesoscale eddies, these structures persist as dominant coherent features governing fluid dynamic behavior. Analyzing their stability mechanisms reveals the underlying principles of how such structures destabilize, break down, and ultimately transition to turbulence.	 Here are typical examples of columnar vortices that are often considered in the literature:
\begin{itemize}
	\item The Rankine vortex
	\begin{equation}\nonumber
		U(r)=\left\{\begin{array}{l}
			r\quad r\le 1,\\
			\frac 1 r\quad r\ge 1,
		\end{array}\right.\quad  W(r)=\left\{\begin{array}{l}
			2\quad r<1,\\
			0\quad r>1.
		\end{array}\right.\end{equation}
	
	\item The Kaufmann-Scully vortex
	\begin{equation}\nonumber
		U(r)=\frac r  {1+r^2},\quad W(r)=\frac 2 {(1+r^2)^2}.
	\end{equation}

	\item The Lamb-Oseen  vortex
	\begin{equation}\nonumber
		U(r)=\frac1 r(1-e^{-r^2}),\quad W(r)=2e^{-r^2}.
	\end{equation}
\end{itemize}
Kelvin was the first to investigate the stability of columnar vortices \cite{Kel}. In two recent works \cite{GS1, GS2}, Gallay and Smets systematically studied the linear stability of columnar vortices for general three-dimensional perturbations. See \cite{GS1} for a detailed historical account of columnar vortices.

\subsection{Main result}

This paper presents a theoretical study of a class of columnar vortices driven by the Coriolis force $\textbf{F}=-\textbf{e}_z\times \textbf{v}$, whose physical essence is captured by the geostrophic balance in a rotating system. In cylindrical coordinates, the specific form of the radial momentum equation characterizing this global balance is as follows:
\begin{equation}
	-\frac {U(r)^2} r+\partial_r P(r)=U(r).
\end{equation}
Although the global stability of 3D incompressible Euler equation is challenging, recently Guo, Pausader and Widmayer made a breakthrough by proving (see \cite{GPW}) the asymptotic stability of uniform rotation (namely, $U(r)=r$) for axisymmetric perturbations. This work demonstrates that uniform rotation can prevent the formation of solution singularities. See \cite{BMN, Du, EW, GR, KLT, RT, PW, Tak, Wan} and the references therein for more relevant works. In the absence of rotation, singularities can develop in finite time, see the recent breakthrough \cite{Elg, EGM}  for $C^{1,\al}$ data and \cite{Chenhou1, Chenhou2} for smooth data.

Motivated by \cite{GPW, Guo},  this paper aims to investigate the asymptotic stability of more general non-uniform rotation (namely,  a general $U(r)$). Let
\[ \textbf{u}_s=U(r) \textbf{e}_{\theta}, \quad p_s=P(r),\]
which is a steady solution of the Euler equations
\begin{align}\label{eq:Eule}
	\partial_t{\textbf{v}}+{\textbf{v}}\cdot
	\nabla {\textbf{v}}+\nabla {p}=0,\quad  \nabla
	\cdot {\textbf{v}}=0.
\end{align}
Let $\textbf{u}={\textbf{v}}-\textbf{u}_s$ and $\pi={p}-p_s$. Then we have
\begin{align*}
	\partial_t\textbf{u}+\textbf{u}\cdot\nabla \textbf{u}+\frac{U(r)}{r}\textbf{e}_z\times \textbf{u}+\frac{U(r)}{r}\partial_{\theta}\textbf{u}
	+r\left(\frac{U(r)}{r}\right)'u_r\textbf{e}_{\theta}+\nabla \pi=0,\quad \nabla\cdot \textbf{u}=0.
\end{align*}
For an axisymmetric  solution
$$\textbf{u}=u^r(t,r,z)\textbf{e}_r+
u^{\theta}(t,r,z)\textbf{e}_{\theta}+u^z(t,r,z)\textbf{e}_z,$$
the system can be rewritten as
\begin{align}\label{nonlinear-Euler-1}
	\partial_t\textbf{u}-2uu^{\theta}\textbf{e}_r+(2u+ru')u^r\textbf{e}_{\theta}
	+\textbf{u}\cdot\nabla \textbf{u}+\nabla \pi=0,
\end{align}
where $u(r)=\frac{U(r)}{r}$.  We denote
\begin{align}\label{steady-potential}
	\Omega(r)=2u+ru',\quad  V(r)=2u\Omega.
\end{align}
We then impose the following assumptions on the potential $V(r)$:\smallskip

\begin{itemize}
	\item [\textbf{(A1)}] $V'(r)>0$, $r\in\mathbb{R}^+$;
	
	\item [\textbf{(A2)}] $\mathrm{Ran} V=[V(0),1)\subset (0,1]$;
	
	\item [\textbf{(A3)}] $V(r)=1-a_0r^{-3}+O(r^{-4})$, $r\to+\infty$, \; $a_0>0$.
\end{itemize}
An example satisfying (A1)-(A3) is that $u(r)=\f12\left(1-\f13\langle r\rangle ^{-3}\right)$, where $\langle r\rangle=r+1$.

As a crucial step towards nonlinear asymptotic stability, the aim of this paper is to study the asymptotic linear stability of non-uniform rotational flows.  The linearized axisymmetric Euler-Coriolis system of \eqref{nonlinear-Euler-1} takes the form
\begin{align}\label{eq:vel-linearize-gernal0}
	\left\{\begin{array}{l}
		\pa_tu^r-2uu^{\theta}+\pa_rp=0,\\
		\pa_tu^{\theta}+(2u+ru')u^r=0,\\
		\pa_tu^{z}+\pa_z p=0,\\
		(\pa_r+r^{-1})u^r+\pa_{z}u^z=0.
	\end{array}\right.
\end{align}
In terms of the angular vorticity $\omega^\theta$ and the angular velocity $u^\theta$, the velocity $\textbf{u}$ is given by
\begin{align*}
	\textbf{u}=\left(\pa_z\Delta_1^{-1}\omega^{\theta}\right)\;\textbf{e}_r
	+u^{\theta}\;\textbf{e}_{\theta}
	+\left(-(\pa_r+r^{-1})\Delta_1^{-1}\omega^{\theta}\right)\;\textbf{e}_{z},
\end{align*}
where
\begin{align}\label{def Delta1}
	\Delta_1=\pa_r^2+r^{-1}\pa_r-r^{-2}+\pa_z^2.
\end{align}

Our main result is stated as follows.

\begin{theorem}\label{Thm:dipersive-decay}
	Let $m,n\in\mathbb{N}$ and $0<\delta\ll 1$ be a fixed number. Under the assumptions \textbf{(A1)}-\textbf{(A3)}, the solution of the linearized Euler-Coriolis system \eqref{eq:vel-linearize-gernal0} admits the following uniform decay bound for $t\gtrsim 1$,
	\begin{align}\label{t-1decay}
		\|\pa_z^n\pa_r^m(u^r, u^{\theta}, u^z)(t)\|_{L^{\infty}(\mathbb{R}^+\times\mathbb{R})}
		&\lesssim t^{-1}\left(\|\omega_0^{\theta}\|_{W^{2+m+n
				,1
			}_{r^{-\delta}+r^{\delta}}(\mathbb{R}^+\times\mathbb{R})}
		+\|u_0^{\theta}\|_{W^{3+m+n
				,1
			}_{r^{-\delta}+r^{\delta}}(\mathbb{R}^+\times\mathbb{R})}\right).
	\end{align}
	Here the norm of the weighted Sobolev space $W^{k,1}_{\omega}$ is defined by
	\begin{align*}
		\|F\|_{W^{k,1}_{\omega}(\mathbb{R}^+\times\mathbb{R})}
		=\sum_{0\leq i+j\leq k}\int_{\mathbb{R}^+\times\mathbb{R}}|\omega(r)\pa^{i}_r\pa^{j}_zF(r,z)|rdrdz.
	\end{align*}	
\end{theorem}

\begin{remark}
	The decay rate obtained is expected to be optimal. Even for the case of uniform rotation ($u\equiv \f12$),  the second author \cite{Ren} established the same decay rate via direct calculations based on the Fourier transform:
	\begin{align*}
		\|
		\pa_z^n\pa_r^m(u_r,u_{\theta},u_z)(t)\|_{L^{\infty}(\mathbb{R}^+\times\mathbb{R})}
		&\leq Ct^{-1}\|\textbf{u}_0\|_{W^{3+m+n
				,1}(\mathbb{R}^+\times\mathbb{R})}.
	\end{align*}
	The core objective of this paper is to extend the aforementioned concept to the scenario of non-uniform rotation; however, this extension has presented numerous intricate analytical challenges (see Section 1.3).	
	
\end{remark}

\subsection{Reformulation of the linearized system}
In view of the last equation in \eqref{eq:vel-linearize-gernal0}, we introduce the stream function $\psi$ such that
\begin{align*}
	\pa_r\psi=-ru^z,\qquad \pa_z\psi=ru^r,
\end{align*}
which gives
\begin{align*}
	\Delta_1(r^{-1}\psi)=\omega^{\theta}.
\end{align*}
We introduce new unknowns
\begin{align*}
	h=V\Delta_1^{-1}(\omega^{\theta})=Vr^{-1}\psi,\quad g=2uu^{\theta}.
\end{align*}
In terms of $\om^\theta$ and $g$, the velocity can be recovered by the following  relationship
\begin{align}
	\begin{split}\label{eq:vel-linearize-gernal} u^r(t,r,z)&=\pa_z\Delta_1^{-1}\omega^{\theta}=\pa_z(V^{-1}h),\\
		u^z(t,r,z)&=-(\pa_r+r^{-1})\Delta_1^{-1}\omega^{\theta}=
		-(\pa_r+r^{-1})(V^{-1}h),\\
		u^{\theta}(t,r,z)&=(2u)^{-1}g.	
	\end{split}
\end{align}
It follows from  \eqref{eq:vel-linearize-gernal0} that
\begin{align*}
	\left\{
	\begin{array}{l}
		\notag\pa_th-\pa_zV\Delta_1^{-1}(g)=0,\\
		\pa_tg+\pa_zh=0,\\
		h|_{t=0}=V\Delta_1^{-1}(\omega^{\theta}_0) ,\quad g|_{t=0}=2uu^{\theta}_{0}.\\
	\end{array}\right.
\end{align*}
Taking the Fourier transform $z\to k$ for $k\neq 0$, the system above is reduced to
\begin{align}\label{eq:vel-linearize-gernal'}
	\left\{
	\begin{array}{l}
		\pa_t\hat{h}+ik\mathcal{A}\hat{g}=0,\\
		\pa_t\hat{g}+ik\hat{h}=0,\\
		\hat{h}|_{t=0}=V\Delta_{1,k}^{-1}(\hat{\omega}^{\theta}_0) ,\quad \hat{g}|_{t=0}=2u\hat{u}^{\theta}_{0},\end{array}\right.
\end{align}
where
\beno
\mathcal{A}=-V\Delta_{1,k}^{-1},\quad  \Delta_{1,k}=\pa_r^2+r^{-1}\pa_r-r^{-2}-k^2.
\eeno

As we shall see, the operator $\calA$ is self-adjoint in $L^{2}(\bbR_{+}; V(r)^{-1}rdr)$ and it only has continuous spectrum which is contained in $[0,\frac{1}{k^{2}}]$. In fact the spectrum fills the entire interval $[0,k^{-2}]$. The system \eqref{eq:vel-linearize-gernal'} is a dispersive-type system and its solution can be formally written as
\begin{align}\begin{split}
		\hat{h}(t,r,k)&=\cos (k\sqrt{\mathcal{A}} t)\left[V\Delta_1^{-1}(\hat{\omega}^{\theta}_0) \right]-i\sqrt{\mathcal{A}}\sin (k\sqrt{\mathcal{A}} t)
		\left[2u\hat{u}^{\theta}_{0}\right],\\
		\label{fm:hg}
		\hat{g}(t,r,k)&=\cos (k\sqrt{\mathcal{A}} t) \left[2u\hat{u}^{\theta}_{0}\right]-i\f{\sin (k\sqrt{\mathcal{A}} t)}{\sqrt{\mathcal{A}}}\left[V\Delta_1^{-1}(\hat{\omega}^{\theta}_0) \right].\end{split}
\end{align}
The pointwise decay of $(h,g)$ relies on the spectral property of $\calA$, in particular its distorted Fourier basis. Let $c\in(0,1)$ be the spectrum of $\calA$. The distorted Fourier basis $\phi(r,c,k)$ associated to $\calA$ satisfies the following Schr\"odinger equation with two parameters $(c,k)$
\begin{align}\label{main ODE}
	\calA\phi=\frac{c}{k^{2}}\phi\quad \Longleftrightarrow\quad \phi^{\prime\prime}(r)+r^{-1}\phi^{\prime}(r)-r^{-2}\phi(r)-k^{2}\phi(r)+\frac{k^{2}V(r)}{c}\phi(r)=0.
\end{align}
To see the dispersive nature of \eqref{eq:vel-linearize-gernal'} and \eqref{fm:hg}, for instance we look at a typical term in the explicit expression (see Lemma \ref{lem:distortedF}) for the solution to \eqref{eq:vel-linearize-gernal0}:
\begin{align}\label{kernel formal}
	\begin{split}
		&\mathcal{F}^{-1}_{k\mapsto z}\left(kV(r)^{-1}\calA^{2}f(\calA)\hat{u}(s,k)\right)\\
		&=\frac{1}{2\pi^{2}}\int_{\mathbb{R}\times \mathbb{R}^{+}}\left(\int_{0}^{1}f\left(\frac{c}{k^{2}}\right)\left(\int_{\mathbb{R}}\frac{\phi(r,c,k)\phi(s,c,k)}{|W(c,k)|^{2}}\frac{e^{ik(z-y)}}{k}dk\right)dc\right)\cdot u(s,y)sdsdy,
	\end{split}
\end{align}
where $W$ is the Wronskian between $\phi$ and the oscillatory solution to \eqref{main ODE} at $r=+\infty$, and $f(x)=\sqrt{x}\left(1-\chi(\delta\sqrt{x})\right)\frac{\cos(k\sqrt{x}t)}{k\sqrt{x}}$ with $\chi(\cdot)$ being a cutoff function defined in \eqref{def:chi-notation}.
Using the relation above and the formula \eqref{eq:vel-linearize-gernal'}, one sees that a main contribution to $u^{z}$, for instance, is given by
\begin{align}\label{vz part}
	I^{L}=\iint_{\bbR^{+}\times \bbR}\int_{0}^{1}\frac{\cos(\sqrt{c}t)}{\sqrt{c}}K_{1,0}^{L}(c,r,s,z-y)\hat{v}_{0}sdcdsdy,
\end{align}
where
\begin{align}\label{kernel example}
	K^{L}_{1,0}(c,r,s,z)=\text{p.v.}\int_{\bbR}\left(1-\chi\left(\frac{\delta\sqrt{c}}{k}\right)\right)\frac{\frac{\sqrt{c}(\partial_{r}+r^{-1})}{k}\phi(r,c,k)\phi(s,c,k)}{|W(c,k)|^{2}}\frac{e^{ikz}}{k}dk.
\end{align}
Here $\delta>0$ is a small constant, and $\chi(\cdot)$ is a smooth cut-off function on $[0,\infty)$ supported near $0$. Then a dispersive estimate on \eqref{vz part} would follow from the fact
\begin{align*}
	\cos\left(\sqrt{c}t\right)=\partial_{\sqrt{c}}\left(\frac{\sin(\sqrt{c}t)}{t}\right),
\end{align*}
as well as an appropriate estimate on
\begin{align}\label{kernel deri example}
	\int_{0}^{1}\left|\partial_{\sqrt{c}}K_{1,0}^{L}(r,s,z,c)\right|dc.
\end{align}
The argument above relies on the global behavior of $\phi(r,c,k)$. In fact, the main ingredients of this paper are devoted to the meticulous construction of $\phi(r,c,k)$ and the analysis of its behavior across various regimes of the parameter pair $(c,k)$.

\subsection{Distorted Fourier transform}
\subsubsection{Overview}
The distorted Fourier transform for self-adjoint Schr\"odinger operators can be derived via Stone's formula. Gesztesy-Zinchenko \cite{GZ} adopt this approach for operators on the half-line, which entails a detailed justification of passing to the limit in the resolvents from the upper and lower half-planes onto the spectrum. For $L^{1}_{\text{loc}}$ potentials, one can leverage the Herglotz property of the Weyl $m$-function. In that case, $\Im\, m(\lambda+\sqrt{-1}\epsilon)$ converges in the weak-star sense to the Herglotz measure, which then serves as the spectral measure. In contrast, for strongly singular potentials, such as the inverse square potentials arising in the present work, the $m$-function is no longer Herglotz and \cite{GZ} carefully deduce the existence of the limit of the resolvents onto the spectrum directly via the spectral theorem for self-adjoint operators.

The space-time resonances method based on the distorted Fourier transform for $\calL_{c,k}$ emerges as a promising approach to address the full asymptotic stability problem for the Euler-Coriolis equation with non-uniform rotation. In recent years, this method has established itself as a powerful tool in the study of the asymptotic stability of solitons when dispersion is not strong. The application of the distorted Fourier transform in this context was pioneered in \cite{KS3}, see \cite{KS2} for an overview of this technique. We refer to the review article \cite{Ger} and the works \cite{Chen, CL, GP, KP2, DM, GPZ, LL, LS} for a sample of recent advances. It is also worth noting that \cite{CL, CG, CGP, KS3, LS1, Li, LL, LSS} involve the development of the spectral and distorted Fourier theory for non self-adjoint matrix operators.\smallskip

\begin{figure}[!h]\label{pic}
	\centering
	\includegraphics[width=0.6\textwidth]{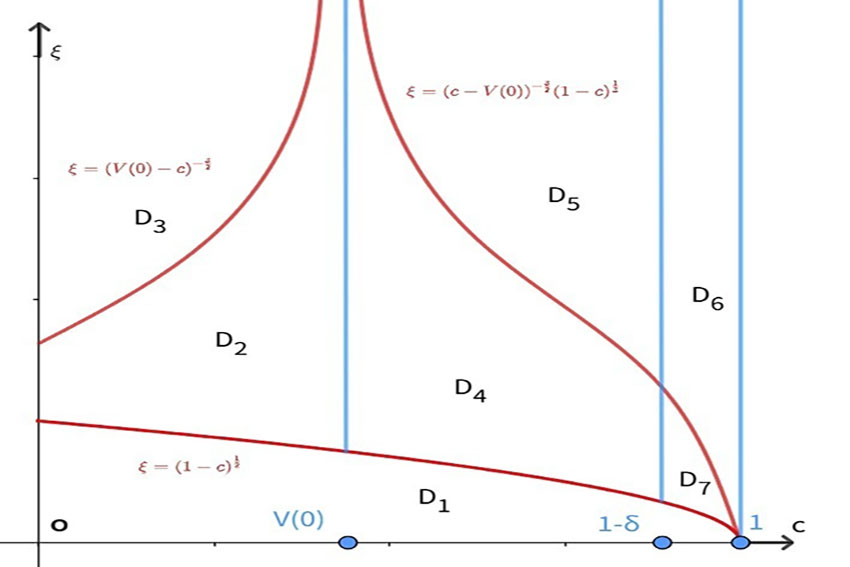}
	\caption{ }
\label{fig2}
\end{figure}

The construction of the distorted Fourier basis associated with $\calA$, which constitutes a major part of this paper, builds on techniques and insights from \cite{BP, CDST, CSST, GZ, KM, KMS, KST}. In particular, to construct $\phi(r,c,k)$ and analyze its global behavior, we carefully analyze the relation between the two parameters of the Schr\"odinger operator 
\begin{align}\label{def: main oper}
	\calL_{c,k}:=-\partial_{r}^{2}-r^{-1}\partial_{r}+r^{-2}+k^{2}\left(1-\frac{V(r)}{c}\right).
\end{align}
But in contrast to these works-particularly \cite{KMS}, we need to decompose the $c-\xi$ plane (where $\xi:=k\sqrt{c^{-1}-1}$) into distinct regions (see Figure \ref{fig2}), with the partitioning guided by the leading-order profiles of \eqref{main ODE} in different regions.

While in the case of uniform rotation considered in \cite{GPW}, the Fourier basis $\phi$ can be explicitly written out using Bessel functions, the construction of distorted Fourier basis is much more sophisticated in non-uniform rotation case.
\subsubsection{Construction of distorted Fourier basis}
We now outline the {\bf main ideas} behind the construction of the distorted Fourier basis.

Starting from $r=0$, we determine the maximal interval $(0, r_0(\xi,c)]$, on which the precise behavior of the real-valued  $\phi(r,\xi,c)$ (defined below) is manifested
\begin{align}  \label{eigen-equation1}
\phi''+r^{-1}\phi'-r^{-2}\phi+\f{\xi^2(V-c)}{1-c}\phi=0,\quad\phi(r,\xi,c)\sim \xi r\quad (r\to 0).
\end{align}
Starting from $r=+\infty$, we determine the maximal interval $[r_{\infty} (\xi,c),+\infty)$, on which the precise behavior of $f_+(r,\xi,c)$(defined below) is manifested
\begin{align}  \label{eigen-equation2}
f''_++r^{-1}f'_+-r^{-2}f_++\f{\xi^2(V-c)}{1-c}f_+=0,\quad f_+(r,\xi,c)\sim (\xi r)^{-\f12}e^{i\xi r}\quad (r\to+\infty).
\end{align}
The lifespan $r_{0} (\xi,c)$ and $r_{\infty} (\xi,c)$ vary with different choices of the parameters, but we expect them to coincide in some sense. As a consequence, the Wronskian
\begin{align}\label{def: Wronskian}
W(\xi,c)=\mathcal{W}(\phi,f_+)(\xi,c)=r\big(\phi'(r,\xi,c)f_+(r,\xi,c)-\phi (r,\xi,c)f_+'(r,\xi,c)\big)
\end{align}
can be evaluated at $r=r_{0} (\xi,c)\simeq r_{\infty} (\xi,c)$ for any $\xi>0$ and $c\in(0,1)$. 
Then for the parameters $(\xi,c)\in(0,+\infty)\times (0,1)$, we obtain the precise behavior for $\widetilde{\phi}(r,\xi,c)=\frac {\phi(r,\xi,c)} {|W(\xi,c)|}$ for all $r>0$ in the following way:
\begin{itemize}

\item for $r\le r_{0} (\xi,c)$, we use the information of $\phi(r,\xi,c)$ by \eqref{eigen-equation1} on $(0, r_0(\xi,c)]$;

\item for $r\gtrsim  r_{\infty} (\xi,c)$, we use the formula $$\widetilde{\phi}(r,\xi,c)=
-\mathrm{Re}\left(\f{W(\xi,c)f_+(r,\xi,c)}
{\left|W(\xi,c)\right|}\right),
$$
and the information of $f_+(r,\xi,c)$  by \eqref{eigen-equation2} for $r\ge r_{\infty}(\xi,c)$.
\end{itemize}
We now briefly discuss how to solve \eqref{eigen-equation1} and \eqref{eigen-equation2} over the maximal interval where all types of precise behaviors are manifested. Indeed, after carrying out extensive calculations, the parameter space $(0,+\infty)\times (0,1)$ (corresponding to $(\xi,c)$) is delicately partitioned into seven distinct regions $D_1-D_7$. The subsequent formal analysis in these regions is based on the assumptions on $V(r)$ stated in \textbf{(A1)}-\textbf{(A3)}.
\begin{enumerate}
\item  In the low-frequency region
$$D_1=\big\{c\in(0,1),0<\xi\lesssim M^2(1-c)^{\f12},\;M\gg 1\;\text{is}\;\text{fixed}\big\},$$
we decompose $r$ into two parts $r\in(0,r_{0}(\xi,c)]$ and $r\in [r_{\infty}(\xi,c),\infty)$.  For the first part (see Lemma \ref{lem:phi-1*} and Lemma \ref{lem:phi-1.5}), we solve \eqref{main ODE} on $r\in (0,\delta \xi^{-1}]$ with $\delta\ll 1$ fixed via the following reformulation
\begin{align}  \label{eigen-equation11}
	\phi''+r^{-1}\phi'-r^{-2}\phi+\xi^2\phi
	=\f{\xi^2(1-V)}{1-c}\phi,\quad\phi(r,\xi,c)\sim \xi r\quad (r\to 0).
\end{align}
Here, the left-hand side corresponds to the main part whose homogeneous solution is exactly the Bessel function of the first kind $J_1(\xi r)$; and the right-hand side $\f{\xi^2(1-V)}{1-c}\phi\simeq \f{\xi^2\langle r\rangle^{-3}}{1-c}\phi$ is a remainder term, provided that $\f{\xi^2\langle r\rangle^{-3}}{1-c}\lesssim r^{-2}$. Indeed, by $\xi\lesssim (1-c)^{\f12}$, if $r\lesssim 1$, then  $\f{\xi^2}{1-c}\lesssim 1\lesssim r^{-2}$; if $1\lesssim r\lesssim \xi^{-1}$, then  $\f{\xi^2r^{-3}}{1-c}\lesssim \f{\xi^2r^{-2}}{1-c} \lesssim r^{-2}$. For the second part (Lemma \ref{lem:k<1-f+}), we solve \eqref{eigen-equation2} on $r\in [\xi^{-1},+\infty)$ via the same formulation
\begin{align}  \label{eigen-equation21}
	f_+''+r^{-1}f_+'-r^{-2}f_++\xi^2f_+=\f{\xi^2(1-V)}{1-c}f_+,\quad f_+(r,\xi,c)\sim (\xi r)^{-\f12} e^{i\xi r}\quad (r\to +\infty),
\end{align}
where the left-hand side corresponds to the main part whose homogeneous solution is exactly the Hankel's function of the first kind $H_1(\xi r)=J_1(\xi r)+iY_1(\xi r)$; and formally, the right-hand side $\f{\xi^2(1-V)}{1-c}f_+$ is a perturbation term provided that $\f{\xi^2(1-V)}{1-c}\sim \f{\xi^2r^{-3}}{1-c}\lesssim \xi^2$, i.e., $r\gtrsim (1-c)^{-\f13}$, which satisfies $r\gtrsim \xi^{-1}$ by $\xi\lesssim (1-c)^{\f12}\leq (1-c)^{\f13}$.
\vspace{1.5mm}

\item  In the small-spectrum and medium-frequency region
$$D_2=\big\{c\in(0,V(0)),1\lesssim \xi\lesssim (V(0)-c)^{-\f32}\big\},$$
we again decompose $r$ into two parts. For the first part (see Lemma \ref{lem: phi-2}), we solve \eqref{eigen-equation1} on $r\in (0,\xi^{-\f23}]$ via the following formulation:
\begin{align}  \label{eigen-equation12}
	\phi''+r^{-1}\phi'-r^{-2}\phi =\f{\xi^2(c-V)}{1-c}\phi,\quad \phi(r,\xi,c)\sim \xi r\quad  (r\to 0),
\end{align}
where the left-hand side corresponds to the main part whose homogeneous solution is exactly the linear function $\xi r$; and the right-hand side $\f{\xi^2(c-V)}{1-c}\phi$ is viewed as a perturbation,  since $1\lesssim\xi\lesssim (V(0)-c)^{-\f23}$ and $r\lesssim \xi^{-\f23}$ imply $$\f{\xi^2|c-V(r)|}{1-c}\lesssim \xi^2\left(|V(0)-c|+|V(r)-V(0)|\right)\lesssim \xi^2\left(|V(0)-c|+r\right)\lesssim \xi^{\f43}\lesssim r^{-2}.$$ For the second part (see Lemma \ref{lem:f+2}), we solve \eqref{eigen-equation2} on $r\in [\xi^{-\f23},+\infty)$ via the  formulation
\begin{align}  \label{eigen-equation22}
	\begin{split}
	&f_+''+r^{-1}f_+' +\f{\xi^2(V-c)}{1-c}f_+=r^{-2}f_+,\\ &f_+(r,\xi,c)\sim \left(\xi \int_0^r\sqrt{\f{V(s)-c}{1-c}}ds\right)^{-\f12}e^{i\xi \int_0^r\sqrt{\f{V(s)-c}{1-c}}ds}\quad  (r\to +\infty),
	\end{split}
\end{align}
where, noticing $\lim_{r\to+\infty}\f{ \int_0^r\sqrt{\f{V(s)-c}{1-c}}ds}{r}=1$,  the $e^{i\xi r}$ oscillation at infinity is replaced by  $e^{i\xi \int_0^r\sqrt{\f{V(s)-c}{1-c}}ds}$, which is the far-field behavior of  the oscillatory Airy function as follows
\begin{align} \label{def:O-Airy-intro}
	\mathrm{Ai}(-\xi^{\f23}\tau)-i\mathrm{Bi}(-\xi^{\f23}\tau)\quad\text{with}\;\;\tau(r,c)
	=\left(\f32\int_0^r\sqrt{\f{V(s)-c}{1-c}}ds\right)^{\f23}.
\end{align}
The behavior \eqref{def:O-Airy-intro} dominates the left-hand side, which can be seen upon applying the Langer transform (to flatten the medium potential term $\f{\xi^2(v-c)}{1-c}f_+$):
\begin{align} \label{langer-intro}
	\begin{split}
		\mathrm{f}_+(\tau,\xi,c)
		&:=C\xi^{\f13}q^{\f14}r^{\f12}f_+(r,\xi,c),\quad q(r,c)=\f{Q}{\tau},\quad Q(r,c)=\f{V(r)-c}{1-c}.
	\end{split}
\end{align}
The right-hand side is a perturbation, which can be verified as follows under the condition $(c,\xi)\in D_2$, $r\gtrsim \xi^{-\f23}$ such that  $r^{-2}\lesssim \f{\xi^2(V-c)}{1-c}$. Indeed,  if $\xi^{-\f23}\lesssim r\lesssim 1$, then $r^{-2}\lesssim \xi^2 r\simeq \f{\xi^2(V(r)-V(0))}{1-c} \leq \f{\xi^2\left((V(r)-c)\right)}{1-c}$; if $r\gtrsim 1$, then $r^{-2}\lesssim 1\lesssim \xi^2 \sim \f{\xi^2(V(r)-V(0))}{1-c} \leq \f{\xi^2(V(r)-c)}{1-c}$.	
\vspace{1.5mm}

\item  In the small-spectrum and high-frequency region
$$D_3=\big\{c\in(0,V(0)),\xi\gtrsim (V(0)-c)^{-\f32}\big\},$$
we decompose $r$ into two parts. For the first part (see Lemma \ref{lem:phi-3} and Lemma  \ref{lem:phi-4.5}), we solve \eqref{eigen-equation1} on $r\in (0,\xi^{-1}(V(0)-c)^{-\f12}]$ via the formulation \eqref{eigen-equation12}, where the left-hand side corresponds to the main part whose homogeneous solution is $\xi r$. The right-hand side of \eqref{eigen-equation12} is viewed as a perturbation, since $$\f{\xi^2|c-V(r)|}{1-c}\lesssim \xi^2\left(|V(0)-c|+r\right)\lesssim \xi^{2}|V(0)-c|\lesssim r^{-2}.$$
For the second part (see Lemma  \ref{lem:f+2}), we solve  \eqref{eigen-equation2} on $r\in [\xi^{-1}(V(0)-c)^{-\f12},+\infty)$ via the formulation \eqref{eigen-equation22}.
The oscillatory Airy behavior \eqref{def:O-Airy-intro} dominates on the left-hand side, upon applying the Langer transform \eqref{langer-intro}, while the right-hand side of \eqref{eigen-equation22} is a perturbation compared to the term $r^{-2}f_+$. Indeed,  if $r\gtrsim \xi^{-1}(V(0)-c)^{-\f12}$, then $r^{-2}\lesssim \xi^2 (V(0)-c) \leq \f{\xi^2(V(r)-c)}{1-c}$.	
\vspace{1.5mm}

\item  In the  medium-spectrum and  medium-frequency region
$$D_4=\big\{c\in(V(0),1-\delta),1\lesssim \xi\lesssim  (c-V(0))^{-\f32},\;\delta\ll 1\;\text{is}\;\text{fixed}\big\},$$
we decompose $r$ into two parts. For the first part (see Lemma\ref{lem: phi-2}), we solve \eqref{eigen-equation1} on $r\in (0,\xi^{-\f23}]$ with the  formulation \eqref{eigen-equation12}.
The main part is again the left-hand side whose homogeneous solution is $\xi r$; while the right-hand side is viewed as a perturbation compared with the term $r^{-2}\phi$. The reason is as follows,
$$\f{\xi^2|c-V(r)|}{1-c}\lesssim \xi^2\left(|V(0)-c|+|V(r)-V(0)|\right)\lesssim \xi^2\left(|V(0)-c|+r\right)\lesssim \xi^{\f43}\lesssim r^{-2}.$$ For the second part (see Lemma  \ref{lem:f+4}), we solve  \eqref{eigen-equation2} on $r\in [\xi^{-\f23},+\infty)$ with the formulation:
\begin{align}  \label{eigen-equation23}
	\begin{split}
	&f_+''+r^{-1}f_+' +\f{\xi^2(V-c)}{1-c}f_+=r^{-2}f_+,\\ &f_+(r,\xi,c)\sim \left(\xi \int_{r_c}^r\sqrt{\Big|\f{V(s)-c}{1-c}\Big|}ds\right)^{-\f12}e^{i\xi \int_{r_c}^r\sqrt{\left|\f{V(s)-c}{1-c}\right|}ds}\;(r\to +\infty).
	\end{split}
\end{align}
The oscillatory Airy behavior
\begin{align} \label{def:O-Airy-intro2}
	\mathrm{Ai}(-\xi^{\f23}\tau)-i\mathrm{Bi}(-\xi^{\f23}\tau)\quad\text{with}\;\tau(r,c)
	=\mathrm{sgn}(r-r_c)\left(\f32\int_{r_c}^r\sqrt{\left|\f{V(s)-c}{1-c}\right|}ds\right)
	^{\f23},\;r_c=V^{-1}(c), \end{align}
dominates the left-hand side, which can be seen upon applying the Langer  transform \eqref{langer-intro}, while the right-hand side of \eqref{eigen-equation22} is a perturbation compared to the term $r^{-2}f_+$. Indeed, if $\xi^{-\f23}\lesssim r\lesssim 1$, then $r^{-2}\lesssim \xi^2 r\sim \f{\xi^2(V(r)-V(0))}{1-c} \leq \f{\xi^2\left((V(r)-c)\right)}{1-c}$; if $r\gtrsim 1$, then
$r^{-2}\lesssim 1\lesssim \xi^2 \sim \f{\xi^2(V(r)-V(0))}{1-c} \leq \f{\xi^2(V(r)-c)}{1-c}$.	
\vspace{1.5mm}

\item  The following two cases can be addressed jointly: the medium-spectrum and high-frequency region
$$D_5=\big\{c\in(V(0),1-\delta),\xi\gtrsim M(c-V(0))^{-\f32},\;\delta\ll 1,\;M\gg 1\;\text{are}\;\text{fixed}\big\},$$
and the large-spectrum and high-frequency region
$$D_6=\big\{c\in(1-\delta,1),\xi\gtrsim M(1-c)^{\f13},\;\delta\ll 1,\;M\gg 1\;\text{are}\;\text{fixed}\big\}.$$
Observing the behavior $c-V(0)\sim \f{r_c}{\langle r_c\rangle}$, $(1-c)^{\f13}\simeq \langle r_c\rangle^{-1}$, the parameters in these cases can be interpreted as
\begin{align}\label{big-xi!}
	r_c>0,\quad \xi r_c^{\f32}/\langle r_c\rangle^{\f12}\gtrsim M,\;\;\;M\gg 1\;\text{fixed}.
\end{align}
In each of these cases, we identify four intervals on the $r$-axis, as we must apply the Langer transform to flatten the potential term  $\f{\xi^2(V(r)-c)}{1-c}f_+$ in this high-frequency case, which makes the coefficient $\f{\xi^2(V(r)-c)}{1-c}$ appear in the denominator. In contrast to the other high-frequency region $D_3$, when $c\in(V(0),1)$, the potential has a unique zero $r_c=V^{-1}(c)$. This gives rise to a size competition between $r-r_c$ and $\xi$, leading to further partitions of the $r$-axis where the solution exhibits distinct behaviors.

For the first part (see Lemma \ref{lem:phi-3}),  we solve  \eqref{eigen-equation1} on $r\in (0,\xi^{-1}(c-V(0))^{-\f12}]$  via the  formulation
\eqref{eigen-equation12}. The main part is the left-hand side with homogeneous solution $\xi r$, while the right-hand side is viewed as a perturbation compared with the term $r^{-2}\phi$. Indeed, we denote $Q(r,c)=\f{V(r)-c}{1-c}$, $r_c=V^{-1}(c)$, $r_*=M^{\f12}\xi^{-1}\left(\f{1-c}{c-V(0)}\right)^{\f12}$ ($M\gg1$ fixed),  and it holds that
\begin{align}\label{relation:xi2Q-r-2}
	\begin{split}
		&\xi^2 |Q(r,c)|\lesssim r^{-2},\quad\quad \text{if}\;\;r\lesssim r_*\big(\sim \xi^{-1}r_c^{-\f12}\langle r_c\rangle^{-1}\ll \min\{1,r_c\}\big),\\
		& r^{-2}\lesssim \xi^2 |Q(r,c)|,\quad\quad \text{if}\;\;r_*\lesssim r\leq r_c-\xi^{-\f23}\langle r_c\rangle^{\f13},\\
		&r^{-2}\lesssim \xi^2 |Q(r,c)|,\quad\quad \text{if}\;\; r\geq r_c+\xi^{-\f23}\langle r_c\rangle^{\f13},\\
		&\xi^2 |Q(r,c)|\lesssim r^{-2},\quad\quad \text{if}\;\;  r_{c}-\xi^{-\f23}\langle r_c\rangle^{\f13}\leq r\leq r_c+\xi^{-\f23}\langle r_c\rangle^{\f13}.
	\end{split}
\end{align}
Here we used the following behaviors in Lemma \ref{lem:behave-Q}:
\begin{align*}
	\begin{split}
		&Q(r,c)\sim \f{(r-r_c)\langle r_c\rangle^2}{\langle r\rangle ^3},\;r\lesssim r_c;\quad\quad Q(r,c)\sim \f{r-r_c}{\langle r\rangle},\;r\gtrsim r_c. 																				\end{split}
\end{align*}
For the second part, due to the second line in \eqref{relation:xi2Q-r-2}, we solve $\phi$ on $r\in[r_*, r_c-C\xi^{-\f23}\langle r_c\rangle^{\f13}]$ via the following formulation of the initial value (given by the first part) problem:
\begin{align}  \label{eigen-equation13}
	\phi''+r^{-1}\phi'+\f{\xi^2(V-c)}{1-c}\phi=r^{-2}\phi,\quad \phi(r)|_{r_*}=\phi(r_*),\quad \phi'(r)|_{r_*}=\phi'(r_*). 		
\end{align}
Upon applying the Langer transform, one finds that the dominant behavior on the left-hand side corresponds to the Airy function
$$C_{61}\xi^{-\f13}q^{-\f14}r^{-\f12}\mathrm{Ai}(-\xi^{\f23}\tau)
\left(1+\widetilde{C}_{61}\int_{\tau(r_*,c)}
^{\tau(r,c)}\textrm{Ai}^{-2}(-\xi^{\f23}t')dt'\right),$$ where $\tau(r,c)
=\mathrm{sgn}(r-r_c)\left(\f32\int_{r_c}^r\sqrt{\sgn(s-r_{c})\f{V(s)-c}{1-c}}ds\right)
^{\f23}$, $C_{61}(\xi,c)$ and $\widetilde{C}_{61}(\xi,c)$ are constants fitting the initial datum. See Lemma \ref{lem:f+3}  for the
details. In the third part, we solve $f_+$ from infinity for $r\in[r_c+C\xi^{-\f23}\langle r_c\rangle^{\f13},+\infty)$ via the formulation \eqref{eigen-equation23}. The right-hand side of \eqref{eigen-equation23} is a perturbation compared with $\frac{\xi^{2}|V(r)-c|}{1-c}f_+$ due to the third line in \eqref{relation:xi2Q-r-2}.  The leading order behavior on the left-hand side is also the oscillatory Airy function in \eqref{def:O-Airy-intro2}. See (2-1) in Lemma \ref{lem:f+4} for the details.
In the fourth part, we solves $f_+$ starting from the third part for $r\in[r_c-C\xi^{-\f23}\langle r_c\rangle^{\f13}, r_c+C\xi^{-\f23}\langle r_c\rangle^{\f13}]$ via the same formulation \eqref{eigen-equation23}.  From \eqref{big-xi!}, we see that the length of this interval is smaller than $M^{-1}r_c$ ($M\gg 1$). Upon applying the Langer transform, one finds that the right-hand side of \eqref{eigen-equation23} is still a perturbation.  The leading-order behavior on the left-hand side remains in the form of the Airy function given in \eqref{def:O-Airy-intro2}, but this time without oscillation, since $\xi^{\f23}|\tau|\lesssim 1\Longleftrightarrow r\in[r_c-C\xi^{-\f23}\langle r_c\rangle^{\f13}, r_c+C\xi^{-\f23}\langle r_c\rangle^{\f13}]$.
See (2-2) in Lemma \ref{lem:f+4} for the details.
\vspace{1.5mm}

\item  In the large-spectrum and medium-frequency region
$$D_7=\big\{c\in(1-\delta,1),M^2(1-c)^{\f12}\lesssim \xi\lesssim M(1-c)^{\f13},\;\delta\ll 1,\;M\gg 1\;\text{are}\;\text{fixed}\big\},$$
where the parameter in the potential, $\f{\xi^2}{1-c}\in [M^4, M^2(1-c)^{-\f13}]$ is of large but not excessively large magnitude,  we again identify four intervals in $r$-axis to classify the distinct behaviors of the solution, where in three of them we solve $\phi$ and in the remaining one we solve $f_+$ from infinity.

In the first part, we solve $\phi$ via the formulation \eqref{eigen-equation11} for $r\in (0, M\f{(1-c)^{\f12}}{\xi}]$, where the leading-order behavior is still $\xi r$ and the {right}-hand side is still a perturbation compared with the $r^{-2}\phi$. See  Lemma \ref{lem:phi-sqrt(1-c)/xi}  for the details. In the second part, we solve $\phi$ for $r\in [ M\f{(1-c)^{\f12}}{\xi}, M^{-\f72}\f{\xi^2}{1-c}]$ via the following formulation
\begin{align}  \label{eigen-equation17}		
	\phi''+r^{-1}\phi'-r^{-2}\phi-\f{\xi^2(1-V)}{1-c}\phi=-\xi^2\phi ,\quad \phi(r)|_{r_{*}}=\phi(r_{*}),\quad \phi'(r)|_{r_{*}}=\phi'(r_{*}),		
\end{align}
where  the initial datum is given by the first part at $r_{*}:=M\f{(1-c)^{\f12}}{\xi}$. The right-hand side of \eqref{eigen-equation17} can be viewed as a perturbation due to  $\xi^2\lesssim r^{-2}\ll \f{\xi^2(1-V(r))}{1-c}$,								where we used
\begin{align*}
	\f{\xi^2(1-V(r))}{1-c}\sim
	\left\{
	\begin{aligned}
		& \f{\xi^2}{1-c}\gg r^{-2},\quad M \f{\sqrt{1-c}}{\xi}\lesssim r\lesssim 1,\\
		& \f{\xi^2r^{-3}}{1-c}\gg r^{-2},\quad 1\lesssim r\lesssim M^{-3.5}\f{\xi^2}{1-c}.
	\end{aligned}
	\right.
\end{align*}
For the right-hand side of \eqref{eigen-equation17}, to flatten the coefficient $\f{\xi^2(1-V(r))}{1-c}$, we apply the Langer transform
\begin{align}  \label{langer-intro7}
	\Phi(v,\xi,c)=\left(1-V(r)\right)^{\f14}r^{\f12}\phi(r,\xi,c),\quad v(r)=\int_r^{+\infty}(1-V(s))^{\f12}ds. \end{align}
Then the leading-order behavior is $C_{71}(\xi,c)\left(1-V(r)\right)^{-\f14}r^{-\f12}\cosh\f{\xi (v(r_{*})-v(r)))}{\sqrt{1-c}}$, where $C_{71}(\xi,c)$, $\widetilde{C}_{71}(\xi,c)$ are constants fitting the initial datum. See Lemma \ref{lem:Phi81} for the details. In the third part, we solve $\phi$ for $r\in [ M^{-\f72}\f{\xi^2}{1-c}, M^{-\f12}\xi^{-1}]$ via the following formulation
\begin{align}  \label{eigen-equation172}		
	\phi''+r^{-1}\phi'-r^{-2}\phi+\xi^2\phi
	=\f{\xi^2(1-V)}{1-c}\phi ,\;\phi(r)|_{r_{*2}}=\phi(r_{*2}),\;\phi'(r)|_{r_{*2}}=\phi'(r_{*2}),		
\end{align}
where  the initial datum is given by the second part at $r_{*2}:=M^{-\f72}\f{\xi^2}{1-c}$. We view $\f{\xi^2(1-V(r))}{1-c}$ as a bounded perturbation which is similar to \eqref{eigen-equation11} when $\xi\lesssim (1-c)^{\f12}$, but the reason is different since we no longer have a natural upper bound for $\f{\xi^2}{1-c}$ in $D_7$. Indeed, since $M^2(1-c)^{\f12}\lesssim \xi$ implies $M^{-\f72}\f{\xi^2}{1-c}\gtrsim M^{\f12}$, it follows for $r\in [ M^{-\f72}\f{\xi^2}{1-c}, M^{-\f12}\xi^{-1}]$ that  $$\f{\xi^2(1-V(r))}{1-c}\sim \f{\xi^2r^{-3}}{1-c} \lesssim M^{3.5}r^{-2}\lesssim M^{2.5},$$
which is rather large but uniformly (in $r,\xi,c$) bounded. Moreover, its positivity helps us to get the lower bound of $\phi$.
Then the leading-order behavior for the left-hand side is $ C_{72}(\xi,c)J_1(\xi r)$, where $C_{72}(\xi,c)$ is a constant fitting the initial datum. See Lemma \ref{lem:Phi82} for the details. In the fourth part (see Lemma  \ref{lem:k<1-f+}), we solve $f_+$ on $r\in [\xi^{-1},+\infty)$ via the formulation \eqref{eigen-equation21}. The right-hand side $\f{\xi^2(1-V)}{1-c}f_+$ is uniformly bounded  by $\xi^2f_+$, since $\f{\xi^2(1-V)}{1-c}\simeq \f{\xi^2r^{-3}}{1-c}\lesssim \f{\xi^2\xi^3}{1-c}\lesssim M^3\xi^2$, using $\xi\lesssim M(1-c)^{\f13}$. The leading-order behavior for the left-hand side of \eqref{eigen-equation21} is again the Hankel's function of the first kind $H_1(\xi r)=J_1(\xi r)+iY_1(\xi r)$.
\end{enumerate}

\subsubsection{The resolvent estimate}
Once we obtain the global behavior of $\phi(r,c,k)$, we can plug into the expression (for instance \eqref{kernel example}) and estimate the $L^{1}_{dc}$-norm \eqref{kernel deri example} of its derivative. More precisely, we plug in the profiles for $\phi(r,c,k)$ and $\phi(s,c,k)$ in the 7 regimes displayed by Figure \ref{fig2} to the expression \eqref{kernel example} and then differentiate in $\sqrt{c}$. One challenge in this estimate is that differentiating in $\sqrt{c}$ may generate factors which are not integrable in $c$. To overcome this, we trade off some decay in the Fourier basis $\phi(r,c,k)$ to compensate the divergence in $c$. Let us consider a prototypical kernel estimate (see Proposition \ref{lem:Boundness-K-pieces1}) in the regime $D_{1}:=\{(c,\xi)|c\in (0,1), 0<|\xi|\lesssim (1-c)^{\frac12}\}$:
\begin{align}\label{kernel esti exam}
\int_0^{1}\left|\pa_c\left(p.v.\int_{\mathbb{R}} \left(1-\chi\left(\f{M^2(1-c)^{\f12}}{\xi}\right)\right)
\widetilde{\phi}_{j}(s,\xi,c)\widetilde{\phi}_{j'}(r,\xi,c)\f{e^{i\xi \sqrt{\f{c}{1-c}}z}}{\xi}d\xi\right)\right|dc.
\end{align}
Here $\widetilde{\phi}_{j}(\xi,c,s)$ and $\widetilde{\phi}_{j'}(\xi,c,r)$ are re-normalized Fourier basis and we consider the case where both $r$ and $s$ are close to $0$. After changing the variables 
$$\xi\mapsto \eta_{1}:=\xi s,\quad y_{1}(c,\lambda):=\sqrt{\frac{c}{1-c}}\frac{z}{s}$$ and denoting the parameters $(z,s,r)$ by $\lambda$, the integral in \eqref{kernel deri example} reduces to
\begin{align}\label{kernel esti exam 1}
I_{1_{0},1_{0}}:=\int_{0}^{1}\left|\partial_{c}\left(\text{p.v.}\int_{\bbR}a_{1_{0},1_{0}}(\eta_{1},c,\lambda)\frac{e^{i\eta_{1}y_{1}(c,\lambda)}}{\eta_{1}}d\eta_{1}\right)\right|dc
\end{align}
where $a_{1_{0},1_{0}}$ behaves like
\begin{align*}
\left((1-c)^{\ell}\right)\left(\partial_{c}^{\ell}a_{1_{0},1_{0}}\right)(\eta_{1},c,\lambda)=\chi(\eta_{1})\cdot O\left(\eta_{1}^{1-\frac{\delta}{3}}\right),\quad \text{for}\quad 0<\delta\ll1 \quad \text{and}\quad \ell=0,1.
\end{align*}
Then the integral $I_{1_{0},1_{0}}$ is bounded by
\begin{align}\label{kernel esti exam 2}
\begin{split}
	I_{1_{0},1_{0}}\leq& \int_{0}^{1}\left|\text{p.v.}\int_{\bbR}(\partial_{c}a_{1_{0},1_{0}})\left(\eta_{1},c,\lambda\right)\frac{e^{i\eta_{1}y_{1}(c,\lambda)}}{\eta_{1}}d\eta_{1}\right|dc\\
	&+\int_{0}^{1}\left|\left(\text{p.v.}\int_{\bbR}a_{1_{0},1_{0}}(\eta_{1},c,\lambda)e^{i\eta_{1}y_{1}(c,\lambda)}d\eta_{1}\right)\right|\cdot \left|\frac{\partial y_{1}(c,\lambda)}{\partial c}\right|dc
\end{split}
\end{align}
The 2nd term on the RHS of \eqref{kernel esti exam 2} can be controlled using the oscillatory effect of the integrand (see Lemma \ref{lem: pifi-Linfty} and Lemma \ref{lem:symbol-chi}). The main challenge comes from the 1st term, which, by using the behavior $(1-c)\partial_{c}a_{1_{0},1_{0}}=\chi\left(\eta_{1}\right)O\left(\eta_{1}^{1-\frac{\delta}{3}}\right)$, can be bounded by
\begin{align*}
\int_{0}^{1}(1-c)^{-1}\left|\text{p.v.}\int_{\bbR}\chi(\eta_{1})O(\eta_{1})\frac{e^{i\eta_{1}y_{1}(c,\lambda)}}{\eta_{1}}d\eta_{1}\right|dc,
\end{align*}
whose integral in $c$ is divergent. To handle this issue, we write
\begin{align*}
(1-c)^{-1}=&(1-c)^{-1+\frac{\delta}{6}}s^{\frac{\delta}{3}}\cdot |\xi|^{\frac{\delta}{3}}(1-c)^{-\frac{\delta}{6}}\cdot (|\xi|s)^{-\frac{\delta}{3}}\\
=&(1-c)^{-1+\frac{\delta}{6}}s^{\frac{\delta}{3}}\cdot\left(\frac{|\xi|}{\sqrt{1-c}}\right)^{\frac{\delta}{3}}|\eta_{1}|^{-\frac{\delta}{3}}\\
\lesssim &(1-c)^{-1+\frac{\delta}{6}}s^{\frac{\delta}{3}}|\eta_{1}|^{-\frac{\delta}{3}}.
\end{align*}
Here we have used the fact that $|\xi|\lesssim \sqrt{1-c}$ in $D_{1}$. Therefore the 1st term on the RHS of \eqref{kernel esti exam 2} can be bounded by
\begin{align*}
\int_{0}^{1}(1-c)^{-1+\frac{\delta}{6}}s^{\frac{\delta}{3}}\left|\text{p.v.}\int_{\bbR}\chi(\eta_{1})|\eta_{1}|^{-\frac{\delta}{3}}d\eta_{1}\right|dc\lesssim s^{\frac{\delta}{3}},
\end{align*}
which suffices. When $(c,\xi)$ belongs to other regimes the non-integrable factor generated upon differentiating in $c$ are different, but as above we compensate by trading off decay of the Fourier basis. The details are given in Section 7 and Section 8.
\subsection{Notations}
We define the norms
\begin{align*}
&\|F(\cdot,\cdot)\|_{L^1(\mathbb{R}^+\times\mathbb{R})}
:=\int_{\mathbb{R}}\int_{\mathbb{R}^+}|F(r,z)|rdrdz,\\
&\|F(\cdot,\cdot)\|_{W^{k,1}_{\omega}(\mathbb{R}^+\times\mathbb{R})}:
=\sum_{i+j=0}^k\|\omega(r)\pa^{i}_r\pa^{j}_zF(r,z)\|
_{L^1(\mathbb{R}^+\times\mathbb{R})},
\end{align*}
where $\omega(r)$ denotes a weight function.  For $1\leq p,q\leq \infty$, we define the iterated  norm as
\begin{align*}
&\|F(c,x)\|_{L^p_cL^q_x}:= \|\|F(c,x)\|_{L^p_c}\|_{L^q_x}. \end{align*}

The Fourier transform in $z$ and its inverse are defined by
\begin{align*}
\hat{u}(k)=\mathcal{F}u=\int_{\mathbb{R}}u(z)e^{-izk}dz,\quad
\mathcal{F}^{-1}v=\f{1}{2\pi}\int_{\mathbb{R}}v(k)e^{izk}dk.
\end{align*}

We use the notation
\begin{align*}
O_{\xi,X,c}^{\xi,\rho^X,\rho^c}\big(f_0(\xi,X,c)\big),\,\,\,
X\in\{r,s,\tau,t,x,y...\}
\end{align*}
to represent a function $f(X,\xi,c)$ satisfying
\begin{align}\label{symbol:f+}
\left|(\rho^c)^l(\rho^X)^i\pa_c^l\pa_X^{i}(\xi\pa_{\xi})^jf(X,\xi,c)\right|\lesssim f_0(X,\xi,c),\;\;\;i,j\in\mathbb{N}
,\;l\in\{0,1\},
\end{align}
for some weight functions $\xi$, $\rho^c$,  $\rho^X$.  For simplicity, we sometimes omit the superscripts and subscripts.
The notation $O\left(f_0(\xi,X,c)\right) $ denotes $\left|f\right|\lesssim f_0$, and $f\sim g$ denotes
\begin{align*}
C^{-1}|g|\leq |f|\leq C|g|. \end{align*}

We introduce smooth cut-off functions $\chi(\xi), \chi_+(\xi)$ on $\mathbb{R}$ as follows
\begin{align}\label{def:chi-notation}
\begin{split}
	\chi(\xi)=
	\left\{
	\begin{aligned}
		&1\quad\quad\quad\quad\quad \;\;|\xi|\leq 1,\\
		&\text{smooth}\quad 1\leq |\xi|\leq 2,\\
		&0\quad\quad\quad\quad \quad \;\;|\xi|\geq 2.
	\end{aligned}
	\right.
	\quad\quad\quad
	\chi_+(\xi)=
	\left\{
	\begin{aligned}
		&1\quad\quad\quad\quad\quad \;\;\xi\leq  1,\\
		&\text{smooth}\quad 1\leq \xi\leq 2,\\
		&0\quad\quad\quad\quad \quad \;\;\xi\geq 1.
	\end{aligned}
	\right.
\end{split}
\end{align}
\medskip

This paper is organized as follows.
In Section 2, we derive the spectral properties of the operator of $\calA$, establish the formal formulas in \eqref{fm:hg} rigorously, and reduce the proof of the main theorem to estimates for kernels of the form \eqref{kernel formal}. Sections 3-6 are devoted to constructing a distorted Fourier basis, which we then employ in Sections 7-8 to derive uniform estimates for kernels of the form \eqref{kernel formal}.


\section{Proof of main result assuming estimates on the resolvent kernel}

\subsection{Spectral calculus}

We know that
\begin{align*}
\Delta_{1,k}^{-1}u&=\f{\pi}{2}K_1(kr)\int_0^rI_1(ks)su(s)ds+
	\f{\pi}{2}I_1(kr)\int_r^{+\infty}K_1(ks)su(s)ds,
\end{align*}
where $I_1(z)$, $K_1(z)$ are the modified Bessel functions of the first and second kind, respectively.  Using the asymptotic properties of  $I_1(z), K_1(z)$ near $z=0$ and $+\infty$, we can deduce that
$\|\Delta_{1,k}^{-1}u\|_{L^2(rdr)} \le C\|u\|_{L^2(rdr)}.$ Notice that  $V^{-1}$ is bounded from below and above. Thus,  for $k\neq 0$, we define the domain of $\calA$ as
$\mathcal{D}(\mathcal{A})= L^2(\mathbb{R}^+;rdr).$

\begin{lemma}\label{lem: spectrum calA}
	Let $k\neq 0$. It holds that
\begin{itemize}
\item[(1)]  $\mathcal{A}$ is a self-adjoint operator on $L^{2}(\mathbb{R}^+;V^{-1}rdr):=\mathcal{H}$;

\item[(2)]   $\sigma(\mathcal{A})\subset [0,\f{1}{k^2}]$ and $\mathrm{ker}(\mathcal{A})=\mathrm{ker}(\mathcal{A}-k^{-2})=\{0\}$;

\item[(3)]  $\calA$ has only continuous spectrum.
\end{itemize}
\end{lemma}

\begin{proof}
	The self-adjointness can be directly verified under the  inner product $\langle\cdot,\cdot\rangle_{\mathcal{H}}$:
	\begin{align}\label{main inner product}
		\langle f,g\rangle_{\mathcal{H}}:=\int_{0}^{\infty}f(r)\,g(r)\,V(r)^{-1}rdr.
	\end{align}
	
	The proof of $\sigma(\mathcal{A})\subset [0,k^{-2}]$ reduces to showing that for any $\phi\in D(\mathcal{A})$, we have
	\begin{align}\label{spectrum ineqs}
		\langle \mathcal{A}\phi,\phi\rangle_{\mathcal{H}}\geq 0\quad   \text{and}\quad \big\langle(1-k^{2}\mathcal{A})\phi,\phi\big\rangle_{\mathcal{H}}\geq 0.
	\end{align}
	For the first inequality in \eqref{spectrum ineqs}, we set
	\begin{align*}
		\psi=\mathcal{A}\phi,\quad \Psi=V^{-1}\psi,
	\end{align*}
	then
	\begin{align*}
		\langle\mathcal{A}\phi,\phi\rangle_{\mathcal{H}}=&-\int_{0}^{\infty}\Psi(r)\Delta_{1}\Psi(r)rdr\\
		=&\int_{0}^{\infty}\Psi(r)\left(-\partial_{r}^{2}\Psi-r^{-1}\partial_{r}\Psi+r^{-2}\Psi+k^{2}\Psi\right)rdr\\
		=&\int_{0}^{\infty}\Big((\partial_{r}\Psi)^{2}+
		(\Psi/r)^{2}+k^{2}\Psi^{2}\Big)rdr\geq 0,
	\end{align*}
	which also gives $\mathrm{ker}\big(k^{2}\mathcal{A}\big)=\{0\}$. For the second inequality in \eqref{spectrum ineqs}, we set
	\begin{align*}
		\psi=-\Delta_{1}^{-1}\phi,\quad \phi=-\Delta_{1}\psi,
	\end{align*}
	then
	\begin{align*}
		\langle(1-k^{2}\mathcal{A})\phi,\phi\rangle_{\mathcal{H}}=&\int_{0}^{\infty}
		k^{2}\psi\,\Delta_{1}\psi\,rdr+\int_{0}^{\infty}\Delta_{1}\psi\,\Delta_{1}\psi\,V^{-1}rdr\\
		=&\int_{0}^{\infty}\Delta_{1}\psi\left(\Delta_{1}\psi
		+k^{2}V\psi\right)V^{-1}rdr\\
		=&\int_{0}^{\infty}\left(\Delta_{1}\psi+k^{2}V(r)
		\psi\right)^{2}V^{-1}rdr+\int_{0}^{\infty}-\left(\Delta_{1}\psi
		+k^{2}V(r)\psi\right)k^{2}\psi\,rdr, \end{align*}
	where the second term on the RHS above is also non-negative. Indeed,
	\begin{align*}
		&\int_{0}^{\infty}\left(-\partial_{r}^{2}\psi
		-r^{-1}\partial_{r}\psi+r^{-2}\psi+k^{2}\psi-k^{2}V(r)\psi\right)k^{2}\psi\,rdr\\
		&=\int_{0}^{\infty}k^2\left((\partial_{r}\psi)^2
		+(\psi/r)^2+k^2(1-V(r))
		\psi^2 \right)\,rdr\\
		&\geq \int_{0}^{\infty}k^2\left(
		(\partial_{r}\psi)^2+(\psi/r)^2\right)\,rdr, \end{align*}
	which also gives $\mathrm{ker}\big(\mathcal{A}-k^{-2}\big)=\{0\}$.
	
To prove that $\calA$ has only continuous spectrum, it suffices to show that $\calA$ has no pure point spectrum, since $\calA$ is self-adjoint. To this end, we argue by contradiction. Assume there exists a non-trivial solution $\psi\in L^{2}(\bbR^+;rdr)$ to the equation $\Delta_{1,k}(V(r)^{-1}\psi)+\lambda \psi=0$. Set $f=V(r)^{-1}\psi$, the above equation becomes $\Delta_{1,k}f+\lambda V(r) f=0$. When $\lambda\neq k^{2}$, as $r\rightarrow\infty$, the equation behaves like
\beno
f^{\prime\prime}+\frac{1}{r}f^{\prime}-\mu^{2}f=0, \quad \mu^{2}=k^{2}-\lambda.
\eeno
When $\mu^{2}<0$, this asymptotic equation admits two oscillatory Bessel solutions with asymptotic behavior $\frac{1}{\sqrt{r}}$, which is not in $L^{2}(\bbR^+; rdr)$. When $\mu=0$, the solution to $f^{\prime\prime}+r^{-1}f^{\prime}=0$ behaves like $r^{2}$, which is again not in $L^{2}(\bbR^+; rdr)$. When $\mu^{2}>0$, there is an exponentially growing solution and a exponentially decaying solution as $r\rightarrow\infty$. In this case, we   rewrite the equation of $f$ as follows
\beno
f^{\prime\prime}+\frac{1}{r}f^{\prime}-\left(\frac{1}{r^{2}}+k^{2}-\lambda V(r)\right)f=0.
\eeno
 But now $\frac{1}{r^{2}}+k^{2}-\lambda V(r)>\frac{1}{r^{2}}$ is strictly positive, which contradicts the fact that $-f^{\prime\prime}-\frac{1}{r}f^{\prime}$ being positive if $f\in L^{2}((0,\infty),rdr)$.
 When $\mu=0$, the asymptotic equation as $r\rightarrow\infty$ becomes
 \begin{align*}
 	f^{\prime\prime}+\frac{1}{r}f^{\prime}-\frac{1}{r^{2}}f=0,
 \end{align*}
 which admits two solutions $r$ and $r^{-1}$, neither of which is in $L^{2}(rdr)$. This together with the sign of the potential as in the case $\mu^{2}>0$ excludes the resonance.
\end{proof}
\begin{remark}\label{rmk: spectral interval}
	As we shall see, for any $c\in (0,1)$ and fixed $k\in \bbR^{2}\setminus\{0\}$ we solve the equation $\calA\phi=\frac{c}{k^{2}}\phi$ and the above lemma implies that the solution $\phi$ is neither an eigenfunction nor a resonance. This in turn implies that the continuous spectrum in fact exhausts the entire interval $[0,k^{-2}]$.
\end{remark}
Using the spectral calculus of self-adjoint operator $\mathcal{A}$,  the formula \eqref{fm:hg} of the solution to \eqref{eq:vel-linearize-gernal'} becomes rigorous.
This along with the relation \eqref{eq:vel-linearize-gernal} yields
\begin{align}
	\label{fm:urthetaz}  \begin{split}
		\hat{u}^r(t,r,k)&=ikV^{-1}\cos (k\sqrt{\mathcal{A}} t)\left[V\Delta_1^{-1}(\hat{\omega}^{\theta}_0) \right]+kV^{-1}\sqrt{\mathcal{A}}\sin (k\sqrt{\mathcal{A}} t)
		\left[2u\hat{u}^{\theta}_{0}\right],\\
		\hat{u}^{z}(t,r,k)&=-(\pa_r+r^{-1})\circ V^{-1}\cos (k\sqrt{\mathcal{A}} t)\left[V\Delta_1^{-1}(\hat{\omega}^{\theta}_0) \right]\\
		&\qquad+i(\pa_r+r^{-1})\circ V^{-1}\sqrt{\mathcal{A}}\sin (k\sqrt{\mathcal{A}} t)[2uu_{\theta}^0],\\
		\hat{u}_{\theta}(t,r,k)&= (2u)^{-1}\cos (k\sqrt{\mathcal{A}} t) [2u\hat{u}^{\theta}_{0}]-i(2u)^{-1}\f{\sin (k\sqrt{\mathcal{A}} t)}{\sqrt{\mathcal{A}}}\left[V\Delta_1^{-1}(\hat{\omega}^{\theta}_0) \right].\end{split}
\end{align}

\begin{lemma}  \label{lem:spectral-fm}
     Let $g:[0,\infty)\to \mathbb{C}$ be a continuous function  and $u,v\in C_c^{\infty} \left((0,+\infty)\right)$ be complex-valued functions, $\ep,\delta>0$. It holds that  \begin{align} \label{spectal-formula}\left(g(\mathcal{A})u,v\right)_{\mathcal{H}}&= \lim_{\delta\to0}\lim_{\epsilon\to 0}\f{1}{\pi }\int_{\delta}^{1-\delta}g\left(\f{c}{k^2}\right)\mathrm{Im}\left(\left(
   (k^2\mathcal{A}-c-i\ep)^{-1}u,v\right)_{\mathcal{H}}\right)
   dc, \end{align}
   provided that the  limit exists. Moreover, the limit has an alternative form
   \begin{align}   \label{spectal-formula-0}
  \left(g(\mathcal{A})u,u\right)_{\mathcal{H}}= \lim_{\delta\to0}
   \lim_{\epsilon\to 0}\f{1}{2\pi i}\int_{\delta}^{1-\delta}g\left(\f{c}{k^2}\right)\left(
   \left(k^2\mathcal{A}-c-i\epsilon\right)^{-1}u-
   \left(k^2\mathcal{A}-c+i\epsilon\right)^{-1}u,u\right)_{\mathcal{H}}dc. \end{align}
 \end{lemma}
      \begin{proof}
       The  Polarization identity writes
 \begin{align*}
4\left(
   \mathcal{L}u,v\right)_{\mathcal{H}}
= &\left( \mathcal{L}(u+v),(u+v)\right)_{\mathcal{H}}- \left( \mathcal{L}(u-v),(u-v)\right)_{\mathcal{H}}\\
&+i\left( \mathcal{L}(u+iv),(u+iv)\right)_{\mathcal{H}}
-i\left( \mathcal{L}(u-iv),(u-iv)\right)_{\mathcal{H}}.\end{align*}
We take $ \mathcal{L}=g(\mathcal{A}), (k^2\mathcal{A}-c-i\ep)^{-1}$ in the above identity, then \eqref{spectal-formula-0} implies \eqref{spectal-formula}.  It remains to prove
       \eqref{spectal-formula-0}.

       Let $P_{\mathcal{A}}(\cdot)$ be the spectral projector associated with the self-adjoint operator $\mathcal{A}$. We notice by Lemma \ref{lem: spectrum calA} that $\mathrm{Ran}\big(P_{\mathcal{A}}(\{0\})\big)
 =\mathrm{ker}\big(\mathcal{A}\big)=\{0\}$,
   $\mathrm{Ran}\big(P_{\mathcal{A}}(\{k^{-2}\})\big)
 =\mathrm{ker}\big(\mathcal{A}-k^{-2}\big)=\{0\}$. Moreover, by $\sigma(\mathcal{A})\subset [0,k^{-2}]$, it follows that
 \begin{align*} \mathrm{1}=P_{\mathcal{A}}([0,k^{-2}])=P_{\mathcal{A}}((0,k^{-2}))=
 P_{\mathcal{A}}((0,k^{-2}]).
  \end{align*}
  For the left-hand side of  \eqref{spectal-formula-0}, it follows from the spectral theorem and the above facts that
  \begin{align*}
\Big(g(\mathcal{A})u,u\Big)_{\mathcal{H}} &=\int_{(0,k^{-2})}g(x)
 d\mu_{u,u}(x). \end{align*}
Here, the spectral measure $\mu_{u,v}$ is defined as
 \begin{align*}
\mu_{u,v}(I):= \big(P_{\mathcal{A}}(I)u,v\big)_{\mathcal{H}},\;
\;\mu_{u,v}(x)=\mu_{u,v}\big((-\infty,x)\big),\end{align*}
 which is supported in $(0,k^{-2}]$ and right-continuous, with $\mu_{u,v}(\{0\})=\mu_{u,v}(\{k^{-2}\})=0$.
  For the right-hand side of  \eqref{spectal-formula-0},
   using the spectral theorem again, we write
  \begin{align*}
\notag&\left(
   \big(\mathcal{A}-\lambda-i\eta\big)^{-1}u,u\right)_{\mathcal{H}}
=\int_{\mathbb{R}}\f{1}{x-\lambda-i\eta}d\mu_{u,u}(x)
:=F_+(\lambda+i\eta).
\end{align*}
One can verify that $F_+(\lambda+i\eta)$ is a Herglotz function, i.e., $F_+(\cdot)$ is analytic on $\mathbb{C}^+$ and $\mathrm{Im}(F_+)>0$. In the same way, we write
 \begin{align*}
\notag&\left(
   \big(\mathcal{A}-\lambda+i\eta\big)^{-1}u,u\right)_{\mathcal{H}}
=\int_{\mathbb{R}}\f{1}{x-\lambda+i\eta}
d\mu_{u,u}(x):=F_-(\lambda-i\eta)=
\overline{F_+(\lambda+i\eta)}, \end{align*}
which gives
\begin{align*}
 &\f{1}{2\pi i}\left(\big(\big(\mathcal{A}-\lambda-i\eta\big)^{-1} -\big(\mathcal{A}-\lambda+i\eta\big)^{-1}\big)u,u\right)_{\mathcal{H}}= \f{1}{\pi}\mathrm{Im}\big(F_+(\lambda+i\eta)
   \big).
\end{align*}
Therefore,  it follows from the property of Herglotz  function (see the proof of \cite[Cor 3.24, p.124]{G-book})  that
\begin{align*}
		&\lim_{\delta\to0}\lim_{\epsilon\to 0}\f{1}{2\pi i }\int_{\delta}^{1-\delta}g\left(\f{c}{k^2}\right)
\left(\big(\big(\mathcal{A}-\lambda-i\eta\big)^{-1} -\big(\mathcal{A}-\lambda+i\eta\big)^{-1}\big)u,u\right)_{\mathcal{H}}dc\\
		&=_{\lambda=\f{c}{k^2},\;\eta=\f{\ep}{k^2}}\lim_{\delta\to 0}\lim_{\ep\to0}\f{1}{\pi}\int_{\delta}^{\f{1-\delta}{k^2}}
		g(\lambda)\mathrm{Im}\left(F_+(\lambda+i\eta)
		\right)d\lambda\\
&=\int_{(0,k^{-2})}g(x)
		d\mu_{u,u}(x), \end{align*}
   which matches the left-hand side.
    \end{proof}

 \subsection{The limiting absorption}

 We first introduce the following definitions.

 \begin{definition}\label{def:useful-fm}
	For any $0<\delta\ll 1$, we denote
	$$S_{\delta}=\big\{\mathrm{Re}(c)\in[\delta,1-\delta],\;\;|\mathrm{Im}(c)|\leq 1\big\}.$$
	For $c\in\mathbb{C}/\{0\}$, $k\in\mathbb{R}/\{0\}$, let $\xi(c,k):=k\sqrt{\f{1}{c}-1}$.
	Then $\xi(\cdot,k)$ is analytic in $S_{\delta}$.
	
	
	Let $\phi(r,c,k)$, $\theta(r,c,k)$, $f_+(r,c,k)$ be smooth solutions of the ODE  with $\mathrm{Re}(c)\in(0,1), k\in\mathbb{R}/\{0\}, r\in\mathbb{R}^+$:
	\begin{align}  \label{eigen-equation}
		f''+r^{-1}f'-r^{-2}f-k^2f+\f{k^2V}{c}f=0
		\end{align}
	with
	\begin{align}
		\phi(r,c,k)\sim \xi(c,k) r(r\to 0),\quad \theta(r,c,k)\sim \f12(\xi(c,k) r)^{-1}\quad\text{as}\quad r\to 0,
		\end{align}
	and as $r\to +\infty$,
	\begin{align}
		\begin{split}
			&f_+(r,c,k)\sim \left\{
			\begin{aligned}
				&(\xi r)^{-\f12}e^{i\xi r},\;\quad\quad \quad\quad |\xi|\lesssim (1-\mathrm{Re}(c))^{\f13},\\
				&(\xi r)^{-\f12}e^{i\xi\int_{0}^r
					\sqrt{\f{V(s)-\mathrm{Re}(c)}{1-\mathrm{Re}(c)}}ds},\;
				|\xi|\gtrsim (1-\mathrm{Re}(c))^{\f13},\;\mathrm{Re}(c)\in(0,V(0)],\\
				&(\xi r)^{-\f12}e^{i\xi\int_{r_c}^r
					\sqrt{\left|\f{V(s)-\mathrm{Re}(c)}{1-\mathrm{Re}(c)}\right|}ds},\;
				|\xi|\gtrsim (1-\mathrm{Re}(c))^{\f13},\;\mathrm{Re}(c)\in(V(0),1).			
				\end{aligned}
			\right.\label{eq:f+behave-complex}
			\end{split}
	\end{align}
Here, we use the notation $f(r)\sim f_0(r)$ to represent $\lim\f{f_0}{f}=\lim \f{f_0'}{f'}=1$.  Then $\phi(c)$, $\theta(c)$ are analytic in $S_{\delta}$ and real-valued if $\mathrm{Im}(c)=0$; $f_+(c)$ is continuous on $S_{\delta}$.
	
	The Wronskian of \eqref{eigen-equation}  is defined by
	\begin{align}\label{fm:wrons-op}\mathcal{W}(f,g)
		:=r(f'g-fg')=(r^{\f12}f)'(r^{\f12}g)-(r^{\f12}g)'(r^{\f12}f).
	\end{align}
	Then we can deduce that $\mathcal{W}(\phi,\theta)\equiv 1$, and $\mathcal{W}(f_+,\theta)$, $\mathcal{W}(\phi,f_+)$ are complex-valued functions which only depend on $c,k$.
	We adopt the following notations throughout the remainder of the paper:
	\begin{align}\label{fm:wrons}
		W(c,k):= \mathcal{W}(\phi,f_+)\;\;\text{or}\;\;W(\xi,c):= \mathcal{W}(\phi,f_+).
	\end{align}
\end{definition}

We notice by the definition that $\phi(c), \theta(c)$ are analytic in  $\mathbb{C}/\{0\}$, $\bar{\phi}(c)=\phi(\bar{c})$, $\bar{\theta}(c)=\theta(\bar{c})$;
$\bar{f_+}(c)=f_-(\bar{c})$. Let $$\textbf{c}:=c+i\ep\in S_{\delta},\; \text{i.e.}, \;c\in [\delta,1-\delta],\; \ep\in[0,1].$$ Then $\mathrm{\xi(c,k)}\geq 0$ and $\xi(c)$ is analytic in $S_{\delta}$,
and $|\xi(c)|$ has positive upper and lower bounds in $S_{\delta}$.  It follows from \eqref{eq:f+behave-complex} that for $r\geq 1$ and  $\textbf{c}\in S_{\delta}$, $|f_+(r,c,k)|\lesssim r^{-1} e^{-\mathrm{Im}(\xi)r}\in L^2((1,+\infty);r^Ndr)$ for any $N\geq 1$.  For $r\leq 1$ and  $\textbf{c}\in S_{\delta}$, $|\phi(r,c,k)|\lesssim |\xi| r\in  L^2((0,1);rdr)$. Then for $\textbf{c}\in S_{\delta}$, $u\in L^2(\mathbb{R}^+;rdr)$, it can be directly verified that
\begin{align*}
	\Phi(r,\textbf{c},k):=\left(\Delta_1+\f{k^2V}{\textbf{c}}\right)^{-1}u
	=\int_0^r
	\phi(s,\textbf{c})\tilde{f}_+(r,\textbf{c})su(s)ds
	+\int_r^{+\infty}
	\phi(r,\textbf{c})\tilde{f}_+(s,\textbf{c})su(s)ds,
	\end{align*}
where  $\tilde{f}_+(r,\textbf{c},k):=\f{f_+(r,\textbf{c},k)}{\mathcal{W}(\phi,f_+)(\textbf{c},k)}$   has the linear expression
\begin{align}\label{f+=aphi+btheta}
	\tilde{f}_+(r,\textbf{c},k)=\theta(r,\textbf{c},k)
	+m(\textbf{c},k)\phi(r,\textbf{c},k),\;\; m(\textbf{c})=\f{\mathcal{W}(f_+,\theta)(\textbf{c})}
	{\mathcal{W}(\phi,f_+)(\textbf{c})},
	\end{align}
due to $\mathcal{W}(\phi,\theta)\equiv 1$.
Let us claim that for $c\in(0,1)$,
\begin{align}    \label{m=1/|W|^2}
	\mathrm{Im}(m(c))= \f{1}{|\mathcal{W}(\phi,f_+)(c,k)|^2}=\f{1}{|W(\phi,f_+)(c,k)|^2},
	\end{align}
where $\mathcal{W}(\phi,f_+)(c)\neq 0$, $c\in(0,1)$.  Indeed, by the asymptotic at $r=+\infty$ of $f_+$ in \eqref{eq:f+behave-complex}, we infer that for $c\in(0,1)$,
\begin{align}\label{W:2i}
	\mathcal{W}(f_+,\bar{f}_+)(c)=2i,
\end{align}
which along with $f_+(c,r)=\mathcal{W}(\phi,f_+)(c)\theta(c,r)+\mathcal{W}(f_+,\theta)(c)\phi(c,r)$(by \eqref{f+=aphi+btheta}), gives
\begin{align} \label{ImW-relation}
	\notag2=\mathrm{Im}\left(\mathcal{W}(f_+,\bar{f}_+)(c)\right)
	&=\mathrm{Im}\left(\mathcal{W}(\phi,f_+)(c)\mathcal{W}(\theta,\bar{f}_+)(c) +\mathcal{W}(f_+,\theta)(c)\mathcal{W}(\phi,\bar{f}_+)(c)\right)\\
	&=2\mathrm{Im}\left(\mathcal{W}(f_+,\theta)(c)\mathcal{W}(\phi,\bar{f}_+)(c)\right),
	\end{align}
which yields $\mathrm{Im}\left(\mathcal{W}(f_+,\theta)(c)\mathcal{W}(\phi,\bar{f}_+)(c)\right)=1$. Here, in the last line of \eqref{ImW-relation}, we used that $\phi(c)$, $\theta(c)$ are real-valued, and $$\overline{\mathcal{W}(f_+,\theta)(c)\mathcal{W}(\phi,\bar{f}_+)(c)}
=-\mathcal{W}(\theta,\bar{f}_+)(c)\mathcal{W}(\phi,f_+)(c)\quad \text{for}\,\, c\in(0,1).$$
Therefore, we obtain
\begin{align*}
	&\mathrm{Im}\big(m(c)\big)=\mathrm{Im}\left(\f{\mathcal{W}(f_+,\theta)}{\mathcal{W}(\phi,f_+)} \right)=\f{\mathrm{Im}\left(\mathcal{W}(f_+,\theta)(c)
		\mathcal{W}(\phi,\bar{f}_+)(c)\right)}
	{|\mathcal{W}(\phi,f_+)(c)|^2}=\f{1}
	{|\mathcal{W}(\phi,f_+)(c)|^2}.
\end{align*}

Now we are in a position to establish the following crucial representation formulas of spectral calculus.

\begin{lemma}\label{lem:distortedF}
	Let $k\in\mathbb{R}/\{0\}$, $f\in C([0,1])$, $\phi(r,c,k)$ and $W(c,k)$ be given by Definition \ref{def:useful-fm}.
For  $\hat{v}(s,k)\in C_c^{\infty}\big((0,\infty)\times\mathbb{R})\big)$,  it holds that
	\begin{align} \label{spectal-formula-result-Fz}
		&\notag\mathcal{F}^{-1}_{k\to z}\left(kV(r)^{-1}\mathcal{A}^2f(\mathcal{A})\hat{v}(s,k)\right)\\
		&= \f{1}{2\pi^2 }\int_{\mathbb{R}\times\mathbb{R}^+}\left(\int_{0}^1f\left(\f{c}{k^2}\right)
		\left(\int_{\mathbb{R}}\f{\phi(r,c,k)
			\phi(s,c,k)}{|W(c,k)|^2}\f{e^{ik(z-y)}}{k}dk \right)dc\right)v(s,y)sdsdy.
	\end{align}
	For a vector-valued function $\textbf{v}=\begin{pmatrix}
		v_1\\v_2\end{pmatrix}$,  $\hat{v}_1(s,k),\hat{v}_2(s,k)\in C_c^{\infty}((0,\infty)\times \mathbb{R})$, it holds that
	\begin{align} \label{spectal-formula-result-vector-Fz}
	\begin{split}
		&\mathcal{F}_{k\to z}^{-1}\left(kV(r)^{-1}\mathcal{A}^{\f{5}{2}}f(\mathcal{A})(\pa_s,ik)\cdot \hat{\textbf{v}}(s,k) \right)\\
		&= \f{-1}{2\pi^2 }\int_{\mathbb{R}\times\mathbb{R}^+}\left(\int_{0}^1
		\sqrt{\f{c}{k^2}}f\left(\f{c}{k^2}\right)
		\left(\int_{\mathbb{R}}\f{\phi(r,c,k)
			\left(\pa_s+s^{-1},ik\right)\phi(s,c,k)}{|W(c,k)|^2}\f{e^{ik(z-y)}}{k}dk\right)
		dc\right)\\&\qquad\qquad\qquad
		\cdot\textbf{v}(s,y)sdsdy,
		\end{split}
	\end{align}
	provided that the integrals on the right-hand side exist.
\end{lemma}

\begin{proof}
Let us claim that for $v\in C_c^{\infty}((0,\infty))$,
	\begin{align} \label{spectal-formula-result}
		kV(r)^{-1}\mathcal{A}^2f(\mathcal{A})v= \f{1}{\pi k }\int_{\mathbb{R}^+}\left(\int_{0}^1f\left(\f{c}{k^2}\right)\f{\phi(r,c,k)
			\phi(s,c,k)}{|W(c,k)|^2}dc \right)v(s)sds,
	\end{align}
and for a vector-valued function $\textbf{v}=\begin{pmatrix}
		v_1\\v_2\end{pmatrix}$ with $v_1,v_2\in C_c^{\infty}((0,\infty))$,
		\begin{align} \label{spectal-formula-result-vector}
		\notag &kV(r)^{-1}\mathcal{A}^{\f52}
		f(\mathcal{A})(\pa_s,ik)\cdot \textbf{v} \\
		& = \f{-1}{k\pi }\int_{\mathbb{R}^+}\left(\int_{0}^1
		\sqrt{\f{c}{k^2}}f\left(\f{c}{k^2}\right)\f{\phi(r,c,k) \left(\pa_s+s^{-1},ik\right)\phi
			(s,c,k)}{|W(c,k)|^2}dc\right)\cdot\textbf{v}(s)sds.
		\end{align}
Then the desired formulas   \eqref{spectal-formula-result-Fz} and \eqref{spectal-formula-result-vector-Fz} can be deduced by taking the inverse Fourier transform $\mathcal{F}^{-1}_{k\to z}$ on  \eqref{spectal-formula-result} and \eqref{spectal-formula-result-vector}, and using Fubini's theorem.

Next, we only present the proof of \eqref{spectal-formula-result}, and \eqref{spectal-formula-result-vector}  follows
by  taking $f(x)=\sqrt{x}f$, $v=(\pa_s,ik)\cdot \bf{v}$ in \eqref{spectal-formula-result} and integration by parts once in $s$.

	Since $C_{c}^{\infty}((0,\infty))$ is dense in $\mathcal{H}=L^2(\mathbb{R}^+; V^{-1}rdr)$, \eqref{spectal-formula-result} can be deduced by the following identity  for $u,v\in C_{c}^{\infty}((0,\infty))$:
	 \begin{align} \label{spectal-formula-result-inner}
		\left(V^{-1}\mathcal{A}^2f(\mathcal{A})u,v\right)_{\mathcal{H}}= \f{k^{-2}}{\pi }\int_{(\mathbb{R}^+)^2}\left(\int_{0}^1f\left(\f{c}{k^2}\right)\phi(r,c,k)
		\phi(s,c,k)\f{dc}{|W(c,k)|^2}\right)v(s)su(r)rdsdr.
	\end{align}
	Notice that for $\mathcal{A}=-V\Delta_1^{-1}$,
	\begin{align} \label{resovlent-identity}
		\notag(k^2\mathcal{A}-c-i\epsilon)^{-1}&=
		-\left((c+i\ep)\left(\Delta_1+\f{k^2V}{c+i\ep}\right)\Delta_1^{-1}\right)^{-1}\\
		\notag&= -(c+i\ep)^{-1}\Delta_1 \left(\Delta_1+\f{k^2V}{c+i\ep}\right)^{-1}\\
		&= -(c+i\ep)^{-1}+ (c+i\ep)^{-2} k^2V  \left(\Delta_1+\f{k^2V}{c+i\ep}\right)^{-1}.
	\end{align}
	Then it follows from \eqref{spectal-formula} that the left-hand side of \eqref{spectal-formula-result-inner}  can be computed via
	\begin{align}
		&\lim_{\delta\to0}\lim_{\epsilon\to 0}\f{1}{\pi }\int_{\delta}^{1-\delta}\left(\f{c}{k^2}\right)^2f\left(\f{c}{k^2}\right)\mathrm{Im}\left(
		\left(
		V^{-1}(k^2\mathcal{A}-c-i\epsilon)^{-1}u,v\right)_{\mathcal{H}}\right)dc\notag\\ &=\lim_{\delta\to0}\lim_{\epsilon\to 0}\f{k^{-2}}{\pi }\int_{\delta}^{1-\delta}f\left(\f{c}{k^2}\right)\mathrm{Im}
		\left(\f{c^2}{(c+i\ep)^2}\left(
		\left(\Delta_1+\f{k^{2}V}{c+i\ep}\right)^{-1}u,v\right)_{\mathcal{H}}\right)dc \notag\\
		&\quad
		- \lim_{\delta\to0}\lim_{\epsilon\to 0}\f{k^{-4}}{\pi }\int_{\delta}^{1-\delta}f\left(\f{c}{k^2}\right)\mathrm{Im}
		\left(\f{c^2}{c+i\ep}\right)\left(V^{-1}u,v\right)_{\mathcal{H}}dc\notag\\
		&=\lim_{\delta\to0}\lim_{\epsilon\to 0}\f{k^{-2}}{\pi }\int_{\delta}^{1-\delta}f\left(\f{c}{k^2}\right)\mathrm{Im}
		\left(\f{c^2}{(c+i\ep)^2}\left(
		\left(\Delta_1+\f{k^2V}{c+i\ep}\right)^{-1}u,v\right)_{\mathcal{H}}\right)dc.   \label{fm:lim-delta-epsilon-Im}
	\end{align}
	
	Since $\phi(\textbf{c})$, $\theta(\textbf{c})$, $f_+(\textbf{c})$ are all continuous in $S_{\delta}$, we know that $m(\textbf{c})$ is  continuous in $S_{\delta}$.  Thus, the left-hand side of \eqref{fm:lim-delta-epsilon-Im} can be computed using the dominated convergence theorem that for $u,v\in C_c^{\infty}[r_1,r_2]$,
	\begin{align*}\notag&\lim_{\delta\to0}\lim_{\epsilon\to 0}\f{k^{-2}}{\pi }\int_{\delta}^{1-\delta}f\left(\f{c}{k^2}\right)\mathrm{Im}
	\left(\f{c^2m(c+i\ep)}{(c+i\ep)^2} \int_{(\mathbb{R}^+)^2}\phi(r,c+i\ep)\phi(s,c+i\ep)su(s)rv(r)dsdr\right)dc\\
	\notag&\quad+ \lim_{\delta\to0}\lim_{\epsilon\to 0}
	\f{k^{-2}}{\pi }\int_{\delta}^{1-\delta}f\left(\f{c}{k^2}\right)\mathrm{Im}
	\left(\f{c^2}{(c+i\ep)^2}
	\int_{\mathbb{R}^+}\int_{0}^r\theta(r,c+i\ep)\phi(s,c+i\ep)
	su(s)rv(r)dsdr\right)dc\\
	&\quad+\lim_{\delta\to0}\lim_{\epsilon\to 0}\f{k^{-2}}{\pi }\int_{\delta}^{1-\delta}f\left(\f{c}{k^2}\right)\mathrm{Im}
	\left(\f{c^2}{(c+i\ep)^2}\int_{\mathbb{R}^+}
	\int_r^{+\infty}\phi(r,c+i\ep)\theta(s,c+i\ep)su(s)rv(r)dsdr
	\right)dc\notag\\
	& =\lim_{\delta\to0}\f{k^{-2}}{\pi }\int_{\delta}^{1-\delta}f\left(\f{c}{k^2}\right)\mathrm{Im}
	\left(m(c) \int_{[r_1,r_2]^2}\phi(r,c)\phi(s,c)su(s)rv(r)dsdr \right) dc\\
	&\quad+ \lim_{\delta\to0}\f{k^{-2}}{\pi }\int_{\delta}^{1-\delta}f\left(\f{c}{k^2}\right)\mathrm{Im}
	\left(
	\int_{r_1}^{r_2}\int_{r_1}^r\theta(r,c)\phi(s,c)
	su(s)rv(r)dsdr\right)dc\\
	&\quad+\lim_{\delta\to0}\f{k^{-2}}{\pi }\int_{\delta}^{1-\delta}f\left(\f{c}{k^2}\right)\mathrm{Im}
	\left(\int_{r_1}^{r_2}
	\int_r^{r_2}\phi(r,c)\theta(s,c)su(s)rv(r)dsdr
	\right)dc\notag\\
	& =\lim_{\delta\to0}\f{k^{-2}}{\pi }\int_{\delta}^{1-\delta}f\left(\f{c}{k^2}\right) \mathrm{Im}\big(m(c)\big) \int_{[r_1,r_2]^2}\phi(r,c)\phi(s,c)su(s)rv(r)dsdrdc\\
	& =\lim_{\delta\to0}\f{k^{-2}}{\pi }\int_{\delta}^{1-\delta}f\left(\f{c}{k^2}\right)\f{1}{|W(\phi,f_+)(c,k)|^2} \int_{[r_1,r_2]^2}\phi(r,c)\phi(s,c)su(s)rv(r)dsdrdc.
	\notag\end{align*}
Here, we have used that the functions in the integrals are supported in compact domains and continuous, and $\phi(c), \theta(c)$ are real-valued.
\end{proof}

\begin{remark} \label{rmk:finite-integral}
		As we shall see, the integrals on the right-hand sides of  \eqref{spectal-formula-result} and \eqref{spectal-formula-result-vector} are finite due to the uniform estimates \eqref{bd:phi/|W|}, Lemma \ref{prop:summery} (which is quite non-trivial), i.e.,
		\begin{align*}
			\left|\f{\phi(r,c,k)}{W(c,k)}\right|\leq C,\quad
			\left|\f{\left(\pa_r+r^{-1}\right)\phi(r,c,k)}{W(c,k)}\right|\leq C.
		\end{align*}
	\end{remark}


 \begin{remark}\label{rmk:D-fourier}
	We observe in \eqref{spectal-formula-result} that there is a lose of $\mathcal{A}^2$ or $c^2$ on the right-hand side.
	Indeed, taking $u=v$, $f\equiv 1$, $\widetilde{\phi}=\phi/|W|$ in  \eqref{spectal-formula-result-inner}, we can obtain a Plancherel's formula corresponding to the operator $\mathcal{A}=-V\left(\pa_r^2+r^{-1}\pa_r-r^{-2}-k^2\right)^{-1}$:
	\begin{align} \label{Plancherel-identity}
		\|\mathcal{A}v\|_{L^2(\mathbb{R}^+;V^{-1}rdr)}^2=\f{1}{k^2\pi} \left\|\int_{\mathbb{R}^+}\phi(r,c,k)v(r)rdr\right\|_{L^2\left((0,1);\f{dc}{|W(c,k)|^2}\right)}^2.
	\end{align}
	Then for $k\neq 0$, we  define the so-called distorted Fourier transform from $\mathcal{D}(\mathcal{A})=\{v|v,V\Delta^{-1}_1v\in L^2(\mathbb{R}^+;V^{-1}rdr)\}$ to $L^2((0,1);|W(c,k)|^{-2}dc)$ by
	\begin{align*} \left(\mathcal{F}_{D}v\right)(c)&:=\int
		_{\mathbb{R}^+}\phi(r,c,k)
		v(r)rdr.\end{align*}
	We also obtain a Plancherel's formula corresponding to $\mathcal{A}=-V\left(\pa_r^2+r^{-1}\pa_r-r^{-2}+\pa_z^2\right)^{-1}$:
	\begin{align*}
		& \|\mathcal{A}v\|_{L^2(V^{-1}rdrdz)}^2=\f{1}{\pi}\left\|\int
		_{\mathbb{R}\times\mathbb{R}^+}\left(\int
		_{\mathbb{R}}\f{\phi(r,k,c)}{|W(c,k)|}
		\f{e^{ik(z-y)}}{k}dk\right)v(r,y)rdrdy\right\|_{L^2((0,1)\times \mathbb{R};dcdz)}^2. \end{align*}
	Thus, we also define a 2D distorted Fourier transform from
	$\mathcal{D}(\mathcal{A})=\big\{u=v|v,V\Delta^{-1}_1u\in L^2(\mathbb{R}^+\times\mathbb{R};V^{-1}rdrdz)\big\}$ to $L^2((0,1)\times \mathbb{R};dcdz)$ by
	\begin{align}\label{def:DF} \big(\mathcal{F}_{\textbf{D}}u\big)(c,z)&:=\int
		_{\mathbb{R}\times\mathbb{R}^+}\widetilde{\Phi}(r,z-y,c)v(r,y)rdrdy,
	\end{align}
	where \begin{align}
		\begin{split}\label{def:nomalize-phi}
			&a(c)=\sqrt{\f{c}{1-c}},\quad \xi(c,k)=k/a(c),\quad
			\tilde{\phi}(r,\xi,c):=\f{\phi(r,c,k)}{W(c,k)},\\ &\widetilde{\Phi}(r,z,c):= \int
			_{\mathbb{R}}\widetilde{\phi}(r,\xi,c)
			\f{e^{i\xi a(c)z}}{\xi}d\xi.
		\end{split}
	\end{align}  Indeed, the following estimate holds due to the  delicate analysis  established in this paper that \begin{align}\label{est:DFL1-Linfty} \|\left(\mathcal{F}_{\textbf{D}}v\right)(c,z)\|_{L^{\infty}_{c,z}}\leq C\|v\|_{L^1_{r,z}}.\end{align}
	And    for some $0<\delta \ll 1$,
	\begin{align} \label{est:pa_cDFL1-L1Linfty} \|\pa_c\left(\mathcal{F}_{\textbf{D}}v\right)(c,z)\|_{L^1_cL^{\infty}_z}
		\leq C\|v\|_{L^1_{r^{-\delta}+r^{\delta}}L^1_z}.\end{align}
\end{remark}

\subsection{Proof of Theorem \ref{Thm:dipersive-decay}}

Applying Lemma \ref{lem:distortedF} with $f(x)=\f{\sin \sqrt{x}t}{\sqrt{x}},\;\f{\cos \sqrt{x}t}{\sqrt{x}}$, we rewrite \eqref{fm:urthetaz}  as
 \begin{align}
		\label{fm:urthetaz-good}  \begin{split}
			\hat{u}^r(t,r,k)&=\left(V^{-1}k\mathcal{A}^{\f52}\f{\cos (k\sqrt{\mathcal{A}} t)}{k\sqrt{\mathcal{A}}}\right)\left[ik\mathcal{A}^{-2}\circ V\Delta_1^{-1}(\hat{\omega}^{\theta}_0) \right]\\
			&\qquad+\left(kV^{-1}\mathcal{A}^2\f{\sin (k\sqrt{\mathcal{A}} t) }{k\sqrt{\mathcal{A}}}\right)
			\left[2k\mathcal{A}^{-1}\circ u\hat{u}^{\theta}_{0}\right]\\
			&=\left(V^{-1}\mathcal{A}^{2}\cos (k\sqrt{\mathcal{A}} t)\right)\left[-ik\mathcal{A}^{-1}(\hat{\omega}^{\theta}_0) \right]+\left(V^{-1}\mathcal{A}^{\f32}\sin (k\sqrt{\mathcal{A}} t) \right)
			\left[2k\mathcal{A}^{-1}(u\hat{u}^{\theta}_{0})\right],\\
			\hat{u}^{z}(t,r,k)&=(\pa_r+r^{-1})\circ V^{-1}k\mathcal{A}^{\f52}\f{\cos (k\sqrt{\mathcal{A}} t)}{k\sqrt{\mathcal{A}}}
			\left[-\mathcal{A}^{-2}\circ V\Delta_1^{-1}(\hat{\omega}^{\theta}_0) \right]\\
			&\quad +(\pa_r+r^{-1})\circ V^{-1}k\mathcal{A}^{2}\f{\sin (k\sqrt{\mathcal{A}} t)}{k\sqrt{\mathcal{A}}}\left[2i\mathcal{A}^{-1}\circ uu_{\theta}^0\right]\\
			&=\left((\pa_r+r^{-1})\circ V^{-1}\mathcal{A}^{2}\cos (k\sqrt{\mathcal{A}} t) \right)
			\left[\mathcal{A}^{-1}(\hat{\omega}^{\theta}_0) \right]\\&\qquad+\left((\pa_r+r^{-1})\circ V^{-1}\mathcal{A}^{\f32}\sin (k\sqrt{\mathcal{A}} t)\right)\left[2i\mathcal{A}^{-1}( uu_{\theta}^0)\right],\\
			\hat{u}^{\theta}(t,r,k)&= 2\Omega\left(V^{-1}k\mathcal{A}^{\f52}\f{\cos (k\sqrt{\mathcal{A}} t)}{k\sqrt{\mathcal{A}}}\right) \left[\mathcal{A}^{-2}\circ u\hat{u}^{\theta}_{0}\right]\\&\qquad+\Omega\left(V^{-1}k\mathcal{A}^2\f{\sin (k\sqrt{\mathcal{A}} t) }{k\sqrt{\mathcal{A}}}\right)
			\left[-i\mathcal{A}^{-2}\circ V\Delta_1^{-1}(\hat{\omega}^{\theta}_0) \right]\\
			&= 2\Omega\left(V^{-1}\mathcal{A}^{2}\cos (k\sqrt{\mathcal{A}} t)\right) \left[\mathcal{A}^{-2}( u\hat{u}^{\theta}_{0})\right]+\Omega\left(V^{-1}\mathcal{A}^{\f32}\sin (k\sqrt{\mathcal{A}} t) \right)
			\left[i\mathcal{A}^{-1}(\hat{\omega}^{\theta}_0) \right].\end{split}
	\end{align}
Then Theorem \ref{Thm:dipersive-decay} is a direct consequence of  \eqref{fm:urthetaz-good} and the following  Proposition \ref{prop:solution-fm}.

		\begin{proposition} \label{prop:solution-fm}
			Let $m,n\in\mathbb{N}$, $t\gtrsim 1$. It holds uniformly in $t,r,z$  that for $\hat{v}(r,k)\in C_c^{\infty}\big((0,\infty)\times \mathbb{R}/\{0\})$,			
			\begin{align}\label{bd:sin}
				&\|\pa_z^n\pa_r^m\mathcal{F}_{k\to z}^{-1}\Big( V^{-1} \mathcal{A}^{\f32}\sin (k\sqrt{\mathcal{A}}t)\big[\hat{v}\big]\Big)\|
				_{L^{\infty}(\mathbb{R}^+\times \mathbb{R})}
				\lesssim t^{-1}\|v\|_{W^{1,m+n}_{r^{-\delta}+r^{\delta}}(\mathbb{R}^+\times \mathbb{R})},\\
			\label{bd:cos-left}
				&\|\pa_z^n\pa_r^m\mathcal{F}_{k\to z}^{-1}\left( V^{-1} \mathcal{A}^2\cos (k\sqrt{\mathcal{A}}t)\big[\mathcal{A}^{-1}(\hat{v})\big]\right)\|_{L^{\infty}
					(\mathbb{R}^+\times \mathbb{R})}
				\lesssim t^{-1}\|v\|_{W^{m+n+1,1}_{r^{-\delta}+r^{\delta}}(\mathbb{R}^+\times \mathbb{R}) },\\
			\label{bd:cos-right}
					&\|\pa_z^n\pa_r^m(\pa_r+r^{-1})\mathcal{F}_{k\to z}^{-1}\left( V^{-1} \mathcal{A}^2\cos (k\sqrt{\mathcal{A}}t)\big[\mathcal{A}^{-1}(\hat{v})\big]\right)\|_{L^{\infty}
						(\mathbb{R}^+\times \mathbb{R})}
					\lesssim  t^{-1}\|v\|_{W^{m+n+2,1}_{r^{-\delta}+r^{\delta}}
						(\mathbb{R}^+\times \mathbb{R})},\\
			\label{bd:sin-RL}
				&\|\pa_z^n\pa_r^m(\pa_r+r^{-1})\mathcal{F}_{k\to z}^{-1}\Big(V^{-1} \mathcal{A}^{\f32}\sin (k\sqrt{\mathcal{A}}t)\big[\hat{v}\big]\Big)\|_{L^{\infty}(\mathbb{R}^+\times \mathbb{R})}\lesssim t^{-1}\|v\|_{W^{m+n+1,1}_{r^{-\delta}+r^{\delta}}(\mathbb{R}^+\times \mathbb{R})}. \end{align}
			\end{proposition}
			
		To prove Proposition \ref{prop:solution-fm}, we first  introduce the following integral kernels.
		\begin{definition} \label{def:kernels}
			Let $c\in(0,1)$, $k\in\mathbb{R}/\{0\}$, $r,s\in\mathbb{R}^+$, and $\phi(r,c,k)$, $f_+(r,c,k)$, $W(c,k)$ be as in Definition \ref{def:useful-fm}. We define
			\begin{align}\label{keyK}
				\begin{split}
					& K(c,r,s,z):=p.v.\int_{\mathbb{R}}\f{\phi(r,c,k)
						\phi(s,c,k)}{\left|
						W(c,k)
						\right|^2}\f{e^{ikz}}{k}dk.
			\end{split}\end{align}
			Let $l',l\in\{0,1\}$, $m\in\mathbb{N}$. The kernel with oscillation for low/high $k-$frequency is defined as follows
			\begin{align}\label{keyK-LH}
				\begin{split}
					& K_{l',l,m}^L(c,r,s,z)\\
					&:=p.v.\int_{\mathbb{R}}\left(1-\chi\Big(\f{M^2 c^{\f12}}{k}\Big)\right)\f{
						\pa_r^m(\pa_r+r^{-1})^{l'}\phi(r,c,k)
						(\f{c^{\f12}(\pa_s+s^{-1})}{k})^{l}\phi(s,c,k)}{\left|
						W(c,k)
						\right|^2}\f{e^{ikz}}{k}dk,\\
					&K_{l',l,m}^H(c,r,s,z)\\
					&:=p.v.\int_{\mathbb{R}}\chi\Big(\f{M^2 c^{\f12}}{k}\Big)\f{
						\left(\f{c^{\f12}\pa_r}{k}\right)^m\left(\f{c^{\f12}(\pa_r+r^{-1})}{k}\right)^{l'}\phi(r,c,k)
						\left(\f{c^{\f12}(\pa_s+s^{-1})}{k}\right)^{l}\phi(s,c,k)}{\left|
						W(c,k)
						\right|^2}\f{e^{ikz}}{k}dk.
			\end{split}\end{align}
			\end{definition}
			
			\begin{proposition}\label{lem:Boundness-K}
		Let $\delta\ll 1$ be fixed, $l',l\in\{0,1\}$, $m\in\mathbb{N}$, $F\in\{L,H\}$.
			The following estimates are uniform in $(r,s,z)\in\mathbb{R}^+\times \mathbb{R}^+\times \mathbb{R}$:
			\begin{align}
				\int_0^1\left|\pa_cK(c,r,s,z)\right|dc\leq C\big(s^{-\delta}+s^{\delta}\big)\label{op1}.
			\end{align}
			For $l'+l+m\geq 1$, there holds
			\begin{align}
				\int_0^1\left|\pa_cK_{l',l,m}^{F}(c,r,s,z)\right|dc\leq C\big(s^{-\delta}+s^{\delta}\big)\label{op2}.
			\end{align}
			\end{proposition}
			
			Now we  prove  Proposition \ref{prop:solution-fm} by using Proposition \ref{lem:Boundness-K}.
			
			\begin{proof} [Proof of Proposition \ref{prop:solution-fm}.]  It suffices  to consider the case $k>0$.
				For \eqref{bd:sin}, we claim that  if $m$ is even, then
				\begin{align}  \label{keyfm-sin}
				\begin{split}
						&\pa_z^n\pa_r^m\mathcal{F}_{k\to z}^{-1}\left( V^{-1} \mathcal{A}^{\f32}\sin (k\sqrt{\mathcal{A}}t)\big[\hat{v}(r,k)\big]\right)\\&=\f{1}{2\pi^2 }\iint_{\mathbb{R}^+\times\mathbb{R}}\left(\int_0^1\f{\sin (c^{\f12}t)}{c^{\f12}}
							K_{0,0,m}^L(c,r,s,z-y)dc\right)\pa_y^nv(s,y)sdsdy \\
							&\quad+\f{1}{2\pi^2 }\iint_{\mathbb{R}^+\times\mathbb{R}}\left(\int_0^1\f{\sin (c^{\f12}t)}{c^{\f12}}K_{0,0,m}^H(c,r,s,z-y)dc\right) \pa_y^n\mathcal{A}^{-\f{m}{2}}v(s,y)sdsdy,	
							\end{split}				\end{align}
and if $m$ is odd, 				
\begin{align}  \label{keyfm-sin-odd}
\begin{split}
						&\pa_z^n\pa_r^m\mathcal{F}_{k\to z}^{-1}\left( V^{-1} \mathcal{A}^{\f32}\sin (k\sqrt{\mathcal{A}}t)\big[\hat{v}(r,k)\big]\right)\\&=\f{1}{2\pi^2 }\iint_{\mathbb{R}^+\times\mathbb{R}}\left(\int_0^1\f{\sin (c^{\f12}t)}{c^{\f12}}
							K_{0,0,m}^L(c,r,s,z-y)dc\right)\pa_y^nv(s,y)sdsdy\\
							&\quad+\f{1}{2\pi^2 }\iint_{\mathbb{R}^+\times\mathbb{R}}\left(\int_0^1\f{\sin (c^{\f12}t)}{c^{\f12}}\left(K_{0,1,m}^H,
							ic^{\f12}K_{0,0,m}^H\right)(c,r,s,z-y)dc\right)\\&\qquad\qquad\cdot\begin{pmatrix}
								\pa_s+s^{-1}\\
								\pa_y\end{pmatrix}\left(-V^{-1}\pa_y^n\mathcal{A}^{-\f{m-1}{2}}v\right)(s,y)sdsdy.
								\end{split}
						\end{align}
We integrate by parts in $c$ by noticing that
\begin{align} \label{integrate-by-parts-in-c}
					\f{\sin(c^{\f12}t)}{c^{\f12}}dc=-\f{2}{t}d\sin(c^{\f12}t),\quad
					\f{\cos(c^{\f12}t)}{c^{\f12}}dc=\f{2}{t}d\sin(c^{\f12}t),
\end{align}
and then we obtain \eqref{bd:sin} by Proposition \ref{lem:Boundness-K}.

				It remains to prove \eqref{keyfm-sin} and \eqref{keyfm-sin-odd}. We decompose
				\begin{align*} 
					\notag  &V^{-1} \mathcal{A}^{\f32}\sin (k\sqrt{\mathcal{A}}t) \big[\hat{v}(r,k)\big]\\
					&=k
					V^{-1}\mathcal{A}^{2}\big(1-\chi(M^2\sqrt{\mathcal{A}})\big)
					\f{\sin (k\sqrt{\mathcal{A}}t) }{k\sqrt{\mathcal{A}}}\big[\hat{v}\big] +k
					V^{-1}\mathcal{A}^{2}\chi(M^2\sqrt{\mathcal{A}}) \f{\sin (k\sqrt{\mathcal{A}}t) }{k\sqrt{\mathcal{A}}} \big[\hat{v}\big]\\ \notag &:=I^L+I^H.\notag \end{align*}
				For $I^L$,  we apply \eqref{spectal-formula-result-Fz} with $f(x)=\big(1-\chi(M^2\sqrt{x})\big)\f{\sin(k\sqrt{x}t)}{k\sqrt{x}}$.
				For $I^H$, we will decompose further to recover the loss while taking the derivative $\pa_r^m$. If $m$ is even,  using $1=\mathcal{A}^{\f{m}{2}}\mathcal{A}^{-\f{m}{2}}$, we decompose
				\begin{align} 
					\notag  I^H=k V^{-1}\mathcal{A}^{2}\chi(M^2\sqrt{\mathcal{A}}) \mathcal{A}^{\f{m}{2}}\f{\sin (k\sqrt{\mathcal{A}}t) }{k\sqrt{\mathcal{A}}} \mathcal{A}^{-\f{m}{2}}\hat{v}(r,k).\notag \end{align}
				Then we apply \eqref{spectal-formula-result-Fz} with $f(x)=\chi(M^2\sqrt{x})x^{\f{m}{2}}\f{\sin(k\sqrt{x}t)}{k\sqrt{x}}$,
				$v= \mathcal{A}^{-\f{m}{2}}v$.
				If $m$ is odd,  using $1
				=\mathcal{A}^{\f{m+1}{2}}\mathcal{A}^{-1}\mathcal{A}^{-\f{m-1}{2}}$ and the observation
				 \begin{align} \label{decom:A-odd} & \mathcal{A}^{-1}=-\Delta_1\circ V^{-1}=(\pa_s,ik)\cdot \begin{pmatrix}\pa_s+s^{-1}\\
						ik\end{pmatrix}\circ (-V^{-1}),
				\end{align}
				we decompose
				\begin{align}
					\notag  I^H=k V^{-1}\mathcal{A}^{2}\chi(M^2\sqrt{\mathcal{A}})\mathcal{A}^{\f{m+1}{2}} \f{\sin (k\sqrt{\mathcal{A}}t) }{k\sqrt{\mathcal{A}}} (\pa_r,ik)\cdot
					\begin{pmatrix}\pa_r+r^{-1}\\
						ik\end{pmatrix}\big(-V^{-1} \mathcal{A}^{-\f{m-1}{2}}\hat{u}_0\big).\notag
				\end{align}
				Then we apply \eqref{spectal-formula-result-vector-Fz} with $f(x)=\chi(M^2 x^{\f12})x^{\f{m+1}{2}}\f{\sin (k\sqrt{x}t)}{k\sqrt{x}}$, $\textbf{u}=\begin{pmatrix}
					\pa_s+s^{-1}\\
					ik\end{pmatrix}\big(-V^{-1} \mathcal{A}^{-\f{m-1}{2}}\hat{v}\big)$.

				 For \eqref{bd:cos-left}, we claim that if $m$ is even, then
				\begin{align}\label{keyfm-cos-right}
				\begin{split}
						&\pa_z^n\pa_r^m\Big(\mathcal{F}_{k\to z}^{-1}\Big(V^{-1} \mathcal{A}^2\cos (k\sqrt{\mathcal{A}t})\big[(\mathcal{A}^{-1}\hat{v})(r,k)\big]\Big)\\
						&=\f{1}{2\pi^2 }\iint_{\mathbb{R}^+\times\mathbb{R}}\int_0^1\f{\cos (c^{\f12}t)}{c^{\f12}}
							\Big(\small K_{0,1,m}^L,\;
							ic^{\f12}K_{0,0,m}^L\Big)(c,r,s,z-y)\\&\qquad\qquad\qquad\qquad\cdot
							\begin{pmatrix}\pa_s+s^{-1}\\
								\pa_y\end{pmatrix}\small  (-V^{-1}\pa_y^nv)(s,y)sdcdsdy\\
							&\quad+\f{1}{2\pi^2 }\iint_{\mathbb{R}^+\times\mathbb{R}}\int_0^1\f{\cos (c^{\f12}t)}{c^{\f12}}\Big(K_{0,1,m}^H,\;
							ic^{\f12}K_{0,0,m}^H\Big)(c,r,s,z-y)\\&\qquad\qquad\qquad\qquad\cdot
							\begin{pmatrix}\pa_s+s^{-1}\\
								\pa_y\end{pmatrix}\big(-V^{-1}\pa_y^n
							\mathcal{A}^{-\f{m}{2}}u_0\big)(s,y)sdcdsdy,		
							\end{split}
							\end{align}
				and if $m$ is odd, then
				  \begin{align}\label{keyfm-cos-right-odd}
				 \begin{split}
						&\pa_z^n\pa_r^m\Big(\mathcal{F}_{k\to z}^{-1}\Big(V^{-1} \mathcal{A}^2\cos (k\sqrt{\mathcal{A}t})\big[(\mathcal{A}^{-1}\hat{v})(r,k)\big]\Big)\\
						&= \f{1}{2\pi^2 }\iint_{\mathbb{R}^+\times\mathbb{R}}\int_0^1\f{\cos (c^{\f12}t)}{c^{\f12}}
							\Big(K_{0,1,m}^L,\;
							ic^{\f12}K_{0,0,m}^L\Big)(c,r,s,z-y)\\&\qquad\qquad\qquad\cdot
							\begin{pmatrix}\pa_s+s^{-1}\\
								\pa_y\end{pmatrix}(-V^{-1}\pa_y^nv)(s,y)sdcdsdy\\
							&\quad+ \f{1}{2\pi^2 }\iint_{\mathbb{R}^+\times\mathbb{R}}\int_0^1\f{\cos (c^{\f12}t)}{c^{\f12}}K_{0,0,m}^H(c,r,s,z-y)\pa_y^n\mathcal{A}^{-\f{m+1}{2}}v(s,y)sdcdsdy.
							\end{split}				\end{align}
	Then we can conclude \eqref{bd:cos-left} by integration by parts in $c$ by \eqref{integrate-by-parts-in-c} and Proposition \ref{lem:Boundness-K}.
	To prove \eqref{keyfm-cos-right} and \eqref{keyfm-cos-right-odd}, we decompose by using \eqref{decom:A-odd} that
				\begin{align*} 
					\notag  & V^{-1} \mathcal{A}\cos (k\sqrt{\mathcal{A}t})\hat{v}(r,k)\\
					&=k
					V^{-1}\mathcal{A}^{\f52}\left(1-\chi(\delta\sqrt{\mathcal{A}})\right)
					\f{\cos (k\sqrt{\mathcal{A}}t) }{k\sqrt{\mathcal{A}}} (\pa_r,ik)\cdot\begin{pmatrix}
						\pa_r+r^{-1}\\
						ik\end{pmatrix}(-V^{-1} \hat{v}(r,k))\\
					&\quad+k
					V^{-1}\mathcal{A}^{\f32}\chi(\delta\sqrt{\mathcal{A}}) \f{\cos (k\sqrt{\mathcal{A}}t) }{k\sqrt{\mathcal{A}}}\hat{v}(r,k):=I^L+I^H.\notag \end{align*}
				Then we apply \eqref{spectal-formula-result-vector-Fz}  on $I^L$ with $f(x)=(1-\chi(M^2\sqrt{x}))\f{\cos(k\sqrt{x}t)}{k\sqrt{x}}$,
				$\hat{u}=\begin{pmatrix}
					\pa_r+r^{-1}\\
					ik\end{pmatrix}(-V^{-1} \hat{v})$. For  $I^H$, we will decompose further to recover the loss while taking the derivative $\pa_r^m$. If $m$ is even,  using
				$\mathcal{A}^{\f32}=\mathcal{A}^{\f52} \mathcal{A}^{\f{m}{2}}
				\mathcal{A}^{-1}\mathcal{A}^{-\f{m}{2}}$ and \eqref{decom:A-odd}, then we decompose
				\begin{align} 
					\notag  I^H=k V^{-1}\mathcal{A}^{\f52}\chi(M^2\sqrt{\mathcal{A}}) \mathcal{A}^{\f{m}{2}}\f{\cos (k\sqrt{\mathcal{A}}t) }{k\sqrt{\mathcal{A}}}(\pa_r,ik)\cdot\begin{pmatrix}
						\pa_r+r^{-1}\\
						ik\end{pmatrix}(-V^{-1} \mathcal{A}^{-\f{m}{2}}\hat{v})
					,\notag \end{align}
				then we apply \eqref{spectal-formula-result-vector-Fz} with $f(x)=\chi(M^2\sqrt{x})x^{\f{m}{2}}\f{\cos(k\sqrt{x}t)}{k\sqrt{x}}$,
				$\hat{\bf{u}}=\begin{pmatrix}
					\pa_r+r^{-1}\\
					ik\end{pmatrix}(-V^{-1}  \mathcal{A}^{-\f{m}{2}}\hat{v})
				$. If $m$ is odd,  using  \eqref{decom:A-odd} and $\mathcal{A}^{\f32}=\mathcal{A}^{2}
				\mathcal{A}^{\f{m}{2}}\mathcal{A}^{-\f{m+1}{2}}$,
				then we decompose
				\begin{align}
					\notag  I^H=k V^{-1}\mathcal{A}^{2}\chi(M^2\sqrt{\mathcal{A}})\mathcal{A}^{\f{m}{2}} \f{\cos (k\sqrt{\mathcal{A}}t) }{k\sqrt{\mathcal{A}}} \big(\big(-\Delta_1\circ V^{-1}\big)^{\f{m+1}{2}}\hat{v}\big),\notag
				\end{align}
				then we apply \eqref{spectal-formula-result-Fz} with $f(x)=\chi(M^2\sqrt{x})x^{\f{m}{2}}\f{\cos(k\sqrt{x}t)}{k\sqrt{x}}$,
				$\hat{u}=\big(-\Delta_1\circ V^{-1}\big)^{\f{m+1}{2}}\hat{v}
				$.
				
				 In a similar way as in \eqref{keyfm-cos-right} and \eqref{keyfm-cos-right-odd}, we can prove \eqref{bd:cos-right} by using the following formulas: if $m$ is even, then
				 \begin{align}
						\notag &\pa_z^n\pa_r^m(\pa_r+r^{-1})\mathcal{F}_{k\to z}^{-1}\Big(V^{-1} \mathcal{A}^{2}\cos (k\sqrt{\mathcal{A}t})\big[\mathcal{A}^{-1}\hat{v}(r,k)\big]\Big)\\
						&=\f{1}{2\pi^2 }\iint_{\mathbb{R}^+\times\mathbb{R}}\int_0^1\f{\cos (c^{\f12}t)}{c^{\f12}}
							\Big(\small K_{1,1,m}^L,\;
							ic^{\f12}K_{1,0,m}^L\Big)(c,r,s,z-y)\\&\qquad\qquad\qquad\cdot
							\begin{pmatrix}\pa_s+s^{-1}\nonumber\\
								\pa_y\end{pmatrix}\small  (-V^{-1}\pa_y^nv)(s,y)sdcdsdy\\
							&\quad+ \f{1}{2\pi^2 }\iint_{\mathbb{R}^+\times\mathbb{R}}\int_0^1\f{\cos (c^{\f12}t)}{c^{\f12}}K_{1,0,m}^H(c,r,s,z-y)\pa_y^n\mathcal{A}^{-\f{m}{2}-1}v(s,y)sdcdsdy,\nonumber
					\end{align}
			and if $m$ is odd, then
				 \begin{align}
				 \begin{split}
						&\pa_z^n\pa_r^m(\pa_r+r^{-1})\mathcal{F}_{k\to z}^{-1}\Big(V^{-1} \mathcal{A}^{2}\cos (k\sqrt{\mathcal{A}t})\big[\mathcal{A}^{-1}\hat{v}(r,k)\big]\Big)\\
						&=\f{1}{2\pi^2 }\iint_{\mathbb{R}^+\times\mathbb{R}}\int_0^1\f{\cos (c^{\f12}t)}{c^{\f12}}
							\Big(K_{1,1,m}^L,\;
							ic^{\f12}K_{1,0,m}^L\Big)(c,r,s,z-y)\\&\qquad\qquad\qquad\cdot
							\begin{pmatrix}\pa_s+s^{-1}\\
								\pa_y\end{pmatrix}(-V^{-1}\pa_y^nv)(s,y)sdcdsdy\\
							&\quad+\f{1}{2\pi^2 }\iint_{\mathbb{R}^+\times\mathbb{R}}\int_0^1\f{\cos (c^{\f12}t)}{c^{\f12}}\Big(K_{1,1,m}^H,\;
							ic^{\f12}K_{1,0,m}^H\Big)(c,r,s,z-y)\\&\qquad\qquad\qquad\cdot
							\begin{pmatrix}\pa_s+s^{-1}\\
								\pa_y\end{pmatrix}\big(-V^{-1}\pa_y^n
							\mathcal{A}^{-\f{m+1}{2}}v\big)(s,y)sdcdsdy.
							\end{split}
					\end{align}

				 In a similar way as in \eqref{keyfm-sin} and \eqref{keyfm-sin-odd}, we can prove  \eqref{bd:sin-RL}  by using the following formulas: if $m$ is even, then				
				\begin{align}  
				\notag &\pa_z^n\pa_r^m(\pa_r+r^{-1})\mathcal{F}_{k\to z}^{-1}\Big( V^{-1} \mathcal{A}^{\f32}\sin (k\sqrt{\mathcal{A}}t)\big[\hat{v}(r,k)\big]\Big)\\&= \f{1}{2\pi^2 }\iint_{\mathbb{R}^+\times\mathbb{R}}\Big(\int_0^1\f{\sin (c^{\f12}t)}{c^{\f12}}
							K_{1,0,m}^L(c,r,s,z-y)dc\Big)\pa_y^nv(s,y)sdsdy \nonumber\\
							&\quad+\f{1}{2\pi^2 }\iint_{\mathbb{R}^+\times\mathbb{R}}\Big(\int_0^1\f{\sin (c^{\f12}t)}{c^{\f12}}\big(K_{1,1,m}^H,
							ic^{\f12}K_{1,0,m}^H\big)(c,r,s,z-y)dc\Big)\\&\qquad\qquad\qquad\cdot\begin{pmatrix}
								\pa_s+s^{-1}\nonumber\\
								\pa_y\end{pmatrix}\big(-V^{-1}\pa_y^n\mathcal{A}^{-\f{m}{2}}v\big)
							(s,y)sdsdy,\nonumber					\end{align}
				and if $m$ is odd, then
				 \begin{align} 
				  \notag &\pa_z^n\pa_r^m(\pa_r+r^{-1})\mathcal{F}_{k\to z}^{-1}\Big( V^{-1} \mathcal{A}^{\f32}\sin (k\sqrt{\mathcal{A}}t)\big[\hat{v}(r,k)\big]\Big)\\
				  &=\f{1}{2\pi^2 }\iint_{\mathbb{R}^+\times\mathbb{R}}\Big(\int_0^1\f{\sin (c^{\f12}t)}{c^{\f12}}
							K_{1,0,m}^L(c,r,s,z-y)dc\Big)\pa_y^nv(s,y)sdsdy\nonumber\\
							&\quad+\f{1}{2\pi^2 }\iint_{\mathbb{R}^+\times\mathbb{R}}\Big(\int_0^1\f{\sin (c^{\f12}t)}{c^{\f12}}K_{1,0,m}^H(c,r,s,z-y)dc\Big) \pa_y^n\mathcal{A}^{-\f{m+1}{2}}v(s,y)sdsdy.\nonumber					\end{align}

				This finishes the proof of the proposition.
				\end{proof}
				
				 \subsection{Revisiting the homogenous case}
				
				 Using the framework of this paper, we provide a simple proof of  asymptotic linear stability in the homogeneous case.
				
				\begin{proposition}\label{prop:homogenous-decay}
					Let $u\equiv\f12$. The solution of the linearized system \eqref{eq:vel-linearize-gernal0} admits the bound
					\begin{align}\label{homogenous-t-1-dacay}
						\|
						\pa_z^n\pa_r^m(u_r,u_{\theta},u_z)(t)\|_{L^{\infty}(\mathbb{R}^+\times\mathbb{R})}
						&\leq Ct^{-1}\|\textbf{u}_0\|_{W^{3+m+n
								,1
							}(\mathbb{R}^+\times\mathbb{R})}.
					\end{align}
				\end{proposition}
				\begin{proof}
					In this case, $2u=V=\Omega\equiv 1$, $\mathcal{A}=-\Delta^{-1}$, and there holds
					\begin{align*}
&u_{r}(t,r,z)=\mathcal{F}_{k\to z}^{-1}\left(\mathcal{A}\cos (k\sqrt{\mathcal{A}} t)\left(\mathcal{A}^{-1}(\hat{u}_r^0)\right)\right)
							+\mathcal{F}_{k\to z}^{-1}\left(V^{-1}\mathcal{A}^{\f32}\sin (k\sqrt{\mathcal{A}} t)\left(k\mathcal{A}^{-1}(\hat{u}_{\theta}^0)\right)\right), \\
&u_{z}(t,r,z)=(\pa_r+r^{-1})
							\mathcal{F}_{k\to z}^{-1}\left(\mathcal{A}^2\cos (k\sqrt{\mathcal{A}} t)\left(ik^{-1}\mathcal{A}^{-2}(\hat{u}_{r}^0)\right)\right)\\&\qquad\qquad\qquad+(\pa_r+r^{-1})\mathcal{F}_{k\to z}^{-1}\left(\mathcal{A}^{\f32}\sin (k\sqrt{\mathcal{A}} t)\left(i\mathcal{A}^{-1}(\hat{u}_{\theta}^0)\right)\right), \\ &u_{\theta}(t,r,z)
							=\f12\left(\mathcal{F}_{k\to z}^{-1}\left(\mathcal{A}\cos (k\sqrt{\mathcal{A}} t)\left(\mathcal{A}^{-1}(\hat{u}_{\theta}^0)\right)\right)\right)
							-\f12\left(\mathcal{F}_{k\to z}^{-1}\left(\mathcal{A}^{\f32}\sin (k\sqrt{\mathcal{A}} t)\left(k^{-1}\mathcal{A}^{-2}(\hat{u}_r^0)\right)\right)\right).
\end{align*}
As in the proof of Proposition \ref{prop:solution-fm}(without introducing cut-off functions to deal with low/high-frequency in $k$), we have
					\begin{align*} 
							\notag &\pa_z^n\pa_r^m\mathcal{F}_{k\to z}^{-1}\left( \mathcal{A}^{\f32}\sin (k\sqrt{\mathcal{A}}t)\hat{v}(r,k)\right)\\&=  \left\{
							\begin{array}{l}
								\f{1}{2\pi^2 }\iint_{\mathbb{R}^+\times\mathbb{R}}\left(\int_0^1\f{\sin (c^{\f12}t)}{c^{\f12}}K_{0,0,m}(c,r,s,z-y)dc\right) \pa_y^n\mathcal{A}^{-\f{m}{2}}v(s,y)sdsdy,\; m\;\text{even};\\ \f{1}{2\pi^2 }\iint_{\mathbb{R}^+\times\mathbb{R}}\left(\int_0^1\f{\sin (c^{\f12}t)}{c^{\f12}}\left(K_{0,1,m},
								ic^{\f12}K_{0,0,m}\right)(c,r,s,z-y)dc\right)\\\qquad\qquad\qquad\cdot\begin{pmatrix}
									\pa_s+s^{-1}\\
									\pa_y\end{pmatrix}\left(-\pa_y^n\mathcal{A}^{-\f{m-1}{2}}u_0\right)(s,y)sdsdy,\;m\;\text{odd};
							\end{array}\right.
						\end{align*}
						and
						\begin{align*}
							\notag&\pa_z^n\pa_r^m\mathcal{F}_{k\to z}^{-1}\left( \mathcal{A}\cos (k\sqrt{\mathcal{A}t})\hat{v}(r,k)\right)\\
							&=  \left\{
							\begin{array}{l}
								\f{1}{2\pi^2 }\iint_{\mathbb{R}^+\times\mathbb{R}}\int_0^1\f{\cos (c^{\f12}t)}{c^{\f12}}\left(K_{0,1,m}^H,\;
								ic^{\f12}K_{0,0,m}\right)(c,r,s,z-y)\\ \qquad\qquad\qquad\cdot
								\begin{pmatrix}\pa_s+s^{-1}\\
									\pa_y\end{pmatrix}\left(-\pa_y^n
								\mathcal{A}^{-\f{m}{2}}u_0\right)(s,y)sdcdsdy,\;m\;\text{even};\\
								\f{1}{2\pi^2 }\iint_{\mathbb{R}^+\times\mathbb{R}}\int_0^1\f{\cos (c^{\f12}t)}{c^{\f12}}K_{0,0,m}(c,r,s,z-y)\pa_y^n\mathcal{A}^{-\f{m+1}{2}}v(s,y)sdcdsdy,\;m\;\text{odd};
								\end{array}\right.
						\end{align*}
						and
						\begin{align*}
							\notag &\pa_z^n\pa_r^m(\pa_r+r^{-1})\mathcal{F}_{k\to z}^{-1}\left(\mathcal{A}^{2}\cos (k\sqrt{\mathcal{A}t})\hat{v}(r,k)\right)\\
							&=  \left\{
							\begin{array}{l} \f{1}{2\pi^2 }\iint_{\mathbb{R}^+\times\mathbb{R}}\int_0^1\f{\cos (c^{\f12}t)}{c^{\f12}}K_{1,0,m}(c,r,s,z-y)\mathcal{A}
								^{-\f{m}{2}}\pa_y^nv(s,y)sdcdsdy,\;\;m\;\text{even};\\
								\f{1}{2\pi^2 }\iint_{\mathbb{R}^+\times\mathbb{R}}\int_0^1\f{\cos (c^{\f12}t)}{c^{\f12}}\left(K_{1,1,m},\;
								ic^{\f12}K_{1,0,m}\right)(c,r,s,z-y)\\\qquad\qquad\qquad\cdot
								\begin{pmatrix}\pa_r+r^{-1}\\
									\pa_y\end{pmatrix}\left(-\pa_y^n\mathcal{A}^{-\f{m-1}{2}}\right)v(s,y)sdcdsdy,\;m\;\text{odd}; \end{array}\right.
						\end{align*}
						and
						\begin{align*}
							&\pa_y^n\pa_r^m(\pa_r+r^{-1}) \mathcal{F}_{k\to z}^{-1}\left(\mathcal{A}^{\f32}\sin (k\sqrt{\mathcal{A}t}) \hat{v}(r,k)\right)\\&=  \left\{
							\begin{array}{l}
								\f{1}{2\pi^2 }\iint_{\mathbb{R}^+\times\mathbb{R}}\int_0^1\f{\sin (c^{\f12}t)}{c^{\f12}}\left(K_{1,1,m},\;
								ic^{\f12}K_{1,0,m}\right)(c,r,s,z-y)\\\qquad\qquad\qquad\cdot
								\begin{pmatrix}\pa_s+s^{-1}  \\
									\pa_y\end{pmatrix}\left(-\pa_y^n\mathcal{A}^{-\f{m}{2}}u_0\right)(s,y)sdcdsdy,\;m\;\text{even};\\ \f{1}{2\pi^2 }\iint_{\mathbb{R}^+\times\mathbb{R}}\int_0^1\f{\sin (c^{\f12}t)}{c^{\f12}}K_{1,0,m}(c,r,s,z-y)
								\pa_y^n\mathcal{A}^{-\f{m+1}{2}}v(s,y)sdcdsdy,\;m\;\text{odd}.
							\end{array}\right.
						\end{align*}
						Here the integral kernel is defined by
					\begin{align}\label{keyK-homo}
						&K_{l',l,m}(c,r,s,z):=p.v.\int_{\mathbb{R}}\f{
							\left(\f{c^{\f12}\pa_r}{k}\right)^m\left(\f{c^{\f12}(\pa_r+r^{-1})}{k}\right)^{l'}\phi(r,c,k)
							\left(\f{c^{\f12}(\pa_s+s^{-1})}{k}\right)^{l}\phi(s,c,k)}{\left|
							W(c,k)
							\right|^2}\f{e^{ikz}}{k}dk,
					\end{align}
					where $\phi(r,c,k)\sim \xi r(r\to 0)$, $f_+(r,c,k)\sim (\xi r)^{-\f12}e^{i\xi r}(r\to +\infty)$, $\xi=k\sqrt{\f{1-c}{c}}$, $W(c,k)=\mathcal{W}(\phi,f_+)$.
					Then the decay estimate \eqref{homogenous-t-1-dacay} can be reduced to the following kernel estimate (Lemma \ref{lem:Boundness-trival}) via integration by parts  in $c$.
					\end{proof}
					
					\begin{lemma}\label{lem:Boundness-trival}
						Let $l',l\in\{0,1\}$, $m\in\mathbb{N}$.
						For the kernel defined in \eqref{keyK-homo}, there exists $C$ independent of $(r,s,z)\in\mathbb{R}^+\times \mathbb{R}^+\times \mathbb{R}$ such that
						\begin{align}
							\int_0^1\left|\pa_cK_{l',l,m}(r,s,z,c)\right|dc\leq C\label{op1-trival}.
						\end{align}
					\end{lemma}
					
					\begin{proof}
					In this case, we can write the functions $\phi,f_+,W$, explicitly: $$\phi(r,c,k)=J_1(\xi r),\quad f_+(r,c,k)=H_+(\xi r),\quad \xi=k\sqrt{\f{1-c}{c}},\quad W(c,k)=i,$$
					where $J_1(\cdot)$ is the Bessel function of the first kind, and $H_+(\cdot)=J_1(\cdot)+iY_1(\cdot)$ is the Hankel's function.
					Making the transform $\f{c^{\f12}}{k}=\f{(1-c)^{\f12}}{\xi}$ and $\f{e^{ikz}}{k}dk=
						\f{e^{i\xi x}}{\xi}d\xi$, $x=\sqrt{\f{c}{1-c}}z$,  the kernel can be rewritten as
					\begin{align}
						\notag &K_{l',l,m}(c,r,s,z)\\
						&=p.v.\int_{\mathbb{R}}
						\left(\f{(1-c)^{\f12}\pa_r}{\xi}\right)^m
						\left(\f{(1-c)^{\f12}(\pa_r+r^{-1})}{\xi}\right)^{i_1}J_1(\xi r)
						\left(\f{(1-c)^{\f12}(\pa_s+s^{-1})}{\xi}\right)^{i_2}J_1(\xi s)\f{e^{i\xi x(c,z) }}{\xi}d\xi,\nonumber
					\end{align}
					where $x(c,z)=\sqrt{\f{c}{1-c}}z$. Notice that					
					\begin{align}\nonumber
						&\pa_c\left((1-c)^{\f{l'+l+m}{2}}p.v.\int_{\mathbb{R}} f_1(\xi,r)f_2(\xi,s)\f{e^{i\xi x(c,z)}}{\xi}d\xi\right)\\
						\notag &=(1-c)^{\f{l'+l+m}{2}}\cdot i\f{\pa x(c,z)}{\pa c}\int_{\mathbb{R}} f_1(\xi,r)f_2(\xi,s)e^{i\xi x(c,z)}d\xi    \\
						&\quad+\pa_c\left((1-c)^{\f{l'+l+m}{2}}\right)p.v.\int_{\mathbb{R}} f_1(\xi,r)f_2(\xi,s)\f{ e^{i\xi x(c,z)}}{\xi}d\xi, \notag
					\end{align}
					where $f_1(\xi,r)=(\f{\pa_r}{\xi})^m
					(\f{(\pa_r+r^{-1})}{\xi})^{l'}J_1(\xi r)
					$ and $f_2(\xi,s)=
					(\f{(\pa_r+r^{-1})}{\xi})^{l}J_1(\xi s)
					$. By the Young's inequality of the convolution, we get
					\begin{align}
						\begin{split}\label{pifi-Linfty-homo} &\left\|\int_{\mathbb{R}}f_1(\xi)f_2(\xi)e^{i\xi x}d\xi\right\|
							_{L^{1}_x(\mathbb{R})}\lesssim \|\mathcal{F}^{-1}_{\xi\to x}(f_1)(x)\|_{L^1_x(\mathbb{R})}
							\|\mathcal{F}^{-1}_{\xi\to x}(f_2)(x)\|_{L^{1}_x(\mathbb{R})},\\
							&\left\|p.v.\int_{\mathbb{R}}f(\xi)\f{e^{i\xi x}}{\xi}d\xi\right\|
							_{L^{\infty}_x(\mathbb{R})}\lesssim \|\mathcal{F}^{-1}_{\xi\to x}(f)(x)\|_{L^1_x(\mathbb{R})},
						\end{split}
					\end{align}
					where we used $\mathcal{F}^{-1}_{\xi\to x}(p.v.(\f{1}{\xi}))=i\pi \mathrm{sgn}(x)$ for the second inequality. Therefore,
					we deduce
					\begin{align}\nonumber
						\begin{split}
							&\left\|\f{\pa x(c,z)}{\pa c}\int_{\mathbb{R}}f_1(\xi)f_2(\xi)e^{i\xi x(c,z)}d\xi\right\|
							_{L^{1}_{c}(0,1)L^{\infty}_{z}(\mathbb{R})}\\
							&= \left\|\int_{\mathbb{R}}f_1(\xi)f_2(\xi)e^{i\xi x}d\xi\right\|
							_{L^{1}_{x}(\mathbb{R})L^{\infty}_z(\mathbb{R})}
							\lesssim
							\prod_{i=1}^2\|\mathcal{F}^{-1}_{\xi\to x}\left(f_i(\xi)\right)(x)\|_{L^1_x(\mathbb{R})L^{\infty}_c(0,1)}.
						\end{split}
					\end{align}
					Thus, \eqref{op1-trival} can be reduced to showing that for $m\in\mathbb{N}$, $l'\in\{0,1\}$,
					\begin{align*}
						&p.v.\int_{\mathbb{R}}
						\left(\f{\pa_r}{\xi}\right)^m
						\left(\f{(\pa_r+r^{-1})}{\xi}\right)^{l'}J_1(\xi r)
						e^{i\xi x}d\xi\in L^1_x(\mathbb{R}),  \quad\f{d}{dc}\left((1-c)^{\f{m+l'}{2}}\right) \in L^1_c((0,1)).
					\end{align*}
					The second one is obvious, and the first one can be written by taking $\eta=\xi r$, $y=x/r$ as
					\begin{align*}
						&\int_{\mathbb{R}}
						\pa_{\eta}^m
						\left(\pa_{\eta}+\eta^{-1}\right)^{l'}J_1(\eta)
						e^{i\eta y}d\eta\in L^1_y(\mathbb{R}),
					\end{align*}
					which follows from the fact that
					\begin{align*}
						\int_{\mathbb{R}}J_1(\eta)e^{i\eta y}d\eta =\f{y}{\sqrt{1-y^2}} \textbf{1}_{(-1,1)}(y),\quad
						\int_{\mathbb{R}}\frac {J_1(\eta)} {\eta}e^{i\eta y}d\eta =\sqrt{1-y^2}\textbf{1}_{(-1,1)}(y).					
						\end{align*}
					\end{proof}
					
					\subsection{Reformulation of Proposition \ref{lem:Boundness-K}}
					
				Motivated by the proof of the homogeneous case, we reformulate the definition of  the kernels in Definition \ref{def:kernels} in terms of variable $\xi=k\sqrt{\f{1-c}{c}}$ as follows.
				
					\begin{definition}\label{def:kernels-new}
					 Let
					\begin{align}\label{def:tphi}
							\widetilde{\phi}(r,\xi,c):= \f{\phi(r,\xi,c)}{|W(\xi,c)|},
						\end{align}
			where for $c\in(0,1)$, $\xi\in\mathbb{R}/\{0\}$, $r\in\mathbb{R}^+$, $\phi(r,\xi,c)$ is a real-valued smooth solution of
						\begin{align} \nonumber 
							\phi''+r^{-1}\phi'-r^{-2}\phi+\f{\xi^2(V(r)-c)}{1-c}\phi=0,\;\;\phi(r,\xi,c)\sim \xi r\;\text{as}\;\;r\to 0,
							\end{align}
						while, $f_+(r,\xi,c)$ is a complex-valued smooth solution of
						\begin{align}\nonumber
							\begin{split}  
								&f_+''+r^{-1}f_+'-r^{-2}f_+
								+\f{\xi^2(V(r)-c)}{1-c}f_+=0\quad \text{with}\\
								&f_+(r,\xi,c)\sim \left\{
								\begin{aligned}
									&(\xi r)^{-\f12}e^{i\xi r},\;\quad\quad \quad\quad\quad |\xi|\lesssim (1-c)^{\f13},\;c\in(0,1)\\
									&(\xi r)^{-\f12}e^{i\xi\int_{0}^r
										\sqrt{\f{V(s)-c}{1-c}}ds},\;
									|\xi|\gtrsim (1-c)^{\f13},\;c\in(0,V(0)] \\
									&(\xi r)^{-\f12}e^{i\xi\int_{r_c}^r
										\sqrt{\f{V(s)-c}{1-c}}ds},\;
									|\xi|\gtrsim (1-c)^{\f13},\;c\in(V(0),1)
								\end{aligned}
								\right. \quad \text{as}\;\;r\to +\infty.\end{split}
						\end{align}
						The Wronskian is defined by
						\begin{align}\label{fm:wrons-op-new}\mathcal{W}(f,g)
							&:=r(f'g-fg')=(r^{\f12}f)'(r^{\f12}g)-(r^{\f12}g)'(r^{\f12}f), \\\label{fm:wrons-new}
						W(\xi,c)&:= \mathcal{W}(\phi,f_+)(\xi,c).
						\end{align}

					Let  $l',l\in\{0,1\}$, $m\in\mathbb{N}$ and
						\begin{align}  \label{def:derivative}
							\mathcal{D}_{l',m}:=\left(\f{(1-c)^{\f12}\pa_r}{\xi}\right)^{m}
							\left(\f{(1-c)^{\f12}(\pa_r+r^{-1})}{\xi}\right)^{l'},\quad
							\mathcal{D}_{l,0}^L:=\left(\f{(1-c)^{\f12}(\pa_s+s^{-1})}{\xi}\right)^{l}.
							\end{align}						
					Then we define
						\begin{align}\label{keyK-new}
							\begin{split}
								K(r,s,z,c)&:=p.v.\int_{\mathbb{R}} \widetilde{\phi}(r,\xi,c)
								\widetilde{\phi}(s,\xi,c)\f{e^{i\xi \sqrt{\f{c}{1-c}}z}}{\xi}d\xi,
								\end{split}
						\end{align}
						and
						\begin{align}\label{keyK-L-new}
							 K_{l',l,m}^L(r,s,z,c)
							:=p.v.\int_{\mathbb{R}}&\left(1-\chi\left(\f{M^2 (1-c)^{\f12}}{\xi}\right)\right)\\
							&\left(\mathcal{D}_{l',m}^L
							\widetilde{\phi}\right)(r,\xi,c)
							\left(\mathcal{D}_{l,0}^L
							\widetilde{\phi}\right)(s,\xi,c)
							\f{e^{i\xi \sqrt{\f{c}{1-c}}z}}{\xi}d\xi,\nonumber
						\end{align}
						\begin{align} \label{keyK-H-new}
							K_{l',l,m}^H(r,s,z,c)
							:=\int_{\mathbb{R}}\chi\left(\f{M^2 (1-c)^{\f12}}{\xi}\right)
							\left(\mathcal{D}_{l',m}
							\widetilde{\phi}\right)(r,\xi,c)
							\left(\mathcal{D}_{l,0}
							\widetilde{\phi}\right)(s,\xi,c)
							\f{e^{i\xi \sqrt{\f{c}{1-c}}z}}{\xi}d\xi.
						\end{align} 				
						\end{definition}

					The kernel estimates in Proposition \ref{lem:Boundness-K} can be reduced to
					
					\begin{proposition}\label{prop:Boundness-K-D}
						Let $\delta\ll 1$ be fixed,  $l',l\in\{0,1\}$, $m\in\mathbb{N}$.
						For the kernels \eqref{keyK-new}, \eqref{keyK-L-new}, \eqref{keyK-H-new} in  Definition \ref{def:kernels-new}, the following estimates hold uniformly in $(r,s,z)\in\mathbb{R}^+\times \mathbb{R}^+\times \mathbb{R}$:
						\begin{align}
							\int_0^1\left|\pa_cK(r,s,z,c)\right|dc\leq C\big(s^{-\delta}+s^{\delta}\big)\label{op-K},
						\end{align}
						and for $l'+l+m\geq 1$, $F\in\{L,H\}$,
						\begin{align}
							\int_0^1\left|\pa_cK_{l',l,m}^{F}(r,s,z,c)\right|dc\leq C\big(s^{-\delta}+s^{\delta}\big).
							\label{op-K-D}
					\end{align}
					\end{proposition}

The proof of Proposition \ref{prop:Boundness-K-D} is deferred to other parts of this paper.
The key ingredient lies in analyzing the quantitative behavior with respect to $(r,\xi,c)$ of $\phi(r,\xi,c)$, $f_+(r,\xi,c)$ and $W(\xi,c)$ as defined in Definition \ref{def:kernels-new}.
					\section{Behavior of Schr\"odinger equation near the origin}
					
					In this section, we consider the Schr\"odinger equation with the prescribed behavior at the origin:
					\begin{align}  \label{eq:phi}
							\phi''+r^{-1}\phi'-r^{-2}\phi+\f{\xi^2(V(r)-c)}{1-c}\phi=0,\;\;\phi(r,\xi,c)\sim \xi r\;\text{as}\;\;r\to 0.
					\end{align}					
					
				Let's begin with a technical lemma.
					
					\begin{lemma}\label{lem:vorteraa-0}  Let the parameters $(\xi,c)\in D$, $r_1(\xi,c)\geq r_0(\xi,c)\geq 0$, $r\in I_1=(r_0,r_1)$,
						$i,j\in\mathbb{N}$, $l\in\{0,1\}$. Let $f(r,\xi,c)$ solve the integral equation
						\begin{align}\label{eq:vorteraa}
							f(r,\xi,c)=g(r,\xi,c)+\int_{r_0}^r K(r,s,\xi,c)f(s,\xi,c)ds,
						\end{align}
						with $K(r,r,\xi,c)=0$ and $g(r,\xi,c)\neq 0$. Let
						\begin{align}\label{def:K1}
							K_1(r,s,\xi,c):= g^{-1}(r,\xi,c)K(r,s,\xi,c)g(s,\xi,c).
						\end{align}
						\begin{itemize}
							
							\item [(1)]
							If there exists a bounded positive function $\kappa(r,\xi,c)$ such that
							\begin{align}
								\int_{r_0}^r
								\sup_{\tilde{r}\in[s,r]}\left|
								K_1(\tilde{r},s,\xi,c)
								\right|ds\leq C
								\kappa(r,\xi,c),\;r\in I_1,\label{condition:K'}
							\end{align}
							then we have
							\begin{align*}
								f(r,\xi,c)&=g(r,\xi,c)\left(1+
								f^{Rem}(r,\xi,c)\right),
							\end{align*}
							with the remainder term satisfying the following bound: for $r\in I_1$ \begin{align}\label{solution-asy-0r}
								|f^{Rem}(r,\xi,c)|\lesssim \kappa(r,\xi,c). \end{align}
							\item
							[(2)]
							If there exist functions $\lambda^r(r,c)$ and  $\lambda^c(r,c,\xi)$ such that for $r\in (2r_0(\xi,c),r_1(\xi,c)]$, $i,j\in\mathbb{N}$, $l\in\{0,1\}$, 
							$|\left(\lambda^c(r,c)\pa_c\right)^l(\xi\pa_{\xi})^j r_0(\xi,c)|  \lesssim r_0(\xi,c)$, 
							and
							\begin{align}\nonumber \int_{0}^1 \sup_{\tilde{t}\in[t,1]}\Big|
								(\lambda^c(r,c))^{l}&(\lambda^r(r,c))^{i}\pa_c^{l}\pa_r^{i}
								(\xi\pa_{\xi})^{j} \\&\Big(
								(r-r_0(\xi,c))K_1(r_0+\tilde{t}(r-r_0),
								r_0+t(r-r_0),\xi,c)\Big)
								\Big|dt\lesssim
								\kappa(r,\xi,c)
								,\label{condition:K'pac-growth}
								\end{align}
							then   it holds that for $r\in( 2r_0(\xi,c),\;r_1(\xi,c)]$,
							\begin{align*}
								|(\lambda^c(r,c))^l(\lambda^r(r,c))^i\pa_c^l\pa_r^i(\xi\pa_{\xi})^j
								f^{Rem}(r,\xi,c)|\lesssim 
								\kappa(r,\xi,c).
							\end{align*}
							In particular, if $r_0=0$,\;$\lambda^r=r$, 
							and  \begin{align}
								\int_{0}^1 \sup_{\tilde{t}\in[t,1]}\left|
								(\lambda^c\pa_c)^{l'}
								(r\pa_{r})^{i_1} (\xi\pa_{\xi})^{j_1}\left(rK_1(\tilde{t}r,tr,\xi,c)  \right)
								\right|dt\lesssim
								\kappa(r,\xi,c),\label{condition:K'pac}
							\end{align}
							then it holds that for $r\in (0,r_1(\xi,c)]$,
							\begin{align}\label{fRes1}
								\left|
								(\lambda^c\pa_c)^l(r\pa_{r})^i(\xi\pa_{\xi})^jf^{Rem}(r,\xi,c)\right|\lesssim \kappa(r,\xi,c).
							\end{align}
							\item    [(3)]
							If
							\begin{align}\label{integral:phi-sgn}
								\begin{split}
									&K_1(r,s,\xi,c)>0,\;\pa_rK_1(r,s,\xi,c)>0,\;\;
									\text{for}\;\;r_0\leq s\leq r\leq r_1,\end{split}
							\end{align}
							then   we have $$f^{Rem}(r,\xi,c)>0,\;\pa_rf^{Rem}(r,\xi,c)>0\quad \text{for} \,\,r\in (r_0,r_1).$$
						\end{itemize}
					\end{lemma}
					
					\begin{proof}
						We construct an infinite series that solves \eqref{eq:vorteraa}:
						\begin{align*}
							f(r,\xi,c)=g(r,\xi,c)\Big(1+\sum_{n=1}^{\infty}h_n(r,\xi,c)\Big):=
							g\big(1+f^{Rem}\big),
						\end{align*}
						where $$h_1(r,\xi,c)=\int_{r_0}^rK_1(r,s,\xi,c)
						ds,\quad
						h_n(r,\xi,c)=\int_{r_0}^rK_1(r,s,\xi,c)h_{n-1}(s,\xi,c)ds (
						n\geq 2).$$
						It follows from the above inductive relation that
						\begin{align}\label{fm:hn}
							&\notag\small h_n(r,\xi,c)\\
							&=\int_{r_0}^rK_1(r,s_1,\xi,c)\int_{r_0}^{s_1}
							K_1(s_1,s_2,\xi,c)...\int_{r_0}^{s_{n-1}}
							K_1(s_{n-1},s_n,\xi,c)ds_n...ds_2ds_1.
						\end{align}
						
						Our aim is to show the convergence of $\sum_{n=1}^{\infty}h_n$. 
						Letting $r\in I_1$ be fixed, we define \begin{align*}\varpi(s,r,\xi,c)=\sup_{\tilde{r}\in[s,r]}
							\left|K_1(\tilde{r},s,\xi,c)\right|.
						\end{align*}
						The condition \eqref{condition:K'} implies
						\begin{align*}
							\int_{r_0}^r\varpi(s,r,\xi,c)ds\leq C\kappa(r,\xi,c)\leq C.
						\end{align*}
						Then we get  by \eqref{fm:hn} that
						\begin{align*}
							\notag|h_n(r,\xi,c)|&\leq \int_{r_0}^r|K(r,s_1,\xi,c)|\int_{r_0}^{s_1}
							|K(s_1,s_2,\xi,c)|...\int_{r_0}^{s_{n-1}}
							|K(s_{n-1},s_n,\xi,c)|ds_n...ds_2ds_1\\
							&\leq \int_{r_0}^r\varpi(s_1,r,\xi,c)\int_{r_0}^{s_1}
							\varpi(s_2,r,\xi,c)...\int_{r_0}^{s_{n-1}}
							\varpi(s_n,r,\xi,c)ds_n...ds_2ds_1\\
							&= \f{\left(\int_{r_0}^r\varpi(s,r,\xi,c)ds\right)^n}{n!}\leq \f{C^n\kappa^n(r,\xi,c)}{n!}\leq \f{C^n}{n!} ,
						\end{align*}
						which along with $f^{Rem}=\sum_{n=1}^{\infty}h_n$
						gives \eqref{solution-asy-0r}.
						
						 To prove \eqref{condition:K'pac-growth}, for fixed $r\in I_1$ and
						$i,j\in \mathbb{N}$, $l\in\{0,1\}$, we define
						\begin{align}\notag
							\varpi_{i,j,l}(t,r,\xi,c):=&\sum_{
									0\leq l'\leq l,\;0\leq i_1\leq i,\;0\leq j_1\leq j}
								\sup_{\tilde{t}\in[t,1]}
								\left|\lambda^c(r,c)^{l'} \lambda(r,c)^{i_1} \pa_c^{l'}
								\pa_{r}^{i_1}(\xi\pa_{\xi})^{j_1}\right.\\
								&\quad \left.\left(
								(r-r_0(\xi,c))K_1(r_0+\tilde{t}(r-r_0),r_0
								+t(r-r_0),\xi,c)\right)\right|.
							\end{align}
							By \eqref{condition:K'pac}, we have $\int_{0}^1\varpi_{i,j,l}(t,r,\xi,c)dt\leq C_{i,j,l} \kappa(r,\xi,c)$. We rewrite \eqref{fm:hn} as
							\begin{align}
								&\notag\small h_n(r,\xi,c)\\
								\label{fm:hn-change}&=\int_{0}^1\int_{0}^{t_1}...
								\int_{0}^{t_{n-1}}(r-r_0(\xi,c))\left(K_1(r,(r-r_0)t_1+r_0,\xi,c) \right) \cdot\\
								\notag&\quad\left((r-r_0(\xi,c))K_1((r-r_0)t_1+r_0,(r-r_0)t_2+r_0,\xi,c)\right)
								\cdot\\
								&\notag\quad...
								\cdot\left((r-r_0(\xi,c))K_1((r-r_0)t_{n-1}+r_0,(r-r_0)t_n+r_0,\xi,c)\right)dt_n...dt_2dt_1.
							\end{align}
							We apply the differential operator $(\lambda^c\pa_c)^{l'}
							(\lambda^r\pa_{r})^{i_1}(\xi\pa_{\xi})^{j_1}$ to the aforementioned formula, employing the Leibniz rule. Then there exist constants $C_{i,j,l}$ independent of $n$ such that
							\begin{align*}
								&\left|(\lambda^c\pa_c)^l(\lambda^r\pa_{r})^i(\xi\pa_{\xi})^jh_n(r,\xi,c)\right|\\
								&\leq C_{i,j,l}
								C^n(r-r_0)^n
								\int_{0}^1\int_{0}^{t_1}...\int_{0}^{t_{n-1}}
								\varpi_{i,j,l}(t_1)
								\varpi_{i,j,l}(t_2)... \varpi_{i,j,l}(t_n)dt_n...dt_2dt_1\\
								&=\f{C_{i,j,l}
									\left(\int_{0}^1\varpi_{i,j,l}(t,r,\xi,c)dt\right)^n}{n!}
								\leq \f{C_{i,j,l}^{n+1}
									\kappa^n(r,\xi,c)}{n!},\notag
							\end{align*}
							which along with $f^{Rem}=\sum_{n=1}^{\infty}h_n$
							gives \eqref{condition:K'pac-growth}.
							
							The proof of \eqref{fRes1} is simpler since $r_0=0$. Finally, the positivity \eqref{integral:phi-sgn} follows directly from \eqref{fm:hn}.
					\end{proof}

					The following lemmas characterize the asymptotic behavior of $\phi$ in the different regimes of the parameters $(c,\xi)$.

					\begin{lemma}\label{lem:phi-sqrt(1-c)/xi}
						Let $c\in(0,1)$, $\xi>0$. The solution $\phi$ of \eqref{eq:phi} may write as
						\begin{align}
							\begin{split}
								&\phi(r,c,\xi)=\xi r\left(1+\phi^{Rem}(r,\xi,c)\right)\;\;\text{for}\;\;r\lesssim \f{\sqrt{1-c}}{\xi},\quad \text{with}\\
								&\quad\left|\left((1-c)\pa_c\right)
								^l(r\pa_r)^i(\xi\pa_{\xi})^j\phi^{Rem}(r,\xi,c)\right|\lesssim \f{\xi^2r^2}{1-c},\;i,j\in\mathbb{N},l\in\{0,1\}.\label{serious:phi}
							\end{split}
						\end{align}
						Moreover, if $c\in(V(0),1)$, $r_c=V^{-1}(c)$, it holds that
						\begin{align}\label{sgn:phi1}
							\phi^{Rem}(r),\quad \pa_r\phi^{Rem}(r)>0\;\;\text{for}\;
							\;r\lesssim \f{\sqrt{1-c}}{\xi}\;\text{and}\;\;r<r_c.
						\end{align}
					\end{lemma}
					\begin{proof}
						We first notice that \eqref{eq:phi}  is a linear second order ODE, which has a unique smooth solution for given initial data. Our aim is to describe the precise behavior of the solution.
						For this, we rewrite the equation as
						\begin{align}\label{Op:L0}
							\phi''+ r^{-1}\phi'-r^{-2}\phi=\f{\xi^2}{1-c}(c-V)\phi,\;\phi(r,\xi,c)\sim \xi r\;\;(r\to 0),
						\end{align}
						which is equivalent to the integral equation
						\begin{align}
							\label{phi:1inter}
							\begin{split}
								\phi(r,\xi,c)&=\xi r+\int_0^rK(r,s,\xi,c)\phi(s,\xi,c)ds,\\
								K(r,s,\xi,c)&:=\f12\left(r-\f{s^2}{r}\right)\f{\xi^2}{1-c}(c-V(s))\;\;\text{with}\;K(s,s,\xi,c)=0.
							\end{split}
						\end{align}
						We apply Lemma \ref{lem:vorteraa-0} to \eqref{phi:1inter} with $D=\big\{(\xi,c)|c\in(0,1),\xi> 0\big\}$, $I_1=\Big\{r|0< r\lesssim \f{\sqrt{1-c}}{|\xi|}\Big\}$, $f=\phi$, $g=\xi r$ and $K_1$ in \eqref{def:K1}:
						\begin{align}\label{def:tK-phi1}
							K_1(r,s,\xi,c)=\f{\xi^2}{2(1-c)}\left(1-\f{s^2}{r^2}\right)s(c-V(s)).
						\end{align}
						It remains to verify conditions \eqref{condition:K'} and \eqref{condition:K'pac}. For \eqref{condition:K'} on $K_1$, taking $\kappa=\f{\xi^2r^2}{1-c}$,
						it holds that for $r\lesssim \f{\sqrt{1-c}}{|\xi|}$,
						\begin{align}\label{est:Kerphi}
							&\int_{0}^r
							\sup_{\tilde{r}\in[s,r]}\left|
							\f{\xi^2}{2(1-c)}\left(1-\f{s^2}{\tilde{r}^2}\right)
							s(c-V(s))\xi s\right|ds\leq \int_{0}^r\f{\xi^2s|V(s)-c|}{1-c}ds\\&\leq \int_{0}^r\f{\xi^2s}{1-c}ds\sim\f{\xi^2r^2}{1-c}.\nonumber
						\end{align}
						For \eqref{condition:K'pac}, taking $\lambda^c=1-c$, 
						it follows that for $i,j\in\mathbb{N},l\in\{0,1\}$ and  $r\lesssim \f{\sqrt{1-c}}{|\xi|}$,
						\begin{align*}
							&\int_{0}^1 \sup_{\tilde{t}\in[t,1]}\left|
							\left((1-c)\pa_c\right)^{l}(r \pa_{r})^{i} \left(\xi\pa_{\xi}\right)^{j}
							\left(\f{\xi^2r^2}{1-c}\left(1-\f{t^2}{\tilde{t}^2}
							\right)t(c-V(tr))\right)\right|dt\\&\lesssim \int_{0}^1\f{\xi^2r^2t}{1-c}dt\lesssim \f{\xi^2r^2}{1-c}.
						\end{align*}
						
						Finally, \eqref{sgn:phi1} follows by applying \eqref{integral:phi-sgn}. Indeed,  $K_1(r,s,\xi,c)> 0$ and $\pa_rK_1(r,s,\xi,c)> 0$ for $0\leq s<r<r_c$.
					\end{proof}
					
					\begin{lemma}\label{lem:phi-1*}
						Let $\delta\ll 1$ be fixed,  $c\in(0,1)$ and $\xi\lesssim (1-c)^{\f12}$. 
						The solution $\phi$ of \eqref{eq:phi} may write as $$\phi(r,c,\xi)=J_1(\xi r)\left(1+\phi^{Rem}(r,\xi,c)\right)$$ with the remainder term satisfying
						\begin{align}   \begin{split}
								\label{serious:phi-1*}
								&\quad 0<\phi^{Rem}(r,\xi,c) \lesssim \f{\xi^2}{1-c}, \quad r\leq \delta\xi^{-1}. \end{split}
						\end{align}
					\end{lemma}
					\begin{proof}
					We rewrite
					\eqref{eq:phi} as
					\begin{align}
						\phi''+r^{-1}\phi'-r^{-2}\phi+\xi^2\phi=\f{\xi^2(1-V)}{1-c}\phi,\;\phi\sim \xi r(r\to 0),
					\end{align}
					which has the equivalent integral form
					\begin{align}
						\label{phi:1inter*}
						\begin{split}
							&\phi(r,\xi,c)=J_1(\xi r)+\int_0^rK(r,s,\xi,c)\phi(s,\xi,c)ds,\\
							&K(r,s,\xi,c):=\f{s\pi}{2}\left(J_1(\xi r)Y_1(\xi s)-J_1(\xi s)Y_1(\xi r)\right)\f{\xi^2(1-V(s))}{1-c}\;\;\text{with}\;K(s,s,\xi,c)=0.
						\end{split}
					\end{align}
					We denote
					\begin{align}\label{def:tK-phi*}
						K_1(r,s,\xi,c)&:=J_1(\xi r)^{-1}\f{s\pi}{2}\left(J_1(\xi r)Y_1(\xi s)-J_1(\xi s)Y_1(\xi r)\right)
						J_1(\xi s)\f{\xi^2(1-V(s))}{1-c}. \end{align}
					Thanks to $J_1'(z)Y_1(z)-J_1(z)Y_1'(z)=\f{1}{\pi z}$, we know that $\f{J_1(z)}{Y_1(z)}$ is monotonically increasing for $z>0$.
					This along with $J_1(z)>0$ for $z<\delta\ll 1$ and $1-V>0$ gives
					\begin{align}\label{def:tK1>0}
						K_1(r,s,\xi,c)>0,\;\;\text{if}\;\; s<r\leq \delta \xi^{-1}.\end{align}
					Then we apply Lemma \ref{lem:vorteraa-0}
					to \eqref{phi:1inter*} with $D=\big\{(\xi,c)|c\in(0,1),\; \xi\lesssim
					(1-c)^{\f12}$, $I_1=\{r|0< r\leq \delta \xi^{-1}\big\}$, $f=\phi$, $g=J_1(\xi r)$ and  $K_1$ as in \eqref{phi:1inter*}.
					We first verify the condition \eqref{condition:K'} on $K_1$. Indeed, taking $\kappa=\f{\xi^2}{1-c}$, it holds that for $r\lesssim \delta \xi^{-1}$,
					\begin{align*} &\quad\int_{0}^r
						\sup_{\tilde{r}\in[s,r]}\left|
						K_1(\tilde{r},s,\xi,c)\right|ds\leq \int_{0}^r\f{\xi^2s|1-V(s)|}{1-c}ds\leq \int_{0}^r\f{\xi^2s\langle s\rangle^{-3}}{1-c}ds\lesssim \f{\xi^2}{1-c}.
					\end{align*}
					This shows $\phi^{Rem}(r,\xi,c) \lesssim \f{\xi^2}{1-c}$.  Finally, we can verify  \eqref{integral:phi-sgn} with \eqref{def:tK-phi*} to obtain $\phi^{Rem}(r)>0$ for $r\leq \delta \xi^{-1}$.
					\end{proof}
					
					\begin{lemma}\label{lem:phi-1.5} 						Let $c\in(0,1)$, $\xi\lesssim (1-c)^{\f12}$. 
						The  solution $\phi$ of \eqref{eq:phi} may write as$$\phi(r,c,\xi)=\xi r\left(1+\phi^{Rem}(r,\xi,c)\right)$$
						with the remainder term satisfying
						\begin{align}   \begin{split}
								\label{serious:phi(1-c)1/2}
								&\left|\left((1-c)\pa_c\right)
								^l(r\pa_r)^i(\xi\pa_{\xi})^j\phi^{Rem}(r,\xi,c)\right|
								\\&\quad\lesssim \xi^2r^2+\f{\xi^2}{1-c}\left(
								r^2\mathbf{1}_{r\lesssim 1}(r)+ \mathbf{1}_{1\lesssim r\lesssim \xi^{-1} }(r)\right) \quad\text{for}\,\, r\lesssim  \xi^{-1},
							\end{split}
						\end{align}
						where $i,j\in\mathbb{N},l\in\{0,1\}$. Moreover,
						the remainder also admits the bound
						\begin{align} \label{serious:phi(1-c)1/2-2}
							\begin{split}
								\left|\left((1-c)\pa_c\right)
								^l\pa_r^i(\xi\pa_{\xi})^j\phi^{Rem}(r,\xi,c)\right|
								\lesssim \f{\xi^2}{1-c}\quad\text{for}\quad r\lesssim 1. \end{split}
						\end{align}
					\end{lemma}
					\begin{proof}
						We notice that in this case, we have $1\lesssim \xi^{-1}$. We rewrite \eqref{eq:phi} as \eqref{Op:L0} and use the integral form \eqref{phi:1inter}.
						Then we apply Lemma \ref{lem:vorteraa-0} to \eqref{phi:1inter} with $D=\{(\xi,c)|c\in(0,1),\; \xi\lesssim (1-c)^{\f12}$, $I_1=\big\{r|0< r\lesssim \xi^{-1}\big\}$, $f=\phi$, $g=\xi r$ and  $K_1$ as in \eqref{def:tK-phi1}:
						\begin{align*}
						K_1(r,s,\xi,c)=\f{\xi^2}{2(1-c)}\left(1-\f{s^2}{r^2}\right)s(c-V(s)).
						\end{align*}
						To obtain the desired results, we need to verify conditions \eqref{condition:K'} and \eqref{condition:K'pac}.
						A key observation in this case is that  for $1\lesssim \xi\lesssim |V(0)-c|^{-\f32}$ and $s\lesssim \xi^{-\f23}$,
						\begin{align}\label{ineq:c-V3}
							|c-V(s)|\leq |1-c|+|1-V(s)|\lesssim |1-c|
							+\langle s\rangle^{-3}.
						\end{align}
						To verify \eqref{condition:K'}, we take $\kappa(r,\xi,c)=\xi^2r^2+\f{\xi^2}{1-c}\left(r^2\mathbf{1}_{r\lesssim 1}(r)+
						r^{-1}\mathbf{1}_{r\gtrsim 1 }(r)\right)$ which is bounded for $r\lesssim \xi^{-1}$. We  bound \eqref{est:Kerphi} by using \eqref{ineq:c-V3} as
						\begin{align*}
							\int_{0}^r
							\sup_{\tilde{r}\in[s,r]}\left|
							\f{\xi^2}{2(1-c)}\left(1-\f{s^2}{\tilde{r}^2}\right)s(c-V(s))\right|ds
							&\lesssim \int_{0}^r\xi^2s+\f{\xi^2}{1-c}\f{s}{\langle s\rangle^3}ds \\
							&\sim  \xi^2r^2+\f{\xi^2}{1-c}\left(
							r^2\mathbf{1}_{r\lesssim 1}(r)+
							\mathbf{1}_{r\gtrsim 1 }(r)\right)
							\notag.
						\end{align*}
						To verify the condition \eqref{condition:K'pac}, we observe that for $s=rt, i\geq 1$,
						\begin{align*}
							|(1-c)\pa_c(c-V(rt))|\leq 1-c,\;
							\left|(r\pa_r)^i\left(c-V(rt)\right)\right|\lesssim \langle rt\rangle^{-3}, \end{align*}
						which along with \eqref{ineq:c-V3}  gives
						\begin{align*}
							&\int_{0}^1 \sup_{\tilde{t}\in[t,1]}\left|
							((1-c)\pa_{c})^{l'}(r \pa_{r})^{i_1}(\xi\pa_{\xi})^{j_1}
							\left(\f{\xi^2r^2}{2(1-c)}\left(1-\f{t^2}{\tilde{t}^2}\right) t(c-V(rt))\right)
							\right|dt\\
							&\lesssim\int_{0}^1 \f{\xi^2r^2}{1-c} t\left((1-c)+\langle rt\rangle^{-3}\right)
							dt\lesssim \xi^2r^2+\f{\xi^2}{1-c}\left(
							r^2\mathbf{1}_{r\lesssim 1}(r)+ \mathbf{1}_{r\gtrsim 1 }(r)\right).
						\end{align*}
						
						Finally, we verify the condition \eqref{condition:K'pac-growth} with $\kappa, \lambda^r\equiv 1$, $\lambda^c=1-c$. We observe that for $r\lesssim 1, i\geq 1$,
						 \begin{align*}
							|((1-c)\pa_{c})^{l}\pa_{r}^{i}\left(r^2(c-V(rt))\right)| \lesssim 1.
						\end{align*}
						Then it follows that for  $r\lesssim 1$, $\xi\lesssim (1-c)^{\f12}$,
						\begin{align*}
							&\int_{0}^1 \sup_{\tilde{t}\in[t,1]}\left|
							((1-c)\pa_{c})^{l}\pa_{r}^{i}(\xi\pa_{\xi})^{j}
							\left(\f{\xi^2r^2}{2(1-c)}\left(1-\f{t^2}{\tilde{t}^2}\right) t(c-V(rt))\right)
							\right|dt\lesssim \f{\xi^2}{1-c}\int_{0}^1  t
							dt\lesssim 1.
						\end{align*}
					\end{proof}
					
					\begin{lemma}\label{lem: phi-2}
						Let $c\in(0,1-\delta)\setminus \{V(0)\}$ for some fixed $\delta\in(0,1)$, $1\lesssim \xi\lesssim |V(0)-c|^{-\f32}$, and $i,j\in\mathbb{N},l\in\{0,1\}$.
						Then the solution $\phi$ of \eqref{eq:phi} may write as
						\begin{align}
							\begin{split}\label{behave:phi2}
								&\phi(r,\xi,c)=\xi r\left(1+\phi^{Rem}(r,\xi,c))\right)\quad \text{for}\;\;r\lesssim \xi^{-\f23}\quad\text{with}\\
								&\quad\left|\left(\xi^{-\f23}\pa_c\right)^l(r\pa_r)^i(\xi\pa_{\xi})^j\phi^{Rem}(r,\xi,c)\right|\lesssim \xi^{\f43} r^2 .
							\end{split}
						\end{align}
						The remainder also admits the bound: for $r\lesssim \xi^{-\f23}$,
						\begin{align}
							\begin{split}\label{behave:phi2-2} &\quad\left|\left(\xi^{-\f23}\pa_c\right)^l(\xi^{-\f23}
								\pa_r)^i(\xi\pa_{\xi})^j\phi^{Rem}(r,\xi,c)\right|\lesssim 1.
							\end{split}
						\end{align}
					\end{lemma}
					
					\begin{proof}
						We rewrite \eqref{eq:phi} as \eqref{Op:L0} and use the integral form \eqref{phi:1inter}. Then we apply Lemma \ref{lem:vorteraa-0} to \eqref{phi:1inter} with
						$D=\big\{(\xi,c)|c\in(0,1-\delta),\;1\lesssim \xi\lesssim |V(0)-c|^{-\f32}\big\}$, $I_1=\{r|0< r\lesssim \xi^{-\f23}\}$, $f=\phi$, $g=\xi r$ and  $K_1$ as in \eqref{def:tK-phi1}: \begin{align*}
							K_1(r,s,\xi,c)=\f{\xi^2}{2(1-c)}\left(1-\f{s^2}{r^2}\right)s(c-V(s)).
						\end{align*}
						
						It suffices to verify conditions \eqref{condition:K'} and \eqref{condition:K'pac} for $K_1$. The key observation in this case is that  for $1\lesssim \xi\lesssim |V(0)-c|^{-\f32}$ and $s\lesssim \xi^{-\f23}$,
						\begin{align}\label{ineq:c-V1}
							|c-V(s)|\leq |c-V(0)|+|V(0)-V(s)|\lesssim |c-V(0)|+s\lesssim \xi^{-\f23}.
						\end{align}
						To verify \eqref{condition:K'}, we take $\kappa(r,\xi,c)=\xi^{\f43}r^2$ and bound \eqref{est:Kerphi} by using \eqref{ineq:c-V1} as
						\begin{align*}
							\int_{0}^r
							\sup_{\tilde{r}\in[s,r]}\left|
							\f{\xi^2}{2(1-c)}\left(1-\f{s^2}{\tilde{r}^2}\right)s(c-V(s))\right|ds\lesssim \int_{0}^r\xi^2s|c-V(s)|ds\lesssim \int_{0}^r\xi^{\f43}sds\sim\xi^{\f43}r^2\notag.
						\end{align*}
						We observe that for $c\in(0,1-\delta)$, $r\lesssim \xi^{-\f23}$, $i\geq 1$, $t\in(0,1)$,
						\begin{align*}
							\left|\xi^{-\f23}\pa_c\left(\f{c-V(rt)}{1-c}\right)\right|+
							\left|(r\pa_r)^i\left(c-V(rt)\right)\right|\lesssim \xi^{-\f23}.
						\end{align*}
						Then we can verify the condition \eqref{condition:K'pac} by taking $\lambda^c=\xi^{-\f23}$:
						\begin{align*}
							&\int_{0}^1 \sup_{\tilde{t}\in[t,1]}\left|
							(\xi^{-\f23}\pa_{c})^{l}(r
							\pa_{r})^{i}(\xi\pa_{\xi})^{j}
							\left(\f{\xi^2r^2}{2(1-c)}\left(1-\f{t^2}{\tilde{t}^2}\right)t (c-V(rt))\right)\right|dt\lesssim \int_{0}^1\xi^{\f43}r^2tdt\lesssim \xi^{\f43}r^2.
						\end{align*}
						
						By using \eqref{ineq:c-V1}, we also observe that for $c\in(0,1-\delta)$, $r\lesssim \xi^{-\f23}$, $t\in(0,1)$, $l=0,1$, $i\geq 1$,
						\begin{align*}
							\left|(\xi^{-\f23}\pa_r)^i\left(\f{\xi^2 r^2(c-V(rt))}{1-c}\right)\right| \lesssim \xi^2(\xi^{-\f23})^3\lesssim 1. 
						\end{align*}
						Then we can verify the condition \eqref{condition:K'pac-growth} with $\kappa=1, \lambda^r=\lambda^c=\xi^{-\f23}$ that
						\begin{align*}
							&\quad\int_{0}^1 \sup_{\tilde{t}\in[t,1]}\left|
							(\xi^{-\f23}\pa_{c})^{l}(\xi^{-\f23}
							\pa_{r})^{i}(\xi\pa_{\xi})^{j}
							\left(\f{\xi^2r^2(c-V(rt))}{2(1-c)}\left(1-\f{t^2}{\tilde{t}^2}\right)t \right)\right|dt\lesssim \int_{0}^1
							tdt
							\sim 1.
						\end{align*}
					\end{proof}
					
					\begin{lemma}\label{lem:phi-3}
						Let $c\in(0,1)\setminus\{V(0)\}$, $\xi\gtrsim  (1-c)^{\f12}|c-V(0)|^{-\f32}$ and $i,j\in\mathbb{N},l\in\{0,1\}$. Then the solution $\phi$ of \eqref{eq:phi} may write as
						\begin{align}
							\begin{split}\label{behave:phi3}
								&\phi(r,\xi,c)=\xi r\left(1+\phi^{Rem}(r,\xi,c)\right)\;\;\text{for}\;\;r\lesssim \xi^{-1}(1-c)^{\f12}|c-V(0)|^{-\f12}\quad\text{with}\\
								&\quad\left|\left((1-c)|c-V(0)|\pa_c\right)^l(r\pa_r)^i(\xi\pa_{\xi})^j
								\phi^{Rem}(r,\xi,c)\right|\lesssim \f{|c-V(0)|\xi^2r^2}{1-c}.
							\end{split}
						\end{align}
						The remainder also admits the bound: for $r\lesssim \xi^{-1}(1-c)^{\f12}|c-V(0)|^{-\f12}$,
						\begin{align}\label{behave:phi3-2} &\quad\left|\left((1-c)|c-V(0)|\pa_c\right)^l
							(\xi^{-1}(1-c)^{\f12}|c-V(0)|^{-\f12}\pa_r)^i(\xi\pa_{\xi})^j
							\phi^{Rem}(r,\xi,c)\right|\lesssim 1.  \end{align}
						Moreover, if $c\in(V(0),1)$, $r_c=V^{-1}(c)$, it holds that
						\begin{align}\label{sgn:phi9}
							\phi^{Rem}(r,\xi,c),\;\pa_r\phi^{Rem}(r,\xi,c)>0,\;r\lesssim \xi^{-1}(1-c)^{\f12}|c-V(0)|^{-\f12}\;\text{and}\;\;r<r_c.
						\end{align}
					\end{lemma}
					\begin{proof}
					 We rewrite \eqref{eq:phi} as \eqref{Op:L0} and use the integral form \eqref{phi:1inter}. Then we apply Lemma \ref{lem:vorteraa-0} to \eqref{phi:1inter} with
					 $D=\big\{(\xi,c)|c\in(0,1-\delta),\; \xi\gtrsim|V(0)-c|^{-\f32}\big\}$, $I_1=\{r|0< r\lesssim \xi^{-1}|V(0)-c|^{-\f12}\}$, $f=\phi$, $g=\xi r$ and  $K_1$ as in \eqref{def:tK-phi1}:
					 \begin{align*}
					 	K_1(r,s,\xi,c)=\f{\xi^2}{2(1-c)}\left(1-\f{s^2}{r^2}\right)s(c-V(s)).
					 \end{align*}
					
					Next, we verify conditions \eqref{condition:K'}, \eqref{condition:K'pac} and \eqref{integral:phi-sgn}. The key observation in this case is that for $ \xi\gtrsim(1-c)^{\f12}|c-V(0)|^{-\f32}$,
					it holds that
					 \begin{align*}
					 	\xi^{-1}(1-c)^{\f12}|c-V(0)|^{-\f12}\lesssim |c-V(0)|.
					 \end{align*}
					As in \eqref{ineq:c-V1},  for $s\leq r\lesssim \xi^{-1}\left(\f{|c-V(0)|}{1-c}\right)^{-\f12}$, we have
					 \begin{align}\label{ineq:c-V2}
					 	|c-V(s)|\leq |c-V(0)|+|V(0)-V(s)|\lesssim |c-V(0)|+s\lesssim |c-V(0)|.
					 \end{align}
					 To verify \eqref{condition:K'}, we take $\kappa(r,\xi,c)=|V(0)-c|\xi^{2}r^2$ and bound \eqref{est:Kerphi} by \eqref{ineq:c-V2} as
					 \begin{align*}
					 	&\quad\int_{0}^r
					 	\sup_{\tilde{r}\in[s,r]}\left|
					 	\left(\left(1-\f{s^2}{\tilde{r}^2}\right)
					 	\cdot \f{\xi^2(c-V(s))s}{1-c}\right)\right|ds\lesssim \int_{0}^r\f{\xi^2|c-V(0)|s}{1-c}ds\sim \f{\xi^2|c-V(0)|r^2}{1-c}\notag.
					 \end{align*}
					 Using \eqref{ineq:c-V2} with $s=rt$, we observe that
					 \begin{align*}
					 	\left|(1-c)|c-V(0)|\pa_c\left(\f{c-V(rt)}{1-c}\right)\right|\lesssim \f{|c-V(0)|}{1-c},\;|(r\pa_r)^i(c-V(rt))|\lesssim  tr\lesssim |c-V(0)|.
					 \end{align*}
					 Then we verify the condition \eqref{condition:K'pac} with $\lambda^c=(1-c)|c-V(0)|$. Indeed, for $r\lesssim \xi^{-1}(1-c)^{\f12}|c-V(0)|^{\f12}$,
					 \begin{align*}
					 	&\int_{0}^1 \sup_{\tilde{t}\in[t,1]}\left|
					 	\left((1-c)|c-V(0)|\pa_c\right)^{l}(r \pa_{r})^{i}\left(\xi\pa_{\xi}\right)^{j}
					 	\left(\left(1-\f{t^2}{\tilde{t}^2}\right)
					 	\cdot
					 	\f{\xi^2r^2(c-V(rt))t}{1-c}\right)\right|dt\\
					 	&\lesssim \int_{0}^1\f{\xi^2|c-V(0)|r^2t}{1-c}dt\sim\f{\xi^2|c-V(0)|r^2}{1-c}.
					 \end{align*}
					
					 Using \eqref{ineq:c-V2}, we also observe that for $r\lesssim \xi^{-1}(1-c)^{\f12}|c-V(0)|^{-\f12}$, $i\geq 1$,
					 \begin{align*}
					 	\left|\left(\xi^{-1}(1-c)^{\f12}|c-V(0)|^{-\f12}\right)^i
					 	\pa_r^i(r^2(c-V(rt)))\right|\lesssim   \f{1-c}{\xi^2}.
					 \end{align*}
					 Then we verify the condition \eqref{condition:K'pac-growth} with $\lambda^r=\xi^{-1}(1-c)^{\f12}|c-V(0)|^{-\f12}, \lambda^c=(1-c)|c-V(0)|$. Indeed,
					  for $r\lesssim \xi^{-1}(1-c)^{\f12}|c-V(0)|^{-\f12}$,
					 \begin{align*}
					 	&\int_{0}^1 \sup_{\tilde{t}\in[t,1]}\left|
					 	\left((1-c)|c-V(0)|\pa_c\right)^{l}\left(\xi^{-1}\Big(\f{1-c}{|c-V(0)|}\Big)^{\f12}
					 	\pa_{r}\right)^{i}\left(\xi\pa_{\xi}\right)^{j}
					 	\left(\left(1-\f{t^2}{\tilde{t}^2}\right)
					 	\cdot
					 	\f{\xi^2r^2(c-V(rt))t}{1-c}\right)\right|dt\\
					 	&\quad\lesssim \int_{0}^1tdt\sim 1.
					 \end{align*}
					
					 For the positivity of $\phi^{Rem}$, it suffices to verify the condition \eqref{integral:phi-sgn}.
					 Indeed, by the definition, $K_1(r,s,\xi,c)> 0$ and $\pa_rK_1(r,s,\xi,c)=\f{\xi^2(c-V(s))}{1-c}\f{s^3}{r^3}>0$ for $0\leq s<r<r_c$.
					\end{proof}
					
					  \begin{lemma}\label{lem:phi-4.5} 						
					  Let $c\in(0,V(0))$, $\xi\gtrsim(V(0)-c)^{-\f32}$, $i,j\in\mathbb{N},l\in\{0,1\}$. 
						The solution $\phi$ of \eqref{eq:phi} may write as: for $ r\leq \xi^{-1}(V(0)-c)^{-\f12}$,
						\begin{align}
							\notag&\phi(r,c,\xi)=2\sqrt{\f{1-c}{V(0)-c}} J_1\left(\xi\sqrt{\f{V(0)-c}{1-c}} r\right)\left(1+\phi^{Rem}(r,\xi,c))\right) \quad \text{with}\\
							& \left|\left((V(0)-c)\pa_c\right)^l(r\pa_r)^i(\xi\pa_{\xi})^j
							\phi^{Rem}(r,\xi,c)\right|\lesssim \xi^2r^3 ,\label{serious:phi-4.5} \\
							&\left|\left((V(0)-c)\pa_c\right)^l(\xi^{-1}(V(0)-c)^{-\f12}
							\pa_r)^i(\xi\pa_{\xi})^j
							\phi^{Rem}(r,\xi,c)\right|\lesssim \xi^{-1}(V(0)-c)^{-\f32}\label{serious:phi-4.5-2}  .
						\end{align}
					\end{lemma}
					\begin{proof}
						We rewrite the equation \eqref{eq:phi} as
						\begin{align}\label{Op:L0-5}
							\phi''+ r^{-1}\phi'-r^{-2}\phi+\f{(V(0)-c)\xi^2}{1-c}
							\phi=\f{\xi^2(V(0)-V(r))}{1-c}\phi,\;\phi(r,\xi,c)\sim \xi r(r\to 0).
							\end{align}
							Letting $\tilde{\xi}=\xi\sqrt{\f{V(0)-c}{1-c}}$, the equation is equivalent to the integral form
							\begin{align*}
								\begin{split}
									&\phi(r,\xi,c)=2\sqrt{\f{1-c}{V(0)-c}} J_1\left(\xi\sqrt{\f{V(0)-c}{1-c}} r\right)J_1\left(\tilde{\xi} r\right)+\int_0^rK(r,s,\xi,c)\phi(s,\xi,c)ds,\\
									&K(r,s,\xi,c):=\f{s\pi}{2}\left(J_1\left(\tilde{\xi} r\right)Y_1\left(\tilde{\xi} s\right)-J_1\left(\tilde{\xi} s\right)Y_1\left(\tilde{\xi} r\right)\right)\f{\xi^2(V(0)-V(s))}{1-c}\;\;\text{with}\;K(s,s,\xi,c)=0.
								\end{split}
							\end{align*}
							We denote
							\begin{align}\label{def:tK-phi4}
								\notag K_1(r,s,\xi,c)&=J_1\left(\tilde{\xi} r\right)^{-1}\f{s\pi}{2}\left(J_1(\tilde{\xi}r)Y_1(\tilde{\xi} s)-J_1(\tilde{\xi} s)Y_1(\tilde{\xi} r)\right)
								\f{\xi^2(V(0)-V(s))}{1-c} J_1(\tilde{\xi} s)\\
								&=\left(J_1(\tilde{\xi} s)Y_1(\tilde{\xi}s)-\f{s^2}{r^2}\f{J_1^2(\tilde{\xi} s)}{(\tilde{\xi} s)^2}\f{(\tilde{\xi}r)^2Y_1(\tilde{\xi} r)}{J_1(\tilde{\xi} r)}\right)\f{\pi\xi^2(V(0)-V(s))s}{2(1-c)}.
								\end{align}
								Then we apply Lemma \ref{lem:vorteraa-0} to \eqref{phi:1inter} with $D=\big\{(\xi,c)|c\in(0,V(0)),\xi\gtrsim 1\big\}$,
								$I_1=\{r|0< r\lesssim \xi^{-1}(V(0)-c)^{-\f12}\}$, $f=\phi$, $g=2\sqrt{\f{1-c}{V(0)-c}} J_1(\xi\sqrt{\f{V(0)-c}{1-c}} r)$.
								To prove \eqref{serious:phi-4.5}, we verify the conditions \eqref{condition:K'} and  \eqref{condition:K'pac}  for $K_1$, with $\kappa=\xi^2r^3\lesssim 1 $.
								For $i=j=l=0$, we verify the condition \eqref{condition:K'} for $r\leq \xi^{-1}(V(0)-c)^{-\f12}(\lesssim 1)$ that
								\begin{align*}
									\notag&\int_{0}^r
									\sup_{\tilde{r}\in[s,r]}\left|
									\left((J_1Y_1)(\tilde{\xi}s)-\f{s^2}{\tilde{r}^2}\f{J_1^2(\tilde{\xi} s)}{(\tilde{\xi} s)^2}\f{(\tilde{\xi}\tilde{r})^2Y_1(\tilde{\xi} \tilde{r})}{J_1(\tilde{\xi} \tilde{r})}\right)\f{\pi\xi^2(V(0)-V(s))s}{2(1-c)}\right|ds\\
									&\lesssim \int^{r}_0\xi^2s^2ds\lesssim \xi^2r^3.
								\end{align*}
								Here we used
								\begin{align}\label{ineq:c-V4}|V(0)-V(s)|\lesssim s,\quad
									\;s\lesssim 1,
								\end{align}
								and
								\begin{align}\label{ineq:c-V4-Bessel}
									\left|\left(z\f{d}{dz}\right)^m\left((J _1Y_1)(z)\right)\right|+\left|\left(z\f{d}{dz}\right)^m\left(\f{J_1(z)}{z}\right)\right|
									+\left|\left(z\f{d}{dz}\right)^m\left(\f{z^2Y_1(z)}{J_1(z)}\right)\right|\lesssim 1,\;\text{for}\;|z|\lesssim 1,\;m\geq 0.
									\end{align}
								For $i+j+l\geq 1$, we verify the condition \eqref{condition:K'pac} for $r\leq \xi^{-1}(V(0)-c)^{-\f12}(\lesssim 1)$ that
								\begin{align*}
									\notag\quad&\int_{0}^1 \sup_{\tilde{t}\in[t,1]}\left|
									\left((V(0)-c)\pa_c\right)^l(r\pa_r)^i(\xi\pa_{\xi})^j\right.\\
									&\qquad\left.\left(
									\f{\pi\xi^2(V(0)-V(rt))r^2t}{2(1-c)}\left((J_1Y_1)(\tilde{\xi}rt)-\f{t^2}{\tilde{t}^2}\f{J_1^2(\tilde{\xi} rt)}{(\tilde{\xi} rt)^2}\f{(\tilde{\xi}r\tilde{t})^2Y_1(\tilde{\xi} r\tilde{t})}{J_1(\tilde{\xi} r\tilde{t})}\right)\right)\right|dt\\
									&\lesssim \int^{1}_0\xi^2r^3dt
									\sim \xi^2r^3, \end{align*}
								where we used \eqref{ineq:c-V4-Bessel} and
								\begin{align*}
									\left|(r\pa_r)^i\left(\f{\xi^2(V(0)-V(rt))r^2}{1-c}\right)\right|\lesssim \xi^2r^3,  \;r\lesssim 1,\,0\leq t\leq 1.
									\end{align*}
									
									To prove \eqref{serious:phi-4.5-2}, we verify the condition \eqref{condition:K'pac-growth} for $K_1$, with $\kappa=\xi^{-1}(V(0)-c)^{-\f32}(\lesssim 1) $. Indeed, for $r\leq \xi^{-1}(V(0)-c)^{-\f12}(\lesssim 1)$,
									\begin{align*}
										\notag\quad&\int_{0}^1 \sup_{\tilde{t}\in[t,1]}\left|
										\left((V(0)-c)\pa_c\right)^l(\xi^{-1}(V(0)-c)^{-\f12}\pa_r)^i(\xi\pa_{\xi})^j\right.\\
										&\qquad\quad\left.\left(
										\f{\pi\xi^2(V(0)-V(rt))r^2t}{2(1-c)}\left((J_1Y_1)(\tilde{\xi}rt)-\f{t^2}{\tilde{t}^2}\f{J_1^2(\tilde{\xi} rt)}{(\tilde{\xi} rt)^2}\f{(\tilde{\xi}r\tilde{t})^2Y_1(\tilde{\xi} r\tilde{t})}{J_1(\tilde{\xi} r\tilde{t})}\right)\right)\right|dt\\
										&\lesssim \int^{1}_0\xi^{-1}(V(0)-c)^{-\f32}tdt
										\sim \xi^{-1}(V(0)-c)^{-\f32}, \end{align*}
									where we used
									\begin{align}\label{ineq:c-V4-Bessel-2}
										\left|\left(\f{d}{dz}\right)^m\left((J _1Y_1)(z)\right)\right|+\left|\left(\f{d}{dz}\right)^m\left(\f{J_1(z)}{z}\right)\right|
										+\left|\left(\f{d}{dz}\right)^m\left(\f{z^2Y_1(z)}{J_1(z)}\right)\right|\lesssim 1,\;\text{for}\;|z|\lesssim 1,\;m\geq 0, \end{align}
									$\xi^{-1}(V(0)-c)^{-\f12}\tilde{\xi}=(1-c)^{-\f12}\lesssim 1$
									, \eqref{ineq:c-V4} and for $r\lesssim \xi^{-1}(V(0)-c)^{-\f12},\,0\leq t\leq 1,\,i\geq 1$,
									\begin{align*}
										\left|\left(\xi^{-1}(V(0)-c)^{-\f12}\pa_r\right)^i\left(\f{\xi^2(V(0)-V(rt))r^2}{1-c}\right)\right|\lesssim \xi^{-1}(V(0)-c)^{-\f32}. \end{align*}
					\end{proof}

					\section{Behavior of Schr\"odinger equation at the infinity}
					
					In this section, we consider the Schr\"odinger equation with the prescribed behavior at the infinity:   					
					\begin{align}
							\begin{split}  \label{eq:f+}
								&f_+''+r^{-1}f_+'-r^{-2}f_+
								+\f{\xi^2(V(r)-c)}{1-c}f_+=0\quad \text{with}\\
								&f_+(r,\xi,c)\sim \left\{
								\begin{aligned}
									&(\xi r)^{-\f12}e^{i\xi r},\;\quad\quad \quad\quad\quad |\xi|\lesssim (1-c)^{\f13},\;c\in(0,1)\\
									&(\xi r)^{-\f12}e^{i\xi\int_{0}^r
										\sqrt{\f{V(s)-c}{1-c}}ds},\;
									|\xi|\gtrsim (1-c)^{\f13},\;c\in(0,V(0)] \\
									&(\xi r)^{-\f12}e^{i\xi\int_{r_c}^r
										\sqrt{\f{V(s)-c}{1-c}}ds},\;
									|\xi|\gtrsim (1-c)^{\f13},\;c\in(V(0),1)
								\end{aligned}
								\right. \quad \text{as}\;\;r\to +\infty.\end{split}
						\end{align}

					\subsection{Volterra type integral equation}	
					
									We again start with a technical lemma solving the Volterra type integral equation with the oscillation kernel.
					\begin{lemma}\label{lem:vorterra-infty}
						Let $i,j\in\mathbb{N}$, $l\in\{0,1\}$, and the parameters $(\xi,c)\in D$, $r\in I_1=(r_0(\xi,c),+\infty)$. Let $f(r,\xi,c)$ solve the integral equation
						\begin{align*}
							f(r,\xi,c)=g(r,\xi,c)+\int_{r}^{+\infty}K(r,s,\xi,c)f(s,\xi,c)ds,
						\end{align*}
						where $K(r,r,\xi,c)=K(r,+\infty,\xi,c)=0$ and $g(r,\xi,c)\neq 0$.
						\begin{itemize}
							\item[(1)]
							If there exists a positive bounded function $\kappa(r,\xi,c)$
							such that for
							\begin{align*}K_1(r,s,\xi,c)=
								g^{-1}(r,\xi,c)
								K(r,s,\xi,c)g(s,\xi,c),\end{align*}
							it holds
							\begin{align}
								&\int_r^{+\infty}
								\sup_{\tilde{r}\in[r,s]}\left|
								K_1(\tilde{r},s,\xi,c)\right|ds\leq C\kappa(r,\xi,c)
								,\label{condition:Kinfty'}\;\text{for}\;\;r\in I_1,
							\end{align}
							then the integral equation 
							admits an unique solution  \begin{align*}f(r,\xi,c)=g(r,\xi,c)\left(1+
								f^{Rem}(r,\xi,c)\right),\end{align*} where the remainder term admits the following bound  on $r\in I_1$:
							\begin{align*}
								|f^{Rem}(r,\xi,c)|\lesssim \kappa(r,\xi,c).
							\end{align*}
							\item [(2) ]  
							Assume  that there exists a smooth invertible transform
							$(r,s,\xi,c)\to \left(x(r,\xi,c),y(s,\xi,c),\xi,c\right)$ satisfying $y(s,\xi,c)|_{s=r}=x(r,\xi,c),\;\lim_{s\to+\infty}y(s,\xi,c)
							=+\infty,\;y_s\neq 0$. We denote the kernel under this transform as
							\begin{align}\label{def:tK1}
								\widetilde{K}_1(x,y,\xi,c)=K_1(r,s,\xi,c)
								y_s^{-1}(s,\xi,c).
							\end{align}
							Let 
							$\lambda^r(r,c)$,\;$\lambda^c(r,c)$ be the weight functions. Assume that for $r\in I_{1}$,
							\begin{align} &
								\xi \pa_{\xi}x=x,\quad
								\left|(\lambda^c)^l(\lambda^r)^i\pa_c^l\pa_r^ix\right|\lesssim |x|,\label{integral2-trans-condition2}
							\end{align}
						 \begin{align}
								\label{condition-trans-tK}
								&\left((\lambda^c)^lx^i\right)(r,\xi,c)\sup_{\tilde{x}\in[x,x+\tilde{y}]}
								\left|\pa_c^l\pa_{\tilde{x}}^i(\xi\pa_{\xi})^j\tilde{K}_1(\tilde{x},
								\tilde{x}+\tilde{y},\xi,c)\right|\lesssim
								C_{i,j,l}\varpi(x+\tilde{y},r,\xi,c)\\
								\notag &\text{with}\;\int_{0}^{+\infty}
								\varpi(x+\tilde{y},r,\xi,c)d\tilde{y}\lesssim \kappa(r,\xi,c).
							\end{align}
							Then the remainder term admits the following derivative estimates:  for $r\in I_{1}$, 							
							\begin{align*}
								\left|(\lambda^c)^l(\lambda^r)^i\pa_c^l\pa_{r}^i(\xi\pa_{\xi})^j
								f^{Rem}(r,\xi,c\right|\lesssim \kappa(r,\xi,c).
							\end{align*}
						\end{itemize}
					\end{lemma}
					
					\begin{proof}
						The proof of 
						(1) is identical as the proof of Lemma \ref{lem:vorteraa-0}. We will prove the derivative bounds under the condition \eqref{condition-trans-tK}.
						Here $x=x(r,\xi,c)$ is a technical  transform to deal with the oscillation in the kernel $K_1$. First, we construct the same infinite series that
						$f(r,\xi,c)=g(r,\xi,c)\left(1+f^{Rem}(r,\xi,c)\right)$, where $f^{Rem}(r,\xi,c)=\sum_{n=1}^{\infty}h_n(r,\xi,c)$ with
						\begin{align*}
							\begin{split}
								h_n(r,\xi,c)&=\int_r^{+\infty}K_1(r,s_1,\xi,c)\int_{s_1}^{+\infty}
								K_1(s_1,s_2,\xi,c)...\int_{s_{n-1}}^{+\infty}
								K_1(s_{n-1},s_n,\xi,c)ds_n...ds_2ds_1.
							\end{split}
						\end{align*}
						Considering the transform $\left(r,s_1,...s_n,\xi,c\right) \to \left(x(r,\xi,c),y_1(s_1,\xi,c),...(y_n(s_n,\xi,c),\xi,c\right)$
						satisfying  \eqref{integral2-trans-condition2}. By the definition of $\widetilde{K}_1$, we can rewrite $h_n$ as
						\begin{align}
							\notag
							h_n(r,\xi,c)&=\tilde{h}_n(x(r,\xi,c),\xi,c)\\ &:=\int_{x}^{+\infty}\widetilde{K}_1(x,y_1,\xi,c)\int_{y_1}^{+\infty}
							\widetilde{K}_1(y_1,y_2,\xi,c)...\int_{y_{n-1}}^{+\infty}\widetilde{K}_1(y_{n-1},y_n,\xi,c)dy_n...dy_2dy_1.\label{eq:t-hn}
						\end{align}
						We claim that under the condition \eqref{integral2-trans-condition2}, for $j,i\in \mathbb{N}$, $l\in\{0,1\}$,
						\begin{align}\label{bd:t-hn}
							&\notag \left| (\lambda^c(r,c))^{l}(\lambda^r(r,c))^i\pa_c^l\pa_r^{i}(\xi\pa_{\xi})^{j}
							h_n(r,\xi,c)\right|\\
							&\lesssim \sum_{\substack{0\leq j_1\leq j\\\mathrm{sgn}(i)\leq i_1\leq i+j-j_1}}\left|x^{i_1}
							(\lambda^c\pa_c)^l\pa_x^{i_1}(\xi\pa_{\xi})^{j_1}
							\tilde{h}_n(x,\xi,c)\right|+\left| x^{i_1+l}\pa_x^{i_1+l}(\xi\pa_{\xi})^{j_1} \tilde{h}_n(x,\xi,c)\right|
							\\
							&\quad+\sum_{\substack{0\leq j_1\leq j\\\mathrm{sgn}(i+l)\leq i_1\leq i+j-j_1}}\left|  x^{i_1}\pa_x^{i_1}(\xi\pa_{\xi})^{j_1} \tilde{h}_n(x,\xi,c)\right|.\notag
						\end{align}
						In fact, we recall that $\xi\pa_{\xi}x=x$ and
						\begin{align*}
							\pa_{\xi} h=\pa_{\xi}x\pa_x\tilde{h}+\pa_{\xi}\tilde{h}
							,\;  \pa rh(r,\xi,c)=\pa_rx\pa_{x}\tilde{h}(x,\xi,c),\;\pa_c
							h=x_{c}\pa_x\tilde{h}+\pa_{c}\tilde{h},
						\end{align*}
						which together with $[\xi\pa_{\xi},x\pa_x]\tilde{h}
						=(x\pa_x-\xi\pa_{\xi})\tilde{h}$, Leibniz's rule and Fa$\grave{\text{a}}$di Bruno's formula, yields that for $i,j\in\mathbb{N}$, $l\in\{0,1\}$,
						\begin{align*}
							&\left|
							(\lambda^c(r))^{l}(\lambda^r(r))^{i}\pa_c^l
							\pa_r^i(\xi\pa \xi)^j
							h_n(r,\xi,c)\right|\\
							&=
							\left|
							(\lambda^c)^{l}(\lambda^r)^{i}\pa_c^l
							\pa_r^i(x\pa_{x}+\xi\pa_{\xi})^{j}\tilde{h}_n(x,\xi,c)\right|
							\\
							&\lesssim
							\sum_{0\leq j_1\leq j,\;0\leq j_2\leq j-j_1
							}\left|(\lambda^c)^{l}(\lambda^r)^{i}(\pa_c^l
							\pa_r^i\circ x^{j_2}\pa_x^{j_2}(\xi\pa_{\xi})^{j_1}
							\tilde{h}_n(x,\xi,c)\right|\\
							&=
							\sum_{0\leq j_1\leq j,\;0\leq j_2\leq j-j_1,\;
								0\leq s\leq i}\left|(\lambda^c)^{l}(\lambda^r)^{i}\pa_c^l\pa_r^{i-s}(x^{j_2})
							\pa_r^s\pa_x^{j_2}(\xi\pa_{\xi})^{j_1}
							\tilde{h}_n(x,\xi,c)\right|\\
							&\lesssim
							\sum_{\substack{0\leq j_1\leq j,\;0\leq j_2\leq j-j_1,\;
									0\leq s\leq i\\k_1,...,k_s\geq 0,\;
									k_1+2k_2+...+sk_s=s\\
									k= k_1+...+k_s}}\left|(\lambda^c)^{l}(\lambda^r)^{i}\pa_c^l\circ
							\pa_r^{i-s}(x^{j_2})
							(\pa_rx)^{k_1}...(\pa_r^sx)^{k_s}\pa_x^k
							\pa_x^{j_2}(\xi\pa_{\xi})^{j_1}
							\tilde{h}_n(x,\xi,c)\right|\\
							&\lesssim\sum_{\substack{0\leq j_1\leq j,\;0\leq j_2\leq j-j_1,\;0\leq s\leq i\\ k'_0,...k'_{i-s}\geq 0,\;k_1,...k_s\geq 0\\k'_0+...+k'_{i-s}=j_2,\;k'_1+...+(i-s)k'_{i-s}=i-s
									\\k_1+...+k_s=k,\;k_1+2k_2+..sk_s=s}}
							\left|(\lambda^c)^{l}(\lambda^r)^{i}\pa_c^l\circ (x)^{k_0'}(\pa_rx)^{k'_1}\right.\\
							&\left.\quad\quad\quad\quad\quad\quad\quad\quad
							\quad\quad\quad\quad\quad\quad\quad\quad...(\pa_r^{i-s}x)^{k'_{i-s}}
							(\pa_rx)^{k_1}...(\pa_r^sx)^{k_s}
							\pa_x^{j_2+k}(\xi\pa_{\xi})^{j_1}
							\tilde{h}_n\right|.
						\end{align*}
						Notice that for integers satisfying $j_2+k\leq (i-s)+s\leq(i-s)j_2+sk $, we have $j_2+k\leq i\leq\max\{i-s,s\}(j_2+k)$, consequently, $j_2+k\geq \mathrm{sgn}(i)$.
						Together with the derivative bounds of $x$ in \eqref{integral2-trans-condition2},
						we may rewrite the above summation with new index $i_1:=j_2+k$, $\alpha_0=k_0'$, $\alpha'=k_1'+k_1$..., as follows
						 \begin{align*}
							&\lesssim
							\sum_{0\leq j_1\leq j,\;\mathrm{sgn}(i)\leq i_1\leq i+j-j_1,\;\atop
								\alpha_0,...,\alpha_i\geq 0,\;
								\alpha_0+\alpha'+..+\alpha_i= i_1,\;\alpha'+2\alpha+..+i\alpha_i=i}
							\left|(\lambda^c)^{l}(\lambda^r)^{i}\pa_c^l\circ x^{\alpha_0}(\pa_rx)^{\alpha'}...(\pa_r^{i}x)^{\alpha_{i}}
							\pa_x^{i_1}(\xi\pa_{\xi})^{j_1}
							\tilde{h}_n(x,\xi,c)\right|\\
							&\lesssim \sum_{0\leq j_1\leq j,\;\mathrm{sgn}(i)\leq i_1\leq i+j-j_1,\;\atop
								\alpha_0,...,\alpha_i\geq 0,\;\alpha_0+\alpha'+..+\alpha_i= i_1,\;\alpha'+2\alpha+..+i\alpha_i=i}\left(\left|(\lambda^c)^{l}(\lambda^r)^{i}x^{\alpha_0}(\pa_rx)^{\alpha'}...(\pa_r^{i}x)^{\alpha_{i}}
							\pa_c^l\pa_x^{i_1}(\xi\pa_{\xi})^{j_1}
							\tilde{h}_n\right|
							\right.\\
							&\left.\quad+ \left|(\lambda^c)^{l}(\lambda^r)^{i}x^{\alpha_0}(\pa_rx)^{\alpha'}
							...(\pa_r^{i}x)^{\alpha_{i}}
							(\pa_cx)^l\pa_x^{i_1+l}(\xi\pa_{\xi})^{j_1}
							\tilde{h}_n\right|\right.\\
							&\left.\quad+\left|(\lambda^c)^{l}(\lambda^r)^{i}\pa_c^l\left(x^{\alpha_0}
							(\pa_rx)^{\alpha'}...(\pa_r^{i}x)^{\alpha_{i}}
							\right)
							\pa_x^{i_1}(\xi\pa_{\xi})^{j_1}
							\tilde{h}_n\right|\right)\\
							&\lesssim\sum_{\substack{0\leq j_1\leq j,\;\mathrm{sgn}(i)\leq i_1\leq i+j-j_1}}\left|(\lambda^c)^{l} x^{i_1}\pa_c^l\pa_x^{i_1}(\xi\pa_{\xi})^{j_1}\tilde{h}_n\right|+\left| x^{i_1+l}\pa_x^{i_1+l}(\xi\pa_{\xi})^{j_1} \tilde{h}_n\right|\\
							&\quad+\sum_{\substack{0\leq j_1\leq j,\;\mathrm{sgn}(i+l)\leq i_1\leq i+j-j_1}}\left|  x^{i_1}\pa_x^{i_1}(\xi\pa_{\xi})^{j_1} \tilde{h}_n\right|.
						\end{align*}
						Now we are in a position to prove the derivative bounds.
						We apply the transform $(y_1,...y_n)\to (y_1',...y_n') $: $y_1'=y_1-x$, $y_2'=y_2-y_1=y_2-y_1'-x$,...$y_n'=y_n-y_{n-1}=y_n-y_{n-1}'-...-y_1'-x$.
						Then \eqref{eq:t-hn} can be rewritten as
						\begin{align}\label{eq:t-hn alt}
							h_n(r,\xi,c)&=\tilde{h}_n(x,\xi,c)\notag\\
							&=\small\int_{(\mathbb{R}^+)^n}\widetilde{K}_1(x,x+y_1',\xi,c)
							...
							\widetilde{K}_1(x+y_1'+...+y_{n-1}',x+y_1'+...+y_{n}',\xi,c)dy_n'...dy_2'dy_1'.
						\end{align}
						Therefore, thanks to the condition \eqref{condition-trans-tK}, for fixed $r\in I_{11}$(or $x$) and $j,i,l$, we get by \eqref{bd:t-hn} and \eqref{eq:t-hn}  that
						\begin{align*}
							&\left|(\lambda^c)^l(\lambda^r)^i\pa_c^l\pa_r^{i}(\xi\pa_{\xi})^{j}
							h_n(r,\xi,c)\right|\\
							&\lesssim \sum_{\substack{0\leq j_1\leq j\\0\leq i_1\leq i+j-j_1}}\left|(\lambda^c)^lx^{i_1}
							\pa_c^l\pa_x^{i_1}(\xi\pa_{\xi})^{j_1}
							\tilde{h}_n(x,\xi,c)\right|
							+ \sum_{\substack{0\leq j_1\leq j\\0\leq i_1\leq i+j-j_1}}\left|
							x^{i_1+l}\pa_x^{i_1+l}(\xi\pa_{\xi})^{j_1}
							\tilde{h}_n(x,\xi,c)\right|\\
							&\leq C^n\int_{(\mathbb{R}^+)^n}\varpi(x+y_1')...
							\varpi(x+...+y_n')dy_n'...dy_1'
							=\f{\left(C\int_x^{+\infty}\varpi(y,r,\xi,c)dy\right)^n}{n!}
							\leq \f{C^n\kappa(r,\xi,c)^n}{n!},
						\end{align*}
						where in the second inequality, we used the condition \eqref{condition-trans-tK} to obtain for $y'\geq 0$,
						\begin{align*}
							&\sum_{\substack{0\leq i'\leq i+j+l,\\0\leq j\leq j,\;0\leq l'\leq l}}\left((\lambda^c)^{l'}x^{i'}\right)(r,\xi,c)\sup_{\tilde{x}\in[x,x+y']}
							\left| \pa_c^{l'}\pa_{\tilde{x}}^{i'}
							\left(\xi\pa_{\xi}\right)^{j'}
							\widetilde{K}_1(\tilde{x},\tilde{x}+y',\xi,c)\right|\leq
							C\varpi(x+y',r,\xi,c).
						\end{align*}
						Summing up $f^{Rem}=\sum_{n=1}^{\infty}h_n$, we arrive at   $|(\lambda^c)^l(\lambda^r)^i\pa_c^l\pa_{r}^i(\xi\pa_{\xi})^j
						f^{Rem}(r,\xi,c)|\leq C\kappa(r,\xi,c)$.
					\end{proof}
					
					\subsection{Behavior of  $f_+$ as the Bessel function}
					We start with introducing some asymptotic behaviors of the Bessel functions as follows
					\begin{align}\label{Bessel:z>1}
						\begin{split}
							Y_{1}(z)&=O(z^{-\f12})\sin z,\;
							J_{1}(z)=O(z^{-\f12})\cos z\quad\text{for}\;z\gtrsim 1,\\
							H_{+}(z)&:=J_1(z)+iY_1(z)=O(z^{-\f12})e^{iz}\quad \text{for}\;z\gtrsim 1,\\
							J_{1}(z)&=O(z),\;
							Y_{1}(z)=O(z^{-1}),\;H_{+}(z):=O(z^{-1})
							\quad\text{for}\;0\leq z\lesssim 1.
						\end{split}
					\end{align}
					Here
					$O\left(f_0(z)\right)$  represents a function $f(z)$ with the bounds $|(z\pa_z)^{i}f(z)|\lesssim f_0(z)$, $i\in\mathbb{N}$.
					The function $H_+(\cdot)=J_{1}(\cdot)+iY_{1}(\cdot)$ is so called  the Hankel's function.

					\begin{lemma}\label{lem:k<1-f+} 	
					Let $c\in(0,1)$, $\xi\lesssim (1-c)^{\f13}$, and $i,j\in\mathbb{N}, l\in\{0,1\}$. The solution $f_+$ of \eqref{eq:f+} may write as
						\begin{align*}
							f_+(r,c,\xi)&=H_{+}(\xi r)\left(1+f^{Rem}(r,\xi,c)\right),
						\end{align*}
						with the remainder term satisfying  \begin{align}\label{Behave:f+1-1}
							\left|\left((1-c)\pa_c\right)^l(r\pa_r)^i(\xi\pa_{\xi})^j
							f^{Rem}(r,\xi,c)\right|\lesssim \f{\xi r^{-2}}{1-c}(\lesssim \f{\xi^3}{1-c}\lesssim 1),\;\;\text{for}\;\;r\gtrsim  \xi^{-1}. \end{align}
							Moreover, if $\xi\lesssim (1-c)^{\f12}$, it holds that
							\begin{align}\label{Behave:f+1-2}
							\left|\left((1-c)\pa_c\right)^l(r\pa_r)^i(\xi\pa_{\xi})^j f^{Rem}(r,\xi,c)\right|\lesssim \f{\xi^2r^{-1}}{1-c}(\lesssim \f{\xi^2}{1-c}\lesssim 1),\;\;\text{for}\;1\lesssim r\lesssim \xi^{-1}; \end{align}
							if $(1-c)^{\f12}\lesssim \xi\lesssim (1-c)^{\f13}$, it holds that
							\begin{align}\label{Behave:f+1-3}
							\left|\left((1-c)\pa_c\right)^l(r\pa_r)^i(\xi\pa_{\xi})^j f^{Rem}(r,\xi,c)\right|\lesssim \f{\xi^2r^{-1}}{1-c}(\lesssim \f{\xi^2}{(1-c)^{\f23}}\lesssim1),\;\;\text{for}\;(1-c)^{-\f13}
							\lesssim r\lesssim \xi^{-1}.
							\end{align}
							\end{lemma}
							\begin{proof}
								We rewrite \eqref{eq:f+} as
								\begin{align} \label{eq:f+Bessel}
									&
									f_+''+r^{-1}f_+'-r^{-2}f_++\xi^2f_+
									=\f{\xi^2(1-V)}{1-c}f_+,\;\;f_+(r,c,\xi)\sim (\xi r)^{-\f12}e^{i\xi r}(r\to +\infty),
									\end{align}
								which is equivalent to the integral equation
								\begin{align}
									\begin{split}
										&f_+(r,\xi,c)=H_+(\xi r)+
										\int_r^{+\infty}K(r,s,\xi,c)f_+(s,\xi,c)ds,\\
										&K(r,s,\xi,c)=
										\f{\pi\xi^2s(1-V(s))}{2(1-c)}\left(J_{1}(\xi r)Y_{1}(\xi s)-J_{1}(\xi s)Y_{1}(\xi r)\right)\;\;\text{with}\;\;K(s,s,\xi,c)=0\label{K:f+rgeq1}.
									\end{split}
								\end{align}
								Note that the assumption $\xi\lesssim(1-c)^{\f13}$ implies
								\begin{align}\label{ineq:xi-1geq 1}
									\xi^{-1}\gtrsim (1-c)^{-\f13}\gtrsim 1.
								\end{align}
								
								 We first prove \eqref{Behave:f+1-1} by applying Lemma \ref{lem:vorterra-infty} with
								 $D=\big\{(\xi,c)|c\in(0,1),\;\xi \lesssim (1-c)^{\f13}\big\}$, $I_1=\{r|r\gtrsim \xi^{-1}\}$,
								$f=f_+$, $g=H_+(\xi r)$, $\kappa=\f{\xi r^{-2}}{1-c}$.  It suffices to verify the condition \eqref{condition:Kinfty'} and \eqref{condition-trans-tK}.
								Notice that $H_+(\xi r)$ has no zero point. Thanks to \eqref{Bessel:z>1} and $1-V(s)=O(s^{-3})(s\gtrsim 1)$,  we infer that for $s\geq r\geq \xi^{-1}$,
								\begin{align}
									\notag K_1(r,s,\xi,c)&:= H_+^{-1}(\xi r)K(r,s,\xi,c)H_+(\xi s)\\
									\notag &=\f{\xi^2 sO(s^{-3})}{1-c}\cdot O\left((\xi r)^{\f12}\right)\cdot\left(O\left((\xi s)^{-\f12}\cdot(\xi r)^{-\f12}\right)\right)\cdot O\left((\xi s)^{-\f12}\right)\\
									\notag\quad&\;\;\;\;\cdot e^{-i\xi r}\left(\cos \xi r\sin \xi s-
									\cos \xi s\sin \xi r\right)e^{i\xi s}\\
									\notag&=\f{\xi^2O(s^{-2})}{1-c}O_{r,s,\xi}
									^{r,s,\xi}\left((\xi s)^{-1}\right)e^{i\xi(s-r)}\sin \xi(s-r),\\
									&=\f{\xi}{1-c}O_{r,s,\xi}^{r,s,\xi}(s^{-3})
									\left(e^{2i\xi(s-r)}-1\right),\label{est:H-1KH}
								\end{align}
								where $K$ is in \eqref{K:f+rgeq1}.
								 Therefore, the condition \eqref{condition:Kinfty'} holds as shown below(with $\kappa=\f{\xi r^{-1}}{1-c}$),
								\begin{align*}
									&\quad\int_r^{+\infty}
									\sup_{\tilde{r}\in[r,s)}\left|
									H_+^{-1}(\xi\tilde{r})K(\tilde{r},s,\xi,c)H_+(\xi s)\right|ds\lesssim\int_r^{+\infty}\f{\xi}{1-c} s^{-3}ds=\f{\xi}{1-c} r^{-2},\;\;
									r\gtrsim \xi^{-1}.
								\end{align*}
								
								To  verify the condition \eqref{condition-trans-tK}, we take $\lambda^c=1-c$, $\lambda^r=r$, $x=\xi r$, $y=\xi s$(satisfying \eqref{integral2-trans-condition2}),
								$\widetilde{K}_1=K_1\cdot\xi^{-1}$ and $\kappa=\f{\xi r^{-2}}{1-c}$. For $i,j\in\mathbb{N}$, $j\in\{0,1\}$, it follows from \eqref{est:H-1KH} that
								\begin{align*}
									\widetilde{K}_1(x,y,\xi,c)=H_+^{-1}(\xi r)K(r,s,\xi,c)H_+(\xi s)\cdot\xi^{-1}=\f{\xi^3}{1-c}O_{x,y}^{x,y}(y^{-3})
									\left(e^{2i(y-x)}-1\right),
								\end{align*}
								which shows that  for $x\gtrsim 1$, $\tilde{y}\geq 0$ and  $\tilde{x}\in[x,x+\tilde{y}]$,
								\begin{align*}
									\tilde{K}(\tilde{x},\tilde{x}+\tilde{y},\xi,c)=\f{\xi^3}{1-c}
									O_{\tilde{x}}^{\tilde{x}}((\tilde{x}+\tilde{y})^{-3})\left(e^{2i\tilde{y}}-1\right).
								\end{align*}
								Therefore, the condition holds as shown below,
								\begin{align*}
									&(1-c)^lx^i\sup_{\tilde{x}\in[x,x+\tilde{y}]}
									\left|(\tilde{x}\pa_{\tilde{x}})^i(\xi\pa_{\xi})^j\pa_c^l\tilde{K}(\tilde{x},
									\tilde{x}+\tilde{y},\xi,c)\right|\lesssim \f{\xi^3(x+\tilde{y})^{-3}}{1-c},
								\end{align*}
								and
								\begin{align*}
									\int_0^{+\infty}\f{\xi^3(x+\tilde{y})^{-3}}{1-c}d\tilde{y}\lesssim \f{\xi^3x^{-2}}{1-c}=\f{\xi r^{-2}}{1-c}.
								\end{align*}
								
								Next, we prove \eqref{Behave:f+1-2} and \eqref{Behave:f+1-3}. We will use another integral equation instead of  \eqref{K:f+rgeq1}.
								We plug $f_+(r,c,\xi)=H_{+}(\xi r)\left(1+f^{Rem}(r,\xi,c)\right)$ into \eqref{eq:f+Bessel} to obtain
								\begin{align*}
									H_{+}(\xi r)(f^{Rem})''
									+\left(2\left(H_{+}(\xi r)\right)'+r^{-1}H_{+}(\xi r)\right)(f^{Rem})' =\f{\xi^2(1-V)}{1-c} H_{+}(\xi r)\left(1+f^{Rem}\right),\end{align*}
								or
								\begin{align*}
									\left(rH_{+}^2(\xi r)(f^{Rem})'\right)'
									=\f{\xi^2(1-V)}{1-c}rH_{+}^2(\xi r))\left(1+f^{Rem}\right).\end{align*}
								For $r\leq C\xi^{-1}$(WLOG, taking $C=1$), we integrate  the above equation on $[r,\xi^{-1}]$  twice to obtain
								\begin{align*}
									f^{Rem}(r)&=f^{Rem}|_{r=\xi^{-1}}
									-\xi^{-1}(f^{Rem})'|_{r=\xi^{-1}}H_{+}^2(1)
									\int_{r}^{\xi^{-1}} s^{-1}H_{+}^{-2}(\xi s) ds \\
									&\quad+  \int_{r}^{\xi^{-1}} w^{-1}H_{+}^{-2}(\xi w)\int_{w}^{\xi^{-1}} sH_{+}^{2}(\xi s) \f{\xi^2(1-V(s))}{1-c}\left(1+f^{Rem}(s)\right) dsdw.
								\end{align*}
								We rewrite the above integral equation as: for $r\leq \xi^{-1}$,
								\begin{align}  \label{itegral:f+}
									\begin{split}
										f^{Rem}(r,\xi,c)&=f_0(r,\xi,c)+  \int_{r}^{\xi^{-1}} \mathrm{K}_1(r,s,\xi,c)\left(1+f^{Rem}(s,\xi,c)\right)ds\quad \text{with}\\
										K_{H_+}(r,s,\xi)&:=  sH_{+}^{2}(\xi s)\int_{r}^{s} w^{-1}H_{+}^{-2}(\xi w) dw,\;\;r\leq s\leq \xi^{-1},\\
										f_0(r,\xi,c)&:=f^{Rem}|_{r=\xi^{-1}}
										-(f^{Rem})'|_{r=\xi^{-1}} K_{H_+}(r,\xi^{-1},\xi), \\
										\mathrm{K}_1(r,s,\xi,c)&:=  K_{H_+}(r,s,\xi)\f{\xi^2(1-V(s))}{1-c}.
								\end{split}\end{align}
								It follows from the behavior of the Bessel function \eqref{Bessel:z>1}($H_+(z)=O_{z}^z(z^{-1})$,\; $z\lesssim 1$) that
								\begin{align} \label{est:KH+}
									K_{H_+}(r,s,\xi)= O_{s,r,\xi}^{s,r,\xi}(s),\;\;r\leq s\leq \xi^{-1}, \end{align}
								which along with $1-V(s)=O(s^{-3})(s\gtrsim 1)$ shows that for $ 1\lesssim r\leq s\lesssim \xi^{-1}$,
								\begin{align}  \label{est:K1f+}
									\mathrm{K}_1(r,s,\xi,c)&=  K_{H_+}(r,s,\xi)\f{\xi^2(1-V(s))}{1-c}=\f{\xi^2}{1-c}
									O_{s,r,\xi}^{s,r,\xi}(s^{-2}).
								\end{align}
								The previous bounds \eqref{Behave:f+1-1} and \eqref{est:KH+} ($ K_{H_+}(r,\xi^{-1},\xi)= O_{r,\xi}^{r,\xi}(\xi^{-1})$) ensure that
								\begin{align}\label{est:f0-f+} |\left((1-c)\pa_c\right)^l(r\pa_r)^i(\xi\pa_{\xi})^j
									f_0(r,\xi,c)|\lesssim \f{\xi^3}{1-c}. \end{align}
									Formally, the solution of equation \eqref{itegral:f+} takes the form of
									\begin{align*}
										&f^{Rem}(r,\xi,c)=\sum_{n=1}^{\infty}h_n(\tau,\xi,c)
										,\;\; \text{where}\;\;
										h_1=f_0(r,\xi,c)+
										\int_{r}^{\xi^{-1}}\mathrm{K}_1(r,s,\xi,c)ds,\\
										&h_{n}=
										\int_{r}^{\xi^{-1}}\mathrm{K}_1(r,s,\xi,c)h_{n-1}(s,\xi,c)ds(n\geq 2). \end{align*}
									 As in the proof of Lemma \ref{lem:vorteraa-0} with condition   \eqref{condition:K'pac-growth},
									 if we verify the following conditions on $f_0$ and $\mathrm{K}_1$:									
									 for $r\leq \f12\xi^{-1}$,
									\begin{align}
										\label{condtion:rf0}
										&\sup_{\tilde{t}\in[0,1]}\left|
										\left((1-c)\pa_c\right)^l \left(r\pa_{r}\right)^i
										(\xi\pa_{\xi})^jf_0(\tilde{t}(\xi^{-1}-r)+r,\xi,c)\right|
										\lesssim \f{\xi^3}{1-c}\leq \f{\xi^2r^{-1}}{1-c},\\
										&\int_0^1\sup_{\tilde{t}\in[0,t]}\left|
										\left((1-c)\pa_c\right)^l \left(r\pa_{r}\right)^i
										(\xi\pa_{\xi})^j\left((\xi^{-1}-r)
										\mathrm{K}_1(r+\tilde{t}(\xi^{-1}-r),
										r+t(\xi^{-1}-r),\xi,c)\right)\right|dt \notag\\
										&\quad\lesssim \f{\xi^2r^{-1}}{1-c},\label{condtion:rK1}
									\end{align}
									then we have
									\begin{align*}
										\left|
										\left((1-c)\pa_c\right)^l \left(r\pa_{r}\right)^i
										(\xi\pa_{\xi})^jh_n(r,\xi,c)\right|
										\leq  \f{C^n_{i,j,l}\left( \f{\xi^2r^{-1}}{1-c}\right)^{n}}{(n-1)!}\;\;n\geq 1, \end{align*}
									and then we infer that for $r\leq \f12\xi^{-1}$,
									\begin{align*} \left|\left((1-c)\pa_c\right)^l(r\pa_r)^i(\xi\pa_{\xi})^j
										f^{Rem}(r,\xi,c)\right|\lesssim \f{\xi^2 r^{-1}}{1-c}. \end{align*}
									If  we further assume $\xi\lesssim (1-c)^{\f12}$ and $ r\gtrsim 1$, then $ \f{\xi^2 r^{-1}}{1-c}\lesssim \f{\xi^2}{1-c}\lesssim 1$, the above bound gives \eqref{Behave:f+1-2}.
									If  we further assume $(1-c)^{\f12}\lesssim \xi\lesssim (1-c)^{\f13}$ and $r\gtrsim (1-c)^{-\f13}$, then $ \f{\xi^2 r^{-1}}{1-c}\lesssim \f{\xi^2}{(1-c)^{\f23}}\lesssim 1$,
									the above bound gives \eqref{Behave:f+1-3}.
									
									It remains to verify the above two conditions. Indeed, \eqref{condtion:rf0} can be deduced directly from \eqref{est:f0-f+}. For \eqref{condtion:rK1},
									notice that for $r\leq \f12\xi^{-1}$,
									$$ |(\xi\pa_{\xi})\left(\xi^{-1}-r\right)|\leq 2(\xi^{-1}-r),\;\;\;|\left(r\pa_{r}\right)
									\left(\xi^{-1}-r\right)|\leq (\xi^{-1}-r).$$ Then it follows from \eqref{est:K1f+}  that  for $r\leq \f12\xi^{-1}$,
									\begin{align*}
										&\int_0^1\sup_{\tilde{t}\in[0,t]}\left|
										\left((1-c)\pa_c\right)^l \left(r\pa_{r}\right)^i
										(\xi\pa_{\xi})^j\left((\xi^{-1}-r)
										\mathrm{K}_1(r+\tilde{t}(\xi^{-1}-r),
										r+t(\xi^{-1}-r),\xi,c)\right)\right|dt \notag\\
										&\lesssim \int_0^1(\xi^{-1}-r)\cdot \f{\xi^2}{1-c}\left(r+t(\xi^{-1}-r)\right)^{-2}dt
										\sim \f{\xi^2r^{-1}}{1-c}. \end{align*}
							\end{proof}

							\subsection{Behavior of  $f_+$ as the Airy function}
							
							\begin{lemma}\label{lem:f+2}
								Let  $c\in(0,V(0)]$, $\xi\gtrsim 1$, $Q(r,c)$ be defined in \eqref{def:Q}. Let $f_+$ be the unique smooth solution to the equation \eqref{eq:f+}, i.e.,
								\begin{align*}
									f''_++r^{-1}f'_+-r^{-2}f_+ +\f{\xi^2(V(r)-c)}{1-c}  f_+=0,\quad f_+(r,\xi,c)\sim (\xi r)^{-\f12}e^{i\xi \int_0^r\sqrt{\f{V(s)-c}{1-c}}ds}.
								\end{align*}
								Then for  $r\gtrsim \xi^{-1}\min\{\xi^{\f13},|c-V(0)|^{-\f12}\}$, there exists a complex number $C$ such that
								\begin{align}
									\label{Behave:f+3-1}
									&f_+(r,\xi,c)
									=C\left(\xi r Q^ {\f12}(r,c)\right)^{-\f12}e^{i\xi\int_0^r
										\sqrt{\f{V(s)-c}{1-c}}ds}\notag
									\left(1+f^{Rem}(r,\xi,c)\right),\\
									&\text{with}\;\;\left|(V(r)-c)r^{i}\pa_c^l\pa_r^i(\xi\pa_{\xi})^jf^{Rem}
									(r,\xi,c)\right|
									\lesssim \left(\xi rQ^ {\f12}(r,c)\right)^{-1}\lesssim 1.
								\end{align}
								Here $i,j\in\mathbb{N},
								\;l\in\{0,1\}$.
							\end{lemma}
							
							\begin{proof}
								Let $\tau=\tau(r,c)$ be defined as \eqref{def:tau} and $q(r,c)=\f{Q(r,c)}{\tau(r,c)}$. We perform the Langer transform that
								\begin{align}\label{def:2-g}
									\mathrm{f}_+(\tau,\xi,c)
									:=C\xi^{\f13}q^{\f14}r^{\f12}f_+(r,\xi,c)
								\end{align}
								for some complex number $C$ and parameters $(\xi,c)$,  such that  $\mathrm{f}_+(\tau,\xi,c)$ satisfies
								\begin{align*}
									\begin{split}
										\pa_{\tau}^2\mathrm{f}_++\xi^2\tau \mathrm{f}_+
										&=W\mathrm{f}_+,\;\;\mathrm{f}_+(\tau,\xi)\sim \mathrm{Ai}(-\xi^{\f23}\tau)-i
										\mathrm{Bi}(-\xi^{\f23}\tau)\;\; \tau\to+\infty,\\
										W&=q^{-\f14}\pa_{\tau}^2(q^{\f14})+\f{3}{4r^2q}.
									\end{split}
								\end{align*}
								Here we have used
								$\pa_r\tau(r,c)=q^{\f12}(r,c), \pa_{r}=q^{\f12}\pa_{\tau}$ and $W$ takes the same as \eqref{W:fm0}.
								
								Let $F_+(r,\xi,c)=\mathrm{f}_+(\tau(r,c),\xi,c)$. Then the equivalent integral form writes as
								\begin{align}
									\label{integral:f+langer1}
									F_+(r,\xi,c)&=\mathrm{Oi}(-\xi^{\f23}\tau(r,c))
									+\int_{\tau(r,c)}^{+\infty}K(r,s,\xi,c)F_+(s,\xi,c)ds,
								\end{align}
								where the oscillation Airy function $\mathrm{Oi}(-z)=\mathrm{Ai}(-z)-i\mathrm{Bi}(-z)$ behaves as \eqref{behave:Oi}  and
								\begin{align*}
									K(r,s,\xi,c)&=\left(\mathrm{Oi}(-\xi^{\f23}\tau(r,c))
									\mathrm{Oi}(-\xi^{\f23}t(s,c))\int_{\tau(r,c)}^{t(s,c)}
									\mathrm{Oi}^{-2}(-\xi^{\f23}t')dt'\right)\cdot W(s,c)\pa_st(s,c).
								\end{align*}
								
								We claim that $F_+(r,\xi,c)=\mathrm{Oi}(-\xi^{\f23}\tau(r,c))\left(1
								+F^{rem}(r,\xi,c)\right)$ for $r\gtrsim\xi^{-\f23}$ with
								\begin{align*}
									\left|(V(r)-c)r^{i}\pa_c^l\pa_r^i(\xi\pa_{\xi})^jF^{rem}(r,\xi,c)\right|\lesssim \left(\xi rQ^{\f12}\right)^{-1},\;i,j\in\mathbb{N},l\in\{0,1\}.
								\end{align*}
								Consequently, using the transform \eqref{def:2-g} and the behavior
								$\mathrm{Oi}(-z)=Cz^{-\f14}e^{\f{2i}{3}z^{\f32}}
								\left(1+O(z^{-\f32})\right)(z\gtrsim 1)$, we have
								\begin{align*}
									f_+&=C\xi^{-\f13}q^{-\f14}r^{-\f12}(\xi^{\f23}\tau)^{-\f14}
									e^{\f{2i}{3}\xi\tau^{\f32}}\left(1+O\left((\xi \tau^{\f32})^{-1}\right)\right)
									\left(1
									+F^{rem}(r,\xi,c)\right)\\
									&=C\left(\xi r Q^ {\f12}(r,c)\right)^{-\f12}e^{i\xi\int_0^r
										\sqrt{\f{V(s)-c}{1-c}}ds}(1+f^{rem}(r,\xi,c))
								\end{align*}
								where $f^{Rem}=\left(1+O\left((\xi \tau^{\f32})^{-1}\right)\right)
								\left(1
								+F^{rem}\right)-1$ admits the derivative bounds \eqref{Behave:f+3-1} due to the estimates of $F^{Rem}$ and $x=\xi\tau^{\f32}$ in \eqref{estx1-3}.
								
								To prove the claim for $F^{Rem}$, we apply Lemma \ref{lem:vorterra-infty}  to the integral equation \eqref{integral:f+langer1}
								with $f=F_+$, $g=Oi(-\xi^{-\f23} \tau)$, $D=\big\{(\xi,c)|c\in(0,V(0)],\;\xi \gtrsim  1\big\}$, $I_1=[C\xi^{-\f23},+\infty)$. Next, we verify two conditions in the lemma.
								
								\emph{Step 1.} We verify the condition \eqref{condition:Kinfty'}($\kappa=\left(\xi rQ^{\f12}\right)^{-1}$)) for $K_1$,
								\begin{align*}
									K_1(r,s,\xi,c)&=\mathrm{Oi}^{-1}\left(-\xi^{\f23}\tau(r,c)\right)
									K(r,s,\xi,c)\mathrm{Oi}\left(-\xi^{\f23}t(s,c)\right)\\
									&=K_{Oi}(\tau(r,c),t(s,c),\xi)
									W(s,c)\pa_st(s,c),\;\;r\leq s
								\end{align*}
								where $K_{Oi}$ is defined in \eqref{def:Knu}.
								It follows from Lemma \ref{lem:behave-W} (a) and \eqref{est:KOi} that for $r\gtrsim \xi^{-1}\min\{\xi^{\f13},|c-V(0)|^{-\f12}\}$,
								\begin{align}\label{bd:K1-f+3}
									\notag&\int_{r}^{+\infty}
									\sup_{\tilde{r}\in[r,s]}|K_1(\tilde{r},s,\xi,c)|ds\lesssim \xi^{-1}\int_{\tau(r)}^{+\infty}t^{-\f12}|W(s(t))|dt\\
									&\qquad\lesssim \xi^{-1}\int_{\tau(r)}^{+\infty}t^{-\f52}dt
									\lesssim (\xi\tau^{\f32}(r))^{-1}\sim \left(\xi rQ^{\f12}(r,c)\right)^{-1}\lesssim 1.
								\end{align}
								Here we used Lemma \ref{lem:behave-Q}, $Q(r,c)\sim \f{r}{\langle r\rangle}+(V(0)-c)$ to obtain
								\begin{align*}
									\xi rQ^{\f12}(r,c)\sim \xi r\left(\f{r^{\f12}}{\langle r\rangle^{\f12}}+(V(0)-c)^{\f12}\right)\gtrsim 
									1,\;\text{for}\;r\gtrsim \xi^{-1}\min\{\xi^{\f13},
									|c-V(0)|^{-\f12}\}.
								\end{align*}
								
								\emph{Step 2.} We verify the condition \eqref{condition-trans-tK}($\lambda^r=r$, $\lambda^c(r,c)=V(r)-c$, $\varpi(y)=y^{-2}$)
								corresponding to $\widetilde{K}_1$,  which is calculated as the formula \eqref{def:tK1}:
								\begin{align}
									&\quad\quad\quad\notag\widetilde{K}_1(x,y,\xi,c)
									:=K_{Oi}(\tau(r,c),t(s,c),\xi)\cdot
									W(s,c)\cdot\f23\xi^{-1}t^{-\f12},
								\end{align}
								where $\f23\xi^{-1}t^{-\f12}=\f{\pa_st}{\pa_sy}$ was used. The transform $(r,s,\xi,c)\to \left(x(r,\xi,c),y(s,\xi,c),\xi,c\right)$ takes for $s\geq r(\text{or}\;\;t\geq \tau,\;y\geq x)$:
								\begin{align*}
									x=\xi\tau^{\f32}(r,c),\;y=\xi t^{\f32}(s,c),
								\end{align*}
								where $\tau$ and $t$ are defined as \eqref{def:tau}, and \eqref{estx1-3} matches the bounds \eqref{integral2-trans-condition2} with  $\lambda^r=r$, $\lambda^c(r,c)=V(r)-c$ for $r>0$:
								\begin{align*}
									\left|(V(r)-c)^lr^i\pa_c^l\pa_r^ix\right|\lesssim x.
								\end{align*}
								
								Since the monotonicity $\pa_rx,\pa_sy>0$, the inverse function $r=r(x,\xi,c)$, $s=s(y,\xi,c)$ exist.
								We conclude by \eqref{asy:tKOi} that for $y\geq x\gtrsim 1$,
								$$\widetilde{K}_{Oi}(x,y,\xi):=K_{Oi}(\tau(r,c),t(s,c),\xi)
								=\f23\xi^{-\f23}y^{-\f13}\left(1+O(y^{-1})\right)
								\int_{0}^{y-x}e^{\f{2i}{3}\tilde{y}}\left(1+O((y-\tilde{y})^{-1})\right)d\tilde{y}.$$
								We also notice that $\f23\xi^{-1}t^{-\f12}=\f23\xi^{-\f23}y^{-\f13}$. Consequently, letting $y=\tilde{x}+\tilde{y}$ and  $\tilde{x}\geq x\gtrsim 1$ (which gives $\tilde{y}\geq 0$), we have
								\begin{align*}
									&\quad\f23\xi^{-\f23}(\tilde{x}+\tilde{y})^{-\f13}
									\widetilde{K}_{Oi}(\tilde{x},\tilde{x}+\tilde{y},\xi)\\
									&=\f43\xi^{-\f43}(\tilde{x}+\tilde{y})^{-\f23}\left(1+O((\tilde{x}+\tilde{y})^{-1})\right)
									\int_{0}^{\tilde{y}}e^{\f{2i}{3}\tilde{w}}\left(1+O\left((\tilde{x}+\tilde{y}
									-\tilde{w})^{-1}\right)\right)
									d\tilde{w},\\
									&\quad\quad\quad\quad\quad\quad\text{independent}\;\text{of}\;c.
								\end{align*}
								Integration by parts gives
								\begin{align*}
									\left|\int_{0}^{\tilde{y}}e^{\f{2i}{3}\tilde{w}}\left(1+O\left((\tilde{x}+\tilde{y}
									-\tilde{w})^{-1}\right)\right)
									d\tilde{w}\right|\lesssim 1.
								\end{align*}
								Therefore, it holds that for $\tilde{x}\geq x$ and $\tilde{y}\geq 0$,
								\begin{align*}
									\left|(V(s)-c)^l(\tilde{x}+\tilde{y})^i\pa_c^l\pa_{\tilde{x}}^i(\xi\pa_{\xi})^j
									\left(\xi^{-\f23}(\tilde{x}+\tilde{y})^{-\f13}
									\widetilde{K}_{Oi}(\tilde{x},\tilde{x}+\tilde{y},\xi)\right)\right|
									\lesssim \xi^{-\f43}(x+\tilde{y})^{-\f23}.
								\end{align*}
								We also conclude by \eqref{est:Dy-tW} that for $\widetilde{W}(y,\xi,c):=W(s(y,\xi,c),c)$, $y\geq 0$,
								\begin{align*}
									\left|(V(s)-c)^ly^i\pa_c^l\pa_y^i(\xi\pa_{\xi})^j\widetilde{W}(y,\xi,c)\right|
									\lesssim \xi^{\f43}y^{-\f43}.
								\end{align*}
								Then letting  $y=\tilde{x}+\tilde{y}$,  we infer that for $y\geq \tilde{x}\geq x\gtrsim 1$,
								\begin{align*}
									\left|(V(s)-c)^l(\tilde{x}+\tilde{y})^i\pa_c^l\pa_{\tilde{x}}^i
									(\xi\pa_{\xi})^j
									\left(\widetilde{W}(\tilde{x}+\tilde{y},\xi,c)\right)\right|
									\lesssim \xi^{\f43}(x+\tilde{y})^{-\f43},
								\end{align*}
								here the inverse function takes $s=s(y)=s(\tilde{x}+\tilde{y})$.
								Collecting the above two derivative bounds and using $V(r)-c\leq V(s)-c$(due to $r(x)\leq \tilde{r}(\tilde{x})\leq s(\tilde{x}+\tilde{y})$ and $V'>0$), we finally arrive at
								\begin{align*}
									&(V(r)-c)^l\sup_{\tilde{x}\in[x,x+\tilde{y}]}\notag
									\left|(\tilde{x})^i\pa_c^l\pa_{\tilde{x}}^i(\xi\pa_{\xi})^j\tilde{K}_1(\tilde{x},
									\tilde{x}+\tilde{y},\xi,c)\right|\lesssim (x+\tilde{y})^{-2},\\
									&\int_{0}^{+\infty}(x+\tilde{y})^{-2}d\tilde{y}\lesssim x^{-1}= \left(\xi \tau^{\f32}\right)^{-1}\sim (\xi rQ^{\f12})^{-1}.
								\end{align*}
								Thus, the condition \eqref{condition-trans-tK} holds.
							\end{proof}
							
							\begin{lemma}\label{lem:f+4}
								Let  $c\in(V(0),1)$, $r_c=V^{-1}(c)$. 	Let $f_+$ be the solution to the equation \eqref{eq:f+}, i.e.,
								\begin{align*}
									f''_+ +r^{-1}f'_+-r^{-2}f_+ +\f{\xi^2(V-c) }{1-c} f_+=0,\quad f_+(r,\xi,c)\sim (\xi r)^{-\f12}e^{i\xi \int_{r_c}^r\sqrt{\left|\f{V(s)-c}{1-c}\right|}ds}(r\to+\infty).
								\end{align*}
								\begin{itemize}
									\item[(1)] If $\xi r_c^{\f32}/\langle r_c\rangle^{\f12}\lesssim 1 $ (which is equivalent to $\xi\lesssim   (1-c)^{\f13}|c-V(0)|^{-\f32}$),
									 then  there exists a complex constant $C$ such that								
									\begin{itemize}

									\item[\textbf{(1.1)}] for $r\geq C_0\xi^{-\f23}\langle \xi^{-\f13}\rangle\geq  2r_c$($C_0$ large enough),
									\begin{align}\label{est:f+6(1+)}
										\notag &f_+(r,\xi,c)
										=C\left(\xi r Q^ {\f12}(r,c)\right)^{-\f12}e^{i\xi\int_{r_c}^r
											\sqrt{\left|\f{V(s)-c}{1-c}\right|}ds}
										\left(1+f^{Rem}(r,\xi,c)\right)\quad\text{with}\\
										&\left|
										\left(\f{r}{\langle r_c\rangle ^{4}}\right)^lr^i
										\pa_c^l\pa_r^i(\xi\pa_{\xi})^jf^{Rem}(r,\xi,c)\right|
										\lesssim (\xi r^{\f32}/\langle r\rangle^{\f12})^{-1},
									\end{align}
									where  $Q$ is defined in \eqref{def:Q}.
									
									\item[\textbf{(1.2)}] for $C_0\xi^{-\f23}\langle \xi^{-\f13}\rangle\leq  r\leq 2r_c $,
									\begin{align}\label{est:f+6(1)}
										f_+(r,\xi,c)
										&=C\xi^{-\f13}q^{-\f14}(r,c)r^{-\f12}\mathrm{Oi}\left(-\xi^{\f23}\tau(r,c)\right)
										\left(1+f^{Rem}(r,\xi,c)\right),\quad\left|
										f^{Rem}\right|
										\lesssim 1,
									\end{align}
									where $\mathrm{Oi}(-z)$ ($z\gtrsim -1$) is defined in \eqref{behave:Oi},  $\tau=\tau(r,c)$ is defined in \eqref{def:tau}, $q=\f{Q}{\tau}$.
									\end{itemize}
									
									\item[(2)] If $\xi r_c^{\f32}/\langle r_c\rangle^{\f12}\gtrsim 1$
									(which is equivalent to $\xi^{-\f23}\langle r_c\rangle^{\f13}\leq Cr_c$ and $\xi\gtrsim  (1-c)^{\f13}|c-V(0)|^{-\f32}$), then  there exists a complex constant  $C$ such that
									\begin{itemize}

									\item[\textbf{(2.1)}] for $r\geq r_c+C\xi^{-\f23}\langle r_c\rangle^{\f13}$,
									\begin{align}\label{est:f+6(2a)}
										\notag &f_+(r,\xi,c)
										=C\left(\xi r Q^ {\f12}(r,c)\right)^{-\f12}e^{i\xi\int_{r_c}^r
											\sqrt{\left|\f{V(s)-c}{1-c}\right|}ds}\notag
										\left(1+f^{Rem}(r,\xi,c)\right)\quad\text{with}\\
										&\left|
										\left(\f{(r-r_c)}{\langle r_c\rangle ^{4}}\right)^l(r-r_c)^i
										\pa_c^l\pa_r^i(\xi\pa_{\xi})^jf^{Rem}(r,\xi,c)\right|
										\lesssim \left(\xi (r-r_c)^{\f32}/\langle r\rangle^{\f12}\right)^{-1}.
									\end{align}

							\item[\textbf{(2.2)}] for  $r_c-C\xi^{-\f23}\langle r_c\rangle^{\f13}\leq r\leq r_c+C\xi^{-\f23}\langle r_c\rangle^{\f13}$, and $\xi r_c^{\f32}/\langle r_c\rangle^{\f12}\gtrsim  M\gg 1$\big(such that $[r_c-C\xi^{-\f23}\langle r_c\rangle^{\f13},r_c+C\xi^{-\f23}\langle r_c\rangle^{\f13}]\subset [r_c/2,2r_c]$\big),
									\begin{align}\label{est:f+6(2b)}
										\notag &f_+(r,\xi,c)
										=C\xi^{-\f13}q(r,c)^{-\f14}r^{-\f12}\text{Oi}\left(-\xi^{\f23}\tau(r,c)\right)
										\left(1+f^{Rem}(r,\xi,c)\right)\quad\text{with}\\
										&\left|\left(\f{\xi^{-\f23}}{\langle r_c\rangle ^{\f{11}{3}}}\right)^l\left(\xi^{-\f23}\langle r_c\rangle^{\f13}\right)^i
										\pa_c^l\pa_r^i(\xi\pa_{\xi})^j
										f^{Rem}(r,\xi,c)\right|
										\lesssim   \left(\xi r_c^{\f32}/\langle r_c\rangle^{\f12}\right)^{-1}.
									\end{align}
								\end{itemize}	
								\end{itemize}
							\end{lemma}
							
							\begin{proof}
								Letting $q(r,c)=\f{Q(r,c)}{\tau(r,c)}$, we again perform Langer transform on the solution $\mathrm{F}_+(\tau,\xi,c)
								=C\xi^{\f13}q^{\f14}r^{\f12}f_+(r,\xi,c)$ for some complex number $C$ to get
								\begin{align} \label{eq:mathrmF+}
									\begin{split}
										\pa_{\tau}^2\mathrm{F}_++\xi^2\tau \mathrm{F}_+
										&=W\mathrm{F}_+,\;\;\mathrm{F}_+(\tau,\xi)\sim  \mathrm{Oi}\left(-\xi^{\f23}\tau\right)\quad \tau\to+\infty,
									\end{split}
								\end{align}
								Here  $W$ is the same as \eqref{W:fm0}.
								Let $F_+(r,\xi,c)=\mathrm{F}_+(\tau(r,c),\xi,c)$. Then the equivalent integral form writes as
								\begin{align}
									\label{integral:f+langer2}
									F_+(r,\xi,c)&=\mathrm{Oi}\left(-\xi^{\f23}\tau(r,c)\right)
									+\int_{\tau(r,c)}^{+\infty}K(r,s,\xi,c)F_+(s,\xi,c)ds,
								\end{align}
								where $$K(r,s,\xi,c)=\left(\mathrm{Oi}\left(-\xi^{\f23}\tau(r,c)\right)
								\mathrm{Oi}\left(-\xi^{\f23}t(s,c)\right)\int_{\tau(r,c)}^{t(s,c)}
								\mathrm{Oi}^{-2}\left(-\xi^{\f23}t'\right)dt'\right)\cdot W(s,c)\pa_st(s,c).$$
								
								We claim that
								\begin{align}
									\label{def:F+case6inf+}F_+(r,\xi,c)=\mathrm{Oi}\left(-\xi^{\f23}\tau(r,c)\right)\left(1+F^{rem}(r,\xi,c)\right),
								\end{align}
								where the remainder holds that for $\xi\gtrsim (1-c)^{\f13}|c-V(0)|^{-\f32}\sim r_c^{-\f32}\langle r_c\rangle^{\f12}$, $c\in(V(0),1)$ and $i,j\in\mathbb{N},l\in\{0,1\}$,
								
								\begin{itemize}
									\item[(a)] if $r\gtrsim r_c$, 
									then
									\begin{align*}
										\left|F^{Rem}(r,\xi,c)\right|
										\lesssim (\xi \max\{r,2r_c\}^{\f32}/\langle \max\{r,2r_c\}\rangle^{\f12})^{-1}.
									\end{align*}
									
									\item[(b)] if $r\geq r_c+C\xi^{-\f23}\langle r_c\rangle^{\f13}$, then
									\begin{align*}
										\left|
										\left(\f{r-r_c}{\langle r_c\rangle ^{4}}\right)^l(r-r_c)^i
										\pa_c^l\pa_r^i(\xi\pa_{\xi})^jF^{Rem}(r,\xi,c)\right|
										\lesssim \left(\xi \max\{r,2r_c\}^{\f32}/\langle \max\{r,2r_c\}\rangle^{\f12}\right)^{-1}.
									\end{align*}
									
									\item[(c)] if $(r_c/2)\leq r_c+C\xi^{-\f23}\langle r_c\rangle^{\f13}\leq r\leq r_c-C\xi^{-\f23}\langle r_c\rangle^{\f13}(\leq 2r_c)$, 
									then
									\begin{align*}
										\left|
										\left(\f{\xi^{-\f23}}{\langle r_c\rangle ^{\f{11}{3}}}\pa_c\right)^l\left(\xi^{-\f23}\langle r_c\rangle^{\f13}\pa_r\right)^i
										(\xi\pa_{\xi})^jF^{Rem}(r,\xi,c)\right|
										\lesssim  \left(\xi r_c^{\f32}/\langle r_c\rangle^{\f12}\right)^{-1}.
									\end{align*}
									\end{itemize}
									Moreover, we apply $F_+(r,\xi,c)
									=C\xi^{\f13}q^{\f14}r^{\f12}f_+$, \eqref{estx4-6-G} and the behavior
									$$\mathrm{Oi}(-z)=Cz^{-\f14}e^{\f{2i}{3}z^{\f32}}
									\left(1+O(z^{-\f32})\right)(z\gtrsim  1)$$ to deduce for  $r\geq r_c+C\xi^{-\f23}\langle r_c\rangle^{\f13}$, or
									\begin{align}\label{est:x-1}
										r\geq r_c\;\text{and}\;x(r,\xi,c)=\xi \mathrm{sgn}(r-r_c)\tau^{\f32}(r,c)\sim \xi (r-r_c)^{\f32}/\langle r\rangle ^{\f12}\gtrsim 1
									\end{align}
									that
									\begin{align*}
										f_+&=C\xi^{-\f13}q^{-\f14}r^{-\f12}(\xi^{\f23}|\tau|)^{-\f14}
										e^{\f{2i}{3}\xi|\tau|^{\f32}}\left(1+O\left(x^{-1}\right)\right)
										\left(1
										+F^{rem}(r,\xi,c)\right)\\
										&=C\left(\xi r Q^ {\f12}(r,c)\right)^{-\f12}e^{i\xi\int_{r_c}^r
											\sqrt{\f{V(s)-c}{1-c}}ds}\left(1+f^{rem}(r,\xi,c)\right)
									\end{align*}
									where \begin{align} \label{def:fRem=FREm}
										f^{Rem}=F^{rem}+O\left( x^{-1}\right)+O\left(x^{-1}\right)F^{rem}. \end{align}
									We also observe that $r^{\f32}/\langle r\rangle ^{\f12}$ is monotonically increasing in $r$.
									Consequently, we will illustrate that (1), (2) can be deduced by (a), (b), (c) and \eqref{est:x-1}, \eqref{def:fRem=FREm}:
									\smallskip
									
									{\bf  (a),(b), \eqref{def:fRem=FREm}, \eqref{est:x-1} $\Rightarrow$ (1.1)}. It suffices to check for $r\geq C\xi^{-\f23}\langle \xi^{-\f13}\rangle\geq  2r_c$ that $r-r_c\sim r$, and
									\begin{align*}
									\xi (r-r_c)^{\f32}/\langle r\rangle^{\f12}&\sim \xi r^{\f32}/\langle r\rangle^{\f12}\geq \xi r^{\f32}/\langle r\rangle^{\f12}|_{r=C\xi^{-\f23}\langle \xi^{-\f13}\rangle }
									\sim \f{\langle \xi^{-\f13}\rangle^{\f32}}{\langle \xi^{-\f23}\langle \xi^{-\f13}\rangle\rangle^{\f12}}\sim 1,
									\end{align*}
									and  \begin{align*}
										\xi \max\{r,2r_c\}^{\f32}/\langle \max\{r,2r_c\}\rangle^{\f12}=\xi r^{\f32}/\langle r\rangle^{\f12}.
									\end{align*}
									
									{\bf (a) $\Rightarrow$ \textbf{(1.2)}}. It suffices to check for $ C\xi^{-\f23}\langle \xi^{-\f13}\rangle\leq r\leq 2r_c$,
									$\xi r_c^{\f32}/\langle r_c\rangle^{\f12}\lesssim 1 $ and
									$\xi\gtrsim  \langle r_c\rangle^{-1}$ that  $\xi\gtrsim 1$, $\xi r_c^{\f32}\sim 1$ for $r_c\lesssim 1$ and  $ \xi\sim r_c^{-1}$  for $r_c\gtrsim 1$, which imply
									\begin{align*}
									&\xi r_c^{\f32}/\langle r_c\rangle^{\f12}\sim  1,\quad r\geq C\xi^{-\f23}\langle \xi^{-\f13}\rangle\gtrsim r_c,\\									
									&\xi \max\{r,2r_c\}^{\f32}/\langle \max\{r,2r_c\}\rangle^{\f12}= \xi (2r_c)^{\f32}/\langle 2r_c\rangle^{\f12}\sim 1.
									\end{align*}
									
									{\bf (b), \eqref{def:fRem=FREm}, \eqref{est:x-1} $\Rightarrow$ \textbf{(2.1)}}.
									It suffices to check for $\xi r_c^{\f32}/\langle r_c\rangle^{\f12}\gtrsim 1$, $r\geq r_c+C^{-1}\xi^{-\f23}\langle r_c\rangle^{\f12}$ that $r\geq r_c$,
									$\xi (r-r_c)^{\f32}/\langle r\rangle^{\f12}\gtrsim 1$ and obviously, it holds that for $r\geq r_c$,
									\begin{align*}
										\xi \max\{r,2r_c\}^{\f32}/\langle \max\{r,2r_c\}\rangle^{\f12}\geq  \xi (r-r_c)^{\f32}/\langle 2r\rangle^{\f12}.
										\end{align*}
										
										{\bf (c) $\Rightarrow$ \textbf{(2.2)}}.
										It suffices to check for $\xi r_c^{\f32}/\langle r_c\rangle^{\f12}\gtrsim 1 $, $r_c-C\xi^{-\f23}\langle r_c\rangle^{\f12}\leq r\leq r_c+C\xi^{-\f23}\langle r_c\rangle^{\f12}$ that
										$\xi (r-r_c)^{\f32}/\langle r\rangle^{\f12}\lesssim 1$ and $r\sim r_c$.
										
										\medskip
										
										Now we prove the claim for the remainder term $F^{Rem}$ in three cases
										via applying Lemma \ref{lem:vorterra-infty} to the integral equation \eqref{integral:f+langer2} with $f=F_+$, $g=\mathrm{Oi}\left(-\xi^{-\f23} \tau\right)$,
										$D=\big\{(\xi,c)|c\in(V(0),1),\;\xi \gtrsim  (1-c)^{\f13}\big\}$.
										We verify the conditions \eqref{condition:Kinfty'} and \eqref{condition-trans-tK} 
										in the lemma.\medskip
										
										\noindent\emph{Proof of (a)}.  We take
										$I_{1}=\{r|r\gtrsim r_c\}$ in  Lemma \ref{lem:vorterra-infty}.
										It suffices to verify the condition \eqref{condition:Kinfty'} ($\kappa=\xi \max\{r,2r_c\}^{\f32}/\langle \max\{r,2r_c\}\rangle^{\f12}$) for $K_1$:
										\begin{align}\label{def:f+6-K1}
											\notag K_1(r,s,\xi,c)&=\mathrm{Oi}^{-1}\left(-\xi^{\f23}\tau(r,c)\right)
											K(r,s,\xi,c)\mathrm{Oi}\left(-\xi^{\f23}t(s,c)\right)\\
											&=K_{Oi}(\tau(r,c),t(s,c),\xi)
											W(s,c)\pa_st(s,c) \notag \\
											&=\mathrm{Oi}\left(-\xi^{\f23}t(s,c)\right)
											\int_{\tau(r,c)}^{t(s,c)}Oi^{-2}\left(-\xi^{\f23}t'\right) dt' \cdot (Wq^{\f12})(s,c),\;\;r\leq s,
										\end{align}
										where $K_{Oi}$ is defined in \eqref{def:Knu}.
										It follows from  \eqref{est:KOi}, Lemma \ref{lem:behave-W}, \eqref{estx4-6-G}
										that for $r\gtrsim r_c$,
										\begin{align}\label{bd:K1-f+6}
											\notag\int_{r}^{+\infty}
											\sup_{\tilde{r}\in[r,s]}|K_1(\tilde{r},s,\xi,c)|ds&\lesssim \xi^{-1}\int_{\tau(r)}^{+\infty}|t|^{-\f12}|W(s(t))|dt\\
											&\lesssim \xi^{-1}\int_{\tau(\max\{r,2r_c\})}^{+\infty}
											|t|^{-\f52}dt
											+\xi^{-1}\int_{r}^{\max\{r,2r_c\}}\f{1}{s^2|Q(s,c)|^{\f12}}ds\\
											&\sim (\xi\tau^{\f32}(2r_c))^{-1}+\xi^{-1}r_c^{-2}\langle r_c\rangle^{\f12}
											\int_{r}^{\max\{r,2r_c\}}\f{1}{|r_c-s|^{\f12}}ds\notag\\
											&\lesssim\left(\xi \max\{r,2r_c\}^{\f32}/\langle \max\{r,2r_c\}\rangle^{\f12}\right)^{-1}.\notag
										\end{align}
										
										\noindent\emph{Proof of (b)}. We take $I_{1}=\{r|r\geq r_c+C^{-1}\xi^{-\f23}\langle r_c\rangle^{\f12}\}$
										and $\kappa=\xi \max\{r,2r_c\}^{\f32}/\langle \max\{r,2r_c\}\rangle^{\f12}$ in Lemma \ref{lem:vorterra-infty}. If $i=j=l=0$, we follow the same as \eqref{bd:K1-f+6}.
										If $i+j+l\geq 1$, we verify the condition \eqref{condition-trans-tK}
								               ($\lambda^r=r-r_c$, $\lambda^c(r,c)=\f{r-r_c}{\langle r_c\rangle^4}$, $\varpi(y)= y^{-2}\mathbf{1}_{y\geq y(\max\{r,2r_c\})}+r_c^{-2}\langle r_c\rangle^{\f23}\mathbf{1}_{y\leq y(\max\{r,2r_c\})}$) corresponding to $\widetilde{K}_1$ (calculated as formula \eqref{def:tK1}):
										\begin{align}\label{def:tK1-6}
											\notag\widetilde{K}_1(x,y,\xi,c)
											&:=\f23\xi^{-1}t^{-\f12}\cdot K_{Oi}(\tau(r,c),t(s,c),\xi)\cdot
											W(s,c),\\
											&=\f23\xi^{-\f23}y^{-\f13}\widetilde{K}_{Oi}(x,y,\xi)\widetilde{W}
											(y,\xi,c)
										\end{align}
										where $\f23\xi^{-1}t^{-\f12}=\f{\pa_st}{\pa_sy}$ was used, $\widetilde{K}_{Oi}(x,y,\xi)=
										K_{Oi}(\tau(r,c),t(s,c),\xi)$ is defined in \eqref{def:Knu}, $\widetilde{W}$ is defined in Lemma \ref{lema:W-trans},
										and the corresponding transform $(r,s,\xi,c)\to \left(x(r,\xi,c),y(s,\xi,c),\xi,c\right)$ takes for $s\geq r(\text{or}\;\;t\geq \tau,\;y\geq x)$:
										\begin{align*}
											x=\xi\tau^{\f32}(r,c),\;y=\xi t^{\f32}(s,c),
										\end{align*}
										where $\tau$ and $t$ are defined as \eqref{def:tau}, and the bounds \eqref{estx4-6-G} matches the bounds \eqref{integral2-trans-condition2} with  $\lambda^r=r-r_c$, $\lambda^c(r,c)=\f{r-r_c}{\langle r_c\rangle^4}$ for $r\geq r_c$:
										\begin{align*}
											\left|\left(\f{r-r_c}{\langle r_c\rangle^4}\right)^l|r-r_c|^i\pa_c^l\pa_r^ix\right|\lesssim x.
										\end{align*}
										By the monotonicity $\pa_rx,\pa_sy>0$, the inverse function $r=r(x,\xi,c)$, $s=s(y,\xi,c)$ exists.
										We conclude by \eqref{convert:Knv-tKnv} and \eqref{asy:tKOi} that for $y\geq x\gtrsim 1$,
										$$\widetilde{K}_{Oi}(x,y,\xi):=K_{Oi}(\tau(r,c),t(s,c),\xi)
										=\f23\xi^{-\f23}y^{-\f13}\left(1+O(y^{-1})\right)
										\int_{0}^{y-x}e^{\f{2i}{3}\tilde{w}}\left(1+O\left((y-\tilde{w})^{-1}\right)\right)d\tilde{w}.$$
										We also notice that $\f23\xi^{-1}t^{-\f12}=\f23\xi^{-\f23}y^{-\f13}$. Identically as the step 2 in the proof of Lemma \ref{lem:f+2}, letting $y=\tilde{x}+\tilde{y}$ and $\tilde{x}\geq x\gtrsim 1$ (which gives $\tilde{y}\geq 0$), we have
										\begin{align*}
											&\f23\xi^{-\f23}(\tilde{x}+\tilde{y})^{-\f13}
											\widetilde{K}_{Oi}(\tilde{x},\tilde{x}+\tilde{y},\xi)\\
											&=\f43\xi^{-\f43}(\tilde{x}+\tilde{y})^{-\f23}\left(1+O\left((\tilde{x}+\tilde{y})^{-1}\right)\right)
											\int_{0}^{\tilde{y}}e^{\f{2i}{3}\tilde{y}}\left(1+O\left((\tilde{x}+\tilde{y}-\tilde{w})^{-1}\right)\right)
											d\tilde{w},
										\end{align*}
										which is independent of $c$ and
										\begin{align*}
											\left|\int_{0}^{\tilde{y}}e^{\f{2i}{3}\tilde{w}}\left(1+O\left((\tilde{x}+\tilde{y}
											-\tilde{w})^{-1}\right)\right)
											d\tilde{w}\right|\lesssim 1.
										\end{align*}
										Therefore, it holds that for $\tilde{x}\geq x$ and $\tilde{y}\geq 0$,
										\begin{align}\label{est:Koi-case6}
											\left|\left(\f{s-r_c}{\langle r_c\rangle^4}\right)^l(\tilde{x}+\tilde{y})^i\pa_c^l
											\pa_{\tilde{x}}^i(\xi\pa_{\xi})^j
											\left(\xi^{-\f23}(\tilde{x}+\tilde{y})^{-\f13}
											\widetilde{K}_{Oi}(\tilde{x},\tilde{x}+\tilde{y},\xi)\right)\right|
											\lesssim \xi^{-\f43}(x+\tilde{y})^{-\f23}.
										\end{align}
										We also conclude by 
										\eqref{est:Dy-tW} in
										Lemma \eqref{lema:W-trans} and \eqref{bd:q} that for $\widetilde{W}(y,\xi,c):=W(s(y,\xi,c),c)$,
										\begin{align}\notag
											&\left|\left(\f{s-r_c}{\langle r_c\rangle^4}\right)^ly^i\pa_c^l\pa_y^i(\xi\pa_{\xi})^j
											\widetilde{W}(y,\xi,c)\right|\\
											\notag &\lesssim
											\xi^{\f43}y^{-\f43}\mathbf{1}_{y\geq y(\max\{r,2r_c\})}(y)+\f{1}{s^2(y)q(s(y),c)}
											\mathbf{1}_{y(\max\{r,r_c\})
												\leq y\leq y(\max\{r,2r_c\})}(y)\\
											&\lesssim \xi^{\f43}y^{-\f43}\mathbf{1}_{y\geq y(\max\{r,2r_c\})}(y)+r_c^{-2}\langle r_c\rangle^{\f23}\mathbf{1}_{y(\max\{r,r_c\})
												\leq y\leq y(\max\{r,2r_c\})}(y).\label{est:tWcase6}
										\end{align}
										Then letting  $y=\tilde{x}+\tilde{y}$,  we infer that for $y\geq \tilde{x}\geq x\gtrsim 1$,
										\begin{align*}
											&\left|\left(\f{s-r_c}{\langle r_c\rangle^4}\right)^l(\tilde{x}+\tilde{y})^i
											\pa_c^l\pa_{x'}^i(\xi\pa_{\xi})^j
											\left(\widetilde{W}(\tilde{x}+\tilde{y},\xi,c)\right)\right|\\
											&\lesssim \xi^{\f43}(x+\tilde{y})^{-\f43}\mathbf{1}_{x+\tilde{y}\geq y(\max\{r,2r_c\})}(\tilde{y})+r_c^{-2}\langle r_c\rangle^{\f23}\mathbf{1}_{y(\max\{r,r_c\})
												\leq x+\tilde{y}\leq y(\max\{r,2r_c\})}(\tilde{y}),
										\end{align*}
										Collecting the above two derivative bounds and
										using $\f{r-r_c}{\langle r_c\rangle^4}\leq \f{s-r_c}{\langle r_c\rangle^4}$ (due to $r_c\leq r(x)\leq r'(\tilde{x})\leq s(\tilde{x}+\tilde{y})$),
										we finally conclude that for $r\geq r_c$
										and $\xi (r-r_c)^{\f32}/\langle r\rangle^{\f12}\gtrsim 1$,
										\begin{align*}
											&\left(\f{r-r_c}{\langle r_c\rangle^4}\right)^lx^i(r)\sup_{\tilde{x}\in[x,x+\tilde{y}]}\notag
											\left|\pa_c^l\pa_{\tilde{x}}^i(\xi\pa_{\xi})^j\tilde{K}_1(\tilde{x},
											\tilde{x}+\tilde{y},\xi,c)\right|\\
											&\lesssim (x+\tilde{y})^{-2}\mathbf{1}_{x+\tilde{y}\geq y(\max\{r,2r_c\})}(\tilde{y})
											+\xi^{-\f43}r_c^{-2}\langle r_c\rangle^{\f23}(x+\tilde{y})^{-\f23}\mathbf{1}_{x+\tilde{y}\geq y(\max\{r,2r_c\})}(\tilde{y}),
										\end{align*}
										which gives
										\begin{align*}
											&\int_{0}^{+\infty}\left(\f{r-r_c}{\langle r_c\rangle^4}\right)^lx^i(r)\sup_{\tilde{x}\in[x,x+\tilde{y}]}\notag
											\left|\pa_c^l\pa_{\tilde{x}}^i(\xi\pa_{\xi})^j\tilde{K}_1(\tilde{x},
											\tilde{x}+\tilde{y},\xi,c)\right|d\tilde{y}\\
											&\lesssim x^{-1}(\max\{r,2r_c\})+(\xi r_c^{\f32}/\langle r_c\rangle^{\f12})^{-\f43}x^{\f13}(\max\{r,2r_c\})
											\sim \left(\xi \max\{r,2r_c\}^{\f32}/\langle \max\{r,2r_c\}\rangle^{\f12}\right)^{-1}.
										\end{align*}
										Thus, we finish the verification of condition \eqref{condition-trans-tK}.\smallskip
										
										\noindent\emph{Proof of (c)}.
										We take  \begin{align*}I_{1}:=[ r_c-C\xi^{-\f23}\langle r_c\rangle^{\f12}, r_c+C\xi^{-\f23}\langle r_c\rangle^{\f12}]\subset\big[\f{r_c}{2},2r_c\big], \end{align*}
										in Lemma \ref{lem:vorterra-infty}.
										We first verify the condition \eqref{condition:Kinfty'} with $\kappa=\left(\xi r_c^{\f32}/\langle r_c\rangle^{\f12}\right)^{-1}$, for $K_1$ defined in
										\eqref{def:f+6-K1}. Indeed, the estimate for $r\in I_{1}$ is similar as that in \eqref{bd:K1-f+6}:
										\begin{align*}
											\notag&\quad\int_{r}^{+\infty}
											\sup_{\tilde{r}\in[r,s]}|K_1(\tilde{r},s,\xi,c)|ds\lesssim \xi^{-1}\left(\int_{\tau(C_1r_c)}^{+\infty}
											+\int_{\tau(r)}^{\tau(C_1r_c)}\right)|t|^{-\f12}|W(s(t))|dt\lesssim
											\left(\xi r_c^{\f32}/\langle r_c\rangle^{\f12}\right)^{-1}.\notag
										\end{align*}
										For the derivative bounds for the interval $[\f{r_c}{2}, 2r_c]$, we use the variable $\tau=\tau(r,c)$.
										We  conclude by \eqref{bd:q} and \eqref{bd:q-Derivative} that
										\begin{align}\label{bd:q-Derivative-in-f+}
											\left|\left(\f{1}{\langle r_c\rangle^3 }\right)^{l}\langle r_c\rangle ^i\pa_c^l\pa_r^i(q^{\f12}(r,c))\right|\lesssim q^{\f12}(r,c)
											\sim \langle r_c\rangle^{-\f13}, \quad \text{for}\;\; r\in[\f{r_c}{2},2r_c],
											\end{align}
										which along with $\pa_r\tau=q^{\f12}$, $\pa_c\tau=q^{\f12}\f{\int_{r_c}^r\pa_c(Q^{\f12})ds}{Q^{\f12}}$, \eqref{estQ-G} and \eqref{behave:Q<rc0} gives
										\begin{align}\label{est:tau-G}
											&
											\left|\left(\f{1}{\langle r_c\rangle^{3}}\right)^l
											\langle r_c\rangle^i\pa_c^l\pa_r^i\tau(r,c)\right|\lesssim \langle r_c\rangle^{\f23}\quad \text{for}\;\; r\in[\f{r_c}{2},2r_c].
										\end{align}
										We denote $\mathrm{F}^{Rem}(\tau,\xi,c)$ such that $\mathrm{F}^{Rem}(\tau(r,c),\xi,c)=F^{Rem}(r,\xi,c)$. Then
										\begin{align*}
											\pa_{r}F(r,\xi,c)=q^{\f12}
											\pa_{\tau}\mathrm{F}(\tau(r),\xi,c),\;
											\pa_{c}F(r,\xi,c)=\pa_{c}\mathrm{F}(\tau,\xi,c)
											+\pa_c\tau \pa_{\tau}\mathrm{F}(\tau,\xi,c).
											\end{align*}
										Using the above facts, and $\xi^{-\f23}\langle r_c\rangle^{\f13}\lesssim \langle r_c\rangle $,
										we can deduce  that for $r\in I_1(\subset[\f{r_c}{2},2r_c])$ (which is equivalent to $|\tau|\leq C \xi^{-\f23}$), \begin{align*}
											&\left|
											\left(\f{\xi^{-\f23}}{\langle r_c\rangle ^{\f{11}{3}}}\pa_c\right)^l\left(\xi^{-\f23}\langle r_c\rangle^{\f13}\pa_r\right)^i
											(\xi\pa_{\xi})^jF^{Rem}(r,\xi,c)\right|\\
											&\lesssim \sum_{\substack{0\leq i_1\leq i}}\left|
											\left(\f{\xi^{-\f23}}{\langle r_c\rangle ^{\f{11}{3}}}\pa_c\right)^{l}\left(\xi^{-\f23}\pa_{\tau}\right)^{i_1}
											(\xi\pa_{\xi})^j\mathrm{F}^{Rem}(\tau,\xi,c)\right|+\left|
											\left(\xi^{-\f23}\pa_{\tau}\right)^{i_1+l}(\xi\pa_{\xi})^j\mathrm{F}^{Rem}(\tau,\xi,c)\right|.
											\end{align*}
										Therefore, the  derivative bounds  in $(c)$ is reduced to that for $|\tau|\leq C \xi^{-\f23}$ (WLOG, we take $C=1$) and $i+j+l\geq 1$,
										\begin{align}\label{bd:(c)}
											\left|
											\left(\f{\xi^{-\f23}}{\langle r_c\rangle ^{\f{11}{3}}}\pa_c\right)^l\left(\xi^{-\f23}\pa_{\tau}\right)^i
											(\xi\pa_{\xi})^j\mathrm{F}^{Rem}(\tau,\xi,c)\right|
											\lesssim  \left(\xi r_c^{\f32}/\langle r_c\rangle^{\f12}\right)^{-1}.
										\end{align}
										
										To this end, we plug $\mathrm{F}_+(\tau,\xi,c)=\mathrm{Oi}(-\xi^{\f23}\tau)\left(1+\mathrm{F}
										_+^{Rem}(\tau,\xi,c)\right)$ into \eqref{eq:mathrmF+} to obtain
										\begin{align*}
											\pa_{\tau}\left(
											\mathrm{Oi}^2\left(-\xi^{\f23}\tau\right)\pa_{\tau}\mathrm{F}
											^{Rem}\right)=W\mathrm{Oi}^2(-\xi^{\f23}\tau)
											\left(1+\mathrm{F}^{Rem}\right),
										\end{align*}
										which gives the integral form
										\begin{align*}
											\mathrm{F}^{Rem}(\tau,\xi,c)&=
											\mathrm{F}^{Rem}|_{\tau=\xi^{-\f23}}
											-\pa_{\tau}\mathrm{F}^{Rem}|_{\tau=\xi^{-\f23}}
											\mathrm{Oi}^2(-1)\int_{\tau}^{\xi^{-\f23}}
											\mathrm{Oi}^{-2}(-\xi^{\f23}t')dt'\\
											&\quad+ \int_{\tau}^{\xi^{-\f23}}\left(\mathrm{Oi}^2(-\xi^{\f23}t)
											\int_{\tau}^{t}
											\mathrm{Oi}^{-2}(-\xi^{\f23}t')dt'\right)\cdot W\left(s(t,c),c\right)\left( 1+\mathrm{F}^{Rem}(t,\xi,c)\right)dt.
											\end{align*}
										We will use this integral form to estimate the derivative bounds \eqref{bd:(c)} near $r_c$ rather than using \eqref{def:F+case6inf+}. For simplicity, we rewrite it as
										\begin{align} \label{new:integral }
											\mathrm{F}^{Rem}(\tau,\xi,c)&=\mathrm{F}_0(\tau,\xi,c)+
											\int_{\tau}^{\xi^{-\f23}}\mathrm{K}_1(\tau,t,\xi,c)dt+
											\int_{\tau}^{\xi^{-\f23}}\mathrm{K}_1(\tau,t,\xi,c)
											\mathrm{F}^{Rem}(t,\xi,c)dt, \end{align}
										where   $|\tau|\leq \xi^{-\f23}$ and
										\begin{align*}
											& \mathrm{F}_0(\tau,\xi,c)
											:=\mathrm{F}^{Rem}|_{\tau=\xi^{-\f23}}
											-\pa_{\tau}\mathrm{F}^{Rem}|_{\tau=\xi^{-\f23}}K_{{Oi}}(
											\tau,-\xi^{-\f23},\xi),\;\;
											\mathrm{K}_1(\tau,t,\xi,c):=K_{{Oi}}(
											\tau,t,\xi)\widetilde{\widetilde{W}}(t,c),\\
											&K_{{Oi}}(
											\tau,t,\xi)=  \mathrm{Oi}^2(-\xi^{\f23}t)
											\int_{\tau}^{t}
											\mathrm{Oi}^{-2}(-\xi^{\f23}t')dt',\;\widetilde{\widetilde{W}}(t,c)=W\left(s(t,c),c\right)\;\;-\xi^{-\f23}\leq \tau\leq t\leq \xi^{-\f23}. \end{align*}
										Formally, the solution of   \eqref{new:integral } takes
										\begin{align*}
											&\mathrm{F}^{Rem}(\tau,\xi,c)=\sum_{n=1}^{\infty}h_n(\tau,\xi,c)
											,\;\; \text{where}\;\;
											h_1=\mathrm{F}_0(\tau,\xi,c)+
											\int_{\tau}^{\xi^{-\f23}}\mathrm{K}_1(\tau,t,\xi,c)dt,\\
											&h_{n}=
											\int_{\tau}^{\xi^{-\f23}}\mathrm{K}_1(\tau,t,\xi,c)h_{n-1}(t,\xi,c)dt\,\,(n\geq 2).
											\end{align*}
										Similar as the proof of Lemma \ref{lem:vorteraa-0} with condition \eqref{condition:K'pac-growth},
										if we verify the following conditions on $\mathrm{F}_0$ and $\mathrm{K}_1$ in \eqref{new:integral }, for $|\tau|\leq \xi^{-\f23}$:
										\begin{align}
											\label{condtion:tauF0}&\sup_{\tilde{s}\in[0,1]}\left|
											\left(\f{\xi^{-\f23}}{\langle r_c\rangle ^{\f{11}{3}}}\pa_c\right)^l \left(\xi^{-\f23}\pa_{\tau}\right)^i
											(\xi\pa_{\xi})^j\mathrm{F}_0(\tilde{s}(\xi^{-\f23}-\tau)+\tau,\xi,c)\right|
											\lesssim \left(\xi r_c^{\f32}/\langle r_c\rangle^{\f12}\right)^{-1},\\
											&\int_0^1\sup_{\tilde{s}\in[0,s]}\left|
											\left(\f{r_c}{\langle r_c\rangle ^{4}}\pa_c\right)^l\left(\xi^{-\f23}\pa_{\tau}\right)^i
											(\xi\pa_{\xi})^j\left((\xi^{-\f23}-\tau)
											\mathrm{K}_1(\tau+\tilde{s}(\xi^{-\f23}-\tau),
											\tau+s(\xi^{-\f23}-\tau),\xi,c)\right)\right|ds\notag \\
											&\lesssim \left(\xi r_c^{\f32}/\langle r_c\rangle^{\f12}\right)^{-\f43},\label{condtion:tauK1}
											\end{align}
											 where $\f{\xi^{-\f23}}{\langle r_c\rangle ^{\f{11}{3}}}\lesssim \f{r_c}{\langle r_c\rangle ^{4}}$,
											 $\left(\xi r_c^{\f32}/\langle r_c\rangle^{\f12}\right)^{-\f43}\lesssim \left(\xi r_c^{\f32}/\langle r_c\rangle^{\f12}\right)^{-1}$,  then we can deduce that
										\begin{align*}
											\left|
											\left(\f{\xi^{-\f23}}{\langle r_c\rangle ^{\f{11}{3}}}\pa_c\right)^l\left(\xi^{-\f23}\pa_{\tau}\right)^i
											(\xi\pa_{\xi})^jh_n(\tau,\xi,c)\right|
											\leq  \f{C^n_{i,j,l}\left(\xi r_c^{\f32}/\langle r_c\rangle^{\f12}\right)^{-n}}{(n-1)!}\;\;n\geq 1,
											\end{align*}
										and finally obtain \eqref{bd:(c)}.
										
										It remains to prove \eqref{condtion:tauF0} and  \eqref{condtion:tauK1}. To prove \eqref{condtion:tauF0}, we let $r=r_c+C\xi^{-\f23}\langle r_c\rangle^{\f13}$ (equivalently, $\tau=\xi^{-\f23}$) in (b) to obtain
										\begin{align*}
											\left|
											\left(\left(\f{\xi^{-\f23}\langle r_c\rangle^{\f13}}{\langle r_c\rangle ^{4}}\right)^l(\xi^{-\f23}\langle r_c\rangle^{\f13})^i
											\pa_c^l\pa_r^i(\xi\pa_{\xi})^jF^{Rem}\right)(r_c+C\xi^{-\f23}\langle r_c\rangle^{\f13},\xi,c)\right|
											\lesssim (\xi r_c^{\f32}/\langle r_c\rangle^{\f12})^{-1},
											\end{align*}
										which along with $\mathrm{F}^{Rem}(\tau,\xi,c)=F^{Rem}(r(\tau,c),\xi,c)$, $\pa_{\tau}\mathrm{F}(\tau)=q^{-\f12}\pa_{r}F(r),\;
										\pa_{c}\mathrm{F}(\tau,c)=\pa_{c}F(r,c)
										-q^{-\f12}\pa_c\tau \pa_{r}F(r,c)$, \eqref{bd:q-Derivative-in-f+}, \eqref{est:tau-G} and $\xi^{-\f23}\langle r_c\rangle^{\f13}\lesssim \langle r_c\rangle$, shows
										\begin{align*}
											\left|
											\left(\left(\f{\xi^{-\f23}}{\langle r_c\rangle ^{\f{11}{3}}}\right)^l(\xi^{-\f23})^i
											\pa_c^l\pa_{\tau}^i(\xi\pa_{\xi})^j\mathrm{F}^{Rem}\right)(\xi^{-\f23},\xi,c)\right|
											\lesssim (\xi r_c^{\f32}/\langle r_c\rangle^{\f12})^{-1}.
											\end{align*}
										We also use \eqref{asy:tKOi-small} to obtain $K_{\mathrm{Oi}}(
										\tau,-\xi^{-\f23},\xi)= O_{\xi,\tau}^{\xi,\xi^{-\f23}}(\xi^{-\f23})$ for $|\tau|\leq  \xi^{-\f23}$. Applying the above two bounds in the definition of $\mathrm{F}_0$, we arrive at \eqref{condtion:tauF0}. To prove \eqref{condtion:tauK1}, we first notice that
										\begin{align*}
											\left|\left(\xi^{-\f23}\pa_{\tau}\right)^i
											(\xi\pa_{\xi})^j\left(\tau+s'(\xi^{-\f23}-\tau)\right)\right|\leq C \xi^{-\f23},
										\end{align*}
										for $C$ independent of $s'\in[0,1]$. Then \eqref{condtion:tauK1} can be reduced to the following:
										\eqref{asy:tKOi-small} in Lemma \ref{lem:tK},
										\begin{align*}
											&K_{{Oi}}(\tau,t,\xi)= O_{\xi,\tau,t}^{\xi,\xi^{-\f23},\xi^{-\f23}}(\xi^{-\f23})\;\;\text{for}\;\;-\xi^{-\f23}\lesssim \tau\leq t\lesssim \xi^{-\f23},
										\end{align*}
										and
										\eqref{est:Dy-ttW} in  Lemma \ref{lema:W-trans},
										\begin{align*} &\left|\left(\f{r_c}{\langle r_c\rangle ^4}\pa_c\right)^l(\xi^{-\f23}\pa_{t})^i\widetilde{\widetilde{W}}
											(t,c)\right|
											\lesssim r_c^{-2}\langle r_c\rangle^{\f23}\mathbf{1}_{|t|\lesssim \xi^{-\f23}}(t).
										\end{align*}
										
										This finishes the proof of the lemma.
							\end{proof}
							\section{Behavior of Schr\"odinger equation on the intermediate region}

							To characterize the behavior of the Wronskian $\mathcal{W}(\phi,f_+)(\xi,c)$, we need to investigate the behavior of the solution $\phi$ to \eqref{eq:phi} on an intermediate region
							that overlaps with the maximal interval $[r_{\infty} (\xi,c),+\infty)$, where the precise behavior of $f_+(r,\xi,c)$ is manifested.

							\subsection{Case of $c\in(V(0),1)$ and $\xi\gtrsim M(1-c)^{\f13}|V(0)-c|^{-\f32}$}

							Let's begin with some technical estimates, which are necessary for handling the matching of initial conditions for the ODE in Lemma \ref{lem:f+3}.
							
							\begin{lemma}\label{lem:b2d2}
								Let $M\gg 1$ be fixed, $\xi r_c^{\f32}/\langle r_c\rangle^{\f12}\gtrsim M$,
								$r_*(\xi,c)=M^{\f12}\xi^{-\f12}\left(\f{1-c}{c-V(0)}\right)^{\f12}$, $x=\xi(-\tau)^{\f32}(r,c)$($\tau$ defined as \eqref{def:tau}), $x_*(\xi,c)=x|_{r=r_*}$.
								Then for $r_*\leq \min\{r_c/2,1\}$, it holds that
								\begin{align}
									\begin{split}
										\label{est:r*x*}
										&\left|\left(\f{r_c}{\langle r_c\rangle^4}\pa_c\right)^l(\xi\pa_{\xi})^jr_*\right|\lesssim r_*\sim M^{\f12}\xi^{-1}\langle r_c\rangle^{-1}r_c^{-\f12},\\
										&\left|\left(\f{r_c}{\langle r_c\rangle^4}\pa_c\right)^l(\xi\pa_{\xi})^jx_*\right|\lesssim x_*\sim \xi r_c^{\f32}/\langle r_c\rangle^{\f12}.
									\end{split}
								\end{align}
								For $r_*\leq r\leq r_c$,
								\begin{align}\label{est:x*-x}
									&
									\left|\left(\f{r_c}{\langle r_c\rangle^4}\pa_c\right)^l(\xi\pa_{\xi})^j\left(x_*-x\right)\right|\lesssim (x_*-x).
								\end{align}
								Moreover, let
								\begin{align}\label{def:B2D2}
									\begin{split}
										&B_1(\xi,c)=-x_*^{-\f13}
										\left(\f{\mathrm{Ai}'(x_*^{\f23})}{\mathrm{Ai}(x_*^{\f23})}
										+\f{\xi^{-\f23}q^{-\f12}\pa_{r}\left(q^{\f14}r^{\f12}\phi(r)\right)
										}{q^{\f14}r^{\f12}\phi(r)}\Big|_{r=r_*}\right),\\
										&D_1(z,x_*)=\f23x_*^{\f13} \mathrm{Ai}^{2}(x_*^{\f23})
										\int^{x_*}_{x_*-z}\mathrm{Ai}^{-2}(w)w^{-\f13}dw.
									\end{split}
								\end{align}
								Then it holds  that
								\begin{align}
									\begin{split}
										\label{est:B2}
										&\f34\leq B_1(\xi,c)\leq \f54\quad\text{if}\;\;M\gg 1,\\
										&\left|\left(\f{r_c}{\langle r_c\rangle^4}\pa_c\right)^l
										(\xi\pa_{\xi})^jB_1(\xi,c)\right|\lesssim M^{-
											\f12},
									\end{split}
								\end{align}
								and for $z\in [0,x_*]$,
								\begin{align}
									\begin{split}\label{est:positivetD2}
										&0< D_1(z,x_*)\leq 
										\f34\;\;\;\text{if}\;\;M\gg 1,\\
										&\left| (x_*\pa_{x_*})^j\left(\langle z\rangle\pa_z\right)^{i}
										D_1(z,x_*)\right|\lesssim 1,\;x_*-z\gtrsim 1.
									\end{split}
								\end{align}
								As a consequence, letting $r_*\leq s\leq r\leq r_c$,  it holds that
								\begin{align}
										&\notag1+\f23\left(\xi^{-\f23}\f{q^{-\f12}\pa_{r}
											\left(q^{\f14}r^{\f12}\phi(r)\right)
										}{q^{\f14}r^{\f12}\phi(r)}\Big|_{r=r_*}+
										\f{\mathrm{Ai}'(x_*^{\f23})}{\mathrm{Ai}(x_*^{\f23})}\right)\mathrm{Ai}^2(x_*^{\f23})\int^{x_*}_{x}\mathrm{Ai}^{-2}(w^{\f23})w^{-\f13}dw\\
										&=1-B_1(\xi,c)D_1(x_*-x,x_*).\label{equal:1+B2D2}
									\end{align}
									Moreover, letting   $1\lesssim x\lesssim x_*$  and $x_*\gtrsim \xi r_c^{\f32}/\langle r_c\rangle^{\f12}\gtrsim M\gg 1$, it holds that
									\begin{align}\label{est1+B2D2-1/16} 1-B_1(\xi,c)D_1\big(x_*-x,x_*\big)\in(\f{1}{16},1),
									\end{align}
									and
									\begin{align}\label{est1+B2D2-x}
										\left|\left(\f{r_c-r}{\langle r_c\rangle^4}\right)^l\left(\f{(r_c-r)r}{\ r_c}\right)^i\pa_c^l\pa_r^i(\xi\pa_{\xi})^j
										\left(1-B_1(\xi,c)D_1\left((x_*(\xi,c)-x(r,\xi,c))
										,x_*(\xi,c)\right)\right)\right|
										\lesssim 1.
									\end{align}
					\end{lemma}
							
							\begin{proof}
								Noticing that $V(r)=1-a_0r^{-3}+O(r^{-4}) (r\to+\infty)$ and $V'(r)\geq c_0>0$, it holds that for $c\in(V(0),1)$,
								\begin{align*}
									1-c\sim \langle r_c\rangle ^{-3}, \;c-V(0)\sim \f{r_c}{\langle r_c\rangle},
								\end{align*}
								and by $\xi r_c^{\f32}/\langle r_c\rangle^{\f12}\gtrsim M$, we can choose $M\gg 1$ such that
								\begin{align}\label{r*<rc/2}
									M^{\f12}r_c^{-\f12}\langle r_c\rangle^{-1}\sim r_*\in(0,\min\{r_c/2,1\}).
								\end{align}
								Then the first bound in \eqref{est:r*x*} follows directly.
								The second bound in \eqref{est:r*x*} follows from the first bound, \eqref{r*<rc/2} and \eqref{estx-G}. The first bound in \eqref{est:x*-x} follows from \eqref{estx-G}, and for the second one, using the definition of $x$, \eqref{est:r*x*}, \eqref{r*<rc/2} and Lemma \ref{lem:behave-Q}, it suffices to check that for $l=1,j=0$, $r\in[r_*,r_c]$,
								\begin{align*}
									\f{\langle r_c\rangle(r_c-r_*)^{\f12}}{\langle r_*\rangle^{\f32}}
									+\int_{r_*}^r\f{r_c\langle r_c\rangle}{\langle s\rangle^{\f32}
										(r_c-s)^{\f12}}ds&\sim |Q(r_*,c)|^{\f12}+\f{r_c}{\langle r_c\rangle^4}\int_{r_*(c)}^r|\pa_cQ(s,c)|\cdot |Q(s,c)|^{-\f12}ds\\
									&\lesssim \int_{r_*(c)}^r|Q(s,c)|^{\f12}ds\sim
									\int_{r_*}^r\f{(r_c-s)^{\f12}\langle r_c\rangle}{\langle s\rangle^{\f32}}ds\\
									&\sim \langle r_c\rangle r_c^{\f12}+
									\langle r_c\rangle ^{-\f12}r_c^{\f32}\sim
									\langle r_c\rangle r_c^{\f12},
								\end{align*}
								where the integrals can be estimated separately by $\int_{r_*}^{r_c/2}ds$ and
								$\int^{r}_{r_c/2}ds$.
								
								To prove \eqref{est:B2}, noticing $x_*\gtrsim M$, it follows from  \eqref{behave:Ai} and \eqref{behave:Ai^-1} that
								\begin{align*}
									-x_*^{-\f13}\f{\textrm{Ai}'(x_*^{\f23})}{\textrm{Ai}(x_*^{\f23})}=1+O(x_*^{-1}).
								\end{align*}
								A direct calculation yields that
								\begin{align}
									&
									-x_*^{-\f13}\f{\xi^{-\f23}q^{-\f12}\pa_{r}\left(q^{\f14}r^{\f12}\phi(r)\right)
									}{q^{\f14}r^{\f12}\phi(r)}\Big|_{r=r_*}=
									\label{fm:tauPhi/Phi} -\left(\xi r (-Q)^{\f12}|_{r=r_*}\right)^{-1}
									\left(\f12
									+\f{r\pa_r\phi}{\phi}\big|_{r=r_*}+\f{r\pa_rq}{4q^{\f12}}\big|_{r=r_*}
									\right).
								\end{align}
								It follows from \eqref{r*<rc/2} and Lemma \ref{lem:phi-3} that
								\begin{align}\label{phir*}
									\begin{split}
										\xi r_*<\phi(r_*),r_*\pa_r\phi(r_*,\xi,c)\lesssim \xi r_*,\;\left|\left(\f{r_c}{\langle r_c\rangle^4}\pa_c\right)^l(r\pa_r)^i
										(\xi\pa_{\xi})^j\phi(r,\xi,c)|_{r=r_*}\right|\lesssim \xi r_*,
									\end{split}
								\end{align}
								which along with \eqref{bd:q-Derivative} and \eqref{est:r*x*} gives
								\begin{align*}
									&\left|\f12
									+\f{r_*\pa_r\phi(r_*,\xi,c)}{\phi(r_*,\xi,c)}+\f{r\pa_rq(r_*,c)}
									{4q^{\f12}(r_*,c)}
									\right|\lesssim \f12+1+\f14\lesssim 1,\\
									&\left|\left(\f{r_c}{\langle r_c\rangle^4}\pa_c\right)^l
									(\xi\pa_{\xi})^j\left(\f12
									+\f{r_*\pa_r\phi(r_*,\xi,c)}{\phi(r_*,\xi,c)}+\f{r\pa_rq(r_*,c)}
									{4q^{\f12}(r_*,c)}\right)
									\right|\lesssim1.
								\end{align*}
								The fact $(-Q)^{\f12}(s,c)\sim \f{\langle r_c\rangle(r_c-s)^{\f12}}{\langle r\rangle^{\f32}}(r\leq r_c)$ and \eqref{r*<rc/2} lead to
								\begin{align*}
									\left(\xi r (-Q)^{\f12}|_{r=r_*}\right)^{-1}\sim
									\left(\xi \langle r_c\rangle \f{r_*}{\langle r_*\rangle^{\f32}} (r_c-r_*)^{\f12}\right)^{-1}\sim \left(\xi \langle r_c\rangle r_* r_c^{\f12}\right)^{-1}\sim M^{-\f12}.
								\end{align*}
								We conclude the above estimate from the definition \eqref{def:B2D2} that for $x_*\gtrsim M\gg 1$,
								\begin{align*}
									B_1(\xi,c)=1+O(x_*^{-1})+O(M^{-\f12})\in \big[\f34,\f54\big].
								\end{align*}
								Thanks to \eqref{behave:Ai}, \eqref{behave:Ai^-1}, \eqref{est:r*x*}, \eqref{bd:q-Derivative} and  \eqref{est-der-Q}, the derivative bounds follow from
								\begin{align*}
									&\left|\left(\f{r_c}{\langle r_c\rangle^4}\pa_c\right)^l
									(\xi\pa_{\xi})^j\left(x_*^{-\f13}\f{\mathrm{Ai}'(x_*^{\f23})}{\mathrm{Ai}(x_*^{\f23})},\;
									\left(\xi r_* (-Q(r_*,c)^{\f12})\right)^{-1}\right)
									\right|\lesssim M^{-
										\f12}.
								\end{align*}
								Thus, we finish the proof of \eqref{est:B2}. Now we prove \eqref{est:positivetD2}. $D_1>0$ is obvious.
								It follows from the definition and behavior of $\mathrm{Ai}$ in \eqref{behave:Ai}, \eqref{behave:Ai^-1} that
								\begin{align*}
									D_1(z,x_*)&=\f23\left(1+O(x_*^{-1})\right)
									\int_{x_*-z}^{x_*}e^{-\f43(x_*-w)}\left(\f{\langle w\rangle}{w}\right)^{\f13}
									\left(1+O\left(\langle w\rangle^{-1}\right)\right)dw.
								\end{align*}
								Therefore, for $z\in[0,x_*)$, if $x_*\geq C^{-1}M$ for some $M\gg 1$, it follows that
								\begin{align*}
									& D_1(z,x_*)\leq D_1(x_*,x_*)\\
									&\leq \f23\left(1+CM^{-1}\right)
									\int_{x_*/2}^{x_*}e^{-\f43(x_*-w)}
									\left(1+CM^{-1}\right)^{\f43}dw+
									C
									\int_{0}^{x_*/2}e^{-\f43(x_*-x_*/2)}\left(w^{-\f13}+1\right)dw\\
									&\leq \f12\left(1+CM^{-1}\right)^{\f73}\left(1-e^{-\f23x_*}\right)
									+Cx_*e^{-\f23x_*}\leq \f12\left(1+CM^{-1}\right)^{\f73}
									+CMe^{-C^{-1}M}\leq \f34.
								\end{align*}
								For the derivative bounds in \eqref{est:positivetD2}, we take $\tilde{w}=x_*-w$ to get
								\begin{align*}
									D_1(z,x_*)=\f23\left(1+O(x_*^{-1})\right)
									\int_{0}^{z}e^{-\f43\tilde{w}}\left(\f{\langle x_*-\tilde{w}\rangle}{x_*-\tilde{w}}\right)^{\f13}\left(1+O\left(\langle x_*
									-\tilde{w}\rangle^{-1}\right)\right)d\tilde{w}.
								\end{align*}
								If $i=0$, we first apply the derivative $(x_*\pa_{x_*})^j$ to the formula above, and then use $x_*\leq (x_*-\tilde{w})+\tilde{w}$. Substituting back $w=x_*-\tilde{w}$,
								the estimate follows similarly as before. While for $i\geq 1 $, using  $x_*\leq (x_*-z)+z$,  it follows that for $x_*-z\gtrsim 1$,
								\begin{align*}
									&\left|(x_*\pa_{x_*})^j \left(\langle z\rangle\pa_z\right)^{i}
									D_1(z,x_*)\right|\\
									&=
									\left|(x_*\pa_{x_*})^j \left(\langle z\rangle\pa_z\right)^{i-1}
									\left((z+1)\cdot
									\f23\left(1+O(x_*^{-1})\right)
									e^{-\f43z}\left(\f{\langle x_*-z\rangle}{x_*-z}\right)^{\f13}\left(1+O(\langle x_*
									-z\rangle^{-1})\right)\right)\right|\\
									&\lesssim \sum_{1\leq k\leq i+j}\langle z\rangle^{k}e^{-\f43z}\left(1+(x_*-z)^{\f23-k}\right)\lesssim 1.
									\end{align*}
									This finishes the proof of \eqref{est:positivetD2}, and
									\eqref{est1+B2D2-x} is the consequence of \eqref{est:x*-x} and \eqref{est:positivetD2}.
							\end{proof}

							\begin{lemma}\label{lem:f+3}
								Let  $c\in(V(0),1)$, $\xi r_c^{\f32}/\langle r_c\rangle^{\f12}\geq M$. Let $M\gg 1$ such that  $r_*=M^{\f12}\xi^{-1}\left(\f{1-c}{c-V(0)}\right)^{\f12}
								\in(0,\min\{r_c/2,1\})$. Let $\phi$ be the solution to \eqref{eq:phi}, i.e.,
								\begin{align*}
									\phi'' +r^{-1}\phi'-r^{-2}\phi+\f{\xi^2(V-c)}{1-c}  \phi=0,\; \phi(r,\xi,c)\sim \xi r,
								\end{align*}
								which has been solved on $(0, r_*)$ in Lemma \ref{lem: phi-2}.
								Let $Q(r,c)$ be defined in \eqref{def:Q}, $\tau=\tau(r,c)$ be defined as \eqref{def:tau}, $q=\f{Q}{\tau}$,
								$\tau_*=\tau(r_*,c)$.
								Then there exist  complex numbers $C_{61}(\xi,c)$, $\widetilde{C}_{61}(\xi,c)$ defined in \eqref{coefficient6}
								such that for
								$M^{\f12}\xi^{-1}\langle r_c\rangle^{-1}r_c^{-\f12}\sim r_*\leq r\leq r_c-C\xi^{-\f23}\langle r_c\rangle^{\f13}$,  it holds that
								\begin{align}
									\notag
									&\phi(r,\xi,c)=C_{61}\xi^{-\f13}q^{-\f14}r^{-\f12}\mathrm{Ai}(-\xi^{\f23}\tau)
									\left(1+\widetilde{C}_{61}\int_{\tau_*(\xi,c)}
									^{\tau(r,c)}\mathrm{Ai}^{-2}(-\xi^{\f23}t')dt'\right)
									\left(1+\phi^{Rem}(r,\xi,c)\right)\\
									&\text{with}
									\;\;  |\phi^{Rem}|
									\lesssim M^{-\f12}\;\;\text{and}\;\; \left|\left(\f{r_c-r}{\langle r_c\rangle^4}\right)^l\left(\f{(r_c-r) r}{ r_c}\right)^i\pa_c^l\pa_r^i(\xi\pa_{\xi})^j\phi^{Rem}(r,\xi,c)\right|
									\lesssim 1\label{behave:m-phi+1},
								\end{align}
								and
								\begin{align}  \label{1+tC61}
									\left(1+\widetilde{C}_{61}\int_{\tau_*(\xi,c)}
									^{\tau(r,c)}\mathrm{Ai}^{-2}(-\xi^{\f23}t')dt'\right)\in \big(\f{1}{16},1\big),
								\end{align}
								and
								\begin{align} 0<C_{61}^{-1}(\xi,c)= O_{\xi,c}^{\xi,\f{r_c}{\langle r_c\rangle^4}}\left(r_c^{\f12}\langle
									r_c\rangle^{2}\cdot(\xi r_c^{\f32}/\langle r_c\rangle^{\f12})^{-N}
									\right),\quad N\geq 1.\label{est:C1-1}
								\end{align}
								Moreover, for
								$M^{\f12}\xi^{-1}\langle r_c\rangle^{-1}r_c^{-\f12}\sim r_*\leq r
								\leq r_c-C\xi^{-\f23}\langle r_c\rangle^{\f13}$, it holds that
								\begin{align} \;\;C_{61}^{-1}\phi(r,\xi,c) &=O_{\xi,r,c}^{\xi,\f{(r_c-r)r}{r_c},
										\f{r_c-r}{\langle r_c\rangle^4}}\left(\left(\xi \langle r_c\rangle \f{r}{\langle r\rangle^{\f32}}(r_c-r)^{\f12}\right)^{-\f12}\right)
									\label{behave:phi-6-itself},
									\end{align}
								and
								\begin{align} \;\;C_{61}^{-1}\phi(r,\xi,c) &\gtrsim \left(\xi \langle r_c\rangle \f{r}{\langle r\rangle^{\f32}}(r_c-r)^{\f12}\right)^{-\f12}. \label{behave:phi-6-itself-lower}
								\end{align}
							\end{lemma}
							\begin{proof}
							Noticing $\pa_{\tau}=q^{-\f12}\pa_r$, we want to extend $\phi$ smoothly  for $r\geq r_*$. It  suffices to match the initial condition at $r=r_*$. We perform the Langer transform $\widetilde{\Phi}(\tau,\xi,c)=\Phi(r,\xi,c):=\xi^{\f13}q^{\f14}r^{\f12}\phi(r,\xi,c)$, which satisfies the equation
							\begin{align}\label{eq:Phi1}
							\left\{
							\begin{aligned}
								&\pa_{\tau}^2\widetilde{\Phi}+\xi^2\tau\widetilde{\Phi}
								=W\widetilde{\Phi},\;\;W
								=q^{-\f14}\pa_{\tau}^2(q^{\f14})+\f{3}{4r^2q},\\
								&\widetilde{\Phi}(\tau_*)=\xi^{\f13}q^{\f14}r^{\f12}\phi|_{r=r_*},\;
								\pa_{\tau}\widetilde{\Phi}(\tau_*)=q^{-\f12}
								\f{d}{dr}\left(\xi^{\f13}q^{\f14}r^{\f12}\phi(r)\right)|_{r=r_*}.
							\end{aligned}
							\right.
						\end{align}
						Direct calculations  give the equivalent integral form of $\Phi(r)$:
						\begin{align}
							\label{integral:Philanger1}
							\Phi(r,\xi,c)&=\tilde{\textrm{Ai}}\left(\tau,\xi,c\right)
							+\int_{r_*(\xi,c)}^{r}K(r,s,\xi,c)\Phi(s,\xi,c)ds,
						\end{align}
						where
						$$K(r,s,\xi,c)=\left(\tilde{\mathrm{Ai}}\left(\tau(r,c),\xi,c\right)
						\tilde{\textrm{Ai}}\left(t(s),\xi,c\right)\int_{t(s,c)}^{\tau(r,c)}
						\tilde{\textrm{Ai}}^{-2}(-\xi^{\f23}t')dt'\right)W(s,c)\pa_st(s,c).$$
						Here $\tilde{\textrm{Ai}}(\tau)$ is the solution of homogenous
						\eqref{eq:Phi1}($W\equiv 0$):
						\begin{align}
							\label{tAi2}&\tilde{\textrm{Ai}}(\tau,\xi,c)=C_{61}(\xi,c)\textrm{Ai}(-\xi^{\f23}\tau)\left(
							1+\widetilde{C}_{61}(\xi,c)\int_{\tau_*}^{\tau}\textrm{Ai}^{-2}(-\xi^{\f23}t)dt\right)\\
							&=_{x=\xi(-\tau)^{\f32},w=\xi(-t')^{\f32}}C_{61}(\xi,c)\textrm{Ai}(x^{\f23})\left(
							1+\f23\xi^{-\f23}\widetilde{C}_{61}(\xi,c)\int_{x}^{x_*(\xi,c)}\textrm{Ai}^{-2}(w^{\f23})
							w^{-\f13}dw\right)\notag,
						\end{align}
						with the coefficients
						\begin{align}
							\label{coefficient6}
							\begin{split}
							C_{61}(\xi,c)&=\f{\xi^{\f13}q^{\f14}r^{\f12}\phi(r)\big|_{r=r_*}}
							{\textrm{Ai}(x_*^{\f23})},\\
							\widetilde{C}_{61}(\xi,c)&=
							\xi^{\f23}\left(\xi^{-\f23}\f{q^{-\f12}\pa_{r}\left(q^{\f14}r^{\f12}
								\phi(r)\right)
							}{q^{\f14}r^{\f12}\phi(r)}\Big|_{r=r_*}+
							\f{\textrm{Ai}'(x_*^{\f23})}{\textrm{Ai}(x_*^{\f23})}\right)\textrm{Ai}^2(x_*^{\f23}).
							\end{split}
						\end{align}
						
						We claim that
						\begin{align}\notag
							&\Phi(r,\xi,c)=\tilde{\textrm{Ai}}\left(\tau(r,c),\xi,c\right)
							\left(1+\Phi^{Rem}(r,\xi,c)\right)\;\;\text{for}\;\;r_*\leq r\leq r_c-C\xi^{-\f23}\langle r_c\rangle^{\f13},\;\text{with}
					\end{align}
					\begin{align}
							\label{tAi4}
							\begin{split}
							&|\Phi^{Rem}|
							\lesssim  M^{-\f12},\quad \left|\left(\f{r_c-r}{\langle r_c\rangle^4}\right)^l\left(\f{(r_c-r) r}{ r_c}\right)^i\pa_c^l\pa_r^i(\xi\pa_{\xi})^j\Phi^{Rem}(r,\xi,c)\right|
							\lesssim 1,\\
							&\Phi^{Rem}|_{r=r_*}=\pa_{r}\Phi^{Rem}|_{r=r_*}
							=0.
							\end{split}
						\end{align}
						Thus, \eqref{behave:m-phi+1} follows directly from $\phi=\xi^{-\f13}q^{-\f14}r^{-\f12}\Phi$; \eqref{1+tC61} follows from \eqref{equal:1+B2D2} and \eqref{est1+B2D2-1/16};
						\eqref{est:C1-1} follows from the definition of $q$, $x$, $\mathrm{Ai}(z)=CO(z^{-\f14})e^{-\f23z^{\f32}}(z\gtrsim 1)$,  \eqref{est:r*x*}, \eqref{behave:Q<rc} ($|Q(r_*)|\sim \f{(r_c-r_*)\langle r_c\rangle^2}{\langle r_*\rangle^3}\sim r_c\langle r_c\rangle^2$), $x_*\sim \xi r_c^{\f32}/\langle
						r_c\rangle^{\f12}\gg 1$ and \eqref{estQ-G}  that
						\begin{align*} \notag C_{61}^{-1}
							&=\f{C\textrm{Ai}(x_*^{\f23})}{\xi^{\f13}q^{\f14}(r)r^{\f12}\xi r|_{r=r_*}}
							=C\xi^{-\f43}r_*^{-\f32}|Q(r_*)|^{-\f14}|\tau(r_*)|^{\f14}O(x_*^{-\f16})e^{-\f23x_*}\\
							&=CO_{\xi,c}^{\xi,\f{r_c}{\langle r_c\rangle^4}} (\xi^{-\f32}r_*^{-\f32}|Q(r_*)|^{-\f14})e^{-\f23x_*}=
							O_{\xi,c}^{\xi,\f{r_c}{\langle r_c\rangle^4}} (r_c^{\f12}\langle r_c\rangle)\cdot O_{\xi,c}^{\xi,\f{r_c}{\langle r_c\rangle^4}} \left((\xi r_c^{\f32}/\langle r_c\rangle^{\f12})^{-N}\right).
							\end{align*}
							To show \eqref{behave:phi-6-itself}, we calculate by the definitions of $q$, $B_1$, $D_1$ and $\mathrm{Ai}(z)=Cz^{-\f14}e^{-\f23z^{\f32}}(z\gtrsim 1)$ that
							\begin{align*}
								\;\;C_{61}^{-1}\phi(r,\xi,c)&=\left(\xi r(-Q)^{\f12}\right)^{-\f12}
								e^{-\f23x(r,\xi,c)}
								\left(1+O(x^{-1})\right)\left(1-B_1(\xi,c)D_1(x_*-x,x_*)\right)
								\left(1+\Phi^{Rem}\right),
							\end{align*}
							where $1\lesssim x(r,\xi,c)=\xi(-\tau)^{\f32}\gtrsim \xi (r_c-r)^{\f32}/\langle r\rangle^{\f12}$ is defined in \eqref{def:x}. Notice by \eqref{behave:Q<rc}, \eqref{estQ-G} that for $r\leq r_c$, it holds  that
							\begin{align*}
								\left(\xi r(-Q)^{\f12}\right)^{-\f12}=
								O_{\xi,c,r}^{\xi,
									\f{r_c-r}{\langle r_c\rangle^4}, \f{(r_c-r)r}{r_c}}\left(\left(\xi \langle r_c\rangle \f{r}{\langle r\rangle^{\f32}}(r_c-r)^{\f12}\right)^{-\f12}\right).
								\end{align*}
								Then the bound \eqref{behave:phi-6-itself} is a consequence of \eqref{estQ-G},  \eqref{estx-G}, \eqref{est1+B2D2-1/16}, \eqref{est1+B2D2-x}, $x^ie^{-\f23x}\lesssim 1$
							and \eqref{behave:m-phi+1}.  The lower bound \eqref{behave:phi-6-itself-lower} follows further by \eqref{1+tC61} and $|\phi^{Rem}|\lesssim M^{-\f12}$.
							
							It remains to prove \eqref{tAi4}. We apply Lemma \ref{lem:vorteraa-0} with $f=\Phi$,
							$g=\tilde{\mathrm{Ai}}$, $D=\big\{(\xi,c)|c\in(V(0),1),\;\xi r_c^{\f32}/\langle r_c\rangle^{\f12}\gtrsim M\big\}$, $I_1=[r_0,r_1]:=[r_*,r_c-C\xi^{-\f23}\langle r_c\rangle^{\f13}]$,
								and
								\begin{align*}&K_1(r,s,\xi,c)=\tilde{\mathrm{Ai}}^{-1}(x(r,\xi,c))
									K(r,s,\xi,c)\tilde{\mathrm{Ai}}(x(s,\xi,c))\\
									&=(Wq^{\f12})(s,c)
									\cdot\f23\xi^{-\f23}\textrm{Ai}^{2}\left(x^{\f23}(s)\right)\int
									^{x(s)}
									_{x(r)}
									\f{\textrm{Ai}^{-2}(y^{\f23})}{y^{\f13}}\cdot \left(\f{1+\f23\xi^{-\f23}\widetilde{C}_{61}
										\int^{x_*(\xi,c)}_{x(s)}y'^{-\f13}\mathrm{Ai}^{-2}(y'^{\f23})dy'}{
										1+\f23\xi^{-\f23}\widetilde{C}_{61}\int^{x_*(\xi,c)}_{y}y'^{-\f13}\mathrm{Ai}^{-2}
										(y'^{\f23})dy'}\right)^{2} dy
									\end{align*}
								with $W$ is defined in \eqref{W:fm0}, $x=\xi(-\tau)^{\f32}(r,c)$. Thanks to the definition and \eqref{equal:1+B2D2}, we rewrite $K_1$ as
								\begin{align}\label{def:Phi1-K1}
								\begin{split}
									K_1(r,s,\xi,c)
									& =(Wq^{\f12})(s,c)
									\cdot\f23\xi^{-\f23}\textrm{Ai}^{2}\left(x^{\f23}(s,\xi,c)\right)\\&\quad\times\int
									^{x(s,\xi,c)}
									_{x(r,\xi,c)}
									\f{\textrm{Ai}^{-2}(y^{\f23})}{y^{\f13}}\left(\f{1-B_1(\xi,c)D_1(x_*-x(s,\xi,c),x_*)}{
										1-B_1(\xi,c)D_1(x_*-y,x_*)}\right)^{2} dy.
										\end{split}
									\end{align}
								Here, $x_*=x(r_*(\xi,c),\xi,c)$, $B_1,D_1$ are defined in \eqref{def:B2D2}. It suffices to verify the conditions in Lemma \ref{lem:vorteraa-0}
								for $M^{\f12}r_c^{-\f12}\langle r_c\rangle^{-1}\sim r_*\leq r\leq r_c-C\xi^{-\f23}\langle r_c\rangle^{\f12}$. The condition \eqref{condition:K'}(with $\kappa=M^{-\f12}$) writes
								\begin{align*} \notag&\int_{r_*}^{r}
									\sup_{\tilde{r}\in[s,r]}|K_1(\tilde{r},s,\xi,c)|ds\lesssim M^{-\f12},\notag
								\end{align*}
								and the condition \eqref{condition:K'pac-growth}  writes
								\begin{align}
									\begin{split}\label{bd:K1-Phi12}&\int_{0}^{1}
										\sup_{\tilde{t}\in[t,1]}\left|\left(\f{r_c-r}{\langle r_c\rangle^4}
										\pa_c\right)^l
										\left(\f{(r_c-r)r}{ r_c}\pa_{r}\right)^i (\xi\pa_{\xi})^j \left((r-r_*)
										K_1\left(\bar{r}(\tilde{t},r,\xi,c),\bar{r}(t,r,\xi,c),\xi,c\right)\right)\right|dt\lesssim 1,\\ &\quad\quad\quad\quad\quad\text{with}\;\;\bar{r}(t,r,\xi,c)=tr+(1-t)r_*(\xi,c),\;\;\text{for}\;r\geq 2r_*.
								\end{split}
							\end{align}
							The first bound follows from \eqref{est:KAi}, \eqref{lem:behave-W} and $Q(s,c)\sim \f{(s-r_c)\langle r_c\rangle^2}{\langle s\rangle ^3}(s\lesssim r_c)$. Indeed, 
							\begin{align}\label{bd:K1-Phi11-pf}
								\notag&\int_{r_*}^{r}
								\sup_{\tilde{r}\in[s,r]}|K_1(\tilde{r},s,\xi,c)|ds
								\lesssim
								\int_{r_*}^{r}\sup_{\tilde{r}\in[s,r]}|
								\xi^{-\f23}\widetilde{K}_{Ai}(x(\tilde{r}),y(s))|\cdot
								|W(s,c)|\cdot|q^{\f12}(s,c)|ds\\
								&\lesssim
								\int_{r_*}^{r}\xi^{-1}|Q(s,c)|^{-\f12}s^{-2}ds\sim  \xi^{-1}\int_{r_*}^{r}\langle r_c\rangle^{-1}|r_c-s|^{-\f12}\langle s\rangle^{\f32}s^{-2}ds\\
								\notag &
								\lesssim \xi^{-1}\langle r_c\rangle^{-1}r_c^{-\f12}
								\int_{r_*}^{\min\{r,r_c/2\}}\langle s\rangle^{\f32}s^{-2}ds
								+\xi^{-1}\langle r_c\rangle^{-1}\langle r_c\rangle^{\f32}r_c^{-2}
								\int^{r}_{\min\{r,r_c/2\}}|r_c-s|^{-\f12}ds\notag \\
								&\lesssim
								\xi^{-1}r_c^{-\f12}\langle r_c\rangle^{-1} r_*^{-1}
								+\left(\xi r_c^{\f32}/\langle r_c\rangle^{\f12}\right)^{-1}\lesssim M^{-\f12}.\notag
							\end{align}
							Thanks to Leibniz law, $\left|\left(\f{r_c-r}{\langle r_c\rangle^4}
							\pa_c\right)^l
							\left(\f{(r_c-r)r}{ r_c}\pa_{r}\right)^i (\xi\pa_{\xi})^j \left((r-r_*(\xi,c))\right)\right|\lesssim r\leq 2(r-r_*)$ for $r\geq 2r_*$, and the following estimate by \eqref{bd:K1-Phi11-pf}, \begin{align*}
								\int_0^1(r-r_*) \cdot
								\xi^{-1}|Q(\bar{r}(t),c)|^{-\f12}\bar{r}(t)^{-2}dt
								=\int_{r_*}^{r}\xi^{-1}|Q(s,c)|^{-\f12}s^{-2}ds\lesssim 1,\; \text{with}\;
								\bar{r}(t)=r_*+t(r-r_*),\end{align*}
								for the bound \eqref{bd:K1-Phi12}, it suffices to prove that
								\begin{align}
									\label{bd:K1-Phi13}&\left|\left(\f{r_c-r}{\langle r_c\rangle^4}
									\pa_c\right)^l
									\left(\f{(r_c-r)r}{ r_c}\pa_{r}\right)^i (\xi\pa_{\xi})^j \left(
									K_1(\bar{r}(\tilde{t},r,\xi,c),\bar{r}(t,r,\xi,c),\xi,c)\right)\right|
									\lesssim \f{\xi^{-\f23}x^{-\f13}(s)}{s^2q^{\f12}(s,c)}\Big|_{s=\bar{r}(t)}\\
									&\sim \xi^{-1}|Q(\bar{r}(t),c)|^{-\f12}\bar{r}(t)^{-2}. \notag
								\end{align}
								Using the behavior of $\mathrm{Ai}$ in \eqref{behave:Ai}, \eqref{behave:Ai^-1} and letting $v=\f{x(s)-y}{x(s)-x(r)}=\f{x(\bar{r}(t))-y}{x(\bar{r}(t))-x(\bar{r}
									(\tilde{t}))}$ in    \eqref{def:Phi1-K1}, the desired form of $K_1$ behaves as follows
									\begin{align}\label{def:Phi1-K1 alt}
									\notag & K_1(\bar{r}(\tilde{t},r,\xi,c),\bar{r}(t,r,\xi,c),\xi,c) \\
									&=C(W q^{\f12})(\bar{r}(t),c)\cdot\xi^{-\f23}x^{-\f13}
									(\bar{r}(t))\left(1+O(x^{-1} (\bar{r}(t)))\right)\left(
									x(\bar{r}(t)- x(\bar{r}(\tilde{t}))\right)\cdot\\
									\notag&\qquad\int
									^{1}   _{0}    e^{-\f43v\left(
										x(\bar{r}(t))- x(\bar{r}(\tilde{t}))\right)}
									\left(1+O\left(\left(x(\bar{r}(t))+v(x(\bar{r}(t))-x(\bar{r}(\tilde{t})))\right)^{-1}\right)
									\right)\cdot\\
									\notag &\qquad\left(\f{1-B_1(\xi,c)D_1(x_*-x(\bar{r}(t,r,\xi,c),\xi,c),x_*)}{
										1-B_1(\xi,c)D_1(x_*-x(\bar{r}(t,r,\xi,c),\xi,c)+v\left(
										x(\bar{r}(t,r,\xi,c))-x(\bar{r}(\tilde{t},r,\xi,c))\right),x_*)}\right)^{2}dw. \end{align}
									Notice  that for $0\leq t\leq \tilde{t}\leq 1$,
									\begin{align*}
										&r_*=\bar{r}(0)\leq \bar{r}(t)\leq  \bar{r}(\tilde{t})\leq \bar{r}(1)=r,\\
										&x_*=x(\bar{r}(0))\geq x(\bar{r}(t))\leq  x(\bar{r}(\tilde{t}))\geq x(\bar{r}(1))=x(r)\gtrsim 1.
										\end{align*}
										We apply Leibniz law to \eqref{def:Phi1-K1}, then \eqref{bd:K1-Phi13} can be reduced to the following bounds for
										$r\leq r_c-C\xi^{-\f23}\langle r_c\rangle^{\f13}$, \;$\bar{r}(t)=\bar{r}(t,\xi,r,c):=(1-t)r_*(\xi,c)+tr$, $t\in [0,1]$:
										 \begin{align}
											\label{est:r*-1}& \left|\left(\f{r_c-r}{\langle r_c\rangle^4}\right)^l
											\left(\f{(r_c-r)r}{ r_c}\right)^i\pa_c^l\pa_r^i
											(\xi\pa_{\xi})^j\left( x(\bar{r}(t,\xi,r,c),\xi,c )
											\right)\right| \lesssim x(\bar{r}(t)),
										\end{align}
										\begin{align}
											\label{est:r*-2} & \left|\left(\f{r_c-r}{\langle r_c\rangle^4}\right)^l
											\left(\f{(r_c-r)r}{ r_c}\right)^i\pa_c^l\pa_r^i
											(\xi\pa_{\xi})^j \left((Wq^{\f12})(\bar{r}(t,r,\xi,c),c) \right)\right| \lesssim \f{1}{s^2q^{\f12}(s,c)|_{s=r(t)}},
										\end{align}
										and for every $0\leq t\leq \tilde{t}\leq 1$,
										\begin{align}
											&\left|\left(\f{r_c-r}{\langle r_c\rangle^4}\right)^l
											\left(\f{(r_c-r)r}{ r_c}\right)^i\pa_c^l\pa_r^i
											(\xi\pa_{\xi})^j \left( x(\bar{r}(t,r,\xi,c))-x(\bar{r}(\tilde{t},r,\xi,c))
											\right)\right| \lesssim x(\bar{r}(t))-x(\bar{r}(\tilde{t})).   \label{est:r*-3}
										\end{align}
										Letting $t=0$, $\tilde{t}=t$ in \eqref{est:r*-3} gives \begin{align*}
											&\left|\left(\f{r_c-r}{\langle r_c\rangle^4}\right)^l
											\left(\f{(r_c-r)r}{ r_c}\right)^i\pa_c^l\pa_r^i  (\xi\pa_{\xi})^j
											\left( x_*-x(\bar{r}(t) )
											\right)\right| \lesssim x_*-x(\bar{r}(t)).
										\end{align*}
										If the derivative $\left(\f{r_c-r}{\langle r_c\rangle^4}\right)^{l'} \left(\f{(r_c-r)r}{ r_c}\right)^{i'}\pa_c^{l'}(\xi\pa_{\xi})^{j'}\pa_r^{i'} $
										acts on the most ``complicated" term $\left(\f{1-B_1D_1}{1-B_1D_1}\right)^2$, we use the bounds of $B_1$, $D_1$
										obtained in \eqref{est:B2} and \eqref{est:positivetD2} and the following observations:
										\begin{align*}
											&(\xi\pa_{\xi})\left(D_1(x_*-x(\bar{r}(t)),\xi,c),x_*\right)
											=  (\xi\pa_{\xi})x_*\pa_{x_*}D_1(z,x_*)+ (\xi\pa_{\xi})(x_*-x(\bar{r}(t)))\pa_zD_1(z,x_*)
											|_{z=x_*-x(\bar{r}(t))},\\
											&\pa_{r}\left(D_1(x_*-x(\bar{r}(t)),\xi,c),x_*\right)
											= \xi\pa_{r}(x_*-x(\bar{r}(t)))\pa_zD_1(z,x_*)
											|_{z=x_*-x(\bar{r}(t))},\\
											&\pa_{c}\left(D_1(x_*-x(\bar{r}(t),\xi,c),x_*)\right)
											=\pa_{c}x_*\pa_{x_*}D_1(z,x_*)+ \pa_{c}(x_*-x(\bar{r}(t)))\pa_zD_1(z,x_*)
											|_{z=x_*-x(\bar{r}(t))}. \end{align*}
											
											It remains to prove \eqref{est:r*-1}-\eqref{est:r*-3}. We denote
											\begin{align*}
												\tilde{x}(\xi,r,c)=\tilde{x}(\xi,r,c)(t):=x(\bar{r}(t,\xi,r,c),\xi,c ), \end{align*}
											and notice that
											\begin{align*}
												&(\xi\pa_{\xi} ) (\tilde{x})=x+\xi\pa_{\xi}\bar{r}\pa_{\bar{r}}x=x-(1-t)r_*
												\pa_{\bar{r}}x ,\;  \pa_r(\tilde{x})=t\pa_{\bar{ r}}x,\;\pa_c(\tilde{x})=  \pa_c\bar{ r}\pa_{\bar{ r}}x+\pa_cx.
											\end{align*}
											To prove \eqref{est:r*-1}, using Leibniz's law, the Fa$\grave{\text{a}}$di Bruno's Formula, the facts above and $\left|\f{r_c-r}{\langle r_c\rangle^4}\pa_c\bar{r}\right|\lesssim (1-t)r_*$, we deduce similarly as  \eqref{bd:t-hn} that
											\begin{align*}
												&\notag \left|\left(\f{r_c-r}{\langle r_c\rangle^4}\right)^l
												\left(\f{(r_c-r)r}{ r_c}\right)^i\pa_c^l\pa_r^i
												(\xi\pa_{\xi})^j\left( x(\bar{r}(t,\xi,r,c),\xi,c )
												\right)\right|
												\\
												&\lesssim \sum_{0\leq j_2\leq j}\left|\left(\f{(r_c-r)rt}{ r_c}\right)^i\left((1-t)r_*\right)^{j_2+l}\pa_{\bar{r}}
												^{i+j_2+l}x(\bar{r}(t),\xi,c)\right|\\
												&\qquad+\left|\left(\f{r_c-r}{\langle r_c\rangle^4}\right)^l\left(\f{(r_c-r)rt}{ r_c}\right)^i\left((1-t)r_*\right)^{j_2}\pa_c^l\pa_{\bar{r}}
												^{i+j_2}x(\bar{r}(t),\xi,c)\right|.\notag
											\end{align*}
											Recalling that \eqref{estx-G} gives us  $\left|\left(\f{(r_c-\bar{r})}{\langle r_c\rangle^{4}}\right)^l
											\left(\f{(r_c-\bar{r}) r}{ r_c}\right)^i\pa_c^l\pa_{\bar{r}}^ix(\bar{r})\right|\lesssim |x(\bar{r})|$,  the above can be bounded by
											\begin{align*}&\lesssim \sum_{0\leq j_2\leq j}
												\left(\f{(r_c-r)rt}{r_c}/\f{(r_c-\bar{r})\bar{r}}{r_c}\right)^i
												\left((1-t)r_*/\f{
													(r_c-\bar{r})\bar{r}}{r_c}\right)^{j_2+l}\\
												&\qquad+ \left(\f{r_c-r}{\langle r_c\rangle^4}/\f{r_c-\bar{r}}{\langle r_c\rangle^4}\right)^l
												\left(\f{(r_c-r)rt}{r_c}/\f{(r_c-\bar{r})\bar{r}}{r_c}\right)^i
												\left((1-t)r_*/\f{(r_c-\bar{r})\bar{r}}{r_c}\right)^{j_2}\sim 1,
												\end{align*}
											where we used the fact that
											\begin{align}
												\begin{split}\label{case6-useful-bd-r}
													&r t\leq \bar{r}(t)\leq r\leq r_c,\;\quad \bar{r}(0)=r_*\leq C\xi^{-\f23}\langle r_c\rangle^{\f13}\leq r_c-r=r-\bar{r}(1),\\
													&\text{and}\; \quad\f{r_*r_c}{(r_c-\bar{r}(t))\bar{r}(t)}\leq \f{r_*r_c}{(r_c-\bar{r}(0))\bar{r}(0)}\lesssim 1.
												\end{split}
											\end{align}
											For \eqref{est:r*-2},  a similar process shows
											\begin{align*}
												&\notag \left|\left(\f{r_c-r}{\langle r_c\rangle^4}\right)^l
												\left(\f{(r_c-r)r}{ r_c}\right)^i\pa_c^l\pa_r^i
												(\xi\pa_{\xi})^j\left( (Wq^{\f12})(\bar{r}(t,\xi,r,c),c )\right)\right|
												\\
												&\lesssim \sum_{0\leq j_2\leq j}\left|\left(\f{(r_c-r)rt}{ r_c}\right)^i\left((1-t)r_*\right)^{j_2+l}\pa_{\bar{r}}
												^{i+j_2+l}(Wq^{\f12})(\bar{r}(t),c )\right|\\
												&\qquad+\left|\left(\f{r_c-r}{\langle r_c\rangle^4}\right)^l\left(\f{(r_c-r)rt}{ r_c}\right)^i\left((1-t)r_*\right)^{j_2}\pa_c^l\pa_{\bar{r}}
												^{i+j_2}(Wq^{\f12})(\bar{r}(t),c )\right|,\notag
											\end{align*}
											and we apply \eqref{bd:q-Derivative} and \eqref{est:W<rc} to obtain that
											$\left|(\f{1}{\langle r_c\rangle ^{3}} )^l \bar{r}^{i}\pa_c^l\pa_{\bar{r}}^i(W q^{\f12})(\bar{r},c)\right|\lesssim
											\f{1}{\bar{r}^2
												q(\bar{r},c)},\;\bar{r}\leq r_c. $
												 Then  for \eqref{est:r*-3}, it remains to check that
												\begin{align*}\left|\f{(r_c-r)rt}{r_c}/\bar{r}(t)\right|\lesssim 1,\;
													\left|(1-t)r_*/\bar{r}(t)\right|\lesssim 1,\; \left|\f{r_c-r}{\langle r_c\rangle^4}/\f{1}{\langle r_c\rangle^3}\right|\lesssim 1,
													\end{align*}
													which are obvious due to $\max\{r t,(1-t)r_*\}\leq \bar{r}(t)\leq r\leq r_c$.
													For \eqref{est:r*-3}, we apply Leibniz's law and  Fa$\grave{\text{a}}$di Bruno's Formula to deduce that
													\begin{align*}
													&\notag \left|\left(\f{r_c-r}{\langle r_c\rangle^4}\right)^l
													\left(\f{(r_c-r)r}{ r_c}\right)^i\pa_c^l\pa_r^i
													(\xi\pa_{\xi})^j\left( x(\bar{r}(t))- x(\bar{r}(\tilde{t}))
													\right)\right|
													\\
													&\lesssim \small\sum_{0\leq j_2\leq j}\left|\left(\f{(r_c-r)r}{ r_c}\right)^ir_*^{j_2+l}\left(t^i(1-t)^{j_2+l}\pa_{\bar{r}}
													^{i+j_2+l}x(\bar{r}(t))-\tilde{t}^i(1-\tilde{t})^{j_2+l}\pa_{\bar{r}}
													^{i+j_2+l}x(\bar{r}(\tilde{t}))\right)\right|\\
													&\qquad+\left|\left(\f{r_c-r}{\langle r_c\rangle^4}\right)^l\left(\f{(r_c-r)r}{ r_c}\right)^ir_*^{j_2}\left(t^i(1-t)^{j_2+l}\pa_c^l\pa_{\bar{r}}
													^{i+j_2}x(\bar{r}(t))-\tilde{t}^i(1-\tilde{t})^{j_2+l}
													\pa_c^l\pa_{\bar{r}}
													^{i+j_2}x(\bar{r}(\tilde{t}))\right)\right|\notag \\
													&:=S_1(\tilde{t})+S_2(\tilde{t}).\notag
												\end{align*}
												
												Fix $t\in[0,1]$. To prove \eqref{est:r*-3}, it remains to check for $ \tilde{t}\in[t,1]$ that $S_k(\tilde{t})\leq C \left(x(\bar{r}(t))-x(\bar{r}(\tilde{t}))\right)$, $k=1,2$. Noticing $S_k(t)= \left(x(\bar{r}(t))-x(\bar{r}(\tilde{t}))\right)|_{\tilde{t}=t}=0$, it suffices to prove the following $\tilde{t}-$derivative bounds for any $i+j+l\geq 1$ and $r\geq 2r_*$:
												\begin{align*}
													\left|\left(\f{r_c-r}{\langle r_c\rangle^4}\right)^l\left(\f{(r_c-r)r}{ r_c}\right)^ir_*^{j}\f{d}{d\tilde{t}}\left(\tilde{t}^i(1-\tilde{t})^{j}\pa_c^l\pa_{\bar{r}}
													^{i+j}x(\bar{r}(\tilde{t}))\right)\right|
													\leq C\left|\f{d}{d\tilde{t}}(x(\bar{r}(\tilde{t})))\right|.
												\end{align*}
												Notice that $\f{d}{d\tilde{t}}(x(\bar{r}(\tilde{t})))
												=\f32\xi(r-r_*)(-Q)^{\f12}(\bar{r}(\tilde{t}))$, $\f{d\bar{r}(\tilde{t})}{d\tilde{t}}=r-r_*$.
												It follows from \eqref{estQ-G} that for $\bar{r}(\tilde{t})\leq r\leq r_c$,
												\begin{align}\label{useful-Q1/2}\left|\left(\f{(r_c-\bar{r}(\tilde{t}))}{\langle r_c\rangle^{4}}\right)^l
													\left(\f{(r_c-\bar{r}(t)) \bar{r}(\tilde{t})}{r_c}\right)^i\pa_c^l\pa_{\bar{r}}^i((-Q)^{\f12})(\bar{r}(t),c)\right|\lesssim |(-Q)^{\f12}(\bar{r}(t),c)|.
												\end{align}
												If $i+j=0$, then $l=1$, and the bound above writes as
												\begin{align*}
													\left|\left(\f{r_c-r}{\langle r_c\rangle^4}\right)(r-r_*)\pa_c
													((-Q)^{\f12})(\bar{r}(\tilde{t}),c)\right|
													\leq C(r-r_*)(-Q)^{\f12}(\bar{r}(\tilde{t}),c),
													\end{align*}
												which is a consequence of \eqref{useful-Q1/2} and $\f{r_c-r}{\langle r_c\rangle^4}/\f{r_c-\bar{r}(\tilde{t})}{\langle r_c\rangle^4}\leq 1$.
												If $i+j\geq 1$,  we need to prove
												\begin{align*}
													\left|\left(\f{r_c-r}{\langle r_c\rangle^4}\right)^l\left(\f{(r_c-r)r}{ r_c}\right)^ir_*^{j}\f{d}{d\tilde{t}}\left(\tilde{t}^i
													(1-\tilde{t})^{j}\pa_c^l\pa_{\bar{r}}
													^{i+j-1}((-Q)^{\f12})(\bar{r}(\tilde{t}),c)\right)\right|
													\leq C(r-r_*)(-Q)^{\f12}(\bar{r}(\tilde{t}),c).
													\end{align*}
												Using \eqref{useful-Q1/2}, $\f{d\bar{r}(\tilde{t})}{d\tilde{t}}=r-r_*$ and \eqref{case6-useful-bd-r}, we only need to check $\f{(r_c-r)r}{r_c}+r_*+(r-r_*)\leq r-r_*$. This is due to $2r_*\leq r\leq r_c$. This finishes the proof of the lemma.
							\end{proof}
							
							\subsection{Case of $c\in(1-\delta,1)$ and $M^2(1-c)^{\f12}\lesssim \xi \lesssim M(1-c)^{\f13}$}
							\begin{lemma}
								\label{lem:Phi81}
								Let $\delta\ll 1$, $M\gg 1$ be fixed. Let  $c\in(1-\delta,1)$ and $M^2(1-c)^{\f12}\lesssim  \xi\lesssim M (1-c)^{\f13}$. Let
								$\phi$ be the solution to \eqref{eq:phi}, i.e.,
								$\phi''+r^{-1}\phi'-r^{-2}\phi+\f{\xi^2(V-c)}{1-c}  \phi=0,\; \phi(r,\xi,c)\sim \xi r$,
								which has been solved on $(0, M\f{\sqrt{1-c}}{\xi})$ in Lemma \ref{lem:phi-sqrt(1-c)/xi}. Let $r_{*}:=M\f{\sqrt{1-c}}{\xi}$, $r_{*2}:=M^{-3.5}\f{\xi^2}{1-c}$ satisfying  $r_{*}\ll1 \ll r_{*2}$. Then there exists a constant
								\begin{align}\label{def:C5}
									\begin{split}
										&C_{71}(\xi,c):=\left(1-V(r)\right)^{\f14}r^{\f12}\phi(r)\Big|_{r=r_{*}},
									\end{split}
								\end{align}
								such that for $M\f{\sqrt{1-c}}{\xi}\lesssim r\lesssim M^{-3.5}\f{\xi^2}{1-c}$, there holds that
								\begin{align}\label{Phi8-1}
									\begin{split}
										&\phi(r,\xi,c)\sim C_{71}(\xi,c)\left(1-V(r)\right)^{-\f14}r^{-\f12}\cosh\f{\xi (v(r_{*})-v(r)))}{\sqrt{1-c}},\\
										&\phi(r,\xi,c)=C_{71}(\xi,c)\left(1-V(r)\right)^{-\f14}r^{-\f12}\cosh\f{\xi (v(r_{*})-v(r))}{\sqrt{1-c}}O_{r,\xi,c}^{r,\xi,1-c}(1) \quad \text{for}\;\;r\geq 2r_*,
								\end{split} \end{align}
								where   $v(r)=\int_r^{+\infty}(1-V(s))^{\f12}ds$, which satisfies $v(r)\sim \langle r\rangle ^{-\f12}$ and   \eqref{bd:v_*-v}. More precisely, the following bound holds for the constant \begin{align}\label{bd:C7}
									\begin{split}
										&\f{1}{C_{71}(\xi,c)}
										=\f{1}{\sqrt{1-c}}O_{\xi,c}^{\xi,1-c}\left(\left(\f{\xi}{(1-c)^{\f12}}\right)^{\f12}\right).
									\end{split}
								\end{align}
							\end{lemma}
							
							\begin{proof}
								This case is different from other cases due to the following observation. For $c\in(1-\delta,1)$, $M^2(1-c)^{\f12}\lesssim \xi \lesssim M(1-c)^{\f13}$, and $M\f{\sqrt{1-c}}{\xi}(:=r_{*})\lesssim r\lesssim M^{-3.5}\f{\xi^2}{1-c}(:=r_{*2})$, where $r_{*}\ll 1$, $1\ll r_{*2}\ll \f{\xi^2}{1-c}$, we observe that the term $-\f{\xi^2(1-V)}{1-c}$ is the main term in the equation due to
								\begin{align*}
									\f{\xi^2(1-V(r))}{1-c}  \gg r^{-2}\gtrsim  \xi^2.
								\end{align*}
								Here, in the first inequality we used                 								\begin{align*}
									\f{\xi^2(1-V(r))}{1-c}\sim
									\left\{
									\begin{aligned}
										& \f{\xi^2}{1-c}\gg r^{-2},\quad M \f{\sqrt{1-c}}{\xi}\lesssim r\lesssim 1,\\
										& \f{\xi^2r^{-3}}{1-c}\gg r^{-2},\quad 1\lesssim r\lesssim M^{-3.5}\f{\xi^2}{1-c}.
									\end{aligned}
									\right.
								\end{align*}
								Therefore, we rewrite the equation as
								\begin{align} \label{rewrite:phi7}
									\left(r^{\f12}\phi\right)''
									-\f{\xi^2(1-V)}{1-c}\left(r^{\f12}\phi\right)
									=\left(\f{3}{4}r^{-2}-\xi^2\right)\left(r^{\f12}\phi\right).
								\end{align}
								We take the new coordinate $r\to v$ such that $v(r)=\int_r^{+\infty}(1-V(s))^{\f12}ds$, which satisfies $v(r)\sim \langle r\rangle^{-\f12}$  by $1-V(r)\sim \langle r\rangle^{-3}$.
								Then we plug  $\Phi(v,\xi,c)=\left(1-V(r)\right)^{\f14}r^{\f12}\phi(r,\xi,c)$ into the above equation to find that
								\begin{align}\label{eq:Phi7}
									&\pa_v^2\Phi-\f{\xi^2}{1-c}\Phi=W(r,\xi)\Phi,
								\end{align}
								where
								\begin{align*}
									W(r,\xi)=\left(\f{3}{4}r^{-2}-\xi^2\right)(1-V(r))^{-1}
									-(1-V(r))^{-\f34}\left((1-V(r))^{-\f14}\right)''
									\end{align*}
								admits the bound
								\begin{align}\label{bd:WPhi7}
									\left|(r\pa_{r})^i(\xi\pa_{\xi})^jW(r,\xi)\right|\lesssim \f{\langle r\rangle^3}{r^2}\quad\text{for}\quad M \f{(1-c)^{\f12}}{\xi}\lesssim r\lesssim M^{-3.5}\f{\xi^2}{1-c}.
									\end{align}
									To smoothly  extend the solution solved in $(0, r_*)$ in Lemma \ref{lem: phi-2}, we match the initial data as follows
									\begin{small}\begin{align}\label{eq:Phi8initial}
											\Phi(v_*)=\left(1-V(r)\right)^{\f14}r^{\f12}\phi(r)\big|_{r=r_*},\;\;
											(\pa_v\Phi)(v_*)=-\left(1-V(r)\right)^{\f12}
											\pa_r\left(\left(1-V(r)\right)^{\f14}r^{\f12}\phi(r)\right)\big|_{r=r_*},  \end{align}
									\end{small}
									where we denote $v_*=v(r_*)$.
									Taking $W\equiv 0$ in \eqref{eq:Phi7}, the solution of the homogeneous equation with initial data \eqref{eq:Phi8initial} writes
									\begin{align} \label{def:Phi0}
										\Phi_0(v,\xi,c)=C_{71}(\xi,c)\cosh \f{\xi(v_*-v)}{\sqrt{1-c}} \left(1-\f{\f{\sqrt{1-c}}{\xi}(\pa_v\Phi)(v_*)}{\Phi(v_*)}\tanh \f{\xi(v_*-v)}{\sqrt{1-c}} \right)
										\end{align}
										with $C_{71}(\xi,c)= \Phi(v_*)$.
										Then a direct calculation gives the equivalent integral form of \eqref{eq:Phi7} as
										\begin{align}\label{eq:Phi5}
											\Phi(v,\xi,c)&=\Phi_0(v,\xi,c)+
											\int_{v_*(\xi,c)}^{v}K(v,w,\xi,c)\Phi(w,\xi,c)dw,
										\end{align}
										where
										$$K(v,w,\xi,c)=\Phi_0(v,\xi,c)\Phi_0(w,\xi,c)
										\int_{w}^{v}
										\Phi_0(\tilde{v},\xi,c)^{-2}d\tilde{v}\cdot W(s(w),\xi).$$
										
										We claim that for $r_*\leq r\leq M^{-3.5}\f{\xi^2}{1-c}$, $l\in\{0,1\}$, $i,j\in\mathbb{N}$, it holds uniformly in $r,\xi,c$ that
										\begin{align} \label{est:Phi7derivative}
											\begin{split}
												&\left|\f{\f{\sqrt{1-c}}{\xi}(\pa_v\Phi)(v_*)}{\Phi(v_*)}\tanh \f{\xi(v-v_*)}{\sqrt{1-c}}\right|\lesssim M^{-1},\\
												&\left|\left((1-c)\pa_c\right)^l(\xi\pa_{\xi})^j(r\pa_r)^i
												\left(\f{\f{\sqrt{1-c}}{\xi}(\pa_v\Phi)(v_*)}{\Phi(v_*)} \right)\right|\lesssim 1,\\
												&\left|\left((1-c)\pa_c\right)^l(\xi\pa_{\xi})^j(r\pa_r)^i
												\left(\tanh \f{\xi(v_*-v(r))}{\sqrt{1-c}} \right)\right|\lesssim 1,
											\end{split}
										\end{align}
										where   $v_*(\xi,c)=v(r_*(\xi,c))$ and
										\begin{align*}
											&(\pa_{v}\Phi)(v_*,\xi,c)=\left(1-V(r_*(\xi,c))\right)^{-\f12}\left(
											\left(1-V(r_*(\xi,c))\right)^{\f14}r_*(\xi,c)^{\f12}\phi(r_*(\xi,c),\xi,c)
											\right)',\\
											&\Phi(v_*,\xi,c)=\left(1-V(r_*(\xi,c))\right)^{\f14}r_*(\xi,c)^{\f12}\phi(r_*(\xi,c),\xi,c),
											\end{align*}
										and the constant for the last bound may depends on $M$; and
										\begin{align}
											\begin{split}\label{claim:Phi7} &\Phi(r,\xi,c)=\Phi_0(v(r),\xi,c)\left(1+\phi^{Rem}(r,\xi,c)\right),\\
												&\quad\text{with}
												\;\;\phi^{Rem}(r_*)=\pa_r\phi^{Rem}(r_*)=0,\;\;|\phi^{Rem}|\lesssim M^{-1}; \end{split}
										\end{align}
										and
										\begin{align}\label{est:phi7-Derivative} &\left|\left((1-c)\pa_c\right)^l\left(r\pa_r\right)^i(\xi\pa_{\xi})^j \phi^{Rem}(r,\xi,c)\right|\lesssim 1\;\;\text{for}\;\;r\geq 2r_*.
										\end{align}
										Noticing  $M\gg 1$, \eqref{def:Phi0}, \eqref{est:Phi7derivative}, \eqref{claim:Phi7} and \eqref{est:phi7-Derivative} give the result \eqref{Phi8-1}, where $C_{71}=\Phi(v_*)$.   It remains to prove
										\eqref{est:Phi7derivative}, \eqref{claim:Phi7} and \eqref{est:phi7-Derivative}. The first bound in \eqref{est:Phi7derivative} follows  from $|\tanh z|\lesssim 1$ and
	\begin{align*}
		\left|\f{\f{\sqrt{1-c}}{\xi}(\pa_v\Phi)(v_*)}{\Phi(v_*)}\right|&=
		\left|\f{\f{\sqrt{1-c}}{\xi}\left(1-V(r)\right)^{-\f12}\pa_r
			\left(\left(1-V(r)\right)^{\f14}r^{\f12}\phi(r)\right)}
		{\left(1-V(r)\right)^{\f14}r^{\f12}\phi(r)}\Big|_{r=r_*}\right|\\
		&=
		M^{-1}\left(1-V(r_*)\right)^{-\f12}\left|\f{r\pa_r
			\left(\left(1-V(r)\right)^{\f14}r^{\f12}\phi(r)\right)}
		{\left(1-V(r)\right)^{\f14}r^{\f12}\phi(r)}\Big|_{r=r_*}\right|\\
		&\lesssim M^{-1}\left(\f{r|V'(r)|}{1-V(r)}\Big|_{r=r_*}+1+\f{r\phi'(r)}{\phi(r)}\Big|_{r=r_*}\right)\lesssim M^{-1},
		\end{align*}
	where we used $M\f{\sqrt{1-c}}{\xi}=r_*\ll 1\ll r_c$ and $\phi(r_*)\geq \xi r_*$, $|r_*\phi'(r_*)|\lesssim \xi r_*$ (by Lemma \ref{lem: phi-2}).  The second bound in \eqref{est:Phi7derivative} follows from $\left|r^k\f{d^k}{dr^k}(1-V(r)) \right|\lesssim (1-V(r)),\; \left|r^k\f{d^k}{dr^k}(v(r)) \right|\lesssim v(r),\;\;k\geq 1$ and
	\begin{align}\label{est:r*7}
		\begin{split}
			&\left|\left((1-c)\pa_c\right)^l(\xi\pa_{\xi})^jr_*(\xi,c)\right|\lesssim r_*\ll 1\ll r_c,\;\;\left|\left((1-c)\pa_c\right)^l(\xi\pa_{\xi})^jv_*(\xi,c)\right|\lesssim r_*,\quad l+j\geq 1,
		\end{split}
	\end{align}
	and by using  Lemma \ref{lem:phi-sqrt(1-c)/xi},
	\begin{align}
		&\left|\left((1-c)\pa_c\right)^l\left(r\pa_r\right)^i(\xi\pa_{\xi})^j \phi (r,\xi,c)\right|\Big|_{r=r_*}\lesssim \xi r_*.
	\end{align}
	The last bound in \eqref{est:Phi7derivative} follows from $\left|\f{d^k}{dz^k}(\tanh z) \right| \lesssim e^{-2z}(z\geq 0,\;k\geq1)$, 
	and the observation that
	\begin{align*}
		\left|\f{\xi r^kv^{(k)}(r)}{\sqrt{1-c}}\right|\lesssim \f{\xi r(1-V(r))^{\f12}}{\sqrt{1-c}}\sim
		\left\{
		\begin{aligned}
			& \f{\xi r}{\sqrt{1-c}},\quad r_*\lesssim r\lesssim 1,\\
			& \f{\xi r^{-\f12}}{\sqrt{1-c}},\quad 1\ll r\lesssim M^{-3.5}\f{\xi^2}{1-c},
		\end{aligned}
		\right.
	\end{align*}
	and
	\begin{align}\label{bd:v_*-v}
		\f{\xi(v_*-v(r))}{\sqrt{1-c}}\sim
		\left\{
		\begin{aligned}
			& \f{\xi (r-r_*)}{\sqrt{1-c}}\sim \f{\xi r}{\sqrt{1-c}}-M,\quad r_*\lesssim r\lesssim 1,\\
			& \f{\xi }{\sqrt{1-c}}\gtrsim \f{\xi r^{-\f12}}{\sqrt{1-c}}\gg M^2\gg 1 ,\quad 1\ll r\lesssim M^{-3.5}\f{\xi^2}{1-c}. \end{aligned}
		\right.
	\end{align}
	
	To prove \eqref{claim:Phi7},  we apply Lemma \ref{lem:vorteraa-0} to the integral equation \eqref{eq:Phi5}. It suffices to verify the condition \eqref{condition:K'} with  $\kappa=M^{-1}$ and $K_1$ as follows
	\begin{align}\label{def:Phi5-K1}
		\notag & K_1(r,s,\xi,c):=\Phi_0^{-1}(v(r))
		K(r,s)\Phi_0(v(s))\f{dv}{ds}\\
		\notag &=\Phi_0(v(s),\xi,c)^2
		\int_{v(s)}^{v(r)}
		\Phi_0(\tilde{v},\xi,c)^{-2}
		d\tilde{v}\cdot W(s,\xi)(1-V(s))^{\f12}\\
		&=
		\int_{v(s)}^{v(r)}\f{\cosh^2 \f{\xi(v(s)-v_*)}{\sqrt{1-c}} }{\cosh^2 \f{\xi(\tilde{v}-v_*)}{\sqrt{1-c}}} K_2(v(s),\tilde{v},\xi,c)d\tilde{v}\cdot W(s,\xi)(1-V(s))^{\f12},
	\end{align}
	where $s\leq r$, $v(r)\leq \tilde{v}\leq v(s)$, $W$ is defined below \eqref{eq:Phi7}, and $K_2$ takes the form
	\begin{align}\label{def:K2-5}
		K_2(v,\tilde{v},\xi,c)&:=\left(\f{1+\f{\f{\sqrt{1-c}}{\xi}(\pa_v\Phi)(v_*)}{\Phi(v_*)}\tanh \f{\xi(v-v_*)}{\sqrt{1-c}}}
		{1+\f{\f{\sqrt{1-c}}{\xi}(\pa_v\Phi)(v_*)}{\Phi(v_*)}\tanh \f{\xi(\tilde{v}-v_*)}{\sqrt{1-c}}}\right)^2.
	\end{align}
	Observing that $\f{\cosh y}{\cosh z}\leq \f{1}{\cosh(z-y)}(z\geq y\geq 0)$ and $v_*-\tilde{v}\geq v_*-v(s)\geq 0$, we have
	\begin{align}   \label{bd:cosh}
		\f{\cosh^2 \f{\xi(v_*-v(s))}{\sqrt{1-c}} }{\cosh^2 \f{\xi(v_*-\tilde{v})}{\sqrt{1-c}}}  \leq  \cosh^{-2} \f{\xi(v(s)-\tilde{v})}{\sqrt{1-c}}\sim e^{-2\f{\xi(v(s)-\tilde{v})}{\sqrt{1-c}}}.
	\end{align}
	Therefore, using the definitions, \eqref{bd:WPhi7} and  \eqref{est:Phi7derivative}, we infer that for $r_*\leq s\leq \tilde{r}\leq r$,
	\begin{align}\label{phi5K1}
		|K_1(\tilde{r},s,\xi,c)|&\lesssim \int_{v(s)}^{v(r)}e^{-2\f{\xi(v(s)-\tilde{v})}{\sqrt{1-c}}} d\tilde{v}\cdot \f{\langle s\rangle^{\f32}}{s^2}\lesssim  \f{\sqrt{1-c}}{\xi}\f{\langle s\rangle^{\f32}}{s^2}. \end{align}
	This shows that for
	$M\f{\sqrt{1-c}}{\xi}\sim r_*\lesssim r\leq M^{-3.5}\f{\xi^2}{1-c} $,
	\begin{align}\label{phi5K1I}
		\notag \int^{r}_{r_*}\sup_{\tilde{r}\in[s,r] }|K_1(\tilde{r},s,\xi,c)|ds
		&\lesssim \f{\sqrt{1-c}}{\xi}\int_{r_*}^{r}\f{\langle s\rangle^{\f32}}{s^2}ds\\
		&\lesssim \f{\sqrt{1-c}}{\xi}\left(\int_{r_*}^{1}s^{-2}ds
		+\int_{1}^{M^{-3.5}\f{\xi^2}{1-c}}s^{-\f12}ds\right)\lesssim M^{-1},
	\end{align}
	i.e., the  condition \eqref{condition:K'} holds.
	
	To prove \eqref{est:phi7-Derivative}, it suffices to verify the condition \eqref{condition:K'pac-growth} in Lemma \ref{lem:vorteraa-0}, i.e., for  $r\geq 2r_*$,
	\begin{small}
		\begin{align} &\int_{0}^1 \sup_{\tilde{t}\in[t,1]}\left|
			((1-c)\pa_c)^{l} (\xi\pa_{\xi})^{j} (r\pa_r)^{i}\left(
			(r-r_*(\xi,c))K_1(r_*+\tilde{t}(r-r_*),
			r_*+t(r-r_*),\xi,c)\right)\right|dt\lesssim 1, \label{condition:K'pac-growth71}
			\end{align}
	\end{small}
	Indeed, we write by \eqref{def:Phi5-K1} as
	\begin{small}\begin{align*} \notag &K_1(r_*+\tilde{t}(r-r_*),r_*+t(r-r_*),\xi,c)\\
			&=-(r-r_*)\int_{t}^{\tilde{t}}\f{\cosh^2 \f{\xi(v(r_*+t(r-r_*))-v(r_*))}{\sqrt{1-c}} }{\cosh^2 \f{\xi(v(r_*+a(r-r_*))-v(r_*))}{\sqrt{1-c}}} \\ &\qquad \cdot K_2\big(v(r_*+t(r-r_*)),v(r_*+a(r-r_*)),\xi,c\big)
			\big(1-V(r_*+a(r-r_*))\big)^{\f12}da\\
			&\quad\cdot W\big(r_*+t(r-r_*),\xi\big)\big(1-V(r_*+t(r-r_*))\big)^{\f12},
		\end{align*}
	\end{small}
	where $2r_*\leq r$, $t\leq a\leq \tilde{t}\leq 1$, $W$ is defined below \eqref{eq:Phi7}, and  $K_2$ defined in \eqref{def:K2-5} satisfies (using \eqref{est:Phi7derivative})
	\begin{align}\label{bd:K2-7}
		\left|\left((1-c)\pa_c\right)^l(\xi\pa_{\xi})^j(s\pa_s)^{i_2}(r\pa_r)^i
		\left(K_2(v(r),v(s),\xi,c)\right)\right| \lesssim 1,
	\end{align}
	where the implied constant is independent of $r,s,\xi,c$ and may depend on $M$.
	Observing that $$\left|\pa_z^{k_2}\pa_y^{k_1}\left(\f{\cosh y}{\cosh z}\right)\right|\leq \f{1}{\cosh(z-y)}(z\geq y\geq 0,\;k_1,k_2\geq 0),$$
	it follows from \eqref{bd:K2-7}, \eqref{bd:cosh}, \eqref{est:r*7}, \eqref{bd:WPhi7}, \eqref{est:Phi7derivative}, $|r^kv^{(k)}(r)|\lesssim v(r),\;r^k|V^{(k)}(r)|\lesssim 1-V(r)(k\geq 1)$, $r\leq 2(r-r_*)(\text{for}\;r\geq 2r_*) $, \eqref{phi5K1} and \eqref{phi5K1I} that the left-hand side of \eqref{condition:K'pac-growth71}  can be bounded by
	
	\begin{small}
		\begin{align*} & &\lesssim \int_{0}^1(r-r_*)\f{\sqrt{1-c}}{\xi}\f{\langle (r_*+t(r-r_*))\rangle^{\f32}}{(r_*+t(r-r_*))^2}dt\lesssim  \f{\sqrt{1-c}}{\xi}\int_{r_*}^{r}\f{\langle s\rangle^{\f32}}{s^2}ds\lesssim 1, \end{align*}
	\end{small}
	for $M\f{\sqrt{1-c}}{\xi}\sim r_*\lesssim r\leq M^{-3.5}\f{\xi^2}{1-c}  $.
	
	In the end, we verify \eqref{bd:C7}. It follows from Lemma \ref{lem: phi-2} and the first bound in \eqref{est:r*7} that $\phi(r_*)\sim  \xi r_*$ and $\phi(r_*)=O_{\xi,c}^{\xi,1-c}(\xi r_*)$. Then we deduce by the  first bound in \eqref{est:r*7} that
	\begin{align*} &\f{1}{C_{71}(\xi,c)}=O_{\xi,c}^{\xi,1-c}\left(\left(\xi r_*^{\f32}\right)^{-1}\right)
		=\f{1}{\sqrt{1-c}}O_{\xi,c}^{\xi,1-c}\left(\left(\f{\xi}{(1-c)^{\f12}}\right)^{\f12}
		\right). \end{align*}
		
	Thus, we finish the proof of the lemma.
							\end{proof}
							
							\begin{lemma}\label{lem:Phi82}
								
								Let $M\gg1$, $\delta\ll 1$ be fixed. Let $c\in(1-\delta,1)$, $M^2(1-c)^{\f12}\lesssim \xi \lesssim M(1-c)^{\f13}$, and $\phi$ be the solution to \eqref{eq:phi},  $\phi''+r^{-1}\phi'-r^{-2}\phi+\f{\xi^2 (V-c)}{1-c} \phi=0,\; \phi(r,\xi,c)\sim \xi r,$
								which has been solved on $(0,r_*)$ and $[r_*, r_{*2}]$ by  in Lemma \ref{lem:phi-sqrt(1-c)/xi} and Lemma \ref{lem:Phi81}, respectively,  with $r_{*}=M\f{\sqrt{1-c}}{\xi}$, $r_{*2}=M^{-3.5}\f{\xi^2}{1-c}$, $r_{*}\ll1 \ll r_{*2}\ll \xi^{-1}$.
								Then there exists a constant $C_{72}(\xi,c)=\f{\phi({r_{*2})}}{J_1(\xi r_{*2})}$ such that for $M^{-3.5}\f{\xi^2}{1-c}\lesssim r\lesssim M^{-\f12}\xi^{-1}$,
								there holds that
								\begin{align}
									\begin{split}\label{phi72}
										\phi(r,\xi,c)&\sim C_{72}(\xi,c)J_1(\xi r),\\
										\phi(r,\xi,c)&= C_{72}(\xi,c)J_1(\xi r)O_{r,\xi,c}^{r,\xi,1-c}(1).
									\end{split}
								\end{align}
								More precisely, the following bounds hold for the constants \begin{align}
									\begin{split}\label{est:C72}
										&\f{1}{C_{72}(\xi,c)}=\f{1}{\sqrt{1-c}}
										O_{\xi,c}^{\xi,1-c}\left(\left(\f{\xi}{(1-c)^{\f13}}\right)^3\right) \cdot O_{\xi,c}^{\xi,1-c}\left(\left(\f{\xi}{(1-c)^{\f12}}\right)^{-N}\right)\;\;\;\;\;N\geq 1,\\
										&\f{C_{71}(\xi,c)}{C_{72}(\xi,c)}=
										\f{O_{\xi,c}^{\xi,1-c}\left(\left(\f{\xi}{(1-c)^{\f13}}\right)^3\right)\cdot O_{\xi,c}^{\xi,1-c}\left(\left(\f{\xi}{(1-c)^{\f12}}\right)^{-\f12}\right)}{ \cosh\f{\xi(v(r_*)-v(r_{*2}))}{\sqrt{1-c}}},
									\end{split}
								\end{align}
								where 
								$C_{71}$ is defined in \eqref{def:C5},  Lemma \ref{lem:Phi81}.
							\end{lemma}
							
							\begin{proof}
								The proof is similar as the one of Lemma \ref{lem:phi-1*}, where we viewed $\f{\xi^2(1-V(r))}{1-c}$ as a perturbation and used both its positivity and the positivity of the Bessel's function $J_1(\xi r)$  when $r\ll \xi^{-1}$. While a difference in this case lies in that in the range $M^2(1-c)^{\f12}\lesssim \xi$, we fail to get a natural upper bound of $\f{\xi^2}{1-c}$. Fortunately, we have \emph{additionally $r\gtrsim r_{*2}:=M^{-3.5}\f{\xi^2}{1-c}\gg 1$} to get
								\begin{align} \label{bd:xi2(1-V)/1-c}
									\f{\xi^2(1-V(r))}{1-c}\sim \f{\xi^2r^{-3}}{1-c} \lesssim M^{3.5}r^{-2}\lesssim M^{2.5}.
								\end{align}
								The solution solved in $(0, r_{*2}]$ in Lemma \ref{lem:Phi81} is \begin{align}\label{use:Phi07}
									\begin{split}
										\phi(r)&=\left(1-V(r)\right)^{-\f14}r^{-\f12}\Phi_0(v(r))\left(1+\phi^{Rem}(r)\right), \\
										\Phi_0(v)&=C_{71}(\xi,c)\cosh \f{\xi(v_*-v)}{\sqrt{1-c}} \left(1-\f{\f{\sqrt{1-c}}{\xi}(\pa_v\Phi)(v_*)}{\Phi(v_*)}\tanh \f{\xi(v_*-v)}{\sqrt{1-c}} \right).
									\end{split}
								\end{align}
								For the new interval $[r_{*2}, M^{-\f12}\xi^{-1}]$, we rewrite the equation as \begin{align}\label{eq:phi82}
									\phi''+r^{-1}\phi'-r^{-2}\phi+\xi^2\phi=\f{\xi^2(1-V)}{1-c}\phi.
								\end{align}
								To smoothly extend the solution on new interval, we match the initial data at $r_{*2}$ as follows
								\begin{small}\begin{align}
										\begin{split}\label{eq:Phi7initial2}
											\phi(r_{*2})&=\phi(r)\big|_{r=r_{*2}}
											=\left(1-V(r)\right)^{-\f14}r^{-\f12}\Phi_0\left(v(r)\right)\left(1+\phi^{Rem}(r)\right)
											\big|_{r=r_{*2}},\\
											\phi'(r_{*2})&=\phi'(r)\big|_{r=r_{*2}}=\left(\left(1-V(r)\right)^{-\f14}
											r^{-\f12}\Phi_0\left(v(r)\right)\left(1+\phi^{Rem}(r)\right)
											\right)'\big|_{r=r_{*2}}. \end{split}
										\end{align}
								\end{small}
								Taking the right-hand side to be zero in \eqref{eq:phi82}, the solution of the homogeneous equation with initial data \eqref{eq:Phi7initial2} writes
								\begin{align}
									\begin{split}\label{def:Phi0 alt}
										\phi_0(r,\xi,c)=&C_{72}(\xi,c)J_1(\xi r) \left(1+\left(\f{\phi'(r_{*2})}{\phi(r_{*2})}- \f{\xi J_1'(\xi r_{*2})}{J_1(\xi r_{*2})}\right)r_{*2}J_1^2(\xi r_{*2})\int_{r_{*2}}^rs^{-1}J_1(\xi s)^{-2}ds \right),\\ \text{with}\quad &C_{72}(\xi,c)=\f{\phi(r_{*2})}{J_1(\xi r_{*2})}.
								\end{split}\end{align}
							Now we give an illustration for \eqref{est:C72}. Indeed, the second bound follows from \eqref{eq:Phi7initial2}, \eqref{use:Phi07}, \eqref{est:Phi7derivative},  \eqref{claim:Phi7} and \eqref{est:phi7-Derivative} that
							 \begin{align*}
								\f{C_{71}(\xi,c)}{C_{72}(\xi,c)}&=\f{J_1(\xi r_{*2})C _{7_1}}{\phi(r_{*2})}\\
								&=\f{J_1(\xi r_{*2})\left(1-V(r_{*2})\right)^{\f14} r_{*2}^{\f12}}{ \left(1-\f{\f{\sqrt{1-c}}{\xi}(\pa_v\Phi)(v_*)}{\Phi(v_*)}\tanh \f{\xi(v(r_*)-v(r_{*2}))}{\sqrt{1-c}} \right)\left(1+\phi^{Rem}(r_{*2})\right)}
								\cdot \cosh^{-1} \f{\xi(v(r_*)-v(r_{*2}))}{\sqrt{1-c}}\\
								&=O_{\xi,c}^{\xi,1-c}(\xi r_{*2}\cdot r_{*2}^{-\f14}) \cdot \cosh^{-1} \f{\xi(v(r_*)-v(r_{*2}))}{\sqrt{1-c}}\\
								&=O_{\xi,c}^{\xi,1-c}\left(\left(\f{\xi}{(1-c)^{\f13}}\right)^3\right)\cdot O_{\xi,c}^{\xi,1-c}\left(\left(\f{\xi}{(1-c)^{\f12}}\right)^{-\f12}\right) \cdot \cosh^{-1} \f{\xi(v(r_*)-v(r_{*2}))}{\sqrt{1-c}}.
								\end{align*}
								Then the first bound in \eqref{est:C72}  follows from  \eqref{bd:C7} and  \begin{align}\label{key:v_*-v_*2}
									\f{\xi(v(r_*)-v(r_{*2}))}{\sqrt{1-c}}\sim  \f{\xi}{\sqrt{1-c}}
								\end{align}
								that
								\begin{small}
									\begin{align*}
										&\quad\f{1}{C_{72}(\xi,c)}=\f{C_{71}}{C_{72}}\cdot\f{1}{C_{71}}\\
										&=O_{\xi,c}^{\xi,1-c}\left(\left(\f{\xi}{(1-c)^{\f13}}\right)^3\right)\cdot O_{\xi,c}^{\xi,1-c}\left(\left(\f{\xi}{(1-c)^{\f12}}\right)^{-\f12}\right) \cdot \cosh^{-1} \f{\xi(v(r_*)-v(r_{*2}))}{\sqrt{1-c}}\cdot  \f{1}{\sqrt{1-c}}O_{\xi,c}^{\xi,1-c}\left(\left(\f{\xi}{(1-c)^{\f12}}\right)^{\f12}
										\right)\\
										&=\f{1}{\sqrt{1-c}}
										O_{\xi,c}^{\xi,1-c}\left(\left(\f{\xi}{(1-c)^{\f13}}\right)^3\right) \cdot \cosh^{-1} \f{\xi(v(r_*)-v(r_{*2}))}{\sqrt{1-c}}\\
										&=\f{1}{\sqrt{1-c}}
										O_{\xi,c}^{\xi,1-c}\left(\left(\f{\xi}{(1-c)^{\f13}}\right)^3\right) \cdot O_{\xi,c}^{\xi,1-c}\left(\left(\f{\xi}{(1-c)^{\f12}}\right)^{-N}\right),
										\end{align*}
								\end{small}
								where $N\geq 1$.
								
								 Now we prove \eqref{phi72}. Using the form \eqref{def:Phi0}, we rewrite the equation \eqref{eq:phi82} with the matched initial data  into an integral form as
								\begin{align}  \label{integral:phi72}
									\phi(r,\xi,c)=\phi_0(r,\xi,c)+\int_{r_{*2}}^r
									K(r,s,\xi,c)\phi(s,\xi,c)ds,
								\end{align}
								where
								\begin{align*}
									K(r,s,\xi,c)=\phi_0(r)\phi_0(s)
									\int_s^r\phi_0(\tilde{s})^{-2}d\tilde{s}\cdot\f{\xi^2(1-V(s))}{1-c}\;\;\text{with}\;\;K(r,r)=0.
									\end{align*}
								For $r_{*2}\leq r\leq M^{-\f12}\xi^{-1}$, we claim that
								\begin{align}
									&\left(\f{\phi'(r_{*2})}{\phi(r_{*2})}- \f{\xi J_1'(\xi r_{*2})}{J_1(\xi r_{*2})}\right)r_{*2}J_1^2(\xi r_{*2})\int_{r_{*2}}^rs^{-1}J_1(\xi s)^{-2}ds\sim M^{\f74}.\label{est:1+C84-phi82}
									\end{align}
								For $i,j\in\mathbb{N}$, $l\in\{0,1\}$,
								\begin{align}
									&\left|(\xi\pa_{\xi})^j\left((1-c)\pa_c\right)^l
									\left(r\pa_r\right)^i\left(\left(\f{\phi'(r_{*2})}{\phi(r_{*2})}- \f{\xi J_1'(\xi r_{*2})}{J_1(\xi r_{*2})}\right)r_{*2}J_1^2(\xi r_{*2})\int_{r_{*2}}^rs^{-1}J_1(\xi s)^{-2}ds\right)\right|\lesssim 1,\label{est:1+C84-phi82-d}
								\end{align}
								and
								\begin{align}
									\phi(r,\xi,c)=\phi_0(r,\xi,c)
									\left(1+\phi^{Rem}(r,\xi,c)\right),\;0<\phi^{Rem}(r)\lesssim M^{\f72}\label{phi82-claim}.
								\end{align}
								\begin{align}
									&\left|\left((1-c)\pa_c\right)^l\left(r\pa_r\right)^i(\xi\pa_{\xi})^j\phi^{Rem}_2(r,\xi,c)\right|
									\lesssim 1,\quad\text{for}\quad 2r_{*2}\leq r\leq M^{-\f12}\xi^{-1}.\label{est:phirem-phi82}
								\end{align}
								Here, the upper implied bounds in  \eqref{est:1+C84-phi82-d} and     \eqref{est:phirem-phi82}   may depend on $M$.
								Finally, \eqref{est:1+C84-phi82}, \eqref{est:1+C84-phi82-d}, \eqref{phi82-claim} and \eqref{phi82-claim} together with the form of $\phi_0$ \eqref{def:Phi0} imply \eqref{phi72}. \eqref{est:1+C84-phi82} can be reduced to the following
								\begin{align*}
									\f{r\pa_r\phi}{\phi}\big|_{r=r_{*2}}\sim M^{\f74},\quad
									\left|\f{\xi rJ_1'(\xi r)}{J_1(\xi r)}\big|_{r=r_{*2}}\right|\lesssim 1,\quad \left|J_1^2(\xi r_{*2})\int_{r_{*2}}^rs^{-1}J_1(\xi s)^{-2}ds\right|\lesssim 1.
								\end{align*}
								Noticing
								\begin{align}  \label{bd:J_1}
									J_1(z)\sim z\;\;\;\text{for}\;\;z\ll 1,\;\;\text{and}\;\; J_1(z)=O_z^z(1)\;\;\;\text{for}\;\;z\lesssim 1. \end{align}
								the second and third bound follows directly by $\xi r_{*2}=M^{-3.5}\f{\xi^3}{1-c}\lesssim M^{-\f12}\ll 1$.  For the first bound, we use the formulation \eqref{eq:Phi7initial2} and \eqref{use:Phi07} to get
								\begin{align*}
									\f{r\pa_r\phi}{\phi}\big|_{r=r_{*2}}
									&=\f{r\pa_r\left(\left(1-V(r)\right)^{-\f14}\right)}{
										\left(1-V(r)\right)^{-\f14}}\big|_{r=r_{*2}}
									+\f{r\pa_r\left(r^{-\f12}\right)}{
										r^{-\f12}}\big|_{r=r_{*2}}
									+\f{r\pa_r\left(\cosh \f{\xi(v_*-v(r))}{\sqrt{1-c}} \right)}{
										\cosh \f{\xi(v_*-v(r))}{\sqrt{1-c}} }\big|_{r=r_{*2}}\\
									&\quad+\f{r\pa_r\left( 1-\f{\f{\sqrt{1-c}}{\xi}(\pa_v\Phi)(v_*)}{\Phi(v_*)}\tanh \f{\xi(v_*-v(r))}{\sqrt{1-c}} \right)}{
										1-\f{\f{\sqrt{1-c}}{\xi}(\pa_v\Phi)(v_*)}{\Phi(v_*)}\tanh \f{\xi(v_*-v(r))}{\sqrt{1-c}}}\big|_{r=r_{*2}}
									+\f{r\pa_r\phi^{Rem}(r)}{
										1+\phi^{Rem}(r)}\big|_{r=r_{*2}}\\
									&:=I_1+...+I_5.
								\end{align*}
								 Indeed, it follows that
								 \begin{align*}
								 |I_3|\sim M^{\f74},\quad  |I_1|+|I_2|\lesssim 1,\quad |I_4|+|I_5|\lesssim M^{\f34}.
								 \end{align*}
					By \eqref{key:v_*-v_*2}, $rv'(r)\sim r^{-\f12}(r\gtrsim 1)$ and $\tanh z\sim 1 (z\sim 1)$, we deduce that
					\begin{align*}
					|I_3|\sim \f{\xi r_{*2}^{-\f12}}{\sqrt{1-c}}\cdot \f{\sinh \f{\xi(v(r_*)-v(r_{*2}))}{\sqrt{1-c}} }{
									\cosh \f{\xi((v(r_*)-v(r_{*2}))}{\sqrt{1-c}} }\sim \f{\xi r_{*2}^{-\f12}}{\sqrt{1-c}}\sim M^{\f74}.
									\end{align*}
									 $|I_1|+|I_2|\lesssim 1$ is obvious. $|I_4|\lesssim M^{\f34}$ follows further from \eqref{est:Phi7derivative}. $|I_5|\lesssim M^{\f34}$ follows by tracking the $M$ dependence in the prove of \eqref{condition:K'pac-growth}. \eqref{est:1+C84-phi82-d} can be reduced to the following
									\begin{align} \label{bd:phi72K2}
									\begin{split}
										&\left|(\xi\pa_{\xi})^j\left((1-c)\pa_c\right)^l
										\left(\f{r_{*2}\phi'(r_{*2})}{\phi(r_{*2})}\right)\right|\lesssim 1,\\
										&\left|(\xi\pa_{\xi})^j\left((1-c)\pa_c\right)^l
										\left( \f{\xi r_{*2} J_1'(\xi r_{*2})}{J_1(\xi r_{*2})}\right)\right|\lesssim 1, \\
										&\left|(\xi\pa_{\xi})^j\left((1-c)\pa_c\right)^l
										\left(r\pa_r\right)^i\left(J_1^2(\xi r_{*2})\int_{r_{*2}}^rs^{-1}J_1(\xi s)^{-2}ds\right)\right|\lesssim 1.
									\end{split}
									\end{align}
									Here, the bound is independent of $r,\xi,c$, and may depend on $M$. We leave the proof to the readers.
		                           To prove \eqref{phi82-claim}, we apply Lemma \ref{lem:vorteraa-0} to the integral equation \eqref{integral:phi72}. It suffices to check that
									for $r_{*2}\leq s\leq r\leq M^{-\f12}\xi^{-1}$,
									\begin{align*}
										\int_{r_{*2}}^r\sup_{\tilde{r}\in [r_{*2},r]} |K_1(\tilde{r},s,\xi,c)|ds\lesssim M^{\f72}\quad\text{and}\quad\;K_1(r,s,\xi,c)>0,
									\end{align*}
									where \begin{align}\label{def:Phi52-K1}
										\notag K_1(r,s,\xi,c)&:=\phi_0^{-1}(r)
										K(r,s)\phi_0(s)\\
										\notag &=\phi_0(s,\xi,c)^2
										\int_{s}^{r}
										\phi_0(\tilde{r},\xi,c)^{-2}
										d\tilde{r}\cdot\f{\xi^2(1-V(s))}{1-c}\\
										&=\int_{s}^{r}\f{J_1(\xi s)^2}{J_1(\xi \tilde{r})^2} K_2(s,\tilde{r},\xi,c)d\tilde{r}\cdot\f{\xi^2(1-V(s))}{1-c}>0,
									\end{align}
									where  $\phi_0$ is defined in \eqref{def:Phi0} and
									\begin{align}\label{def:K2-52}
										K_2(s,\tilde{r},\xi,c)&:=\left(\f{1+\left(\f{\phi'(r_{*2})}{\phi(r_{*2})}- \f{\xi J_1'(\xi r_{*2})}{J_1(\xi r_{*2})}\right)r_{*2}J_1^2(\xi r_{*2})\int_{r_{*2}}^s\tilde{s}^{-1}J_1(\xi\tilde{ s})^{-2}d\tilde{s} }{1+\left(\f{\phi'(r_{*2})}{\phi(r_{*2})}- \f{\xi J_1'(\xi r_{*2})}{J_1(\xi r_{*2})}\right)r_{*2}J_1^2(\xi r_{*2})\int_{r_{*2}}^{\tilde{r}}\tilde{s}^{-1}J_1(\xi \tilde{s})^{-2}d\tilde{s } }\right)^2. \end{align}
									It follows from \eqref{est:1+C84-phi82} that for $r_{*2}\leq s\leq  \tilde{r}\leq M^{-\f12}\xi^{-1}$,
									\begin{align}\label{bd:K2-52}
										0<K_2(s,\tilde{r},\xi,c)\sim 1.
										\end{align}
										By \eqref{bd:J_1},  for  $r_{*2}\leq s\leq  \tilde{r}\leq M^{-\f12}\xi^{-1}\ll \xi^{-1}$, we have
										\begin{align}\label{bd:J/J-72}
											0<\f{J_1(\xi s)^2}{J_1(\xi \tilde{r})^2} \sim \f{s^2}{ \tilde{r}^2} .
											\end{align}
										Then for $1\ll r_{*2}\leq s\leq  \tilde{r}\leq M^{-\f12}\xi^{-1}$, we have
										\begin{align}\label{positive:K1-52}
											0<K_1(\tilde{r},s,\xi,c)\lesssim \f{\xi^2s^{-2}}{1-c}.
											 \end{align}
										Therefore, we infer  that for   $M^{-3.5}\f{\xi^2}{1-c}\leq  r\leq M^{-\f12}\xi^{-1}$,
										\begin{align}\label{intK-phi82}
											\int_{r_{*2}}^r\sup_{\tilde{r}\in [r_{*2},r]} |K_1(\tilde{r},s,\xi,c)|ds\lesssim
											\int_{r_{*2}}^r\f{\xi^2s ^{-2}}{1-c}ds\leq \f{\xi^2r_{*2} ^{-1}}{1-c}
											=M^{\f72}.
										\end{align}
										
										To prove \eqref{est:phirem-phi82},  it suffices to verify the condition  \eqref{condition:K'pac-growth} in Lemma \ref{lem:vorteraa-0},
										i.e., for $2r_{*2}\leq r\leq M^{-\f12}\xi^{-1}$:
										\begin{small}
											\begin{align}
												&\int_{0}^1 \sup_{\tilde{t}\in[t,1]}\left|
												((1-c)\pa_c)^{l} (\xi\pa_{\xi})^{j} (r\pa_r)^{i}\left(
												(r-r_{*2}(\xi,c))K_1(r_{*2}+\tilde{t}(r-r_{*2}),
												r_{*2}+t(r-r_{*2}),\xi,c)  \right)\right|dt\lesssim 1,\label{condition:K'pac-growth72}
												\end{align}
										\end{small}
										Indeed, we write by \eqref{def:K2-52} that
										\begin{small}
											\begin{align*} \notag & K_1(r_{*2}+\tilde{t}(r-r_{*2}),r_{*2}+t(r-r_{*2}),\xi,c)\\
												&=
												-(r-r_{*2})\int_{t}^{\tilde{t}}\f{J_1\left(\xi(r_{*2}+t(r-r_{*2}))\right)^2}{
													J_1\left(\xi (r_{*2}+a(r-r_{*2}))\right)^2} \cdot K_2\left(r_{*2}+t(r-r_{*2}),r_{*2}+a(r-r_{*2}),\xi,c\right)
												da\\
												&\quad\cdot \f{\xi^2(1-V(r_{*2}+t(r-r_{*2})))}{1-c},
											\end{align*}
										\end{small}
										where $2r_{*2}\leq r$, $t\leq a\leq \tilde{t}\leq 1$, and $K_2$ defined in \eqref{def:K2-52} satisfies (using \eqref{bd:phi72K2})
										\begin{align}\label{bd:K2-72}
											\left|\left((1-c)\pa_c\right)^l(\xi\pa_{\xi})^j(s\pa_s)^{i_2}(r\pa_r)^i
											\left(K_2(s,r,\xi,c)\right)\right| \lesssim 1,
										\end{align}
										where the implied constant is independent of $r,s,\xi,c$ and may depends on $M$.  Finally, the left-hand side of \eqref{condition:K'pac-growth72}  can be bounded by
										\begin{small}
											\begin{align*} &\lesssim \int_{0}^1(r-r_{*2})\f{\xi^2\left(r_{*2}+t(r-r_{*2})\right) ^{-2}}{1-c}dt\sim \int_{r_{*2}}^r\f{\xi^2s ^{-2}}{1-c}ds\lesssim 1
											\end{align*}
										\end{small}
										for $2r_{*2}\leq r\leq M^{-\f12}\xi^{-1} $.
							\end{proof}
							
							\section{Global behavior of Schr\"odinger equation}

							In this section, by combining the behavior of the Schr\"odinger equation in the different regions, we establish its global behavior.

							\begin{definition} \label{def:W0rhoc}
						        Let $\delta \ll 1$ and $M\gg 1$ be fixed. For $\xi>0$, we define
								\begin{align}\label{behave:Wron0}
									W_0(\xi,c)\sim
									\left\{
									\begin{array}{l}
										1\;\;\; \;\;\quad\quad \quad\quad \; c\in(0,1),\;\xi\lesssim M^2(1-c)^{\f12},\\
										\xi^{\f13}\;\;\; \quad\quad \quad\quad\; c\in(0,V(0)),\;1\lesssim \xi\lesssim |c-V(0)|^{-\f32},\\
										|c-V(0)|^{-\f12}\;\;\; c\in(0,V(0)),\;\xi\gtrsim |c-V(0)|^{-\f32},\\
										\xi^{\f13}\;\;\; \quad\quad \quad\quad  \; c\in(V(0),1-\delta),\;1\lesssim \xi\lesssim |c-V(0)|^{-\f32},\\
										C_{62}(\xi,c)\;\;\quad\quad c\in(V(0),1),\;\xi\gtrsim M(1-c)^{\f13}|c-V(0)|^{-\f32},\\
										C_{72}(\xi,c)\;\;\quad\quad c\in(1-\delta,1),\;M^2(1-c)^{\f12}\lesssim \xi\lesssim M|1-c|^{\f13},
										\end{array}\right.
								\end{align}
								 and $W_0(\xi)=W_0(-\xi)$ for $\xi<0$.
								 Here $C_{61}>0$ is defined in Lemma \ref{lem:f+3} and behaves as \eqref{est:C1-1},  and $C_{72}>0$ is defined in Lemma \ref{lem:Phi82} and behaves as \eqref{est:C72}.
								
								We also define a weight function as follows: for $\xi>0$,
								\begin{align}\label{behave:rho}
									\rho_0(\xi,c)=
									\left\{
									\begin{array}{l}
										1-c\;\;\; \quad\;\;\quad\quad \quad\quad\quad c\in(0,1),\;\xi\lesssim (1-c)^{\f12},\\
										\xi^{-\f23}\;\;\;\quad\quad\quad\quad \quad\quad\quad  c\in(0,V(0)),\;1\lesssim \xi\lesssim |c-V(0)|^{-\f32},\\
										c-V(0)\;\;\;\quad\quad\quad\quad\quad  c\in(0,V(0)),\;\xi\gtrsim |c-V(0)|^{-\f32},\\
										\xi^{-\f23}\;\;\;\quad\quad\quad\quad \quad\quad\quad c\in(V(0),1-\delta),\;1\lesssim \xi\lesssim |c-V(0)|^{-\f32},\\
										(1-c)|c-V(0)|\;\;\; \quad\quad  c\in(V(0),1),\;\xi\gtrsim M(1-c)^{\f13}|c-V(0)|^{-\f32},\\
										1-c\;\;\;\quad\quad\quad\quad\quad\quad\quad  c\in(1-\delta,1),\;\delta\ll 1,\;M^2(1-c)^{\f12}\lesssim \xi\lesssim |1-c|^{\f13},
									\end{array}\right.
								\end{align}
								and $\rho_0(\xi)=\rho_0(-\xi)$ for $\xi<0$.
								\end{definition}
It is obvious that $W_0(\xi,c)=O_{\xi,c}^{\xi,\rho_0(\xi,c)}
									\left(W_0(\xi,c)\right)$.								Recall the Definition \eqref{def:derivative}:
									 \begin{align*}
										\mathcal{D}_{l,m}=\left(\f{(1-c)^{\f12}\pa_r}{\xi}\right)^{m}
										\left(\f{(1-c)^{\f12}(\pa_r+r^{-1})}{\xi}\right)^{l},\quad
										\mathcal{D}_{l,m}^L=\pa_r^{m}
										\left(\f{(1-c)^{\f12}(\pa_r+r^{-1})}{\xi}\right)^{l},\;\,\, l\in\{0,1\},m\in\mathbb{N}.
										\end{align*}

								\begin{proposition}  \label{prop:summery}
									Let $c\in(0,1)\setminus\{V(0)\}$, $\xi\in\mathbb{R}\setminus\{0\}$ be parameters, and  $\phi(r,\xi,c)$ be the real-valued solution of \eqref{eq:phi},
									$f_+(r,\xi,c)$ be the complex-valued solution of \eqref{eq:f+}. The Wronskian $W(\xi,c)$ is as \eqref{fm:wrons-new}, and $W_0(\xi,c)$ and $\rho_0(\xi,c)$ are defined in Definition \ref{def:W0rhoc}. We  list the bound for $\xi>0$, the bounds for $\xi<0$ are identical with respects to $-\xi$ since the functions are defined symmetrically:
\begin{align*}
										\phi(\xi):=-\phi(-\xi),\;f_+(\xi):=-f_+(-\xi),\;W(\xi)=W(-\xi).
\end{align*}
 
It holds that
									\begin{align} \label{behaveW1}
										W(\xi,c)=O_{\xi,c}^{\xi,\rho_0(\xi,c)}\big(W_0(\xi,c)\big)\quad\text{and}\quad  W(\xi,c)\sim W_0(\xi,c),
									\end{align}
									and
									\begin{align} \label{behaveW2}
									\f{W(\xi,c)}{|W(\xi,c)|} =O_{\xi,c}^{\xi,\rho_0(\xi,c)}(1).
									\end{align}
									
									Let
									$\widetilde{\phi}(r,\xi,c)=\f{\phi(r,\xi,c)}{|W(\xi,c)|}$ be as  in \eqref{def:tphi}. Then the following bounds hold.
									
									\begin{itemize}
										\item[(1)] $c\in(0,1)$,\;$\xi\lesssim M^2(1-c)^{\f12}$. It holds that for $r\lesssim \xi^{-1}$,
										\begin{align}
											\begin{split}\label{1-phi}
												&\widetilde{\phi}(\xi,c,r)=
												O_{\xi,c,r}^{\xi,1-c,r}(\xi r);\\ &(\mathcal{D}^L_{l,m}\widetilde{\phi})(\xi,c,r)=O_{\xi,c}^{\xi,1-c}
												(\xi)\;(m+l\geq 1);\\
												&(\mathcal{D}_{1,0}\widetilde{\phi})(\xi,c,r)=(1-c)^{\f12}O_{\xi,c}^{\xi,1-c}
												(1).
										\end{split}\end{align}
										It holds that for $r\gtrsim \xi^{-1}$,
										\begin{align}
											\begin{split}\label{1-f+}
												&f_+(\xi,c,r)=O_{\xi,c,r}^{\xi,1-c,r}\left((\xi r)^{-\f12}\right)e^{i\xi r};\\
												&(\mathcal{D}^L_{l,m}f_+)(\xi,c,r)=\xi O_{\xi,c}^{\xi,1-c}\left((\xi r)^{-\f12}\right)e^{i\xi r}\;(m+l\geq 1);\\
												& (\mathcal{D}_{1,0}f_+)(\xi,c,r)=O_{\xi,c}^{\xi,1-c}\left((\xi r)^{-\f12}\right)e^{i\xi r}.
											\end{split}
										\end{align}
In the following cases, $\xi$ is rather large, we only need to consider the derivatives of corresponding to the higher frequencies $\mathcal{D}_{l,m}=\left(\f{(1-c)^{\f12}\pa_r}{\xi}\right)^{m}
										\left(\f{(1-c)^{\f12}(\pa_r+r^{-1})}{\xi}\right)^{l}$.										
										\item[(2)] $c\in(0,V(0))$,\;$1\lesssim \xi\lesssim |c-V(0)|^{-\f32}$. It holds that for $r\lesssim \xi^{-\f23}$,
										\begin{align}
											\begin{split}\label{2-phi}
												&\widetilde{\phi}(\xi,c,r)=
												O_{\xi,c,r}^{\xi,\xi^{-\f23},r}\left(\xi^{\f23} r\right);\\
												&(\mathcal{D}_{l,m}
												\widetilde{\phi})(\xi,c,r)=O_{\xi,c}^{\xi,\xi^{-\f23}}(1 ).
											\end{split}
										\end{align}
										It holds that for $r\gtrsim  \xi^{-\f23}$,
										\begin{align} \begin{split}
												\label{2-f+}
												&f_+(\xi,c,r)=
												O_{\xi,c,r}^{\xi,\xi^{-\f23},r}\left(\left(\xi\int_0^r \sqrt{\f{V(s)-c}{1-c}}ds\right)^{-\f12}\right)e^{i\xi\int_0^r \sqrt{\f{V(s)-c}{1-c}}ds};\\
												&\left(\xi^{-1}\pa_r\right)^n
												(\mathcal{D}_{l,m}f_+)(\xi,c,r)=O_{\xi,c}^{\xi,\xi^{-\f23}}
												\left(\left(\xi^{\f23} r\right)^{-\f12}\right)e^{i\xi\int_0^r \sqrt{\f{V(s)-c}{1-c}}ds}.
											\end{split}
										\end{align}
										\item[(3)] $c\in(0,V(0))$,\;$ \xi\gtrsim |c-V(0)|^{-\f32}$.  It holds that for $r\lesssim \xi^{-1}|c-V(0)|^{-\f12}$,
										\begin{align}
											\begin{split}\label{4-phi}
												&\widetilde{\phi}(\xi,c,r)=
												O_{\xi,c,r}^{\xi, V(0)-c,r}\left( (V(0)-c)^{\f12}\xi r\right);\\
												&(\mathcal{D}_{l,m}
												\widetilde{\phi})(\xi,c,r)=O_{\xi,c}^{\xi,V(0)-c}(1 ). \end{split}
										\end{align}
										It holds that for $r\gtrsim  \xi^{-1}|c-V(0)|^{-\f12}$,
										\begin{align}
											\begin{split}\label{4-f+}
												&f_+(\xi,c,r)=
												O_{\xi,c,r}^{\xi,c-V(0),r}\left(\left(\xi\int_0^r \sqrt{\f{V(s)-c}{1-c}}ds\right)^{-\f12}\right)e^{i\xi\int_0^r \sqrt{\f{V(s)-c}{1-c}}ds}; \\
												&(\mathcal{D}_{l,m}f_+)(\xi,c,r)=
												O_{\xi,c}^{\xi,c-V(0)}\left(\left(\xi\int_0^r \sqrt{\f{V(s)-c}{1-c}}ds\right)^{-\f12}\right)e^{i\xi\int_0^r \sqrt{\f{V(s)-c}{1-c}}ds}.
										\end{split} \end{align}
										\item[(4)] $c\in(V(0),1-\delta)$,\;$1\lesssim \xi\lesssim |c-V(0)|^{-\f32}$. It holds that for $r\lesssim \xi^{-\f23}$,
										\begin{align}
											\begin{split}\label{3-phi}
												&\widetilde{\phi}(\xi,c,r)=
												O_{\xi,c,r}^{\xi,\xi^{-\f23},r}\left(\xi^{\f23} r\right);\\
												&(\mathcal{D}_{l,m}\widetilde{\phi})(\xi,c,r)=
												O_{\xi,c}^{\xi,\xi^{-\f23}}(1).
										\end{split} \end{align}
										In this case, $r_c:=V^{-1}(c)\sim c-V(0)$, $\langle r_c\rangle \sim 1$. It holds that for $\xi^{-\f23}\lesssim r\leq 2r_c$,
										\begin{align}\label{3-f+supplement}
											&f_+(\xi,c,r)=
											O((\xi r^{\f32})^{-\f13}).
										\end{align}
										For $r\geq M\xi^{-\f23}\geq 2r_c$, it holds that
										\begin{align}
											\label{3-f+}
											\begin{split}
												&f_+(\xi,c,r)=
												O_{\xi,c,r}^{\xi,\xi^{-\f23},r}\left(\left(\xi\int_{r_c}^r \sqrt{\f{V(s)-c}{1-c}}ds\right)^{-\f12}\right)e^{i\xi\int_{r_c}^r \sqrt{\f{V(s)-c}{1-c}}ds};\\
												&(\mathcal{D}_{l,m}f_+)(\xi,c,r)=
												O_{\xi,c}^{\xi,\xi^{-\f23}}\left(\left(\xi\int_{r_c}^r \sqrt{\f{V(s)-c}{1-c}}ds\right)^{-\f12}\right)e^{i\xi\int_{r_c}^r \sqrt{\f{V(s)-c}{1-c}}ds}. \end{split}
										\end{align}
										
										\item[(5),(6)]  $c\in(V(0),1)$,\;$ \xi\gtrsim M(1-c)^{\f13}|c-V(0)|^{-\f32}$. In this case, for $r_c=V^{-1}(c)$,\;it holds that
										\begin{align}\label{6-phi-basic}
											(1-c)^{\f13}\sim \langle r_c\rangle^{-1},\; |V(0)-c|\sim \f{r_c}{\langle r_c\rangle},\;\; \xi r_c^{\f32}/\langle r_c\rangle^{\f12}\gtrsim M.
										\end{align}
										It holds that for $ r\lesssim M^{\f12}\xi^{-1}r_c^{-\f12}\langle r_c\rangle^{-1}$,
										\begin{align}
											\begin{split}\label{6-phi-1}
												&\widetilde{\phi}(\xi,c,r)= O_{\xi,c,r}^{\xi,\f{r_c}{\langle r_c\rangle^4},r}\left(\xi r_c^{\f12}\langle r_c\rangle r\right)\cdot
												O_{\xi,c}^{\xi,\f{r_c}{\langle r_c\rangle^4}}\left(\left(\xi r_c^{\f32}/\langle r_c\rangle^{\f12}\right)^{-1}\right);\\
												&(\mathcal{D}_{l,m}\widetilde{\phi})(\xi,c,r)= \left(\f{r_c}{\langle r_c\rangle}\right)^{\f12}
												O_{\xi,c}^{\xi,\f{r_c}{\langle r_c\rangle^4}}\left(\left(\xi r_c^{\f32}/\langle r_c\rangle^{\f12}\right)^{-1}\right).
										\end{split}\end{align}
										If $r_c\lesssim 1$, it holds that for  $M^{\f12}\xi^{-1}r_c^{-\f12}\langle r_c\rangle^{-1}\lesssim r\leq \f{r_c}{2}$,
										\begin{align}\label{6-phi-2}
											\begin{split}
												&\widetilde{\phi}(\xi,c,r)= O_{\xi,c}^{\xi,
													r_c}\left(\left(\xi r_c^{\f12} r\right)^{-\f12}\right);\\
												&(\mathcal{D}_{l,m}\widetilde{\phi})(\xi,c,r)=
												r_c^{\f12}O_{\xi,c}^{\xi,
													r_c}\left(\left(\xi r_c^{\f12} r\right)^{-\f12}\right).
													\end{split}
										\end{align}
										If $r_c\gtrsim 1$, it holds that for $M^{\f12}\xi^{-1}r_c^{-\f12}\langle r_c\rangle^{-1}\lesssim r\leq \f{1}{2}$,
										\begin{align}\label{6-phi-2'}
											\begin{split}
												&\widetilde{\phi}(\xi,c,r)= O_{\xi,c}^{\xi,
													r_c^{-3}}\left(\left(\xi r_c^{\f32} r\right)^{-\f12}\right);\\
												&(\mathcal{D}_{l,m}\widetilde{\phi})(\xi,c,r)= O_{\xi,c}^{\xi,
													r_c^{-3}}\left(\left(\xi r_c^{\f32} r\right)^{-\f12}\right). \end{split}
										\end{align}
										If $r_c\geq \f12$, it holds that for $\f14\leq r\leq \f{r_c}{2}$,
										\begin{align}\label{6-phi-3}
											\begin{split}
												&\widetilde{\phi}(\xi,c,r)= O_{\xi,c}^{\xi, r_c^{-3}}\left(\left( \xi r_c^{\f32}/ r^{\f12}\right)^{-\f12}\right);\\
												&(\mathcal{D}_{l,m}\widetilde{\phi})(\xi,c,r)=
												r_c^{-\f12}O_{\xi,c}^{\xi,
													r_c^{-3}}\left(\left( \xi r_c^{\f32}/ r^{\f12}\right)^{-\f12}\right). \end{split}
										\end{align}
										It holds that for $\f{r_c}{4}\leq r\leq r_c- C\xi^{-\f23}\langle r_c\rangle^{\f13}$,
										\begin{align}\label{6-phi-4} \begin{split}
												&\widetilde{\phi}(\xi,c,r)= O_{\xi,c,r}^{\xi,
													\f{r_c-r}{\langle r_c\rangle^4},\f{(r_c-r)r}{r_c}}\left(\left(\xi r_c(r_c-r)^{\f12}/\langle r_c\rangle^{\f12}\right)^{-\f12}\right);\\
												&(\mathcal{D}_{l,m}\widetilde{\phi})(\xi,c,r)=\f{r_c^{\f12}}{\langle r_c\rangle^2} O_{\xi,c}^{\xi,
													\f{r_c-r}{\langle r_c\rangle^4}}\left(\left(\xi r_c(r_c-r)^{\f12}/\langle r_c\rangle^{\f12}\right)^{-\f12}\right). \end{split}
										\end{align}
										It holds that for $r_c-C\xi^{-\f23}\langle r_c\rangle^{\f13}\leq r\leq r_c+C\xi^{-\f23}\langle r_c\rangle^{\f13}$,
										\begin{align}
											\begin{split}\label{6-f+a1} &f_+(\xi,c,r)=
												O_{\xi,c,r}^{\xi,\xi^{-\f23}\langle r_c\rangle^{-\f{11}{3}},\xi^{-\f23}\langle r_c\rangle^{\f13}}\left(\left(\xi r_c^{\f32}/\langle r_c\rangle^{\f12}\right)^{-\f13}\right);\\
												&(\mathcal{D}_{l,m}f_+)(\xi,c,r)=
												\f{r_c^{\f12}}{\langle r_c\rangle^2} O_{\xi,c}^{\xi,\xi^{-\f23}\langle r_c\rangle^{-\f{11}{3}}}\left(\left(\xi r_c^{\f32}/\langle r_c\rangle^{\f12}\right)^{-\f13}\right).
											\end{split}
										\end{align}
										It holds that for $r\geq r_c+C\xi^{-\f23}\langle r_c\rangle^{\f13}(\geq r_c)$,
										\begin{align}
											\begin{split}\label{6-f+}
												&f_+(\xi,c,r)=
												O_{\xi,c}^{\xi,\f{r-r_c}{\langle r_c\rangle^4}}\left(\left(\xi\int_{r_c}^r \sqrt{\f{V(s)-c}{1-c}}ds\right)^{-\f16}\cdot \left(\xi r_c^{\f32}/\langle r_c\rangle^{\f12}\right)^{-\f13}\right)e^{i\xi\int_{r_c}^r \sqrt{\f{V(s)-c}{1-c}}ds};\\
												&(\mathcal{D}_{l,m}f_+)(\xi,c,r)=
												O_{\xi,c}^{\xi,\f{r-r_c}{\langle r_c\rangle^4}}\left(\left(\xi\int_{r_c}^r \sqrt{\f{V(s)-c}{1-c}}ds\right)^{-\f16}\cdot \left(\xi r_c^{\f32}/\langle r_c\rangle^{\f12}\right)^{-\f13}\right)e^{i\xi\int_{r_c}^r \sqrt{\f{V(s)-c}{1-c}}ds}. \end{split}
										\end{align}

										\item[(7)] For $c\in(1-\delta,1)$,\;$M^2(1-c)^{\f12}\lesssim \xi \lesssim M (1-c)^{\f13}$,  $N\geq 1$, it holds that
										\begin{align}
											\label{7-phi}
											&\widetilde{\phi}(r,\xi,c)
											=\left\{
											\begin{array}{l}
												O_{\xi,c,r}^{\xi,1-c,r}\left(\f{\xi r}{\sqrt{1-c}}\right)\cdot O_{\xi,c}^{\xi,1-c}\left( \left(\f{\xi}{(1-c)^{\f13}}\right)^3\right)\cdot O_{\xi,c}^{\xi,1-c}\left(\left(\f{\xi}{(1-c)^{\f12}}\right)^{-N}\right)\;\;\;\;
												r\lesssim M\xi^{-1}\sqrt{1-c},\\
												O_{\xi,c,r}^{\xi,1-c,r}\left(\left(\f{\xi r}{\sqrt{1-c}}\right)^{-\f12}\right)\cdot O_{\xi,c,r}^{\xi,1-c,r}\left( \left(\f{\xi}{(1-c)^{\f13}}\right)^3\right)\cdot O_{\xi,c}^{\xi,1-c}\left(\left(\f{\xi}{(1-c)^{\f12}}\right)^{-N}\right)\;\;\; M\xi^{-1}\sqrt{1-c}\lesssim r\lesssim 1,\\
												O_{\xi,c,r}^{\xi,1-c,r}\left( \left(\f{\xi}{(1-c)^{\f13}}\right)^3\right)\cdot O_{\xi,c,r}^{\xi,1-c,r}\left(\left(\f{\xi r^{-\f12}}{(1-c)^{\f12}}\right)^{-N-\f12}\right)\;\;\; 1\ll r\lesssim M^{-\f72}\f{\xi^2}{1-c},\\
												O_{\xi,c,r}^{\xi,1-c,r}( \xi r)\;\;\; M^{-\f72}\f{\xi^2}{1-c}\lesssim r\lesssim M^{-\f12}\xi^{-1}, \end{array}\right.
										\end{align}
										and
										\begin{align}\label{7-phi-D}
											&(\mathcal{D}_{l,m}\widetilde{\phi})(r,\xi,c)=
											\left\{
											\begin{array}{l}
												O_{\xi,c}^{\xi,1-c}\left( \left(\f{\xi}{(1-c)^{\f13}}\right)^3\right)\cdot O_{\xi,c}^{\xi,1-c}\left(\left(\f{\xi}{(1-c)^{\f12}}\right)^{-N}\right)\;\;\;\;
												r\lesssim M\xi^{-1}\sqrt{1-c},\\
												O_{\xi,c}^{\xi,1-c}\left(\left(\f{\xi r}{\sqrt{1-c}}\right)^{-\f12}\right)\cdot O_{\xi,c}^{\xi,1-c}\left( \left(\f{\xi}{(1-c)^{\f13}}\right)^3\right)\cdot O_{\xi,c}^{\xi,1-c}\left(\left(\f{\xi}{(1-c)^{\f12}}\right)^{-N}\right)\quad M\xi^{-1}\sqrt{1-c}\lesssim r\lesssim 1,\\
												O_{\xi,c,r}^{\xi,1-c,r}\left( \left(\f{\xi}{(1-c)^{\f13}}\right)^3\right)\cdot O_{\xi,c,r}^{\xi,1-c,r}\left(\left(\f{\xi r^{-\f12}}{(1-c)^{\f12}}\right)^{-N-\f12}\right)\;\;\; 1\ll r\lesssim M^{-\f72}\f{\xi^2}{1-c},\\
												O_{\xi,c,r}^{\xi,1-c,r}( \xi r)\;\;\; M^{-\f72}\f{\xi^2}{1-c}\lesssim r\lesssim M^{-\f12}\xi^{-1},
											\end{array}\right. \end{align}
										and
										\begin{align}\label{7-f+}
											\begin{split}
												&f_+(r,\xi,c)=
												O_{\xi,c,r}^{\xi,1-c,r}\left((\xi r)^{-\f12}\right)e^{i\xi r},\\ &(\mathcal{D}_{l,m}f_+)(r,\xi,c)=
												O_{\xi,c}^{\xi,1-c}\left((\xi r)^{-\f12}\right)e^{i\xi r},\;\;r\gtrsim  M^{-\f12}\xi^{-1}.
											\end{split}
										\end{align}
										\end{itemize}
										
									A direct consequence of the above estimates are, it holds uniformly for $(c,\xi)\in (0,1)\setminus\{V(0)\}\times \mathbb{R}\setminus\{0\}$ and $r>0$ that
									\begin{align} \label{bd:phi/|W|}
										\left|\f{\phi(r,\xi,c)}{W(\xi,c)}\right|+\left|\f{\mathcal{D}\phi(r,\xi,c)}{W(\xi,c)}\right|\leq C,
										\end{align}
									
								\end{proposition}
								
								\begin{proof} The proof proceeds as follows. In the first step, we define $\widetilde{\phi}_0 (r,\xi,c)=\f{\phi (r,\xi,c)}{W_0(\xi,c)}$, where $W_0$ is as in \eqref{behave:Wron0}. We  prove \eqref{1-phi}- \eqref{6-f+} for $\widetilde{\phi}_0$ (instead of $\tilde{\phi}$, which has the identical bounds) and $f_+$. In the second step, we show that
									\begin{align} \label{behaveW3}
										\mathcal{W}(\widetilde{\phi}_0,f_+)(\xi,c)=O_{\xi,c}
										^{\xi,\rho_0(\xi,c)}(1),\quad\text{and}\quad \mathcal{W}(\widetilde{\phi}_0,f_+)(\xi,c)\sim 1, \end{align}
									which gives  \eqref{behaveW1}, since  $\f{W}{W_0}
									=\f{\mathcal{W}(\phi,f_+)}{W_0}=\mathcal{W}(\widetilde{\phi}_0,f_+)$. And \eqref{behaveW2}  follows directly from \eqref{behaveW1}.
								In the last step, it is straightforward to get  by \eqref{behaveW1} that
									\begin{align*}
										\widetilde{\phi}=\widetilde{\phi}_0\cdot \f{W_0}{|W|}
										=\widetilde{\phi}_0\cdot O_{\xi,c}
										^{|\xi|,\rho_0(\xi,c)}(1).\end{align*}
									Then using the bounds obtained for $\widetilde{\phi}_0$ in step 1, we obtain the same desired bounds for $\widetilde{\phi}$. Lastly,  \eqref{bd:phi/|W|} follows directly.

									Notice that in our proof,  $\delta \ll 1$ and $M\gg 1$ are fixed.\smallskip
									
									To prove \eqref{1-phi} (for $\widetilde{ \phi}_0$), it follows from \eqref{serious:phi(1-c)1/2} and Lemma \ref{lem:phi-1.5} that for $r\lesssim \xi^{-1}$,
									\begin{align}\label{summerize-phi0-case1}
										\widetilde{\phi}_0(r,c,\xi)&=\xi r\left(1+O_{\xi,c,r}^{\xi,1-c,r}
										\left(\xi^2r^2+\f{\xi^2}{1-c}(r^2\mathbf{1}_{r\lesssim 1}+ \mathbf{1}_{r\gtrsim 1 })\right)\right),
										\end{align}
									which gives $\widetilde{\phi}(r,\xi,c)=O_{\xi,c,r}^{\xi,1-c,r}(\xi r)$. To derive the estimate for the derivatives in $r$, it follows from \eqref{serious:phi(1-c)1/2}
									and \eqref{serious:phi(1-c)1/2-2} in Lemma \ref{lem:phi-1.5} that for $r\lesssim \xi^{-1}$,
									\begin{align*}
										\widetilde{ \phi}(r,c,\xi)&=\xi r\left(1+O_{\xi,c,r}^{\xi,1-c,1}
										\left(\f{\xi^2}{1-c}\right)\mathbf{1}_{r\lesssim 1}+ O_{\xi,c,r}^{\xi,1-c,r}
										\left(\xi^2r^2+\f{\xi^2}{1-c}\right)\mathbf{1}_{1\lesssim r\lesssim\xi^{-1} }\right)\\
										&=\xi r+\xi rO_{\xi,c,r}^{\xi,1-c,1}
										\left(\f{\xi^2}{1-c}\right)\mathbf{1}_{r\lesssim 1}+ O_{\xi,c,r}^{\xi,1-c,r}
										\left(\xi r\right)\mathbf{1}_{1\lesssim r\lesssim\xi^{-1} },
									\end{align*}
									which gives
									\begin{align*}
										r^{-1}\widetilde{\phi}(r,c,\xi),\quad \pa_r\widetilde{ \phi}(r,c,\xi)&=\xi +\xi O_{\xi,c,r}^{\xi,1-c,1}
										(1)\mathbf{1}_{r\lesssim 1}+ O_{\xi,c,r}^{\xi,1-c,r}
										(\xi)\mathbf{1}_{1\lesssim r\lesssim\xi^{-1} }=O_{\xi,c}^{\xi,1-c}(\xi).
									\end{align*}
									From two formulas above, we deduce that  for $n\geq 1$, $r\lesssim \xi^{-1}$,
									\begin{align*}
										\pa_r^n\widetilde{ \phi}(r,c,\xi)&=
										O_{\xi,c}^{\xi,1-c}(\xi), \;\pa_r^n(\pa_r+r^{-1})\widetilde{ \phi}(r,c,\xi)=
										O_{\xi,c}^{\xi,1-c}(\xi).
										\end{align*}
										
										To prove \eqref{1-f+}, it follows from Lemma \ref{lem:k<1-f+}, $H_{+}(z)=Cz^{-\f12}e^{iz}\left(1+O\left(z^{-1}\right)\right)(z\gtrsim 1)$ and \eqref{estQ-G-0} that
										\begin{align*}
											f_+(r,c,\xi)&=C(\xi r)^{-\f12}e^{i\xi r}\left(1+
											O_{\xi,r}^{\xi,r}\left((\xi r)^{-1}\right)\right)
											\left(1+O_{\xi,c,r}^{\xi,1-c,r}\left(\f{\xi^2}{1-c}\right)\right)\\
											&=  C(\xi r)^{-\f12}e^{i\xi r} O_{\xi,c,r}^{\xi,1-c,r}(1),
											\end{align*}
										which implies the first bound. For the second bound, it suffices to observe that when $r\gtrsim \xi^{-1}$,
										\begin{align*}
										\left|\pa_r\left((\xi r)^{-\f12}\right)\right|\lesssim r^{-1}(\xi r)^{-\f12}\lesssim \xi,\quad r^{-1}\lesssim \xi,\quad \pa_r\left(e^{i\xi r}\right)= i\xi e^{i\xi r}.
										\end{align*}
										
										To prove \eqref{2-phi}, it follows from \eqref{behave:phi2-2} in Lemma \ref{lem: phi-2} that for $r\lesssim \xi^{-\f23}$,
										\begin{align*}
											\widetilde{ \phi}(r,c,\xi)&=\xi^{\f23} r\left(1+O_{\xi,c,r}^{\xi,\xi^{-\f23},\xi^{-\f23}}
											(1)\right).
										\end{align*}
										Then the first line in \eqref{2-phi} follows directly. For the second line in \eqref{2-phi}, we deduce the following stronger bounds:
										for $c\in(0,V(0))$, $\xi\gtrsim 1$, $r\lesssim \xi^{-\f23}$,
										\begin{align*}
											(\xi r)^{-1}\widetilde{\phi}(r,\xi,c)&=
											O_{\xi,c}^{\xi,\xi^{-\f23}}
											\left(\xi^{-\f{1}{3}}\right),\\
											\left(\xi^{-1}\pa_r\right)^m\widetilde{\phi}(r,\xi,c)&=
											\xi^{-m}  m\xi^{\f23}\pa_r^{m-1}\left(O_{\xi,c,r}^{\xi,\xi^{-\f23},\xi^{-\f23}}
											(1)\right)+\xi^{-m}\xi^{\f23} r\pa_r^{m}\left(O_{\xi,c,r}^{\xi,\xi^{-\f23},\xi^{-\f23}}(
											1)\right)\\
											&=
											O_{\xi,c}^{\xi,\xi^{-\f23}}\left(\xi^{-\f{m}{3}}\right).
											\end{align*}
											
											To prove \eqref{2-f+},  we get by Lemma \ref{lem:f+2} that for $r\gtrsim \xi^{-\f23}$,
											\begin{align*} f_+(r,\xi,c)
												&=C\left(\xi r Q^ {\f12}(r,c)\right)^{-\f12}e^{i\xi\int_0^r
													\sqrt{\f{V(s')-c}{1-c}}ds'}\notag
												\left(1+O_{\xi,c,r}^{\xi,V(r)-c,r}\left(\left(\xi r Q^ {\f12}(r,c)\right)^{-\f12}\right)\right).
												\end{align*}
											Noticing that for $r\gtrsim \xi^{-\f23}\gtrsim c-V(0)$, we get by \eqref{behave:Q-V(0)} and \eqref{estx1-3} that
											\begin{align*}
												&\xi r Q(r,c)^{\f12} \sim \xi\int_0^r Q(s',c)^{\f12}ds'\sim \xi\left(\f{r^{\f32}}{\langle r\rangle^{\f12}}+(V(0)-c)^{\f12}r\right)\gtrsim \xi r^{\f23}\gtrsim 1,  \\
												& V(r)-c\gtrsim \f{r}{\langle r\rangle} \gtrsim \xi^{-\f23},
											\end{align*}
											which along with \eqref{estQ-G-0} gives
											\begin{align*}
												\left(\xi r Q^ {\f12}(r,c)\right)^{-\f12}
												=O_{\xi,c,r}^{\xi,V(r)-c,r}\left(\left(\xi\int_0^r
												\sqrt{\f{V(s')-c}{1-c}}ds' \right)^{-\f12}\right)=O_{\xi,c,r}^{\xi,\xi^{-\f23},r}
												\left(\left(\xi\int_0^r
												\sqrt{\f{V(s')-c}{1-c}}ds' \right)^{-\f12}\right).
												\end{align*}
											Therefore, for $r\gtrsim \xi^{-\f23}$, using $V(r)-c\gtrsim \xi^{-\f23}$ again, we get
											\begin{align*}f_+(r,\xi,c)&
												=O_{\xi,c,r}^{\xi,V(r)-c,r}\left(\left(\xi\int_0^r
												\sqrt{\f{V(s')-c}{1-c}}ds'\right)^{-\f12}\right)e^{i\xi\int_0^r
													\sqrt{\f{V(s')-c}{1-c}}ds'}\notag
												\left(1+O_{\xi,c,r}^{\xi,V(r)-c,r}(1)\right)\\
												&=O_{\xi,c,r}^{\xi,\xi^{-\f23},r}\left(\left(\xi\int_0^r
												\sqrt{\f{V(s')-c}{1-c}}ds'\right)^{-\f12}\right)e^{i\xi\int_0^r\sqrt{\f{V(s)-c}{1-c}}ds}.
											\end{align*}
											Then the first bound in \eqref{2-f+} follows. Noticing that when $r\gtrsim \xi^{-\f23}$,
											it follows $(\xi r)^{-1}\lesssim \xi^{-\f13}$ and $\sqrt{\f{V(r)-c}{1-c}}\lesssim 1$ (coming from the oscillation part). Then we get the second bound in \eqref{2-f+}.
											
											To prove the first line in \eqref{4-phi}, it follows from \eqref{behave:phi3}, Lemma \ref{lem:phi-4.5} and $J_1(z)=\f12z\left(1+O_{z}^z(z^2)\right)(|z|\lesssim 1)$ that  for $c\in(0,V(0))$, $r\lesssim \xi^{-1}(V(0)-c)^{-\f12}(\lesssim V(0)-c)$,
											\begin{align*}
												\widetilde{\phi}(r,\xi,c)
												&=(V(0)-c)^{\f12}\xi r\left(1+O_{\xi,c,r}^{\xi,|c-V(0)|,r}\left((V(0)-c)\xi^2 r^2\right)\right)\left(1+O_{\xi,c,r}^{\xi,|c-V(0)|,r}\left(\xi^2r^3\right)\right)\\
												&=(V(0)-c)^{\f12}\xi r\left(1+O_{\xi,c,r}^{\xi,|c-V(0)|,r}\left((V(0)-c)\xi^2 r^2\right)\right),
												\end{align*}
											which gives a stronger bounds for $(l,n)=(1,0), (0,1)$ in the second line of \eqref{4-phi}:
											\begin{align*}
												(\mathcal{D}_{1,0}\widetilde{\phi})(r,\xi,c), (\mathcal{D}_{0,1}\widetilde{\phi})(r,\xi,c)=  (V(0)-c)^{\f12}+
												O_{\xi,c}^{\xi,V(0)-c}\left(\left( (V(0)-c)^{\f12}\xi r\right)^2\right).
												\end{align*}
											If $l+n\geq 2$ in the second line of \eqref{4-phi}, we also have stronger bounds:
											 \begin{align*}
												(\mathcal{D}_{0,n}\widetilde{\phi})(r,\xi,c)=&C_n(1-c)^{\f{1-n}{2}}(V(0)-c)
												^{\f{n}{2}}+(V(0)-c)^{\f{n}{2}}O_{\xi,c}^{\xi,V(0)-c}\left( (V(0)-c)^{\f{1}{2}}
												\xi r \right)\\
												&+(V(0)-c)^{\f{n}{2}}O_{\xi,c}^{\xi,V(0)-c}\left(\xi^{-1} (V(0)-c)^{-\f{3}{2}}\right),\\
												(\mathcal{D}_{1,n}\widetilde{\phi})(r,\xi,c)=&C_n(1-c)^{\f{-n}{2}}(V(0)-c)
												^{\f{n+1}{2}}+(V(0)-c)^{\f{n+1}{2}}O_{\xi,c}^{\xi,V(0)-c}\left( (V(0)-c)^{\f{1}{2}}
												\xi r \right)\\
												&+(V(0)-c)^{\f{n+1}{2}}O_{\xi,c}^{\xi,V(0)-c}\left(\xi^{-1} (V(0)-c)^{-\f{3}{2}}\right),
												\end{align*}
											where $C_n$ are constants only depending on $m$, in particular, $C_n=0$ if $n$ is even. The details are as follows. We first recall for the Bessel function $J_1(\cdot)$ that
											\begin{align}\label{behave:Bessel-derivative}
												J_1^{(m)}(z)=S_m(z):= \left\{
												\begin{array}{l}
													C_m+O_{z}^z(z^2),\;m\geq 1\;\text{odd},\\
													O_z^z(z),\;m\geq 2\;\text{even}. \end{array}\right.
											\end{align}
											To prove the last  line in \eqref{4-phi}, it follows from \eqref{behave:phi3}, \eqref{serious:phi-4.5} in Lemma \ref{lem:phi-4.5} and \eqref{behave:Bessel-derivative} that for $c\in(0,V(0))$, $r\lesssim \xi^{-1}((0)-c)^{-\f12} (\lesssim V(0)-c)$ and $m\geq 1$,
											\begin{align*}
												&(\xi^{-1}\pa_r)^m\widetilde{\phi}(r,\xi,c)\\
												&=2(1-c)^{\f12}(\xi^{-1}\pa_r)^m\left(J_1\left(\xi
												\sqrt{\f{c-V(0)}{1-c}}r\right)\right)\left(1+O_{\xi,c,r}^{\xi,V(0)-c,r}
												\left(\xi^2r^3\right)\right)\\
												&\quad+\sum_{m_1=1}^m(\xi^{-1}\pa_r)^{m-m_1}\left(J_1\left(\xi
												\sqrt{\f{c-V(0)}{1-c}}r\right)\right) \cdot
												(\xi^{-1}\pa_r)^{m_1}\left(O_{\xi,c,r}^{\xi,V(0)-c,\xi^{-1}
													(c-V(0))^{-\f12}}
												\left(\xi^{-1}(c-V(0))^{-\f32}\right)\right)
												\\
												&=2(1-c)^{\f12}\left(\f{c-V(0)}{1-c}\right)^{\f{m}{2}}
												S_m\left(\xi
												\sqrt{\f{c-V(0)}{1-c}}r\right)\left(1+
												O_{\xi,c}^{\xi,V(0)-c}
												\left((c-V(0))\xi^2r^2\right)\right)\\
												&\quad+\sum_{m_1=1}^m\left(\f{c-V(0)}{1-c}\right)^{\f{m-m_1}{2}}
												S_{m_1}\left(\xi
												\sqrt{\f{c-V(0)}{1-c}}r\right) \cdot
												(V(0)-c)^{\f{m_1}{2}}\left(O_{\xi,c}^{\xi,V(0)-c}
												\left(\xi^{-1}(c-V(0))^{-\f32}\right)\right)\\
												&=
												\left\{
												\begin{array}{l}
													(1-c)^{\f12}\left(\f{c-V(0)}{1-c}\right)^{\f{m}{2}}\left(C_m+O_{\xi,c}^{\xi,V(0)-c}
													\left((V(0)-c)\xi^2r^2\right)+O_{\xi,c}^{\xi,V(0)-c}
													\left(\xi^{-1}(c-V(0))^{-\f32}\right)\right),\;m\;\text{odd},\\
													(1-c)^{\f12}\left(\f{c-V(0)}{1-c}\right)^{\f{m}{2}}\left(O_{\xi,c}^{\xi,V(0)-c}
													\left((V(0)-c)^{\f12}\xi r \right)+O_{\xi,c}^{\xi,V(0)-c}
													\left(\xi^{-1}(c-V(0))^{-\f32}\right)\right),\;m\;\text{even}. \end{array}\right.
											\end{align*}
											
											To prove the first line in \eqref{4-f+}, we apply  $(c-V(0))^{-\f12}\lesssim \xi^{\f13}$ and Lemma \ref{lem:f+2} to deduce that for $r\gtrsim  \xi^{-1}|c-V(0)|^{-\f12}$,
											\begin{align*}
												f_+(r,\xi,c)
												&=C\left(\xi rQ(r,c)^{\f12}\right)^{-\f12}e^{i\xi\int_0^r
													\sqrt{\f{V(s)-c}{1-c}}ds}\notag
												\left(1+O_{\xi,c,r}^{\xi,V(r)-c,r}\left(\left(\xi rQ (r,c)^{\f12}\right)^{-1}\right)\right)\\
												&=O_{\xi,c,r}^{\xi,V(r)-c,r}\left(\left(\xi\int_0^r
												\sqrt{\f{V(s)-c}{1-c}}ds\right)^{-\f12}\right)e^{i\xi\int_0^r
													\sqrt{\f{V(s)-c}{1-c}}ds},
											\end{align*}
											where we used  \eqref{behave:Q-V(0)}, \eqref{estQ-G-0} to get
											$Q(r,c)\geq V(0)-c$, $\xi rQ(r,c)^ {\f12}\gtrsim \xi\int_0^r
											\sqrt{\f{V(s)-c}{1-c}}ds$ and  for $\alpha=\f12,1$,
											\begin{align*}
												O_{\xi,c,r}^{\xi,V(r)-c,r}\left(\left(\xi rQ^ {\f12}(r,c)\right)^{-\alpha}\right)&=
												O_{\xi,c,r}^{\xi,V(0)-c,r}\left(\left(\xi rQ^ {\f12}(r,c)\right)^{-\alpha}\right)\\
												&=  O_{\xi,c,r}^{\xi,V(0)-c,r}\left(\left(\xi r(V(0)-c)^ {\f12}\right)^{-\alpha}\right)=O_{\xi,c,r}^{\xi,V(0)-c,r}(1).
											\end{align*}
											 Noticing that $\xi^{-1}r^{-1}\lesssim  (V(0)-c)^{\f12}\lesssim 1$ and $ \sqrt{\f{V(r)-c}{1-c}}\lesssim 1$, the second line in \eqref{4-f+} follows.
											
											 The proof of \eqref{3-phi} follows identically as \eqref{2-phi}, where the last line has a stronger bound for $l+m\geq 1$:
											\begin{align*}
												(\mathcal{D}_{l,m} \widetilde{\phi})(r,\xi,c)&=O_{\xi,c}^{\xi,\xi^{-\f23}}\left(\xi^{-\f{l+n}{3}}\right).
											\end{align*}
											\eqref{3-f+supplement} follows from \eqref{est:f+6(1)} and Lemma \ref{lem:f+4}, where we used $q\gtrsim \langle r_c\rangle^{-\f23}\sim 1$ by \eqref{bd:q}, $\left|-\xi^{\f23}\tau(r,c)=x^{\f23}\right|\lesssim 1$ for $C^{-1}r_c\lesssim \xi^{-\f23}\lesssim r\leq 2r_c$ by \eqref{estx4-6-G}, and $|\mathrm{Oi}(z)|\lesssim 1 (|z|\lesssim 1)$.
											
											To prove \eqref{3-f+}, we get by \eqref{est:f+6(1+)} and Lemma \ref{lem:f+4} that for $r\geq M\xi^{-\f23}\geq  2r_c$,
											\begin{align*}
												f_+(r,\xi,c)
												&=C\left(\xi r Q^ {\f12}(r,c)\right)^{-\f12}e^{i\xi\int_{r_c}^r
													\sqrt{\f{V(s')-c}{1-c}}ds'}
												\left(1+O_{\xi,c,r}^{\xi,\f{r}{\langle r_c\rangle ^{4}},r}\left(\left(\xi r^{\f32}/\langle r\rangle^{\f12}\right)^{-1}\right)\right)\\
												&=O_{\xi,c,r}^{\xi,r,r}\left(\left(\xi\int_{r_c}^r
												\sqrt{\f{V(s')-c}{1-c}}ds'\right)^{-\f12}\right)e^{i\xi\int_{r_c}^r
													\sqrt{\f{V(s')-c}{1-c}}ds'}
												\left(1+O_{\xi,c,r}^{\xi,r,r}\left(\left(\xi r^{\f32}/\langle r\rangle^{\f12}\right)^{-1}\right)\right) \\
												&=O_{\xi,c,r}^{\xi,\xi^{-\f23},r}\left(\left(\xi\int_{r_c}^r
												\sqrt{\f{V(s')-c}{1-c}}ds'\right)^{-\f12}\right)e^{i\xi\int_{r_c}^r
													\sqrt{\f{V(s')-c}{1-c}}ds'}.
											\end{align*}
											Here we have used  \eqref{behave:Q<rc0}, \eqref{estQ-G},  $r-r_c\geq \f12 r$, $r^{\f12}/\langle r\rangle^{\f12}\gtrsim \xi^{-\f13}$ and $\xi^{-\f23}\lesssim r$ to obtain
											\begin{align*}
												\left(\xi r Q^ {\f12}(r,c)\right)^{-\f12}&=O_{\xi,c,r}^{\xi,\f{r-r_c}{\langle r_c\rangle ^{4}},\f{(r_c-r)\langle r\rangle}{\langle r_c\rangle}}\left(\left(\xi\int_{r_c}^r
												\sqrt{\f{V(s')-c}{1-c}}ds'\right)^{-\f12}\right)\\
												&=O_{\xi,c,r}^{\xi,\xi^{-\f23},r}\left(\left(\xi\int_{r_c}^r
												\sqrt{\f{V(s')-c}{1-c}}ds'\right)^{-\f12}\right),
												\end{align*}
											and
											\begin{align*}
												\left(\xi r Q^ {\f12}(r,c)\right)^{-\f12}&=O_{\xi,c,r}^{\xi,\f{r-r_c}{\langle r_c\rangle ^{4}},\f{(r_c-r)\langle r\rangle}{\langle r_c\rangle}}\left(\left(\xi r(r-r_c)^{\f12}/\langle r\rangle^{\f12}\right)^{-\f12}\right)\\
												&=O_{\xi,c,r}^{\xi,r,r}\left(\left(\xi r^{\f32}/\langle r\rangle^{\f12}\right)^{-\f12}\right)
												=O_{\xi,c,r}^{\xi,\xi^{-\f23},r}\left(\left(\xi^{\f23} r\right)^{-\f12}\right).
											\end{align*}
											Then the first line in \eqref{3-f+} follows.
											Noticing that $(\xi r)^{-1}\lesssim \xi^{-\f13}$ and  $\sqrt{\f{V(r)-c}{1-c}}\lesssim 1$, we also obtain the second line in \eqref{3-f+}.
											
											To prove \eqref{6-phi-1}, using \eqref{behave:phi3} in Lemma \ref{lem:phi-3}
											\big(noticing $\xi^{-1}(1-c)^{\f12}|c-V(0)|^{-\f12}\sim \xi^{-1}r_c^{-\f12}\langle r_c\rangle^{-1}$\big)
											and \eqref{est:C1-1}\Big($C_{61}^{-1}=O_{\xi,c}^{\xi,\f{r_c}{\langle r_c\rangle^4}}\Big(r_c^{\f12}\langle
											r_c\rangle^{2}\cdot(\xi r_c^{\f32}/\langle r_c\rangle^{\f12})^{-N}
											\Big)$\Big) in Lemma \ref{lem:f+3}, we infer that for $r\lesssim M^{-\f12}\xi^{-1}r_c^{-\f12}\langle r_c\rangle^{-1}$,
											\begin{align*}
												\widetilde{\phi}(r,\xi,c)
												&=C_{61}^{-1}(\xi,c)\phi(r,\xi,c)\\
												&=
												O_{\xi,c}^{\xi,\f{r_c}{\langle r_c\rangle^4}}\left(\left(\xi r_c^{\f32}/\langle r_c\rangle^{\f12}\right)^{-1}\right)\cdot O_{\xi,r,c}^{\xi,r,\f{r_c}{\langle r_c\rangle^4}}\left(\xi r_c^{\f12}\langle r_c\rangle r\right)\left(1+O_{\xi,c,r}^{\xi,\f{r_c}{\langle r_c\rangle^4},\xi^{-1} r_c^{-\f12} \langle r_c\rangle^{-1}}(1)\right)\\
												&=
												O_{\xi,c}^{\xi,\f{r_c}{\langle r_c\rangle^4}}\left(\left(\xi r_c^{\f32}/\langle r_c\rangle^{\f12}\right)^{-1}\right)\cdot O_{\xi,r,c}^{\xi,\xi^{-1} r_c^{-\f12} \langle r_c\rangle^{-1},\f{r_c}{\langle r_c\rangle^4}}\left(\xi r_c^{\f12}\langle r_c\rangle r\right) .
												\end{align*}
											This gives the first line, and together with $(1-c)^{\f12}\cdot r_c^{\f12}\langle r_c\rangle \sim \f{r_c}{\langle r_c\rangle}\lesssim 1$,  we get the second line of \eqref{6-phi-1}.
											We  recall by \eqref{behave:phi-6-itself} in Lemma \ref{lem:f+3} that for
											$M^{\f12}\xi^{-1}\langle r_c\rangle^{-1}r_c^{-\f12}\lesssim  r
											\leq r_c-C\xi^{-\f23}\langle r_c\rangle^{\f13}$, \begin{align}\label{6-phi-important} \;\;\widetilde{\phi}(r,\xi,c) &=O_{\xi,c,r}^{\xi,
													\f{r_c-r}{\langle r_c\rangle^4},\f{(r_c-r)r}{r_c}}\left(\left(\xi (r_c-r)^{\f12}\langle r_c\rangle  \f{r}{\langle r\rangle^{\f32}}\right)^{-\f12}\right). \end{align}
												Then the first lines in  \eqref{6-phi-2}-\eqref{6-phi-4} follow directly by their ranges, respectively.  Indeed, for $M^{\f12}\xi^{-1}r_c^{-\f12}\langle r_c\rangle^{-1}\lesssim r\leq \min\{\f12,\f{r_c}{2}\}$, it suffices to use $(r_c-r)^{\f12}\sim r_c^{\f12}$ and $\f{1}{\langle r\rangle^{\f32}}\gtrsim 1$; for $\min\{\f12,\f{r_c}{2}\}\leq r\leq \f{r_c}{2}$, the non-trivial case is when $r_c\geq 1$ and $\f12\leq r\leq \f{r_c}{2}$, which gives  $\xi (r_c-r)^{\f12}\langle r_c\rangle  \f{r}{\langle r\rangle^{\f32}}\sim \xi r_c^{\f32}/\langle r\rangle^{\f12}$ ; for $\f{r_c}{2}\leq r\leq r_c- C\xi^{-\f23}\langle r_c\rangle^{\f13}$, it suffices to notice $r\sim r_c$. Noticing the $r-$weight in \eqref{6-phi-important}, to show the second lines in  \eqref{6-phi-2}-\eqref{6-phi-4}, it suffices to verify  for
												$M^{\f12}\xi^{-1}\langle r_c\rangle^{-1}r_c^{-\f12}\leq r \leq r_c-C\xi^{-\f23}\langle r_c\rangle^{\f13}$:
												\begin{align*}
													\f{(1-c)^{\f12}}{\xi r}\lesssim  1\;\;\text{and}\;\;\f{(1-c)^{\f12}r_c}{\xi r(r_c-r)} \lesssim 1.
												\end{align*}
												Noticing that $(1-c)^{\f12}\sim \langle r_c\rangle^{-\f32}$, the first one follows directly. For the second inequality, when $
												M^{\f12}\xi^{-1}\langle r_c\rangle^{-1}r_c^{-\f12}\leq r \leq r_c/2$, it follows that $\f{(1-c)^{\f12}r_c}{\xi r(r_c-r)}\lesssim
												\f{\langle r_c\rangle^{-\f32}r_c}{\langle r_c\rangle^{-1}r_c^{-\f12}r_c}\lesssim 1 $; when
												$r_c/2\leq r
												\leq r_c-C\xi^{-\f23}\langle r_c\rangle^{\f13}$, it follows from \eqref{6-phi-basic} that
												$\f{(1-c)^{\f12}r_c}{\xi r(r_c-r)}\lesssim
												\f{\langle r_c\rangle^{-\f32}r_c}{\xi r_c\cdot\xi^{-\f23}\langle r_c\rangle^{\f13}}\lesssim (\xi r_c^{\f32}/\langle r_c\rangle^{\f12})^{-\f13}\cdot \f{r_c^{\f12}}{\langle r_c\rangle^2}\lesssim 1 $.
												
												To prove \eqref{6-f+a1}, we notice that when $\xi r_c^{\f32}/\langle r_c\rangle^{\f12}\gtrsim M\gg 1$, the range $r_c-C\xi^{-\f23}\langle r_c\rangle^{\f13}\leq r\leq r_c+C\xi^{-\f23}\langle r_c\rangle^{\f13}$  falls in $[\f{r_c}{2},2r_c]$.  We  conclude by \eqref{bd:q} and \eqref{bd:q-Derivative} that
												\begin{align*}
													\left|\left(\f{1}{\langle r_c\rangle^3 }\right)^{l}\langle r_c\rangle ^i\pa_c^l\pa_r^i\left(q^{\f12}(r,c)\right)
													\right|\lesssim q^{\f12}(r,c)\sim \langle r_c\rangle^{-\f13}, \quad \text{for}\;\; r\in[\f{r_c}{2},2r_c],
													\end{align*}
												which along with $\pa_r\tau=q^{\f12}$, $\pa_c\tau=q^{\f12}\f{\int_{r_c}^r\pa_c(Q^{\f12})ds}{Q^{\f12}}$, \eqref{estQ-G} and \eqref{behave:Q<rc0} gives \begin{align*}
													&
													\left|\left(\f{1}{\langle r_c\rangle^{3}}\right)^l
													\langle r_c\rangle^i\pa_c^l\pa_r^i\tau(r,c)\right|\lesssim \langle r_c\rangle^{\f23}\quad \text{for}\;\; r\in \big[\f{r_c}{2},2r_c\big].
												\end{align*}
												Therefore, it follows from   \eqref{est:f+6(2b)}
												\; in Lemma  \ref{lem:f+4}   and $\xi \langle r_c\rangle \gtrsim \xi r_c^{\f32}/\langle r_c\rangle^{\f12}\gtrsim 1$  that
												\begin{align*}
													\notag f_+(r,\xi,c)
													&=C\xi^{-\f13}q^{-\f14}(r,c)r^{-\f12}\mathrm{Oi}\left(-\xi^{\f23}\tau(r,c)\right)
													O_{\xi,c,r}^{\xi,\xi^{-\f23}\langle r_c\rangle^{-\f{11}{3}},\;\xi^{-\f23}\langle r_c\rangle^{\f13}}\big(1\big)\\
													&=O_{\xi,c,r}^{\xi,\xi^{-\f23}\langle r_c\rangle^{-\f{11}{3}},\;\xi^{-\f23}\langle r_c\rangle^{\f13}}\left(\left(\xi r_c^{\f32}/\langle r_c\rangle^{\f12}\right)^{-\f13}\right),
													\end{align*}
												which gives the first line in \eqref{6-f+a1}. Recall that $(1-c)^{\f12}\sim \langle r_c\rangle^{-\f32}$, the second line follows from $(\xi r_c)^{-1}\cdot\langle r_c\rangle^{-\f32}\lesssim 1$ and   $\xi^{-1} \cdot (\xi^{\f23}\langle r_c\rangle^{-\f13})\cdot \langle r_c\rangle^{-\f32}\lesssim 1$.
												
												To prove \eqref{6-f+},  it follows from \eqref{est:f+6(2a)}, $Q(r,c)\sim \f{r-r_c}{\langle r\rangle}$ (by \eqref{behave:Q<rc0}), \eqref{estQ-G}, $\xi\int_{r_c}^r
												\sqrt{\f{V(s)-c}{1-c}}ds\sim \xi (r-r_c)^{\f32}/\langle r\rangle^{\f12}$ (by \eqref{estx4-6-G}) that for $r\geq r_c+C\xi^{-\f23}\langle r_c\rangle^{\f13}\geq r_c$, it holds
												\begin{align} \notag f_+(r,\xi,c)
													&=C\left(\xi r Q^ {\f12}(r,c)\right)^{-\f12}e^{i\xi\int_{r_c}^r
														\sqrt{\f{V(s)-c}{1-c}}ds}\notag
													O_{\xi,c,r}^{\xi,\f{r-r_c}{\langle r_c\rangle^4},r-r_c}(1)
													,\notag\\
													&= O_{\xi,c,r}^{\xi,\f{r-r_c}{\langle r_c\rangle^4},r-r_c}\left(\left(\xi r(\f{r-r_c}{\langle r\rangle})^{\f12}\right)^{-\f12}\right)
													e^{i\xi\int_{r_c}^r
														\sqrt{\f{V(s)-c}{1-c}}ds}\cdot O_{\xi,c,r}^{\xi,\f{r-r_c}{\langle r_c\rangle^4},r-r_c}(1)\notag \\
													&= O_{\xi,c,r}^{\xi,\f{r-r_c}{\langle r_c\rangle^4},r-r_c}\left(\left(\xi^{\f13}(r-r_c)^{\f12}/\langle r\rangle^{\f16}\cdot \xi^{\f23}r/\langle r\rangle^{\f13}\right)^{-\f12}\right)
													e^{i\xi\int_{r_c}^r
														\sqrt{\f{V(s)-c}{1-c}}ds}\notag\\
													&= O_{\xi,c,r}^{\xi,\f{r-r_c}{\langle r_c\rangle^4},r-r_c}\left(\left(\left(\xi\int_{r_c}^r \sqrt{\f{V(s)-c}{1-c}}ds\right)^{\f13}\cdot \xi^{\f23}r_c/\langle r_c\rangle^{\f13}\right)^{-\f12}\right)
													e^{i\xi\int_{r_c}^r
														\sqrt{\f{V(s)-c}{1-c}}ds}\notag\\
													&= O_{\xi,c,r}^{\xi,\f{r-r_c}{\langle r_c\rangle^4},r-r_c}\left(\left(\xi\int_{r_c}^r \sqrt{\f{V(s)-c}{1-c}}ds\right)^{-\f16}\cdot \left(\xi r_c^{\f32}/\langle r_c\rangle^{\f12}\right)^{-\f13}\right)
													e^{i\xi\int_{r_c}^r
														\sqrt{\f{V(s)-c}{1-c}}ds},\label{6-f+-important} \end{align}
												where we used $\f{|r_c-r|\langle r\rangle}{\langle r_c\rangle}\geq r-r_c$ and $\f{r}{\langle r\rangle^{\f13}}\geq \f{r_c}{\langle r_c\rangle^{\f13}}$, for $r\geq r_c$.
												This gives the first line of \eqref{6-f+} . Then the other bounds in \eqref{6-f+} follow from
												\begin{align*}\f{(1-c)^{\f12}}{\xi(r-r_c)}\sim \xi^{-1}(r-r_c)^{-1}\langle r_c\rangle^{-\f32}\lesssim \left(\xi (r-r_c)^{\f32}/\langle r_c\rangle^{\f12}\right)^{-\f23}\cdot \left(\xi r_c^{\f32}/\langle r_c\rangle^{\f12}\right)^{-\f13}\lesssim 1.\end{align*}
												
												To prove \eqref{7-phi}, we use \eqref{serious:phi} and \eqref{est:C72} to get the first bound;
												we use \eqref{Phi8-1},  \eqref{est:C72}, $\f{\cosh y}{\cosh z}\lesssim e^{-(z-y)}(z\geq y\geq 0)$ and $\int_{r}^{r_{*2}}(1-V(s))ds\sim \f{1}{\sqrt{\langle r\rangle}}(r_*\lesssim r\leq \f12r_{*2}$ to get the second and third bound;
												we use \eqref{phi72} to get the fourth bound. The derivative bound follows directly from \eqref{7-phi}. \eqref{7-f+} follows from
												\eqref{Behave:f+1-1} in Lemma \ref{lem:k<1-f+}.\smallskip

												Next, we prove \eqref{behaveW3} for $\xi>0$, while for $\xi<0$, we use $W(\xi)=-\overline{W(-\xi)}$.\smallskip
												
												To prove the first bound in \eqref{behaveW3}, it follows from  \eqref{summerize-phi0-case1} (for $\widetilde{\phi}_0$) and \eqref{1-f+};
												\eqref{2-phi} (for $\widetilde{\phi}_0$) and \eqref{2-f+};
												\eqref{3-phi} (for $\widetilde{\phi}_0$) and \eqref{3-f+};
												\eqref{4-phi} (for $\widetilde{\phi}_0$) and \eqref{4-f+};
												\eqref{6-phi-4} (for $\widetilde{\phi}_0$) and \eqref{6-f+a1};
												\eqref{7-phi} (for $\widetilde{\phi}_0$), \eqref{7-phi-D} and \eqref{7-f+}; that
												\begin{align*}
														&\mathcal{W}\left(\widetilde{\phi}_0,f_+\right)(\xi,c)
														=\left(r\widetilde{\phi}_0'(r)\cdot f_+(r)-rf_+'(r)\cdot\widetilde{\phi}_0(r)\right)
														\\
														&=
														\left\{
														\begin{array}{l}
															|_{r=\xi^{-1}}=O_{\xi,c}^{\xi,1-c}
															\left(1\right)\;\;\; \;\;\quad c\in(0,1),\;\xi\lesssim M^2 (1-c)^{\f12}.\\
															|_{r=\xi^{-\f23}}=O_{\xi,c}^{\xi,\xi^{-\f23}}
															\left(1\right)\;\;\; \quad c\in(0,V(0)),\;1\lesssim \xi\lesssim |c-V(0)|^{-\f32},\\
															|_{r=\xi^{-\f23}}=O_{\xi,c}^{\xi,\xi^{-\f23}}
															\left(1\right)\;\;\; \quad  c\in(V(0),1-\delta),\;1\lesssim \xi\lesssim |c-V(0)|^{-\f32},\\ |_{r=r_c-C\xi^{-\f23}\langle r_c\rangle^{\f13}}=O_{\xi,c}^{\xi,\xi^{-\f23}\langle r_c\rangle^{-\f{11}{3}}}
															\left(r_c\xi^{\f23}\langle r_c\rangle^{-\f13}\cdot\left(\xi^{\f23} r_c/\langle r_c\rangle^{\f13}\right)^{-\f12}\cdot \left(\xi r_c^{\f32}/\langle r_c\rangle^{\f12}\right)^{-\f13}\right)\\
															\quad\quad\quad\quad\quad\quad=O_{\xi,c}^{\xi,\xi^{-\f23}\langle r_c\rangle^{-\f{11}{3}}}
															\left(1\right)\;\;\; \quad \quad c\in(V(0),1),\;\xi\gtrsim M(1-c)^{\f13}|c-V(0)|^{-\f32},\\
															|_{r=\xi^{-1}(V(0)-c)^{-\f12}}=O_{\xi,c}^{\xi,V(0)-c}
															\left(1\right)\;\; \quad  c\in(0,V(0)),\;\xi\gtrsim (V(0)-c)^{-\f32},\\
															|_{r\sim \xi^{-1}}=
															O_{\xi,c}^{\xi,1-c}\left(1\right)\quad  c\in(1-\delta,1),\;M^2(1-c)^{\f12}\lesssim \xi\lesssim M|1-c|^{\f13}.
														\end{array}\right.
													\end{align*}
												
												To prove the second bound in \eqref{behaveW3}, it follows from \eqref{W:2i} that
												\begin{align*}
													&\mathrm{Im}\left(r(f_+'\bar{f}_+)(r,\xi,c)\right)=1.
													\end{align*}
												Then using $\phi$ being real-valued, we have
												\begin{align*}
													\mathrm{Im}\left(\mathcal{W}\left(\widetilde{\phi}_0,f_+\right)(\xi,c)\bar{f}_+(r,\xi,c)\right)=
													\mathrm{Im}\left(\left(r\widetilde{\phi}_0'f_+
													-rf_+'\widetilde{\phi}_0\right)\bar{f}_+\right)= -\mathrm{Im}\left(rf_+'\bar{f}_+\right)\tilde{\phi}_0(r)
													=-\widetilde{\phi}_0(r,\xi,c),
													\end{align*}
												which shows that for any $r\in\mathbb{R}^+$,
												\begin{align*}
													\left|\mathcal{W}\left(\widetilde{\phi}_0,f_+\right)(\xi,c)\right|\geq \f{|\widetilde{\phi}_0(r,\xi,c)|}{|f_+(r,\xi,c)|}.
													\end{align*}
													Therefore, letting $0<\delta\ll 1$ be fixed, it follows from Lemma \ref{lem:phi-1*} and \eqref{1-f+};
												\eqref{behave:phi2} and \eqref{2-f+};
												\eqref{behave:phi2} and \eqref{3-f+};
												\eqref{serious:phi-4.5} and \eqref{4-f+};
												\eqref{behave:phi-6-itself-lower}  and \eqref{6-f+a1};
												\eqref{phi72} and \eqref{7-f+} that
												\begin{align*}
													& \left|\mathcal{W}\left(\widetilde{\phi}_0,f_+\right)(\xi,c)\right|\geq \f{|\widetilde{\phi}_0(r,\xi,c)|}{|f_+(r,\xi,c)|} \\
													&\gtrsim  \left\{
													\begin{array}{l} \f{J_1(\xi r)}{(\xi r)^{-\f12}}\Big|_{r=\delta\xi^{-1}}\sim 1\;\;\; \;\;\quad \; c\in(0,1),\;\xi\lesssim M^2(1-c)^{\f12},\\
														\f{\xi^{\f23}r}{(\xi^{\f23}r)^{-\f12}}
														\Big|_{r=\delta\xi^{-\f23}}\sim  1\;\;\; \quad\; c\in(0,V(0)),\;1\lesssim \xi\lesssim |c-V(0)|^{-\f32},\\
														\f{\xi^{\f23}r}{(\xi^{\f23}r)^{-\f12}}
														\Big|_{r=\delta\xi^{-\f23}}\sim  1\;\;\; \quad  \; c\in(V(0),1-\delta),\;1\lesssim \xi\lesssim |c-V(0)|^{-\f32},\\
														\f{(V(0)-c)^{\f12}\xi r}{((V(0)-c)^{\f12}\xi r)^{-\f12}}
														\Big|_{r=\delta\xi^{-1}(V(0)-c)^{-\f12}}\sim 1\;\; \quad c\in(0,V(0)),\;\xi\gtrsim M(V(0)-c)^{-\f32},\\
														\f{\left(\xi \langle r_c\rangle \f{r}{\langle r\rangle^{\f32}}(r_c-r)^{\f12}\right)^{-\f12}}{\left(\xi r(r_c-r)^{\f12}/\langle r_c\rangle^{\f12}\right)^{-\f12}} \Big|_{r=r_c-C\xi^{-\f23}\langle r_c\rangle^{\f13}\sim r_c}\sim 1\;\;\; \quad  \;c\in(V(0),1),\;\xi\gtrsim M(1-c)^{\f13}|c-V(0)|^{-\f32},\\
														\f{J_1(\xi r)}{(\xi r)^{-\f12}}\Big|_{r=M^{-\f12}\xi^{-1}}\sim  1\quad c\in(1-\delta,1),\;M^2(1-c)^{\f12}\lesssim \xi\lesssim M|1-c|^{\f13}.
													\end{array}\right.
												\end{align*}
												This finishes the proof of \eqref{behaveW3}.
								\end{proof}
							\section{Uniform estimates of the integral kernel: Part I}
							
							In this section, we prove the uniform  estimate of the integral kernel \eqref{op-K} in Proposition \ref{prop:Boundness-K-D},
							whose framework could be applied to the derivative estimate \eqref{op-K-D} of the kernel.
							
							We  will first prove $l'=l=m=0$ in Proposition \ref{prop:Boundness-K-D}, which will be divided into seven parts. The sketch of the proof is as follows. Firstly, we change the coordinate in the $\xi-$ integral; Secondly, we take $\pa_c$ for the functions which have the the oscillation part and the non-oscillation part, and we expect the factor $y'(c)$ that emerges from the oscillation part will help us to take a transform in $c$, while for the  non-oscillation part, some  non-integrable factor $\rho^{-1}(c)$ will emerge, we will use the decomposition \eqref{ineq:trho} to correct it to be integrable in $c$; Finally, we find that the bounds of the integrated functions fit the conditions of Lemma \ref{lem:symbol-chi} and obtain the desired uniform bound. This framework can be applied to the higher derivatives estimates in the next section.
														
							\subsection{Estimates of oscillatory integral}
						
							\begin{lemma}\label{lem: pifi-Linfty}
								Let $\lambda\in\mathbb{R}^n$ be parameter, $y(\lambda)$ be a function. It holds uniformly in $\lambda$ that
								\begin{align}
									\begin{split}\label{pifi-Linfty-1} &\left|p.v.\int_{\mathbb{R}}\f{e^{i\eta y(\lambda)}}{\eta}d\eta\right|
										\leq \pi ,
									\end{split}
								\end{align}
								and
								\begin{align}
									\begin{split}\label{pifi-Linfty} &\left|p.v.\int_{\mathbb{R}}f(\eta,\lambda)\f{e^{i\eta y(\lambda)}}{\eta}d\eta\right|
										\leq \pi \left\|\int_{\mathbb{R}}f(\eta,\lambda)e^{i\eta x}d\eta\right\|_{L^1_x(\mathbb{R})}. \end{split}
								\end{align}
							\end{lemma}
							
							\begin{proof}
								Using
								\begin{align*}
									&p.v.\int_{\mathbb{R}}\f{e^{i\eta y(\lambda)}}{\eta}d\eta
									=\pi\mathrm{sgn}(y(\lambda)),
									\end{align*}
								we get \eqref{pifi-Linfty-1} and
								write the left hand side of \eqref{pifi-Linfty} as
								\begin{align*}
									&p.v.\int_{\mathbb{R}}f(\eta,\lambda)\f{e^{i\eta y(\lambda)}}{\eta}d\eta
									=\int_{\mathbb{R}}\pi\mathrm{sgn}\left(y(\lambda)-x\right)\left(\int_{\mathbb{R}}f(\eta,\lambda)e^{i\eta x}d\eta \right)dx. \end{align*}
								Then the desired inequality \eqref{pifi-Linfty} follows.
							\end{proof}
							
							\begin{lemma}\label{lem:symbol-chi} Let $\lambda\in \mathbb{R}^n$ be a parameter, $A$ be some fixed positive number, $a(\cdot,x,\lambda)$ be a smooth function of $x$ supported in $(-A,A)$ or $\mathbb{R}\setminus[-A,A]$. For $0\leq k\leq 2$, there exist constants $C_k$  independent of $\eta,x,\lambda$ such that
								\begin{align}\label{bd:a}
									\left|(\eta\pa_{\eta})^ka(\eta,x,\lambda)\right|\leq C_k|\eta|^{\alpha}.
									\end{align}
								\begin{itemize}
									\item If $\alpha>0$, $a(\cdot,x,\lambda)\in C_0^2((-A,A))$, then it holds uniformly in $\lambda$ that
									\begin{align}\label{I_A}
										\left\|\int_{\mathbb{R}}a(\eta,x,\lambda)e^{i\eta x}d\eta\right\|_{L^1_x(\mathbb{R})}\lesssim 1.
									\end{align}
									\item If $\alpha<0$, $a(\cdot)\in C_0^2(\mathbb{R}\setminus[-A,A])$, then it holds uniformly in $\lambda$ that
									\begin{align}
										\left\|\int_{\mathbb{R}}a(\eta,x,\lambda)e^{i\eta x}d\eta\right\|_{L^1_x(\mathbb{R})}\lesssim 1\label{I_B}.
										\end{align}
								\end{itemize}
							\end{lemma}
							
							\begin{proof}
								For \eqref{I_A}, it suffices to consider $\alpha\in(0,1)$. Denote $I(x,\lambda)=\int_{\mathbb{R}}a(\eta,x,\lambda)e^{i\eta x}d\eta$.
								Indeed, we will prove $|I(x)|\lesssim 1$ and $|I(x)|\lesssim |x|^{-1-\alpha}$, which gives
								\begin{align*}
									\|I\|_{L^1_x(\mathbb{R})}\lesssim \int_{|x|\leq 1}dx+
									\int_{|x|\geq 1}|x|^{-1-\alpha}dx\lesssim 1.
								\end{align*}
								The former is oblivious using \eqref{bd:a}($k=0$), and the latter follows by integrating by parts twice and using \eqref{bd:a} ($\alpha\in(0,1)$),
								\begin{align}
									\label{x-al-1}
									\left|I(x,\lambda)\right|
									&=\left|\int_{|\eta|\leq |x|^{-1}}a(\eta,x,\lambda)e^{i\eta x}d\eta+\int_{|\eta|\geq |x|^{-1}}a(\eta,x,\lambda)e^{i\eta x}d\eta\right|\\
									\notag &\lesssim \int_{|\eta|\leq |x|^{-1}}|\eta|^{\alpha}d\eta+\left|x^{-1}\int_{|\eta|\geq |x|^{-1}}a(\eta,x,\lambda)\f{de^{i\eta x}}{d\eta}d\eta\right|\\
									\notag &\lesssim |x|^{-\alpha-1}+\left|x^{-1}a(\eta,x,\lambda)|_{|\eta|=|x|^{-1}}\right|
									+\left|x^{-1}\int_{|\eta|\geq |x|^{-1}}a'(\eta,x,\lambda)e^{i\eta x}d\eta\right|\\
									&\lesssim |x|^{-\alpha-1}+\left|x^{-2}a'(\eta,x,\lambda)|_{|\eta|=|x|^{-1}}\right|
									+\left|x^{-2}\int_{|\eta|\geq |x|^{-1}}a''(\eta,x,\lambda)e^{i\eta x}d\eta\right|\lesssim |x|^{-\alpha-1}.  \notag
								\end{align}
								
								For \eqref{I_B}, it suffices to consider $\alpha\in(-1,0)$. Indeed, we will prove $|I(x,\lambda)|\lesssim |x|^{-2}$ and $|I(x,\lambda)|\lesssim |x|^{-1-\alpha}$, which gives
								\begin{align*}
									\|I\|_{L^1_x(\mathbb{R})}\lesssim \int_{|x|\leq 1}|x|^{-1-\alpha}dx+
									\int_{|x|\geq 1}|x|^{-2}dx\lesssim 1.
								\end{align*}
								Using $\alpha<0$ and \eqref{bd:a}, the former follows by integrating by parts twice directly, while the latter follows identically as \eqref{x-al-1}, with $\alpha\in(-1,0)$.
							\end{proof}
							
							In what follows, we use the notation
								$f(x,c,r,s)=m(c,r,s)\mathcal{L}^1_x(\mathbb{R})$ to represent a function $f(x,c,r,s)$ with $\|f/m\|_{L^1_x(\mathbb{R})}\leq C$, where $C$ is a universal constant, independent of $r,s,c$.
								
								\subsection{A decomposition of $\widetilde{\phi}$ by its behavior}
								
								Using $\mathcal{W}(f_+,\bar{f}_+)=2i$, we have an alternative representation of $\widetilde{\phi}$ via $f_+$, used when $r$ large:
								\begin{align} \label{phi-f+-linear}
									\notag\widetilde{\phi}&=
									\f{\mathcal{W}(\widetilde{\phi},\bar{f}_+)}{2i}
									f_+- \f{\mathcal{W}(\widetilde{\phi},f_+)}{2i}
									\bar{f}_+=
									-\mathrm{Re}\left(\mathcal{W}(\widetilde{\phi},\bar{f}_+)f_+\right)\\
									&=  \mathrm{Im}\left(\mathcal{W}(\widetilde{\phi},f_+)\right) \mathrm{Re}(f_+) -
									\mathrm{Re}\left(\mathcal{W}(\widetilde{\phi},f_+)\right) \mathrm{Im}(f_+)\notag\\
									\notag&= \f{\mathrm{Im}(W)(\xi,c)}{|W(\xi,c)|}\mathrm{Re}(f_+)(\xi,c,r) -
									\f{\mathrm{Re}(W)(\xi,c)}{|W(\xi,c)|}\mathrm{Im}(f_+)(\xi,c,r)\\
									&=  \f{i\overline{W}(\xi,c)}{2|W(\xi,c)|}f_+(\xi,c,r)+
									\f{-iW(\xi,c)}{2|W(\xi,c)|}\bar{f}_+(\xi,c,r).
								\end{align}
								
								Recall that the smooth cut-off functions $\chi, \chi_+$ defined in \eqref{def:chi-notation} satisfy
								\begin{align} \label{identity:chi}
									\begin{split}
										\chi_+(\xi)\chi(\xi)&=\chi(\xi);\\ f\left(g_1(\eta)/M\right)f\left(g_2(\eta)\right)&=f\left(g_2(\eta)\right),\;\text{for}\;\;
										M\gg 1,\;g_1(\eta)\sim g_2(\eta),\;f\in\{\chi,\chi_+\}.
								\end{split}\end{align}
							
							\begin{definition} \label{def:decom-phi}
								Let $\delta\ll 1$, $M\gg 1$ be fixed. According to the behavior of $\widetilde{\phi}$  in Proposition \ref{prop:summery},
								we decompose for $\xi\in\mathbb{R}/\{0\}$, $c\in(0,1)/\{V(0)\}$, $r>0$  as
								\begin{small}
									\begin{align}
										\begin{split}
											& \widetilde{\phi}(\xi,c,r)=\sum_{j\in J^L}\widetilde{\phi}_{j},\quad J^L=\{1_0,1_{\infty}\},\;\text{if}\;\;0<|\xi|\lesssim (1-c)^{\f12},\\
											& \widetilde{\phi}(\xi,c,r)=\sum_{j\in J^H_{2}}\widetilde{\phi}_{j},\quad J^H_{2}
											=\{2_0,2_{\infty}\},\;\text{if}\;\;(1-c)^{\f12}\lesssim |\xi|\lesssim  (V(0)-c)^{-\f32},\;c\in(0,V(0)),\\
											& \widetilde{\phi}(\xi,c,r)=\sum_{j\in J^H_{3}}\widetilde{\phi}_{j},\quad J^H_{3}
											=\{3_0,3_{\infty}\},\;\text{if}\;\;|\xi|\gtrsim  (V(0)-c)^{-\f32},\;c\in(0,V(0)),\\
											& \widetilde{\phi}(\xi,c,r)=\sum_{j\in  J^H_{4}}\widetilde{\phi}_{j},\;\;J^H_{4}=\{
											4_0,4_{\infty}\},\;\text{if}\;\;(1-c)^{\f12}\lesssim |\xi|\lesssim (c-V(0))^{-\f32}  ,\;c\in(V(0),1-\delta),\\
											&\widetilde{\phi}(\xi,c,r)=\sum_{j\in J^H_{5}}\widetilde{\phi}_{j},\;\; J^H_{5}=\{5_0,5_1,5_2,5_3,5_{\infty}\},\;
											\text{if}\;\;|\xi|\gtrsim   M (c-V(0))^{-\f32},\;c\in(V(0),1-\delta),\\
											& \widetilde{\phi}(\xi,c,r)=\sum_{j\in  J^H_{6}}\widetilde{\phi}_{j} ,\;J^H_{6}=\{6_0,6_1,6_2,6_3,6_4,6_{\infty}\},\;
											\text{if}\;\;|\xi|\gtrsim   M(1-c)^{\f13},\;c\in(1-\delta,1),\\
											& \widetilde{\phi}(\xi,c,r)=\sum_{j\in  J^H_{7}}\widetilde{\phi}_{j} ,\; J^H_{7}
											=\{7_0,7_1,7_2,7_3, 7_{\infty}\},\;
											\text{if}\;\;M^2(1-c)^{\f12}\lesssim |\xi|\lesssim  M (1-c)^{\f13},\;c\in(1-\delta,1).
											 \end{split}
									\end{align}
								\end{small}
								Here for $\xi<0$, $\widetilde{\phi}_j(\xi)=-\widetilde{\phi}_j(-\xi)$; and for $\xi>0$,
								\begin{align}
										\begin{split}\label{def:tphi-1}
											&\widetilde{\phi}_{1_0}(\xi,c,r)
											=\chi\left(\xi r\right)\widetilde{\phi},\\
											&\widetilde{\phi}_{1_{\infty}}(\xi,c,r)
											=\widetilde{\phi_{1_{\infty}}}(\xi,c,r)e^{i\xi r}+c.c,\\
											&\text{with}\;\;\widetilde{\phi_{1_{\infty}}}(\xi,c,r)=
											\left(1-\chi(\xi r)\right)\f{i\overline{W}(\xi,c)}{2|W(\xi,c)|}\left(f_+e^{-i\xi r}\right);
										\end{split}
									\end{align}
								and
									\begin{align} \label{def:tphi-2}
										\begin{split}
											&\widetilde{\phi}_{2_0}(\xi,c,r)
											=\chi\left(\xi^{\f23} r\right)\widetilde{\phi},\\
											&\widetilde{\phi}_{2_{\infty}}(\xi,c,r)
											= \widetilde{\phi_{2_{\infty}}}(\xi,c,r) e^{i\xi \int_{0}^r\sqrt{\f{V(s')-c}{1-c}}ds'}+c.c, \\
											&\text{with}\;\;\widetilde{\phi_{2_{\infty}}}(\xi,c,r)=\left(1-\chi(\xi^{\f23} r)\right)\f{i\overline{W}(\xi,c)}{2|W(\xi,c)|}\left(f_+e^{-i\xi \int_{0}^r\sqrt{\f{V(s')-c}{1-c}}ds'}\right);
										\end{split}
									\end{align}
								and
									\begin{align}  \label{def:tphi-3}
										\begin{split}
											&\widetilde{\phi}_{3_0}(\xi,c,r)
											=\chi\left(\xi (c-V(0))^{\f12}r\right)\widetilde{\phi},\\
											&\widetilde{\phi}_{3_{\infty}}(\xi,c,r)
											= \widetilde{\phi_{3_{\infty}}}(\xi,c,r) e^{i\xi \int_{0}^r\sqrt{\f{V(s')-c}{1-c}}ds'}+c.c, \\
											&\text{with}\;\; \widetilde{\phi_{3_{\infty}}}(\xi,c,r)
											=\left(1-\chi(\xi (c-V(0))^{\f12} r)\right)\f{i\overline{W}}{2|W|}\left(f_+e^{-i\xi \int_{0}^r\sqrt{\f{V(s')-c}{1-c}}ds'}\right);
										\end{split}
									\end{align}
								and
									\begin{align} \label{def:tphi-4}
										\begin{split}
											&\widetilde{\phi}_{4_0}(\xi,c,r)
											=\chi\left(\xi^{\f23} r\right)\widetilde{\phi}_{4_0},\\
											&\widetilde{\phi}_{4_{\infty}}(\xi,c,r)
											=\widetilde{\phi_{4_{\infty}}}(\xi,c,r) e^{i\xi \int_{r_c}^r\sqrt{\f{V(s')-c}{1-c}}ds'}+c.c, \\
											&\text{with}\;\; \widetilde{\phi_{4_{\infty}}}(\xi,c,r)
											= \left(1-\chi(\xi^{\f23} r)\right)\f{i\overline{W}}{2|W|}\left(f_+e^{-i\xi \int_{r_c}^r\sqrt{\f{V(s')-c}{1-c}}ds'}\right);
									\end{split} \end{align}
								and
								\begin{align}\label{def:tphi-5}
										\begin{split}
											\widetilde{\phi}_{5_0}(\xi,c,r)
											&=
											\chi\left(\xi r_c^{\f12} r/M^{\f12}\right)\widetilde{\phi},\\ \widetilde{\phi}_{5_1}(\xi,c,r)
											&=
											\left(\chi_+(4r/r_c)-\chi\left(\xi r_c^{\f12} r/M^{\f12}\right)\right)\widetilde{\phi},\\
											\phi_{5_2}(\xi,c,r)
											&= \left(1- \chi_+(4r/r_c)\right)\left(1-\chi_+(\xi^{\f23} (r_c-r))\right)\widetilde{\phi},\\
											\widetilde{\phi}_{5_{3}}(\xi,c,r)
											&= \left(1- \chi_+(4r/r_c)\right)\chi_+(\xi^{\f23} (r_c-r))\chi(\xi^{\f23}(r-r_c))\left(\f{i\overline{W}}{2|W|}f_+-
											\f{iW}{2|W|}\bar{f}_+\right),\\
											\widetilde{\phi}_{5_{\infty}}(\xi,c,r)
											&=\widetilde{\phi_{5_{\infty}}}(\xi,c,r) e^{i\xi \int_{r_c}^r\sqrt{\f{V(s')-c}{1-c}}ds'}+c.c,										
											\end{split}
									\end{align}
								where $$\widetilde{\phi _{5_{\infty}}}(\xi,c,r)
											=\f{i\overline{W}}{2|W|}\left(1- \chi_+(4r/r_c)\right)
											\chi_+\left(\xi^{\f23}(r_c-r)\right)
											\left(1-\chi(\xi^{\f23}(r-r_c))\right)\left(f_+e^{-i\xi \int_{r_c}^r\sqrt{\f{V(s')-c}{1-c}}ds'}\right);$$
								and
									\begin{align}\label{def:tphi-6}
										\begin{split}
											\widetilde{\phi}_{6_0}(\xi,c,r)
											&=
											\chi\left(\xi r_c^{\f32} r/M^{\f12}\right)\widetilde{\phi},
											\\ \widetilde{\phi}_{6_1}(\xi,c,r)&=
											\left(\chi_+(4r)-\chi(\xi r_c^{\f32} r/M^{\f12})\right)\widetilde{\phi},\\
											\widetilde{\phi}_{6_2}(\xi,c,r)
											&=
											\left(\chi_+(4r/r_c)-\chi_+(4r)\right)\widetilde{\phi},\\
											\phi_{6_3}(\xi,c,r)
											&= \left(1- \chi_+(4r/r_c)\right)\left(1-\chi_+\left(\xi^{\f23} (r_c-r)/r_c^{\f13}\right)\right)\widetilde{\phi},\\
											\widetilde{\phi}_{6_4}(\xi,c,r)
											&= \left(1- \chi_+(4r/r_c)\right)\chi_+\left(\xi^{\f23} (r_c-r)/r_c^{\f13}\right)\chi(\xi^{\f23}(r-r_c)/r_c^{\f13})\left(\f{i\overline{W}}{2|W|}f_+-
											\f{iW}{2|W|}\bar{f}_+\right),\\
											\widetilde{\phi}_{6_{\infty}}(\xi,c,r)
											&=\widetilde{\phi_{6_{\infty}}}(\xi,c,r) e^{i\xi \int_{r_c}^r\sqrt{\f{V(s')-c}{1-c}}ds'}+c.c,										
											\end{split}
									\end{align}
								where
								$$ \widetilde{\phi _{6_{\infty}}}(\xi,c,r)
											=\f{i\overline{W}}{2|W|}\left(1- \chi_+(4r/r_c)\right)
											\chi_+\left(\xi^{\f23} (r_c-r)/r_c^{\f13}\right)\left(1-\chi\left(\xi^{\f23}(r-r_c)/r_c^{\f13}\right)\right)
											\left(f_+e^{-i\xi \int_{r_c}^r\sqrt{\f{V(s')-c}{1-c}}ds'}\right);$$
											and
									\begin{align}
										\begin{split}\label{def:tphi-7}
											\widetilde{\phi}_{7_0}(\xi,c,r)
											&=\chi\left(\f{\xi r}{M (1-c)^{\f12}}\right)\widetilde{\phi},\\
											\widetilde{\phi}_{7_1}(\xi,c,r)
											&= \left(\chi\left(r/M^{\f12}\right)-\chi\left(\f{\xi r}{M(1-c)^{\f12}}\right)\right)\widetilde{\phi},\\
											\widetilde{\phi}_{7_2}(\xi,c,r)
											&= \left(\chi\left(\f{M^{\f72}(1-c) r}{\xi^2}\right)-\chi(r/M)\right)\widetilde{\phi},\\
											\widetilde{\phi}_{7_3}(\xi,c,r)
											&=
											\left(\chi\left(M^{\f12}\xi r\right)-\chi\left(\f{M^{\f72}(1-c) r}{\xi^2}\right)\right)\widetilde{\phi},\\
											\widetilde{\phi}_{7_{\infty}}(\xi,c,r)
											&= \left(1-\chi(M^{\f12}\xi r)\right)\f{i\overline{W}}{2|W|}f_+(\xi,c,r) +c.c.
											\end{split}
									\end{align}
								\end{definition}
								
					\subsection{The analysis for the weight functions.}
			
								The following lemma shows more precise behavior of the weight defined in Definition \ref{def:W0rhoc}.
								
								\begin{lemma}  \label{lem:rho}
									Let $\delta\ll 1$, $M\gg 1$ be fixed. Let $r_c=V^{-1}(c)$ for $c\in(V(0),1)$, and $\rho(c,\xi,s)$ be defined as follows
									\begin{small}
										\begin{align}\label{def:rho-1}
											\begin{split}
												\rho(c,\xi,s)^{-1}:= \left\{
												\begin{array}{l}
													(1-c)^{-1}\triangleq \rho_{1}^{-1}\;\;\quad\quad \quad\;\;\text{for}\quad c\in(0,1),\;\xi\lesssim (1-c)^{\f12};\\
													\xi^{\f23}\triangleq \rho_{2}^{-1}\;\;\; \quad\quad\quad\quad\quad\quad\text{for}\quad c\in(0,V(0)),\;1\lesssim \xi\lesssim |c-V(0)|^{-\f32};\\
													(V(0)-c)^{-1}\triangleq \rho_{3}^{-1}\;\; \quad\quad \text{for}\quad c\in(0,V(0)),\;\xi\gtrsim |c-V(0)|^{-\f32};\\
													\xi^{\f23}\triangleq \rho_{4}^{-1}\;\;\;\quad\quad\quad\quad\quad\quad\text{for}\quad  c\in(V(0),1-\delta),\;1\lesssim \xi\lesssim |c-V(0)|^{-\f32};\\
													\left\{
													\begin{array}{l}
														r_c^{-1}\triangleq \rho_{5_0}^{-1}\;\;\; \;\;\quad\quad\quad\quad\quad s\lesssim \xi^{-1}r_c^{-\f32}\\
														r_c^{-1}\triangleq \rho_{5_1}^{-1}\;\;\; \;\;\quad\quad\quad\quad\quad \xi^{-1}r_c^{-\f32}\lesssim s\leq r_c/2 \\
														(r_c-s)^{-1}\triangleq \rho_{5_2}^{-1}\;\;\; \;\quad\quad r_c/2\leq s\leq r_c-C\xi^{-\f23} \quad\text{for}\quad c\in(V(0),1-\delta),\;\xi r_c^{\f32}\gtrsim M\\
														\xi^{\f23}\triangleq \rho_{5_3}^{-1}\;\;\; \;\;\;\;\;\quad\quad\quad\quad r_c-C\xi^{-\f23}\leq s\leq r_c+C\xi^{-\f23} \\
														(s-r_c)^{-1}\triangleq \rho_{5_{\infty}}^{-1}\;\;\; \quad\quad\; s\geq r_c+C\xi^{-\f23} .
													\end{array}\right.\\ \left\{
													\begin{array}{l}
														r_c^{3}\triangleq \rho_{6_0}^{-1}\;\;\; \;\;\quad\quad\quad\quad\quad\quad s\lesssim \xi^{-1}r_c^{-\f12}\\
														r_c^{3}\triangleq \rho_{6_1}^{-1}\;\;\; \;\;\quad\quad\quad\quad\quad\quad \xi^{-1}r_c^{-\f12}\lesssim s\leq 1/2 \\
														r_c^{3}\triangleq \rho_{6_2}^{-1}\;\;\; \;\;\quad\quad\quad\quad\quad\quad \f12\leq s\leq r_c/2 \\
														(r_c-s)^{-1}r_c^4\triangleq \rho_{6_3}^{-1}\;\;\; \;\quad\quad\quad r_c/2\leq s\leq r_c-C\xi^{-\f23}r_c^{\f13} \quad\quad\text{for}\quad c\in(1-\delta,1),\;\xi r_c\gtrsim M\\
														\xi^{\f23}r_c^{\f{11}{3}}\triangleq \rho_{6_4}^{-1}\;\;\; \;\;\;\;\quad\quad\quad\quad r_c-C\xi^{-\f23}r_c^{\f13}\leq s\leq r_c+C\xi^{-\f23}r_c^{\f13} \\
														(s-r_c)^{-1}r_c^4\triangleq \rho_{6_{\infty}}^{-1}\;\;\; \quad\quad\quad\; s\geq r_c+C\xi^{-\f23}r_c^{\f13} .
													\end{array}\right.   \\
													(1-c)^{-1}\;\;\;\quad\quad\quad\text{for}\quad c\in(1-\delta,1),\;\delta\ll 1,\;M(1-c)^{\f12}\lesssim \xi\lesssim |1-c|^{\f13}.
												\end{array}\right.
											\end{split}
										\end{align}
									\end{small}
									Then for ${\rho_0}$ in Definition \ref{def:W0rhoc}, it holds ${\rho_0}^{-1}\lesssim \rho^{-1}$. Moreover, there exist a corrector $\sigma(c,s,\xi)$ and a multiplier  $\mathrm{m}(\xi,c,s)$ such that  for the corresponding range, their components admit the following decomposition
									\begin{align} \label{ineq:trho} \rho_j^{-1}=\mathrm{m}_j(\xi,c,s)\cdot \sigma_j(c,s,\xi) ,  \quad j\in\{1,2,...,7_1,...7_{\infty}\},
									\end{align}
									where
								\begin{small}
										\begin{align}\label{def:m-eta}
											&\mathrm{m}(\xi,c,s):= \\\nonumber
											&\left\{
											\begin{array}{l}
												(\f{\xi}{(1-c)^{\f12}})^{\delta/3}(\xi s)^{-\delta/3}\triangleq \mathrm{m}_{1}(\xi,c,s)\;\;\; \quad\text{for}\;\;c\in(0,1),\;\xi\lesssim (1-c)^{\f12},\\
												(\xi(V(0)-c)^{\f32})^{\f23-\delta/9}(\xi s^{\f32})^{\delta/9}\triangleq \mathrm{m}_{2}(\xi,c,s)\;\;\; \quad\text{for}\;\;c\in(0,V(0)),\;1\lesssim \xi\lesssim |c-V(0)|^{-\f32},\\
												(\xi(V(0)-c)^{\f32})^{-\delta/6} (\xi(V(0)-c)^{\f12}s)^{\delta/6}\triangleq \mathrm{m}_{3}(\xi,c,s)\;\;\; \quad\text{for}\;\;c\in(0,V(0)),\;\xi\gtrsim |c-V(0)|^{-\f32},\\
												(\xi(c-V(0))^{\f32})^{\f23-\delta/9}(\xi s^{\f32})^{\delta/9}\triangleq \mathrm{m}_{4}(\xi,c,s)\;\;\; \quad\text{for}\;\;c\in(V(0),1-\delta),\;1\lesssim \xi\lesssim |c-V(0)|^{-\f32},\\
												(1-c)^{-1}\;\;\;\text{for}\;\; c\in(1-\delta,1),\;\delta\ll 1,\;M(1-c)^{\f12}\lesssim \xi\lesssim |1-c|^{\f13}(\xi r_c^{\f32}\gtrsim M),
												\end{array}\right.
									\end{align}\end{small}
									and   for $c\in(V(0),1-\delta),\; \xi r_c^{\f32}\gtrsim M$,
																		\begin{align}\nonumber											=\left\{
												\begin{array}{l}
													(\xi r_c^{\f32})^{-\delta/6}
													(\xi r_c^{\f12}r)^{\delta/6}\textbf{1}_{ s\lesssim \xi^{-1}r_c^{-\f12}}(s)\triangleq \mathrm{m}_{5_1}(\xi,c,s),\\
													(\xi r_c^{\f32})^{-\delta/6}
													(\xi r_c^{\f12}r)^{\delta/6}\textbf{1}_{\xi^{-1}r_c^{-\f12}
														\lesssim s\lesssim r_c/2}(s)\triangleq \mathrm{m}_{5_2}(\xi,c,s),\\
													(\xi (r_c-s)^{\f32})^{-\delta/9}
													(\xi r_c^{\f32})^{\delta/9}\textbf{1}_{r_c/2\leq s\leq r_c-C\xi^{-\f23}}(s)\triangleq \mathrm{m}_{5_3}(\xi,c,s),\\
													(\xi |s-r_c|^{\f32})^{\f23-\delta/9}
													(\xi r_c^{\f32})^{\delta/9}\textbf{1}_{r_c-C\xi^{-\f23}\leq s\leq r_c+C\xi^{-\f23}}(s)\triangleq \mathrm{m}_{5_4}(\xi,c,s),\\
													(\xi (s-r_c)^{\f32})^{-\delta/9}
													(\xi r_c^{\f32})^{\delta/9}\textbf{1}_{ s\geq r_c+C\xi^{-\f23}}(s)
													\triangleq \mathrm{m}_{5_{\infty}}(\xi,c,s),
												\end{array}\right.
									\end{align}	
									and for $c\in(1-\delta,1),\; \xi r_c\gtrsim  M$,
									\begin{align}\nonumber
											=\left\{
												\begin{array}{l}
													(\xi r_c)^{-\delta}(\xi r_c^{\f32}r)^{\delta}\textbf{1}_{s\lesssim  \xi^{-1}r_c^{-\f32}}(s)\triangleq \mathrm{m}_{7_1}(\xi,c,s),\\
													(\xi r_c)^{-\delta}(\xi r_c^{\f32}r)^{\delta}\textbf{1}_{ \xi^{-1}r_c^{-\f32}\lesssim s\leq \f12 }(s)\triangleq \mathrm{m}_{7_2}(\xi,c,s),\\
													(\f{r_c}{r})^{\delta/12}\textbf{1}_{\f12\leq s\leq r_c/2 }(s)\triangleq \mathrm{m}_{7_3}(\xi,c,s),\\
													(\xi( r_c-s)^{\f32}/r_c^{\f12})^{-\delta/9}(\xi r_c)^{\delta/9}\textbf{1}_{r_c/2\leq s\leq r_c-C\xi^{-\f23}r_c^{\f13} }(s)\triangleq \mathrm{m}_{7_4}(\xi,c,s),\\
													(\xi|r_c-s|^{\f32}/r_c^{\f12})^{\f23-\delta/9}(\xi r_c)^{\delta/9}\textbf{1}_{r_c-C\xi^{-\f23}r_c^{\f13}\leq s\leq r_c+C\xi^{-\f23}r_c^{\f13} }(s)\triangleq  \mathrm{m}_{7_5}(\xi,c,s),\\
													(\xi( s-r_c)^{\f32}/r_c^{\f12})^{-\delta/9}(\xi r_c)^{\delta/9}\textbf{1}_{s\geq r_c+C\xi^{-\f23}r_c^{\f13} }(s)\triangleq \mathrm{m}_{7_{\infty}}(\xi,c,s);
												\end{array}\right.
									\end{align}									
		and
						\begin{align}\label{def:rho*} \sigma(c,s,\xi)= \left\{
											\begin{array}{l}
												(1-c)^{-1+\delta/6}s^{\delta/3}\triangleq  \sigma_{1}(c,s),\;\;\; \quad\text{for}\;\;c\in(0,1),\;\xi\lesssim (1-c)^{\f12};\\
												(V(0)-c)^{-1+\delta/6}s^{-\delta/6}\triangleq \sigma_{2}(c,s),\;\;\; \quad\text{for}\;\;c\in(0,V(0)),\;1\lesssim \xi\lesssim |c-V(0)|^{-\f32};\\
												(V(0)-c)^{-1+\delta/6}s^{-\delta/6}\triangleq \sigma_{3}(c,s),\;\;\; \quad\text{for}\;\;c\in(0,V(0)),\;\xi\gtrsim |c-V(0)|^{-\f32};\\
												(c-V(0))^{-1+\delta/6}s^{-\delta/6}\triangleq \sigma_{4}(c,s),\;\;\; \quad\text{for}\;\;c\in(V(0),1-\delta),\;1\lesssim \xi\lesssim |c-V(0)|^{-\f32};\\
												 (1-c)^{-1}\;\;\;\text{for}\;\; c\in(1-\delta,1),\;\delta\ll 1,\;M(1-c)^{\f12}\lesssim \xi\lesssim |1-c|^{\f13}(\xi r_c^{\f32}\gtrsim M);
												 \end{array}\right.
										\end{align}
									and 	for $c\in(V(0),1-\delta),\; \xi r_c^{\f32}\gtrsim M$,
																\begin{align}\nonumber
												=\left\{
												\begin{array}{l}  r_c^{-1+\delta/6}s^{-\delta/6}\textbf{1}_{ s\leq r_c/2 }(s)\triangleq \sigma_{5_1}(c,s,\xi),\\
													r_c^{-1+\delta/6}s^{-\delta/6}\textbf{1}_{s\leq r_c/2}(s)\triangleq \sigma_{5_2}(c,s),\\
													(r_c-s)^{-1+\delta/6}r_c^{-\delta/6}\textbf{1}_{r_c/2\leq s\leq r_c}(s)\triangleq \sigma_{5_3}(c,s),\\
													|s-r_c|^{-1+\delta/6}r_c^{-\delta/6}\textbf{1}_{r_c/2\leq s\leq 2r_c}(s)\triangleq \sigma_{5_4}(c,s),\\
													(s-r_c)^{-1+\delta/6}r_c^{-\delta/6}\textbf{1}_{ s\geq r_c}(s)\triangleq \sigma_{5_{\infty}}(c,s),
												\end{array}\right.
										\end{align}
									and for $c\in(1-\delta,1),\; \xi r_c\gtrsim  M$,									
									\begin{align}\nonumber
												=\left\{
												\begin{array}{l}
													r_c^{3-\delta/2}s^{-\delta}\textbf{1}_{s\lesssim  \xi^{-1}r_c^{-\f32}}(s)\triangleq \sigma_{7_1}(c,s),\\
													r_c^{3-\delta/2}s^{-\delta}\textbf{1}_{ \xi^{-1}r_c^{-\f32}\lesssim s\leq \f12 }(s)\triangleq \sigma_{7_2}(c,s),\\
													r_c^{3-\delta/12}s^{\delta/12}\textbf{1}_{\f12\leq s\leq r_c/2 }(s)\triangleq \sigma_{7_3}(c,s),\\
													(r_c-s)^{-1+\delta/6}r_c^{4-\delta/6} \textbf{1}_{r_c/2\leq s\leq r_c-C\xi^{-\f23}r_c^{\f13} }(s)\triangleq \sigma_{7_4}(c,s),\\
													|s-r_c|^{-1+\delta/6}r_c^{4-\delta/6} \textbf{1}_{r_c-C\xi^{-\f23}r_c^{\f13}\leq s\leq r_c+C\xi^{-\f23}r_c^{\f13} }(s)\triangleq \sigma_{7_5}(c,s),\\
													(s-r_c)^{-1+\delta/6}r_c^{4-\delta/6}\textbf{1}_{s\geq r_c+C\xi^{-\f23}r_c^{\f13} }(s)\triangleq \sigma_{7_{\infty}}(c,s). \end{array}\right.
																							\end{align}							
									For the components of $\sigma(c,s)$, it holds uniformly in $s$ that
									\begin{align}\label{int:sigma}
											\begin{split}
												&\int_{0}^1\sigma_1(c,s)dc\lesssim s^{\delta},\;j=1\\
												& \int_{0}^{V(0)}\sigma_j(c,s)dc\lesssim s^{-\delta},\;\;j=2,3\\
												& \int_{V(0)}^{1-\delta}\sigma_j(c,s,\xi)dc\lesssim s^{-\delta},\;\;j=4,5_1,5_2;\\
												&\int_{V(0)}^{1-\delta}\sigma_j(c,s)dc\lesssim 1,\;\;j=5_3,5_4;5_{\infty};\\
												& \int_{1-\delta}^1\sigma_j(c,s)dc\lesssim s^{-\delta},\;\;j=7_1,7_2;\;\; \int_{1-\delta}^1\sigma_{7_3}(c,s)dc\lesssim s^{\delta};\\
												&\int_{1-\delta}^1\sigma_j(c,s)dc\lesssim 1,\;\;j=7_4,7_5,7_{\infty}.
											\end{split}
										\end{align}
									\end{lemma}
									
								\begin{proof}
									We first decompose $\rho^{-1}$ directly that
									\begin{align}\label{decomp:rho2}
											\begin{split}
												\rho(c,s,\xi)^{-1}= \left\{
												\begin{array}{l}
													1\cdot (1-c)^{-1}\;\;\; \text{for}\;\;c\in(0,1),\;\xi\lesssim (1-c)^{\f12}.\\
													(\xi^{\f23}(V(0)-c))\cdot (V(0)-c)^{-1}\;\;\; \text{for}\;\;c\in(0,V(0)),\;1\lesssim \xi\lesssim |c-V(0)|^{-\f32}.\\
													1\cdot (V(0)-c)^{-1}\;\;\; \text{for}\;\;c\in(0,V(0)),\;\xi\gtrsim |c-V(0)|^{-\f32}.\\
													(\xi^{\f23}(c-V(0)))\cdot (c-V(0))^{-1}\;\;\; \text{for}\;\;c\in(V(0),1-\delta),\;1\lesssim \xi\lesssim |c-V(0)|^{-\f32}.\\
													1\cdot(1-c)^{-1}\;\;\;\text{for}\;\; c\in(1-\delta,1),\;\delta\ll 1,\;M(1-c)^{\f12}\lesssim \xi\lesssim |1-c|^{\f13}(\xi r_c^{\f32}\gtrsim M).												\end{array}\right.
											\end{split}
										\end{align}
										and   for $c\in(V(0),1-\delta),\;\xi |c-V(0)|^{\f32}\sim \xi r_c^{\f32}\gtrsim M$,
																				\begin{align}\nonumber
											\left\{
													\begin{array}{l}
														1\cdot r_c^{-1}\;\;\; r\lesssim  \xi^{-1}r_c^{-\f12},\\
														1\cdot r_c^{-1}\;\;\; \xi^{-1}r_c^{-\f12}\lesssim s\leq r_c/2,\\
														1\cdot (r_c-s)^{-1}\;\;\;\quad\quad \quad\quad  r_c/2\leq s\leq r_c-C\xi^{-\f23},\\
														\xi^{\f23}|s-r_c|\cdot |s-r_c|^{-1}\;\;\;\quad\quad   \quad\quad r_c-C\xi^{-\f23}\leq s\leq r_c+C\xi^{-\f23},\\
														1\cdot
														(r_c-s)^{-1}\;\;\; s\geq r_c+C\xi^{-\f23},
													\end{array}\right.										\end{align}
and for $c\in(1-\delta,1),\;\f{\xi}{(1-c)^{\f13}}\sim \xi r_c\gtrsim M$,
\begin{align}		\nonumber											\left\{
													\begin{array}{l}
														1\cdot r_c^{3}\;\;\; r\lesssim  \xi^{-1}r_c^{-\f32},\\
														1\cdot r_c^{3}\;\;\; \xi^{-1}r_c^{-\f32}\lesssim r\leq 1/2,\\
														1\cdot r_c^{3}\;\;\; \f12\leq r\leq r_c/2,\\
														1\cdot(r_c-s)^{-1}r_c^4\;\;\; r_c/2\leq r\leq r_c-C\xi^{-\f23}r_c^{\f13},\\
														\big(\xi^{\f23}|s-r_c|/r_c^{\f12}\big)\cdot |s-r_c|^{-1}r_c^4\;\;\; r_c-C\xi^{-\f23}r_c^{\f13}\leq r\leq r_c+C\xi^{-\f23}r_c^{\f13},\\
														1\cdot(s-r_c)^{-1}r_c^4\;\;\; r\geq  r_c+C\xi^{-\f23}r_c^{\f13}.
													\end{array}\right.
																						\end{align}
According to the behaviors in \eqref{1-phi}-\eqref{6-f+}, we make a further decomposition as
								\begin{align}\label{decomp:rho}
											\begin{split}
												\rho(c,\xi,s)^{-1}\sim \left\{
												\begin{array}{l}
													\f{\xi^2}{1-c}(\xi r)^{-2}\cdot r^2\;\;\; \text{for}\;\;c\in(0,1),\;\xi\lesssim (1-c)^{\f12},\\
													(\xi^{\f23}r)\cdot r^{-1}\;\;\; \text{for}\;\;c\in(0,V(0)),\;1\lesssim \xi\lesssim |c-V(0)|^{-\f32},\\
													\left(\xi|c-V(0)|^{\f32} \right)^{-1}
													\cdot\left(\xi|c-V(0)|^{\f12} r\right)\cdot r^{-1}\;\;\; \text{for}\;\;c\in(0,V(0)),\;\xi\gtrsim |c-V(0)|^{-\f32},\\
													(\xi^{\f23}r)\cdot r^{-1}\;\;\; \text{for}\;\;c\in(V(0),1-\delta),\;1\lesssim \xi\lesssim |c-V(0)|^{-\f32},\\
													 (1-c)^{-1}\;\;\;\text{for}\;\; c\in(1-\delta,1),\;\delta\ll 1,\;M(1-c)^{\f12}\lesssim \xi\lesssim |1-c|^{\f13}(\xi r_c^{\f32}\gtrsim M),\\
																								\end{array}\right.
										\end{split}
							\end{align}
and for $c\in(V(0),1-\delta),\; \xi r_c^{\f32}\gtrsim M$,								\begin{align}\nonumber
											\left\{
													\begin{array}{l}
														\left(\xi r_c^{\f32} \right)^{-1}\left(\xi r_c^{\f12} r\right)\cdot r^{-1}\;\;\; r\lesssim  \xi^{-1}r_c^{-\f12},\\
														\left(\xi r_c^{\f32} \right)^{-1}\left(\xi r_c^{\f12} r\right)\cdot r^{-1}\;\;\; \xi^{-1}r_c^{-\f12}\lesssim r\leq r_c/2,\\
														\left(\xi^{\f23}(r_c-s)\right)^{-1}(\xi^{\f23}r_c)\cdot r_c^{-1}\;\;\;\quad\quad   r_c/2\leq r\leq r_c-C\xi^{-\f23},\\
														\xi^{\f23}r_c\cdot r_c^{-1}\;\;\;\quad\quad   \quad\quad r_c-C\xi^{-\f23}\leq r\leq r_c+C\xi^{-\f23},\\
														\left(\xi^{\f23}(s-r_c)\right)^{-1}\left(\xi^{\f23}r_c\right)\cdot r_c^{-1}\;\;\; r\geq r_c+C\xi^{-\f23},
													\end{array}\right.							
													\end{align}
and  for $c\in(1-\delta,1),\;\f{\xi}{(1-c)^{\f13}}\sim \xi r_c\gtrsim M$,									
														\begin{align}\nonumber
													\left\{
											\begin{array}{l}
											\left(\xi r_c\right)^{-6}
											\left(\xi r_c^{\f32} r\right)^{6}\cdot r^{-6}\;\;\; r\lesssim  \xi^{-1}r_c^{-\f32},\\
												\left(\xi r_c\right)^{-6}
														\left(\xi r_c^{\f32} r\right)^{6}\cdot r^{-6}\;\;\; \xi^{-1}r_c^{-\f32}\lesssim r\leq 1/2,\\
														\left(\f{r_c}{r}\right)^{\f12}\cdot\left(\f{r}{r_c}\right)^{\f12}r_c^3\;\;\; \f12\leq r\leq r_c/2,\\
														\left(\xi^{\f23}(r_c-s)/r_c^{\f13}\right)^{-1}\left(\xi^{\f23}r_c^{\f23}\right)\cdot r_c^{3}\;\;\; r_c/2\leq r\leq r_c-C\xi^{-\f23}r_c^{\f13},\\
														\left(\xi^{\f23}r_c^{\f23}\right)\cdot r_c^3\;\;\; r_c-C\xi^{-\f23}r_c^{\f13}\leq r\leq r_c+C\xi^{-\f23}r_c^{\f13},\\
														\left(\xi^{\f23}(s-r_c)/r_c^{\f13}\right)^{-1}\left(\xi^{\f23}r_c^{\f23}\right)\cdot r_c^{3}\;\;\; r\geq  r_c+C\xi^{-\f23}r_c^{\f13}.
									\end{array}\right.
																		\end{align}
														Then \eqref{ineq:trho} follows by $ \rho^{-1}=\left(\eqref{decomp:rho2}\right)^{1-\delta/6}\cdot
								\left(\eqref{decomp:rho}\right)^{\delta/6}$. Indeed,
								\begin{align*} \left\{
										\begin{array}{l}
											\int_0^1(1-c)^{-1+\delta/6}s^{\delta/3}dc\lesssim s^{\delta/3},\\
											\int_0^{V(0)}(V(0)-c)^{-1+\delta/6}s^{-\delta/6}dc\lesssim s^{-\delta/6},\\
											\int_0^{V(0)} (V(0)-c)^{-1+\delta/6}s^{-\delta/6}dc \lesssim s^{-\delta/6},\\
											\int_{V(0)}^{1-\delta}(c-V(0))^{-1+\delta/6}s^{-\delta/6}dc
											\lesssim  s^{-\delta/6},\\ (1-c)^{-1}\;\;\;\text{for}\;\; c\in(1-\delta,1),\;\delta\ll 1,\;M(1-c)^{\f12}\lesssim \xi\lesssim |1-c|^{\f13}\left(\xi r_c^{\f32}\gtrsim M\right);\\
											\end{array}\right.
								\end{align*}
							and
							\begin{align*}
											\left\{
											\begin{array}{l}
												\int_0^{V^{-1}(1-\delta)}r_c^{-1+\delta/6}s^{-\delta/6}\textbf{1}_{ s\lesssim \xi^{-1}r_c^{-\f12}}(s)dr_c\\
												\lesssim
												\int_0^{V^{-1}(1-\delta)}r_c^{-1+\delta/6}s^{-\delta/6}\textbf{1}_{ s\lesssim r_c}(s)dr_c\lesssim  \int_s^{V^{-1}(1-\delta)}r_c^{-1+\delta/6}s^{-\delta/6}dr_c\lesssim s^{-\delta/6},\\    \\
												\int_0^{V^{-1}(1-\delta)}
												r_c^{-1+\delta/6}s^{-\delta/6}\textbf{1}_{\xi^{-1}r_c^{-\f12}
													\lesssim s\lesssim r_c/2}(s)dr_c\\
												\lesssim  \int_0^{V^{-1}(1-\delta)}
												r_c^{-1+\delta/6}s^{-\delta/6}dr_c\lesssim s^{-\delta/6},\\  \\
												\int_0^{V^{-1}(1-\delta)}
												(r_c-s)^{-1+\delta/6}r_c^{-\delta/6}\textbf{1}_{r_c/2\leq s\leq r_c-C\xi^{-\f23}(\leq r_c)}(s)dr_c\\
												\lesssim  \int_s^{2s}
												(r_c-s)^{-1+\delta/6}r_c^{-\delta/6}dr_c\lesssim 1,\\  \\
												\int_0^{V^{-1}(1-\delta)}|s-r_c|^{-1+\delta/6}
												r_c^{-\delta/6}\textbf{1}_{(r_c/2\leq )r_c-C\xi^{-\f23}\leq s\leq r_c+C\xi^{-\f23}(\leq 2r_c)}(s)dr_c\\
												\lesssim \int_{s/2}^{2s}|s-r_c|^{-1+\delta/6}
												r_c^{-\delta/6}dr_c\lesssim 1,\\ \\
												\int_0^{V^{-1}(1-\delta)}
												(s-r_c)^{-1+\delta/6}r_c^{-\delta/6}\textbf{1}_{ s\geq r_c+C\xi^{-\f23}(\geq r_c)}(s)dr_c\\
												\lesssim \int_0^{s}
												(s-r_c)^{-1+\delta/6}r_c^{-\delta/6}dr_c\lesssim 1,
											\end{array}\right.
								\end{align*}
							and	
								\begin{align*}
											\left\{
											\begin{array}{l}
												\int_{V^{-1}(1-\delta)}^{+\infty} r_c^{3-\delta/2}s^{-\delta}\textbf{1}_{s\lesssim  \xi^{-1}r_c^{-\f32}}(s)r_c^{-4}dr_c\\
												\lesssim \int_{V^{-1}(1-\delta)}^{+\infty} r_c^{-1-\delta/2}s^{-\delta}\textbf{1}_{s\lesssim  r_c^{-\f12}}(s)dr_c\lesssim \int_{V^{-1}(1-\delta)}^{s^{-2}} r_c^{-1-\delta/2}s^{-\delta}dr_c\lesssim s^{-\delta},\\
												\\
												\int_{V^{-1}(1-\delta)}^{+\infty}
												r_c^{3-\delta/2}s^{-\delta}\textbf{1}_{ \xi^{-1}r_c^{-\f32}\lesssim s\leq \f12 }(s)r_c^{-4}dr_c\\
												\lesssim \int_{V^{-1}(1-\delta)}^{+\infty}
												r_c^{-1-\delta/2}s^{-\delta}dr_c\lesssim  s^{-\delta} ,\\ \\
												\int_{V^{-1}(1-\delta)}^{+\infty}
												r_c^{3-\delta/12}s^{\delta/12}\textbf{1}_{\f12\leq s\leq r_c/2 }(s)r_c^{-4}dr_c\\
												\lesssim   \int_{2s}^{+\infty}
												r_c^{-1-\delta/12}s^{\delta/12}dr_c \lesssim 1,\\ \\
												\int_{V^{-1}(1-\delta)}^{+\infty}
												(r_c-s)^{-1+\delta/6}r_c^{4-\delta/6} \textbf{1}_{r_c/2\leq s\leq r_c-C\xi^{-\f23}r_c^{\f13} (\leq r_c)}(s)r_c^{-4}dr_c\\
												\lesssim   \int_{s}^{2s}
												(r_c-s)^{-1+\delta/6}r_c^{-\delta/6}dr_c\lesssim 1, \\ \\
												\int_{V^{-1}(1-\delta)}^{+\infty}
												|s-r_c|^{-1+\delta/6}r_c^{4-\delta/6} \textbf{1}_{(r_c/2\leq )r_c-C\xi^{-\f23}r_c^{\f13}\leq s\leq r_c+C\xi^{-\f23}r_c^{\f13}(\leq 2r_c) }(s)r_c^{-4}dr_c\\
												\lesssim \int_{s/2}^{2s}|s-r_c|^{-1+\delta/6}r_c^{4-\delta/6} dr_c\lesssim 1,\\       \\
												\int_{V^{-1}(1-\delta)}^{+\infty}
												(s-r_c)^{-1+\delta/6}r_c^{4-\delta/6}\textbf{1}_{s\geq r_c+C\xi^{-\f23}r_c^{\f13} }(s)r_c^{-4}dr_c\\
												\lesssim
												\int_{V^{-1}(1-\delta)}^{s}
												(s-r_c)^{-1+\delta/6}r_c^{-\delta/6}dr_c\lesssim 1, \end{array}\right.
								\end{align*}
								where we used $\int_{V(0)}^{1-\delta}dc\sim \int_0^{V^{-1}(1-\delta)}dr_c$  and   $\int_{1-\delta}^1dc\sim \int_{V^{-1}(1-\delta)}^{+\infty}r_c^{-4}dr_c$.
								\end{proof}
								
								\subsection{Proof of the uniform estimate \eqref{op-K}}
								
								We first have the following observation:\smallskip
								
							 If $\eta$ depends linearly on $\xi$, then we have
									\begin{align} \label{fm:f-tf1}
										\begin{split}
											& (\eta\pa_{\eta})f(\eta,c)=(\xi\pa_{\xi})\tilde{f}(\xi,c),\quad \text{and}\quad \pa_cf(\eta,c)=\pa_c\tilde{f}(\xi,c)
											-\f{\pa_c\eta}{\eta}(\xi\pa_{\xi})\tilde{f}(\xi,c).
										\end{split}
									\end{align}
									Here $\tilde{f}(\xi,c)=f(\eta(\xi),c)$.

									The proof of the uniform estimate \eqref{op-K} is reduced to seven parts:  Propositions \eqref{lem:Boundness-K-pieces1}- Propositions \eqref{lem:Boundness-K-pieces7}.

								\begin{proposition}\label{lem:Boundness-K-pieces1}
									 Let  $\delta\ll 1$ be fixed, $J^{L}=\{1_0,1_{\infty}\}$, and $\widetilde{\phi}_j$($j\in J^L$) be defined as in Definition \ref{def:decom-phi}.
									Then for $(j,j')\in J^L\times J^L$, it holds uniformly for $(r,s,z)\in\mathbb{R}^+\times \mathbb{R}^+\times \mathbb{R}$  that
									\begin{align}
										\int_0^{1}\left|\pa_c\left(p.v.\int_{\mathbb{R}} \left(1-\chi\left(\f{M^2(1-c)^{\f12}}{\xi}\right)\right)
										\widetilde{\phi}_{j}(\xi,c,s)\widetilde{\phi}_{j'}(\xi,c,r)\f{e^{i\xi \sqrt{\f{c}{1-c}}z}}{\xi}d\xi\right)\right|dc\lesssim s^{\delta}+1\label{op-pieces1}.
										\end{align}
								\end{proposition}
								
								\begin{proof}
									 Without further illustration, $c\in(0,1)$, $| \xi| \lesssim (1-c)^{\f12}$.
									Let the new coordinate be defined as
									\begin{align}
										\label{trans:eta1}
										\eta_1(\xi,c,s)=  \xi s.
									\end{align}
									We also recall, by virtue of the weight defined in \eqref{def:rho-1}, that $\rho(c,\xi,s)=1-c\triangleq \rho_{1}$, and by \eqref{def:rho*},
									\begin{align} \label{def:sigma1}
										\sigma_1(c,s)=(1-c)^{-1+\delta/6}s^{\delta/3},
									\end{align}
									which admits the uniform bound in $s$ that
									\begin{align}\label{int:sigma1}
										& \int_0^{\Blue{1}}\sigma_{1}(c,s)dc\lesssim s^{\frac{\delta}{3}}.
										\end{align}
										Under the  coordinate \eqref{trans:eta1},  the multiplier $\mathrm{m}$ in \eqref{def:m-eta} behaves as
										\begin{align}\label{rewriten:m-eta1}
											& O_{\eta_{1_0}}^{\eta_{1}}
											\left(\eta_{1}^{-\delta/3}\right):= \mathrm{m}_{1}(\xi,c,s).
										\end{align}
									By \eqref{ineq:trho},  we can further decompose
										\begin{align}\label{decompose:rho-1}
											\rho_{1}^{-1}&=\mathrm{m}_{1}(\xi,c,s)\sigma_{1}(c,s).
										\end{align}
										We take the transform such that
										\begin{align}
											\begin{split}
												\label{12:1}&\phi_{1_0}(\eta_{1},c,s)
												:=\widetilde{\phi}_{1_0}\left(\eta_{1}s^{-1},c,s\right), \quad \phi_{1_0}^*(\eta_{1},c,s,r):=\widetilde{\phi}_{1_0}
												\left(\eta_{1}s^{-1},c,r\right),   \\
												&\phi_{1_{\infty}}(\eta_{1},c,s)
												:=\widetilde{\phi_{1_{\infty}}}\left(\eta_{1}s^{-1},c,s\right), \quad \phi_{1_{\infty}}^*(\eta_{1},c,s,r):=\widetilde{\phi_{1_{\infty}} }
												\left(\eta_{1}s^{-1},c,r\right),
										\end{split}\end{align}
										with $\widetilde{\phi_{1_{\infty}}}$,  $\widetilde{\phi}_{1_0}$ defined in \eqref{def:tphi-1}.
										We summarize the bounds for \eqref{12:1} as follows.
										
										Using the formula \eqref{fm:f-tf1}, the definitions of $\widetilde{\phi}_{1_0}$, $\widetilde{\phi}_{1_{\infty}}$ in \eqref{def:tphi-1},  the bounds \eqref{1-phi}, \eqref{1-f+} in Proposition \ref{prop:summery}, the decomposition \eqref{decompose:rho-1}  for $\rho_1=1-c$ and the bounds of $\mathrm{m}_{1}$ in \eqref{rewriten:m-eta1}, we obtain
										\begin{small}
											\begin{align}  \label{bd:10}
												\begin{split}
													&\phi_{1_0}(\eta_{1},c,s)=
													\chi(\eta_{1})
													O_{\eta_{1},c}^{\eta_{1},\rho_1}(\eta_{1})
													=\chi(\eta_{1_0})
													O_{\eta_{1},c}^{\eta_{1},(\sigma_{1}\mathrm{m}_{1})^{-1}}(\eta_{1}) =\chi(\eta_{1})
													O_{\eta_{1},c}^{\eta_{1},\sigma_{1}(c,s)^{-1}}\left(\eta_{1}^{1-\delta/3}\right),\\
													&\rho^{-1}_1\phi_{1_0}(\eta_{1},c,s)=
													\sigma_{1}(c,s)\mathrm{m}_{1}(\xi,c,s)\chi(\eta_{1})
													O_{\eta_{1}}^{\eta_{1}}(\eta_{1})=
													\sigma_{1}(c,s)\chi(\eta_{1})
													O_{\eta_{1}}^{\eta_{1}}\left(\eta_{1}^{1-\delta/3}\right), \end{split}
											\end{align}
										\end{small}
										and
										\begin{small}
											\begin{align}  \label{bd:1infty}
												\begin{split}
													&\phi_{1_{\infty}}(\eta_{1},c,s)=
													\left(1-\chi(\eta_{1})\right)
													O_{\eta_{1},c}^{\eta_{1},\rho_1}\left(\eta_{1}^{-\f12}\right)
													=\chi(\eta_{1_0})
													O_{\eta_{1},c}^{\eta_{1},(\sigma_{1}\mathrm{m}_{1})^{-1}}
													\left(\eta_{1}^{-\f12}\right) =\left(1-\chi(\eta_{1})\right) O_{\eta_{1},c}^{\eta_{1},\sigma_{1}(c,s)^{-1}}\left(\eta_{1}^{-\f12}\right),\\
													&\rho^{-1}_1\phi_{1_{\infty}}(\eta_{1},c,s)=\sigma_{1}(c,s)\mathrm{m}_{1}(\xi,c,s)
													\left(1-\chi(\eta_{1})\right)
													O_{\eta_{1}}^{\eta_{1}}\left(\eta_{1}^{-\f12}\right)=
													\sigma_{1}(c,s)\left(1-\chi(\eta_{1})\right)
													O_{\eta_{1}}^{\eta_{1}}\left(\eta_{1}^{-\f12}\right).
													\end{split}
											\end{align}
										\end{small}
										For $\phi^*$, we use  the formula  \eqref{useful:trans-2}
										with $\rho=\rho_1=1-c$, $\f{\pa_c\eta_{1}}{\eta_{1}}=0$, which are independent of $r$. Then, using the bounds \eqref{1-phi}, \eqref{1-f+} in Proposition \ref{prop:summery}, the decomposition \eqref{decompose:rho-1}  for $\rho_1=1-c$ and the bounds of $\mathrm{m}_{1}$ in \eqref{rewriten:m-eta1},
										we  obtain
										\begin{align}&  \label{bi-bd:1infty}
											\begin{split}
												&\phi_{1_0}^*(\eta_{1},c,s,r)= O_{\eta_{1},c}^{\eta_{1},\rho_1}(1),\quad
												\phi_{1_{\infty}}^*(\eta_{1},c,s,r)= O_{\eta_{1},c}^{\eta_{1},\rho_1}(1).
												\end{split}
										\end{align}
										For the cut-off function, by $0<|\xi| \lesssim (1-c)^{1\f12}$,  we infer that
										\begin{align} \label{bd-cutoff:1infty}
											\left(1-\chi\left(\f{(1-c)^{\f12}}{\xi}\right)\right)
											=O_{\eta_{1},c}^{\eta_{1},\rho_1}(1).
										\end{align}
										
										Now we are in a position to prove \eqref{op-pieces1}. Denote $\lambda=(z,s,r)$ as the parameter.
								For $(j,j')=(1_0,1_0)$, We apply the change of variable $\xi\to \eta_{1}=\xi s$ in the $\xi$-integral, thereby reducing \eqref{op-pieces1} to
										\begin{align}
											I_{1_0,1_0}:=\int_0^{1}\left|\pa_c\left(p.v.\int_{\mathbb{R}} a_{1_0,1_0}(\eta_{1},c,\lambda)\f{e^{i\eta_{1} y_{1}(c,\lambda) }}{\eta_{1}}d\eta_{1}\right)\right|dc\lesssim s^{\delta}+1\label{op-pieces1-new-case1},
											\end{align}
										with
										\begin{align}
											\label{def:y1} y_{1}(c,\lambda)= \sqrt{\f{c}{1-c}}z/s\;\;\text{is}\;\;\text{monotonic}\;\;\text{in}\;\;c,
											\end{align}
										and by \eqref{bd:10}, \eqref{bd:1infty} and \eqref{bd-cutoff:1infty},
										\begin{align} \label{def:a_10}
											\begin{split}
												a_{1_0,1_0}(\eta_{1},c,\lambda)&:=
												\left(1-\chi\left(\f{(1-c)^{\f12}}{\xi}\right)\right)\phi_{1_0}(\eta_{1},c,s)
												\phi_{1_0}^* (\eta_{1},c,s,r)\\&=\chi(\eta_{1})
												O_{\eta_{1},c}^{\eta_{1},\sigma_{1}(c,s)^{-1}}\left(\eta_{1}^{1-\delta/3}\right).
												\end{split}
										\end{align}
										Therefore, for the following two parts of $\pa_c$, we apply Lemma \ref{lem: pifi-Linfty} to the first term and  take the transform $c\to y_{1}$ for the second term,
										to obtain
										 \begin{align}
											\begin{split}
												\label{process:integral10}I_{1_0,1_0}&\leq \int_0^{1}\sigma_{1}(c,s)\left|\left(p.v.\int_{\mathbb{R}} \left(\sigma_{1}^{-1}\pa_c\right)a_{1_0,1_0}(\eta_{1},c,\lambda)\f{e^{i\eta_{1} y_{1}(c,\lambda) }}{\eta_{1}}d\eta_{1}\right)\right|dc\\
												&\quad+
												\int_0^{1}\left|\left(p.v.\int_{\mathbb{R}} a_{1_0,1_0}(\eta_{1},c,\lambda)e^{i\eta_{1} y_{1}(c,\lambda) }d\eta_{1}\right)\right|\cdot \left|\f{\pa y_{1}(c,\lambda)}{\pa c}\right|dc\\
												&\lesssim  \int_0^{1}\sigma_{1}(c,s)\left\|\int_{\mathbb{R}} \left(\sigma_{1}^{-1}\pa_c\right)a_{1_0,1_0}(\eta_{1},c,\lambda)e^{i\eta_{1} x}d\eta_{1}\right\|_{L^1_x(\mathbb{R})}dc\\
												&\quad+
												\int_{\mathbb{R}}\left|\int_{\mathbb{R}} a_{1_0,1_0}(\eta_{j},c,\lambda)e^{i\eta_{1} y_{1} }d\eta_{1}\right|d y_{1}.
										\end{split}
									\end{align}
										Then we apply  Lemma \ref{lem:symbol-chi} and the bound in \eqref{def:a_10} to estimate the integral as
										\begin{align*}
											\lesssim  \int_0^{1}\sigma_{1}(c,s)dc+1,
											\end{align*}
										which gives the desired bound \eqref{op-pieces1-new-case1} by \eqref{int:sigma1}.
										
										For $(j,j')=(1_{\infty},1_{0})$, we apply the change of variable  $\xi\to \eta_{1}=\xi s$ in the $\xi$-integral, thereby reducing \eqref{op-pieces1} to
										\begin{align}
											 I_{1_{\infty},1_{0}}:=\int_0^{1}\left|\pa_c\left(p.v.\int_{\mathbb{R}} a_{1_{\infty},1_{0}}(\eta_{1},c,\lambda)\f{e^{i\eta_{1} (y_{1}(c,\lambda) \pm 1)}}{\eta_{1}}d\eta_{1}\right)\right|dc\lesssim s^{\delta}+1\label{op-pieces1-new-case2},
											 \end{align}
										with $y_1$ in \eqref{def:y1}, and by \eqref{bd:10}, \eqref{bd:1infty} and \eqref{bd-cutoff:1infty},
										\begin{align} \label{def:a_1infty}
											\begin{split}
												 a_{1_{\infty},1_{0}}(\eta_{1},c,\lambda)&:=
												\left(1-\chi\left(\f{(1-c)^{\f12}}{\xi}\right)\right)f(\eta_{1},c,s)
												\phi_{1_0}^* (\eta_{1},c,s,r)\\&=\left(1-\chi\left(\f{(1-c)^{\f12}}{\xi}\right)\right)\left(1-\chi(\eta_{1})\right)
												O_{\eta_{1},c}^{\eta_{1},\sigma_{1}(c,s)^{-1}}\left(\eta_{1}^{-\f12}\right)\\
&=\left(1-\chi(\eta_{1})\right)
												O_{\eta_{1},c}^{\eta_{1},\sigma_{1}(c,s)^{-1}}
\left(\eta_{1}^{-\f12}\right),\end{split}
										\end{align}
										where $f\in\{\phi_{1_{\infty}},\bar{\phi}_{1_{\infty}}\}$.
										Therefore, for the following two parts of $\pa_c$, we apply Lemma \ref{lem: pifi-Linfty} to the first term and  take the transform $c\to y_{1}$ for the second term, deducing similarly as in \eqref{process:integral10}. We then apply Lemma \ref{lem:symbol-chi} and the bound given in \eqref{def:a_1infty} to estimate the integral as $\int_0^{1}\sigma_{1}(c,s)dc+1$,  which yields the desired bound \eqref{op-pieces1-new-case2} by virtue of \eqref{int:sigma1}.
																				
										For $(j,j')=(1_{0},1_{\infty})$, we apply the change of variable $\xi\to \eta_{1}=\xi r$ in the $\xi$-integral. The derivation is similar to that for $(j,j')=(1_{0},1_{\infty})$, we leave the detail to the readers.
										
										For $(j,j')=(1_{\infty},1_{\infty})$, we apply the change of variable $\xi\to \eta_{1}=\xi s$ in the $\xi$-integral, thereby  reducing  \eqref{op-pieces1}  to
										\begin{align}
											I_{1_{\infty},1_{\infty}}:=\int_0^{1}\left|\pa_c\left(p.v.\int_{\mathbb{R}} a_{1_{\infty},1_{\infty}}(\eta_{1},c,\lambda)\f{e^{i\eta_{1} (y_{1}(c,\lambda) \pm 1)}}{\eta_{1}}d\eta_{1}\right)\right|dc\lesssim s^{\delta}+1\label{op-pieces1-new-case3},
											\end{align}
										with $y_1$ in \eqref{def:y1}, and by \eqref{bd:10}, \eqref{bd:1infty} and \eqref{bd-cutoff:1infty},
										\begin{align} \label{def:a_1inftyinfty}
											\begin{split}
												 a_{1_{\infty},1_{\infty}}(\eta_{1},c,\lambda)&:=
												\left(1-\chi\left(\f{(1-c)^{\f12}}{\xi}\right)\right)f(\eta_{1},c,s)
												f^* (\eta_{1},c,s,r)\\&=\left(1-\chi(\eta_{1})\right)
												O_{\eta_{1},c}^{\eta_{1},\sigma_{1}(c,s)^{-1}}\left(\eta_{1}^{-\f12}\right),
												\end{split}
										\end{align}
										where $f\in\{\phi_{1_{\infty}},\bar{\phi}_{1_{\infty}}\}$, $f^*\in\{\phi_{1_{\infty}}^*,\bar{\phi}_{1_{\infty}}^*\}$ .
										Therefore, for the two parts of $\pa_c$, we apply Lemma \ref{lem: pifi-Linfty} to the first term and take the transform $c\to y_{1}$ for the second term,
										deducing similarly as in \eqref{process:integral10}. We then apply Lemma \ref{lem:symbol-chi} and the bound in \eqref{def:a_1inftyinfty} to estimate the integral as   $\int_0^{1}\sigma_{1}(c,s)dc+1$, which gives the desired bound \eqref{op-pieces1-new-case2} by \eqref{int:sigma1}.
								\end{proof}
								
								\begin{proposition}\label{lem:Boundness-K-pieces2}
									Let   $\phi_j$ be defined in Definition \ref{def:decom-phi}, $j\in J^H_{2}=\{2_0,2_{\infty}\} $.
									Then for $(j,j')\in J^H_{2}\times J^H_{2}$, it holds uniformly for $(r,s,z)\in\mathbb{R}^+\times \mathbb{R}^+\times \mathbb{R}$  that
									\begin{align}\label{op-pieces2}
																			&\int_0^{V(0)}\left|\pa_c\left(p.v.\int_{\mathbb{R}} \chi\left(\f{M^2(1-c)^{\f12}}{\xi}\right) \chi\left(\xi (V(0)-c)^{\f32}\right)
										\widetilde{\phi}_{j}(\xi,c,s)\widetilde{\phi}_{j'}(\xi,c,r)\f{e^{i\xi z\sqrt{\f{c}{1-c}}}}{\xi}d\xi\right)\right|dc\\
										&\lesssim s^{-\delta}+1.\nonumber
										\end{align}
								\end{proposition}
								
								\begin{proof}
									Without further illustration, we assume $c\in(0,V(0))$, $1\lesssim \xi \lesssim (V(0)-c)^{-\f32}$.
									Let the new coordinate be defined as follows
									\begin{align}
										\label{trans:eta2}\begin{split}
											\eta(\xi,c,s)= \left\{
											\begin{array}{l} \xi s^{\f32} := \eta_{2_0}\;\;\;\quad\quad\quad\quad\quad\quad s\lesssim \xi^{-\f23} ,\\
												\xi \int_{0}^s\sqrt{\f{V(s')-c}{1-c}}ds':= \eta_{2_{\infty}} \;\;\;\quad s\gtrsim \xi^{-\f23}, \end{array}\right.
										\end{split}
									\end{align}
									which satisfies (by \eqref{estx1-3})
									\begin{align} \label{lower:eta_2infty}
										\eta_{2\infty}=\xi \int_{0}^s\sqrt{\f{V(s')-c}{1-c}}ds'\sim \xi \left(\f{s^{\f32}}{\langle s\rangle^{\f12}}+(V(0)-c)^{\f12}s\right)\gtrsim \xi^{\f23}s=\eta_{2_0}^{\f23},\quad \text{if}\quad s\gtrsim \xi^{-\f23},
									\end{align}
									and \begin{align} \label{est:paceta/eta2}
										\f{\pa_c\eta_{2_0}}{\eta_{2_0}}=0,\;\; \quad   \left|\f{\pa_c\eta_{2_{\infty}}}{\eta_{2_{\infty}}}\right|\lesssim \xi^{\f23}=\rho_{2}^{-1},
									\end{align}
									by using the bounds(by  \eqref{behave:Q-V(0)}) \begin{align} \label{lower:v-c}
										V(s)-c\sim  \f{s}{\langle s\rangle}+(V(0)-c)\geq  \f{s}{\langle s\rangle}\gtrsim  \xi^{-\f23},\quad \text{if}\quad s\gtrsim \xi^{-\f23},
									\end{align}
									where the weight  defined  in \eqref{def:rho-1} writes  $\rho(c,\xi,s)=\xi^{-\f23}:= \rho_{2}$. The corrector in \eqref{def:rho*} takes the form
									\begin{align} \label{def:sigma2}
										\sigma_2(c,s)= \left\{
										\begin{array}{l} (V(0)-c)^{-1+\delta/6}s^{-\delta/6}\textbf{1}_{(V(0)-c)\geq 0 }(c):= \sigma_{2_0}(c,s) ,\\
											(V(0)-c)^{-1+\delta/6}s^{-\delta/6}\textbf{1}_{s\gtrsim (V(0)-c)\geq 0 }(c):= \sigma_{2_{\infty}}(c,s), \end{array}\right.
									\end{align}
									which admits the uniform bounds in $s$ that
									\begin{align}\label{int:sigma2}
										& \int_0^{V(0)}\sigma_{2_0}(c,s)dc\lesssim s^{-\delta},\\
										&\int_0^{V(0)}\sigma_{2_{\infty}}(c,s)dc\lesssim s^{-\delta/6}\int_{V(0)-Cs}^{V(0)}(V(0)-c)^{-1+\delta/6}dc\lesssim 1.  \label{int:sigmainfty}
										\end{align}
									Under the coordinates \eqref{trans:eta2},  by \eqref{lower:eta_2infty},  the multiplier $\mathrm{m}$ in  \eqref{def:m-eta} behaves as
									\begin{align}\label{rewriten:m-eta2}
										&\quad\mathrm{m}_2(\xi,c,s)=
										\left\{
										\begin{array}{l}
											O_{\eta_{2_0}}^{\eta_{2_0}}
											\left(\eta_{2_0}^{\delta/9}\right):= \mathrm{m}_{2_0}(\xi,c,s),\\ O_{\eta_{2_{\infty}}}^{\eta_{2_{\infty}}}
											\left(\eta_{2_{\infty}}^{\delta/6}\right):= \mathrm{m}_{2_{\infty}}(\xi,c,s). \end{array}\right.
									\end{align}
									Finally,  by \eqref{ineq:trho},  we can further decompose
									\begin{align}
										\begin{split}\label{decompose:rho-2}
											\rho_{2}^{-1}&=\mathrm{m}_{2_0}(\xi,c,s)\sigma_{2_0}(c,s),\;\text{if}\;\;s\lesssim \xi^{-\f23},\\
											\rho_{2}^{-1}&=\mathrm{m}_{2_{\infty}}(\xi,c,s)\sigma_{2_{\infty}}(c,s),\;\text{if}\;\;s\gtrsim \xi^{-\f23}.
										\end{split}
									\end{align}
									
									\textit{Case 1. }  $(j,j')= (2_0, 2_0)$.
									
									We first apply the change of variable $\xi\to \eta_{2_0}$ in the $\xi$-integral, reducing  \eqref{op-pieces2}  to
									\begin{align}
										I_{2_0,2_0}:=\int_0^{V(0)}\left|\pa_c\left(p.v.\int_{\mathbb{R}} a_{2_0,2_0}(\eta_{2_0},c,\lambda)\f{e^{i\eta_{2_0} y_{2_0}(c,\lambda) }}{\eta_{2_0}}d\eta_{2_0}\right)\right|dc\lesssim s^{-\delta}+1\label{op-pieces2-new-case1}.
									\end{align}
									Here $\lambda=(z,s,r)$ is the parameter. The new functions are defined as follows
									\begin{align}
										\label{def:y20} y_{2_0}(c,\lambda)= \sqrt{\f{c}{1-c}}z/s^{\f32}\;\;\text{is}\;\;\text{monotonic}\;\;\text{in}\;\;c,
									\end{align}
									\begin{align} \label{def:a_20}
										\begin{split}
											&\quad a_{2_0,2_0}(\eta_{2_0},c,\lambda):=\chi\left(\f{(1-c)^{\f12}}{\xi}\right) \chi\left(\xi(V(0)-c)^{\f32}\right)\phi_{2_0}(\eta_{2_0},c,s)\phi_{2_0}^* (\eta_{2_0},c,s,r),
										\end{split}
									\end{align}
									where $\xi(\eta_{2_0},s)=\eta_{2_0}s^{-\f32}$ and
									\begin{align}
										\label{22:1}&\phi_{2_0}(\eta_{2_0},c,s)
										:=\widetilde{\phi}_{2_0}\left(\xi(\eta_{2_0},s),s,c\right), \quad \phi_{2_0}^*(\eta_{2_0},c,s,r):=\widetilde{\phi}_{2_0}\left(\xi(\eta_{2_0},s),r,c\right),
									\end{align}
									with $\widetilde{\phi}_{2_0}$ defined in \eqref{def:tphi-2}.
									Using the formula \eqref{fm:f-tf1}, the definitions of $\widetilde{\phi}_{2_0}$ in \eqref{def:tphi-2},  the bounds \eqref{2-phi}, \eqref{2-f+} in Proposition \ref{prop:summery}, the decomposition \eqref{decompose:rho-2}  for $\rho_2=\xi^{-\f23}$ and the bounds of $\mathrm{m}_{2_0}$ in \eqref{rewriten:m-eta2}, we obtain
									\begin{small}
										\begin{align}  \label{bd:20}
											\begin{split}
												&\phi_{2_0}(\eta_{2_0},c,s)=
												\chi\left(\eta_{2_0}^{\f23}\right)
												O_{\eta_{2_0},c}^{\eta_{2_0},\rho_2}\left(\eta_{2_0}^{\f23}\right)
												=\chi\left(\eta_{2_0}^{\f23}\right)
												O_{\eta_{2_0},c}^{\eta_{2_0},(\sigma_{2_0}\mathrm{m}_{2_0})^{-1}}
												\left(\eta_{2_0}^{\f23}\right)
												=\chi\left(\eta_{2_0}^{\f23}\right)
												O_{\eta_{2_0},c}^{\eta_{2_0},\sigma_{2_0}^{-1}}\left(\eta_{2_0}^{\f23}\right),\\
												&\rho^{-1}_2\phi_{2_0}(\eta_{2_0},c,s)=
												\sigma_{2_0}\mathrm{m}_{2_0}\chi\left(\eta_{2_0}^{\f23}\right)
												O_{\eta_{2_0}}^{\eta_{2_0}}\left(\eta_{2_0}^{\f23}\right)
												=\sigma_{2_0}(c,s)\chi\left(\eta_{2_0}^{\f23}\right)
												O_{\eta_{2_0}}^{\eta_{2_0}}\left(\eta_{2_0}^{\f23}\right).
											\end{split}
										\end{align}
									\end{small}
									For $\phi^*$, we have the formula
									\begin{small} \begin{align}  \notag
											\pa_c\phi^*(\eta,c,s,r)&=
											\pa_c\widetilde{\phi}^*(\xi,c,r)
											-\f{\pa_c\eta(c,s,\xi)}{\eta(c,s,\xi)} (\xi\pa_{\xi})\widetilde{\phi}^*(\xi,c,r)\\
											&=
											O\left(\rho(c,r,\xi)^{-1}\right)(\rho\pa_c)\widetilde{\phi}^*(\xi,c,r)
											+O\left(\f{\pa_c\eta(c,s,\xi)}{\eta(c,s,\xi)}\right) (\xi\pa_{\xi})\widetilde{\phi}^*(\xi,c,r).\label{useful:trans-2}
										\end{align}
									\end{small}
									Note in this case, $\rho=\rho_2=\xi^{-\f23}$, $\f{\pa_c\eta_{2_0}}{\eta_{2_0}}=0$, which along with the bound \eqref{2-phi} in Proposition \ref{prop:summery} gives
									\begin{align}  \label{bi-bd:20}
										\begin{split}
											&\phi_{2_0}^*(\eta_{2_0},c,s,r)= O_{\eta_{2_0},c}^{\eta_{2_0},\rho_2}(1). \end{split}
									\end{align}
									For the cut-off function, by $c\in(0,V(0))$, $1\lesssim |\xi| \lesssim (V(0)-c)^{-\f32}$, $\eta_{2_0}=\xi s^{\f32}$,  we have
									\begin{align} \label{bd-cutoff:20}
										\chi\left(\f{(1-c)^{\f12}}{\xi}\right) \chi\left(\xi (V(0)-c)^{\f32}\right)=O_{\eta_{2_0},c}^{\eta_{2_0},\rho_2}(1).
									\end{align}
									
									Summing up  \eqref{bd:20}, \eqref{bi-bd:20} and  \eqref{bd-cutoff:20}, we obtain
									\begin{align}  \label{bd:a20}
										\begin{split}
											a_{2_0,2_0}(\eta_{2_0},c,\lambda)=\chi\left(\eta_{2_0}^{\f23}\right)
											O_{\eta_{2_0},c}^{\eta_{2_0},\sigma_{2_0}(c,s)^{-1}}\left(\eta_{2_0}^{\f23}\right). \end{split}
									\end{align}
									
									 Now we are in a position to deal with \eqref{op-pieces2-new-case1}.
									We take $\pa_c$ on  the non-oscillation part and the oscillation part. Applying Lemma \ref{lem: pifi-Linfty} to the first term and  taking the transform $c\to y_{2_0}$ for the second term, we infer that  for $(j,j')=(2_0,2_0)$,
									\begin{align}
										\begin{split}
											\label{process:integral20}I_{j,j'}&\leq \int_0^{V(0)}\sigma_{j}(c,s)\left|\left(p.v.\int_{\mathbb{R}} \left(\sigma_{j}^{-1}\pa_c\right)a_{j,j'}(\eta_{j},c,\lambda)\f{e^{i\eta_{j} y_{j}(c,\lambda) }}{\eta_{j}}d\eta_{j}\right)\right|dc\\
											&\quad+
											\int_0^{V(0)}\left|\left(p.v.\int_{\mathbb{R}} a_{j,j'}(\eta_{j},c,\lambda)e^{i\eta_{j} y_{j}(c,\lambda) }d\eta_{j}\right)\right|\cdot \left |\f{\pa y_{j}(c,\lambda)}{\pa c}\right|dc\\
											&\lesssim  \int_0^{V(0)}\sigma_{j}(c,s)\left\|\int_{\mathbb{R}} (\sigma_{j}^{-1}\pa_c)a_{j,j'}(\eta_{j},c,\lambda)e^{i\eta_{j} x}d\eta_{j}\right\|_{L^1_x(\mathbb{R})}dc\\
											&\quad+
											\int_{\mathbb{R}}\left|\int_{\mathbb{R}} a_{j,j'}(\eta_{j},c,\lambda)e^{i\eta_{j} y_{j} }d\eta_{j}\right|d y_{j}.
									\end{split}\end{align}
									Then we apply  Lemma \ref{lem:symbol-chi} and the bound \eqref{bd:a20} to bound the integral as  \begin{align*} \lesssim  \int_0^{V(0)}\sigma_{2_0}(c,s)dc+1, \end{align*}
									which gives the desired bound \eqref{op-pieces2-new-case1} by \eqref{int:sigma2}.
									
									\textit{Case 2. }  $(j,j')=\{2_{\infty}, 2_{0}\}$ or $(j,j')=\{2_0, 2_{\infty}\}$.
									We only prove the first one. For the second one, we change the variable $s$ and $r$, $j$ and $j'$ to reduce it to the first proof.
									For $(j,j')=(2_{\infty},2_0)$, we  apply the change of variable $\xi\to \eta_{2_\infty}$ in the $\xi$-integral, reducing \eqref{op-pieces2} to
									\begin{align}
										I_{2_{\infty},2_0}:=\int_0^{V(0)}\left|\pa_c\left(p.v.\int_{\mathbb{R}} a_{2_{\infty},2_0}(\eta_{2_{\infty}},c,\lambda)\f{e^{i\eta_{2_{\infty}} \left(y_{2_{\infty}}(c,\lambda) \pm 1\right) }}{\eta_{2_{\infty}}}d\eta_{2_{\infty}}\right)\right|dc\lesssim 1\label{op-pieces2-new-case2}.
										\end{align}
									Here, the new coordinate $\eta_{2_{\infty}}$ takes $\xi \int_{0}^s\sqrt{\f{V(s')-c}{1-c}}ds'$, and  $\lambda=(z,s,r)$ is the parameter. The new functions are defined as follows
									\begin{align} \label{def:y2infty} y_{2_{\infty}}(c,\lambda)= \sqrt{\f{c}{1-c}}z/\int_{0}^s\sqrt{\f{V(s')-c}{1-c}}ds'\;\;\text{is}\;\;\text{monotonic}\;\;\text{in}\;\;c.
									\end{align}
									\begin{align} \label{def:a_infty20}
										\begin{split}
											&\quad a_{2_{\infty},2_0}(\eta_{2_0},c,\lambda):=\chi\left(\f{(1-c)^{\f12}}{\xi}\right)
											\chi\left(\xi(V(0)-c)^{\f32}\right)f(\eta_{2_{\infty}},c,s) \phi_{2_0}^* (\eta_{2_{\infty}},c,s,r),
										\end{split}
									\end{align}
									where $\xi(\eta_{2_{\infty}},c,s)=\eta_{2_{\infty}}
									\left(\int_{0}^s\sqrt{\f{V(s')-c}{1-c}}ds'\right)^{-1}$, $f\in\{\phi_{2_{\infty}}, \bar{\phi}_{2_{\infty}}\}$, and
									\begin{align}
										\label{22:2}&\phi_{2_{\infty}}(\eta_{2_{\infty}},c,s)
										:=\widetilde{\phi_{2_{\infty}}}\left(\xi(\eta_{2_{\infty}},c,s),s,c\right), \quad \phi_{2_0}^*(\eta_{2_{\infty}},c,s,r):=\widetilde{\phi}_{2_0}
										\left(\xi(\eta_{2_{\infty}},c,s),r,c\right),
										\end{align}
									with $\widetilde{\phi_{2_{\infty}}}$,$\widetilde{\phi}_{2_0}$ defined in \eqref{def:tphi-2}.
									Noticing by \eqref{lower:eta_2infty}, we have  for some fixed $M\gg 1$,
									\begin{align*}
										1-\chi(\xi^{\f23} s)=\left(1-\chi(M\eta_{2_{\infty}})\right)\left(1-\chi(\xi^{\f23} s)\right).
									\end{align*}
									Therefore, using the definition of $\widetilde{\phi_{2_{\infty}}}$ in \eqref{def:tphi-2}, the formula \eqref{fm:f-tf1},  the bound \eqref{2-f+} in Proposition \ref{prop:summery}, the decomposition \eqref{decompose:rho-2} for $\rho_2=\xi^{-\f23}$ and the bound of $\mathrm{m}_{2_{\infty}}$ in \eqref{rewriten:m-eta2}, we obtain
									\begin{small}
										\begin{align}  \label{bd:2infty}
											\begin{split}
												\phi_{2_{\infty}}(\eta_{2_{\infty}},c,s)&=
												\left(1-\chi\left(\xi^{\f23} s\right)\right)
												O_{\eta_{2_{\infty}},c}^{\eta_{2_{\infty}},\rho_2}\left(\eta_{2_{\infty}}
												^{-\f12}\right)=
												\left(1-\chi(M\eta_{2_{\infty}})\right)\left(1-\chi\left(\xi^{\f23} s\right)\right)
												O_{\eta_{2_{\infty}},c}^{\eta_{2_{\infty}},
													(\sigma_{2_{\infty}}\mathrm{m}_{2_{\infty}})^{-1}}\left(\eta_{2_{\infty}}
												^{-\f12}\right)\\
												&=\left(1-\chi(M\eta_{2_{\infty}})\right)
												O_{\eta_{2_{\infty}},c}^{\eta_{2_{\infty}},
													(\sigma_{2_{\infty}}\mathrm{m}_{2_{\infty}})^{-1}}\left(\eta_{2_{\infty}}
												^{-\f12}\right)
												=\left(1-\chi(M\eta_{2_{\infty}})\right)
												O_{\eta_{2_{\infty}},c}^{\eta_{2_{\infty}},
													\sigma_{2_{\infty}}(c,s)^{-1}}\left(\eta_{2_{\infty}}
												^{-\f12+\delta/6}\right),\\
												\rho^{-1}_2\phi_{2_{\infty}}(\eta_{2_{\infty}},c,s)&=
												\sigma_{2_{\infty}}\mathrm{m}_{2_{\infty}}
												\left(1-\chi(M\eta_{2_{\infty}})\right)
												O_{\eta_{2_{\infty}}}^{\eta_{2_{\infty}}}\left(\eta_{2_{\infty}}
												^{-\f12+\delta/6}\right) =\sigma_{2_{\infty}}(c,s)\left(1-\chi(M\eta_{2_{\infty}})\right)
												O_{\eta_{2_{\infty}}}^{\eta_{2_{\infty}}}\left(\eta_{2_{\infty}}
												^{-\f12+\delta/6}\right).
											\end{split}
										\end{align}
									\end{small}
									For $\phi^*$, we use  the formula  \eqref{useful:trans-2}
									with $\rho=\rho_2=\xi^{-\f23}$, $\f{\pa_c\eta_{2_{\infty}}}{\eta_{2_{\infty}}}=O(\rho_2^{-1})$(by \eqref{est:paceta/eta2}), which are independent of $r$, $s$.
									Then by similar arguments as in the derivation of \eqref{bi-bd:20}, we obtain									
									\begin{align} \label{bi-bd:2infty}
										\begin{split}
											&\phi_{2_0}^*(\eta_{2_{\infty}},c,s,r)= O_{\eta_{2_{\infty}},c}^{\eta_{2_{\infty}},\rho_2}(1).
											\end{split}
									\end{align}
									For the cut-off function, by $c\in(0,V(0))$, $1\lesssim |\xi| \lesssim (V(0)-c)^{-\f32}$, we obtain
									\begin{align} \label{bd-cutoff:2infty}
										\chi\left(\f{(1-c)^{\f12}}{\xi}\right) \chi\left(\xi (V(0)-c)^{\f32}\right)=O_{\eta_{2_{\infty}},c}^{\eta_{2_{\infty}},\rho_2}(1).
									\end{align}
									Summing up \eqref{bd:2infty}, \eqref{bi-bd:2infty} and \eqref{bd-cutoff:2infty}, we obtain
									\begin{align}  \label{bd:a2infty}
										\begin{split}
											a_{2_{\infty},2_0}(\eta_{2_{\infty}},c,\lambda)=\left(1-\chi(M\eta_{2_{\infty}})\right)
											O_{\eta_{2_{\infty}},c}^{\eta_{2_{\infty}},
												\sigma_{2_{\infty}}(c,s)^{-1}}\left(\eta_{2_{\infty}}
											^{-\f12+\delta/6}\right).
											\end{split}
									\end{align}
									
									Now we are in a position to deal with \eqref{op-pieces2-new-case2}.
									Taking $\pa_c$ on  the non-oscillation part and the oscillation part, applying Lemma \ref{lem: pifi-Linfty} to the first term and taking the transform $c\to y_{2_{\infty}}$ for the second term, we deduce \eqref{process:integral20} for $(j,j')=(2_{\infty},2_0)$.
									Then we apply  Lemma \ref{lem:symbol-chi}, the bound \eqref{bd:a2infty}  and \eqref{int:sigmainfty} to estimate the integrals in  \eqref{process:integral20}($(j,j')=(2_{\infty},2_0)$) as  \begin{align*} \lesssim  \int_0^{V(0)}\sigma_{2_{\infty}}(c,s)dc+1\lesssim 1, \end{align*}
									which yields the desired result \eqref{op-pieces2-new-case2}.
									
									\textit{Case 3. }  $(j,j')=\{2_{\infty}, 2_{\infty}\}$.  WLOG, we assume $s\geq r$($r\geq s$ can be treated in a similar manner). We apply the change of variable $\xi\to \eta_{2_\infty}$ in the $\xi$-integral, reducing \eqref{op-pieces2} as
									\begin{align} I_{2_{\infty},2_{\infty}}
										:=\int_0^{V(0)}\left|\pa_c\left(p.v.\int_{\mathbb{R}} a_{2_{\infty},2_{\infty}}
										(\eta_{2_{\infty}},c,\lambda)\f{e^{i\eta_{2_{\infty}} \left(y_{2_{\infty}}(c,\lambda) \pm 1\right) }}{\eta_{2_{\infty}}}d\eta_{2_{\infty}}\right)\right|dc\lesssim 1\label{op-pieces2-new-case3}.
										\end{align}
										Here, the new coordinate $\eta_{2_{\infty}}$ takes $\xi \int_{0}^s\sqrt{\f{V(s')-c}{1-c}}ds'$ and  $\lambda=(z,s,r)$ is the parameter. The new functions are defined as follows.  $y_{2_{\infty}}(c,\lambda)$ takes as in \eqref{def:y2infty} and
										\begin{small} \begin{align} \label{def:a_infty2infty}
												\begin{split}
													&\quad a_{2_{\infty},2_{\infty}}(\eta_{2_{\infty}},c,\lambda):=\chi\left(\f{(1-c)^{\f12}}{\xi}\right)
													\chi\left(\xi(V(0)-c)^{\f32}\right)f(\eta_{2_{\infty}},c,s) f^* (\eta_{2_{\infty}},c,s,r)
													e^{\pm i\f{\int_{0}^r\sqrt{\f{V(s')-c}{1-c}}ds'}{\int_{0}^s\sqrt{\f{V(s')-c}{1-c}}ds'}},
												\end{split}
											\end{align}
										\end{small}
										where $\xi(\eta_{2_{\infty}},c,s)=\eta_{2_{\infty}}
										\left(\int_{0}^s\sqrt{\f{V(s')-c}{1-c}}ds'\right)^{-1}$, $f\in\{\phi_{2_{\infty}}, \bar{\phi}_{2_{\infty}}\}$, $f^*\in\{\phi_{2_{\infty}}^*, \bar{\phi}_{2_{\infty}}^*\}$,  with $\phi_{2_{\infty}}(\eta_{2_{\infty}},c,s)$ as in \eqref{22:2} and
										\begin{align}
											\label{22:3}\phi_{2_{\infty}}^*(\eta_{2_{\infty}},c,s,r)
											:=\widetilde{\phi_{2_{\infty}}}
											\left(\xi(\eta_{2_{\infty}},c,s),r,c\right).
											\end{align}
																				
										For $\phi^*$, we use  the formula  \eqref{useful:trans-2}
										with $\rho=\rho_2=\xi^{-\f23}$, $\f{\pa_c\eta_{2_{\infty}}}{\eta_{2_{\infty}}}=O(\rho_2^{-1})$ (by \eqref{lower:v-c}), which are independent of $r$.
										Then by similar arguments as in the derivation of \eqref{bi-bd:2infty}, we  obtain
										\begin{align}  \label{bi-bd:2inftyinfty}
											&\phi_{2_{\infty}}^*(\eta_{2_{\infty}},c,s,r)= O_{\eta_{2_{\infty}},c}^{\eta_{2_{\infty}},\rho_2}(1).
											\end{align}
										For the oscillation term, it follows from \eqref{estx1-3}  and \eqref{lower:v-c}  that for $s\geq r$, $s\gtrsim \xi^{-\f23}$,
										\begin{align}  \label{bi-bd:2osc}
											\left|\pa_c \left(e^{\pm i\f{\int_{0}^r\sqrt{\f{V(s')-c}{1-c}}ds'}{\int_{0}^s\sqrt{\f{V(s')-c}{1-c}}ds'}}
											\right) \right|\lesssim
											\left((V(r)-c)^{-1}+(V(s)-c)^{-1}\right)\cdot \f{\int_{0}^r\sqrt{\f{V(s')-c}{1-c}}ds'}{\int_{0}^s\sqrt{\f{V(s')-c}{1-c}}ds'}
											\lesssim \xi^{\f23}=\rho_2^{-1}.
											\end{align}
										Summing up \eqref{bd:2infty}, \eqref{bi-bd:2inftyinfty}, \eqref{bd-cutoff:2infty} and \eqref{bi-bd:2osc}, we obtain
										\begin{align}  \label{bd:a2inftyinfty}
											\begin{split}
												a_{2_{\infty},2_{\infty}}(\eta_{2_{\infty}},c,\lambda)
												=\left(1-\chi(M\eta_{2_{\infty}})\right)
												O_{\eta_{2_{\infty}},c}^{\eta_{2_{\infty}},
													\sigma_{2_{\infty}}(c,s)^{-1}}\left(\eta_{2_{\infty}}
												^{-\f12+\delta/6}\right). \end{split}
										\end{align}
										
										Finally, we deal with \eqref{op-pieces2-new-case3}.
										Taking $\pa_c$ on $a$ and the oscillation part, applying Lemma \ref{lem: pifi-Linfty} to the first term and taking the transform $c\to y_{2_{\infty}}$ for the second term, we deduce \eqref{process:integral20} for $(j,j')=(2_{\infty},2_{\infty})$.
										Then we apply  Lemma \ref{lem:symbol-chi}, the bound \eqref{bd:a2infty}  and \eqref{int:sigmainfty} to estimate the integrals in  \eqref{process:integral20} ($(j,j')=(2_{\infty},2_{\infty})$) as
										\begin{align*} \lesssim  \int_0^{V(0)}\sigma_{2_{\infty}}(c,s)dc+1\lesssim 1,
										\end{align*}
										which gives the desired result \eqref{op-pieces2-new-case3}.
								\end{proof}
								
								\begin{proposition}\label{lem:Boundness-K-pieces3}
									Let  $J^H_{3}=\{3_0,3_{\infty}\} $ and $\phi_j$($j\in J^H_{3}$) be defined in Definition \ref{def:decom-phi}.
									Then for $(j,j')\in J^H_{3}\times J^H_{3}$, it holds uniformly for  $(r,s,z)\in\mathbb{R}^+\times \mathbb{R}^+\times \mathbb{R}$  that
									\begin{small}
									\begin{align}\label{op-pieces3}
									&\int_0^{V(0)}\left|\pa_c\left(p.v.\int_{\mathbb{R}} \chi\left(\f{M^2(1-c)^{\f12}}{\xi}\right)\left(1- \chi\left(\xi (V(0)-c)^{\f32}\right) \right)
										\widetilde{\phi}_{j}(\xi,c,s)\widetilde{\phi}_{j'}(\xi,c,r)\f{e^{i\xi \sqrt{\f{c}{1-c}}z}}{\xi}d\xi\right)\right|dc\\&\lesssim s^{-\delta}+1.\nonumber
										\end{align}
										\end{small}
									\end{proposition}
									
									\begin{proof}
										 Without further illustration, assume that  $c\in(0,V(0))$, $ |\xi| \gtrsim  (V(0)-c)^{-\f32}$.
										Let the new coordinate be defined as follows
										\begin{align}
											\label{trans:eta3}\begin{split}
												\eta(\xi,c,s)= \left\{
												\begin{array}{l} \xi (V(0)-c)^{\f12}s := \eta_{3_0}\;\;\;\quad\quad\quad s\lesssim |\xi|^{-1}(V(0)-c)^{-\f12} ,\\
													\xi \int_{0}^s\sqrt{\f{V(s')-c}{1-c}}ds':= \eta_{3_{\infty}} \;\;\;\quad\quad s\gtrsim |\xi|^{-1}(V(0)-c)^{-\f12} , \end{array}\right.
											\end{split}
										\end{align}
										which satisfies( by \eqref{estx1-3} )
										\begin{small}
											\begin{align} \label{lower:eta_3infty}
												\eta_{3\infty}=\xi \int_{0}^s\sqrt{\f{V(s')-c}{1-c}}ds'\sim \xi (\f{s^{\f32}}{\langle s\rangle^{\f12}}+(V(0)-c)^{\f12}s)\gtrsim \xi (V(0)-c)^{\f12}s=\eta_{3_0},
											\end{align}
										\end{small}
										and \begin{align} \label{est:paceta/eta3}
											\left|\f{\pa_c\eta_{j}}{\eta_{j}}\right|\lesssim r_c^{-1}\sim \rho_3^{-1},\;\;j=3_0,3_{\infty},
											\end{align}
										where the weight in \eqref{def:rho-1} writes $\rho(c,\xi,s)=V(0)-c:= \rho_{3}$.  The corrector in  \eqref{def:rho*} takes the form \begin{align} \label{def:sigma3}
											\sigma_3(c,s)= \left\{
											\begin{array}{l} (V(0)-c)^{-1+\delta/6}s^{-\delta/6}\textbf{1}_{(V(0)-c)\geq 0 }(c):= \sigma_{3_0}(c,s) ,\\
												(V(0)-c)^{-1+\delta/6}s^{-\delta/6}\textbf{1}_{(V(0)-c)\geq 0 }(c):= \sigma_{3_{\infty}}(c,s), \end{array}\right.
										\end{align}
										which has the uniform bounds in $s$ that \begin{align}\label{int:sigma3}
											& \int_0^{V(0)}\sigma_{3_0}(c,s)dc\lesssim s^{-\delta},\\
											&\int_0^{V(0)}\sigma_{3_{\infty}}(c,s)dc\lesssim s^{-\delta}.  \label{int:sigma3infty} \end{align}
										Under the  coordinates \eqref{trans:eta3}, by \eqref{lower:eta_3infty}, the multiplier $\mathrm{m}$ in  \eqref{def:m-eta} behaves as
										\begin{align}\label{rewriten:m-eta3}
											&\quad\mathrm{m}_3(\xi,c,s)=
											\left\{
											\begin{array}{l}
												O_{\eta_{3_0}}^{\eta_{3_0}}
												\left(\eta_{3_0}^{\delta/6}\right):= \mathrm{m}_{3_0}(\xi,c,s),\\ O_{\eta_{2_{\infty}}}^{\eta_{2_{\infty}}}
												\left(\eta_{3_{\infty}}^{\delta/6}\right):= \mathrm{m}_{3_{\infty}}(\xi,c,s). \end{array}\right.
										\end{align}
										By \eqref{ineq:trho},  we can further decompose
										\begin{align}
											\begin{split}\label{decompose:rho-3}
												\rho_{3}^{-1}&=\mathrm{m}_{3_0}(\xi,c,s)\sigma_{3_0}(c,s),\;\text{if}\;\;s\lesssim |\xi|^{-1}(V(0)-c)^{-\f12} ,\\
												\rho_{3}^{-1}&=\mathrm{m}_{3_{\infty}}(\xi,c,s)\sigma_{3_{\infty}}(c,s),
												\;\text{if}\;\;s\gtrsim |\xi|^{-1}(V(0)-c)^{-\f12} .
											\end{split}
										\end{align}
										
										\textit{Case 1. }  $(j,j')= (3_0, 3_0)$.
										We first apply the change of varaibel $\xi\to \eta_{3_0}$ in the $\xi$-integral, reducing \eqref{op-pieces3} to
										\begin{align}
											I_{3_0,3_0}:=\int_0^{V(0)}\left|\pa_c\left(p.v.\int_{\mathbb{R}} a_{3_0,3_0}(\eta_{3_0},c,\lambda)\f{e^{i\eta_{3_0} y_{3_0}(c,\lambda) }}{\eta_{3_0}}d\eta_{3_0}\right)\right|dc\lesssim s^{-\delta}+1\label{op-pieces3-new-case1}.
											\end{align}
											Here, $\lambda=(z,s,r)$ is the parameter. The new functions are defined as follows
											\begin{align} \label{def:y30}
												y_{3_0}(c,\lambda)= \sqrt{\f{c}{1-c}}z/(V(0)-c)^{\f12}s,
												\end{align}
											\begin{align} \label{def:a_30}
												& a_{3_0,3_0}(\eta_{3_0},c,\lambda):=\chi\left(\f{(1-c)^{\f12}}{\xi}\right) \left(1-\chi\left(\xi(V(0)-c)^{\f32}\right)\right)\phi_{3_0}(\eta_{3_0},c,s)\phi_{3_0}^* (\eta_{3_0},c,s,r),
											\end{align}
											where $\xi(\eta_{3_0},c,s)=\eta_{3_0}(V(0)-c)^{-\f12}s^{-1}$, and
											\begin{align}
												\label{32:1}&\phi_{3_0}(\eta_{3_0},c,s)
												:=\widetilde{\phi}_{3_0}\left(\xi(\eta_{3_0},c,s),s,c\right), \quad \phi_{3_0}^*(\eta_{3_0},c,s,r):=\widetilde{\phi}_{3_0}
												\left(\xi(\eta_{3_0},c,s),r,c\right).
												\end{align}
					Using the formula \eqref{fm:f-tf1}, the definitions of $\tilde{\phi}_{3_0}$ in \eqref{def:tphi-3},  the bounds \eqref{3-phi}, in Proposition \ref{prop:summery}, the decomposition \eqref{decompose:rho-3}  for $\rho_3=V(0)-c$ and the bounds of $\mathrm{m}_{3_0}$ in \eqref{rewriten:m-eta3}, we obtain
									\begin{small}
										\begin{align}  \label{bd:30}
											\begin{split}
												&\phi_{3_0}(\eta_{3_0},c,s)=
												\chi(\eta_{3_0})
												O_{\eta_{3_0},c}^{\eta_{3_0},\rho_3}(\eta_{3_0})
												=\chi(\eta_{3_0})
												O_{\eta_{3_0},c}^{\eta_{3_0},(\sigma_{3_0}\mathrm{m}_{3_0})^{-1}}
												(\eta_{3_0})
												=\chi(\eta_{3_0})
												O_{\eta_{3_0},c}^{\eta_{3_0},\sigma_{3_0}^{-1}}(\eta_{3_0}),\\
												&\rho^{-1}_3\phi_{3_0}(\eta_{3_0},c,s)=
												\sigma_{3_0}\mathrm{m}_{3_0}\chi(\eta_{3_0})
												O_{\eta_{3_0}}^{\eta_{3_0}}(\eta_{3_0})
												=\sigma_{3_0}(c,s)\chi(\eta_{3_0})
												O_{\eta_{3_0}}^{\eta_{3_0}}(\eta_{3_0}).
											\end{split}
										\end{align}
									\end{small}
									For $\phi^*$, we take $\rho=\rho_3=V(0)-c$, $\f{\pa_c\eta_{3_0}}{\eta_{3_0}}=O(\rho_3^{-1})$ in the formula \eqref{useful:trans-2} and  use the bounds \eqref{3-phi}, \eqref{3-f+} in Proposition \ref{prop:summery},  to  obtain
									\begin{align} \label{bi-bd:30}
										\begin{split}
											&\phi_{3_0}^*(\eta_{3_0},c,s,r)= O_{\eta_{3_0},c}^{\eta_{3_0},\rho_3}(1).
										\end{split}
									\end{align}
									For the cut-off function, by $c\in(0,V(0))$, $|\xi| \gtrsim (V(0)-c)^{-\f32}$, $\eta_{3_0}=\xi (V(0)-c)^{\f12}s$, we have
									\begin{align} \label{bd-cutoff:30}
										\chi\left(\f{(1-c)^{\f12}}{\xi}\right)\left(1- \chi(\xi (V(0)-c)^{\f32})\right)=O_{\eta_{3_0},c}^{\eta_{3_0},\rho_3}(1).
									\end{align}
									Summing up \eqref{bd:30}, \eqref{bi-bd:30} and \eqref{bd-cutoff:30}, we obtain
									\begin{align}  \label{bd:a30}
										\begin{split}
											a_{3_0,3_0}(\eta_{3_0},c,\lambda)=\chi(\eta_{3_0})
											O_{\eta_{3_0},c}^{\eta_{3_0},\sigma_{3_0}(c,s)^{-1}}(\eta_{3_0}).
											\end{split}
									\end{align}
									
									 Now we  prove \eqref{op-pieces3-new-case1}.
									We take $\pa_c$ on  $a_{3_0,3_0}$ and the oscillation part. Applying Lemma \ref{lem: pifi-Linfty} to the first term and taking the transform $c\to y_{3_0}$ (via the Area formula) for the second term, we conclude that  \eqref{process:integral20} holds for $(j,j')=(3_0,3_0)$. We remark that while applying the one dimensional Area formula, we have used that for fixed $\lambda=(z,s,r)$ and $x$, the algebraic equation $y_{3_0}(c,\lambda)=x$ has finite number (independent of $\lambda$, $x$) roots $c$, due to the explicit form in \eqref{def:y30} and the Fundamental theorem of algebra.
									Then we apply  Lemma \ref{lem:symbol-chi} and the bound \eqref{bd:a30} to bound the last line in \eqref{process:integral20} as
									\begin{align*}
										\lesssim  \int_0^{V(0)}\sigma_{3_0}(c,s)dc+1,
									\end{align*}
									which along with  \eqref{int:sigma2} yields the desired bound \eqref{op-pieces2-new-case1}.
									
									\textit{Case 2. }  $(j,j')=\{3_{\infty}, 3_{0}\}$ or $(j,j')=\{3_0, 3_{\infty}\}$.
									We first consider $(j,j')=(3_{\infty},3_0)$. We  apply the change of variable $\xi\to \eta_{3_\infty}$ in the $\xi$-integral, reducing  \eqref{op-pieces3} to
									\begin{align} I_{3_{\infty},3_0}:=\int_0^{V(0)}\left|\pa_c\left(p.v.\int_{\mathbb{R}} a_{3_{\infty},3_0}(\eta_{3_{\infty}},c,\lambda)\f{e^{i\eta_{3_{\infty}} \left(y_{3_{\infty}}(c,\lambda) \pm 1\right) }}{\eta_{3_{\infty}}}d\eta_{3_{\infty}}\right)\right|dc\lesssim 1\label{op-pieces3-new-case2}.
									\end{align}
									Here, the new coordinate $\eta_{3_{\infty}}$ takes $\xi \int_{0}^s\sqrt{\f{V(s')-c}{1-c}}ds'$ and  $\lambda=(z,s,r)$ is the parameter.
									The new functions are defined as follows
									\begin{align} \label{def:y3infty} y_{3_{\infty}}(c,\lambda)= \sqrt{\f{c}{1-c}}z/\int_{0}^s\sqrt{\f{V(s')-c}{1-c}}ds'\;\;\text{is}\;\;\text{monotonic}\;\;\text{in}\;\;c,
									\end{align}
									\begin{align} \label{def:a_infty30}
										\begin{split}
											&\quad a_{3_{\infty},3_0}(\eta_{3_0},c,\lambda):=\chi\left(\f{(1-c)^{\f12}}{\xi}\right)
											\left(1-\chi(\xi(V(0)-c)^{\f32})\right)f(\eta_{3_{\infty}},c,s) \phi_{3_0}^* (\eta_{3_{\infty}},c,s,r),
										\end{split}
									\end{align}
									and $\xi(\eta_{3_{\infty}},c,s)=\eta_{3_{\infty}}
									\left(\int_{0}^s\sqrt{\f{V(s')-c}{1-c}}ds'\right)^{-1}$, $f\in\{\phi_{3_{\infty}}, \bar{\phi}_{3_{\infty}}\}$, where
									\begin{align}
										\label{32:2}&\phi_{3_{\infty}}(\eta_{3_{\infty}},c,s)
										:=\widetilde{\phi_{3_{\infty}}}\left(\xi(\eta_{3_{\infty}},c,s),s,c\right), \quad \phi_{3_0}^*(\eta_{3_{\infty}},c,s,r):=\widetilde{\phi}_{3_0}
										\left(\xi(\eta_{3_{\infty}},c,s),r,c\right),
										\end{align}
									with $\widetilde{\phi}_{3_{\infty}}(\xi,c,r)
									= \widetilde{\phi_{3_{\infty}}}(\xi,c,r) e^{i\xi \int_{0}^r\sqrt{\f{V(s')-c}{1-c}}ds'}+c.c$, $\widetilde{\phi_{3_{\infty}}}$,$\tilde{\phi}_{3_0}$ are defined in \eqref{def:tphi-3}.
									Noticing by  \eqref{lower:eta_3infty}, we have  for some fixed $M\gg 1$
									\begin{align*}
										\left(1-\chi\left(\xi (V(0)-c)^{\f12} s\right)\right)=\left(1-\chi\left(M\eta_{3_{\infty}}\right)\right)\left(1-\chi\left(\xi (V(0)-c)^{\f12} s\right)\right).
									\end{align*}
									 Therefore, using the definitions of $\widetilde{\phi_{3_{\infty}}}$ in \eqref{def:tphi-3}, the formula \eqref{fm:f-tf1},  the bounds \eqref{3-f+} in Proposition \ref{prop:summery}, the decomposition \eqref{decompose:rho-3} for $\rho_3=V(0)-c$ and the bounds of $\mathrm{m}_{3_{\infty}}$ in \eqref{rewriten:m-eta3}, we obtain
									\begin{small}
										\begin{align}  \label{bd:3infty}
											\begin{split}
												\phi_{3_{\infty}}(\eta_{3_{\infty}},c,s)&=
												\left(1-\chi\left(\xi (V(0)-c)^{\f12} s\right)\right)
												O_{\eta_{3_{\infty}},c}^{\eta_{3_{\infty}},\rho_3}\left(\eta_{3_{\infty}}
												^{-\f12}\right)\\
												&=
												\left(1-\chi\left(M\eta_{3_{\infty}}\right)\right)\left(1-\chi\left(\xi (V(0)-c)^{\f12} s\right)\right)
												O_{\eta_{3_{\infty}},c}^{\eta_{2_{\infty}},
													(\sigma_{3_{\infty}}\mathrm{m}_{3_{\infty}})^{-1}}\left(\eta_{3_{\infty}}
												^{-\f12}\right)\\
												&=\left(1-\chi\left(M\eta_{3_{\infty}}\right)\right)
												O_{\eta_{3_{\infty}},c}^{\eta_{3_{\infty}},
													(\sigma_{3_{\infty}}\mathrm{m}_{3_{\infty}})^{-1}}\left(\eta_{3_{\infty}}
												^{-\f12}\right)
												=\left(1-\chi(M\eta_{3_{\infty}})\right)
												O_{\eta_{3_{\infty}},c}^{\eta_{3_{\infty}},
													\sigma_{3_{\infty}}(c,s)^{-1}}\left(\eta_{3_{\infty}}
												^{-\f12+\delta/6}\right),\\
												\rho^{-1}_3\phi_{3_{\infty}}(\eta_{3_{\infty}},c,s)&=
												\sigma_{3_{\infty}}\mathrm{m}_{3_{\infty}}
												\left(1-\chi(M\eta_{3_{\infty}})\right)
												O_{\eta_{3_{\infty}}}^{\eta_{3_{\infty}}}(\eta_{3_{\infty}}
												^{-\f12+\delta/6}) =\sigma_{3_{\infty}}(c,s)\left(1-\chi(M\eta_{3_{\infty}})\right)
												O_{\eta_{3_{\infty}}}^{\eta_{3_{\infty}}}\left(\eta_{3_{\infty}}
												^{-\f12+\delta/6}\right).
											\end{split}
										\end{align}
									\end{small}
									We also get  by  \eqref{behave:Q-V(0)} that
									\begin{align}
										 \label{lower:v-c3} V(s)-c\sim  \f{s}{\langle s\rangle}+(V(0)-c)\geq  (V(0)-c).
									\end{align}
									For $\phi^*$, we use  the formula  \eqref{useful:trans-2}
									with $\rho=\rho_3=V(0)-c$, $\f{\pa_c\eta_{3_{\infty}}}{\eta_{3_{\infty}}}=O\left(\rho_3^{-1}\right)$ (by \eqref{lower:v-c3}), which are independent of $r$, $s$.
									Then by similar arguments as in the derivation of  \eqref{bi-bd:30}, we obtain
									\begin{align}  \label{bi-bd:3infty}
										\begin{split}
											&\phi_{3_0}^*(\eta_{3_{\infty}},c,s,r)= O_{\eta_{3_{\infty}},c}^{\eta_{3_{\infty}},
												\rho_3}(1). \end{split}
									\end{align}
									For the cut-off function, by $c\in(0,V(0))$, $|\xi| \gtrsim (V(0)-c)^{-\f32}$, we infer that
									\begin{align} \label{bd-cutoff:3infty}
										\chi\left(\f{(1-c)^{\f12}}{\xi}\right)\left(1- \chi\left(\xi (V(0)-c)^{\f32}\right)\right)=O_{\eta_{3_{\infty}},c}^{\eta_{3_{\infty}},\rho_3}(1).
									\end{align}
									Summing up \eqref{bd:3infty}, \eqref{bi-bd:3infty} and \eqref{bd-cutoff:3infty}, we obtain
									\begin{align}  \label{bd:a3infty}
										\begin{split}
											a_{3_{\infty},3_0}(\eta_{3_{\infty}},c,\lambda)=
											\left(1-\chi\left(M\eta_{3_{\infty}}\right)\right)
											O_{\eta_{3_{\infty}},c}^{\eta_{3_{\infty}},
												\sigma_{3_{\infty}}(c,s)^{-1}}\left(\eta_{3_{\infty}}
											^{-\f12+\delta/6}\right).
											\end{split}
									\end{align}
									
									Now we are in a position to deal with \eqref{op-pieces3-new-case2}.
									Taking $\pa_c$ on  the non-oscillation part and the oscillation part, applying Lemma \ref{lem: pifi-Linfty} to the first term and then transform $c\to y_{2_{\infty}}$, we then
									apply  Lemma \ref{lem:symbol-chi}, the bound \eqref{bd:a2infty} and \eqref{int:sigmainfty} to estimate the integrals in  \eqref{process:integral20} ($(j,j')=(3_{\infty},3_0)$) as
									\begin{align*} \lesssim  \int_0^{V(0)}\sigma_{3_{\infty}}(c,s)dc+1\lesssim 1,
									\end{align*}
									which gives the desired result \eqref{op-pieces3-new-case2}.
									For \emph{ $(j,j')=(3_0,3_{\infty})$}, the proof needs minor changes. Indeed, we apply the change of variable $\xi\to \eta_{3_\infty}:=\xi \int_{0}^r\sqrt{\f{V(s')-c}{1-c}}ds'$
									 in the $\xi$-integral, reducing \eqref{op-pieces3} to
									\begin{align}
										I_{3_0,3_{\infty}}:=\int_0^{V(0)}\left|\pa_c\left(p.v.\int_{\mathbb{R}} a_{3_0,3_{\infty}}(\eta_{3_{\infty}},c,\lambda)\f{e^{i\eta_{3_{\infty}} \left(y_{3_{\infty}}(c,\lambda) \pm 1\right) }}{\eta_{3_{\infty}}}d\eta_{3_{\infty}}\right)\right|dc\lesssim s^{-\delta}\label{op-pieces3-new-case2'},
										\end{align}
									where $y_{3_{\infty}}(c,\lambda)= \sqrt{\f{c}{1-c}}z/\int_{0}^r\sqrt{\f{V(s')-c}{1-c}}ds'$ is monotonic in $c$, and
									\begin{align*} 
										\begin{split}
											&\quad a_{3_{\infty},3_0}(\eta_{3_{\infty}},c,\lambda):=\chi\left(\f{(1-c)^{\f12}}{\xi}\right)
											\left(1-\chi\left(\xi(V(0)-c)^{\f32}\right)\right)f(\eta_{3_{\infty}},c,r) \phi_{3_0}^* (\eta_{3_{\infty}},c,s,r),
										\end{split}
									\end{align*}
									where $\xi(\eta_{3_{\infty}},c,r)=\eta_{3_{\infty}}
									\left(\int_{0}^r\sqrt{\f{V(s')-c}{1-c}}ds'\right)^{-1}$, $f\in\{\phi_{3_{\infty}}, \bar{\phi}_{3_{\infty}}\}$, and
									\begin{align}
										\label{32:2'}&\phi_{3_{\infty}}(\eta_{3_{\infty}},c,r)
										:=\widetilde{\phi_{3_{\infty}}}\left(\xi(\eta_{3_{\infty}},c,r),c,r\right), \quad \phi_{3_0}^*(\eta_{3_{\infty}},c,r,s):=\widetilde{\phi}_{3_0}
										\left(\xi(\eta_{3_{\infty}},c,r),c,s\right).
										\end{align}
									Using the bounds \eqref{3-phi}, \eqref{3-f+} in Proposition \ref{prop:summery}, the decomposition \eqref{decompose:rho-3} for $\rho_3=V(0)-c$ and the bounds of $\mathrm{m}_{3_{0}}$ in \eqref{rewriten:m-eta3}, we deduce
									\begin{align}  \label{bd:3infty'}
										\begin{split}
											\phi_{3_{\infty}}(\eta_{3_{\infty}},c,r)&
											=\left(1-\chi\left(M\eta_{3_{\infty}}\right)\right)
											O_{\eta_{3_{\infty}},c}^{\eta_{3_{\infty}},
												\rho_3^{-1}}\left(\eta_{3_{\infty}}
											^{-\f12}\right), \end{split}
									\end{align}
									and
									\begin{align} \label{bi-bd:30'}
										&\phi_{3_0}^*(\eta_{3_{\infty}},c,r,s)= O_{\eta_{3_{\infty}},c}^{\eta_{3_{\infty}},
											\sigma_{3_{\infty}}(c,s)^{-1}}(1),\quad  \rho_3^{-1}\phi_{3_0}^*(\eta_{3_{\infty}},c,r,s)=
										\sigma_{3_{\infty}}(c,s)O_{\eta_{3_{\infty}}}^{\eta_{3_{\infty}}}(1).
									\end{align}
									Here, due to $\rho_3^{-1}=(V(0)-c)^{-1}$ independent of $r$ or $s$, we can make the emerging factor $\sigma_{3_{\infty}}$ independent of $r$.
									Then we get by \eqref{bd:3infty'}, \eqref{bi-bd:30'} and \eqref{bd-cutoff:3infty} that
									\begin{align}  \label{bd:a3infty'}
										\begin{split}
											a_{3_0,3_{\infty}}(\eta_{3_{\infty}},c,\lambda)=
											\left(1-\chi\left(M\eta_{3_{\infty}}\right)\right)
											O_{\eta_{3_{\infty}},c}^{\eta_{3_{\infty}},
												\sigma_{3_{\infty}}(c,s)^{-1}}\left(\eta_{3_{\infty}}
											^{-\f12}\right), \end{split}
									\end{align}
									which along with \eqref{int:sigma3} gives \eqref{op-pieces3-new-case2'}, following the previous process.
									
									\textit{Case 3. }  $(j,j')=\{3_{\infty}, 3_{\infty}\}$.  WLOG, we assume $s\geq r$ ($r\geq s$ can be treated in a similar manner). We apply the change of variable $\xi\to \eta_{3_\infty}$ in the $\xi$-integral, reducing \eqref{op-pieces3} to
									\begin{align} I_{3_{\infty},3_{\infty}}
										:=\int_0^{V(0)}\left|\pa_c\left(p.v.\int_{\mathbb{R}} a_{3_{\infty},3_{\infty}}
										(\eta_{3_{\infty}},c,\lambda)\f{e^{i\eta_{3_{\infty}} \left(y_{3_{\infty}}(c,\lambda) \pm 1\right) }}{\eta_{3_{\infty}}}d\eta_{3_{\infty}}\right)\right|dc\lesssim s^{-\delta}\label{op-pieces3-new-case3}.
										\end{align}
										Here, the new coordinate $\eta_{3_{\infty}}$ takes $\xi \int_{0}^s\sqrt{\f{V(s')-c}{1-c}}ds'$ and  $\lambda=(z,s,r)$ is the parameter. The new functions are defined as follows.  $y_{3_{\infty}}(c,\lambda)$ takes as in \eqref{def:y3infty} and
									\begin{small} \begin{align} \label{def:a_infty3infty}
											\begin{split}
													&\quad a_{3_{\infty},3_{\infty}}(\eta_{3_{\infty}},c,\lambda)
												:=\chi\left(\f{(1-c)^{\f12}}{\xi}\right)
												\left(1-\chi\left(\xi(V(0)-c)^{\f32}\right)\right)f(\eta_{3_{\infty}},c,s) f^* (\eta_{3_{\infty}},c,s,r)
													e^{\pm i\f{\int_{0}^r\sqrt{\f{V(s')-c}{1-c}}ds'}{\int_{0}^s
															\sqrt{\f{V(s')-c}{1-c}}ds'}},
												\end{split}
										\end{align}
										\end{small}
										where $\xi(\eta_{3_{\infty}},c,s)=\eta_{3_{\infty}}
										\left(\int_{0}^s\sqrt{\f{V(s')-c}{1-c}}ds'\right)^{-1}$, $f\in\{\phi_{3_{\infty}}, \bar{\phi}_{3_{\infty}}\}$, $f^*\in\{\phi_{3_{\infty}}^*, \bar{\phi}_{3_{\infty}}^*\}$, and $\phi_{3_{\infty}}(\eta_{3_{\infty}},c,s)$ takes as in \eqref{32:2} and
										\begin{align}
											\label{32:3}&\phi_{3_{\infty}}^*(\eta_{3_{\infty}},c,s,r)
											:=\widetilde{\phi_{3_{\infty}}}
											\left(\xi(\eta_{3_{\infty}},c,s),r,c\right).
											\end{align}
										 For $\phi^*$, we use  the formula  \eqref{useful:trans-2}
										with $\rho=\rho_3=V(0)-c$, $\f{\pa_c\eta_{3_{\infty}}}{\eta_{3_{\infty}}}=O\left(\rho_3^{-1}\right)$ (by \eqref{lower:v-c3}), which are independent of $r$.
										Then by similar arguments as in the derivation of  \eqref{bi-bd:30}, we  obtain
										\begin{align}  \label{bi-bd:3inftyinfty}
											&\phi_{3_{\infty}}^*(\eta_{3_{\infty}},c,s,r)= O_{\eta_{3_{\infty}},c}^{\eta_{3_{\infty}},
												\rho_3}(1).
											\end{align}
										For the oscillation term, it follows from \eqref{estx1-3} and \eqref{lower:v-c3}  that for $s\geq r$, $s\gtrsim |\xi|^{-1}(V(0)-c)^{-\f12}$,
										\begin{align}  \label{bi-bd:3osc}
											\left|\pa_c \left(e^{\pm i\f{\int_{0}^r\sqrt{\f{V(s')-c}{1-c}}ds'}{\int_{0}^s\sqrt{\f{V(s')-c}{1-c}}ds'}}
											\right) \right|\lesssim
											(V(0)-c)^{-1}\cdot \f{\int_{0}^r\sqrt{\f{V(s')-c}{1-c}}ds'}{\int_{0}^s\sqrt{\f{V(s')-c}{1-c}}
												ds'}
											\lesssim (V(0)-c)^{-1}=\rho_3^{-1}.
											\end{align}
											
											Summing up  \eqref{bd:3infty}, \eqref{bi-bd:3inftyinfty}, \eqref{bd-cutoff:3infty} and \eqref{bi-bd:3osc}, we obtain
											\begin{align}  \label{bd:a3inftyinfty}
												\begin{split}
													a_{3_{\infty},3_{\infty}}(\eta_{3_{\infty}},c,\lambda)
													=\left(1-\chi\left(M\eta_{3_{\infty}}\right)\right)
													O_{\eta_{3_{\infty}},c}^{\eta_{3_{\infty}},
														\sigma_{3_{\infty}}(c,s)^{-1}}\left(\eta_{3_{\infty}}
													^{-\f12+\delta/6}\right). \end{split}
											\end{align}
											
											Finally, we deal with \eqref{op-pieces3-new-case3}.
											Taking $\pa_c$ on $a_{3_{\infty}, 3_{\infty}}$ and the oscillation part, applying Lemma \ref{lem: pifi-Linfty} to the first term and transform $c\to y_{3_{\infty}}$ for the second term, we deduce \eqref{process:integral20} for $(j,j')=(3_{\infty},3_{\infty})$.
											We then apply  Lemma \ref{lem:symbol-chi}, the bound \eqref{bd:a3infty}  and \eqref{int:sigma3infty} to estimate the integrals in  \eqref{process:integral20}($(j,j')=(3_{\infty},3_{\infty})$) as
											\begin{align*} \lesssim  \int_0^{V(0)}\sigma_{3_{\infty}}(c,s)dc+1\lesssim s^{-\delta}, \end{align*}
											which gives the desired result \eqref{op-pieces3-new-case3}.
									\end{proof}
									
									\begin{proposition}\label{lem:Boundness-K-pieces4}
										Let $\delta\ll 1$ be fixed, $J^H_{2}=\{4_0,4_{\infty}
										\}$ and $\phi_j$ ($j\in J^H_{4}$) be defined in Definition \ref{def:decom-phi}.
										Then for $(j,j')\in J^H_{4}\times J^H_{4}$, it holds uniformly for $(r,s,z)\in\mathbb{R}^+\times \mathbb{R}^+\times \mathbb{R}$ that \begin{align}\label{op-pieces4}											&\int_{V(0)}^{1-\delta}\left|\pa_c\int_{\mathbb{R}} \chi(\f{M^2(1-c)^{\f12}}{\xi}) \chi\left(\xi (c-V(0))^{\f32}\right)
											\widetilde{\phi}_{j}(\xi,c,s)\widetilde{\phi}_{j'}(\xi,c,r)\f{e^{i\xi \sqrt{\f{c}{1-c}}z}}{\xi}d\xi\right|dc\\&\lesssim s^{-\delta}+1.\nonumber
											\end{align}
									\end{proposition}
									
									\begin{proof}
										Let $M\gg 1$ be fixed, the new coordinates be defined as follows
										\begin{align}
											\label{trans:eta4}\begin{split}
												\eta(\xi,c,s)= \left\{
												\begin{array}{l} \xi s^{\f32} := \eta_{4_0}\;\;\;\quad\quad\quad\quad\quad\quad\quad s\lesssim M\xi^{-\f23} ,\\
													\xi \int_{r_c}^s\sqrt{\f{V(s')-c}{1-c}}ds':= \eta_{4_{\infty}} \;\;\;\quad\quad s\gtrsim M\xi^{-\f23}, \end{array}\right.
											\end{split}
										\end{align}
										which satisfies (by \eqref{estx4-6-G})
										\begin{align} \label{lower:eta_4infty}
											\eta_{4_\infty}=\xi \int_{r_c}^s\sqrt{\f{V(s')-c}{1-c}}ds'\sim \xi (s-r_c)^{\f32} \gtrsim \xi s^{\f32}=\eta_{4_0},\quad \text{if}\quad s\geq C M\xi^{-\f23}\geq 2r_c,\;M\gg 1 ,
											\end{align}
										and  \begin{align} \label{est:paceta/eta4}
											\f{\pa_c\eta_{4_0}}{\eta_{4_0}}=0,\;\; \quad   \left|\f{\pa_c\eta_{4_{\infty}}}{\eta_{4_{\infty}}}\right|\lesssim \xi^{\f23}=\rho_{4}^{-1},
										\end{align}
										where the weight in \eqref{def:rho-1} writes $\rho(c,\xi,s)=\xi^{-\f23}:= \rho_{4}$. The corrector in \eqref{def:rho*} takes
										\begin{align} \label{def:sigma4}
											\sigma_2(c,s)= \left\{
											\begin{array}{l} (c-V(0))^{-1+\delta/6}s^{-\delta/6}\textbf{1}_{(c-V(0))\geq 0 }(c):= \sigma_{4_0}(c,s) ,\\
												(c-V(0))^{-1+\delta/6}s^{-\delta/6}\textbf{1}_{s\gtrsim (c-V(0))\geq 0 }(c):= \sigma_{4_{\infty}}(c,s),
												\end{array}\right.
										\end{align}
										which have the uniform bounds in $s$ that \begin{align}\label{int:sigma4}
											& \int_0^{V(0)}\sigma_{4_0}(c,s)dc\lesssim s^{-\delta},\\
											&\int_0^{V(0)}\sigma_{4_{\infty}}(c,s)dc\lesssim 1.  \label{int:sigma4infty}
											\end{align}
										Under the  coordinates \eqref{trans:eta4}, by \eqref{lower:eta_4infty}, the multiplier $\mathrm{m}$ in \eqref{def:m-eta} behaves as
										\begin{align}\label{rewriten:m-eta4}
											&\quad\mathrm{m}_4(\xi,c,s)=
											\left\{
											\begin{array}{l}
												O_{\eta_{4_0}}^{\eta_{4_0}}
												\left(\eta_{4_0}^{\delta/9}\right):= \mathrm{m}_{4_0}(\xi,c,s),\\ O_{\eta_{4_{\infty}}}^{\eta_{4_{\infty}}}
												\left(\eta_{4_{\infty}}^{\delta/9}\right):= \mathrm{m}_{4_{\infty}}(\xi,c,s).
												\end{array}\right.
										\end{align}
										Using  the definitions of $\widetilde{\phi}_{4_0}$ in \eqref{def:tphi-4},  the bounds \eqref{4-phi}, \eqref{4-f+} in Proposition \ref{prop:summery},
										the decomposition \eqref{ineq:trho} for $\rho_4=\xi^{-\f23}$, the proof is divided into three cases,
										whose details are  similar as those in Proposition \ref{lem:Boundness-K-pieces4}. We leave the details to the readers.
									\end{proof}
									
									\begin{proposition}\label{lem:Boundness-K-pieces5}
										Let  $\delta\ll 1$  and $M\gg 1$ be fixed. Let $\tilde{\phi}_j$ be defined in \eqref{def:tphi-5}, Definition \ref{def:decom-phi}, where $j\in J^H_{5}=
										\{5_0,5_1,5_2,5_3,5_4,5_{\infty}
										\}  $.
										Then for $
										(j,j')\in J^H_{5}\times J^H_{5} $, it holds uniformly for  $(r,s,z)\in\mathbb{R}^+\times \mathbb{R}^+\times \mathbb{R}$  that
										\begin{small}
											\begin{align} \label{op-pieces5}
											&\int_{V(0)}^{1-\delta}\left|\pa_c\left(p.v.\int_{\mathbb{R}} \chi\left(\f{M^2(1-c)^{\f12}}{\xi}\right)\left(1- \chi\left(\xi (V(0)-c)^{\f32}/M\right) \right)
												\widetilde{\phi}_{j}(\xi,c,s)\widetilde{\phi}_{j'}(\xi,c,r)\f{e^{i\xi \sqrt{\f{c}{1-c}}z}}{\xi}d\xi\right)\right|dc\\
												&\lesssim s^{-\delta}+1.\nonumber
												\end{align}
											\end{small}
									\end{proposition}
									
									\begin{proof}
									Without further illustration, we assume $c\in(V(0),1-\delta),\;\xi r_c^{\f32}\gtrsim M$, which implies $0<c-V(0)\sim r_c\lesssim 1$. Define the new coordinates as
									\begin{align} \label{trans:eta5}\begin{split}
											\eta(\xi,c,s)= \left\{
											\begin{array}{l} \xi(c-V(0))^{\f12} s := \eta_{5_0}\;\;\;\quad\quad\quad\quad s\lesssim  M^{\f12}\xi^{-1}r_c^{-\f12},\\
												\xi(c-V(0))^{\f12} s:=
												\eta_{5_1}\;\;\quad\quad\quad\quad\; M^{\f12}\xi^{-1}r_c^{-\f12}\lesssim s\leq r_c/2,\\
												\xi (c-V(0))^{\f32}:=  \eta_{5_2}\;\;\;\quad\;\;\quad\quad\;\;\; r_c/2\leq s\leq r_c-C\xi^{-\f23},\\
												\xi (c-V(0))^{\f32}:=  \eta_{5_3}\;\;\;\;\; \quad\quad \quad \quad r_c-C\xi^{-\f23}\leq s\leq r_c+C\xi^{-\f23},\\
												\xi \int_{r_c}^s\sqrt{\f{V(s')-c}{1-c}}ds':= \eta_{5_{\infty}} \;\;\;\;\;\quad\quad s\geq r_c+C\xi^{-\f23},
												\end{array}\right.
										\end{split}
									\end{align}
									which  satisfy
									\begin{align} \label{est:paceta/eta5}
										\left|\f{\pa_c\eta_j}{\eta_j}\right|\lesssim r_c^{-1}=\rho_{j}(c)^{-1},\;j=5_0,...,5_3,\quad \text{and}\quad   \left|\f{\pa_c\eta_{5_{\infty}}}{\eta_{5_{\infty}}}\right|\lesssim (s-r_c)^{-1}=\rho_{5_{\infty}}(c,s)^{-1},
									\end{align}
									where the weights  defined  in \eqref{def:rho-1} writes
									\begin{align}\label{recall:rho-5} \rho(c,\xi,s)=
										\left\{
										\begin{array}{l}
											\rho_{5_0}=\rho_{5_1}=r_c\;\;\; \;\;\quad\quad\quad\quad\quad\quad s\leq r_c/2,\\
											\rho_{5_2}=(r_c-s) \;\;\; \;\quad\quad \quad\quad r_c/4\leq s\leq r_c-C\xi^{-\f23}, \quad\quad\quad\quad\\
											\rho_{5_3}=\xi^{-\f23}\;\;\; \;\;\;\;\;\;\quad\quad\quad\quad r_c-C\xi^{-\f23}\leq s\leq r_c+C\xi^{-\f23}, \\
											\rho_{5_{\infty}}=(s-r_c) \;\;\; \quad\quad\; \quad \quad s\geq r_c+C\xi^{-\f23} .
										\end{array}\right.
									\end{align}
									The corrector in this case \eqref{def:rho*} writes
									\begin{align}\label{def:sigma5} \sigma(c,s)=\left\{
										\begin{array}{l}  \sigma_{5_0}= \sigma_{5_1}=r_c^{-1+\delta/6}s^{-\delta/6}\textbf{1}_{s\leq r_c/2}(s)\quad\quad\quad\quad
											s\leq r_c/2\\
											\sigma_{5_2}=\sigma_{5_3}=  |r_c-s|^{-1+\delta/6}r_c^{-\delta/6}\textbf{1}_{r_c/4\leq s\leq 2r_c}(s)\;\;\;\quad r_c/4\leq s\leq 2r_c\\
											\sigma_{5_{\infty}}=(s-r_c)^{-1+\delta/6}r_c^{-\delta/6}\textbf{1}_{s\geq r_c}(s)\quad\quad\quad\quad s\geq r_c,\\
											\sigma_{5_{\infty}'}:=r_c^{-1+\delta/6}(s-r_c)^{-\delta/6}
											\textbf{1}_{s\geq r_c}(s)\quad\quad\quad\quad s\geq r_c,
										\end{array}\right.
									\end{align}
									which have the uniform bounds in $s$ that
									\begin{align}\label{int:sigma5}
										\begin{split} & \int_{V(0)}^{1-\delta}\sigma_j(c,s,\xi)dc\lesssim s^{-\delta},\;\;j=5_0,5_1;\;\;
											\int_{V(0)}^{1-\delta}\sigma_j(c,s)dc\lesssim 1,\;\;j=5_2,5_3,5_{\infty},5_{\infty}'.
										\end{split}
									\end{align}
									Under the  coordinates \eqref{trans:eta5},  for $\xi r_c^{\f32}\gtrsim M$, $c\in(V(0),1-\delta)$, the multiplier $\mathrm{m}$ in  \eqref{def:m-eta} behaves as
									\begin{align}\label{rewriten:m-eta}
										&\quad\mathrm{m}(\xi,c,s)=
										\left\{
										\begin{array}{l}
											\mathrm{m}_{5_0}(\xi,c,s)=O_{\eta_{5_0}}^{\eta_{5_0}}
											\left(\eta_{5_0}^{\delta/6}\right),\\
											\mathrm{m}_{5_1}(\xi,c,s)=O_{\eta_{5_1}}^{\eta_{5_1}}
											\left(\eta_{5_1}^{\delta/6}\right),\\
											\mathrm{m}_{5_2}(\xi,c,s)=O_{\eta_{5_2}}^{\eta_{5_2}}
											\left(\eta_{5_2}^{\delta/9}\right),\\
											\mathrm{m}_{5_3}(\xi,c,s)=O_{\eta_{5_3}}^{\eta_{5_3}}
											\left(\eta_{5_3}^{\delta/9}\right),\\
											\mathrm{m}_{5_{\infty}}(\xi,c,s)=O_{\eta_{5_{\infty}}}^{\eta_{5_{\infty}}}
											\left(\eta_{5_{\infty}}^{-\delta/9}\cdot
											\left(\xi r_c^{\f32}\right)^{\delta/9}\right).
											\end{array}\right.
									\end{align}
									Finally, by \eqref{ineq:trho},  we can further decompose
									\begin{align}
										\begin{split}\label{decompose:rho-5}
											\rho_{j}^{-1}&=\mathrm{m}_{j}(\xi,c,s)\sigma_{j}(c,s),\;j=5_0,5_1\;\;
											\text{if}\;\;s\leq r_c/2,\\
											\rho_{5_2}^{-1}&=\mathrm{m}_{5_2}(\xi,c,s)\sigma_{5_2}(c,s),\;\;
											\text{if}\;\;r_c/4\leq s\leq r_c-C\xi^{-\f23},\\
											\rho_{5_3}^{-1}&=\mathrm{m}_{5_3}(\xi,c,s)\sigma_{5_3}(c,s),\;\;
											\text{if}\;\;r_c-C\xi^{-\f23}\leq s\leq r_c+C\xi^{-\f23},\\
											\rho_{5_{\infty}}^{-1}&=\mathrm{m}_{5_{\infty}}(\xi,c,s)
											\sigma_{5_{\infty}}(c,s),\;\text{if}\;\;s\leq r_c+C\xi^{-\f23},\\
											r_c^{-1}&=O_{\eta_{5_{\infty}}}^{\eta_{5_{\infty}}}
											\left(\eta_{5_{\infty}}^{\delta/6}\right)
											\sigma_{5_{\infty}'}(c,s),\;\text{if}\;\;s\leq r_c+C\xi^{-\f23},
											\end{split}
									\end{align}
									where the last line follows by using
									\begin{align*}
										r_c^{-1}=r_c^{-1+\delta/6}(s-r_c)^{-\delta/6} \cdot \left(\xi r_c^{\f32}\right)^{-\delta/9}  \left(\xi (s-r_c)^{\f32}\right)^{\delta/9},
									\end{align*}
									and the fact that if $\xi r_c^{\f32}/\langle r_c\rangle ^{\f12}\gtrsim 1$ and $s\geq r_c+C\xi^{-\f23}\langle r_c\rangle^{\f13}$, then(by \eqref{estx4-6-G}) 									\begin{align} \label{observe:eta_infty}
										\xi \int_{r_c}^s\sqrt{\f{V(s')-c}{1-c}}ds'\gtrsim \left(\xi (s-r_c)^{\f32}/\langle r_c\rangle ^{\f12}\right)^{\f23}.
									\end{align}
									 Due to the oscillation behavior of $\widetilde{\phi}_{5_{\infty}}(\xi,c,r) $ and  $\widetilde{\phi}_{5_{\infty}}(\xi,c,s)$,
									 the proof of \eqref{op-pieces5} will be divided into four cases.
									
									 \textit{Case 1. }  $(j,j')\in \{5_0,5_1,5_2,5_3\}\times \{5_0,5_1,5_2,5_3\}$.
									
									We first apply the change of variable \eqref{trans:eta5} to $\tilde{\phi}_j$ defined in \eqref{def:tphi-5} that for $j\in \{5_0,5_1,5_2,5_3\}$,
									\begin{align}
										\label{25:1}&\phi_j(\eta_{j},c,s)
										:=\widetilde{\phi}_j\left(\xi(\eta_{j},c,s),s,c\right),
										\end{align}
									where $\xi(\eta_j,c,s)$ are the corresponding inverse functions of $\eta_j(\xi,c,s)$. And for  $j',j\in \{5_0,5_1,5_2,5_3\}$,
									\begin{align}
										\begin{split}
											&\phi_{j'}(\eta_{j},c,s,r):=\widetilde{\phi}_{j'}
											\label{25:2}\left(\xi(\eta_{j},c,s),r,c\right).
										\end{split}
									\end{align}
									 We summarize the quantitative bounds for $\phi_j$, $\phi_j'$ in variable $\eta_j$ and $c$ as follows. Using \eqref{fm:f-tf1} and the definitions of $\tilde{\phi}_j$ in \eqref{def:tphi-5}, the bounds \eqref{est:paceta/eta5},  the bounds \eqref{6-phi-1}, \eqref{6-phi-2} in Proposition \ref{prop:summery}  and  the decomposition \eqref{decompose:rho-5} and the bounds of $\mathrm{m}_j$ in \eqref{rewriten:m-eta}, we deduce that  for $c\in(V(0),1-\delta)$, $|\xi| r_c^{\f32}\gtrsim M\gg1$,
									\begin{align}  \label{bd:50-53}
											\begin{split}
												\phi_{5_0}(\eta_{5_0},c,s)&=
												\chi\left(\xi (c-V(0))^{\f12}s/M\right)\chi\left(\xi r_c^{\f12}s/M^{\f12}\right) O_{\eta_{5_0},c}^{\eta_{5_0},\rho_{5_0}}(\eta_{5_0})\\
												&=\chi\left(\eta_{5_0}/M\right)
												O_{\eta_{5_0},c}^{\eta_{5_0},\left(\sigma_{5_0}\mathrm{m}_{5_0}\right)^{-1}}
												(\eta_{5_0})\\
												&=\chi\left(\eta_{5_0}/M\right)
												O_{\eta_{5_0},c}^{\eta_{5_0},\sigma_{5_0}^{-1}}
												(\eta_{5_0}),\\
												\phi_{5_1}(\eta_{5_1},c,s)&=
												\left(1-\chi\left(\xi (c-V(0))^{\f12}s\right) \right)
												\left(\chi_+\left(4s/r_c\right)-\chi\left(\xi r_c^{\f12}s/M^{\f12}\right) \right)
												O_{\eta_{5_1},c}^{\eta_{5_1},\rho_{5_1}}\left(\eta_{5_1}^{-\f12}\right)\\
												&=\left(1-\chi\left(\eta_{5_1}\right) \right) O_{\eta_{5_1},c}^{\eta_{5_1},\left(\sigma_{5_1}\mathrm{m}_{5_1}\right)^{-1}}
												\left(\eta_{5_1}^{-\f12}\right)\\
												&=\left(1-\chi\left(\eta_{5_1}\right) \right) O_{\eta_{5_1},c}^{\eta_{5_1},\sigma_{5_1}^{-1}}
												\left(\eta_{5_1}^{-\f12+\delta/6}\right),
											\end{split}
										\end{align}
									
									and using the bounds \eqref{6-phi-4}, \eqref{6-f+a1} in Proposition \ref{prop:summery},
									\begin{align} \begin{split}
												&\quad\phi_{5_2}(\eta_{5_2},c,s)\\ &= \left(1-\chi\left(\xi (c-V(0))^{\f32}/M^{\f12}\right)\right) \left(1-\chi_+\left(4s/r_c\right)\right)\left(1-\chi_+\left(\xi^{\f23} (r_c-s)\right)\right)
												O_{\xi,c}^{\xi,\rho_{5_2}}\left(\left(\xi (r_c-s)^{\f32}\right)^{-\f16}\cdot\left(\xi r_c^{\f32}\right)^{-\f13}\right)\\
												&= \left(1-\chi\left(\xi (c-V(0))^{\f32}/M^{\f12}\right)\right)\left(1-\chi_+\left(\xi^{\f23} (r_c-s)\right)\right) O_{\eta_{5_2},c}^{\eta_{5_2},\rho_{5_2}}\left(\eta_{5_2}^{-\f13}\right)\\
												&= \left(1-\chi\left(\eta_{5_2}/M^{\f12}\right)\right) O_{\eta_{5_2},c}^{\eta_{5_2},(\sigma_{5_2}\mathrm{m}_{5_2})^{-1}}
												\left(\eta_{5_2}^{-\f13}\right)\\
												\label{bd:52-53}
												&=\left(1-\chi\left(\eta_{5_2}/M^{\f12}\right)\right)  O_{\eta_{5_2},c}^{\eta_{5_2},\sigma_{5_2}^{-1}}
												\left(\eta_{5_2}^{-\f16+\delta/9 }\right), \\
												&\quad\phi_{5_3}(\eta_{5_3},c,s)\\
												&=
												\left(1-\chi\left(\xi (c-V(0))^{\f32}/M^{\f12}\right)\right)
												\left(1- \chi_+\left(4s/r_c\right)\right)\chi_+\left(\xi^{\f23} (r_c-s)\right)\chi\left(\xi^{\f23}(s-r_c)\right) O_{\eta_{5_3},c}^{\eta_{5_3},\rho_{5_3}}\left(\eta_{5_3}^{-\f13}\right)\\
												&=
												\left(1-\chi\left(\xi (c-V(0))^{\f32}/M^{\f12}\right)\right) O_{\eta_{5_3},c}^{\eta_{5_3},\rho_{5_3}}\left(\eta_{5_3}^{-\f13}\right)\\
												&= \left(1-\chi\left(\eta_{5_3}/M^{\f12}\right)\right)
												O_{\eta_{5_3},c}^{\eta_{5_3},(\sigma_{5_3}\mathrm{m}_{5_3})^{-1}}
												\left(\eta_{5_3}^{-\f13}\right)\\
												&= \left(1-\chi\left(\eta_{5_3}/M^{\f12}\right)\right)
												O_{\eta_{5_3},c}^{\eta_{5_3},\sigma_{5_3}^{-1}}
												\left(\eta_{5_3}^{-\f13+\delta/9}\right).\end{split}
										\end{align}
								
									Similarly, we obtain
									\begin{align}  \label{bd:5-rho-1}
											\begin{split}
												&\rho_{5_0}(c)^{-1}\phi_{5_0}(\eta_{5_0},c,s)
												=\chi\left(\eta_{5_0}/M\right)
												O_{\eta_{5_0}}^{\eta_{5_0}}
												(\eta_{5_0})\sigma_{5_0}(c,s),\\
												&\rho_{5_1}(c)^{-1}\phi_{5_1}(\eta_{5_1},c,s)
												=\left(1-\chi(\eta_{5_1}) \right) O_{\eta_{5_1}}^{\eta_{5_1}}
												\left(\eta_{5_1}^{-\f12+\delta/6}\right)\sigma_{5_1}(c,s),\\
												&\rho_{5_2}(c,s,\xi)^{-1} \phi_{5_2}(\eta_{5_2},c,s)=
												\left(1-\chi\left(\eta_{5_2}/M^{\f12}\right)\right) O_{\eta_{5_2}}^{\eta_{5_2}}
												\left(\eta_{5_2}^{-\f16+\delta/9}\right)\sigma_{5_2}(c,s), \\
												&\rho_{5_3}(c,s,\xi)^{-1}\phi_{5_3}(\eta_{5_3},c,s)=
												\left(1-\chi\left(\eta_{5_3}/M^{\f12}\right)\right)
												O_{\eta_{5_3}}^{\eta_{5_3}}
												\left(\eta_{5_3}^{-\f13+\delta/9}\right)\sigma_{5_3}(c,s),\\
												&r_c^{-1}\phi_{5_2}(\eta_{5_2},c,s)=
												\left(1-\chi\left(\eta_{5_2}/M^{\f12}\right)\right)O_{\eta_{5_2}}^{\eta_{5_2}}
												\left(\eta_{5_2}^{-\f16+\delta/9}\right)\sigma_{5_2}(c,s),\\
												&r_c^{-1}\phi_{5_3}(\eta_{5_3},c,s)=
												\left(1-\chi\left(\eta_{5_3}/M^{\f12}\right)\right)
												O_{\eta_{5_3}}^{\eta_{5_3}}
												\left(\eta_{5_3}^{-\f13+\delta/9}\right)\sigma_{5_3}(c,s),
												\end{split}
										\end{align}
								where for the last two lines, we use $r_c^{-1}\leq (r_c-s)^{-1}=\rho_{5_2}^{-1}$ for $s\leq r_c-C\xi^{-\f23}\leq r_c$; and $r_c^{-1}\lesssim \xi^{\f23}=\rho_{5_3}^{-1}$ for $r_c-C\xi^{-\f23}\leq s\leq r_c+C\xi^{-\f23}$. For $\phi_{j'}^*$, the following derivative formula is useful:  \begin{small} \begin{align}  \notag
											\pa_c\phi_{j'}^*(\eta_{j},c,s,r)&=
											\pa_c\widetilde{\phi}_{j'}^*(\xi,c,r)
											-\f{\pa_c\eta_j(c,s,\xi)}{\eta_j(c,s,\xi)} (\xi\pa_{\xi})\widetilde{\phi}_{j'}^*(\xi,c,r)\\
											&=
											O\left(\rho_{j'}(c,r)^{-1}\right)(\rho\pa_c)\widetilde{\phi}_{j'}^*(\xi,c,r)
											+O\left(\left(\f{\pa_c\eta}{\eta_j}\right)(c,s)^{-1}\right) (\xi\pa_{\xi})\widetilde{\phi}_{j'}^*(\xi,c,r).\label{useful:trans-5}
										\end{align}
									\end{small}
									Using similar facts as in obtaining \eqref{bd:50-53} and \eqref{bd:52-53}, such as  the definition in \eqref{def:tphi-5}, the coordinates \eqref{trans:eta5},  the bounds \eqref{est:paceta/eta5}, the bounds \eqref{6-phi-1}, \eqref{6-phi-2}, \eqref{6-phi-4}, \eqref{6-f+a1} in Proposition \ref{prop:summery} and the decomposition in \eqref{decompose:rho-5},
									we can deduce that for $j',j\in \{5_0,5_1,5_2,5_3\}$,
									\begin{align}&  \label{bi-bd:50-5infty} \begin{split}
											&\phi_{j'}^*(\eta_{j},c,s,r)= O_{\eta_{j},c}^{\eta_{j},r_c}(1)
											+O_{\eta_{j},c}^{\eta_{j},\rho_j(c,s,\xi)}(1),\quad j'=5_0,5_1,\\
											&\phi_{j'}^*(\eta_{j},c,s,r)= O_{\eta_{j},c}^{\eta_{j},\sigma_{j'}(c,r)^{-1}}(1)
											+O_{\eta_{j},c}^{\eta_{j},\rho_j(c,s,\xi)}(1),\quad j'=5_2,5_3, \end{split}
									\end{align}
									where  $\rho_{5_0}=\rho_{5_1}=r_c$ are independent of $r$,$s$.
									Using the formula \eqref{fm:f-tf1} and the bound in \eqref{est:paceta/eta5},  it follows directly for $j\in \{5_0,5_1,5_2,5_3\} $  that  \begin{align} \label{behave:chi5}
										\chi\left(\f{(1-c)^{\f12}}{\xi}\right)\left(1- \chi\left(\xi (V(0)-c)^{\f32}/M\right)\right)
										=O_{\eta_j,c}^{\eta_j,r_c}(1).
									\end{align}
									 Therefore, using \eqref{bd:50-53}, \eqref{bd:52-53}, \eqref{bd:5-rho-1} and \eqref{bi-bd:50-5infty}, we can reduce
									\eqref{op-pieces5} to
									\begin{align}
										I_{j,j'}:=\int_{V(0)}^{1-\delta}\left|\pa_c\left(p.v.\int_{\mathbb{R}} a_{j,j'}(\eta_j,c,\lambda)\f{e^{i\eta_j y_j(c,\lambda) }}{\eta_j}d\eta_j\right)\right|dc\lesssim s^{-\delta}+1\label{op-pieces5-new-case1},
										\end{align}
									where $j,j'\in\{5_0,5_1,5_2,5_3\} $,  $\lambda=(s,r,z)$  and
									\begin{align} \label{def:a_j,j'}
										\begin{split}
											&a_{j,j'}(\eta_j,c,\lambda):=\chi\left(\f{(1-c)^{\f12}}{\xi}\right)\left(1- \chi\left(\xi(V(0)-c)^{\f32}/M\right) \right)\phi_{j}(\eta_{j},c,s)\phi_{j}^*(\eta_{j},c,s,r)\\
											&=
											\left\{
											\begin{array}{l}
												\chi(\eta_{5_0}/M)
												O_{\eta_{5_0},c}^{\eta_{5_0},\sigma_{5_0,j'}(c,s)^{-1}}
												(\eta_{5_0}),\;\;\;\quad\quad\quad\quad\quad\quad\quad j=5_0, \\
												\left(1-\chi(\eta_{5_1}) \right) O_{\eta_{5_1},c}^{\eta_{5_1},\sigma_{5_1,j'}(c,s)^{-1}}
												\left(\eta_{5_1}^{-\f12+\delta/6}\right),
												\;\;\;\quad\quad \quad\quad  j=5_1,\\
												\left(1-\chi\left(\eta_{5_2}/M^{\f12}\right)\right) O_{\eta_{5_2},c}^{\eta_{5_2},\sigma_{5_2,j'}(c,s,r)^{-1}}
												\left(\eta_{5_2}^{-\f16+\delta/9}\right),\;\;\; \quad\quad j=5_2,\\
												\left(1-\chi\left(\eta_{5_3}/M^{\f12}\right)\right)
												O_{\eta_{5_3},c}^{\eta_{5_3},\sigma_{5_3,j'}(c,s,r)^{-1}}
												\left(\eta_{5_3}^{-\f13+\delta/9}\right),\; \;\;\; \quad j=5_3, \end{array}\right.
										\end{split}
									\end{align}
									and
									\begin{align} \label{def:sigma_j,j'}
										\begin{split}
											\sigma_{j,j'}(c,s,r)= \left\{
											\begin{array}{l}
												\sigma_{j}(c,s) ,\;\quad\quad\quad\quad \quad j'=5_0,5_1,\;j=5_0,...,5_3,\\ \sigma_{j}(c,s)+\sigma_{j'}(c,r) ,\;j'=5_2,5_3,\; j=5_0,5_1,5_2,5_3, \end{array}\right.
										\end{split}
									\end{align}
									where $\sigma_{j}$ are as in \eqref{def:sigma5}. It follows from \eqref{int:sigma5} that the following uniform (in $s,r$) bounds hold,
									\begin{align} \label{sigma-5infty1}
										\int_{V(0)}^{1-\delta}\sigma_{j,j'}(c,s,r)dc\lesssim  s^{-\delta},\;\;j=5_0,5_1;\quad \int_{V(0)}^{1-\delta}\sigma_{j,j'}(c,s,r)dc\lesssim  1,\;\;j=5_2,5_3,
									\end{align}
									and
									\begin{align} \label{def:y_j51}
										\begin{split}
											y_j(c,z,s)= \left\{
											\begin{array}{l} \sqrt{\f{c}{1-c}}z/\sqrt{c-V(0)}s,\;\;\;\quad\quad j=5_0,5_1,\\
												\sqrt{\f{c}{1-c}}z/(c-V(0))^{\f32},\;\;\;
												\quad\quad
												j=5_2, 5_3. \end{array}\right.
										\end{split}
									\end{align}
									Indeed, applying Lemma \ref{lem: pifi-Linfty} to the first term and the transform $c\to y_j$ (via the Area formula) for the second term, we have
									\begin{align}
										\begin{split}
											\label{process:integral}I_{j,j'}&\leq \int_{V(0)}^{1-\delta}\sigma_{j,j'}(c,s,r)\left|\left(p.v.\int_{\mathbb{R}} \left(\sigma_{j,j'}^{-1}\pa_c\right)a_{j,j'}(\eta_j,c,\lambda)\f{e^{i\eta_j y_j(c,\lambda) }}{\eta_j}d\eta_j\right)\right|dc\\
											&\quad+
											\int_{V(0)}^{1-\delta}\left|\left(p.v.\int_{\mathbb{R}} a_{j,j'}(\eta_j,c,\lambda)e^{i\eta_j y_j(c,\lambda) }d\eta_j\right)\right|\cdot \left|\f{\pa y_j(c,\lambda)}{\pa c}\right|dc\\
											&\lesssim  \int_{V(0)}^{1-\delta}\sigma_{j,j'}(c,s,r)\left\|\int_{\mathbb{R}} \left(\sigma_{j,j'}^{-1}\pa_c\right)a_{j,j'}(\eta_j,c,\lambda)e^{i\eta_j x}d\eta_j\right\|_{L^1_x(\mathbb{R})}dc\\
											&\quad+
											\int_{\mathbb{R}}\left|\int_{\mathbb{R}} a_{j,j'}(\eta_j,c,\lambda)e^{i\eta_j y_j }d\eta_j\right|d y_j,
									\end{split}\end{align}
									while applying the one dimensional Area formula, we have used that for fixed $\lambda=(z,s,r)$ and $x$, the algebraic equation $y_j(c,\lambda)=x$ has finite number (independent of $\lambda$, $x$) roots $c$, due to the explicit form in \eqref{def:y_j51} and the Fundamental theorem of algebra. Then we apply Lemma \ref{lem:symbol-chi} and the bounds \eqref{def:a_j,j'} to bound the integrals as  \begin{align*} \lesssim  \int_{V(0)}^{1-\delta}\sigma_{j,j'}(c,s,r)dc+1. \end{align*}
									Using the definition \eqref{def:sigma_j,j'} and the bounds \eqref{int:sigma5},  we obtain the  desired bound \eqref{op-pieces5-new-case1}.
									
									\textit{Case 2. }  $(j,j')\in \{5_{\infty}\}\times \{5_0,5_1,5_2,5_3\}$ and $(j,j')\in  \{5_0,5_1,5_2,5_3\}\times\{5_{\infty}\}$.
									
									We only prove the first case.  The case where $(j,j')\in  \{5_0,5_1,5_2,5_3\}\times\{5_{\infty}\}$ is  symmetric  to the case where $(j,j')\in \{5_{\infty}\}\times \{5_0,5_1,5_2,5_3\}$; we merely swap the indices $j$ and $j'$ and the variable $r$ and $s$, and the proof is identical.
									
									We first apply the change of variable \eqref{trans:eta5} in the $\xi$-integral that
									\begin{align}
										\label{25:3}&\phi_{5_{\infty}}(\eta_{5_{\infty}},c,s)
										:=\widetilde{\phi_{5_{\infty}}}\left(\left(\int_{r_c}^s
										\sqrt{\f{V(s')-c}{1-c}}ds'\right)^{-1}\eta_{5_{\infty}},s,c\right),
										\end{align}
									where $\widetilde{\phi_{5_{\infty}}}$ satisfies $\widetilde{\phi}_{5_{\infty}}
									=\widetilde{\phi_{5_{\infty}}} e^{i\xi \int_{r_c}^s\sqrt{\f{V(s')-c}{1-c}}ds'}
									+c.c$  as  in  \eqref{def:tphi-5}, $\widetilde{\phi}_{5_{\infty}}$ is the function in the target integral \eqref{op-pieces5}. For $j'\in \{5_0,5_1,5_2,5_3\}$, we also take the transform
									\begin{align}
										\begin{split}
											&\phi_{j'}(\eta_{5_{\infty}},c,s,r):=\widetilde{\phi_{j'} }
											\left((\int_{r_c}^s
											\sqrt{\f{V(s')-c}{1-c}}ds')^{-1}\eta_{5_{\infty}},c,r\right).
										\end{split}
									\end{align}
									 Then for $(j,j')\in \{5_{\infty}\}\times \{5_0,5_1,5_2,5_3\}$, we reduce \eqref{op-pieces5} to
									\begin{align}
										I_{5_{\infty},j'}:=\int_{V(0)}^{1-\delta}\left|\pa_c\left(p.v.\int_{\mathbb{R}} a_{5_{\infty},j'}(\eta_{5_{\infty}},c,\lambda)\f{e^{i\eta_{5_{\infty}} (y_{5_{\infty}}(c,s,r)\pm 1) }}{\eta_{5_{\infty}}}d\eta_{5_{\infty}}\right)\right|dc\lesssim 1\label{op-pieces5-new-case2},
										\end{align}
									where   $\lambda=(z,s,r)$  and
									\begin{align} \label{def:y_j52}
										\begin{split}
											y_{5_{\infty}}(c,\lambda)= \sqrt{\f{c}{1-c}}z/\int_{r_c}^s\sqrt{\f{V(s')-c}{1-c}}ds' \quad \text{monotonic}\;\text{in}\;\;c,\end{split}
									\end{align}
									and
									\begin{align} \label{5phi-case2}
										a_{5_{\infty},j'}(\eta_{5_{\infty}},c,\lambda)
										:=\chi\left(\f{(1-c)^{\f12}}{\xi}\right)\left(1- \chi\left(\xi(V(0)-c)^{\f32}/M\right)\right)ff^*,
									\end{align}
									with
									$\xi(\eta_{5_{\infty}},c,s)=\left(\int_{r_c}^s
									\sqrt{\f{V(s')-c}{1-c}}ds'\right)^{-1}\eta_{5_{\infty}}$ and
									\begin{align*}
										f\in \big\{\phi_{5_{\infty}}(\eta_{5_{\infty}},c,s), \overline{\phi_{5_{\infty}}}(\eta_{5_{\infty}},c,s)\big\},\; \;f^*\in \big\{\phi_{j'}^*(\eta_{5_{\infty}},c,s,r), \overline{\phi_{j'}}^*(\eta_{5_{\infty}},c,s,r)\big\}.
									\end{align*}
									We claim that for $j'\in \{5_0,5_1,5_2,5_3\}$,
									\begin{align} \label{def:a_j,j'-52}
										\begin{split}
											a_{5_{\infty},j'}(\eta_{5_{\infty}},c,\lambda)=  \left(1-\chi\left(M^{\f12}\eta_{5_{\infty}}\right)\right) O_{\eta_{5_{\infty}},c}^{\eta_{5_{\infty}},\sigma_{5_{\infty},j'}
												(c,\lambda)
												^{-1}}\left(\eta_{5_{\infty}}^{-\f16}\right),
											\end{split}
									\end{align}
									where
									\begin{align} \label{def:sigma_j,j'52} \begin{split}
											\sigma_{5_{\infty},j'}(c,s,r)= \left\{
											\begin{array}{l}
												\sigma_{5_{\infty}}(c,s)+\sigma_{5_{\infty}'}(c,s) ,\;j'=5_0,5_1, \\
												\sigma_{5_{\infty}}(c,s) +\sigma_{j'}(c,r),\;j'=5_2,5_3, \end{array}\right.
										\end{split}
									\end{align}
									where $\sigma_{5_{\infty}}$, $\sigma_{5_{\infty}'}$, $\sigma_{j'}$ are as in \eqref{def:sigma5}. It follows from \eqref{int:sigma5} that  for $j'\in \{5_0,5_1,5_2,5_3\}$, the following uniform (in $s,r$) bounds hold,
									\begin{align} \label{sigma-5infty}
										\int_{V(0)}^{1-\delta}\sigma_{5_{\infty},j'}(c,s,r)dc\lesssim  1.
									\end{align}
									With \eqref{def:a_j,j'-52} at hand, \eqref{op-pieces5-new-case2} follows by a similar process as \eqref{process:integral}, using Lemma \ref{lem: pifi-Linfty}, the Area formula (with \eqref{def:y_j52}),  Lemma \ref{lem:symbol-chi}, and \eqref{sigma-5infty}.
									
									It remains to prove \eqref{def:a_j,j'-52}.
									To this end, we first notice by  \eqref{estx4-6-G} that if $\xi r_c^{\f32}\gtrsim 1$ and $s\geq r_c+C\xi^{-\f23}$, then
									\begin{align} \label{observe:eta_infty5}
									\xi \int_{r_c}^s\sqrt{\f{V(s')-c}{1-c}}ds'\gtrsim \left(\xi (s-r_c)^{\f32}\right)^{\f23},
								\end{align}
									since
									\begin{align*}
										\notag&\xi \int_{r_c}^s\sqrt{\f{V(s')-c}{1-c}}ds'\sim \xi (s-r_c)^{\f32}/\langle s\rangle ^{\f12}\\
										&\sim \left\{
										\begin{array}{l}
											\xi (s-r_c)^{\f32},\quad \text{if}\quad  s\leq 2r_c,\\
											\left(\xi (s-r_c)^{\f32}\right)^{\f23}\cdot \left(\xi r_c^{\f32}\right)^{\f13}\cdot\left(\f{(s-r_c) }{r_c\langle s\rangle}\right)^{\f12},\quad \text{if}\quad  s\geq 2r_c. \end{array}\right.
										\end{align*}
										Therefore, using the definitions of $\tilde{\phi}_{5_{\infty}}$ in \eqref{def:tphi-5}, the identity \eqref{fm:f-tf1}, combing the bound \eqref{est:paceta/eta5}, \eqref{6-f+} in Proposition \ref{prop:summery}, \eqref{observe:eta_infty5}, and finally using the decomposition \eqref{decompose:rho-5} as well as the bounds of $\mathrm{m}_{5_{\infty}}$ in \eqref{rewriten:m-eta},
										we obtain
										\begin{align}
											\label{bd:54case2} \begin{split}
												\phi_{5_{\infty}}(\eta_{5_{\infty}},c,s)&=
												\left(1-\chi\left(M^{\f12}\xi \int_{r_c}^s\sqrt{\f{V(s')-c}{1-c}}ds'\right)\right) \left(1-\chi\left(\xi^{\f23}(s-r_c)\right)\right) \\
												&\quad \;\;\left(1- \chi_+\left(4s/r_c\right)\right)
												\chi_+\left(\xi^{\f23}(r_c-s)\right)
												O_{\eta_{5_{\infty}},c}^{\eta_{5_{\infty}},\rho_{5_{\infty}}}
												\left(\eta_{5_{\infty}}^{-\f16}\cdot \left(\xi r_c^{\f32}\right)^{-\f13}\right)
												\\
												&= \left(1-\chi\left(M^{\f12}\eta_{5_{\infty}}\right)\right)
												O_{\eta_{5_{\infty}},c}^{\eta_{5_{\infty}},\rho_{5_{\infty}}}
												\left(\eta_{5_{\infty}}^{-\f16}\cdot
												\left(\xi r_c^{\f32}\right)^{-\f13}\right)\\
												&= \left(1-\chi\left(M^{\f12}\eta_{5_{\infty}}\right)\right) O_{\eta_{5_{\infty}},c}^{\eta_{5_{\infty}},(\sigma_{5_{\infty}}
													\mathrm{m}_{5_{\infty}})^{-1}}\left(\eta_{5_{\infty}}^{-\f16}\right)\\
												&= \left(1-\chi\left(M^{\f12}\eta_{5_{\infty}}\right)\right) O_{\eta_{5_{\infty}},c}^{\eta_{5_{\infty}},\sigma_{5_{\infty}}
													^{-1}}\left(\eta_{5_{\infty}}^{-\f16}\right),
										\end{split}\end{align}
										with $M\gg 1$. Similarly, we get
										\begin{align}  \label{bd:5-rho-11case2} \begin{split} \rho_{5_{\infty}}(c,s)^{-1}\phi_{5_{\infty}}(\eta_{5_{\infty}},c,s)
												&=\left(1-\chi\left(M^{\f12}\eta_{5_{\infty}}\right)\right) O_{\eta_{5_{\infty}}}^{\eta_{5_{\infty}}}\left(\eta_{5_{\infty}}^{-\f16}\right) \sigma_{5_{\infty}}(c,s),\\
												r_c^{-1}\phi_{5_{\infty}}(\eta_{5_{\infty}},c,s)
												&=\left(1-\chi\left(M^{\f12}\eta_{5_{\infty}}\right)\right) O_{\eta_{5_{\infty}}}^{\eta_{5_{\infty}}}\left(\eta_{5_{\infty}}
												^{-\f16+\delta/6}\right)\sigma_{5_{\infty}'}(r_c,s),
											\end{split}
										\end{align}
										where  $\sigma_{5_{\infty}}$, $\sigma_{5_{\infty}'}$  are as in \eqref{def:sigma5}, and we used the last line in decomposition \eqref{decompose:rho-5}.
										For $\phi_{j'}^*$, using the formula \eqref{useful:trans-5} and similar facts, we can also obtain
										\begin{align} \label{bi-bd:5inftycase2}
											\begin{split}
												&\phi_{j'}^*(\eta_{5_{\infty}},c,s,r)= O_{\eta_{5_{\infty}},c}^{\eta_{5_{\infty}},r_c}(1)
												+O_{\eta_{5_{\infty}},c}^{\eta_{5_{\infty}},\rho_{5_{\infty}}(c,s)}(1),\quad j'=5_0,5_1,\\
												&\phi_{j'}^*(\eta_{5_{\infty}},c,s,r)= O_{\eta_{5_{\infty}},c}^{\eta_{5_{\infty}},\sigma_{j'}(c,r)^{-1}}(1)
												+O_{\eta_{5_{\infty}},c}^{\eta_{5_{\infty}},\rho_{5_{\infty}}(c,s)}(1),\quad j'=5_2,5_3. \end{split}
										\end{align}
										It follows directly from the identity \eqref{fm:f-tf1} and the last bound in \eqref{est:paceta/eta5} that
										\begin{align} \label{behave:chi5case2}
											\chi\left(\f{(1-c)^{\f12}}{\xi}\right)\left(1- \chi\left(\xi (V(0)-c)^{\f32}/M\right) \right)
											=O_{\eta_{5_{\infty}},c}^{\eta_{5_{\infty}},
												r_c}(1)+O_{\eta_{5_{\infty}},c}^{\eta_{5_{\infty}},
												\rho_{5_{\infty}}(c,s)}(1).
										\end{align}
										Then \eqref{def:a_j,j'-52} follows from \eqref{bd:54case2}, \eqref{bd:5-rho-11case2}, \eqref{bi-bd:5inftycase2} and         \eqref{behave:chi5case2}.

										\textit{Case 3. }  $(j,j')\in  \{5_{\infty}\}\times\{5_{\infty}\}$.
										WLOG, we assume $s\geq r>r_c$ ($r\geq s>r_c$ can be treated in a similar manner).
										We apply the change of variable $\eta_{5_{\infty}}=\xi \int_{r_c}^s
										\sqrt{\f{V(s')-c}{1-c}}ds'$(as in \eqref{trans:eta5}) that
										\begin{align*}
											&\phi_{5_{\infty}}(\eta_{5_{\infty}},c,s)
											:=\widetilde{\phi_{5_{\infty}}}
											\left((\int_{r_c}^s
											\sqrt{\f{V(s')-c}{1-c}}ds')^{-1}\eta_{5_{\infty}},s,c\right),
											\end{align*}
										and
										\begin{align*}
											\begin{split}
												&\phi_{5_{\infty}}^*(\eta_{5_{\infty}},c,s,r):=
												\widetilde{\phi_{5_{\infty}}} \left((\int_{r_c}^s
												\sqrt{\f{V(s')-c}{1-c}}ds')^{-1}\eta_{5_{\infty}},c,r\right).
											\end{split}
										\end{align*}
										In the following, we use the old notation $Q(r,c)=\f{V(r)-c}{1-c}$. By Lemma \ref{lem:behave-Q} and Lemma \ref{lem:behave-x}, we infer that for $s>r_c$,
										 \begin{align}
											\begin{split}\label{Basic:Q}
												&Q(s,c)\sim \f{s-r_c}{\langle s\rangle },\quad \int_{r_c}^sQ(s',c)^{\f12}ds'\sim (s-r_c)^{\f32};\\ &\left|\f{\pa_cQ(s,c)}{Q(s,c)}\right|
												+\left|\f{\pa_c(\int_{r_c}^sQ(s',c)^{\f12}ds')}
												{\int_{r_c}^sQ(s',c)^{\f12}ds'}\right|\lesssim (s-r_c)^{-1}.
											\end{split}
										\end{align}
										Then by   \eqref{def:tphi-5}, we reduce \eqref{op-pieces5} to
									\begin{align}  \label{est:5-infty}
											I_{5_{\infty},5_{\infty}}:=\int_{V(0)}^{1-\delta}
										\left|\pa_c\left(p.v.\int_{\mathbb{R}} a_{5_{\infty},5_{\infty}}(\eta_{5_{\infty}},c,\lambda) \f{e^{i\eta_{5_{\infty}}(y_{5_{\infty}}(c,\lambda)\pm 1)}}{\eta_{5_{\infty}}}d\eta_{5_{\infty}}\right)\right|dc\lesssim  1,
									\end{align}
										where $\lambda=(z,s,r)$, $y_{5_{\infty}}$ as in \eqref{def:y_j52}, and
											\begin{align*}
												a_{5_{\infty},5_{\infty}}(\eta_{5_{\infty}},c,\lambda)
												:=\chi\left(\f{(1-c)^{\f12}}{\xi}\right)\left(1- \chi\left(\xi(V(0)-c)^{\f32}/M\right) \right)f f^* e^{\pm i\f{ \int_{r_c}^rQ(s',c)^{\f12}ds'}{\int_{r_c}^sQ(s',c)^{\f12}ds'}},
											\end{align*}
										with $f\in \big\{\phi_{5_{\infty}}(\eta_{5_{\infty}},c,s), \overline{\phi_{5_{\infty}}}(\eta_{5_{\infty}},c,s)\big\}$, $f^*\in \big\{\phi_{5_{\infty}}^*(\eta_{5_{\infty}},c,s,r), \overline{\phi_{5_{\infty}}^*}(\eta_{5_{\infty}},c,s,r)\big\}$ and
										$\xi=\xi(\eta_{5_{\infty}},c,s)=\left(\int_{r_c}^s
										\sqrt{\f{V(s')-c}{1-c}}ds'\right)^{-1}\eta_{5_{\infty}}$.
										We notice that compared with \eqref{5phi-case2} in Case 2, there is an extra oscillation type term $e^{\pm i\f{ \int_{r_c}^rQ(s',c)^{\f12}ds'}{\int_{r_c}^sQ(s',c)^{\f12}ds'}}$.
										Using the formula \eqref{useful:trans-5}, the bound \eqref{6-f+} in Proposition \ref{prop:summery} and  \eqref{ineq:trho},
										and  the bounds of $\mathrm{m}_{5_{\infty}}$ in \eqref{rewriten:m-eta}, we have
										\begin{align}
											\begin{split} \label{bd:tphi5-case4'}
												&\phi_{5_{\infty}}^*(\eta_{5_{\infty}},c,s,r)= O_{\eta_{5_{\infty}},c}^{\eta_{5_{\infty}},\sigma_{5_{\infty}}(c,r)^{-1}}(1)
												+O_{\eta_{5_{\infty}},c}^{\eta_{5_{\infty}},\rho_{5_{\infty}}(c,s)}(1).
												\end{split}
										\end{align}
										 We notice by \eqref{Basic:Q} and $s-r_c\geq r-r_c>0$ that
										\begin{align}
											\left|\pa_c\Big(e^{\pm i\f{ \int_{r_c}^rQ(s',c)^{\f12}ds'}{\int_{r_c}^sQ(s',c)^{\f12}ds'}}\Big) \right|\lesssim (s-r_c)^{-1}.
											\end{align}
										Then, by using the bounds \eqref{bd:54case2}, \eqref{bd:5-rho-11case2}, \eqref{behave:chi5case2}, \eqref{bd:tphi5-case4'} and \eqref{int:sigma5},  we obtain
										\begin{align} \label{def:a_infty}
											\begin{split}
												a_{5_{\infty},5_{\infty}}(\eta_{5_{\infty}},c,\lambda)=  \left(1-\chi\left(M^{\f12}\eta_{5_{\infty}}\right)\right) O_{\eta_{5_{\infty}},c}^{\eta_{5_{\infty}},\sigma_{5_{\infty},5_{\infty}}
													(c,\lambda)
													^{-1}}\left(\eta_{5_{\infty}}^{-\f16}\right), \end{split}
										\end{align}
										where
										\begin{align} \notag
										&\sigma_{5_{\infty},5_{\infty}}(c,s,r)= \sigma_{5_{\infty}}(c,s)+\sigma_{5_{\infty}'}(c,s) +\sigma_{5_{\infty}}(c,r) ,\\ &\;\text{with}\;\;\int_{V(0)}^{1-\delta}\sigma_{5_{\infty},5_{\infty}}
											(c,s,r)dc\lesssim  1,\label{bd:sigma_5infty}  \end{align}
										where $\sigma_{5_{\infty}}$, $\sigma_{5_{\infty}'}$ are as in \eqref{def:sigma5}.  To prove \eqref{est:5-infty}, we take $\pa_c$ on $a$ and the oscillation, respectively.  Then we first apply Lemma \ref{lem: pifi-Linfty} to the first term and the transform $c\to y$ for the second  term, and then use the bound \eqref{def:a_infty} and Lemma  \ref{lem:symbol-chi}, and finally use \eqref{bd:sigma_5infty}, to obtain
										\begin{align}
											\begin{split}
												\label{process:integral-5case4} I_{5_{\infty},5_{\infty}}&\leq \int_{V(0)}^{1-\delta}\sigma_{5_{\infty},5_{\infty}}(c,\lambda)
												\left|\left(p.v.\int_{\mathbb{R}} (\sigma_{5_{\infty},5_{\infty}}^{-1}\pa_c)
												a_{5_{\infty},5_{\infty}}(\eta_{5_{\infty}},c,\lambda)\f{e^{i\eta_{5_{\infty}} y_{5_{\infty}}(c,\lambda) }}{\eta_{5_{\infty}}}d\eta_{5_{\infty}}\right)\right|dc\\
												&\quad+
												\int_{V(0)}^{1-\delta}\left|\left(p.v.\int_{\mathbb{R}} a_{5_{\infty},5_{\infty}}(\eta_{5_{\infty}},c,\lambda)e^{i\eta_{5_{\infty}} y_{5_{\infty}}(c,\lambda) }d\eta_{5_{\infty}}\right)\right|\cdot \left|\f{\pa y_{5_{\infty}}(c,\lambda)}{\pa c}\right|dc\\
												&\lesssim  \int_{V(0)}^{1-\delta}\sigma_{5_{\infty},5_{\infty}}(c,\lambda)
												\left\|\int_{\mathbb{R}} (\sigma_{5_{\infty},5_{\infty}}^{-1}\pa_c)a_{_{5_{\infty}},5_{\infty}}
												(\eta_{5_{\infty}},c,\lambda)
												e^{i\eta_{5_{\infty}} x}d\eta_{5_{\infty}}\right\|_{L^1_x(\mathbb{R})}dc\\
												&\quad+
												\int_{\mathbb{R}}\left|\int_{\mathbb{R}} a_{5_{\infty},{5_{\infty}}}(\eta_{5_{\infty}},c,\lambda)e^{i\eta_{5_{\infty}} y_{5_{\infty}} }d\eta_{5_{\infty}}\right|d y_{5_{\infty}}\\
												&\lesssim \int_{V(0)}^{1-\delta}\sigma_{5_{\infty},5_{\infty}}(c,\lambda) dc\lesssim 1.\end{split}\end{align}

												This finishes the proof of the proposition.
									\end{proof}
									
									\begin{proposition}\label{lem:Boundness-K-pieces6}
										Let  $\delta\ll 1$ such that $1-\delta\geq V(1)$ and $M\gg 1$ be fixed. Let  $J^H_{6}=
										\{6_0,6_1,6_2,6_3,6_4,6_5,6_{\infty}
										\}  $ and $\tilde{\phi}_j$($j\in J^H_{6}$) be defined in \eqref{def:tphi-6}, Definition \ref{def:decom-phi}.
										Then for $(j,j')\in \tilde{J}^H_{6}\times  \tilde{J}^H_{6} $, it holds uniformly for  $(r,s,z)\in\mathbb{R}^+\times \mathbb{R}^+\times \mathbb{R}$  that
										\begin{small}
											\begin{align}\label{op-pieces6} &\int_{1-\delta}^1\left|\pa_c\left(p.v.\int_{\mathbb{R}} \chi\left(\f{M^2(1-c)^{\f12}}{\xi}\right)\left(1- \chi\left(\xi/(1-c)^{\f13}M\right) \right)
												\widetilde{\phi}_{j}(\xi,c,s)\widetilde{\phi}_{j'}(\xi,c,r)\f{e^{i\xi \sqrt{\f{c}{1-c}}z}}{\xi}d\xi\right)\right|dc\\
												&\lesssim s^{-\delta}+1.\nonumber
												\end{align}
										\end{small}
									\end{proposition}
									
									\begin{proof}
										Without further illustration, we will limit the range $c\in(1-\delta, 1),\;\xi r_c\gtrsim M$. And in this case, it holds that $(1-c)^{-\f13}\sim r_c\gtrsim  1$.
										Let the new coordinate be defined as follows
										\begin{align} \label{trans:eta6}\begin{split}
												\eta(\xi,c,s)= \left\{
												\begin{array}{l} \f{ \xi s}{(1-c)^{\f12}} := \eta_{6_0}=\eta_{6_1}\;\;\;\quad\quad\quad\quad s\leq \f12,\\
													\f{\xi}{(1-c)^{\f13}}:=
													\eta_{6_2}\;\;\quad\quad\quad\quad\; \f12\leq s\leq r_c/2,\\
													\f{\xi}{(1-c)^{\f13}}:=  \eta_{6_3}\;\;\;\quad\;\; \f{r_c}{2}\leq s\leq r_c-C\xi^{-\f23}r_c^{\f13},\\
													\f{\xi}{(1-c)^{\f13}}:=  \eta_{6_4}\;\;\; \quad\quad \quad \quad r_c-C\xi^{-\f23}r_c^{\f13}\leq s\leq r_c+C\xi^{-\f23}r_c^{\f13},\\
													\xi \int_{r_c}^s\sqrt{\f{V(s')-c}{1-c}}ds':= \eta_{6_{\infty}} \;\;\;\quad\quad s\geq r_c+C\xi^{-\f23}r_c^{\f13}, \end{array}\right.
											\end{split}
										\end{align}
										which satisfy
										\begin{align} \label{est:paceta/eta6}
											\left|\f{\pa_c\eta_j}{\eta_j}\right|\lesssim r_c^{3}=\rho_{j}(c)^{-1},\;j=6_0,...,6_3,\;\;
											\text{and}\quad   \left|\f{\pa_c\eta_{6_{\infty}}}{\eta_{6_{\infty}}}\right|\lesssim (s-r_c)^{-1}r_c^4=\rho_{5_{\infty}}(c,s)^{-1},
										\end{align}
										where the weights  defined  in \eqref{def:rho-1} write
										\begin{align}\label{recall:rho-6} \rho(c,\xi,s)=
											\left\{
											\begin{array}{l}
												\rho_{6_0}=\rho_{6_1}=\rho_{6_2}=r_c^{-3},\;\;\; \;\;\quad s\leq \f{r_c}{2},\\
												\rho_{6_3}=\f{r_c-s}{r_c^4}, \;\;\; \;\quad r_c/2\leq s\leq r_c-C\xi^{-\f23}, \\
												\rho_{6_4}=\xi^{-\f23}r_c^{-\frac{11}{3}},\;\;\; \quad r_c-C\xi^{-\f23}r_c^{\f13}\leq s\leq r_c+C\xi^{-\f23}r_c^{\f13}, \\
												\rho_{6_{\infty}}=\f{s-r_c}{r_c^4} \;\;\; \quad s\geq r_c+C\xi^{-\f23}r_c^{\f13},
											\end{array}\right.
										\end{align}
										and
										\begin{align}\label{def:sigma6}
											\sigma(c,s)=\left\{
											\begin{array}{l}  \sigma_{6_0}= \sigma_{6_1}=r_c^{3-\delta/2}s^{-\delta}\textbf{1}_{s\leq \f12}(s),\quad\quad
												s\leq 1/2,\\
												\sigma_{6_2}=r_c^{3-\delta/12}s^{\delta/12}\textbf{1}_{\f12\leq s\leq r_c/2}(s),\quad\quad
												\f14\leq s\leq r_c/2,\\
												\sigma_{6_3}=\sigma_{6_4}=  |r_c-s|^{-1+\delta/6}r_c^{4-\delta/6}\textbf{1}_{r_c/2\leq s\leq 2r_c}(s),\;\;\;\quad r_c/4\leq s\leq 2r_c,\\
												\sigma_{6_{\infty}}=(s-r_c)^{-1+\delta/6}r_c^{4-\delta/6}\textbf{1}_{
													s\geq r_c}(s),
												\quad\quad s\geq r_c+C\xi^{-\f23}r_c^{\f13},\\
												\sigma_{6_{\infty}'}=(s-r_c)^{\delta}r_c^{3-\delta}\textbf{1}_{
													s\geq r_c}(s),
												\quad\quad s\geq r_c+C\xi^{-\f23}r_c^{\f13},
											\end{array}\right.
										\end{align}
										which have the uniform bounds in $s$ that
										\begin{align}\label{int:sigma6}
											\begin{split} & \int^{1}_{1-\delta}\sigma_j(c,s,\xi)dc\lesssim s^{-\delta},\;\;j=6_0,6_1;\;\;
												\int^{1}_{1-\delta}\sigma_j(c,s)dc\lesssim 1,\;\;j=6_2,6_3,6_4,6_{\infty},6_{\infty}'.
											\end{split}
										\end{align}
										Under the  coordinates \eqref{trans:eta6}, using the facts that
										$(1-c)^{-\f13}\sim r_c$, $ \xi \int_{r_c}^s\sqrt{\f{V(s')-c}{1-c}}ds'\sim \xi (s-r_c)^{\f32}/\langle r\rangle^{\f12}$,   for $\xi r_c \gtrsim M$, $c\in(1-\delta,1)$,
										the multiplier $\mathrm{m}$ in \eqref{def:m-eta} behaves as
										\begin{align}\label{rewriten:m-eta6}
											&\quad\mathrm{m}(\xi,c,s)=
											\left\{
											\begin{array}{l}
												\mathrm{m}_{6_0}(\xi,c,s)=O_{\eta_{6_0}}^{\eta_{6_0}}
												\left(\eta_{6_0}^{\delta}\right),\\
												\mathrm{m}_{6_1}(\xi,c,s)=O_{\eta_{6_0}}^{\eta_{6_0}}
												\left(\eta_{6_0}^{\delta}\right),\\
												\mathrm{m}_{6_2}(\xi,c,s)=(\f{s}{r_c})^{-\f{\delta}{12}}
												O_{\eta_{6_2}}^{\eta_{6_2}}
												(1),\\
												\mathrm{m}_{6_3}(\xi,c,s)=O_{\eta_{6_3}}^{\eta_{6_3}}
												\left(\eta_{6_3}^{\delta/9}\right),\\
												\mathrm{m}_{6_4}(\xi,c,s)=O_{\eta_{6_3}}^{\eta_{6_3}}
												\left(\eta_{6_3}^{\delta/9}\right),\\
												\mathrm{m}_{6_{\infty}}(\xi,c,s)
												=O_{\eta_{6_{\infty}}}^{\eta_{6_{\infty}}}
												\left(\eta_{6_{\infty}}^{-\delta/9}\cdot
												(\xi r_c)^{\delta/9}\right). \end{array}\right.
										\end{align}
										Due to the oscillation behavior of $\widetilde{\phi}_{6_{\infty}}(\xi,c,r) $  and   $\widetilde{\phi}_{6_{\infty}}(\xi,c,s) $, the proof of  \eqref{op-pieces6} will be divided into four cases.
										
										 \textit{Case 1. }  $(j,j')\in \{6_0,6_1,6_2,6_3,6_4\}\times \{6_0,6_1,6_2,6_3,6_4\}$.
										
										We first apply the change of variable  \eqref{trans:eta6} to $\tilde{\phi}_j$ defined in \eqref{def:tphi-6}, in the $\xi$-integral that for $j\in \{6_0,6_1,6_2,6_3,6_4\}$,
										\begin{align}
											\label{26:1}&\phi_j(\eta_{j},c,s)
											:=\widetilde{\phi}_j\left(\xi(\eta_{j},c,s),s,c\right),
											\end{align}
										where $\xi(\eta_j,c,s)$ are the corresponding inverse functions of $\eta_j(\xi,c,s)$. For  $j',j\in \{6_0,6_1,6_2,6_3,6_4\},$
										\begin{align}
											\begin{split}
												&\phi_{j'}^*(\eta_{j},c,s,r):=\widetilde{\phi}_{j'}
												\label{26:2}\left(\xi(\eta_{j},c,s),r,c\right). \end{split}
											\end{align}
											The following bounds are deduced in the range  $c\in(1-\delta,1)$, $|\xi| r_c\gtrsim M\gg1$, and use the definitions of $\tilde{\phi}_j$ in \eqref{def:tphi-6},  the decomposition \eqref{ineq:trho} and the bounds of $\mathrm{m}_j$ in \eqref{rewriten:m-eta6}. Using the formula \eqref{fm:f-tf1} with $\rho_j^{-1}=r_c^3$, $\f{\pa_c\eta_j}{\eta_j}=O(r_c^3)$, $j=6_0,6_1$, and the bounds \eqref{6-phi-1}- \eqref{6-phi-2} in Proposition \ref{prop:summery} ($1\lesssim r_c\sim (1-c)^{-\f13}$),
											we obtain
												\begin{align}  \label{bd:60-61}
													\begin{split}
														\phi_{6_0}(\eta_{6_0},c,s)&=
														\chi\left(\xi s/(1-c)^{\f12}M\right)\chi\left(\xi r_c^{\f32}s/M^{\f12}\right) O_{\eta_{6_0},c}^{\eta_{6_0},\rho_6}(\eta_{6_0})\\
														&=\chi(\eta_{6_0}/M)
														O_{\eta_{6_0},c}^{\eta_{6_0},(\sigma_{6_0}\mathrm{m}_{6_0})^{-1}}
														(\eta_{6_0})\\
														&=\chi(\eta_{6_0}/M)
														O_{\eta_{6_0},c}^{\eta_{6_0},\sigma_{6_0}^{-1}}
														(\eta_{6_0}),\\
														\phi_{6_1}(\eta_{6_1},c,s)&=\left(1-\chi\left(\xi s/(1-c)^{\f12}\right) \left(\chi_+(4s)-\chi\left(\xi r_c^{\f12}s/M^{\f12}\right) \right)
														\right)O_{\eta_{6_1},c}^{\eta_{6_1},\rho_6}\left(\eta_{6_1}^{-\f12}\right)\\
														&=\left(1-\chi(\eta_{6_1}) \right) O_{\eta_{6_1},c}^{\eta_{6_1},(\sigma_{6_1}\mathrm{m}_{6_1})^{-1}}
														\left(\eta_{6_1}^{-\f12}\right),\\
														&=\left(1-\chi(\eta_{6_1}) \right) O_{\eta_{6_1},c}^{\eta_{6_1},\sigma_{6_1}^{-1}}
														\left(\eta_{6_1}^{-\f12+\delta/6}\right).
													\end{split}
												\end{align}
Using the formula \eqref{fm:f-tf1} with $\rho_{6_2}^{-1}=r_c^3$, $\f{\pa_c\eta_{6_2}}{\eta_{6_2}}=O(r_c^3)$,  and  the bounds \eqref{6-phi-3} in Proposition \ref{prop:summery}, we deduce that  for $r_c\geq \f12$, $|\xi| r_c\gtrsim M\gg1$,
																						\begin{align}  \label{bd:62}
													\begin{split}
														\phi_{6_2}(\xi,c,s)&= \left(1-\chi\left(\xi r_c/M^{\f12}\right)\right)
														\left(\chi_+(4s/r_c)-\chi_+(4s)\right)\left(\f{r_c}{s}\right)^{-\f14}
														O_{\eta_{6_2},c}^{\eta_{6_2},\rho_{6_2}}
														\left(\eta_{6_2}^{-\f12}\right)\\
														&= \left(1-\chi\left(\eta_{6_2}/M^{\f12}\right)\right)
														\left(\left(\f{r_c}{s}\right)^{-\f14}\textbf{1}_{\f14\leq s\leq r_c/2} (s)\right)
														O_{\eta_{6_2},c}^{\eta_{6_2},(\sigma_{6_2}\mathrm{m}_{6_2})^{-1}}
														\left(\eta_{6_2}^{-\f12}\right)\\
														&= \left(1-\chi\left(\eta_{6_2}/M^{\f12}\right)\right)
														\left(\left(\f{r_c}{s}\right)^{-\f14+\delta/12}\textbf{1}_{\f14\leq s\leq r_c/2} (s)\right)
														O_{\eta_{6_2},c}^{\eta_{6_2},\sigma_{6_2}^{-1}}
														\left(\eta_{6_2}^{-\f12}\right)\\
														&= \left(1-\chi\left(\eta_{6_2}/M^{\f12}\right)\right)
														O_{\eta_{6_2},c}^{\eta_{6_2},\sigma_{6_2}^{-1}}
														\left(\eta_{6_2}^{-\f12}\right).\end{split}
												\end{align}
											
											Using the formula \eqref{fm:f-tf1} with $\rho_{6_3}^{-1}=(r_c-s)^{-1}r_c^4$, $\f{\pa_c\eta_{6_3}}{\eta_{6_3}}=O(r_c^3)=O\left(\rho_{6_3}^{-1}\right)$, and $\rho_{6_4}^{-1}=\xi^{\f23}r_c^{\f{11}{3}}$, $\f{\pa_c\eta_{6_4}}{\eta_{6_4}}=O(r_c^3)=O\left(\rho_{6_4}^{-1}\right)$,  along with the bounds \eqref{6-phi-4}, \eqref{6-f+a1} in Proposition \ref{prop:summery},
											we deduce that  for $r_c\gtrsim 1$ and $|\xi| r_c\gtrsim M\gg1$,
											\begin{small}
												\begin{align} \begin{split}
														&\quad\phi_{6_3}(\eta_{6_3},c,s)\\ &=
														\left(1-\chi\left(\xi r_c/M^{\f12}\right)\right)\left(1- \chi_+(4r/r_c)\right)\left(1-\chi_+\left(\xi^{\f23} (r_c-r)/r_c^{\f13}\right)\right)O_{\xi,c}^{\xi,\rho_{6_3}}\left(\left(\xi (r_c-s)^{\f32}/r_c^{\f12}\right)^{-\f16}\cdot(\xi r_c)^{-\f13}\right)\\
														&= \left(1-\chi\left(\eta_{6_3}/M^{\f12}\right)\right) O_{\eta_{6_3},c}^{\eta_{6_3},\rho_{6_3}}\left(\eta_{6_3}^{-\f13}\right)\\
														&= \left(1-\chi\left(\eta_{6_3}/M^{\f12}\right)\right) O_{\eta_{6_3},c}^{\eta_{6_3},(\sigma_{6_3}\mathrm{m}_{6_3})^{-1}}\left(\eta_{6_3}^{-\f13}\right)\\ \label{bd:63-64}
														&= \left(1-\chi\left(\eta_{6_3}/M^{\f12}\right)\right) O_{\eta_{6_3},c}^{\eta_{6_3},\sigma_{6_3}^{-1}}
														\left(\eta_{6_3}^{-\f13+\delta/9}\right), \\
														&\quad\phi_{6_4}(\eta_{6_4},c,s)\\
														&=
														\left(1-\chi\left(\xi r_c/M^{\f12}\right)\right)
														\left(1- \chi_+(4r/r_c)\right)\chi_+\left(\xi^{\f23} (r_c-r)/r_c^{\f13}\right)\chi\left(\xi^{\f23}(r-r_c)/r_c^{\f13}\right) O_{\eta_{6_4},c}^{\eta_{6_4},\rho_{6_4}}\left(\eta_{6_4}^{-\f13}\right)\\
														&= \left(1-\chi\left(\eta_{6_4}/M^{\f12}\right)\right) O_{\eta_{6_4},c}^{\eta_{6_4},(\sigma_{6_4}\mathrm{m}_{6_4})^{-1}}
														\left(\eta_{6_4}^{-\f13}\right)\\
														&= \left(1-\chi\left(\eta_{6_4}/M^{\f12}\right)\right) O_{\eta_{6_4},c}^{\eta_{6_4},\sigma_{6_4}^{-1}}
														\left(\eta_{6_4}^{-\f13+\delta/9}\right).
														\end{split}
												\end{align}
											\end{small}
											We also have
												\begin{align}  \label{bd:6-rho-1}
												\begin{split}
														&\rho_{6_0}(c)^{-1}\phi_{6_0}(\eta_{6_0},c,s)
														=\chi\left(\eta_{6_0}/M\right)
														O_{\eta_{6_0}}^{\eta_{6_0}}
														(\eta_{6_0})\sigma_{6_0}(c,s),\\
														&\rho_{6_1}(c)^{-1}\phi_{6_1}(\eta_{6_1},c,s)
														=\left(1-\chi(\eta_{6_1}) \right) O_{\eta_{6_1}}^{\eta_{6_1}}
														\left(\eta_{6_1}^{-\f12+\delta/6}\right)\sigma_{6_1}(c,s),\\
														&\rho_{6_2}(c,s,\xi)^{-1} \phi_{6_2}(\eta_{6_2},c,s)=
														\left(1-\chi\left(\eta_{6_2}/M^{\f12}\right)\right) O_{\eta_{6_2}}^{\eta_{6_2}}
														\left(\eta_{6_2}^{-\f12}\right)\sigma_{6_2}(c,s), \\
														&\rho_{6_3}(c,s,\xi)^{-1} \phi_{6_3}(\eta_{6_3},c,s)=
														\left(1-\chi\left(\eta_{6_3}/M^{\f12}\right)\right) O_{\eta_{6_3}}^{\eta_{6_3}}
														\left(\eta_{6_3}^{-\f13+\delta/9}\right)\sigma_{6_3}(c,s), \\
														&\rho_{6_4}(c,s,\xi)^{-1} \phi_{6_4}(\eta_{6_4},c,s)=
														\left(1-\chi\left(\eta_{6_4}/M^{\f12}\right)\right) O_{\eta_{6_4}}^{\eta_{6_4}}
														\left(\eta_{6_4}^{-\f13+\delta/9}\right)\sigma_{6_4}(c,s), \\
														&r_c^{3}\phi_{j}(\eta_{j},c,s)=\left(1-\chi\left(\eta_{j}/M^{\f12}\right)\right) O_{\eta_{j}}^{\eta_{j}}
														\left(\eta_{j}^{-\f13+\delta/9}\right)\sigma_{j}(c,s) ,\;\;j=6_3,6_4,
														\end{split}
												\end{align}
											where we notice that $r_c^{3}=\rho_{j}(c)^{-1}$, $j=6_0,6_1,6_2$.
											
											For $\phi_{j'}^*$, using \eqref{useful:trans-5} and similar facts as before, such as  the definition in \eqref{def:tphi-6}, the coordinates \eqref{trans:eta6},  the bounds in Proposition \ref{prop:summery}   and  the decomposition \eqref{ineq:trho}, we  deduce that for $j',j\in \big\{6_0,6_1,6_2,6_3, 6_4\big\}$,
											\begin{align}&  \label{bi-bd:60-64}
												\begin{split}
													&\phi_{j'}^*(\eta_{j},c,s,r)= O_{\eta_{j},c}^{\eta_{j},r_c^{-3}}(1)+ O_{\eta_{j},c}^{\eta_{j},\rho_{j}(c,s,\xi)}(1)
													,\quad j'=6_0,6_1, 6_2,\\
													&\phi_{j'}^*(\eta_{j},c,s,r)= O_{\eta_{j},c}^{\eta_{j},\sigma_{j'}(c,r)^{-1}}(1)
													+O_{\eta_{j},c}^{\eta_{j},\rho_j(c,s,\xi)}(1),\quad j'=6_3,6_4. \end{split}
											\end{align}
											Using the formula \eqref{fm:f-tf1} and the bound in \eqref{est:paceta/eta6}, it follows directly that  \begin{align} \label{behave:chi6}
												\chi\left(\f{(1-c)^{\f12}}{\xi}\right)\left(1- \chi\left(\xi/(1-c)^{\f13}M\right) \right)
												=O_{\eta_j,c}^{\eta_j,r_c^{-3}}(1),\;\;j\in \{6_0,6_1,6_2,6_3, 6_4\}.
											\end{align}
											Therefore, using \eqref{bd:60-61}, \eqref{bd:62}, \eqref{bd:63-64}, \eqref{bd:6-rho-1} and \eqref{bi-bd:60-64}, we can reduce
											\eqref{op-pieces6} to
											\begin{align} I_{j,j'}:=\int_{V(0)}^{1-\delta}\left|\pa_c\left(p.v.\int_{\mathbb{R}} a_{j,j'}(\eta_j,c,\lambda)\f{e^{i\eta_j y_j(c,\lambda) }}{\eta_j}d\eta_j\right)\right|dc\lesssim s^{-\delta}+1\label{op-pieces6-new-case1},
											\end{align}
											where   $\lambda=(s,r,z)$  and
											\begin{align} \label{def:a_j,j'6}
												\begin{split}
													a_{j,j'}(\eta_j,c,\lambda)&:=\chi\left(\f{(1-c)^{\f12}}{\xi}\right)\left(1- \chi\left(\xi/(1-c)^{\f13}M\right)\right)
													\phi_{{j}}(\eta_{{j}},,c,s)\phi_{{j'}}(\eta_{{j'}},,c,s,r)^*\\
													&=
													\left\{
													\begin{array}{l}
														\chi\left(\eta_{6_0}/M\right)
														O_{\eta_{6_0},c}^{\eta_{6_0},\sigma_{6_0,j'}^{-1}}
														(\eta_{6_0}),\;\;\;\quad\quad\quad\quad\quad\quad\quad j=6_0, \\
														\left(1-\chi\left(\eta_{6_1}\right)\right) O_{\eta_{6_1},c}^{\eta_{6_1},\sigma_{6_1,j'}^{-1}}
														\left(\eta_{6_1}^{-\f12+\delta/6}\right),
														\;\;\;\quad\quad \quad\quad  j=6_1,\\
														\left(1-\chi\left(\eta_{6_2}/M^{\f12}\right)\right)
														O_{\eta_{6_2},c}^{\eta_{6_2,j'},\sigma_{6_2}^{-1}}
														\left(\eta_{6_2}^{-\f12}\right),\;\quad\quad\quad\quad j=6_2,\\
														\left(1-\chi\left(\eta_{6_4}/M^{\f12}\right)\right) O_{\eta_{6_4},c}^{\eta_{6_4},\sigma_{6_4,j'}^{-1}}
														\left(\eta_{6_4}^{-\f13+\delta/9}\right),\;\quad\quad\quad j=6_3,\\
														\left(1-\chi\left(\eta_{6_4}/M^{\f12}\right)\right) O_{\eta_{6_4},c}^{\eta_{6_4},\sigma_{6_4,j'}^{-1}}
														\left(\eta_{6_4}^{-\f13+\delta/9}\right),\;\quad\quad\quad j=6_4, \end{array}\right.
												\end{split}
											\end{align}
											with  \begin{align} \label{def:sigma_j,j'6}
												\begin{split}
													\sigma_{j,j'}(c,s,r)= \left\{
													\begin{array}{l}
														\sigma_{j}(c,s)+\sigma_{j'}(c,s) ,\;j'=6_0,6_1,6_2,\;\;j=6_0,6_1,...,6_4,\\ \sigma_{j}(c,s)+\sigma_{j'}(c,r) ,\;\;j'=6_3, 6_4,\;\; \;\; j=6_0,6_1,...,6_4, \end{array}\right.
												\end{split}
											\end{align}
											where $\sigma_{j}$ are as in \eqref{def:sigma6}. It follows from \eqref{int:sigma6} that the following uniform (in $s,r$) bounds hold, \begin{align} \label{sigma-6infty1}
												\int_{1-\delta}^1\sigma_{j,j'}(c,s,r)dc\lesssim  s^{-\delta},\;\;j=6_0,6_1;\quad \int_{1-\delta}^1\sigma_{j,j'}(c,s,r)dc\lesssim  1,\;\;j=6_2,6_3,6_4,
											\end{align}
											and
											\begin{align} \label{def:y_j61}
												\begin{split}
													y_j(c,z,s,r)= \left\{
													\begin{array}{l} \sqrt{\f{c}{1-c}}z/\f{  s}{(1-c)^{\f12}} ,\;\;\;\quad\quad\quad\quad\quad\quad\quad\quad\quad\quad j=6_0,6_1,\\
														\sqrt{\f{c}{1-c}}z/\f{ 1 }{(1-c)^{\f13}},\;\;\;\quad\quad\quad\quad\quad\quad\quad\quad \quad\quad j=6_2,6_3,6_4. \end{array}\right.
												\end{split}
											\end{align} Indeed, applying Lemma \ref{lem: pifi-Linfty} to the first term and the transform $c\to y_j$(via the Area formula) for the second term, we have
											\begin{align}
												\begin{split}
													\label{process:integra61}I_{j,j'}&\leq \int_{V(0)}^{1-\delta}\sigma_{j,j'}(c,s,r)\left|\left(p.v.\int_{\mathbb{R}} \left(\sigma_{j,j'}^{-1}\pa_c\right)a_{j,j'}(\eta_j,c,\lambda)\f{e^{i\eta_j y_j(c,\lambda) }}{\eta_j}d\eta_j\right)\right|dc\\
													&\quad+
													\int_{V(0)}^{1-\delta}\left|\left(p.v.\int_{\mathbb{R}} a_{j,j'}(\eta_j,c,\lambda)e^{i\eta_j y_j(c,\lambda) }d\eta_j\right)\right|\cdot \left|\f{\pa y_j(c,\lambda)}{\pa c}\right|dc\\
													&\lesssim  \int_{V(0)}^{1-\delta}\sigma_{j,j'}(c,s,r)\left\|\int_{\mathbb{R}} \left(\sigma_{j,j'}^{-1}\pa_c\right)a_{j,j'}(\eta_j,c,\lambda)e^{i\eta_j x}d\eta_j\right\|_{L^1_x(\mathbb{R})}dc\\
													&\quad+
													\int_{\mathbb{R}}\left|\int_{\mathbb{R}} a_{j,j'}(\eta_j,c,\lambda)e^{i\eta_j y_j }d\eta_j\right|d y_j,
											\end{split}\end{align}
											while applying the one dimensional Area formula, we have used that for fixed $\lambda=(z,s,r)$ and $x$,
											the algebraic equation $y_j(c,\lambda)=x$ has finite number(independent of $\lambda$, $x$) roots $c$, due to the explicit form in \eqref{def:y_j61}
											and the Fundamental theorem of algebra.  Then we apply  Lemma \ref{lem:symbol-chi} and the bounds \eqref{def:a_j,j'6} to bound the integrals as
											\begin{align*}
												\lesssim  \int_{V(0)}^{1-\delta}\sigma_{j,j'}(c,s,r)dc+1.
												\end{align*}
											Using \eqref{int:sigma6},  we obtain the desired bound \eqref{op-pieces6-new-case1}.
											
											 \textit{Case 2. }  $(j,j')\in \{6_{\infty}\}\times \{6_0,6_1,6_2,6_3, 6_4\}$ and $(j,j')\in  \{6_0,6_1,6_2,6_3, 6_4\}\times\{6_{\infty}\}$.\smallskip
											
											We only prove the first case. The case where $(j,j')\in   \{6_0,6_1,6_2,6_3, 6_4\}\times \{6_{\infty}\}$ is symmetric  to that where $(j,j')\in \{6_{\infty}\}\times \{6_0,6_1,6_2,6_3, 6_4\}$;
											we merely need to swap  the indices $j$ and $j'$ and the variable $r$ and $s$, and the proof is identical.
											
											We first apply the change of variable   \eqref{trans:eta6} in the $\xi$-integral that
											\begin{align}
												\label{26:3}&\phi_{6_{\infty}}(\eta_{6_{\infty}},c,s)
												:=\widetilde{\phi_{6_{\infty}}}\left(\left(\int_{r_c}^s
												\sqrt{\f{V(s')-c}{1-c}}ds'\right)^{-1}\eta_{6_{\infty}},s,c\right),
												\end{align}
											where $\widetilde{\phi_{6_{\infty}}}$ satisfies $\widetilde{\phi}_{6_{\infty}}
											=\widetilde{\phi_{6_{\infty}}} e^{i\xi \int_{r_c}^s\sqrt{\f{V(s')-c}{1-c}}ds'}
											+c.c$  as  in  \eqref{def:tphi-6}, $\widetilde{\phi}_{6_{\infty}}$ is the function in the target integral \eqref{op-pieces6}.
											For $j'\in \{6_0,6_1,6_2,6_3, 6_4\}$, we also take the transform
											\begin{align}
												\begin{split}
													&\phi_{j'}(\eta_{6_{\infty}},c,s,r):=\widetilde{\phi}_{j'}
													\left(\left(\int_{r_c}^s
													\sqrt{\f{V(s')-c}{1-c}}ds'\right)^{-1}\eta_{6_{\infty}},c,r\right).
												\end{split}
											\end{align}
											Then for $(j,j')\in (j,j')\in \{6_{\infty}\}\times \{6_0,6_1,6_2,6_3, 6_4\}$, we reduce \eqref{op-pieces6} to
											\begin{align} I_{6_{\infty},j'}:=\int_{1-\delta}^1
												\left|\pa_c\left(p.v.\int_{\mathbb{R}} a_{6_{\infty},j'}(\eta_{6_{\infty}},c,\lambda)\f{e^{i\eta_{6_{\infty}} (y_{6_{\infty}}(c,\lambda)\pm 1) }}{\eta_{6_{\infty}}}d\eta_{6_{\infty}}\right)\right|dc\lesssim 1\label{op-pieces6-new-case2},
											\end{align}
											where  $\lambda=(z,s,r)$  and
											\begin{align} \label{def:y_j62}
												\begin{split}
													y_{6_{\infty}}(c,\lambda)= \sqrt{\f{c}{1-c}}z/\int_{r_c}^s\sqrt{\f{V(s')-c}{1-c}}ds', \quad \text{monotonic}\;\text{in}\;\;c,\end{split}
											\end{align}
											and
											\begin{align} \label{6phi-case2}
												a_{6_{\infty},j'}(\eta_{6_{\infty}},c,\lambda)
												:=\chi\left(\f{(1-c)^{\f12}}{\xi}\right)\left(1- \chi\left(\xi/(1-c)^{\f13}M\right) \right)f(\eta_{6_{\infty}},c,s)f^*(\eta_{6_{\infty}},c,s,r),
											\end{align}
											with
											$\xi(\eta_{6_{\infty}},c,s)=\left(\int_{r_c}^s
											\sqrt{\f{V(s')-c}{1-c}}ds'\right)^{-1}\eta_{6_{\infty}}$ and $f\in \{\phi_{6_{\infty}}, \overline{\phi_{6_{\infty}}}\},\; \;f^*\in\{\phi_{j'}^*, \overline{\phi_{j'}}^*\}. $
											We claim that for $j'\in  \{6_0,6_1,6_2,6_3, 6_4\}$,
											\begin{align} \label{def:a_j,j'-62}
												\begin{split}
													a_{6_{\infty},j'}(\eta_{6_{\infty}},c,\lambda)=  \left(1-\chi\left(M^{\f12}\eta_{6_{\infty}}\right)\right) O_{\eta_{6_{\infty}},c}^{\eta_{6_{\infty}},\sigma_{6_{\infty},j'}
														(c,\lambda)
														^{-1}}\left(\eta_{6_{\infty}}^{-\f16}\right),
												\end{split}
											\end{align}
											where
											\begin{align} \label{def:sigma_j,j'62}
												\begin{split}
													\sigma_{6_{\infty},j'}(c,\lambda)= \left\{
													\begin{array}{l}
														\sigma_{6_{\infty}}(c,s)+\sigma_{6_{\infty}'}(c,s) ,\;j'=6_0,6_1, 6_2,\\
														\sigma_{6_{\infty}}(c,s)+\sigma_{6_{\infty}'}(c,s) +\sigma_{j'}(c,r) ,\;j'=6_3,6_4, \end{array}\right.
												\end{split}
											\end{align}
											where $\sigma_{6_{\infty}}$, $\sigma_{6_{\infty}'}$, $\sigma_{j'}$ are as in \eqref{def:sigma6}. And it follows from \eqref{int:sigma6} that  for $j'\in \{6_0,6_1,6_2,6_3, 6_4\}$, the uniform( in $\lambda)$ bounds hold,
											\begin{align} \label{sigma-6infty}
												\int_{1-\delta}^1\sigma_{6_{\infty},j'}(c,\lambda)dc\lesssim  1.
											\end{align}
											With \eqref{def:a_j,j'-62} at hand, \eqref{op-pieces6-new-case2} follows by a similar process as \eqref{process:integral}, using Lemma \ref{lem: pifi-Linfty}, the Area formula (with \eqref{def:y_j62}),  Lemma \ref{lem:symbol-chi}, and \eqref{sigma-6infty}.
											
											 It remains to prove \eqref{def:a_j,j'-62}. To this end, we use \eqref{observe:eta_infty}, and
											since
											\begin{align*}
												\notag& \xi \int_{r_c}^s\sqrt{\f{V(s')-c}{1-c}}ds'\sim \xi (s-r_c)^{\f32}/\langle s\rangle ^{\f12}\\
												&= \left\{
												\begin{array}{l}
													\xi (s-r_c)^{\f32}/\langle r_c\rangle ^{\f12}\cdot \left(\f{\langle r_c\rangle }{\langle s\rangle}\right)^{\f12},\quad \text{if}\quad  s\leq 2r_c,\\
													\left(\xi (s-r_c)^{\f32}/\langle r_c\rangle ^{\f12}\right)^{\f23}\cdot \left(\xi r_c^{\f32}/\langle r_c\rangle ^{\f12}\right)^{\f13}\cdot\left(\f{\langle r_c\rangle(s-r_c) }{r_c\langle s\rangle}\right)^{\f12},\quad \text{if}\quad  s\geq 2r_c, \end{array}\right. \end{align*}
											where we used  $\langle r_c\rangle(s-r_c) \geq \langle r_c\rangle s \geq r_c\langle s\rangle$ for $s\geq 2r_c$.
											Therefore, using the definition  of $\widetilde{\phi_{6_{\infty}}}$ in \eqref{def:tphi-6}, applying the identity \eqref{fm:f-tf1} on \eqref{26:3} with $\rho_6^{-1}=(s-r_c)^{-1}r_c^4$, using the bound \eqref{6-f+} in Proposition \ref{prop:summery} and the decomposition \eqref{ineq:trho} and the bound of $\mathrm{m}_{6_{\infty}}$ in \eqref{rewriten:m-eta6}, we obtain
											\begin{align}
												\label{bd:6inftycase2} \begin{split}
													&\quad\phi_{6_{\infty}}(\eta_{6_{\infty}},c,s)\\
													&=\left(1-\chi\left(M^{\f12}\xi \int_{r_c}^s\sqrt{\f{V(s')-c}{1-c}}ds'\right)\right)\cdot \left(1-\chi(\xi^{\f23}(s-r_c)//r_c^{\f13})\right) \left(1- \chi_+(4s/r_c)\right)
													\chi_+\left(\xi^{\f23}(r_c-s)/r_c^{\f13}\right)
													\\ &\quad\;\;O_{\eta_{6_{\infty}},c}^{\eta_{6_{\infty}},\rho_6}
													\left(\eta_{6_{\infty}}^{-\f16}\cdot (\xi r_c )^{-\f13}\right)
													\\
													&= \left(1-\chi\left(M^{\f12}\eta_{6_{\infty}}\right)\right)
													O_{\eta_{6_{\infty}},c}^{\eta_{6_{\infty}},\rho_6}
													\left(\eta_{6_{\infty}}^{-\f16}\cdot
													(\xi r_c)^{-\f13}\right)\\
													&= \left(1-\chi\left(M^{\f12}\eta_{6_{\infty}}\right)\right) O_{\eta_{6_{\infty}},c}^{\eta_{6_{\infty}},(\sigma_{6_{\infty}} (c,s)
														\mathrm{m}_{6_{\infty}})^{-1}}\left(\eta_{6_{\infty}}^{-\f16}\cdot (\xi r_c)^{-\f13}\right)\\
													&= \left(1-\chi\left(M^{\f12}\eta_{6_{\infty}}\right)\right) O_{\eta_{6_{\infty}},c}^{\eta_{6_{\infty}},\sigma_{6_{\infty}}(c,s)
														^{-1}}\left(\eta_{6_{\infty}}^{-\f16}\right),
											\end{split}\end{align}
											where the first cut-off function is $1$ due to \eqref{observe:eta_infty} provided that $M\gg 1$.
											Similarly, we have
											\begin{align}  \label{bd:6-rho-11case2}
												\begin{split}
													&\rho_6(c,s)^{-1}\phi_{6_{\infty}}(\eta_{6_{\infty}},c,s)
													=\left(1-\chi\left(M^{\f12}\eta_{6_{\infty}}\right)\right) O_{\eta_{6_{\infty}}}^{\eta_{6_{\infty}}}\left(\eta_{6_{\infty}}^{-\f16}\right) \sigma_{6_{\infty}}(c,s),\\
													&r_c^{3}\phi_{6_{\infty}}(\eta_{6_{\infty}},c,s)
													=\left(1-\chi\left(M^{\f12}\eta_{6_{\infty}}\right)\right) O_{\eta_{6_{\infty}}}^{\eta_{6_{\infty}}}(\eta_{6_{\infty}}
													^{-\f16})\sigma_{6_{\infty}'}(c,s),
												\end{split}
											\end{align}
											where $\sigma_{6_{\infty}}$, $\sigma_{6_{\infty}'}$ are defined in \eqref{def:sigma6}, and for the second one, we used the decomposition  by \eqref{observe:eta_infty}:
											\begin{align} \label{decom:rho6'}
												r_c^3=r_c^{3-\delta}(s-r_c)^{\delta}\cdot (\xi r_c)^{\f23\delta}
												\left(\xi(s-r_c)^{\f32}/r_c^{\f12}\right)^{-\f23\delta}= \sigma_{6_{\infty}'}(c,s)\cdot O_{\eta_{6_{\infty}}}^{\eta_{6_{\infty}}}\left(\eta_{6_{\infty}}^{-\delta}\cdot (\xi r_c)^{\f23\delta}\right).
											\end{align}
											For $\phi_{j'}^*$, using the formula \eqref{useful:trans-5} and similar facts, we can also obtain
											\begin{align} \label{bi-bd:6inftycase2}
												\begin{split}
													&\phi_{j'}^*(\eta_{6_{\infty}},c,s,r)= O_{\eta_{6_{\infty}},c}^{\eta_{6_{\infty}},\rho_{j'}(c)}(1)
													+O_{\eta_{6_{\infty}},c}^{\eta_{6_{\infty}},\rho_{6_{\infty}}(c,s)}(1),\quad j'=6_0,6_1,6_2\\
													&\phi_{j'}^*(\eta_{6_{\infty}},c,s,r)= O_{\eta_{6_{\infty}},c}^{\eta_{6_{\infty}},\sigma_{j'}(c,r)^{-1}}(1)
													+O_{\eta_{6_{\infty}},c}^{\eta_{5_{\infty}},\rho_{6_{\infty}}(c,s)}(1),\quad j'=6_3,6_4, \end{split}
											\end{align}
											noting that $\rho_{6_0}=\rho_{6_1}=\rho_{6_2}=r_c^3$ are independent of $r$, $s$.  And for $s\geq r_c+C\xi^{-\f23}r_c^{\f13}$, using the formula \eqref{fm:f-tf1} and the last bound in \eqref{est:paceta/eta6},  it follows directly that  \begin{align} \label{behave:chi6case2}
												\chi\left(\f{(1-c)^{\f12}}{\xi}\right)\left(1- \chi\left(\xi /M(1-c)^{\f13}\right) \right)
												=O_{\eta_{6_{\infty}},c}^{\eta_{6_{\infty}},r_c^{-3}}(1)
												+O_{\eta_{6_{\infty}},c}^{\eta_{6_{\infty}},\rho_{6_{\infty}}(c,s)}(1).
											\end{align}
											Then \eqref{def:a_j,j'-62} follows by  \eqref{bd:6inftycase2}, \eqref{bd:6-rho-11case2}, \eqref{bi-bd:6inftycase2} and        \eqref{behave:chi6case2}.\smallskip
											
											\textit{Case 3. }  $(j,j')\in  \{6_{\infty}\}\times\{6_{\infty}\}$.\smallskip
											
											WLOG, we assume $s\geq r>r_c$($r\geq s>r_c$ can be treated in a similar manner).
											We apply the change of variable $\eta_{6_{\infty}}=\xi \int_{r_c}^s
											\sqrt{\f{V(s')-c}{1-c}}ds'$(as in \eqref{trans:eta6}) that
											\begin{align*}
												&\phi_{6_{\infty}}(\eta_{6_{\infty}},c,s)
												:=\widetilde{\phi_{6_{\infty}}}
												\left((\int_{r_c}^s
												\sqrt{\f{V(s')-c}{1-c}}ds')^{-1}\eta_{6_{\infty}},s,c\right),\\
													&\phi_{6_{\infty}}(\eta_{6_{\infty}},c,s,r):=
													\widetilde{\phi_{6_{\infty}}}\left((\int_{r_c}^s
													\sqrt{\f{V(s')-c}{1-c}}ds')^{-1}\eta_{6_{\infty}},c,r\right).
																							\end{align*}
											In the following, we use the old notation $Q(r,c)=\f{V(r)-c}{1-c}$. By Lemma \ref{lem:behave-Q} and Lemma \ref{lem:behave-x}, we infer that for $s>r_c$,\begin{align}
												\begin{split}\label{Basic:Q6}
													&Q(s,c)\sim \f{s-r_c}{\langle s\rangle },\; \int_{r_c}^sQ(s',c)^{\f12}ds'\sim (s-r_c)^{\f32}/\langle s\rangle;\\ &\left|\f{\pa_cQ(s,c)}{Q(s,c)}\right|
													+\left|\f{\pa_c(\int_{r_c}^sQ(s',c)^{\f12}ds')}
													{\int_{r_c}^sQ(s',c)^{\f12}ds'}\right|\lesssim (s-r_c)^{-1}r_c^4 .
												\end{split}
											\end{align}
											Then by   \eqref{def:tphi-6}, we reduce \eqref{op-pieces6} to
									\begin{align}  \label{est:6-infty}
											I_{6_{\infty}, 6_{\infty}}:=\int_{1-\delta}^1\left|\pa_c\left(p.v.\int_{\mathbb{R}} a_{6_{\infty},6_{\infty}}(\eta_{6_{\infty}},c,\lambda)
											\f{e^{i\eta_{6_{\infty}}(y_{6_{\infty}}(c,\lambda)\pm 1)}} {\eta_{6_{\infty}}}d\eta_{6_{\infty}}\right)\right|dc\lesssim  1 ,
											\end{align}
											where   $\lambda=(z,s,r)$, $y_{6_{\infty}}$  as in  \eqref{def:y_j62} and
											\begin{small}
												\begin{align*}
													a_{6_{\infty},6_{\infty}}(\eta_{6_{\infty}},c,\lambda)
													:=\chi\left(\f{(1-c)^{\f12}}{\xi}\right)\left(1- \chi\left(\xi/M(1-c)^{\f13}\right) \right)f(\eta_{6_{\infty}},c,s) f^*(\eta_{6_{\infty}},c,s,r) e^{\pm i\f{ \int_{r_c}^rQ(s',c)^{\f12}ds'}{\int_{r_c}^sQ(s',c)^{\f12}ds'}},
												\end{align*}
											\end{small}
											with $f\in \{\phi_{6_{\infty}}, \overline{\phi_{5_{\infty}}}\}$, $f^*\in \{\phi_{6_{\infty}}^*, \overline{\phi_{6_{\infty}}}\}$ and
											$\xi(\eta_{6_{\infty}},c,s)=\left(\int_{r_c}^s
											\sqrt{\f{V(s')-c}{1-c}}ds'\right)^{-1}\eta_{6_{\infty}}$. We note that compared with \eqref{6phi-case2} in Case 2, there is an extra oscillation type term $e^{\pm i\f{ \int_{r_c}^rQ(s',c)^{\f12}ds'}{\int_{r_c}^sQ(s',c)^{\f12}ds'}}$.
											Using the formula \eqref{useful:trans-5}, the bound \eqref{6-f+} in Proposition \ref{prop:summery} and  \eqref{ineq:trho}, and the bounds of $\mathrm{m}_{6_{\infty}}$ in \eqref{rewriten:m-eta}, we have
											\begin{align} \begin{split} \label{bd:tphi6-case4'}
													&\phi_{6_{\infty}}(\eta_{6_{\infty}},c,s,r)= O_{\eta_{6_{\infty}},c}^{\eta_{6_{\infty}},
														\sigma_{6_{\infty}}(c,r)^{-1}}(1)
													+O_{\eta_{6_{\infty}},c}^{\eta_{6_{\infty}},\rho_{6_{\infty}}(c,s)}(1). \end{split}
											\end{align}
											Noticing by \eqref{Basic:Q6} and $s-r_c\geq r-r_c>0$ that
											\begin{align}
												\left|\pa_ce^{\pm i\f{ \int_{r_c}^rQ(s',c)^{\f12}ds'}{\int_{r_c}^sQ(s',c)^{\f12}ds'}} \right|\lesssim (s-r_c)^{-1}r_c^4.
												\end{align}
											Then, using the bounds \eqref{bd:6inftycase2}, \eqref{bd:6-rho-11case2}, \eqref{behave:chi6case2}, \eqref{bd:tphi6-case4'} and \eqref{int:sigma6},  we get
											\begin{align} \label{def:a_infty6}
												\begin{split}
													a_{6_{\infty}}(\eta_{6_{\infty}},c,\lambda)=
													\left(1-\chi\left(M^{\f12}\eta_{6_{\infty}}\right)\right) O_{\eta_{6_{\infty}},c}^{\eta_{6_{\infty}},\sigma_{6_{\infty},6_{\infty}}
														(c,\lambda)
														^{-1}}\left(\eta_{6_{\infty}}^{-\f16}\right), \end{split}
											\end{align}
											where
											\begin{align} \notag\sigma_{6_{\infty},6_{\infty}}(c,\lambda)&= \sigma_{6_{\infty}}(c,s)+\sigma_{5_{\infty}'}(c,s) +\sigma_{6_{\infty}}(c,r) ,\\
											&\;\text{with}\;\;\int_{1-\delta}^1\sigma_{6_{\infty},6_{\infty}}
												(c,\lambda)dc\lesssim  1,\label{bd:sigma_6infty}  \end{align}
											by \eqref{int:sigma6}, with $\sigma_{6_{\infty}}$, $\sigma_{6_{\infty}'}$ be as in \eqref{def:sigma6}.  To prove \eqref{est:6-infty},
											we take $\pa_c$ on $a$ and the oscillation, respectively.  Then we first apply Lemma \ref{lem: pifi-Linfty} to the first term
											and the transform $c\to y_{6_{\infty}}$ for the second  term, then use the bound \eqref{def:a_infty} and Lemma  \ref{lem:symbol-chi},
											and finally use \eqref{bd:sigma_6infty}, to obtain similar process as \eqref{process:integral-5case4},
											with $\int_{V(0)}^{1-\delta}$ replaced by $\int_{1-\delta}^1$ and $5_{\infty}$ replaced by $6_{\infty}$.
									\end{proof}
									
									\begin{proposition}\label{lem:Boundness-K-pieces7}
										Let  $\delta\ll 1$ and $M\gg 1$ be fixed. Let  $J^H_{7}=
										\{7_0,7_1,7_2,7_3, 7_{\infty}
										\}  $ and $\tilde{\phi}_j$($j\in J^H_{7}$) be defined in \eqref{def:tphi-7}, Definition \ref{def:decom-phi}.
										Then for $
										(j,j')\in J^H_{7}\times  J^H_{7}$, it holds uniformly for $(r,s,z)\in\mathbb{R}^+\times \mathbb{R}^+\times \mathbb{R}$  that
										\begin{small}
											\begin{align} \int_{1-\delta}^1\left|\pa_c\left(p.v.\int_{\mathbb{R}} \chi\left(\f{M^2(1-c)^{\f12}}{\xi}\right)\chi\left(\xi/(1-c)^{\f13}M\right)
												\widetilde{\phi}_{j}(\xi,c,s)\widetilde{\phi}_{j'}(\xi,c,r)\f{e^{i\xi \sqrt{\f{c}{1-c}}z}}{\xi}d\xi\right)\right|dc\lesssim s^{\delta}+1\label{op-pieces7}.
												\end{align}
											\end{small}
									\end{proposition}
									
									\begin{proof}
										 Without further illustration, $c\in(0,1)$, $ (1-c)^{\f12}\lesssim | \xi| \lesssim (1-c)^{\f13}$.
										Let the new coordinate be defined as follows
										\begin{align} \label{trans:eta7}\begin{split}
												\eta(\xi,c,s)= \left\{
												\begin{array}{l} \f{\xi }{(1-c)^{\f13}} := \eta_{7_0}=\eta_{7_1}=\eta_{7_2}\;\;\;\quad\quad\quad\quad s\lesssim \f{\xi^2}{1-c},\\
													\xi s:= 													\eta_{7_3}= \eta_{7_{\infty}}\;\;\quad\quad\quad\quad\quad\quad\quad \quad \;\; s\gtrsim \f{\xi^2}{1-c} . \end{array}\right.
											\end{split}
										\end{align}
										Recall  the weight in \eqref{def:rho-1} that $\rho(c,\xi,s)=1-c:= \rho_{7}$.  We take the transform such that
										\begin{align}
											\begin{split}
												\label{72:1}&\phi_{j}(\eta_{j},c,s)
												:=\widetilde{\phi}_{j}\left(\eta_{j}(1-c)^{\f13},c,s\right), \quad \phi_{j'}^*(\eta_{j},c,s,r):=\widetilde{\phi}_{j'}
												\left(\eta_{j}(1-c)^{\f13},c,r\right) \quad j=7_0,7_1,7_2,  \\
												&\phi_{j}(\eta_{j},c,s)
												:=\widetilde{\phi}_{j}\left(\eta_{j}s^{-1},c,s\right), \quad \phi_{j'}^*(\eta_{j},c,s,r):=\widetilde{\phi}_{j'}
												\left(\eta_{j}s^{-1},c,r\right) \quad j=7_3,7_{\infty},
										\end{split}\end{align}
										with $\widetilde{\phi}_{j}$ defined in \eqref{def:tphi-7}.
										Observing that
										\begin{align} \label{7-observe}
											\begin{split}
												&(1-c)^{-\delta/3}= \left(\f{\xi}{(1-c)^{\f12}}\right)^{2\delta}\cdot
												\left(\f{\xi}{(1-c)^{\f13}}\right)^{-2\delta}\quad\quad s\lesssim 1,\\
												&\left(\f{\xi s^{-\f12}}{(1-c)^{\f12}}\right)^{2\delta}=
												s^{\delta}\left(\f{\xi}{(1-c)^{\f12}}\right)^{2\delta}\quad\quad \quad\quad \quad 1\lesssim s\lesssim \f{\xi^2}{1-c},\\
												&(1-c)^{-\delta/3}= s^{\delta}(\xi s)^{-\delta}\cdot
												\left(\f{\xi }{(1-c)^{\f13}}\right)^{\delta}\quad\quad  s\gtrsim \f{\xi^2}{1-c}\gtrsim \f{\sqrt{1-c}}{\xi}, \end{split}
										\end{align}
										we get by \eqref{7-phi} and \eqref{7-f+} that
										\begin{align}
											\begin{split}\label{7-phi-use}
												&\widetilde{\phi}(s,\xi,c)= \left\{
												\begin{array}{l}
													O_{\xi,c}^{\xi,(1-c)^{1-\delta/3}}\left(\left(\f{\xi}{(1-c)^{\f13}}\right)^{3-2\delta}\right)\;\;\;\;
													s\lesssim M\xi^{-1}\sqrt{1-c}.\\
													O_{\xi,c}^{\xi,(1-c)^{1-\delta/3}}\left(\left(\f{\xi}{(1-c)^{\f13}}\right)^{3-2\delta}\right)\;\;\;\;\;\;\; M\xi^{-1}\sqrt{1-c}\lesssim s\lesssim 1.\\
													s^{\delta}O_{\xi,c}^{\xi,(1-c)^{1-\delta/3}}\left(\left(\f{\xi}{(1-c)^{\f13}}\right)^{3-2\delta}\right)\;\;\; 1\lesssim  s\lesssim M^{-\f72}\f{\xi^2}{1-c}.\\
													s^{\delta}O_{\xi,c}^{\xi,(1-c)^{1-\delta/3}}\left((\xi s)^{1-\delta}\right)\;\;\; M^{-\f72}\f{\xi^2}{1-c}\lesssim s\lesssim M^{-\f12}\xi^{-1}.
												\end{array}\right.  \\
												&(1-c)^{-1}\widetilde{\phi}(s,\xi,c)=(1-c)^{-1+\delta/3} \left\{
												\begin{array}{l}
													O_{\xi}^{\xi}\left(\left(\f{\xi}{(1-c)^{\f13}}\right)^{3-2\delta}\right)\;\;\;\;
													s\lesssim M\xi^{-1}\sqrt{1-c}.\\
													O_{\xi}^{\xi}\left(\left(\f{\xi}{(1-c)^{\f13}}\right)^{3-2\delta}\right)\;\;\;\;\;\;\; M\xi^{-1}\sqrt{1-c}\lesssim s\lesssim 1.\\
													s^{\delta}O_{\xi}^{\xi}\left(\left(\f{\xi}{(1-c)^{\f13}}\right)^{3-2\delta}\right)\;\;\; 1\lesssim  s\lesssim M^{-\f72}\f{\xi^2}{1-c}.\\
													s^{\delta}O_{\xi}^{\xi}\left((\xi s)^{1-\delta}\right)\;\;\; M^{-\f72}\f{\xi^2}{1-c}\lesssim s\lesssim M^{-\f12}\xi^{-1}.
												\end{array}\right.
											\end{split}
										\end{align}
										\begin{align}
											\begin{split}\label{7-f+-use}
												f_+(\xi,c,s)&=  s^{\delta}O_{\xi,c}^{\xi,(1-c)^{1-\delta/3}}\left((\xi s)^{-\f12}\right)e^{i\xi r}\quad \quad r\gtrsim M^{-\f12}\xi^{-1}. \\
												(1-c)^{-1}f_+(\xi,c,s)&= (1-c)^{-1+\delta/3} s^{\delta}O_{\xi}^{\xi}\left((\xi s)^{-\f12}\right)e^{i\xi r}\quad \quad s\gtrsim M^{-\f12}\xi^{-1}. \\
											\end{split}
										\end{align}
										
										We summarize the bounds for \eqref{72:1} as follows. Using the formula \eqref{fm:f-tf1}, the definitions of $\tilde{\phi}_{j}$ in \eqref{def:tphi-7},
										the bounds \eqref{7-phi-use}, \eqref{7-f+-use}, we obtain
										\begin{align}  \label{bd:70}
												\begin{split}
													&\pa_c\phi_{j}(\eta_{j},c,s),\;
													(1-c)^{-1}\phi_{j}(\eta_{j},c,s)=(1-c)^{-1+\delta/3}\chi\left(\eta_{j}/2M\right)
													O_{\eta_{j}}^{\eta_{j}}\left(\eta_{j}^{3-2\delta}\right)\quad j=7_0,7_1,\\
													&\pa_c\phi_{7_2}(\eta_{7_2},c,s),\;
													(1-c)^{-1}\phi_{7_2}(\eta_{7_2},c,s)=s^{\delta}(1-c)^{-1+\delta/3}
													\chi\left(\eta_{7_2}/2M\right)
													O_{\eta_{7_2}}^{\eta_{7_2}}\left(\eta_{7_2}^{3-2\delta}\right),\\
													&\pa_c\phi_{7_3}(\eta_{7_3},c,s),\;
													(1-c)^{-1}\phi_{7_3}(\eta_{7_3},c,s)=s^{\delta}(1-c)^{-1+\delta/3}
													\chi\left(\eta_{7_3}/M\right)
													O_{\eta_{7_3}}^{\eta_{7_3}}\left(\eta_{7_3}^{1-\delta}\right),
												\end{split}
											\end{align}
										\begin{small}
											\begin{align}  \label{bd:7infty}
												\begin{split}
													\pa_c\phi_{7_{\infty}}(\eta_{\infty},c,s),\;
													(1-c)^{-1}\phi_{7_{\infty}}(\eta_{7_{\infty}},c,s)
													&=\sum_{\iota\in\{-,+\}}s^{\delta}(1-c)^{-1+\delta/3}
													\left(1-\chi(\eta_{7_{\infty}})\right)
													O_{\eta_{7_{\infty}}}^{\eta_{7_{\infty}}}
													\left(\eta_{7_{\infty}}^{-\f12}\right)e^{\iota i\eta_{7_{\infty}}}.\end{split}
											\end{align}
										\end{small}
										For $\phi^*$, we use  the formula  \eqref{useful:trans-2}
										with $\rho=\rho_7=1-c$, $\f{\pa_c\eta_{j}}{\eta_{j}}=O\left((1-c)^{-1}\right)$. Then, using the bounds \eqref{7-phi-use}, \eqref{7-f+-use}  again, we  obtain
											\begin{align}&  \label{bi-bd:70}
												\begin{split}
													&\phi_{j'}^*(\eta_{j},c,s,r)= O_{\eta_{j},c}^{\eta_{j},1-c}(1)
													\quad\quad\quad\quad j'\neq 7_{\infty},\;j=7_0,...,7_{\infty},
													\end{split}
											\end{align}
											and
											\begin{align}  \label{bi-bd:7infty}
												\begin{split}
													&\phi_{7_{\infty}}^*(\eta_{j},c,s,r)= \sum_{\iota\in\{-,+\}}O_{\eta_{j},c}^{\eta_{j},1-c}(1)e^{\iota i\eta_j (1-c)^{\f13}r } \quad j'=7_{\infty},\;j=7_0,7_{1},7_2,\\
													&\phi_{7_{\infty}}^*(\eta_{j},c,s,r)= \sum_{\iota\in\{-,+\}}O_{\eta_{j},c}^{\eta_{j},1-c}(1)e^{\iota i\eta_j \f{r}{s} }
													\quad \quad\quad j'=7_{\infty},\;j=7_3,7_{\infty}.\end{split}
											\end{align}
																				
																For the cut-off function, we have
									\begin{align} \label{bd-cutoff:7infty}
												 \chi\left(\f{M^2(1-c)^{\f12}}{\xi}\right)\chi\left(\xi/(1-c)^{\f13}M\right) =O_{\eta_{j},c}^{\eta_{j},1-c}(1).
											\end{align}

										Now we are in a position to prove \eqref{op-pieces7}. Denote $\lambda=(z,s,r)$ as the parameter. For $(j,j')\in \{7_0,7_1,7_2,7_3\}\times  \{7_0,7_1,7_2,7_3\}$,
										we apply the change of variable $\xi\to \eta_{j}$ \eqref{trans:eta7} in the $\xi$-integral,  reducing \eqref{op-pieces7} to
										\begin{align}
											I_{j,j'}:=\int_{1-\delta}^{1}\left|\pa_c\left(p.v.\int_{\mathbb{R}} a_{j,j'}(\eta_{j},c,\lambda)\f{e^{i\eta_{j} y_{j}(c,\lambda) }}{\eta_{j}}d\eta_{j}\right)\right|dc\lesssim s^{\delta}+1\label{op-pieces7-new-case1},
											\end{align}
										with
										\begin{align} \label{def:y71}
											\begin{split}
												y_j(c,z,s,r)= \left\{
												\begin{array}{l} \sqrt{\f{c}{1-c}}z/(1-c)^{\f13} ,\;\;\;\quad\quad\quad\quad\quad\quad\quad j=7_0,7_1,7_2,\\
													\sqrt{\f{c}{1-c}}z/s,\;\;\;\quad\quad\quad\quad\quad\quad\quad\quad \quad\quad j=7_3. \end{array}\right.
											\end{split}
										\end{align}
										and \begin{align} \label{def:a_71}
											\begin{split}
												&\quad a_{j,j}(\eta_{1},c,\lambda):=
												\chi\left(\f{M^2(1-c)^{\f12}}{\xi}\right)\chi\left(\xi/(1-c)^{\f13}M\right) \phi_{j}(\eta_{j},c,s)
												\phi_{j'}^* (\eta_{j},c,s,r).
												\end{split}
										\end{align}
										It follows from  \eqref{bd:70}, \eqref{bi-bd:70} and \eqref{bd-cutoff:7infty}  that
										\begin{align} \label{bd:a71}
											\begin{split}
												&a_{j,j'}(\eta_j,c,\lambda)
												=
												\left\{
												\begin{array}{l}
													\chi\left(\eta_{7_0}/2M\right)
													O_{\eta_{7_0},c}^{\eta_{7_0},(1-c)^{1-\delta/3}}
													\left(\eta_{7_0}^{3-\delta}\right),\;\;\;\quad\quad j=7_0, \\
													\chi\left(\eta_{7_1}/2M\right)
													O_{\eta_{7_1},c}^{\eta_{7_1},(1-c)^{1-\delta/3}}
													\left(\eta_{7_1}^{3-\delta}\right),\;\;\;\quad\quad j=7_1, \\
													s^{\delta}\chi\left(\eta_{7_2}/2M\right)
													O_{\eta_{7_2},c}^{\eta_{7_2},(1-c)^{1-\delta/3}}
													\left(\eta_{7_2}^{3-\delta}\right),\;\;\;\quad\quad j=7_2, \\
													s^{\delta}\chi(\eta_{7_3}/M)
													O_{\eta_{7_3},c}^{\eta_{7_3},(1-c)^{1-\delta/3}}
													\left(\eta_{7_3}^{1-\delta}\right),\;\;\;\quad\quad j=7_3. \end{array}\right.
											\end{split}
										\end{align}
										Therefore, for the following two parts of $\pa_c$, we apply Lemma \ref{lem: pifi-Linfty} to the first term and  take the transform $c\to y_{j}$ for the second term, to deduce that for  $(j,j')\in \{7_0,7_1,7_2,7_3\}\times  \{7_0,7_1,7_2,7_3\}$,
										\begin{align}
											\begin{split}
												\label{process:integral70}I_{j,j'}&\leq \int_{1-\delta}^{1}(1-c)^{-1+\delta/3}\left|\left(p.v.\int_{\mathbb{R}}
												\left((1-c)^{1-\delta/3}\pa_c\right)a_{j,j'}(\eta_{j},c,\lambda)\f{e^{i\eta_{j} y_{j}(c,\lambda) }}{\eta_{j}}d\eta_{j}\right)\right|dc\\
												&\quad+
												\int_{1-\delta}^{1}\left|\left(p.v.\int_{\mathbb{R}} a_{j,j'}(\eta_{j},c,\lambda)e^{i\eta_{j} y_{j}(c,\lambda) }d\eta_{j}\right)\right|\cdot
												\left|\f{\pa y_{j}(c,\lambda)}{\pa c}\right|dc\\
												&\lesssim  \int_{1-\delta}^{1}(1-c)^{-1+\delta/3}
												\left\|\int_{\mathbb{R}} \left((1-c)^{1-\delta/3}\pa_c\right)a_{j,j'}(\eta_{j},c,\lambda)e^{i\eta_{j} x}d\eta_{j}\right\|_{L^1_x(\mathbb{R})}dc\\
												&\quad+
												\int_{\mathbb{R}}\left|\int_{\mathbb{R}} a_{j,j'}(\eta_{j},c,\lambda)e^{i\eta_{j} y_{j} }d\eta_{1}\right|d y_{j}.
										\end{split}\end{align}
										We then apply  Lemma \ref{lem:symbol-chi} and the estimates in \eqref{bd:a71} to bound the above integrals by
										\begin{align*}
											\lesssim  s^{\delta}+1,
											\end{align*}
										which yields  the desired bound \eqref{op-pieces7-new-case1}.
										
										For $(j,j')=\{7_{\infty}\}\times \{7_0,7_1,7_2,7_3\}$, we apply the change of variable  $\xi\to \eta_{7_{\infty}}=\xi s$ in the $\xi$-integral,  reducing  \eqref{op-pieces7} to
										\begin{align}
											 I_{7_{\infty},j'}:=\int_{1-\delta}^{1}\left|\pa_c\left(p.v.\int_{\mathbb{R}} a_{7_{\infty},j'}(\eta_{7_{\infty}},c,\lambda)\f{e^{i\eta_{7_{\infty}} y_{7_{\infty}}(c,\lambda)}}{\eta_{7_{\infty}}}d\eta_{7_{\infty}}\right)\right|dc\lesssim s^{\delta}+1\label{op-pieces7-new-case2},
											 \end{align}
										with $y_{7_{\infty}}=\sqrt{\f{c}{1-c}}z/s$, and by \eqref{bd:7infty}, \eqref{bi-bd:70} and \eqref{bd-cutoff:7infty},
										\begin{align} \label{def:a_7infty}
											\notag a_{7_{\infty},j'}(\eta_{7_{\infty}},c,\lambda)&:=
											\chi\left(\f{M^2(1-c)^{\f12}}{\xi}\right)\chi\left(\xi/(1-c)^{\f13}M\right) \phi(\eta_{7_{\infty}},c,s)
											\phi_{j'}^* (\eta_{7_{\infty}},c,s,r)\\
											&=\sum_{\iota\in\{-,+\}}s^{\delta}\left(1-\chi(\eta_{7_{\infty}}) \right) O_{\eta_{1},c}^{\eta_{1},(1-c)^{1-\delta/3}}\left(\eta_{7_{\infty}}^{-\f12}\right)
											e^{\iota i \eta_{7_{\infty}}}.
											\end{align}
											Therefore, for the following two parts of $\pa_c$, we apply Lemma \ref{lem: pifi-Linfty} to the first term
											and  take the transform $c\to y_{7_{\infty}}$ for the second term, to deduce the estimate in the same manner as in \eqref{process:integral70}.
											Then we apply  Lemma \ref{lem:symbol-chi} and the bound in \eqref{def:a_1infty} to estimate the integral by  $s^{\delta}+1$,
											which yields the desired bound \eqref{op-pieces7-new-case2}.
											
											 For $(j,j')=\{7_0,7_1,7_2,7_3\}\times \{7_{\infty}\}$, we apply the change of variable $\xi\to \eta_{7_{\infty}}=\xi r$ in the $\xi$-integral,
											 and the argument proceeds similarly for $(j,j')= \{7_{\infty}\}\times\{7_0,7_1,7_2,7_3\}$ by employing the bounds as follows
										\begin{align*}
													&\phi_{7_{\infty}}^*(\eta_{7_{\infty}},c,r)
													:=\widetilde{\phi}_{7_{\infty}}\left(\eta_{7_{\infty}}r^{-1},c,r\right)
													=\sum_{\iota\in\{-,+\}}
													\left(1-\chi(\eta_{7_{\infty}}) \right)
													O_{\eta_{7_{\infty}},c}^{\eta_{7_{\infty}},1-c}
													\left(\eta_{7_{\infty}}^{-\f12}\right)e^{\iota i\eta_{7_{\infty}}},\\
													&\pa_c\phi_{j}(\eta_{7_{\infty}},c,s,r),\;
													(1-c)^{-1}\phi_{j}(\eta_{7_{\infty}},c,s,r)
													=(1-c)^{-1+\delta/3} O_{\eta_{7_{\infty}}}^{\eta_{7_{\infty}}}(1)\quad j=7_0,7_1,\\
													&\pa_c\phi_{j}(\eta_{7_{\infty}},c,s,r),\;
													(1-c)^{-1}\phi_{j}(\eta_{7_{\infty}},c,s,r)
													=s^{\delta}(1-c)^{-1+\delta/3} O_{\eta_{7_{\infty}}}^{\eta_{7_{\infty}}}(1)\quad j=7_2,7_3,\\
													&\text{with}\quad  \phi_{j}(\eta_{7_{\infty}},c,s,r)
													:=\widetilde{\phi}_{j}\left(\eta_{7_{\infty}}r^{-1},c,s\right),\quad \widetilde{\phi}_{j}\quad \text{defined}\quad  \text{in}\quad  \eqref{def:tphi-7}.
												\end{align*}
																						We leave the details to the readers.
											
											 For $(j,j')=(7_{\infty},7_{\infty})$, we apply the change of variable  $\xi\to \eta_{7_{\infty}}=\xi s$ in the $\xi$-integral, reducing \eqref{op-pieces7}  to
											\begin{align}
												I_{7_{\infty},7_{\infty}}:=\int_{1-\delta}^{1}
												\left|\pa_c\left(\int_{\mathbb{R}} a_{7_{\infty},7_{\infty}}(\eta_{7_{\infty}},c,\lambda)
												\f{e^{i\eta_{7_{\infty}} (y_{7_{\infty}}(c,\lambda) \pm 1)}}{\eta_{7_{\infty}}}d\eta_{7_{\infty}}\right)\right|dc\lesssim s^{\delta}+1\label{op-pieces7-new-case3},
												\end{align}
											with $y_{7_{\infty}}=\sqrt{\f{c}{1-c}}z/s$, and by  \eqref{bd:7infty}, \eqref{bi-bd:7infty} and \eqref{bd-cutoff:7infty},
											\begin{align} \label{def:a_7inftyinfty}
												\notag & a_{7_{\infty},7_{\infty}}(\eta_{7_{\infty}},c,\lambda):=
												\chi\left(\f{M^2(1-c)^{\f12}}{\xi}\right)\chi\left(\xi/(1-c)^{\f13}M\right) \phi(\eta_{7_{\infty}},c,s)
												\phi_{j'}^* (\eta_{7_{\infty}},c,s,r)\\
												&=\sum_{\iota\in\{-,+\}}\sum_{\tau\in\{-,+\}}
												s^{\delta}\left(1-\chi(\eta_{7_{\infty}}) \right) O_{\eta_{1},c}^{\eta_{1},(1-c)^{1-\delta/3}}\left(\eta_{7_{\infty}}^{-\f12}\right)
												e^{\iota i \eta_{7_{\infty}}}e^{\tau i \eta_{7_{\infty}}\f{r}{s}}.
												\end{align}
											Therefore, we deduce \eqref{op-pieces7-new-case3} by a similar argument to that in \eqref{process:integral70}.
									\end{proof}

									\section{Uniform estimates of the integral kernel: Part II}
									
									In this section, we prove the uniform  estimate of the integral kernel \eqref{op-K-D} in Proposition \ref{prop:Boundness-K-D}. The proof follows a similar manner in the last section.

									\subsection{A decomposition of $\mathcal{D}_{l,m}\widetilde{\phi}$  by its behavior}
									By \eqref{phi-f+-linear}, we have  \begin{align} \label{phi-f+-linear-D}
										(\mathcal{D}\widetilde{\phi})(\xi,c,r) &=  \f{i\overline{W}(\xi,c)}{2|W(\xi,c)|}
										\left(\mathcal{D}f_+\right)(\xi,c,r)+
										\f{-iW(\xi,c)}{2|W(\xi,c)|}\left(\mathcal{D}\bar{f}_+\right)(\xi,c,r),
									\end{align}
									where by the definition \eqref{def:derivative'},
									\begin{align}  \label{def:derivative'}
										\mathcal{D}\in\left\{\mathcal{D}_{l,m}=\left(\f{(1-c)^{\f12}\pa_r}{\xi}\right)^{m}
										\left(\f{(1-c)^{\f12}(\pa_r+r^{-1})}{\xi}\right)^{l},\quad
										\mathcal{D}_{l,m}^L=\pa_r^{m}
										\left(\pa_r+r^{-1}\right)^{l}\right\}.
										\end{align}
										Recall  that smooth cut-off functions $\chi, \chi_+$ are defined in \eqref{def:chi-notation}.																				
										\begin{definition} \label{def:decom-phi-D}
											Let $\delta\ll 1$, $M\gg 1$ be fixed, $\mathcal{D}\in\{\mathcal{D}_{l,m},\mathcal{D}^L_{l,0}\}$ be as \eqref{def:derivative'}. According to the behaviors in Proposition \ref{prop:summery}, we decompose for $\xi\in\mathbb{R}/\{0\}$, $c\in(0,1)/\{V(0)\}$, $r>0$  such that
											\begin{small}
												\begin{align}
													\begin{split}
														& \left(\mathcal{D}^L_{l,m}\widetilde{\phi}\right)(\xi,c,r)=\sum_{j\in J^L}\left(\mathcal{D}_{l,m}^L\widetilde{\phi}\right)_{j},\quad J^L=\{1_0,1_{\infty}\},\;\text{if}\;\;0<|\xi|\lesssim (1-c)^{\f12},\\
														& \left(\mathcal{D}_{l,m}\widetilde{\phi}\right)(\xi,c,r)=\sum_{j\in J^H_{2}}\left(\mathcal{D}_{l,m}\widetilde{\phi}\right)_{j},\quad J^H_{2}
														=\{2_0,2_{\infty}\},\;\text{if}\;\;(1-c)^{\f12}\lesssim |\xi|\lesssim  (V(0)-c)^{-\f32},\;c\in(0,V(0)),\\
														& \left(\mathcal{D}_{l,m}\widetilde{\phi}\right)(\xi,c,r)=\sum_{j\in J^H_{3}}\left(\mathcal{D}_{l,m}\widetilde{\phi}\right)_{j},\quad J^H_{3}
														=\{3_0,3_{\infty}\},\;\text{if}\;\;|\xi|\gtrsim  (V(0)-c)^{-\f32},\;c\in(0,V(0)),\\& \left(\mathcal{D}_{l,m}\widetilde{\phi}\right)(\xi,c,r)=\sum_{j\in  J^H_{4}}\left(\mathcal{D}_{l,m}\widetilde{\phi}\right)_{j},\;\;J^H_{4}=\{
														4_0,4_{\infty}\},\;\text{if}\;\;(1-c)^{\f12}\lesssim |\xi|\lesssim (c-V(0))^{-\f32}  ,\;c\in(V(0),1-\delta),\\
														&\left(\mathcal{D}_{l,m}\widetilde{\phi}\right)(\xi,c,r)=\sum_{j\in J^H_{5}}\left(\mathcal{D}_{l,m}\widetilde{\phi}\right)_{j},\;\; J^H_{5}=\{5_0,5_1,5_2,5_3,5_{\infty}\},\;
														\text{if}\;\;|\xi|\gtrsim   M (c-V(0))^{-\f32},\;c\in(V(0),1-\delta),\\
														& \left(\mathcal{D}_{l,m}\widetilde{\phi}\right)(\xi,c,r)=\sum_{j\in  J^H_{6}}\left(\mathcal{D}_{l,m}\widetilde{\phi}\right)_{j} ,\;J^H_{6}=\{6_0,6_1,6_2,6_3,6_4,6_{\infty}\},\;
														\text{if}\;\;|\xi|\gtrsim   M(1-c)^{\f13},\;c\in(1-\delta,1),\\
														& \left(\mathcal{D}_{l,m}\widetilde{\phi}\right)(\xi,c,r)=\sum_{j\in  J^H_{7}}\left(\mathcal{D}_{l,m}\widetilde{\phi}\right)_{j} ,\; J^H_{7}
														=\{7_0,7_1,7_2,7_{\infty}\},\;
														\text{if}\;\;M^2(1-c)^{\f12}\lesssim |\xi|\lesssim  M (1-c)^{\f13},\;c\in(1-\delta,1), \end{split}
												\end{align}
											\end{small}
											where for $\xi<0$, $\left(\mathcal{D}\widetilde{\phi}\right)_j(\xi)=-\left(\mathcal{D}\widetilde{\phi}\right)_j(-\xi)$ and for $\xi>0$:
											\begin{align}
													\begin{split}\label{def:tphi-1-D}
														&\left(\mathcal{D}\widetilde{\phi}\right)_{1_0}(\xi,c,r)
														=\chi(\xi r)\left(\mathcal{D}\widetilde{\phi}\right),\\
														&\left(\mathcal{D}\widetilde{\phi}\right)_{1_{\infty}}(\xi,c,r)
														=\left(1-\chi(\xi r)\right)\f{i\overline{W}(\xi,c)}{2|W(\xi,c)|}
														\left(\mathcal{D}f_+\right)(\xi,c,r)+c.c,
														\end{split}
												\end{align}
											where $\mathcal{D}\in\{\mathcal{D}_{l,m},\mathcal{D}^L_{l,0}\}$ and
											\begin{align} \label{def:tphi-2-D}
													\begin{split}
														\left(\mathcal{D}_{l,m}\widetilde{\phi}\right)_{2_0}(\xi,c,r)
														&=\chi\left(\xi^{\f23} r\right) \left(\mathcal{D}_{l,m}\widetilde{\phi}\right),\\
														\left(\mathcal{D}_{l,m}\widetilde{\phi}\right)_{2_{\infty}}(\xi,c,r)
														&= \left(1-\chi\left(\xi^{\f23} r\right)\right)\f{i\overline{W}(\xi,c)}{2|W(\xi,c)|}
														\left(\mathcal{D}_{l,m}f_+\right)(\xi,c,r)+c.c,
														\end{split}
												\end{align}
											and
												\begin{align}  \label{def:tphi-3-D}
													\begin{split}
														\left(\mathcal{D}_{l,m}\widetilde{\phi}\right)_{3_0}(\xi,c,r)
														&=\chi\left(\xi (c-V(0))^{\f12} r\right) \left(\mathcal{D}_{l,m}\widetilde{\phi}\right),\\
														\left(\mathcal{D}_{l,m}\widetilde{\phi}\right)_{3_{\infty}}(\xi,c,r)
														&= \left(1-\chi\left(\xi (c-V(0))^{\f12} r\right)\right)\f{i\overline{W}(\xi,c)}{2|W(\xi,c)|}
														\left(\mathcal{D}_{l,m}f_+\right)(\xi,c,r)+c.c,\end{split}
												\end{align}
													and												
													\begin{align} \label{def:tphi-4-D}
													\begin{split}
														\left(\mathcal{D}_{l,m}\widetilde{\phi}\right)_{4_0}(\xi,c,r)
														&=\chi\left(\xi^{\f23} r/M\right)\left(\mathcal{D}_{l,m}\widetilde{\phi}\right)_{4_0},\\
														\left(\mathcal{D}_{l,m}\widetilde{\phi}\right)_{4_{\infty}}(\xi,c,r)
														&= \left(1-\chi\left(\xi^{\f23} r/M\right)\right)\f{i\overline{W}(\xi,c)}{2|W(\xi,c)|}
														\left(\mathcal{D}_{l,m}f_+\right)(\xi,c,r)+c.c,
												\end{split} \end{align}
											and
											\begin{small}
												\begin{align}\label{def:tphi-5-D}
													\begin{split}
														\left(\mathcal{D}_{l,m}\widetilde{\phi}\right)_{5_0}(\xi,c,r)
														&=
														\chi\left(\xi r_c^{\f12} r/M^{\f12}\right) \left(\mathcal{D}_{l,m}\widetilde{\phi}\right),\\  \left(\mathcal{D}_{l,m}\widetilde{\phi}\right)_{5_1}(\xi,c,r)
														&=
														\left(\chi_+(4r/r_c)-\chi\left(\xi r_c^{\f12} r/M^{\f12}\right)\right) \left(\mathcal{D}_{l,m}\widetilde{\phi}\right),\\
														(\mathcal{D}_{l,m}\widetilde{\phi})_{5_2}(\xi,c,r)
														&= \left(1- \chi_+(4r/r_c)\right)\left(1-\chi_+\left(\xi^{\f23} (r_c-r)\right)\right) \left(\mathcal{D}_{l,m}\widetilde{\phi}\right),\\
														\left(\mathcal{D}_{l,m}\widetilde{\phi}\right)_{5_{3}}(\xi,c,r)
														&= \left(1- \chi_+(4r/r_c)\right)\chi_+\left(\xi^{\f23} (r_c-r)\right)\chi\left(\xi^{\f23}(r-r_c)\right)\left(\f{i\overline{W}}{2|W|}f_+-
														\f{iW}{2|W|}\bar{f}_+\right),\\
														\left(\mathcal{D}_{l,m}\widetilde{\phi}\right)_{5_{\infty}}(\xi,c,r)
														&=\left(1- \chi_+(4r/r_c)\right)
														\chi_+\left(\xi^{\f23}(r_c-r)\right)
														\left(1-\chi\left(\xi^{\f23}(r-r_c)\right)\right)\f{i\overline{W}(\xi,c)}{2|W(\xi,c)|}
														\left(\mathcal{D}_{l,m}f_+\right)(\xi,c,r)\\
														&\quad+c.c,
														\end{split}
												\end{align}
											\end{small}
											and
											\begin{small}
												\begin{align}\label{def:tphi-6-D}
													\begin{split}
														&\left(\mathcal{D}_{l,m}\widetilde{\phi}\right)_{6_0}(\xi,c,r)
														=
														\chi\left(\xi r_c^{\f32} r/M^{\f12}\right)\left(\mathcal{D}_{l,m}\widetilde{\phi}\right),\\
														&\left(\mathcal{D}_{l,m}\widetilde{\phi}\right)_{6_1}(\xi,c,r)
														=
													        \left(\chi_+(4r)-\chi\left(\xi r_c^{\f32} r/M^{\f12}\right)\right)\left(\mathcal{D}_{l,m}\widetilde{\phi}\right),\\
														&\left(\mathcal{D}_{l,m}\widetilde{\phi}\right)_{6_2}(\xi,c,r)
														=
														\left(\chi_+(4r/r_c)-\chi_+(4r)\right)\left(\mathcal{D}_{l,m}\widetilde{\phi}\right),\\
														&\left(\mathcal{D}_{l,m}\widetilde{\phi}\right)_{6_3}(\xi,c,r)
														= \left(1- \chi_+(4r/r_c)\right)\left(1-\chi_+\left(\xi^{\f23} (r_c-r)/r_c^{\f13}\right)\right)\left(\mathcal{D}_{l,m}\widetilde{\phi}\right),\\
														&\left(\mathcal{D}_{l,m}\widetilde{\phi}\right)_{6_4}(\xi,c,r)
														= \left(1- \chi_+(4r/r_c)\right)\chi_+\left(\xi^{\f23} (r_c-r)/r_c^{\f13}\right)\chi\left(\xi^{\f23}(r-r_c)/r_c^{\f13}\right)
														\left(\f{i\overline{W}}{2|W|}f_+-
														\f{iW}{2|W|}\bar{f}_+\right),\\
														&\left(\mathcal{D}_{l,m}\widetilde{\phi}\right)_{6_{\infty}}(\xi,c,r)=\left(1- \chi_+(4r/r_c)\right)
														\chi_+\left(\xi^{\f23}(r_c-r)/r_c^{\f13}\right)
														\left(1-\chi\left(\xi^{\f23}(r-r_c)/r_c^{\f13}\right)\right)\\&\qquad\qquad\qquad\qquad\times
														\f{i\overline{W}(\xi,c)}{2|W(\xi,c)|}
														\left(\mathcal{D}_{l,m}f_+\right)(\xi,c,r)+c.c,
													\end{split}
												\end{align}
											\end{small}
											and
												\begin{align}
													\begin{split}\label{def:tphi-7-D}
														&\left(\mathcal{D}_{l,m}\widetilde{\phi}\right)_{7_0}(\xi,c,r)
														=\chi\left(\f{\xi r}{M (1-c)^{\f12}}\right)\widetilde{\phi},\\
														&\left(\mathcal{D}_{l,m}\widetilde{\phi}\right)_{7_1}(\xi,c,r)
														= \left(\chi(r/M^{\f12})-\chi\left(\f{\xi r}{M(1-c)^{\f12}}\right)\right)\widetilde{\phi},\\
														&\left(\mathcal{D}_{l,m}\widetilde{\phi}\right)_{7_2}(\xi,c,r)
														= \left(\chi\left(\f{M^{\f72}(1-c) r}{\xi^2}\right)-\chi(r/M)\right)\widetilde{\phi},\\
														&\left(\mathcal{D}_{l,m}\widetilde{\phi}\right)_{7_3}(\xi,c,r)
														=
														\left(\chi\left(M^{\f12}\xi r\right)-\chi\left(\f{M^{\f72}(1-c) r}{\xi^2}\right)\right)\widetilde{\phi},\\
														&\left(\mathcal{D}_{l,m}\widetilde{\phi}\right)_{7_{\infty}}(\xi,c,r)
														= \left(1-\chi\left(M^{\f12}\xi r\right)\right)\f{i\overline{W}}{2|W|} \left(\mathcal{D}_{l,m}f_+\right)(\xi,c,r)+c.c.
														\end{split}
												\end{align}
											
										\end{definition}
										
										\subsection{Proof of  the uniform estimate \eqref{op-K-D}}
										
										The proof is reduced to the following propositions.

										\begin{proposition}\label{lem:Boundness-K-pieces1-D}
											Let  $l=l',l\in\{0,1\}$, $m\in \mathbb{N}$,  $\left(\mathcal{D}_{l,m}^L\widetilde{\phi}\right)_j$ be defined in Definition \ref{def:decom-phi-D},
											$j\in J^L_{1}=\{1_0,1_{\infty}\} $.
											Then for $l'+l+m\geq 1$, $(j,j')\in J^L_{1}\times J^L_{1}$, it holds uniformly for  $(r,s,z)\in\mathbb{R}^+\times \mathbb{R}^+\times \mathbb{R}$  that
											\begin{align}
											\begin{split}
												\int_0^{1}&\left|\pa_c\left(p.v.\int_{\mathbb{R}} \left(1-\chi\left(\f{M^2(1-c)^{\f12}}{\xi}\right)\right)\right.\right.\\&\quad\left.\left.
												\left(\mathcal{D}_{l,0}\widetilde{\phi}\right)_{j}(\xi,c,s)
												\left(\mathcal{D}_{l',m}^L\widetilde{\phi}\right)_{j'}(\xi,c,r)\f{e^{i\xi \sqrt{\f{c}{1-c}}z}}{\xi}d\xi\right)\right|dc\lesssim s^{\delta}+1.
												\end{split}\label{op-pieces1-D}
												\end{align}
										\end{proposition}
										
										\begin{proof}
											 In the proof of  this case, we always assume $|\xi|\lesssim (1-c)^{\f12}$.
											The proof  relies on the bounds \eqref{1-phi}, \eqref{1-f+} in Proposition \ref{prop:summery} that for $s,r\lesssim |\xi|^{-1}$,
											\begin{align}
												\begin{split}\label{1-phi'}
													&\widetilde{\phi}(\xi,c,r)=
													O_{\xi,c}^{\xi,1-c}(\xi r),\;\;\widetilde{\phi}(\xi,c,s)=
													O_{\xi,c}^{\xi,1-c}(\xi s),\\
													&\left(\mathcal{D}^L_{l',m}\widetilde{\phi}\right)(\xi,c,r)=O_{\xi,c}^{\xi,1-c}
													(\xi)\;(l'+m\geq 1),\quad \left(\mathcal{D}_{1,0}\widetilde{\phi}\right)(\xi,c,s)=(1-c)^{\f12}O_{\xi,c}^{\xi,1-c}
													(1), \end{split}
												\end{align}
											and  for $s,r\gtrsim |\xi|^{-1}$,
											\begin{align}
												\begin{split}\label{1-f+'}
													&f_+(\xi,c,r)=O_{\xi,c}^{\xi,1-c}\left((\xi r)^{-\f12}\right)e^{i\xi r},\;f_+(\xi,c,s)=O_{\xi,c}^{\xi,1-c}\left((\xi s)^{-\f12}\right)e^{i\xi s},\\
													&\left(\mathcal{D}^L_{l',m}f_+\right)(\xi,c,r)=\xi O_{\xi,c}^{\xi,1-c}\left((\xi r)^{-\f12}\right)e^{i\xi r}=O_{\xi,c}^{\xi,1-c}\left((\xi r)^{-\f12}\right)e^{i\xi r}\;(l'+m\geq 1),\\ &\left(\mathcal{D}_{1,0}f_+\right)(\xi,c,s)=(1-c)^{\f12}O_{\xi,c}^{\xi,1-c}\left((\xi s)^{-\f12}\right)e^{i\xi s}=O_{\xi,c}^{\xi,1-c}\left((\xi s)^{-\f12}\right)e^{i\xi s}.
												\end{split}
											\end{align}
											
											If $(j,j')=(1_{\infty},1_{\infty})$, note that the final equalities for the derivative bounds in \eqref{1-f+'} are identical to the non-derivative ones.
											Thus, we have already completed this case in the proof of Proposition \ref{lem:Boundness-K-pieces1}.			
																			
											If $(j,j')=(1_{0},1_{\infty})$ and $l=0$, the proof is also identical to the corresponding case in Proposition \ref{lem:Boundness-K-pieces1}.
											If $(j,j')=(1_{0},1_{\infty})$ and $l=1$, thanks to the formulation in \eqref{def:tphi-1-D} and the bounds \eqref{1-phi}, \eqref{1-f+},
											\eqref{op-pieces1-D} can be reduced to
											\begin{align} \label{op1D-reduce1}
												\int_0^{1}\left|\pa_c\left(p.v.\int_{\mathbb{R}} \left(1-\chi(\xi r)\right)
												(1-c)^{\f12}O_{\xi,c}^{\xi,1-c}\left((\xi r)^{-\f12}\right)e^{\pm i\xi r}\f{e^{i\xi \sqrt{\f{c}{1-c}}z}}{\xi}d\xi\right)\right|dc\lesssim 1,
												\end{align}
											where, taking the transform $\eta_1=\xi r$ and denoting $y_1=\sqrt{\f{c}{1-c}}\f{z}{r}$, the left-hand side follows from Lemma \ref{lem: pifi-Linfty} and  Lemma \ref{lem:symbol-chi} that
											\begin{align*}
												&\int_0^{1}\left|\pa_c\left(p.v.\int_{\mathbb{R}} \left(1-\chi(\eta_1)\right)
												O_{\eta_1,c}^{\eta_1,(1-c)^{\f12}}\left(\eta_1^{-\f12}\right)\f{e^{i\eta_1 \left(y_1(c,z,r)\pm 1\right)}}{\eta_1}d\eta_1\right)\right|dc\\
												&\lesssim \int_0^{1}(1-c)^{-\f12}\left|\int_{\mathbb{R}} \left(1-\chi(\eta_1)\right)
												O_{\eta_1}^{\eta_1}\left(\eta_1^{-\f12}\right)\f{e^{i\eta_1 \left(y_1(c,z,r)\pm 1\right)}}{\eta_1}d\eta_1\right|dc\\
												&\quad+\int_{\mathbb{R}}\left|\int_{\mathbb{R}} \left(1-\chi(\eta_1)\right)
												O_{\eta_1}^{\eta_1}\left(\eta_1^{-\f12}\right)e^{i\eta_1 \left(y_1\pm 1\right)}d\eta_1\right|d(y_1\pm 1) \\
												&\lesssim 1.
											\end{align*}
											
											If $(j,j')=(1_{\infty},1_{0})$ and $l'=m=0$, the proof is also identical to the corresponding case in Proposition \ref{lem:Boundness-K-pieces1}. If $(j,j')=(1_{\infty},1_{0})$ and $l'+m\geq 1$, thanks to the formulation in \eqref{def:tphi-1-D} and  the bounds \eqref{1-phi} (first of the second line), \eqref{1-f+}(third line), \eqref{op-pieces1-D} can be reduced to
											\begin{align*}
												\int_0^{1}\left|\pa_c\left(p.v.\int_{\mathbb{R}} \left(1-\chi(\xi s)\right)
												\xi O_{\xi,c}^{\xi,1-c}\left((\xi s)^{-\f12}\right)e^{\pm i\xi s}\f{e^{i\xi \sqrt{\f{c}{1-c}}z}}{\xi}d\xi\right)\right|dc\lesssim 1.\end{align*}
											Noticing  $|\xi|\lesssim (1-c)^{\f12}$, the above can be treated similarly as in \eqref{op1D-reduce1} with the transform $\eta_1=\xi s$.
											
											If $(j,j')=(1_{0},1_{0})$ and $l=0$, the proof is similar to the corresponding case in Proposition \ref{lem:Boundness-K-pieces1}.
											If $(j,j')=(1_{0},1_{0})$, $l=1$ and $l'=m=0$,  thanks to the formulation in \eqref{def:tphi-1-D} and the bounds in\eqref{1-phi}, \eqref{op-pieces1-D} can be reduced to
											\begin{align*}
												\int_0^{1}\left|\pa_c\left(p.v.\int_{\mathbb{R}}\chi(\xi r) \xi O_{\xi,c}^{\xi,1-c}(\xi r)\f{e^{i\xi \sqrt{\f{c}{1-c}}z}}{\xi}d\xi\right)\right|dc\lesssim 1,
												\end{align*}
											which follows similarly as in \eqref{op1D-reduce1}.
											
											The case  $(j,j')=(1_{0},1_{0})$, $l=1$, $l'+m\geq 1$ the most worth mention, we reduce \eqref{op-pieces1-D}
											by the formulation in \eqref{def:tphi-1-D} and the bounds  in \eqref{1-phi} as
											\begin{align*}
												\int_0^{1}\left|\pa_c\left(p.v.\int_{\mathbb{R}}
												\left(1-\chi\left(\f{M^2(1-c)^{\f12}}{\xi}\right)\right) (1-c)^{\f12} O_{\xi,c}^{\xi,1-c}(\xi )\f{e^{i\xi \sqrt{\f{c}{1-c}}z}}{\xi}d\xi\right)\right|dc\lesssim 1.
												\end{align*}
											Indeed, we take the transform $\eta=\f{\xi}{(1-c)^{\f12}}$ and denoting $y=c^{\f12}z$, the left-hand side follows from Lemma \ref{lem: pifi-Linfty} and  Lemma \ref{lem:symbol-chi} that
											\begin{align*}
												&\int_0^{1}\left|\pa_c\left(p.v.\int_{\mathbb{R}} \left(1-\chi(M^2/\eta)\right)
												O_{\eta,c}^{\eta,1}(\eta)\f{e^{i\eta y(c,z)}}{\eta}d\eta\right)\right|dc\\
												&\lesssim \int_0^{1}\left|\int_{\mathbb{R}} \left(1-\chi(M^2/\eta)\right)
												O_{\eta}^{\eta}(\eta)\f{e^{i\eta y(c,z)}}{\eta}d\eta\right|dc
												+\int_{0}^1\left|\int_{\mathbb{R}} \left(1-\chi(M^2/\eta)\right)
												O_{\eta}^{\eta}(\eta)e^{i\eta y(c,z)}d\eta \right|\cdot|y'(c)|dc\\
												&\lesssim\int_{\mathbb{R}}\left|\int_{\mathbb{R}} \left(1-\chi(M^2/\eta)\right)
												O_{\eta}^{\eta}(\eta)e^{i\eta y}d\eta \right|dy +\int_{\mathbb{R}}\left|\int_{\mathbb{R}} \left(1-\chi(M^2/\eta)\right)
												O_{\eta}^{\eta}(\eta)e^{i\eta y}d\eta \right|dy\lesssim 1.
												\end{align*}
										\end{proof}
										
										\begin{remark}\label{rmk:1D-stronger}
											Following a similar process, one can prove a stronger estimate for $l+l'+m\geq 1$:
											\begin{align}
											\begin{split}
												\int_0^{1}&\left|\pa_c\left(p.v.\int_{\mathbb{R}} \left(1-\chi\left(\f{M^2(1-c)^{\f12}}{\xi}\right)\right)\right.\right.\\&\quad\left.\left.
												\left(\mathcal{D}_{l,0}\widetilde{\phi}\right)_{j}(\xi,c,s)
												\left(\pa_r^m\mathcal{D}_{l',0}\widetilde{\phi}\right)_{j'}(\xi,c,r)\f{e^{i\xi \sqrt{\f{c}{1-c}}z}}{\xi}d\xi\right)\right|dc\lesssim 1,
												\end{split}\label{op-pieces1-D'}
												\end{align}
											where the $\left(\mathcal{D}\widetilde{\phi}\right)_{j}$ are defined similarly as  \eqref{def:tphi-1-D}, with
											$\mathcal{D}\in\{\mathcal{D}_{l,0},\pa_r^m\mathcal{D}_{l',0}\}$, $\mathcal{D}_{l,0}=\left(\f{(1-c)^{\f12}\left(\pa_s+s^{-1}\right)}{\xi}\right)$, $l,l'\in\{0,1\}$, $m\in \mathbb{N}$.
										\end{remark}
										
										\begin{proposition}\label{lem:Boundness-K-pieces2-D}
											Let  $l=l',l\in\{0,1\}$, $m\in \mathbb{N}$,  $\left(\mathcal{D}_{l,m}\widetilde{\phi}\right)_j$ be defined in Definition \ref{def:decom-phi-D}, $j\in J^H_{2}=\{2_0,2_{\infty}
											\}  $.
											Then for $l'+l+m\geq 1$, $(j,j')\in J^H_{2}\times J^H_{2}$, it holds uniformly for  $(r,s,z)\in\mathbb{R}^+\times \mathbb{R}^+\times \mathbb{R}$  that
												\begin{align}
												\begin{split}
												\int_0^{V(0)}&\left|\pa_c\int_{\mathbb{R}} \chi\left(\f{M^2(1-c)^{\f12}}{\xi}\right) \chi\left(\xi (V(0)-c)^{\f32}\right)\right.\\&\quad\left.
													\left(\mathcal{D}_{l,0}\widetilde{\phi}\right)_{j}(\xi,c,s)
													\left(\mathcal{D}_{l',m}\widetilde{\phi}\right)_{j'}(\xi,c,r)\f{e^{i\xi z\sqrt{\f{c}{1-c}}}}{\xi}d\xi\right|dc\lesssim 1+s^{-\delta}.
													\end{split}\label{op-pieces2-D}
													\end{align}
											
										\end{proposition}
										
										\begin{proof}
											 Without further illustration, we assume $c\in(0,V(0))$, $1\lesssim \xi \lesssim (V(0)-c)^{-\f32}$.
											Let the new coordinates be defined as \eqref{trans:eta2}:
											\begin{align*}
												\begin{split}
													\eta(\xi,c,s)= \left\{
													\begin{array}{l} \xi s^{\f32} := \eta_{2_0}\;\;\;\quad\quad\quad\quad\quad\quad\quad\quad s\lesssim \xi^{-\f23} ,\\
														\xi \int_{0}^s\sqrt{\f{V(s')-c}{1-c}}ds':= \eta_{2_{\infty}} \;\;\;\quad\quad s\gtrsim \xi^{-\f23}, \end{array}\right.
												\end{split}
											\end{align*}
											whose inverse function is denoted as $\xi(\eta_{j},c,s)$,
											where $\eta_j$ satisfies \eqref{lower:eta_2infty} and \eqref{est:paceta/eta2} with the weight $\rho(c,\xi,s)=\xi^{-\f23}:= \rho_{2}$.
											We take the transform such that
											\begin{align}
												\begin{split}
													\label{22:1D}&\left(\mathcal{D}_{l',0}\phi\right)_{j}(\eta_{j},c,s)
													:=\left(\mathcal{D}_{l',0}\widetilde{\phi}\right)_{j}\left(\xi(\eta_{j},c,s),s,c\right), \\&\left(\mathcal{D}_{l,m}\phi\right)^*_{j'}(\eta_{j},c,s,r)
													:=\left(\mathcal{D}_{l,m}\widetilde{\phi}\right)^*_{j'}
													\left(\xi(\eta_{j'},c,s),r,c\right).
												\end{split}
											\end{align}
											For simplicity, we use the notation $\mathcal{D}=\mathcal{D}_{l,m}=\left(\f{(1-c)^{\f12}\pa_r}{\xi}\right)^{m}
											\left(\f{(1-c)^{\f12}(\pa_r+r^{-1})}{\xi}\right)^{l}$, since the estimates are the same for different $l,m$. It follows from the formula \eqref{fm:f-tf1},
											the definitions of $\left(\mathcal{D}\widetilde{\phi}\right)_{j}$ in \eqref{def:tphi-2-D},  the bounds \eqref{2-phi}, \eqref{2-f+} in Proposition \ref{prop:summery} that
											\begin{align}  \label{bd:20-D}
													\begin{split}
														(\mathcal{D}\phi)_{2_0}(\eta_{2_0},c,s)&=
														\chi(\eta_{2_0}^{\f23})
														O_{\eta_{2_0},c}^{\eta_{2_0},\rho_2}(1)
														=\chi\left(\eta_{2_0}^{\f23}\right)
														O_{\eta_{2_0},c}^{\eta_{2_0},\xi^{-\f23}}\left(\eta^{\f23\delta}\cdot \xi^{-\f23}\cdot s^{-\delta}\right)
														\\&=\chi\left(\eta_{2_0}^{\f23}\right)
														O_{\eta_{2_0},c}^{\eta_{2_0},\xi^{-\f23}(1-\delta)}\left(\eta^{\f23\delta}_{2_0} \right)      s^{-\delta}  \\
														& =s^{-\delta}\chi\left(\eta_{2_0}^{\f23}\right)
														O_{\eta_{2_0},c}^{\eta_{2_0},(V(0)-c)^{1-\delta}}\left(\eta^{\f23\delta} \right),\\
														\rho_2^{-1}(\mathcal{D}\phi)_{2_0}(\eta_{2_0},c,s)& =s^{-\delta}\chi\left(\eta_{2_0}^{\f23}\right)
														O_{\eta_{2_0}}^{\eta_{2_0}}\left(\eta^{\f23\delta}_{2_0} \right)(V(0)-c)^{-1+\delta}, \end{split}
												\end{align}
											and similar to \eqref{bd:2infty}, we get by the derivative bounds \eqref{2-f+} in Proposition \ref{prop:summery} that
										 \begin{align}  \label{bd:2infty-D}
													\begin{split}
														(\mathcal{D}\phi)_{2_{\infty}}(\eta_{2_{\infty}},c,s)&=
														\left(1-\chi\left(\xi^{\f23} s\right)\right)
														\left(O_{\eta_{2_{\infty}},c}^{\eta_{2_{\infty}},\rho_2}\left(\eta_{2_{\infty}}
														^{-\f12}\right)e^{i\eta_{2_{\infty}}}
														+O_{\eta_{2_{\infty}},c}^{\eta_{2_{\infty}},\rho_2}\left(\eta_{2_{\infty}}
														^{-\f12}\right)e^{-i\eta_{2_{\infty}}}\right)\\
														&=\sum_{\iota\in\{-,+\}}\left(1-\chi\left(M\eta_{2_{\infty}}\right)\right)
														O_{\eta_{2_{\infty}},c}^{\eta_{2_{\infty}},
															\sigma_{2_{\infty}}(c,s)^{-1}}\left(\eta_{2_{\infty}}
														^{-\f12+\delta/6}\right)e^{\iota i\eta_{2_{\infty}}},\\
														\rho^{-1}_2(\mathcal{D}\phi)_{2_{\infty}}(\eta_{2_{\infty}},c,s)& =\sigma_{2_{\infty}}(c,s)\sum_{\iota\in\{-,+\}}\left(1-\chi\left(M\eta_{2_{\infty}}\right)\right)
														O_{\eta_{2_{\infty}}}^{\eta_{2_{\infty}}}\left(\eta_{2_{\infty}}
														^{-\f12+\delta/6}\right)e^{\iota i\eta_{2_{\infty}}}.
													\end{split}
												\end{align}
											where $\sigma_{2_{\infty}}(c,s)$ is as \eqref{def:sigma2}, satisfying \eqref{int:sigmainfty}.
											Using \eqref{useful:trans-2} with $\rho=\rho_2=\xi^{-\f23}$ and \eqref{est:paceta/eta2},
											we get by  the derivative bounds \eqref{2-phi}, \eqref{2-f+} in Proposition \ref{prop:summery} that \begin{align}  \label{bi-bd:20-D}
												\begin{split}
													& (\mathcal{D}\phi)_{2_0}^*(\eta_{j},c,s,r)= O_{\eta_{2_0},c}^{\eta_{2_0},\rho_2}(1),\;j=2_0,2_{\infty},\\
													& (\mathcal{D}\phi)_{2_{\infty}}^*(\eta_{2_0},c,s,r)= \sum_{\iota\in\{-,+\}}O_{\eta_{2_{\infty}},c}^{\eta_{2_{\infty}},
														\rho_2}(1)e^{\iota i\eta_{2_{\infty}}}.
													\end{split}
											\end{align}
											For the cut-off function, by $c\in(0,V(0))$, $1\lesssim |\xi| \lesssim (V(0)-c)^{-\f32}$, we get by \eqref{est:paceta/eta2} that
											\begin{align} \label{bd-cutoff:20-D}
												\chi\left(\f{(1-c)^{\f12}}{\xi}\right) \chi\left(\xi (V(0)-c)^{\f32}\right)=O_{\eta_{j},c}^{\eta_{j},\rho_2}(1),\;\;j=2_0,2_{\infty}.
											\end{align}
											
											 Now we are in a position to prove \eqref{op-pieces2-D}. If $(j,j')= (2_0, 2_0)$, we reduce it to
											\begin{align*} I_{2_0,2_0,l',l,m}:=\int_0^{V(0)}\left|\pa_c\left(p.v.\int_{\mathbb{R}} a_{2_0,2_0,l,l',m}(\eta_{2_0},c,\lambda)\f{e^{i\eta_{2_0} y_{2_0}(c,\lambda) }}{\eta_{2_0}}d\eta_{2_0}\right)\right|dc\lesssim s^{-\delta}+1,
											\end{align*}
											where $\lambda=(z,s,r)$ is the parameter, $y_{2_0}(c,\lambda)$ as \eqref{def:y20} is monotonic in $c$, and by \eqref{bd:20-D}, \eqref{bi-bd:20-D}, \eqref{bd-cutoff:20-D},
											\begin{align} \label{def:a_20-D}
												\notag
												a_{2_0,2_0,l,l',m}(\eta_{2_0},c,\lambda)&:=\chi\left(\f{(1-c)^{\f12}}{\xi}\right) \chi\left(\xi(V(0)-c)^{\f32}\right)(\mathcal{D}_{l,0}\phi)_{2_0}(\eta_{2_0},c,s)
												(\mathcal{D}_{l',m}\phi)_{2_0}^* (\eta_{2_0},c,s,r),\\
												&= s^{-\delta}\chi\left(\eta_{2_0}^{\f23}\right)
												O_{\eta_{2_0},c}^{\eta_{2_0},(V(0)-c)^{1-\delta}}\left(\eta^{\f23\delta}_{2_0} \right).
											\end{align}
											If $(j,j')= (2_{\infty}, 2_0)$, we take $\eta_{2_{\infty}}=\xi \int_{0}^s\sqrt{\f{V(s')-c}{1-c}}ds' $ and reduce \eqref{op-pieces2-D} to
											\begin{align*} I_{2_{\infty},2_0,l,l',m}:
												=\int_0^{V(0)}\left|\pa_c\left(p.v.\int_{\mathbb{R}} a_{2_{\infty},2_0,l',l,m}(\eta_{2_{\infty}},c,\lambda)
												\f{e^{i\eta_{2_{\infty}} y_{2_{\infty}}(c,\lambda) }}{\eta_{2_{\infty}}}d\eta_{2_{\infty}}\right)\right|dc\lesssim 1,
												\end{align*}
											where  $y_{2_{\infty}}(c,\lambda)$ as \eqref{def:y20} is monotonic in $c$, and by \eqref{bd:20-D}, \eqref{bi-bd:20-D}, \eqref{bd-cutoff:20-D},
											\begin{align} \label{bd:a2infty-D}
												\notag
												a_{2_{\infty},2_0,l',l,m}(\eta_{2_{\infty}},c,\lambda)
												&:=\chi\left(\f{(1-c)^{\f12}}{\xi}\right) \chi\left(\xi(V(0)-c)^{\f32}\right)(\mathcal{D}_{l,0}\phi)_{2_{\infty}}
												(\eta_{2_{\infty}},c,s)
												(\mathcal{D}_{l',m}\phi)_{2_0}^* (\eta_{2_{\infty}},c,s,r),\\
												&=\sum_{\iota\in\{-,+\}} \left(1-\chi\left(M\eta_{2_{\infty}}\right)\right)
												O_{\eta_{2_{\infty}},c}^{\eta_{2_{\infty}},
													\sigma_{2_{\infty}}(c,s)^{-1}}\left(\eta_{2_{\infty}}
												^{-\f12+\delta/6}\right)e^{\iota i\eta_{2_{\infty}}}.
											\end{align}
											Then the oscillation bounds $I_{2_0,2_0,l',l,m}$ and  $I_{2_{\infty},2_0,l',l,m}$ follow in the same manner as in \eqref{process:integral20}, utilizing \eqref{def:a_20-D} and \eqref{bd:a2infty-D}, \eqref{int:sigmainfty}, respectively.		
																				
											If $(j,j')=(2_0,2_{\infty})$, we interchange the variables $s$ and $r$, and swap $j$ and $j'$,
											thereby reducing the case to $(j,j')=(2_{\infty},2_{0})$. For $(j,j')=(2_{\infty},2_{\infty})$,
											the proof is identical to case 3 in the proof of Proposition \ref{lem:Boundness-K-pieces2},
											as the bounds in \eqref{2-f+} (given in Proposition \ref{prop:summery}) are the same for both the derivative and non-derivative versions.										\end{proof}
										
										\begin{proposition}\label{lem:Boundness-K-pieces3-D}
											 Let  $l',l\in\{0,1\}$, $m\in \mathbb{N}$, $j\in J^H_{3}=\{3_0,3_{\infty}
											\}  $. Let $(\mathcal{D}_{l,m}\widetilde{\phi})_j$ be defined in \eqref{def:tphi-3-D}, Definition \ref{def:decom-phi-D} .
											Then for $l'+l+m\geq 1$, $(j,j')\in J^H_{3}\times J^H_{3}$, it holds uniformly for $(r,s,z)\in\mathbb{R}^+\times \mathbb{R}^+\times \mathbb{R}$ that
												\begin{align}
												\begin{split}
												\int_0^{V(0)}&\left|\pa_c\int_{\mathbb{R}} \chi\left(\f{M^2\sqrt{1-c}}{\xi}\right) \left(1-\chi\left(\xi (V(0)-c)^{\f32}\right)\right)\right.\\&\quad\left.
													\left(\mathcal{D}_{l,0}\widetilde{\phi}\right)_{j}(\xi,c,s)
													\left(\mathcal{D}_{l',m}\widetilde{\phi}\right)_{j'}(\xi,c,r)\f{e^{i\xi z\sqrt{\f{c}{1-c}}}}{\xi}d\xi\right|dc\lesssim s^{-
														\delta}+1.
														\end{split} \label{op-pieces3-D}
													\end{align}
											
										\end{proposition}
										
										\begin{proof}
											Without further illustration, we assume $c\in(0,V(0))$, $ \xi \gtrsim  (V(0)-c)^{-\f32}$ in the proof.
											Let the new coordinates be defined as \eqref{trans:eta3}:
											\begin{align*}
												\begin{split}
													\eta(\xi,c,s)= \left\{
													\begin{array}{l} \xi (V(0)-c)^{\f12}s = \eta_{3_0}\;\;\;\quad\quad\quad s\lesssim |\xi|^{-1}(V(0)-c)^{-\f12} ,\\
														\xi \int_{0}^s\sqrt{\f{V(s')-c}{1-c}}ds'= \eta_{3_{\infty}} \;\;\;\quad\quad s\gtrsim |\xi|^{-1}(V(0)-c)^{-\f12} , \end{array}\right.
												\end{split}
											\end{align*}
											whose inverse function is denoted as $\xi(\eta_{j},c,s)$, where $\eta_j$ satisfies \eqref{lower:eta_3infty} and \eqref{est:paceta/eta3} with the weight  $\rho(c,\xi,s)=V(0)-c:= \rho_{3}$. We take the transform  for $(j,j')\in\{3_0,3_{\infty}\}\times \{3_0,3_{\infty}\}$ such that
											\begin{align*}
												\begin{split}
													&(\mathcal{D}_{l,0}\phi)_{j}(\eta_{j},c,s)
													:=\left(\mathcal{D}_{l,0}\widetilde{\phi}\right)_{j}\left(\xi(\eta_{j},c,s),s,c\right), \\
													&(\mathcal{D}_{l',m}\phi)^*_{j'}(\eta_{j},c,s,r)
													:=(\mathcal{D}_{l',m}\widetilde{\phi})^*_{j'}
													\left(\xi(\eta_{j'},c,s),r,c\right).
												\end{split}
											\end{align*}
											For simplicity, we use the notation $\mathcal{D}\in\{\mathcal{D}_{l,m},
											\mathcal{D}_{l',0}\}$, since the estimates are the same for different $l,l',m$. It follows by the formula \eqref{fm:f-tf1},
											the definitions of $(\mathcal{D}\widetilde{\phi})_{j}$ in \eqref{def:tphi-3-D},  the bounds \eqref{3-phi}, \eqref{3-f+} in Proposition \ref{prop:summery} that
											\begin{small}
												\begin{align}  \label{bd:30-D}
													\begin{split}
														(\mathcal{D}\phi)_{3_0}(\eta_{3_0},c,s)&=
														\chi(\eta_{3_0})
														O_{\eta_{3_0},c}^{\eta_{3_0},\rho_3}(1)
														=\chi(\eta_{3_0})
														O_{\eta_{3_0},c}^{\eta_{3_0},V(0)-c}\left(\eta_{3_0}^{\delta}\cdot (V(0)-c)^{\delta}\cdot \left(\xi (V(0)-c)^{\f32}\right)^{-\delta}\cdot s^{-\delta}\right) \\
														&=\chi(\eta_{3_0})
														O_{\eta_{3_0},c}^{\eta_{3_0},(V(0)-c)(1-\delta)}(\eta_{3_0}^{\delta} )      s^{-\delta}  \\
														& =s^{-\delta}\chi(\eta_{3_0})
														O_{\eta_{3_0},c}^{\eta_{3_0},(V(0)-c)^{1-\delta}}\left(\eta^{\delta}_{3_0} \right)      ,\\
														\rho_3^{-1}(\mathcal{D}\phi)_{3_0}(\eta_{3_0},c,s)& =s^{-\delta}\chi(\eta_{3_0})
														O_{\eta_{3_0}}^{\eta_{3_0}}\left(\eta_{3_0}^{\delta} \right)      (V(0)-c)^{-1+\delta},
														\end{split}
												\end{align}
											\end{small}
											and similar to \eqref{bd:3infty}, we get by the derivative bounds of \eqref{3-f+}, in Proposition \ref{prop:summery} and  \eqref{lower:eta_3infty}, \eqref{est:paceta/eta3},  \eqref{decompose:rho-3}, \eqref{rewriten:m-eta3} that
											\begin{small}
												\begin{align}  \label{bd:3infty-D}
													\begin{split}
														(\mathcal{D}\phi)_{3_{\infty}}(\eta_{3_{\infty}},c,s)&=
														\left(1-\chi(M\eta_{3_{\infty}})\right)\left(1-\chi\left(\xi (V(0)-c)^{\f12} s\right)\right)
														\left(O_{\eta_{3_{\infty}},c}^{\eta_{3_{\infty}},\rho_3}\left(\eta_{3_{\infty}}
														^{-\f12}\right)e^{i\eta_{3_{\infty}}}
														+O_{\eta_{3_{\infty}},c}^{\eta_{3_{\infty}},\rho_3}\left(\eta_{3_{\infty}}
														^{-\f12}\right)e^{-i\eta_{2_{\infty}}}\right)\\
														&= \sum_{\iota\in\{-,+\}}\left(1-\chi(M\eta_{3_{\infty}})\right)
														O_{\eta_{3_{\infty}},c}^{\eta_{3_{\infty}},
															\sigma_{3_{\infty}}(c,s)^{-1}}\left(\eta_{3_{\infty}}
														^{-\f12+\delta/6}\right)e^{\iota i\eta_{3_{\infty}}},\\
														\rho^{-1}_3(\mathcal{D}\phi)_{3_{\infty}}(\eta_{2_{\infty}},c,s)& =\sigma_{3_{\infty}}(c,s)\sum_{\iota\in\{-,+\}}\left(1-\chi(M\eta_{3_{\infty}})\right)
														O_{\eta_{3_{\infty}}}^{\eta_{3_{\infty}}}\left(\eta_{3_{\infty}}
														^{-\f12+\delta/6}\right)e^{\iota i\eta_{3_{\infty}}},
													\end{split}
												\end{align}
											\end{small}
											 where $\sigma_{3_{\infty}}(c,s)$ is as \eqref{def:sigma3}, satisfying \eqref{int:sigma3infty}. Using \eqref{useful:trans-2} with $\rho=\rho_3=V(0)-c$ and \eqref{est:paceta/eta3}, we get by the derivative bounds of \eqref{3-phi}, \eqref{3-f+}, in Proposition \ref{prop:summery} that
											 \begin{align}  \label{bi-bd:30-D}
												\begin{split}
													& (\mathcal{D}\phi)_{3_0}^*(\eta_{j},c,s,r)= O_{\eta_{3_0},c}^{\eta_{3_0},\rho_3}(1),\;j=3_0,3_{\infty},\\
													& (\mathcal{D}\phi)_{3_{\infty}}^*(\eta_{3_0},c,s,r)= \sum_{\iota\in\{-,+\}}O_{\eta_{3_{\infty}},c}^{\eta_{3_{\infty}},
														\rho_3}(1)e^{\iota i\eta_{3_{\infty}}}.
													\end{split}
											\end{align}
											For the cut-off function, by $c\in(0,V(0))$, $|\xi| \gtrsim (V(0)-c)^{-\f32}$, we get by \eqref{est:paceta/eta3} that
											\begin{align} \label{bd-cutoff:30-D}
												\chi\left(\f{(1-c)^{\f12}}{\xi}\right)\left(1- \chi\left(\xi (V(0)-c)^{\f32}\right)\right)=O_{\eta_{j},c}^{\eta_{j},\rho_3}(1),
												\;\;j=3_0,3_{\infty}.
											\end{align}
											
											 Now we are in a position to prove \eqref{op-pieces3-D}. If $(j,j')= (3_0, 3_0)$, we reduce it to
											\begin{align*} I_{3_0,3_0,l,l',m}:=\int_0^{V(0)}\left|\int_{\mathbb{R}} a_{3_0,3_0,l,l',m}(\eta_{3_0},c,\lambda)\f{e^{i\eta_{3_0} y_{3_0}(c,\lambda) }}{\eta_{3_0}}d\eta_{3_0}\right|dc\lesssim s^{-\delta},\end{align*}
											where $\lambda=(z,s,r)$ is the parameter, $y_{3_0}(c,\lambda)$ as \eqref{def:y30} is monotonic in $c$, and by \eqref{bd:30-D}, \eqref{bi-bd:30-D}, \eqref{bd-cutoff:30-D},
											\begin{align} \label{def:a_30-D}
												\notag
												a_{3_0,3_0,l,l',m}(\eta_{3_0},c,\lambda)&:=\chi\left(\f{(1-c)^{\f12}}{\xi}\right) \left(1-\chi\left(\xi(V(0)-c)^{\f32}\right)\right)(\mathcal{D}_{l,0}\phi)_{3_0}
												(\eta_{3_0},c,s)
												(\mathcal{D}_{l',m}\phi)_{3_0}^* (\eta_{3_0},c,s,r),\\
												&= s^{-\delta}\chi(\eta_{3_0})
												O_{\eta_{3_0},c}^{\eta_{3_0},(V(0)-c)^{1-\delta}}(\eta^{\delta}_{3_0} ).
											\end{align}
											If $(j,j')= (3_{\infty}, 3_0)$, we take $\eta_{3_{\infty}}=\xi\int_{0}^s\sqrt{\f{V(s')-c}{1-c}}ds' $ and reduce  \eqref{op-pieces3-D} to
											\begin{align*} I_{3_{\infty},3_0,l,l',m}:
												=\int_0^{V(0)}\left|\pa_c\int_{\mathbb{R}} a_{3_{\infty},3_0,l,l',m}(\eta_{3_{\infty}},c,\lambda)\f{e^{i\eta_{3_{\infty}} y_{3_{\infty}}(c,\lambda) }}{\eta_{3_{\infty}}}d\eta_{3_{\infty}}\right|dc\lesssim 1,
												\end{align*}
											where  $y_{3_{\infty}}(c,\lambda)$ as \eqref{def:y30} is monotonic in $c$, and by \eqref{bd:30-D}, \eqref{bi-bd:30-D}, \eqref{bd-cutoff:30-D},
											\begin{small}
												\begin{align} \label{bd:a3infty-D}
													\notag
													a_{2_{\infty},2_0,l,l',m}(\eta_{2_{\infty}},c,\lambda)
													&:=\chi\left(\f{(1-c)^{\f12}}{\xi}\right) \left(1-\chi\left(\xi(V(0)-c)^{\f32}\right)\right)(\mathcal{D}_{l,0}\phi)_{2_{\infty}}
													(\eta_{3_{\infty}},c,s)
													(\mathcal{D}_{l',m}\phi)_{3_0}^* (\eta_{3_{\infty}},c,s,r),\\
													&=\sum_{\iota\in\{-,+\}} \left(1-\chi(M\eta_{3_{\infty}})\right)
													O_{\eta_{3_{\infty}},c}^{\eta_{3_{\infty}},
														\sigma_{3_{\infty}}(c,s)^{-1}}\left(\eta_{3_{\infty}}
													^{-\f12+\delta/6}\right)e^{\iota i\eta_{3_{\infty}}}.
												\end{align}
											\end{small}

											Then the oscillation bounds $I_{3_0,3_0,l,l',m}$ and  $I_{3_{\infty},3_0,l,l',m}$ follow in the same manner as in \eqref{process:integral20},
											utilizing \eqref{def:a_30-D}, and \eqref{bd:a3infty-D}, \eqref{int:sigma3infty}, respectively.
											
											If $(j,j')=(3_0,3_{\infty})$, we interchange the variables $s$ and $r$, and swap $j$ and $j'$, thereby reducing the case to $(j,j')=(3_{\infty},3_{0})$.
											For $(j,j')=(3_{\infty},3_{\infty})$, the proof is identical to case 3 in the proof of Proposition \ref{lem:Boundness-K-pieces3},
											as the bounds in \eqref{3-f+}(given in Proposition \ref{prop:summery})  are the same for both the derivative and non-derivative versions.
										\end{proof}
										
										\begin{proposition}\label{lem:Boundness-K-pieces4-D}
											Let $\delta\ll 1$ be fixed. Let  $l',l\in\{0,1\}$, $m\in \mathbb{N}$, $j\in J^H_{4}=\{4_0,4_{\infty}
											\} $. Let  $(\mathcal{D}_{l,m}\widetilde{\phi})_j$ be defined in  \eqref{def:tphi-4-D}, Definition \ref{def:decom-phi-D} .
											Then for $l'+l+m\geq 1$, $(j,j')\in J^H_{4}\times J^H_{4}$, it holds uniformly for  $(r,s,z)\in\mathbb{R}^+\times \mathbb{R}^+\times \mathbb{R}$  that
											\begin{align}
											\begin{split}
											\int_{V(0)}^{1-\delta}&\left|\pa_c\int_{\mathbb{R}} \chi\left(\f{M^2(1-c)^{\f12}}{\xi}\right) \chi\left(\xi (c-V(0))^{\f32}\right)\right.\\&\quad\left.
													(\mathcal{D}_{l,0}\widetilde{\phi})_{j}(\xi,c,s)
													(\mathcal{D}_{l',m}\widetilde{\phi})_{j'}(\xi,c,r)\f{e^{i\xi \sqrt{\f{c}{1-c}}z}}{\xi}d\xi\right|dc\lesssim s^{-\delta}+1.
													\end{split}
													\label{op-pieces4-D}
													\end{align}
											
										\end{proposition}
										
										\begin{proof}
											The proof is similar to that of Proposition \ref{lem:Boundness-K-pieces4-D}, with the coordinate transformation given in \eqref{trans:eta4}:
											\begin{align*}
												\begin{split}
													\eta(\xi,c,s)= \left\{
													\begin{array}{l} \xi s^{\f32} = \eta_{4_0}\;\;\;\quad\quad\quad\quad\quad\quad\quad\quad s\lesssim \xi^{-\f23} ,\\
														\xi \int_{r_c}^s\sqrt{\f{V(s')-c}{1-c}}ds'=\eta_{4_{\infty}} \;\;\;\quad\quad s\gtrsim \xi^{-\f23}, \end{array}\right.
												\end{split}
											\end{align*}
											whose inverse function is denoted as $\xi(\eta_{j},c,s)$, where $\eta_j$ satisfies \eqref{lower:eta_4infty} and \eqref{est:paceta/eta4} with the weight $\rho(c,\xi,s)=\xi^{-\f23}:= \rho_{4}$. For $(j,j')\in \{4_0,4_{\infty}\}\times \{4_0,4_{\infty}\}$,
											\begin{align*}
												\begin{split}
													&(\mathcal{D}_{l',0}\phi)_{j}(\eta_{j},c,s)
													:=(\mathcal{D}_{l',0}\widetilde{\phi})_{j}\left(\xi(\eta_{j},c,s),s,c\right), \\
													&(\mathcal{D}_{l,m}\phi)^*_{j'}(\eta_{j},c,s,r)
													:=(\mathcal{D}_{l,m}\widetilde{\phi})^*_{j'}
													\left(\xi(\eta_{j'},c,s),r,c\right).
												\end{split}
											\end{align*}
											The details are left to the readers.
										\end{proof}
										
										\begin{proposition}\label{lem:Boundness-K-pieces5-D}
											Let  $\delta\ll 1$  and $M\gg 1$ be fixed. Let  $l',l\in\{0,1\}$, $m\in \mathbb{N}$,  $j\in J^H_{5}=
											\{5_0,5_1,5_2,5_3,5_{\infty}
											\} $. Let  $(\mathcal{D}_{l,m}\widetilde{\phi})_j$ be defined in  \eqref{def:tphi-5-D}, Definition \ref{def:decom-phi-D} .
											Then for $l'+l+m\geq 1$, $(j,j')\in J^H_{5}\times J^H_{5}$,  it holds uniformly for  $(r,s,z)\in\mathbb{R}^+\times \mathbb{R}^+\times \mathbb{R}$  that
											\begin{align}
											\begin{split}
											 \int_{V(0)}^{1-\delta}&\left|\pa_c\int_{\mathbb{R}} \chi(\f{M^2(1-c)^{\f12}}{\xi})\Big(1- \chi(\xi (V(0)-c)^{\f32}/M) \Big)\right.\\
													&\quad\left.(\mathcal{D}_{l,0}\widetilde{\phi})_{j}(\xi,c,s)
													(\mathcal{D}_{l',m}\widetilde{\phi})_{j'}(\xi,c,r)\f{e^{i\xi \sqrt{\f{c}{1-c}}z}}{\xi}d\xi\right|dc\lesssim s^{-\delta}+1.
												\end{split}
													\label{op-pieces5-D}\end{align}
										
										\end{proposition}
										
										\begin{proof}
											Without further illustration, we assume $c\in(V(0),1-\delta),\;\xi r_c^{\f32}\gtrsim M$, which implies $0<c-V(0)\sim r_c\lesssim 1$.
											We only modify the first coordinate $\eta_{5_0}$ among the five coordinates in \eqref{trans:eta5}, such that											                                                                               \begin{align} \label{trans:eta5-D}\begin{split}
													\eta(\xi,c,s)= \left\{
													\begin{array}{l} \xi(c-V(0))^{\f32} = \eta_{5_0}\;\;\;\quad\quad\quad\quad s\lesssim  M^{\f12}\xi^{-1}r_c^{-\f12},\\
														\xi(c-V(0))^{\f12} s= \eta_{5_1}\;\;\quad\quad\quad\quad\; M^{\f12}\xi^{-1}r_c^{-\f12}\lesssim s\leq r_c/2,\\
														\xi (c-V(0))^{\f32}= \eta_{5_2}\;\;\;\quad\;\;\quad\quad\;\;\; r_c/2\leq s\leq r_c-C\xi^{-\f23},\\
														\xi (c-V(0))^{\f32}= \eta_{5_3}\;\;\;\;\; \quad\quad \quad \quad r_c-C\xi^{-\f23}\leq s\leq r_c+C\xi^{-\f23},\\
														\xi \int_{r_c}^s\sqrt{\f{V(s')-c}{1-c}}ds'=\eta_{5_{\infty}} \;\;\;\;\;\quad\quad s\geq r_c+C\xi^{-\f23},
														\end{array}\right.
												\end{split}
											\end{align}
											whose inverse function are denoted as $\xi(\eta_{j},c,s)$. Moreover, $\eta_j$ satisfies  the bounds \eqref{est:paceta/eta5},
											with the same weights as \eqref{recall:rho-5}:
											\begin{align*} \rho(c,\xi,s)=
												\left\{
												\begin{array}{l}
													\rho_{5_0}=\rho_{5_1}=r_c\;\;\; \;\;\quad\quad\quad\quad s\leq r_c/2,\\
													\rho_{5_2}=(r_c-s) \;\;\; \;\quad\quad \quad\quad r_c/4\leq s\leq r_c-C\xi^{-\f23}, \quad\quad\quad\quad\\
													\rho_{5_3}=\xi^{-\f23}\;\;\; \;\;\;\;\;\;\quad\quad\quad\quad r_c-C\xi^{-\f23}\leq s\leq r_c+C\xi^{-\f23}, \\
													\rho_{5_{\infty}}=(s-r_c) \;\;\; \quad\quad\; \quad \quad s\geq r_c+C\xi^{-\f23} .
												\end{array}\right. \end{align*}
											If \emph{$j'\neq 5_{\infty}$}, i.e,   $(j,j')\in  \{5_0,5_1,5_2,5_3,5_{\infty}\}\times \{5_0,5_1,5_2,5_3
											\}$, we make the transform
											\begin{align*}
													&(\mathcal{D}_{l',0}\phi)_{j}(\eta_{j},c,s)
													:=(\mathcal{D}_{l',0}\widetilde{\phi})_{j}\left(\xi(\eta_{j},c,s),s,c\right), \\
													&(\mathcal{D}_{l,m}\phi)^*_{j'}(\eta_{j},c,s,r)
													:=(\mathcal{D}_{l,m}\widetilde{\phi})^*_{j'}
													\left(\xi(\eta_{j'},c,s),r,c\right).
											\end{align*}
											Then we  reduce \eqref{op-pieces5-D} to
											\begin{align}\label{bd:I5-D}
											I_{j,j',l,l',m}:=\int^{1-\delta}_{V(0)}\left|\int_{\mathbb{R}} a_{j,j',l,l',m}(\eta_{j},c,\lambda)\f{e^{i\eta_{j} y_{j}(c,\lambda) }}{\eta_{j}}d\eta_{j}\right|dc\lesssim 1,
											\end{align}
											where $\lambda=(z,s,r)$ is the parameter and
											\begin{align}\label{def:y_j5-D} y_j(c,\lambda)=
												\left\{
												\begin{array}{l}
													z\sqrt{\f{c}{1-c}}/(c-V(0))^{\f32}\;\;\; \;\;\quad\quad j=5_0,5_2,5_3,\\
													z\sqrt{\f{c}{1-c}}/(c-V(0))^{\f12}s\;\;\; \;\;\quad\quad j=5_1,\\
													z\sqrt{\f{c}{1-c}}/\int_{r_c}^s\sqrt{\f{V(s')-c}{1-c}}ds'\;\;\; \;\;\quad j=5_{\infty},\end{array}\right. \end{align}
											and  \begin{align} \label{def:a_50-D}
												&\notag a_{j,j',l,l',m}(\eta_{j},c,\lambda)\\
												&:=\chi\left(\f{(1-c)^{\f12}}{\xi}\right) \left(1-\chi\left(\xi(c-V(0))^{\f32}/M\right)\right)(\mathcal{D}_{l,0}\phi)_{j}
												(\eta_{j},c,s)(\mathcal{D}_{l',m}\phi)_{j'}^* (\eta_{j},c,s,r).
											\end{align}
											We claim that
											\begin{align} \label{bd:a_j,j'5-D}
												\begin{split}
													a_{j,j'}(\eta_j,c,\lambda)&=
													\left\{
													\begin{array}{l}
														\left(1-\chi(\eta_{5_0}/M)\right)
														O_{\eta_{5_0},c}^{\eta_{5_0},r_c^{\f12}}
														\left(\eta_{5_0}^{-1}\right),\;\;\;\quad\quad\quad\quad\quad\quad\quad j=5_0, \\
														\left(1-\chi(\eta_{5_1}) \right) O_{\eta_{5_1},c}^{\eta_{5_1},r_c^{\f12}}
														\left(\eta_{5_1}^{-\f12}\right),
														\;\;\;\quad\quad \quad\quad \quad \quad\quad \quad\quad  j=5_1,\\
														\left(1-\chi\left(\eta_{5_2}/M^{\f12}\right)\right) O_{\eta_{5_2},c}^{\eta_{5_2},\sigma_{5_2,j'}(c,s,r)^{-1}}
														\left(\eta_{5_2}^{-\f16+\delta/9}\right),\;\;\; \quad\quad j=5_2,\\
														\left(1-\chi\left(\eta_{5_3}/M^{\f12}\right)\right)
														O_{\eta_{5_3},c}^{\eta_{5_3},\sigma_{5_3,j'}(c,s,r)^{-1}}
														\left(\eta_{5_3}^{-\f13+\delta/9}\right),\; \;\;\; \quad \quad j=5_3,\\
														\sum_{\iota\in\{-,+\}}\left(1-\chi\left(M^{\f12}\eta_{5_{\infty}}\right)\right) O_{\eta_{5_{\infty}},c}^{\eta_{5_{\infty}},\sigma_{5_{\infty},j'}
															(c,\lambda)
															^{-1}}\left(\eta_{5_{\infty}}^{-\f16}\right)e^{\iota i \eta_{5_{\infty}}},\quad j=5_{\infty}, \end{array}\right.
												\end{split}
											\end{align}
											where for $j=5_2,5_3,5_{\infty}$, $\sigma_{j,j'}(c,s,r)$ is defined in \eqref{def:sigma_j,j'} and \eqref{def:sigma_j,j'52},
											satisfying \eqref{sigma-5infty1} and \eqref{sigma-5infty}. Therefore, for $j=5_2,5_3, 5_4$,
											\eqref{bd:I5-D} follows a similar process as \eqref{process:integral}; for $j=5_0,5_1$,
											it follows from Lemma \ref{lem: pifi-Linfty} and Lemma \ref{lem:symbol-chi}  that
											\begin{align}
												\begin{split}
													\label{process:integral5-D}I_{j,j',l,l',m}
													&\leq \int_{V(0)}^{1-\delta}r_c^{-\f12}\left|\left(p.v.\int_{\mathbb{R}} \left(r_c^{\f12}\pa_c\right)a_{j,j',l,l',m}(\eta_j,c,\lambda)\f{e^{i\eta_j y_j(c,\lambda) }}{\eta_j}d\eta_j\right)\right|dc\\
													&\quad+
													\int_{V(0)}^{1-\delta}\left|\left(p.v.\int_{\mathbb{R}} a_{j,j'}(\eta_j,c,\lambda)e^{i\eta_j y_j(c,\lambda) }d\eta_j\right)\right|\cdot \left|\f{\pa y_j(c,\lambda)}{\pa c}\right|dc\\
													&\lesssim  \int_{V(0)}^{1-\delta}r_c^{-\f12}\left\|\int_{\mathbb{R}} \left(r_c^{\f12}\pa_c\right)a_{j,j',l,l',m}(\eta_j,c,\lambda)e^{i\eta_j x}d\eta_j\right\|_{L^1_x(\mathbb{R})}dc\\
													&\quad+
													\int_{\mathbb{R}}\left|\int_{\mathbb{R}} a_{j,j',l,l',m}(\eta_j,c,\lambda)e^{i\eta_j y_j }d\eta_j\right|d y_j\\
													&\lesssim  \int_{V(0)}^{1-\delta}r_c^{-\f12}dc +1\lesssim 1.
											\end{split}\end{align}
											
										 It remains to prove \eqref{bd:a_j,j'5-D}. For simplicity, we use the notation $\mathcal{D}\in\{\mathcal{D}_{l,m},
										\mathcal{D}_{l',0}\}$, since the estimates are the same for different $l,l',m$.
										We first use the formula \eqref{fm:f-tf1} and  the definitions of $(\mathcal{D}\widetilde{\phi})_{j}$ in \eqref{def:tphi-5-D} to obtain the following bounds.
										Thanks to the bounds \eqref{6-phi-1} in Proposition \ref{prop:summery},
										\begin{align} \label{bd:50-D}
												\begin{split}
													&(\mathcal{D}\phi)_{5_0}(\eta_{5_0},c,s)=
													\left(1-\chi(\eta_{5_0}/M) \right) r_c^{\f12}
													O_{\eta_{5_0},c}^{\eta_{5_0},r_c}\left(\eta_{5_0}^{-1}\right)
													=\left(1-\chi(\eta_{5_0}/M) \right) O_{\eta_{5_0},c}^{\eta_{5_0},r_c^{\f12}}\left(\eta_{5_0}^{-1}\right), \\
													&r_c^{-1}(\mathcal{D}\phi)_{5_0}(\eta_{5_0},c,s) =r_c^{-\f12}
													O_{\eta_{5_0}}^{\eta_{5_0}}\left(\eta_{5_0}^{-1}\right). \end{split}
											\end{align}
									Thanks to the bounds \eqref{6-phi-2} in Proposition \ref{prop:summery},
										\begin{align}  \label{bd:51-D}
												\begin{split}
													&(\mathcal{D}\phi)_{5_1}(\eta_{5_1},c,s)=
													\left(1-\chi(\eta_{5_1}) \right)r_c^{\f12} O_{\eta_{5_1},c}^{\eta_{5_1},r_c}
													\left(\eta_{5_1}^{-\f12}\right)
													=\left(1-\chi(\eta_{5_1}) \right) O_{\eta_{5_1},c}^{\eta_{5_1},r_c^{\f12}}
													\left(\eta_{5_1}^{-\f12}\right), \\  &r_c^{-1}(\mathcal{D}\phi)_{5_1}(\eta_{5_1},c,s)=\left(1-\chi(\eta_{5_1}) \right) O_{\eta_{5_1},c}^{\eta_{5_1},r_c^{\f12}}
													\left(\eta_{5_1}^{-\f12}\right).
													\end{split}
											\end{align}
										
				Since the derivative bounds in \eqref{6-phi-4} and \eqref{6-f+a1} are the same as the non-derivative ones,
				it follows identically as in \eqref{bd:52-53}, \eqref{bd:5-rho-1}, \eqref{bd:54case2} , \eqref{bd:5-rho-11case2}, respectively that
										
										\begin{align}  \label{bd:523-D}
												\begin{split}
													&(\mathcal{D}\phi)_{5_2}(\eta_{5_2},c,s)=
													\left(1-\chi\left(\eta_{5_2}/M^{\f12}\right)\right) O_{\eta_{5_2},c}^{\eta_{5_2},\sigma_{5_2}(c,s)^{-1}}
													\left(\eta_{5_2}^{-\f16+\delta/9}\right), \\
													&\rho_{5_2}(c,s)^{-1}(\mathcal{D}\phi)_{5_2}(\eta_{5_2},c,s)=
													\left(1-\chi\left(\eta_{5_2}/M^{\f12}\right)\right) O_{\eta_{5_2}}^{\eta_{5_2}}
													\left(\eta_{5_2}^{-\f16+\delta/9}\right)\sigma_{5_2,j'}(c,s), \\
													&r_c^{-1}(\mathcal{D}\phi)_{5_2}(\eta_{5_2},c,s)=
													\left(1-\chi(\eta_{5_2}/M^{\f12})\right) O_{\eta_{5_2}}^{\eta_{5_2}}
													\left(\eta_{5_2}^{-\f16+\delta/9}\right)\sigma_{5_2,j'}(c,s), \\ & (\mathcal{D}\phi)_{5_3}(\eta_{5_3},c,s)=
													\left(1-\chi(\eta_{5_3}/M^{\f12})\right)
													O_{\eta_{5_3},c}^{\eta_{5_3},\sigma_{5_3}(c,s)^{-1}}
													\left(\eta_{5_3}^{-\f13+\delta/9}\right), \\
													&\rho_{5_3}(\xi)^{-1}(\mathcal{D}\phi)_{5_3}(\eta_{5_3},c,s)=
													\left(1-\chi\left(\eta_{5_3}/M^{\f12}\right)\right)
													O_{\eta_{5_3}}^{\eta_{5_3}}
													\left(\eta_{5_3}^{-\f13+\delta/9}\right)\sigma_{5_3}(c,s), \\
													&r_c^{-1}(\mathcal{D}\phi)_{5_3}(\eta_{5_3},c,s)=
													\left(1-\chi\left(\eta_{5_3}/M^{\f12}\right)\right)
													O_{\eta_{5_3}}^{\eta_{5_3}}
													\left(\eta_{5_3}^{-\f13+\delta/9}\right)\sigma_{5_3}(c,s),\\
													&(\mathcal{D}\phi)_{5_{\infty}}(\eta_{5_{\infty}},c,s)=
													\sum_{\iota\in\{-,+\}}\left(1-\chi\left(M^{\f12}\eta_{5_{\infty}}\right)\right) O_{\eta_{5_{\infty}},c}^{\eta_{5_{\infty}},\sigma_{5_{\infty}}(c,s)^{-1}}\left(\eta_{5_{\infty}}^{-\f16}\right)e^{\iota i \eta_{5_{\infty}}}, \\
													&\rho_{5_{\infty}}(c,s)^{-1}
													(\mathcal{D}\phi)_{5_{\infty}}(\eta_{5_{\infty}},c,s)=
													\sum_{\iota\in\{-,+\}}\left(1-\chi\left(M^{\f12}\eta_{5_{\infty}}\right)\right) O_{\eta_{5_{\infty}}}^{\eta_{5_{\infty}}}\left(\eta_{5_{\infty}}^{-\f16}\right)
													\sigma_{5_{\infty}}(c,s) e^{\iota i \eta_{5_{\infty}}},\\
													&r_c^{-1}
													(\mathcal{D}\phi)_{5_{\infty}}(\eta_{5_{\infty}},c,s)=
													\sum_{\iota\in\{-,+\}}\left(1-\chi(M^{\f12}\eta_{5_{\infty}})\right) O_{\eta_{5_{\infty}}}^{\eta_{5_{\infty}}}\left(\eta_{5_{\infty}}^{-\f16}\right)
													\sigma_{5_{\infty}'}(c,s) e^{\iota i \eta_{5_{\infty}}}.
												\end{split}
											\end{align}
										
										For $\phi_{j'}^*$, it follows similarly as  in \eqref{bi-bd:5inftycase2}, \eqref{bi-bd:50-5infty}  and \eqref{bd:tphi5-case4'} that for $j\in\{5_0,5_1,5_2,5_3,5_{\infty}\}$,
										\begin{align}&  \label{bi-bd:50-5infty-D}
											\begin{split}
												&\phi_{j'}^*(\eta_{j},c,s,r)= O_{\eta_{j},c}^{\eta_{j},r_c}(1)
												+O_{\eta_{j},c}^{\eta_{j},\rho_j(c,s,\xi)}(1),\quad j'=5_0,5_1,\\
												&\phi_{j'}^*(\eta_{j},c,s,r)= O_{\eta_{j},c}^{\eta_{j},\sigma_{j'}(c,r)^{-1}}(1)
												+O_{\eta_{j},c}^{\eta_{j},\rho_j(c,s,\xi)}(1),\quad j'=5_2,5_3,5_{\infty}. \end{split}
										\end{align}
										For the cut-off function, it follows similarly as in \eqref{behave:chi5} and \eqref{behave:chi5case2} that
										\begin{align} \label{behave:chi5-D}
											\begin{split}
												&\chi\left(\f{(1-c)^{\f12}}{\xi}\right)\left(1- \chi\left(\xi (V(0)-c)^{\f32}/M\right) \right)\\
												&=
												\left\{
												\begin{array}{l}
													O_{\eta_j,c}^{\eta_j,r_c}(1)\;\;\; \;\;\quad\quad \quad\quad\quad\quad\quad \quad j=5_0,5_2,5_3,5_4,\\ O_{\eta_{5_{\infty}},c}^{\eta_{5_{\infty}},
														r_c}(1)+O_{\eta_{5_{\infty}},c}^{\eta_{5_{\infty}},
														\rho_{5_{\infty}}(c,s)}(1)\;\;\; \;\;\quad j=5_{\infty}.\end{array}\right.
											\end{split}
										\end{align}
										Finally, we conclude  the claim \eqref{bd:a_j,j'5-D} by \eqref{bd:50-D}, \eqref{bd:51-D}, \eqref{bd:523-D}, \eqref{bi-bd:50-5infty-D}, \eqref{behave:chi5-D}.
										
										If   \emph{$(j,j')\in\{5_0,5_1,5_2,5_3\}\times\{5_{\infty}\}$}, this case is symmetric to $(j,j')\in\{5_{\infty}\}\times \{5_0,5_1,5_2,5_3\}$.
										We thus interchange the roles $s$ and $r$ in the proof of  \eqref{op-pieces5-D} for $(j,j')\in\{5_{\infty}\}\times \{5_0,5_1,5_2,5_3\}$.
										Specifically, we set $\eta_{5_{\infty}}= \xi \int_{r_c}^r\sqrt{\f{V(s')-c}{1-c}}ds'$.
										The details are left to the readers.
										
										For \emph{$(j,j')\in\{5_{\infty}\}\times\{5_{\infty}\}$}, the  proof of  \eqref{op-pieces5-D} is identical to that of case 3 in Proposition \ref{lem:Boundness-K-pieces5},
										as the derivative bounds coincide with the non-derivative bounds in \eqref{6-f+} (given in Proposition \ref{prop:summery}).										\end{proof}
										
										\begin{proposition}\label{lem:Boundness-K-pieces6-D}
											Let  $\delta\ll 1$  and $M\gg 1$ be fixed. Let  $l',l\in\{0,1\}$, $m\in \mathbb{N}$,  $j\in J^H_{6}=
											\{6_0,6_1,6_2,6_3,6_4,6_{\infty}
											\} $. Let  $(\mathcal{D}_{l,m}\widetilde{\phi})_j$ be defined in  \eqref{def:tphi-6-D}, Definition \ref{def:decom-phi-D} .
											Then for $l'+l+m\geq 1$, $(j,j')\in J^H_{6}\times J^H_{6}$,  it holds uniformly for  $(r,s,z)\in\mathbb{R}^+\times \mathbb{R}^+\times \mathbb{R}$  that
												\begin{align}
												\begin{split}
												 \int_{1-\delta}^1&\left|\pa_c\int_{\mathbb{R}} \chi\left(\f{M^2(1-c)^{\f12}}{\xi}\right)\left(1- \chi\left(\xi /M(1-c)^{\f13}\right) \right)\right.\\&\quad\left.
													(\mathcal{D}_{l,0}\widetilde{\phi})_{j}(\xi,c,s)
													(\mathcal{D}_{l',m}\widetilde{\phi})_{j'}(\xi,c,r)\f{e^{i\xi \sqrt{\f{c}{1-c}}z}}{\xi}d\xi\right|dc\lesssim s^{-\delta}+1.
													\end{split}
													\label{op-pieces6-D}\end{align}
											
										\end{proposition}
										
										\begin{proof}
										Without further illustration, we assume $c\in(1-\delta,1),\;\xi r_c\gtrsim M$, which implies $(1-c)^{-\f13}\sim r_c\gtrsim  1$.
										We take the same coordinate $\eta_{6_0}$ as \eqref{trans:eta6} that   \begin{align*} 
												\begin{split}
													\eta(\xi,c,s)=
													\left\{
													\begin{array}{l} \f{ \xi s}{(1-c)^{\f12}} = \eta_{6_0}\;\;\;\quad\quad s\lesssim M^{\f12}\xi^{-1}r_c^{-\f32},\\
														\f{ \xi s}{(1-c)^{\f12}} =\eta_{6_1}\;\;\;\quad\quad M^{\f12}\xi^{-1}r_c^{-\f32}\lesssim s\leq \f12,\\
														\f{\xi}{(1-c)^{\f13}}:=
														\eta_{6_2}\;\;\quad\quad\quad \f12\leq s\leq r_c/2,\\
														\f{\xi}{(1-c)^{\f13}}=  \eta_{6_3}\;\;\;\quad\;\; \f{r_c}{2}\leq s\leq r_c-C\xi^{-\f23}r_c^{\f13},\\
														\f{\xi}{(1-c)^{\f13}}=  \eta_{6_4}\;\;\; \quad\quad  r_c-C\xi^{-\f23}r_c^{\f13}\leq s\leq r_c+C\xi^{-\f23}r_c^{\f13},\\
														\xi \int_{r_c}^s\sqrt{\f{V(s')-c}{1-c}}ds'= \eta_{6_{\infty}} \;\;\;\quad\quad s\geq r_c+C\xi^{-\f23}r_c^{\f13},
														\end{array}\right.
												\end{split}
											\end{align*}
											whose inverse function is denoted as $\xi(\eta_{j},c,s)$.
											Moreover, $\eta_j$ satisfies  the bounds \eqref{est:paceta/eta6}, with the same weight as \eqref{recall:rho-6}:
											\begin{align*} \rho(c,\xi,s)=
												\left\{
												\begin{array}{l}
													\rho_{6_0}=\rho_{6_1}=\rho_{6_2}=r_c^{-3}\;\;\; \;\;\quad\quad s\leq \f{r_c}{2},\\
													\rho_{6_3}=\f{r_c-s}{r_c^4} \;\;\; \;\quad\quad r_c/2\leq s\leq r_c-C\xi^{-\f23}, \\
													\rho_{6_4}=\xi^{-\f23}r_c^{-\frac{11}{3}}\;\;\; \quad\quad r_c-C\xi^{-\f23}r_c^{\f13}\leq s\leq r_c+C\xi^{-\f23}r_c^{\f13},\\
													\rho_{6_{\infty}}=\f{s-r_c}{r_c^4} \;\;\; \quad\quad s\geq r_c+C\xi^{-\f23}r_c^{\f13} .
												\end{array}\right.
											\end{align*}
									        If \emph{$j'\neq 6_{\infty}$}, i.e,   $(j,j')\in  \big\{6_0,6_1,6_2,6_3,6_4,6_{\infty}
											\big\}\times  \big\{6_0,6_1,6_2,6_3,6_4\big\}$, we make the transform such that
											\begin{align*}  \begin{split}
													&(\mathcal{D}_{l',0}\phi)_{j}(\eta_{j},c,s)
													:=(\mathcal{D}_{l',0}\widetilde{\phi})_{j}\left(\xi(\eta_{j},c,s),s,c\right), \\&(\mathcal{D}_{l,m}\phi)^*_{j'}(\eta_{j},c,s,r)
													:=(\mathcal{D}_{l,m}\widetilde{\phi})^*_{j'}
													\left(\xi(\eta_{j'},c,s),r,c\right).
												\end{split}
											\end{align*}
											Then we  reduce
											\eqref{op-pieces6-D} to
											\begin{align}
												\label{bd:I6-D} I_{j,j',l,l',m}&:=\int_{1-\delta}^{1}\left|\int_{\mathbb{R}} a_{j,j',l,l',m}(\eta_{j},c,\lambda)\f{e^{i\eta_{j} y_{j}(c,\lambda) }}{\eta_{j}}d\eta_{j}\right|dc\lesssim s^{-\delta},\;\;j=6_0,6_1, \\
												I_{j,j',l,l',m}&:=\int_{1-\delta}^{1}\left|\int_{\mathbb{R}} a_{j,j',l,l',m}(\eta_{j},c,\lambda)\f{e^{i\eta_{j} y_{j}(c,\lambda) }}{\eta_{j}}d\eta_{j}\right|dc\lesssim 1,\;\;j=6_2,6_3,6_4,6_{\infty},
												\label{bd:I6-D-34infty}
											\end{align}
											where $\lambda=(z,s,r)$ is the parameter and
											\begin{align}\label{def:y_j6-D}
											y_j(c,\lambda)=
												\left\{
												\begin{array}{l}
													z\sqrt{\f{c}{1-c}}\cdot(1-c)^{\f13}\;\;\; \;\;\quad\quad\quad\quad\quad j=6_2,6_3,6_4,\\
													z\sqrt{\f{c}{1-c}}/\f{s}{(1-c)^{\f12}}\;\;\; \;\;\quad\quad\quad\quad\quad\quad j=6_0,6_1,\\
													z\sqrt{\f{c}{1-c}}/\int_{r_c}^s\sqrt{\f{V(s')-c}{1-c}}ds'\;\;\; \;\;\quad\quad j=6_{\infty},\end{array}
													\right.
												\end{align}
											and  \begin{align*} 
												&\notag a_{j,j',l,l',m}(\eta_{j},c,\lambda)\\
												&:=\chi\left(\f{(1-c)^{\f12}}{\xi}\right) \left(1-\chi\left(\xi/M(1-V(0))^{\f13}\right)\right)(\mathcal{D}_{l,0}\phi)_{j}
												(\eta_{j},c,s)
												(\mathcal{D}_{l',m}\phi)_{j'}^* (\eta_{j},c,s,r).
											\end{align*}
											We claim that
											\begin{align} \label{bd:a_j,j'6-D}
												\begin{split}
													a_{j,j'}(\eta_j,c,\lambda) &=
													\left\{
													\begin{array}{l}
														\chi(\eta_{6_0}/M)
														O_{\eta_{6_0},c}^{\eta_{6_0},r_c^{-3+\f{\delta}{2}}}
														\left(\eta_{6_0}^{\delta}\right)\cdot s^{-\delta},\;\;\;\quad\quad\quad\quad\quad\quad\quad j=6_0, \\
														\left(1-\chi(\eta_{6_1}/M) \right) O_{\eta_{6_1},c}^{\eta_{6_1},r_c^{-3+\f{\delta}{2}}}
														\left(\eta_{6_1}^{-\f12+\delta}\right)\cdot s^{-\delta},
														\;\;\quad\quad \quad\quad  j=6_1,\\
														\left(1-\chi\left(\eta_{6_2}/M^{\f12}\right)\right) O_{\eta_{6_2},c}^{\eta_{6_2},\sigma_{6_2,j'}(c,s,r)^{-1}}
														\left(\eta_{6_2}^{-\f12}\right),\;\;\; \quad\quad\quad j=6_2,\\
														\left(1-\chi\left(\eta_{6_3}/M^{\f12}\right)\right) O_{\eta_{6_3},c}^{\eta_{6_3},\sigma_{6_3,j'}(c,s,r)^{-1}}
														\left(\eta_{6_3}^{-\f13+\delta/9}\right),\;\;\; \quad\quad j=6_3,\\
														\left(1-\chi\left(\eta_{6_4}/M^{\f12}\right)\right) O_{\eta_{6_4},c}^{\eta_{6_4},\sigma_{6_4,j'}(c,s,r)^{-1}}
														\left(\eta_{6_4}^{-\f13+\delta/9}\right),\;\;\; \quad\quad j=6_4,\\
														\sum_{\iota\in\{-,+\}}\left(1-\chi\left(M^{\f12}\eta_{6_{\infty}}\right)\right) O_{\eta_{6_{\infty}},c}^{\eta_{6_{\infty}},\sigma_{6_{\infty},j'}
															(c,\lambda)
															^{-1}}\left(\eta_{6_{\infty}}^{-\f16}\right)e^{\iota i \eta_{6_{\infty}}},\quad j=6_{\infty}, \end{array}\right.
												\end{split}
											\end{align}
											where  for $j=6_2,6_3, 6_4, 6_{\infty}$, $\sigma_{j,j'}(c,s,r)$ are defined in \eqref{def:sigma_j,j'6} and \eqref{def:sigma_j,j'62},
											satisfying \eqref{sigma-6infty1} and \eqref{sigma-6infty}, respectively.
											Therefore, for $j=6_3,6_4,6_{\infty}$, \eqref{bd:I6-D-34infty} follows a similar process as in \eqref{process:integral}.
											For $j=6_0,6_1$, \eqref{bd:I5-D}  follows from Lemma \ref{lem: pifi-Linfty} and Lemma \ref{lem:symbol-chi}  that \begin{align}
												\begin{split}
													\label{process:integral6-D}I_{j,j',l,l',m}&\leq s^{-\delta}\int_{1-\delta}^1r_c^{3-\delta/2}\left|\left(p.v.\int_{\mathbb{R}} \left(r_c^{-3+\f{\delta}{2}}\pa_c\right)a_{j,j',l,l',m}(\eta_j,c,\lambda)\f{e^{i\eta_j y_j(c,\lambda) }}{\eta_j}d\eta_j\right)\right|dc\\
													&\quad+
													s^{-\delta}\int_{1-\delta}^1\left|\left(p.v.\int_{\mathbb{R}} a_{j,j'}(\eta_j,c,\lambda)e^{i\eta_j y_j(c,\lambda) }d\eta_j\right)\right|\cdot \left|\f{\pa y_j(c,\lambda)}{\pa c}\right|dc\\
													&\lesssim s^{-\delta} \int_{1-\delta}^1r_c^{3-\delta/2}\left\|\int_{\mathbb{R}} \left(r_c^{-3+\f{\delta}{2}}\pa_c\right)a_{j,j',l,l',m}(\eta_j,c,\lambda)e^{i\eta_j x}d\eta_j\right\|_{L^1_x(\mathbb{R})}dc\\
													&\quad+
													s^{-\delta}\int_{\mathbb{R}}\left|\int_{\mathbb{R}} a_{j,j',l,l',m}(\eta_j,c,\lambda)e^{i\eta_j y_j }d\eta_j\right|d y_j\\
													&\lesssim s^{-\delta} \int_{1-\delta}^1r_c^{3-\delta/2}dc +s^{-\delta}\lesssim s^{-\delta} \int_{V^{-1}(1-\delta)}^{+\infty}r_c^{-\delta/2-1}dr_c +s^{-\delta} \lesssim s^{-\delta}.
											\end{split}\end{align}
											
										It remains to prove \eqref{bd:a_j,j'6-D}. Indeed, the first and second lines follow a similar argument to that in \eqref{bd:30-D}, based on the following observation for $j=6_0,6_1$:
										\begin{align}  \label{observation:60-61}
											1=(\xi r_c^{\f32}s)^{\delta}(\xi r_c)^{-\delta}r_c^{-\f{\delta}{2}}s^{-\delta}=O_{\eta_j}^{\eta_j}(\eta_j^{\delta})\cdot r_c^{-\f{\delta}{2}}\cdot s^{-\delta}.
										\end{align}
										The remaining cases follow similarly to the bounds established in the proof of Proposition \ref{lem:Boundness-K-pieces6},
										owing to the observation that the derivative bounds coincide with the non-derivative bounds
										in \eqref{6-phi-3}, \eqref{6-phi-4}, \eqref{6-f+a1}, and \eqref{6-f+} (given in Proposition \ref{prop:summery}).  The details are left to the readers.
										
										If   \emph{$(j,j')\in\{5_0,5_1,5_2,5_3\}\times\{5_{\infty}\}$}, this case is symmetric to $(j,j')\in\{5_{\infty}\}\times \{5_0,5_1,5_2,5_3\}$.
										We thus interchange the roles  $s$ and $r$ in the proof of  \eqref{op-pieces5-D} for $(j,j')\in\{5_{\infty}\}\times \{5_0,5_1,5_2,5_3\}$.
										Specifically, we set $\eta_{5_{\infty}}= \xi \int_{r_c}^r\sqrt{\f{V(s')-c}{1-c}}ds'$. The details are left to the readers.

										If   \emph{$(j,j')\in\{6_0,6_1,6_2,6_3,6_4\}\times\{6_{\infty}\}$}, this case is symmetric  to $(j,j')\in\{6_{\infty}\}\times \{6_0,6_1,6_2,6_3,6_4\}$.
										We thus interchange the roles $s$ and $r$ in the proof of  \eqref{op-pieces6-D} with $(j,j')\in\{6_{\infty}\}\times \{6_0,6_1,6_2,6_3,6_4\}$. Specifically, we set $\eta_{6_{\infty}}= \xi \int_{r_c}^r\sqrt{\f{V(s')-c}{1-c}}ds'$. The details are left to the readers.
										
										If \emph{$(j,j')\in\{6_{\infty}\}\times\{6_{\infty}\}$}, the  proof of  \eqref{op-pieces6-D} is identical to that of case 3 in Proposition \ref{lem:Boundness-K-pieces6}, since we observe that the derivative bounds are the same as the  non-derivative bounds in \eqref{6-f+} in Proposition \ref{prop:summery}.
										\end{proof}
										
										\begin{proposition}\label{lem:Boundness-K-pieces7-D}
											Let  $\delta\ll 1$ such that $1-\delta\geq V(1)$ and $M\gg 1$ be fixed, and $J^H_{7}=
											\{7_0,7_1,7_2,7_{\infty}
											\}$. Let  $(\mathcal{D}_{l,m}\widetilde{\phi})_j$ be defined in \eqref{def:tphi-7-D} and Definition \ref{def:decom-phi-D} .
											Then for $l'+l+m\geq 1$, $(j,j')\in J^H_{7}\times J^H_{7}$,  it holds uniformly for  $(r,s,z)\in\mathbb{R}^+\times \mathbb{R}^+\times \mathbb{R}$ that
												\begin{align}
												\begin{split}
												\int_{1-\delta}^1&\left|\pa_c\int_{\mathbb{R}} \chi\left(\f{M^2(1-c)^{\f12}}{\xi}\right)\chi\left(\xi/(1-c)^{\f13}M\right)\right.\\&\quad\left.
													(\mathcal{D}_{l,0}\widetilde{\phi})_{j}(\xi,c,s)(\mathcal{D}_{l',m}\widetilde{\phi})_{j'}(\xi,c,r)\f{e^{i\xi \sqrt{\f{c}{1-c}}z}}{\xi}d\xi\right|dc
													\lesssim s^{\delta}+1.\end{split} \label{op-pieces7-D}.
													\end{align}
										
										\end{proposition}
										
										\begin{proof}
											The proof is similar to that of Proposition \ref{lem:Boundness-K-pieces6}, since the derivative bounds in \eqref{7-phi-D} are similar to those in \eqref{7-phi},
											and the derivative bound of \eqref{7-f+} is identical to the non-derivative one.
										\end{proof}

\appendix

\section{Bounds on potential functions}
							Let $c\in(0,1)$, $r\in\mathbb{R}^+$. We define
							\begin{align}\label{def:Q}
								Q(r,c):=\f{V(r)-c}{1-c}.
							\end{align}
							For  $c\in[V(0),1)$,  we denote $r_c:=V^{-1}(c)$ and define
							\begin{align}\label{def:tau}
								&\tau(r,c):=
								\left\{
								\begin{aligned}
									&\left(\f32\int_0^r\sqrt{Q(s,c)}ds\right)^{\f23},\;c\in(0,V(0)],\\
									&\mathrm{sgn}(r-r_c)\left(\f32\int_{r_c}^r\sqrt{|Q(s,c)|}ds\right)^{\f23},\;c\in(V(0),1),
								\end{aligned}
								\right.
							\end{align}
							and \begin{align*}q(r,c)=\f{Q(r,c)}{\tau(r,c)}.\end{align*}
						 Note that for $c\in (V(0),1)$, due to $V(r)=1-a_0^2 r^{-3}+O(r^{-4})(r\to+\infty)$ and $V'>0$, we infer that
							\begin{align}\label{rc-behave}
								1-c\sim \langle r_c\rangle^{-3}\; \text{and}\; V'(r_c)\sim \langle r_c\rangle^{-4}.\end{align}
							\begin{lemma}\label{lem:behave-Q}
								Let $i\in\mathbb{N}$, $l\in\{0,1\}$, $c\in(0,1)$, $r>0$ . We  have the following precise behavior  for $Q$:
								\begin{itemize}
									\item For $c\in(0,V(0)]$, $r>0$, it holds that
									\begin{align} \label{behave:Q-V(0)}
										\begin{split}
											&Q(r,c)\sim \f{r}{\langle r\rangle}+(V(0)-c),
										\end{split}
									\end{align}
									and
									\begin{align}\label{estQ-G-0} &
										\left|\left(V(r)-c\right)^l
										r^i\pa_c^l\pa_r^iQ\right|\lesssim |Q|.
									\end{align}
									
									\item For $c\in(V(0),1)$, $r\gtrsim r_c$, it holds  that        									
									\begin{align} \label{behave:Q<rc0}
										&Q(r,c)\sim \f{r-r_c}{\langle r\rangle}.
									\end{align}
									\item
									For $c\in(V(0),1)$, $r\lesssim r_c$, it holds that
									\begin{align}\label{behave:Q<rc}
										\begin{split}
											&Q(r,c)\sim \f{(r-r_c)\langle r_c\rangle^2}{\langle r\rangle ^3},\;r\lesssim r_c.
										\end{split}
									\end{align}
									\item
									For $c\in(V(0),1)$ and $i+l\geq 1$, it holds that for $r>0$,
									\begin{align}
										\label{est-der-Q}
										&\left|\langle r_c\rangle ^{-3l}\pa_c^l\pa_r^iQ(r,c)\right|\lesssim \langle r_c\rangle ^{3}\langle r\rangle ^{-3-i}.
									\end{align}
										which along with \eqref{behave:Q<rc0}, \eqref{behave:Q<rc} implies
									\begin{align}\label{estQ-G}  &
										\left|\left(\f{|r_c-r|}{\langle r_c\rangle^{4}}\right)^l
										\left(\f{|r_c-r|\langle r\rangle}{\langle r_c\rangle}\right)^i\pa_c^l\pa_r^iQ\right|\lesssim |Q|.
									\end{align}
								\end{itemize}
								The following facts hold for $\tau$:
								\begin{align}
									\begin{split}
										\label{est-tleqQ}
										&\tau^{\f32}(r,c)=\f32\int_0^{r}Q^{\f12}(r',c)dr'\sim \f32 rQ^{\f12}(r,c),\;\;c\in(0,V(0)].\\
										&\tau^{\f32}(r,c)=\f32\left|\int_{r_c}^{r}|Q|^{\f12}(r',c)dr'\right|\lesssim |r-r_c||Q|^{\f12}(r,c),\;\;c\in(V(0),1).\\
										&\tau^{\f32}(r,c)=\f32\left|\int_{r_c}^{r}|Q|^{\f12}(r',c)dr'\right|\sim |r-r_c|Q^{\f12}(r,c),\;\;c\in(V(0),1),\;r\geq r_c.
									\end{split}
								\end{align}
							\end{lemma}
							\begin{proof}
								For $c\in(0,V(0)]$, since $V'(r)\sim \langle r\rangle ^{-4}$, we deduce that $V(r)-V(0)\sim r (r\lesssim 1)$, $V(r)-V(0)\sim 1 (r\gtrsim 1)$.
								Then \eqref{behave:Q-V(0)} follows from $V(r)-c=V(r)-V(0)+V(0)-c\sim \f{r}{1+r}+(V(0)-c)$, and \eqref{estQ-G-0} follows directly.
								
								For $c\in(V(0),1)$,  we first prove derivative bounds of $Q$ in \eqref{est-der-Q}. We write $Q(r,c)=
								1+\f{V(r)-1}{1-c}$. It follows directly that $|(1-c)\pa_cQ=\f{V(r)-1}{1-c}|\lesssim \f{\langle r\rangle ^{-3}}{1-c}$ and
								$\left|\left((1-c)\pa_c\right)^l\pa_r^iQ=\f{V^{(i)}(r)}{1-c}\right|\lesssim \f{\langle r\rangle ^{-3-i}}{1-c}(i\geq 1)$. We write
								\begin{align}\label{def:g1}
									Q(r,c)&=(r-r_c)g_1(r,c),\;g_1(r,c):=\int_{0}^1\f{V'(r_c+s(r-r_c))}{1-c}ds,
								\end{align}
								then the remaining two situations follow by
								\begin{align}\label{g1-behave}
									g_1(r,c)\sim \langle r_c\rangle^3\left|\f{\langle r\rangle^{-3}-\langle r_c\rangle^{-3}}{r-r_c}\right|\sim \f{\langle r+r_c\rangle^{2}}{\langle r\rangle^{3}}\sim
									\left\{
									\begin{aligned}
										&\f{1}{\langle r\rangle },\;\text{if}\;\;r\gtrsim r_c;\\
										&\f{\langle r_c\rangle^2}{\langle r\rangle^3 },\;\text{if}\;r\lesssim r_c.
									\end{aligned}
									\right.
								\end{align}
								The upper bound of $|\tau|^{\f32}$ in \eqref{est-tleqQ} follows by the monotonicity of $Q$ in $r$. We will prove the lower bounds.  For $c\in(V(0),1)$, it can be obtained easily by the behavior of $Q$. For $c\in(V(0),1)$, using the behavior of $Q$, it suffices to prove that
								\begin{align*}
									f_1(r)=\int_{r_c}^r\f{(s-r_c)^{\f12}}{\langle s\rangle^{\f12}}ds\gtrsim
									\f{(r-r_c)^{\f32}}{\langle r\rangle^{\f12}},\quad r\geq r_c.
								\end{align*}
								In fact, integration by parts yields that for $r\geq r_c$,
								\begin{align*}
									f_1(r)=\f23\f{(r-r_c)^{\f32}}{\langle r\rangle^{\f12}}+\f13\int_{r_c}^r\f{(s-r_c)^{\f32}}{\langle s\rangle^{\f32}}ds\gtrsim
									\f{(r-r_c)^{\f32}}{\langle r\rangle^{\f12}}.
								\end{align*}
								Thus, the lemma is proved.
							\end{proof}
						
							\begin{lemma}\label{lem:behave-x}
								Let  $\tau$ be defined in \eqref{def:tau} and
								\begin{align}\label{def:x}
									x(r,c,\xi):=\xi\mathrm{sgn}(r-r_c)\tau^{\f32}(r,c)=
									\f32\xi\int_{r_c}^r|Q|^{\f12}(s,c)ds.
								\end{align}
								It holds that \begin{align*}
									\xi\pa_{\xi}x=x,\quad \pa_rx=\f32\xi |Q|^{\f12}\geq 0,
									\end{align*}
								and for $i\in\mathbb{N}$, $l\in\{0,1\}$ and $r_c=V^{-1}(c)$:
								
								\begin{itemize}
									\item For $c\in(0,V(0))$, $r>0$, it holds that
									\begin{align}
										\begin{split}
											\label{estx1-3}
											&x\sim \xi\left(\f{r^{\f32}}{\langle r\rangle ^{\f12}}+(V(0)-c)^{\f12}r\right),\quad \left|(V(r)-c)^lr^i\pa_c^l\pa_r^ix\right|
											\lesssim |x|.
										\end{split}
									\end{align}
									
									\item For $c\in(V(0),1)$, $r\geq r_c/4$, it holds that
									\begin{align}\label{estx4-6-G} \begin{split}
											&|x|\sim \xi |r-r_c|^{\f32}/\langle r\rangle ^{\f12}\;\text{and}\;\;
											\left|\left(\f{|r-r_c|}{\langle r_c\rangle^{4}}\right)^l|r-r_c|^{i}\pa_c^l\pa_r^ix\right|\lesssim |x|.
										\end{split}
									\end{align}
									
									\item For $c\in(V(0),1)$, $r\leq r_c$, it holds that
									\begin{align}\label{estx-G}
										\begin{split}
											&|x|\gtrsim \xi |r-r_c|^{\f32}/\langle r\rangle ^{\f12}\;\text{and}\;\;
											\left|\left(\f{(r_c-r)}{\langle r_c\rangle^{4}}\right)^l
											\left(\f{(r_c-r)\langle r\rangle}{\langle r_c\rangle}\right)^i\pa_c^l\pa_r^ix\right|\lesssim |x|.
										\end{split}
									\end{align}
								\end{itemize}
							\end{lemma}
							
							\begin{proof}
								We first prove \eqref{estx1-3}. By the definition,  $x=\f32\xi \int_{0}^r|Q|^{\f12}(s,c)ds\sim \xi r |Q|^{\f12}\sim \xi\left(\f{r^{\f32}}{(1+r)^{\f12}}+(V(0)-c)^{\f12}r\right) $,
								where we used
								\begin{align*}
									0<\f{rV'(r)}{V(r)-c}\lesssim 1,\quad c\in(0,V(0)).
									\end{align*}
								For the derivative bounds, if $i=0$, using \eqref{est-der-Q}, it suffices to  check that
								 \begin{align*}
									f(r)=\int_{0}^r(V(s)-c)^{-\f12}\langle s\rangle ^{-3}ds
									\lesssim r(V(r)-c)^{-\f12}.
									\end{align*}
									Indeed, using $V(r)-c\sim \f{r}{r+1}+V(0)-c$, for $r\lesssim V(0)-c$,
									$f(r)\lesssim \int_0^r(V(0)-c)^{-\f12}ds\lesssim r(V(r)-c)^{-\f12}$;
									for $V(0)-c\lesssim r\lesssim 1$,
									$f(r)\lesssim \int_0^rs^{-\f12}ds\sim r^{\f12}\lesssim r(V(r)-c)^{-\f12}$;
									for $V(0)-c\lesssim r\lesssim 1$,
									$f(r)\lesssim \int_0^1s^{-\f12}ds
									+\int_1^rs^{-3}ds\lesssim 1\lesssim r(V(r)-c)^{-\f12}$.
									For $i\geq 1,\;l\in \{0,1\}$, it follows
									from Fa$\grave{\text{a}}$di Bruno's formula, Leibnitz's rule and $|\pa_c^l\pa_r^iQ|\lesssim \langle r\rangle ^{-3-i}$ that
									\begin{align}
										\notag\left|\pa_c^l\pa_r^ix\right|&=\f32\xi \left|\pa_c^l\pa_r^{i-1}(Q^{\f12})\right|
										\lesssim \xi\sum_{\substack{k_1,...,k_{i-1}\geq 0,\;k_1+...+(i-1)k_{i-1}=i-1\\
												k=k_1+...+k_{i-1}}}\left|\pa_c^l\left(Q^{\f12-k}(\pa_rQ)^{k_1}
										...(\pa_r^{i-1}Q)^{k_{i-1}}\right)\right|\\
										\notag&\lesssim \xi\sum_{\substack{k_1,...,k_{i-1}\geq 0,\;k_1+...+(i-1)k_{i-1}=i-1\\
												k=k_1+...+k_{i-1}}}\left|Q^{\f12-k-l}
										(\pa_c^lQ)^{l}(\pa_rQ)^{k_1}
										...(\pa_r^{i-1}Q)^{k_{i-1}}\right|\\
										&\quad+\xi\sum_{\substack{k_1,...,k_{i-1}\geq 0,\;k_1+...+(i-1)k_{i-1}=i-1\\
												k=k_1+...+k_{i-1}}}\left|Q^{\f12-k}
										\pa_c^{l}\left((\pa_rQ)^{k_1}
										...(\pa_r^{i-1}Q)^{k_{i-1}}\right)\right|\notag\\
										&\lesssim  \xi\sum_{\substack{k_1,...,k_{i-1}\geq 0,\;k_1+...+(i-1)k_{i-1}=i-1\\
												k=k_1+...+k_{i-1}}}\left|Q^{\f12-k-l}\right|\cdot \langle r\rangle ^{-3l-(3+1)k_1-...-(3+i-1)k_{i-1}}k\notag\\
										\notag&\quad+\xi\sum_{\substack{k_1,...,k_{i-1}\geq 0,\;k_1+...+(i-1)k_{i-1}=i-1\\
												k=k_1+...+k_{i-1}}}\left|Q^{\f12-k}\right|\cdot  \mathrm{sgn}(k)\langle r\rangle ^{-(3+1)k_1-...-(3+i-1)k_{i-1}}\\
										&\lesssim \xi\sum_{\substack{0\leq k\leq i-1}}\left|Q^{\f12-k-l}\right|\cdot\left(\langle r\rangle ^{-3l-3k-(i-1)}+\mathrm{sgn}(k)|Q|^l\langle r\rangle ^{-3k-(i-1)}\right)\label{Faadi:x},
									\end{align}
									where in the last line we used the  observation  $i-1=k_1+...+(i-1)k_{i-1}\geq k=k_1+...+k_{i-1}\geq 0$.
									Consequently, using $Q(r,c)\gtrsim  \f{r}{\langle r\rangle }$ and $|x|\sim \xi r|Q|^{\f12}$,  we have
									\begin{align*}
										\left|(V(r)-c)^{l}r^i\pa_c^l\pa_r^ix\right|&\lesssim \xi r|Q|^{\f12}r^{i-1}\sum_{\substack{0\leq k\leq i-1}}|Q|^{-k}\cdot\left(\langle r\rangle ^{-3l-3k-(i-1)}+\mathrm{sgn}(k)|Q|^l\langle r\rangle ^{-3k-(i-1)}\right)\\
										&\lesssim |x|\sum_{\substack{0\leq k\leq i-1}}r^{i-1}\f{\langle r\rangle ^{k}}{r^{k}}\cdot\left(\langle r\rangle ^{-3l-3k-(i-1)}+\mathrm{sgn}(k)\langle r\rangle ^{-3k-(i-1)}\right)\lesssim |x|.
									\end{align*}
									Next, we prove \eqref{estx4-6-G}.
									By the definition, \eqref{g1-behave}, we infer that for $r_c\lesssim 1$, $r>0$ and $r_c\gtrsim 1$, $r\gtrsim r_c$,
									\begin{align}\label{equiv:x1}
										|x|&=\xi |\tau|^{\f32}\sim \xi |r-r_c|^{\f32}\Phi_1(r,c)\sim
										\xi |r-r_c|^{\f32}\int_0^1\f{s^{\f12}}{(1+(1-s)r_c+sr)^{\f12}}ds\sim
										\xi |r-r_c|^{\f32}\langle r\rangle ^{-\f12},
									\end{align}
									where we have defined
									\begin{align}\label{def:Phi1 pre}
										\Phi_1(r,c)=\int_0^1s^{\f12}g_1^{\f12}
										((1-s)r_c+sr,c)ds.
									\end{align}
									and used $I(r):=\int_0^1\f{s^{\f12}}{(1+(1-s)r_c+sr)^{\f12}}ds\lesssim \int_0^1\f{s^{\f12}}{(s+sr)^{\f12}}ds= \langle r\rangle ^{-\f12}$ and
									\begin{align*}
										&I(r)\geq \int_0^1\f{s^{\f12}}{(1+(1-s)C+sr)^{\f12}}ds\gtrsim
										\int_0^1\f{s^{\f12}}{(1+C+r)^{\f12}}ds\sim \langle r\rangle ^{-\f12},\;r_c\leq C\\
										&I(r)\geq \int_0^1\f{s^{\f12}}{(1+(1-s)Cr+sr)^{\f12}}ds\gtrsim
										\int_0^1\f{s^{\f12}}{(1+(C+1)r)^{\f12}}ds\sim \langle r\rangle ^{-\f12},\;1\lesssim r_c\leq C r.
									\end{align*}
									
									Now we prove the derivative bounds in \eqref{estx4-6-G}. For $i=0$, $l=1$,
									thanks to $|\pa_cQ(s,c)|\lesssim \langle r_c\rangle ^{6}\langle r\rangle ^{-3}$ \eqref{est-der-Q} and $|x|\sim \xi |r-r_c|^{\f32}\langle r\rangle ^{-\f12}$,
									it suffices to show that for $r_c\lesssim 1$, $r>0$ and $r_c\gtrsim 1$, $r\gtrsim r_c$,
									\begin{align*}
										f(r)=\langle r_c\rangle ^{2}\left|\int_{r_c}^r|r_c-s|^{-\f12}\langle s\rangle ^{-\f52}ds\right|\lesssim  |r-r_c|^{\f12}\langle r\rangle
										^{-\f12}.\end{align*}
									Indeed, if $r\sim r_c$, then $f(r)\sim \langle r_c\rangle ^{-\f12}\left|\int_{r_c}^r|r_c-s|^{-\f12}ds\right|\sim |r-r_c|^{\f12}\langle r\rangle
									^{-\f12} $; if $r\geq 2r_c$, then
									$f(r)= \langle r_c\rangle ^{2}\int_{r_c}^{2r_c}|r_c-s|^{-\f12}\langle s\rangle ^{-\f52}ds+\langle r_c\rangle ^{2}\int_{2r_c}^{r}|r_c-s|^{-\f12}\langle s\rangle ^{-\f52}ds
									\lesssim 1+\langle r_c\rangle ^{\f32}
									\langle r\rangle ^{-\f32}\lesssim  1\lesssim |r-r_c|^{\f12}\langle r\rangle
									^{-\f12} $;
									if $r\leq r_c/2$ and $r_c\lesssim 1$, then
									$f(r)= \langle r_c\rangle ^{2}\int^{r_c}_{r_c/2}|r_c-s|^{-\f12}\langle s\rangle ^{-\f52}ds+\langle r_c\rangle ^{2}\int^{r_c/2}_{r}|r_c-s|^{-\f12}\langle s\rangle ^{-\f52}ds
									\lesssim \langle r_c\rangle ^{-\f12}r_c^{\f12}
									+\langle r_c\rangle ^{2}
									\langle r\rangle ^{-\f52}|r_c-r|^{\f12}\lesssim r_c^{\f12} \leq
									|r-r_c|^{\f12}\sim |r-r_c|^{\f12}\langle r\rangle
									^{-\f12} $.
									For $i\geq 1,\;l\in \{0,1\}$, using \eqref{Faadi:x},  $Q(r,c)\sim \f{r-r_c}{\langle r\rangle }$ and $|r-r_c|\lesssim \langle r\rangle $($r_c\lesssim 1$), we have
									\begin{align*}
										\left||r-r_c|^{i+l}\pa_c^l\pa_r^ix\right|&\lesssim \xi\sum_{0\leq k\leq i-1}|r-r_c|^{i+l+\f12-k-l}\langle r\rangle^{-\f12+k+l}\left(\langle r\rangle^{-3l-3k-(i-1)}+\mathrm{sgn}(k)
										\langle r\rangle^{-3k-(i-1)}\right)\\
										&\lesssim \xi |r-r_c|^{\f32}\langle r\rangle^{-\f12}\sum_{0\leq k\leq i-1}\langle r\rangle^{-2l-3k}+\mathrm{sgn}(k)
										\langle r\rangle^{l-3k}\lesssim \xi |r-r_c|^{\f32}\langle r\rangle^{-\f12}
										\sim |x|.
									\end{align*}
									Finally, we prove \eqref{estx4-6-G} for $r_c\gtrsim 1 $, $r\leq r_c/2$.  By the definition,
									$|Q(s,c)|\sim r_c^2\f{r_c-s}{\langle s\rangle ^3}$ and $r_c-s\sim r_c(s\leq r_c/2)$, it follows that
									\begin{align} \label{equiv:x}
|x|=\xi r_c\int_{r}^{r_c/2}|r_c-s|^{\f12}\langle s\rangle ^{-\f32}ds
									\gtrsim  \xi r_c^{\f32}\int_{r}^{r_c/2}\langle s\rangle ^{-\f32}ds
									\sim \xi r_c^{\f32}\langle r\rangle ^{-\f12}.\end{align} 									Now we prove the derivative bounds $|\langle r_c\rangle^{-3l}\langle r\rangle ^{i}\pa_c^l\pa_r^ix|\lesssim |x|$ in \eqref{estx4-6-G}. For $i=0$, $l=1$, thanks to
									$|Q(r,c)|\sim r_c^2\f{|r-r_c|}{\langle r\rangle ^3}(r\lesssim r_c)$, $|\pa_cQ(s,c)|\lesssim \langle r_c\rangle ^{6}\langle r\rangle ^{-3}$ and $|x|\gtrsim r_c^{\f32}\langle r\rangle ^{-\f12}$, it suffices to prove that for $r_c\gtrsim 1$ and $r\lesssim r_c$,
									\begin{align*}
										f(r)=\langle r_c\rangle ^{\f12}
										\int_{r}^{r_c}|r_c-s|^{-\f12}\langle s\rangle ^{-\f32}ds\lesssim \langle r\rangle
										^{-\f12}.\end{align*}
										Indeed, $f(r)= \langle r_c\rangle ^{\f12}\int^{r_c}_{r_c/2}|r_c-s|^{-\f12}\langle s\rangle ^{-\f32}ds+\langle r_c\rangle ^{\f12}\int^{r_c/2}_{r}|r_c-s|^{-\f12}\langle s\rangle ^{-\f32}ds
										\sim \langle r_c\rangle ^{\f12}r_c^{\f12}\langle r_c\rangle ^{-\f32}+
										\langle r_c\rangle ^{\f12}r_c^{-\f12}\langle r\rangle ^{-\f12}
										\lesssim  \langle r\rangle ^{-\f12} $.
										For $i\geq 1,\;l\in \{0,1\}$, using \eqref{Faadi:x},  $|Q(r,c)|\sim \f{r_c^3}{\langle r\rangle^3 }(r\leq r_c/2)$, we have
										for $r\leq r_c/2$
										\begin{align*}
											&\notag\left|\langle r_c\rangle ^{-3l}\langle r\rangle^i\pa_c^l\pa_r^ix\right|\\
											&\lesssim \xi \langle r_c\rangle ^{-3l}\langle r\rangle^i
											\sum_{\substack{0\leq k\leq i-1}}(\f{r_c^3}{\langle r\rangle ^{3}})^{\f12-k-l}\cdot\left(\langle r\rangle ^{-3l-3k-(i-1)}+\mathrm{sgn}(k)
											(\f{r_c^3}{\langle r\rangle^3})^l\langle r\rangle ^{-3k-(i-1)}\right)\\
											&= \xi r_c^{\f32}\langle r\rangle^{-\f12}r_c^{-3l}
											\sum_{\substack{0\leq k\leq i-1}}r_c^{-3k}\left(r_c^{-3l}+\mathrm{sgn}(k)\right)\lesssim \xi r_c^{\f32}\langle r\rangle^{-\f12}r_c^{-3l}\lesssim |x|.
										\end{align*}
										\eqref{estx-G} is a summary of the above estimate, \eqref{equiv:x1},\eqref{equiv:x} and 
										\eqref{estx4-6-G}.
										Indeed, we notice that, for $r\leq r_c$, it holds $|r_c-r|\gtrsim
										\f{|r_c-r|\langle r\rangle}{\langle r_c\rangle}$; and for $r\leq r_c/2$, it holds $|r-r_c|\sim r_c$. 							
										\end{proof}
							
							\begin{lemma}\label{lem:behave-q}
								Let $q=\f{Q}{\tau}$, $\tau$ be defined in \eqref{def:tau}.  Then $q>0$, $\pa_r\tau=q^{\f12}$ and $\pa_{r}=q^{\f12}\pa_{\tau}$.
								 Let $c\in(V(0),1)$, i.e., $r_c\in(0,+\infty)$.  For $r\lesssim r_c$, it holds that
									\begin{align}\label{bd:q}
										\f{\langle r_c\rangle^{\f43}}{\langle r\rangle^{2}}\lesssim q(r,c)\lesssim \f{\langle r_c\rangle^{\f73}}{\langle r\rangle^{3}},
									\end{align}
									and
									\begin{align}\label{bd:q-Derivative} \left|\langle r_c\rangle ^{-3l}\langle r\rangle ^i\pa_c^l\pa_r^iq(r,c)
										\right|\lesssim q(r,c).
									\end{align}
								
							\end{lemma}
							
							\begin{proof}
							We use the form
							\begin{align*}q(r,c)=g_1(r,c)\cdot\left(\f32\Phi_1(r,c)\right)^{-\f23}.
							\end{align*} In fact, by the definition, we write $Q=(r-r_c)g_1(r,c)$ and $\tau=(r-r_c)\big(\f32\Phi_1(r,c)\big)^{\f23}$, where
							$g_1$ is defined in \eqref{def:g1}  and
							\begin{align}\label{def:Phi1}
								\Phi_1(r,c)=\int_0^1s^{\f12}g_1^{\f12}
								((1-s)r_c+sr,c)ds.
							\end{align}
							We claim that for $r\lesssim r_c$,
							\begin{align*}
								g_1(r,c)\sim \f{\langle r_c\rangle^{2}}{\langle r\rangle^{3}},\;\;
								\langle r_c\rangle^{-\f12}\lesssim \Phi_1(r,c)\lesssim \f{\langle r_c\rangle}{\langle r\rangle^{\f32}},
							\end{align*}
							and
							\begin{align*}
								\left|\langle r_c\rangle ^{-3l}\langle r\rangle ^i\pa_c^l\pa_r^ig_1(r,c)\right|\lesssim g_1(r,c),\quad
								\left|\langle r_c\rangle ^{-3l}\langle r\rangle ^i\pa_c^l\pa_r^i\Phi_1(r,c)\right|\lesssim \Phi_1(r,c).
							\end{align*}
							The first two bounds give \eqref{bd:q}, and \eqref{bd:q-Derivative} can be obtained via the latter two bounds.
							More precisely, we use Fa$\grave{\text{a}}$di Bruno's formula and Leibniz's rule to obtain
							\begin{align*}
								|\langle r_c\rangle ^{-3l}\langle r\rangle ^i\pa_c^l\pa_r^iq(r,c)|\quad\quad\quad\quad\quad\quad\quad\quad\quad\quad&\\
								\lesssim \langle r_c\rangle ^{-3l}\langle r\rangle ^i\sum_{\substack{0\leq s\leq i,\;k_1,...,k_i\geq 0\\k_1+...+sk_s=s,\;k=k_1+...+k_s}}\Big|\pa_c^l\Big(
								(\pa_r^{i-s}g_1)&\Phi_1^{\f23-k}
								(\pa_r\Phi_1)^{k_1}...(\pa_r^s\Phi_1)^{k_s}\Big)\Big|\\
								\lesssim (g_1\Phi_1^{-\f23})\langle r_c\rangle ^{-3l}\langle r\rangle ^i\sum_{\substack{0\leq s\leq i,\;0\leq t\leq l\\
										0\leq k_1,...,k_i\leq k,\;0\leq l',...,l_s\leq t\\k_1+...+sk_s=s,\;k_1+...+k_s=k\\
										l_0+l'+...+l_s=t}}\Big|&
								(\pa_c^{l-t}\pa_r^{i-s}g_1/g_1)\Phi_1^{-k-l_0}
								(\pa_r\Phi_1)^{k_1-l'}...(\pa_r^s\Phi_1)^{k_s-l_s}\\
								&\cdot(\pa_c\Phi_1)^{l_0}(\pa_c\pa_r\Phi_1)^{l'}...(\pa_c\pa_r^s\Phi_1)
								^{l_s}\Big|,\end{align*}
							which can be rewritten as
							\begin{align*}
								(g_1\Phi_1^{-\f23})(\langle r_c\rangle ^{-3})^l\langle r\rangle ^i\sum_{\substack{0\leq s\leq i,\;0\leq t\leq l\\
										0\leq k_1,...,k_i\leq k,\;0\leq l',...,l_s\leq t\\k_1+...+sk_s=s,\;k_1+...+k_s=k\\
										l_0+l'+...+l_s=t}}\Big|&
								(\pa_c^{l-t}\pa_r^{i-s}g_1/g_1)
								(\pa_r\Phi_1/\Phi)^{k_1-l'}...(\pa_r^s\Phi_1/\Phi_1)^{k_s-l_s}\\
								&\cdot(\pa_c\Phi_1/\Phi_1)^{l_0}...(\pa_c\pa_r^s\Phi_1/\Phi_1)
								^{l_s}\Big|\lesssim q(r,c).
							\end{align*}
							
							Now we prove the  bounds in the claim. The behavior of $g_1$ is obtained in the proof of Lemma \ref{lem:behave-Q}, where for $r\lesssim r_c$,
							 it simplifies to $\f{\langle r+r_c\rangle^{2}}{\langle r\rangle^{3}}\sim \f{\langle r_c\rangle^{2}}{\langle r\rangle^{3}}$.
							Using this and the definition, we obtain the upper/lower bound of  $\Phi_1$.
							Indeed, it suffices to check that the followings hold uniformly in $s\in(0,1)$:
							\begin{align*}
								\f{( r_c+1)}{( r_c+1)^{\f32}}\lesssim \f{(1-s)r_c +s(r+r_c)+1
								}{\left((1-s)r_c +sr+1
									\right)^{\f32}}\lesssim \f{( r_c+1)}{(r+1)^{\f32}},\;\;r\lesssim r_c.
							\end{align*}
							To prove the derivative bounds of $g_1$, we write by the definition as
							\begin{align*}
								\pa_c^l\pa_r^ig_1(r,c):=\int_{0}^1\f{s^iV^{(i+1)}(r_c+s(r-r_c))}
								{(1-c)^{1+l}}ds
								+\int_{0}^1\f{s^i(1-s)^lV^{(i+1+l)}(r_c+s(r-r_c))}{(V'(r_c))^{l}(1-c)}ds,
							\end{align*}
							which along with $|V^{(i+1)}(r)|+|r^lV^{(i+1+l)}(r)|\lesssim \langle r\rangle^{-i-4}$, $(1-c)^{-1}\sim \langle r_c\rangle^{3}$, $
							(V'(r_c))^{-1}\sim \langle r_c\rangle^{4}$, gives
							\begin{align*}
								|\pa_c^l\pa_r^ig_1(r,c)|&\lesssim \langle r_c\rangle^{3+3l}\int_{0}^1\langle r_c+s(r-r_c)\rangle^{-4-i}ds\sim \langle r_c\rangle^{3+3l}\left|\f{\langle r\rangle^{-3-i}-\langle r_c\rangle^{-3-i}}{r-r_c}\right|
								\\
								&\sim\f{ \langle r_c\rangle^{3l-i}\langle r+r_c\rangle^{i+2}}{
									\langle r\rangle^{i+3}}
								\sim \langle r_c\rangle^{3l}\left(
								\f{ 1}{
									\langle r\rangle}\right)^ig_1(r,c).
							\end{align*}
							To prove the derivative bounds of $\Phi_1$ for $r\lesssim r_c$, we apply Fa$\grave{\text{a}}$di Bruno's formula, Leibniz's rule and
							$\left|\langle r_c\rangle ^{-3l'}\langle r\rangle ^{k'}\left(\pa_c^{l'}\pa_r^{k'}g_1\right)((1-s)r_c+sr,c)\right|\lesssim g_1((1-s)r_c+sr,c)$, to obtain
							\begin{align*}
								&\left|\langle r_c\rangle ^{-3l}\langle r\rangle ^i\pa_c^l\pa_r^i\Phi_1(r,c)\right|\\
								&\lesssim \langle r_c\rangle ^{-3l}\langle r\rangle ^i\sum_{\substack{k_1,...,k_i\geq 0\\k_1+...+ik_i=i\\k=k_1+...+k_i}}\left|\int_0^1s^{\f12+i}
								(\f{1-s}{V'(r_c)})^{l}\pa_r^l\left(g_1^{\f12-k}
								(\pa_rg_1)^{k_1}...(\pa_r^ig_1)^{k_i}\right)
								((1-s)r_c+sr,c)ds\right|\\
								&\quad+\langle r_c\rangle ^{-3l}\langle r\rangle ^i\sum_{\substack{k_1,...,k_i\geq 0\\k_1+...+ik_i=i\\k=k_1+...+k_i}}\left|\int_0^1s^{\f12+i}
								\pa_c^l\left(g_1^{\f12-k}
								(\pa_rg_1)^{k_1}...(\pa_r^ig_1)^{k_i}\right)
								((1-s)r_c+sr,c)ds\right|\\
								&\lesssim \left|\int_0^1s^{\f12}\left(\f{\left((1-s)\langle r_c\rangle\right)^{l}\left(s\langle r\rangle \right)^i}{\langle (1-s)r_c+sr\rangle^{l+i}}+\f{\left(s\langle r\rangle\right) ^i}{\langle (1-s)r_c+sr\rangle^{i}}\right)|g_1|^{\f12}
								((1-s)r_c+sr,c)ds\right|\\
								&\lesssim  \left|\int_0^1|g_1|^{\f12}
								((1-s)r_c+sr,c)ds\right|\sim\int_0^1\f{\langle r_c\rangle}{
									\langle(1-s)r_c+sr\rangle^{\f32}}ds\\
								&\leq C  \int_0^1\f{\langle r_c\rangle}{
									\langle(1-s)C^{-1}r+sr\rangle^{\f32}}ds\sim
								\f{\langle r_c\rangle}{
									\langle r\rangle^{\f32}}\sim g_1^{\f12}(r,c).
							\end{align*}
							In the third inequality, we used $
							\f{(1-s)\langle r_c\rangle}{\langle (1-s)r_c+sr\rangle}
							\leq 1$ and $\f{s\langle r\rangle }{\langle (1-s)r_c+sr\rangle}\leq 1$ for $s\in(0,1)$.							
							\end{proof}
							
							\begin{lemma}\label{lem:behave-W}
								Let $c\in(0,1)$, $s>0$ and define
								\begin{align}\label{W:fm0}
									W(s,c):=q^{-\f14}\pa_{t}^2(q^{\f14})+\f{3}{4s^2q},
								\end{align}
								where $t=t(s,c)$ is defined in \eqref{def:tau}, and $q(s,c)=\f{Q(s,c)}{t(s,c)}$  .
								Then $W$ has two alternative forms as follows
								\begin{align}\label{W:fm1}
									W(s,c)&=\f{5}{16 t^2(s,c)}+t(s,c)\left(\f{-\pa_s^2Q(s,c)}{4Q^2(s,c)}+\f{5\pa_sQ^2(s,c)}{16Q(s,c)^3}
									+\f{3}{4s^2Q(s,c)}\right),\\
									W(s,c)
									&=\f{\pa_s^2q(s,c)}{4q^2(s,c)}-\f{5\pa_sq^2(s,c)}{16q^3(s,c)}
									+\f{3}{4s^2q(s,c)}.\label{W:fm2}
								\end{align}
								Moreover, the following estimates hold
								\begin{itemize}
									\item  For $c\in(0,V(0))$,
									\begin{align}\label{est:W<V(0)}
										\left|(V(s)-c)s^i\pa_c^l\pa_s^iW(s,c)\right|\lesssim t^{-2},\;\quad s>0.
									\end{align}
									\item  For $c\in(V(0),1)$ and $r_c=V^{-1}(c)$,
									\begin{align}\label{est:W>rc}
										&\left|\left(\f{s}{\langle r_c\rangle ^{4}}\right)^ls^{i}\pa_c^l\pa_s^iW(s,c)\right|\lesssim t^{-2},\quad s\geq 2r_c,\\
									 \label{est:W<rc}
										&\left|\left(\f{1}{\langle r_c\rangle ^{3}}\right)^l s^{i}\pa_c^l\pa_s^iW(s,c)\right|\lesssim\f{1}{s^2q},\quad s\lesssim r_c.
									\end{align}
								\end{itemize}
							\end{lemma}
							
							\begin{proof}
								\eqref{W:fm2} follows by expanding \eqref{W:fm0} with $\pa_{t}=q^{-\f12}\pa_{s}$; \eqref{W:fm1} follows by expanding \eqref{W:fm2} with $q=\f{Q}{t}$ and $\pa_st=q^{\f12}$.
								
								Concerning \eqref{est:W<V(0)}-\eqref{est:W<rc}, roughly speaking,  for $s\geq 2r_c$, we use \eqref{W:fm1}; for $s\lesssim r_c$, we use \eqref{W:fm2}. Notice for $c\in(0,V(0))$, $s\geq 2r_c$, it holds $s-r_c\sim s$. We conclude by
								Lemma \ref{lem:behave-Q}  that
								\begin{align*}
									&|t|^{\f32}(s,c)\lesssim s|Q|^{\f12}(s,c),\;\text{for}\; c\in(0,V(0)],s>0,\;\text{and}\;c\in(V(0),1),s\geq r_c.\\
									&\left|\langle r_c\rangle ^{-4l}s^{i+l}\pa_c^l\pa_s^iQ(s,c)\right|
									\lesssim |Q(s,c)|,\;\text{for}\; c\in(V(0),1),s\geq 2r_c.\\
									&\left|(V(s)-c)^ls^i\pa_c^l\pa_s^iQ(s,c)\right|
									\lesssim |Q(s,c)|,\;\text{for}\; c\in(0,V(0)],s>0.
								\end{align*}
								We also get by Lemma \ref{lem:behave-x} that
								\begin{align}\label{behave-tau-in-W}
									\begin{split}
										&\left|\langle r_c\rangle ^{-4l}s^{i+l}\pa_c^l\pa_s^it(s,c)\right|
										\lesssim t(s,c),\;\text{for}\; c\in(V(0),1),s\geq 2r_c.\\
										&\left|(V(s)-c)^ls^i\pa_c^l\pa_s^it(s,c)\right|
										\lesssim t(s,c),\;\text{for}\; c\in(0,V(0)],s>0.
									\end{split}
								\end{align}
								To prove \eqref{est:W<V(0)} and \eqref{est:W>rc}, we rewrite \eqref{W:fm1} as
								\begin{align*}
									W(s,c)&=\f{5}{16 t^2(s,c)}+\f{3t(s,c)}{4s^2Q(s,c)}
									\left(\f{-s^2\pa_s^2Q(s,c)}{3Q(s,c)}+\f{5s^2\pa_sQ^2(s,c)}{12Q(s,c)^2}
									+1\right).
								\end{align*}
								Using the above estimates, we obtain $|W|\lesssim t^{-2}$ for $c\in(0,V(0)]$, $s>0$ and $c\in(V(0),1)$, $s\geq 2r_c$. The derivative bounds for $i\in\mathbb{N}$, $l\in\{0,1\}$ follow by taking $\pa_c^l\pa_s^i$ directly in the above formula and using again the estimates of $Q$ and $t$.
								To prove \eqref{est:W<rc}, we conclude by Lemma \eqref{lem:behave-q} that for $s\lesssim r_c$,
								\begin{align*}
									\left|\langle r_c\rangle ^{-3l}\langle s\rangle ^i\pa_c^l\pa_s^iq(s,c)\right|\lesssim q(s,c).
								\end{align*}
								We also need  the formulation
								\begin{align*}W(s,c)
									=\f{3}{4s^2q(s,c)}\left(\f{s^2\pa_s^2q(s,c)}{3q(s,c)}
									-\f{5s^2\pa_sq^2(s,c)}{12q^2(s,c)}
									+1\right).
								\end{align*}
								Utilizing the bounds of $q$, we obtain $|W|\lesssim \f{1}{s^2q}$. Using again the bounds and Leibniz's rule, the derivatives of $W$ admit the same bounds.
							\end{proof}
							
							\begin{lemma}\label{lema:W-trans}
								Let  $c\in(0,1)$, $s>0$, $t=t(s,c)$ be defined in \eqref{def:tau}, $W(s,c)$ be defined in \eqref{W:fm0}, $y(s,\xi,c)=\xi\mathrm{sgn}(s-r_c) t^{\f32}(s,c)$ be defined in \eqref{def:x}, where $s=s(y)$  and $s=s(t)$ are monotonic inverse functions of $y=y(s)$ and $t=t(s)$, respectively, due to $\pa_sy,\pa_st> 0$. We denote
								\begin{align*}\widetilde{W}(y,\xi,c)=W(s(y,\xi,c),c),\;\;\widetilde{\widetilde{W}}(t,c)
									=W(s(t,c),c).
									\end{align*}
								Then for $j,i\in\mathbb{N}$, $l\in\{0,1\}$, it holds that
								\begin{align}
									\begin{split}\label{est:Dy-tW}
										&\left|y^i\left(V(s)-c)\pa_c\right)^l\pa_y^i(\xi\pa_{\xi})^j
										\widetilde{W}(y,\xi,c)\right|
										\lesssim \xi^{\f43}y^{-\f43}\quad \text{for}\;c\in(0,V(0)),\\
										&\left|\left(\f{r_c}{\langle r_c\rangle ^4}\right)^ly^i\pa_c^l\pa_y^i(\xi\pa_{\xi})^j\widetilde{W}(y,\xi,c)\right|\\
										&\quad\lesssim \xi^{\f43}y^{-\f43}\mathbf{1}_{y\geq y(2r_c)}(y)
										+\f{\langle r_c\rangle^{\f23}}{r_c^2}\mathbf{1}_{y(r_c/4)\leq y\leq
											y(2r_c)}(y)\quad \text{for}\;c\in(V(0),1),
									\end{split}
								\end{align}
								and for $c\in(V(0),1)$ and $\xi r_c^{\f32}/\langle r_c\rangle^{\f12}\gtrsim M\gg 1$,							
								\begin{align} \label{est:Dy-ttW}
									&\left|\left(\f{r_c}{\langle r_c\rangle ^4}\pa_c\right)^l(\xi^{-\f23}\pa_{t})^i\widetilde{\widetilde{W}}
									(t,c)\right|
									\lesssim r_c^{-2}\langle r_c\rangle^{\f23}\mathbf{1}_{|t|\lesssim \xi^{-\f23}}(t).
								\end{align}
							\end{lemma}
							
							\begin{proof}
								Observe that $t^{-2}(s)=\xi^{\f43}y^{-\f43}(s)$ for $s\geq r_c$ and $\f{1}{s^2q(s,c)}\sim \f{\langle r_c\rangle^{\f23}}{r_c^2}$ for $s\sim r_c$, due to \eqref{bd:q}. We denote
								\begin{align*}
									\lambda^c(s,c)=
									\left\{
									\begin{aligned}
										&V(s)-c,\;\text{if}\;c\in(0,V(0)]\;\text{and}\;s>0,\\
										&\f{r_c}{\langle r_c\rangle^{4}},\;\text{if}\;c\in(V(0),1)\;\text{and}\;s>0,
									\end{aligned}
									\right.
								\end{align*}
								and
								\begin{align*}
									\mathcal{W}(s,c):=
									\left\{
									\begin{aligned}
										&t^{-2}(s,c),\;\text{if}\;c\in(0,V(0)]\;\text{and}\;s>0,\\
										&t^{-2}(s,c)\mathbf{1}_{s\geq 2r_c}+\f{1}{s^2q(s,c)}\mathbf{1}_{r_c/4\leq s\leq 2r_c}
										,\;\text{if}\;c\in(V(0),1).
									\end{aligned}
									\right.
								\end{align*}
								We also notice that $y^i\pa_y^i=\sum_{1\leq k\leq i}(-1)^kC_n^k(y\pa_y)^k$. For $\widetilde{f}(y,\xi,c)=f(s(y,\xi,c),c)$, it holds that  $\xi\pa_{\xi}\widetilde{f}(y,\xi,c)
								=\f{-\xi\pa_{\xi}y}{\pa_sy}\pa_sf(s,c)
								=\f{-y}{\pa_sy}\pa_sf(s,c)$, $y\pa_{y}\widetilde{f}(y,\xi,c)
								=\f{y}{\pa_sy}\pa_sf(s,c)$ and $\lambda^c\pa_{c}\widetilde{f}(y,\xi,c)
								=-\f{\lambda^c\pa_cy}{\pa_sy}\pa_sf(s,c)+(\lambda^c\pa_c)f(s,c)$.
								Then \eqref{est:Dy-tW} can be deduced by
								\begin{align}\label{est:Dy-tW-2}
									&\notag \left|\left(\lambda^c(s,c)\pa_c\right)^l(y\pa_y)^i
									(\xi\pa_{\xi})^j\widetilde{W}(y,\xi,c)\right|\\
									\notag &\leq \left|\left(\f{\lambda^c\pa_cy}{\pa_sy}\pa_s\right)^l
									\circ\left(\f{y}{\pa_sy}\pa_s\right)^{i+j}
									W(s,c)\right|
									+\left|(\lambda^c\pa_c)^l\circ\left(\f{y}{\pa_sy}\pa_s\right)^{i+j}
									W(s,c)\right|\\
									&:=A+B\lesssim \mathcal{W}(s,c),\;\;s\geq r_c/4.
								\end{align}
								Before estimating $A$ and $B$, we summarize some bounds as follows. By the definition, $y(s)=\f23\xi\int_s^{r_c}|Q(w,c)|^{\f12}dw$, $\pa_sy(s)=\f23\xi|Q(s,c)|^{\f12}$. By \eqref{est-tleqQ}, we have
								\begin{align}\label{est:tWladar}
									&\left|\f{y}{\pa_sy}\right|
									=\left|\f{\int_s^{r_c}|Q(w,c)|^{\f12}dw}{|Q(s,c)|^{\f12}}\right|\sim \lambda^r(s,c)\leq s,\;\\
									\notag&\lambda^r(s,c):=
									\left\{
									\begin{aligned}
										&s,\quad\text{if}\;c\in(0,V(0)],\,\,\text{and}\;s>0;\\
										&|s-r_c|,\quad\text{if}\;c\in(V(0),1),\,\,\text{and}\;s\geq r_c/4.
									\end{aligned}
									\right.
								\end{align}
								By Lemma \ref{lem:behave-Q}, we have for $n\geq 1$, $l\in\{0,1\}$,
								\begin{align*}
									\left|\f{\lambda^r(s,c)\langle s\rangle^{n-1}(\lambda^c\pa_c)^l\pa_s^nQ(s,c)}{Q(s,c)}\right|\lesssim 1,\;s\geq r_c/4,
								\end{align*}
								which along with Leibniz's rule and Fa$\grave{\text{a}}$di Bruno's formula, implies that for $n\geq 1$, $l\in\{0,1\}$,
								\begin{align}\label{est:lamcpasy/spasy}
									\left|\left(\f{y}{\pa_sy}\right)^{n-1} \f{(\lambda^c\pa_c)^l\pa_s^ny}{\pa_sy}\right|\lesssim 1,\;\left|\left(\f{y}{\pa_sy}\right)^{n-1} (\lambda^c\pa_c)^l\pa_s^n\left( \f{y}{\pa_sy}\right)\right|\lesssim 1,\;\;s\geq r_c/4.
								\end{align}
								Moreover, we can check by  Lemma \ref{lem:behave-x} that
								\begin{align}\label{est:lamcpacy/spasy}
									\left|\left(\lambda^c\pa_c\right)^ly\right|\lesssim |s\pa_sy|.
								\end{align}
								Indeed, $l=0$ follows from \eqref{est:tWladar}; for $l=1$, it follows that for $s\gtrsim r_c$,
								\begin{align*}
									\f{r_c}{\langle r_c\rangle^4}
									\left|\int_{r_c}^s\langle r_c\rangle^6\langle w\rangle^{-\f52}|w-r_c|^{-\f12}dw\right|\lesssim
									\f{r_c}{\langle r_c\rangle^{\f12}}|s-r_c|^{\f12}
									\lesssim \f{s}{\langle s\rangle^{\f12}}|s-r_c|^{\f12}.
								\end{align*}
								With these bounds at hand, we are in a position to estimate $A, B$.
								Thanks to \eqref{est:lamcpacy/spasy}($l=1$),  Leibniz's rule, Fa$\grave{\text{a}}$di Bruno's formula, \eqref{est:lamcpasy/spasy}($l=0$), $\left|\f{y}{\pa_sy}\right|\lesssim s$ and Lemma \ref{lem:behave-W}, it follows that for $s\geq r_c/4$,
								\begin{align*}
									A&\lesssim \left|(s\pa_s)^l\left(\f{y}{\pa_sy}\pa_s\right)^{j+i}
									W(s,c)\right|\\
									&\lesssim \sum_{\substack{0\leq m\leq i+j+l-1,\;k_1,...,k_{m}\geq 0\\k_1+...+mk_{m}=m,\;
											k_1+...+k_{m}=k}}\left|\left(\f{y}{\pa_sy}\right)^{i+j+l-k}\cdot
									\left(\pa_s\left(\f{y}{\pa_sy}\right)\right)^{k_1}...
									\left(\pa_s^m\left(\f{y}{\pa_sy}\right)\right)^{k_m}\pa_s^{i+j+l-m}W(s,c)\right|\notag\\
									&\lesssim \sum_{\substack{0\leq m\leq i+j+l-1,\;k_1,...,k_{m}\geq 0\\k_1+...+mk_{m}=m,\;
											k_1+...+k_{m}=k}}\left|\left(\f{y}{\pa_sy}\right)^{i+j+l-k-(k_2+2k_3+...+(m-1)k_m)}\cdot
									\pa_s^{i+j+l-m}W(s,c)\right|\notag\\
									&= \sum_{\substack{0\leq m\leq i+j+l-1}}\left|\left(\f{y}{\pa_sy}\right)^{i+j+l-m}
									\pa_s^{i+j+l-m}W(s,c)\right|\\
									&\lesssim \sum_{\substack{0\leq m\leq i+j+l-1}}\left|s^{i+j+l-m}
									\pa_s^{i+j+l-m}W(s,c)\right|\lesssim \mathcal{W}(s,c).
								\end{align*}
								Similarly, it follows  from 
								Leibniz's rule,  Fa$\grave{\text{a}}$di Bruno's formula, \eqref{est:lamcpasy/spasy}($l=0,1$), $\left|\f{y}{\pa_sy}\right|\lesssim s$ and Lemma \ref{lem:behave-W}
								that for $s\geq r_c/4$,
								\begin{align*}
									B
									&\lesssim \sum_{\substack{0\leq m\leq i+j-1,\\k_0,...,k_{m}\geq 0\\k_1+...+mk_{m}=m,\\
											k_0+k_1...+k_{m}=k}}\left(\left|\left(\f{y}{\pa_sy}\right)^{i+j-k}\cdot(\lambda^c\pa_c)^l
									\left(
									\pa_s\left(\f{y}{\pa_sy}\right)^{k_0}
									\left(\pa_s\left(\f{y}{\pa_sy}\right)\right)^{k_1}...
									\left(\pa_s^m\left(\f{y}{\pa_sy}\right)\right)^{k_m}\right)\pa_s^{i+j-m}W(s,c)\right|\right.\notag\\
									&\left.\quad+ \sum_{\substack{0\leq m\leq i+j-1,\;k_1,...,k_{m}\geq 0\\k_1+...+mk_{m}=m,\;
											k_1+...+k_{m}=k}}\left|
									\left(\f{y}{\pa_sy}\right)^{i+j-k}\cdot \left(\pa_s\left(\f{y}{\pa_sy}\right)\right)^{k_1}...
									\left(\pa_s^m\left(\f{y}{\pa_sy}\right)\right)^{k_m}\left(\lambda^c\pa_c\right)^l
									\pa_s^{i+j-m}W(s,c)\right|\right)\notag\\
									&\lesssim \sum_{\substack{0\leq m\leq i+j-1}}\left(\left|\left(\f{y}{\pa_sy}\right)^{i+j-m}
									\pa_s^{i+j-m}W(s,c)\right|+\left|\left(\f{y}{\pa_sy}\right)^{i+j-m}
									(\lambda^c\pa_c)^l\pa_s^{i+j-m}W(s,c)\right|\right)\\
									&\lesssim \sum_{\substack{0\leq m\leq i+j-1}}\left|s^{i+j-m}
									\pa_s^{i+j-m}W(s,c)\right|+\left|s^{i+j-m}
									(\lambda^c\pa_c)^l\pa_s^{i+j-m}W(s,c)\right|\lesssim \mathcal{W}(s,c).
								\end{align*}
								
								Finally, we prove \eqref{est:Dy-ttW}. It follows from  \eqref{estx4-6-G} ($|x|=\xi\tau^{\f32}\sim \xi |r-r_c|^{\f32}/\langle r\rangle ^{\f12}$) that
								\begin{align*} &|\tau|\lesssim \xi^{-\f23}\Leftrightarrow r_c-C\xi^{-\f23}\langle r_c\rangle^{\f13}\leq r\leq r_c+C\xi^{-\f23}\langle r_c\rangle^{\f13}. \end{align*}
								Using $\xi r_c^{\f32}/\langle r_c\rangle^{\f12}\gtrsim M\gg 1$, we know, $|\tau|\lesssim \xi^{-\f23}$ implies $s=s(t,c)\in[ r_c/2,2r_c]$.
								We  conclude by \eqref{bd:q} and \eqref{bd:q-Derivative} that
								\begin{align*}
									\left|\left(\f{1}{\langle r_c\rangle^3 }\right)^{l}\langle r_c\rangle ^i\pa_c^l\pa_r^i(q^{\f12}(s,c))\right|\lesssim q^{\f12}(s,c)\sim \langle r_c\rangle^{-\f13}, \quad \text{for}\;\; s\in[\f{r_c}{2},2r_c],
									\end{align*}
								which along with $\pa_st=q^{\f12}$, $\pa_ct=q^{\f12}\f{\int_{r_c}^s\pa_c(Q^{\f12})ds'}{Q^{\f12}}$, \eqref{estQ-G} and \eqref{behave:Q<rc0} gives \begin{align*}
									&\left|\left(\f{1}{\langle r_c\rangle^{3}}\right)^l
									\langle r_c\rangle^i\pa_c^l\pa_s^it(s,c)\right|\lesssim \langle r_c\rangle^{\f23}\quad \text{for}\;\; s\in[\f{r_c}{2},2r_c].
								\end{align*}
								On the other hand, \eqref{est:W<rc} takes the form of
								\begin{align*}
									\left|\left(\f{1}{\langle r_c\rangle ^{3}}\right)^l r_c^{i}\pa_c^l\pa_s^iW(s,c)\right|\lesssim\f{\langle r_c\rangle^{\f23}}{r_c^2}\quad \text{for}\;\; s\in[\f{r_c}{2},2r_c].
									\end{align*}
								Notice that  for  $\tilde{\tilde{f}}(t,c)=f(s,c)$,  it holds that $\pa_t\tilde{\tilde{f}}(t)=q^{-\f12}\pa_sf(s)$ and $$\pa_c\tilde{\tilde{f}}(t,c)=-\pa_ct\cdot q^{-\f12}\pa_sf(s(t,c),c)+\pa_cf(s,c).$$  Then \eqref{est:Dy-ttW} follows by the estimates above and $\xi r_c^{\f32}/\langle r_c\rangle^{\f12}\gtrsim 1$.
								
								This finishes the proof of the lemma.
							\end{proof}
							
			\section{Airy functions}
							
					We recall some basic facts on Airy functions $\mathrm{Ai}(z)$. For $z\geq 0$, $\mathrm{Ai}(z)$ is a positive smooth function with the following asymptotic behavior:
							\begin{align}\label{behave:Ai}
								\mathrm{Ai}(z)=C\langle z^{\f32}\rangle^{-\f16}e^{-\f{2}{3}z^{\f32}}\left(1+O\left(\langle z^{\f32}\rangle^{-1}\right)\right)\;\;z\geq 0,
							\end{align}
							where $\langle z^{\f32}\rangle=z^{\f32}+1$ and $f(z)=O(f_0(z))$ represents $\left|\langle z^{\f32}\rangle^{\f23k}\pa_z^kf(z)\right|\lesssim f_0(z)$($k\in\mathbb{N}$).
							This along with $\mathrm{Ai}(z)>0$ implies
							\begin{align}\label{behave:Ai^-1}
								\mathrm{Ai}^{-1}(z)=C\langle z^{\f32}\rangle^{\f16}e^{\f{2}{3}z^{\f32}}\left(1+O(\langle z^{\f32}\rangle^{-1})\right)\;\;z\geq 0.
							\end{align}
							For $z\geq 0$, $ \mathrm{Oi}(-z):= \mathrm{Ai}(-z)-i \mathrm{Bi}(-z)$ is a complex-valued non-zero smooth function with the following asymptotic behavior:
							\begin{align}\label{behave:Oi}
								\begin{split}
									&\mathrm{Oi}(-z)=C\langle z^{\f32}\rangle^{-\f16}e^{\f{2i}{3}z^{\f32}}
									\left(1+O(\langle z^{\f32}\rangle^{-1})\right)\;\;z\geq 0,
								\end{split}
							\end{align}
							which along with $\mathrm{Oi}(-z)\neq 0(z\geq 0)$ implies
							\begin{align}\label{behave:Oi-1}
								\begin{split}
									&\mathrm{Oi}^{-1}(-z)=C\langle z^{\f32}\rangle^{\f16}e^{-\f{2i}{3}z^{\f32}}
									\left(1+O(\langle z^{\f32}\rangle^{-1})\right)\;\;z\geq 0.
								\end{split}
							\end{align}
							The above  $C$ represent different real/complex constants.
							\begin{definition}\label{def:K}
								Let $\nu\in\{\mathrm{Ai},\mathrm{Oi}\}$. We define the following kernels
								\begin{align}\label{def:Knu}
									K_{\nu}(\tau,t,\xi):=
									\left\{
									\begin{array}{l}
										\mathrm{Oi}^2(-\xi^{\f23}t)\int_{\tau}^tOi^{-2}(-\xi^{\f23}s)ds\quad 0\leq \tau\leq t,\; \nu=\mathrm{Oi},\\
										\mathrm{Ai}^2(-\xi^{\f23}t)\int_{\tau}^t\mathrm{Ai}^{-2}(-\xi^{\f23}s)ds\quad t\leq \tau\leq 0,\; \nu=\mathrm{Ai},
									\end{array}\right.
								\end{align}
								and
								\begin{align}\label{def:tKnu}
									\widetilde{K}_{\nu}(x,y):=
									\left\{
									\begin{array}{l}
										\mathrm{Oi}^2(-y^{\f23})\int_{x}^{y}
										\mathrm{Oi}^{-2}(-\tilde{y}^{\f23})\tilde{y}^{-\f13}
										d\tilde{y}\quad 0\leq x\leq y,\;\; \nu=\mathrm{Oi},\\
										\mathrm{Ai}^2(y^{\f23})\int_{y}^{x}\mathrm{Ai}^{-2}(\tilde{y}^{\f23})
										\tilde{y}^{-\f13}d\tilde{y}\quad  0\leq x\leq y,\;\; \nu=\mathrm{Ai}.
									\end{array}\right.
								\end{align}
							\end{definition}
							We notice that for $x=\xi \mathrm{sgn}(r-r_c)\tau^{\f32}$, $y=\xi \mathrm{sgn}(s-r_c)t^{\f32}$, it holds that
							\begin{align}\label{convert:Knv-tKnv}
								K_{\nu}(\tau,t,\xi)=\f23\xi^{-\f23}\widetilde{K}_{\nu}(x,y).
							\end{align}
							
							\begin{lemma}\label{lem:tK}
								The asymptotic behavior of  $K_{\nu}/\widetilde{K}_{\nu}$ is stated as follows
								\begin{align}\label{asy:tKOi}
									\notag&\widetilde{K}_{{Oi}}(x,y)=
									y^{-\f13}
									\int_{x}^{y}e^{\f{4i}{3}(y-w)} \left(1+O(y^{-1})\right)
									\left(1+O(w^{-1})\right)
									dw\\
									&=_{\tilde{w}=y-w}y^{-\f13}
									\int_{0}^{y-x}e^{\f{4i}{3}\tilde{w}} \left(1+O(y^{-1})\right)
									\left(1+O\left( (y-\tilde{w})^{-1}\right)\right)
									d\tilde{w}\;\;1\lesssim x\leq y, \end{align}
								\begin{align}\label{asy:tKAi}
									\notag\widetilde{K}_{{Ai}}(x,y)&=
									y^{-\f13}
									\int_{x}^{y}e^{-\f{4}{3}(y-w)} \left(1+O( y^{-1})\right)
									\left(1+O( w^{-1})\right)
									dw\\
									&=_{\tilde{w}=y-w} y^{-\f13}
									\int_{0}^{y-x}e^{-\f{4}{3}\tilde{w}}
									\left(1+O(y^{-1})\right)
									\left(1+O(( y-\tilde{w})^{-1})\right)
									d\tilde{w}\;\;1\lesssim x\leq y,
								\end{align}
								and
								\begin{align}\label{asy:tKOi-small}
									&K_{{Oi}}(\tau,t,\xi)= O_{\xi,\tau,t}^{\xi,\xi^{-\f23},\xi^{-\f23}}(\xi^{-\f23})\;\;-\xi^{-\f23}\lesssim \tau\leq t\lesssim \xi^{-\f23}, \end{align}
								where $f(\xi,\tau)=O_{\xi,\tau}^{\xi,\rho^{\tau}}(f_0(\xi,\tau))$ represents $\left|(\xi^{-\f23}\pa_{\tau})^i(\xi\pa_{\xi})^jf(\tau,\xi)\right|\lesssim f_0(\tau,\xi)$ and $f(z)=O_{z}(f_0(z))$ represents $\left|z^i\pa_z^{i}f(z)\right|\lesssim f_0(z)$, $i,j\in\mathbb{N}$.
								In particular, it holds that
								\begin{align}
									\left|K_{{Ai}}(\tau,t,\xi)=\f23\xi^{-\f23}\widetilde{K}_{\mathrm{Ai}}(x,y)\right|&\leq C\xi ^{-1}(-t)^{-\f12}\;\; t\leq \tau\leq 0,\label{est:KAi}\\
									\label{est:KOi}
									\left|K_{{Oi}}(\tau,t,\xi)=\f23\xi^{-\f23}\widetilde{K}_{\mathrm{Oi}}(x,y)\right|&\leq C\xi ^{-1}t^{-\f12}\quad0\leq \tau\leq t.
								\end{align}
							\end{lemma}
							
							\begin{proof}
\eqref{asy:tKOi} is the consequence of \eqref{behave:Oi} and \eqref{behave:Oi-1}.
\eqref{asy:tKAi} is the consequence of \eqref{behave:Ai} and \eqref{behave:Ai^-1}.
To prove \eqref{est:KAi} and \eqref{est:KOi}, we claim that for the integral emerged in \eqref{asy:tKOi}, it holds that for $y-z\gtrsim 1$, $z\geq 0$,
								\begin{align}\label{est:IOi}
									\begin{split}
										&\left|I_{\mathrm{Oi}}(y,z)=\int_{0}^{z}e^{\f{2i}{3}\tilde{y}}\left(1+O\left((y-\tilde{y})^{-1}\right)\right)
										d\tilde{y}\right|\lesssim 1,\\
										&\left|I_{\mathrm{Ai}}(y,z)=\int_{0}^{z}e^{-\f{2}{3}\tilde{y}}\left(1+O\left((y-\tilde{y})^{-1}\right)\right)
										d\tilde{y}\right|\lesssim 1.
									\end{split}
								\end{align}
								Indeed, the boundedness follows by integration by parts, here we only show the first one:
								\begin{align*}
									&I_{\mathrm{Oi}}(y,z)=\f{3}{2i}e^{\f{2i}{3}\tilde{y}}\left(1+O\left(\langle y-\tilde{y}\rangle^{-1}\right)\right)
									\Big|_{\tilde{y}=0}^z
									-\int_{0}^{z}e^{\f{2i}{3}\tilde{y}}O\left(\langle y-\tilde{y}\rangle^{-2}\right)
									d\tilde{y}=I_1+I_2.
								\end{align*}
								It follows directly that $|I_1|,|I_2|\lesssim 1$.
								Thus, together with \eqref{convert:Knv-tKnv}, \eqref{est:KAi} can be deduced by \eqref{asy:tKAi}, and \eqref{est:KOi} can be deduced by \eqref{asy:tKOi}.
								To prove \eqref{asy:tKOi-small}, we rewrite the kernel for $-\xi^{-\f23}\lesssim \tau\leq t\lesssim \xi^{-\f23}$ as
								\begin{align*}
									K_{{Oi}}(\tau,t,\xi)=
									(\tau-t)\mathrm{Oi}^2(-\xi^{\f23}t)\int_{0}^1Oi^{-2}\left(-\xi^{\f23}
									(\tau+s(\tau-t))\right)ds. \end{align*}
								It suffices to notice $|\tau-t|\lesssim \xi^{-\f23}$, and for $|z|\lesssim 1$, $\mathrm{Oi}(-z), \mathrm{Oi}^{-1}(-z)=O(1)$, by  \eqref{behave:Oi}, \eqref{behave:Oi-1} .
							\end{proof}

\section*{Acknowledgments}
The authors would like to thank Joachim Krieger and Wilhelm Schlag for sharing their ideas to improve an early draft of this paper. The authors would also like to thank Lydia Bieri, Yan Guo, Sergiu Klainerman, Peter Miller, Dongyi Wei, Sijue Wu for their helpful conversations and comments in the preparation of this article.

S. Miao is partially supported by the National Key R \& D Program of China 2021YFA1001700 and NSF of China under Grants 12426203, 12221001. Z. Zhang is partially supported by NSF of China under Grant 12288101.

\end{document}